**FLORENTIN SMARANDACHE**
**ION PĂTRAŞCU**

# THE GEOMETRY OF HOMOLOGICAL TRIANGLES

**2012**

# Table of Contents

## PREFACE

This book is addressed to students, professors and researchers of geometry, who will find herein many interesting and original results. The originality of the book *The Geometry of Homological Triangles* consists in using the homology of triangles as a “filter” through which remarkable notions and theorems from the geometry of the triangle are unitarily passed.

Our research is structured in seven chapters, the first four are dedicated to the homology of the triangles while the last ones to their applications.

In the first chapter one proves the theorem of homological triangles (Desargues, 1636), one survey the remarkable pairs of homological triangles, making various connections between their homology centers and axes.

Second chapter boards the theorem relative to the triplets of homological triangles. The Veronese Theorem is proved and it is mentioned a remarkable triplet of homological triangles, and then we go on with the study of other pairs of homological triangles.

Third chapter treats the bihomological and trihomological triangles. One proves herein that two bihomological triangles are trihomological (Rosanes, 1870), and the Theorem of D. Barbilian (1930) related to two equilateral triangles that have the same center.

Any study of the geometry of triangle is almost impossible without making connections with the circle. Therefore, in the fourth chapter one does research about the homological triangles inscribed into a circle. Using the duality principle one herein proves several classical theorems of Pascal, Brianchon, Aubert, Alasia.

The fifth chapter contains proposed problems and open problems about the homological triangles, many of them belonging to the book authors.

The sixth chapter presents topics which permit a better understanding of the book, making it self-contained.

The last chapter contains solutions and hints to the 100 proposed problems from the fifth chapter. The book ends with a list of references helpful to the readers.

*The authors*

# Chapter 1

## Remarkable pairs of homological triangles

In this chapter we will define the homological triangles, we'll prove the homological triangles' theorem and it's reciprocal. We will also emphasize on some important pairs of homological triangles establishing important connections between their centers and axes of homology.

### 1.1. Homological triangles' theorem

**Definition 1**

Two triangles $ABC$ and $A_1B_1C_1$ are called homological if the lines $AA_1$, $BB_1$ and $CC_1$ are concurrent. The concurrence point of lines $AA_1$, $BB_1$, $CC_1$ is called the homological center of $ABC$ and $A_1B_1C_1$ triangles.

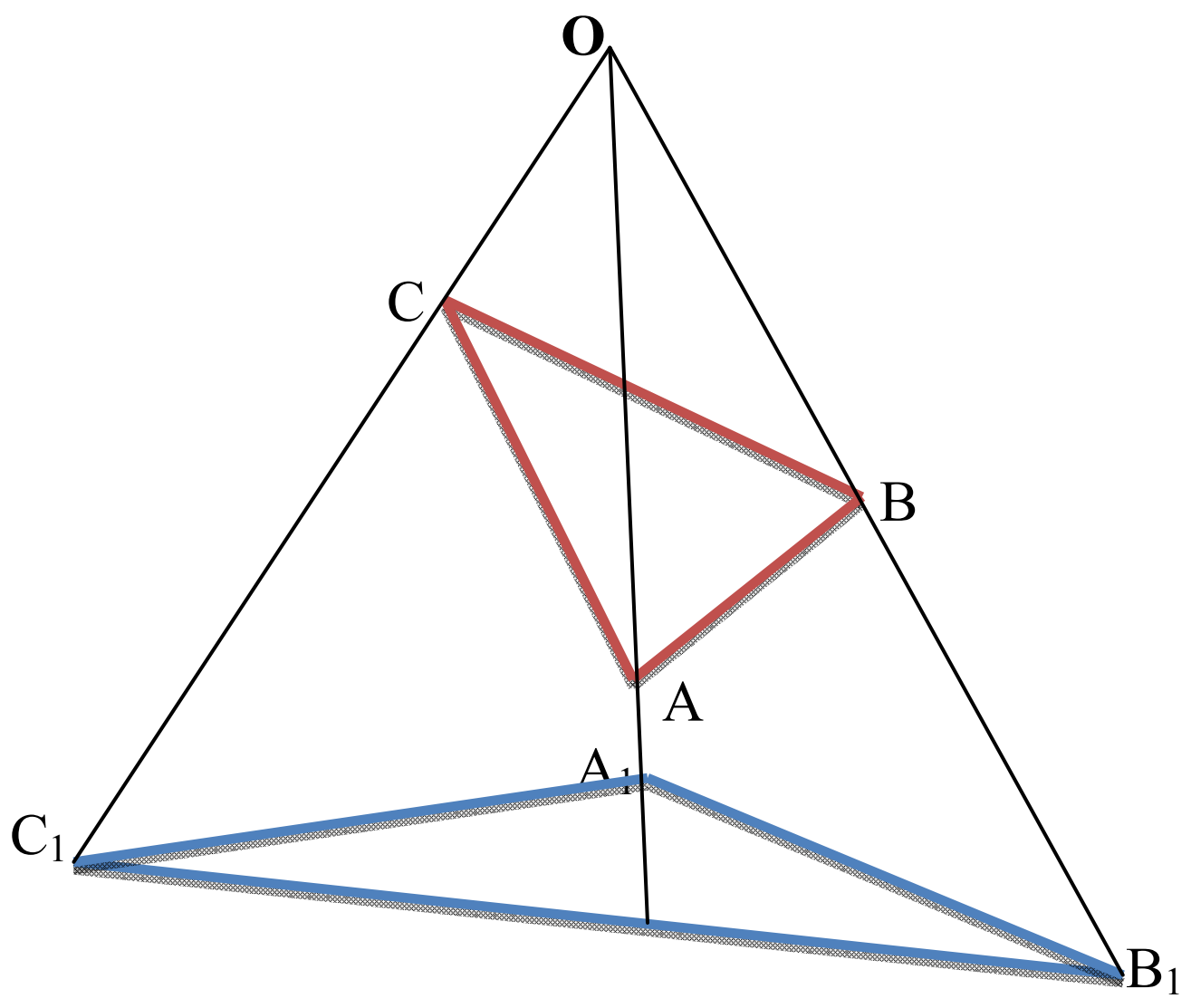

Fig. 1

**Observation 1**

In figure 1 the triangles $ABC$ and $A_1B_1C_1$ are homological. The homology center has been noted with $O$ .

**Definition 2**

If two triangles $ABC$ and $A_1B_1C_1$ are homological and the lines $AB_1$, $BC_1$, $CA_1$ are concurrent, the triangles are called double homological (or bio-homological).

If the lines $AC_1$, $BA_1$, $CB_1$ are also concurrent, we called these triangles triple homological or tri-homological.

**Theorem 1** (G. Desargues – 1636)

1) If $ABC$ and $A_1B_1C_1$ are homological triangles such that:

$$AB \cap A_1B_1 = \{N\},\ BC \cap B_1C_1 = \{M\},\ CA \cap C_1A_1 = \{P\},$$

Then the points $N$, $M$, $P$ are collinear.

2) If $ABC$ and $A_1B_1C_1$ are homological triangles such that:

$$AB \cap A_1B_1 = \{N\},\ BC \cap B_1C_1 = \{M\},\ CA \cap C_1A_1 = \{P\},$$

Then $MN \parallel AC$.

3) If $ABC$ and $A_1B_1C_1$ are homological triangles such that:

$$AB \parallel A_1B_1,\ BC \parallel B_1C_1$$

Then $AC \parallel A_1C_1$.

**Proof**

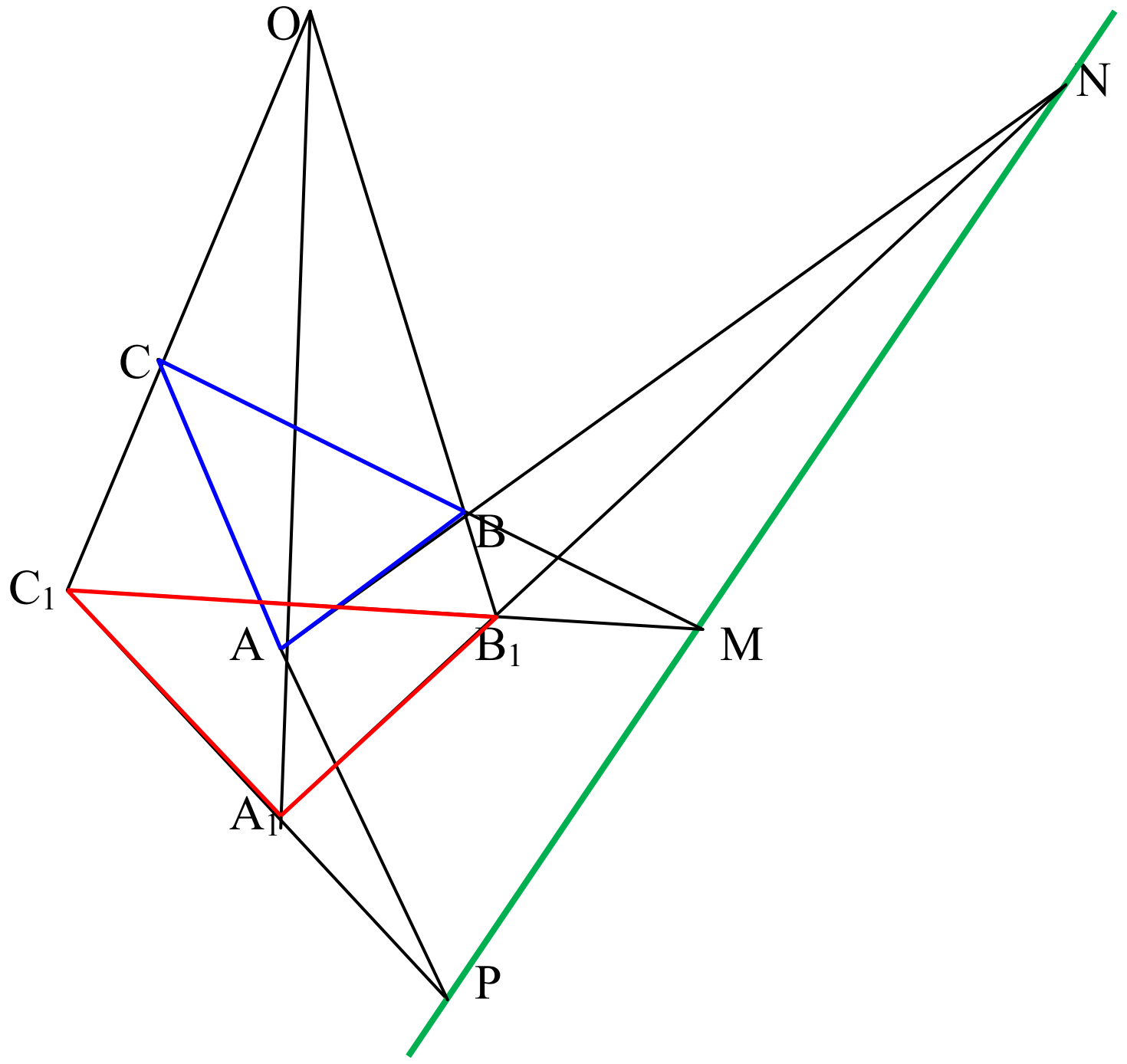

Fig. 2

1) Let $O$ be the homology center of triangles $ABC$ and $A_1B_1C_1$ (see figure 2). We apply the Menelaus' theorem in triangles: $OAC$, $OBC$, $OAB$ for the transversals: $P, A_1, C_1$; $M, B_1, A_1$; $N, B_1, A_1$ respectively. We obtain:

$$\frac{PA}{PC} \cdot \frac{A_1O}{A_1A} \cdot \frac{C_1C}{C_1O} = 1 \qquad (1)$$

$$\frac{MC}{MB} \cdot \frac{B_1B}{B_1O} \cdot \frac{C_1O}{C_1C} = 1 \qquad (2)$$

$$\frac{NB}{NA} \cdot \frac{B_1O}{B_1B} \cdot \frac{A_1A}{A_1O} = 1 \qquad (3)$$

By multiplying side by side the relations (1), (2), and (3) we obtain, after simplification:

$$\frac{PA}{PC}\cdot\frac{MC}{MB}\cdot\frac{NB}{NA}=1 \qquad (4)$$

This relation, in accordance with the Menelaus' reciprocal theorem in triangle $ABC$ implies the collinearity of the points $N$, $M$, $P$.

**Definition 3**

The line determined by the intersection of the pairs of homological lines of two homological triangles is called the triangles' homological axis.

**Observation 2**

The line $M, N, P$ from figure 2 is the homological axis of the homological triangles $ABC$ and $A_1B_1C_1$.

2) Menelaus' theorem applied in triangle $OAB$ for the transversal $N, A_1, B_1$ implies relation (3) $\frac{NB}{NA}\cdot\frac{B_1O}{B_1B}\cdot\frac{A_1A}{A_1O}=1$.

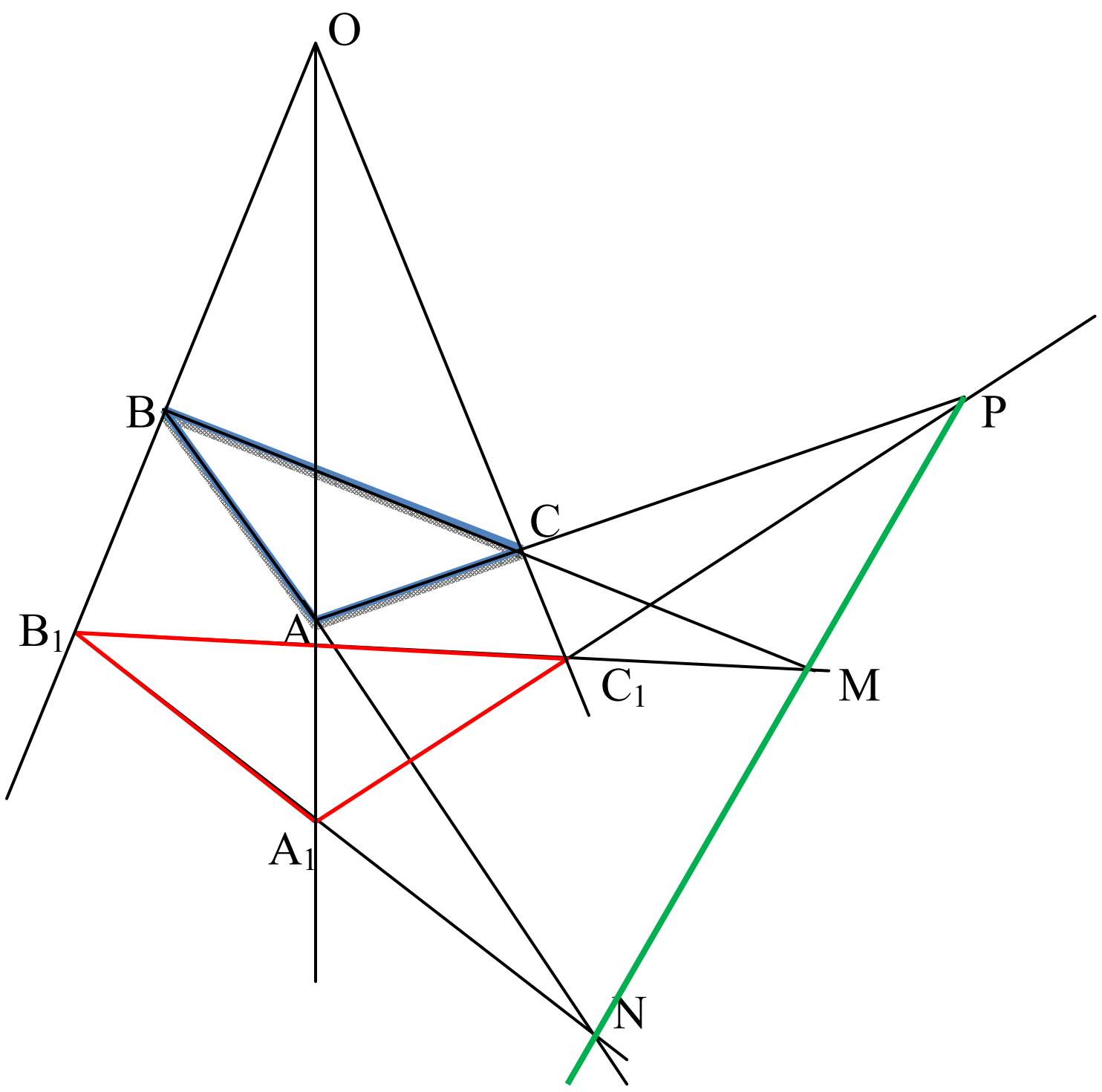

Fig. 3

Menelaus' theorem applied in triangle $OBC$ for the transversal $M, C_1, B_1$ implies relation (2) $\frac{MC}{MB}\cdot\frac{B_1B}{B_1O}\cdot\frac{C_1O}{C_1C}=1$.

By multiplication of side by side of these two relations we find

$$\frac{NB}{NA}\cdot\frac{MC}{MB}\cdot\frac{C_1O}{C_1C}\cdot\frac{A_1A}{A_1O}=1 \qquad (5)$$

Because $AB \parallel A_1C_1$, the Thales' theorem in the triangle $OAC$ gives us that:

$$\frac{C_1O}{C_1C}=\frac{A_1O}{A_1A} \qquad (6)$$

Considering (6), from (5) it results"

$$\frac{NB}{NA}=\frac{MB}{MC}$$

which along with the Thales' reciprocal theorem in triangle $BAC$ gives that $MN \parallel AC$.

**Observation 3.**

The 2) variation of Desargues' theorem tells us that if two triangles are homological and two of their homological lines are parallel and the rest of the pairs of homological sides are concurrent, it results that: the line determined by the intersection of the points of the pairs of homological lines (the homological axis) is parallel with the homological parallel lines.

3) The proof results from Thales' theorem applied in triangles $OAB$ and $OBC$ then by applying the Thales' reciprocal theorem in triangle $OAC$.

**Observation 4.**

The variation 3) of Desargues' theorem is also called the weak form of Desargues' theorem. Two homological triangles which have the homological sides parallel are called homothetic.

**Theorem 2.** (The reciprocal of Desargues' theorem)

1) If two triangles $ABC$ and $A_1B_1C_1$ satisfy the following relations

$$AB\cap A_1B_1=\{N\},\ BC\cap B_1C_1=\{M\},\ CA\cap C_1A_1=\{P\},$$

and the points $N$, $M$, $P$ are collinear then the triangles are collinear.

2) If two triangles $ABC$ and $A_1B_1C_1$ have a pair of parallel lines and the rest of the pairs of lines are concurrent such that the line determined by their concurrency points is parallel with one of the pairs of parallel lines, then the triangles are homological.

3) If two triangles $ABC$ and $A_1B_1C_1$ have

$$AB \parallel A_1B_1,\ AC \parallel A_1C_1, \text{ and } \frac{AB}{A_1B_1}=\frac{AC}{A_1C_1}\neq 1$$

then $BC \parallel B_1C_1$ and the triangles are homological.

Proof.

1) We will be using the method of reduction ad absurdum. Let

$$\{O\}=BB_1\cap AA_1;\ \{O_1\}=AA_1\cap CC_1;\ \{O_2\}=CC_1\cap BB_1\ \ O\neq O_1\neq O_2$$

The Menelaus' theorem applied in triangles $OAB$; $O_1AC$; $O_2BC$ respectively for the transversals $A_1,B_1,N$; $P,A_1,C_1$; $M,B_1,C_1$ respectively provides the following relations:

$$\frac{NB}{NA}\cdot\frac{B_1O}{B_1B}\cdot\frac{A_1A}{A_1O}=1 \quad (7)$$

$$\frac{PA}{PC}\cdot\frac{A_1O_1}{A_1A}\cdot\frac{C_1C}{C_1O_1}=1 \quad (8)$$

$$\frac{MC}{MB}\cdot\frac{B_1B}{B_1O_2}\cdot\frac{C_1O_2}{C_1C}=1 \quad (9)$$

Multiplying side by side these relations and taking into account that the points $M, N, P$ are collinear, that is

$$\frac{PA}{PC}\cdot\frac{MC}{MB}\cdot\frac{NB}{NA}=1 \quad (10)$$

After simplifications we obtain the relation:

$$\frac{A_1O_1}{A_1O}\cdot\frac{B_1O}{B_1O_2}\cdot\frac{C_1O_2}{C_1O_1}=1 \quad (11)$$

The Menelaus' reciprocal theorem applied in triangle $A_1B_1C_1$ and relation (11) shows that the points $O, O_1, O_2$ are collinear. On the other side $O$ and $O_1$ are on $AA_1$, it results that $O_2$ belongs to $AA_1$ also. From

$$\{O_2\}=AA_1\cap BB_1;\ \{O_2\}=AA_1\cap CC_1;\ \{O_2\}=BB_1\cap CC_1,$$

it results that

$$\{O_2\}=AA_1\cap B_1B_2\cap CC_1$$

and therefore $O_2=O_1=O$, which contradicts our assumption.

2) We consider the triangles $ABC$ and $A_1B_1C_1$ such that $AC \parallel A_1C_1$,

$AB\cap A_1B_1=\{N\}$, $BC\cap B_1C_1=\{M\}$ and $MN \parallel AC$.

Let

$$\{O\}=BB_1\cap AA_1;\ \{O_1\}=AA_1\cap CC_1;\ \{O_2\}=CC_1\cap BB_1$$

we suppose that $O\neq O_1\neq O_2\neq 0$. Menelaus' theorem applied in the triangles $OAB$, $O_2BC$ for the transversals $N, A_1, B_1$; $M, C_1, B_1$ respectively leads to the following relations:

$$\frac{NB}{NA}\cdot\frac{B_1O}{B_1B}\cdot\frac{A_1A}{A_1O}=1 \quad (12)$$

$$\frac{MC}{MB}\cdot\frac{B_1B}{B_1O_2}\cdot\frac{C_1O_2}{C_1C}=1 \quad (13)$$

On the other side from $MN \parallel AC$ and $A_1C_1 \parallel MN$ with Thales' theorem, it results

$$\frac{NB}{NA}=\frac{MB}{MC}\cdot \quad (14)$$

$$\frac{A_1A}{A_1O_1}=\frac{C_1C}{C_1O_1} \quad (15)$$

By multiplying side by side the relations (12) and (13) and considering also (14) and (15) we obtain:

$$\frac{A_1O_1}{A_1O}\cdot\frac{B_1O}{B_1O_2}\cdot\frac{C_1O_2}{C_1O_1}=1 \qquad (16)$$

This relation implies the collinearity of the points $O, O_1, O_2$.

Following the same reasoning as in the proof of 1) we will find that $O = O_1 = O_2$, and the theorem is proved.

3) Let $\{O\} = BB_1 \cap AA_1$; $\{O_1\} = AA_1 \cap CC_1$; $\{O_2\} = CC_1 \cap BB_1$, suppose that $O \neq O_1$. Thales' theorem applied in the triangles $OAB$ and $O_1AC$ leads to

$$\frac{OA}{OA_1}=\frac{OB}{OB_1}=\frac{AB}{A_1B_1} \qquad (17)$$

$$\frac{O_1A}{O_1A_1}=\frac{O_1C}{O_1C_1}=\frac{AC}{AC_1} \qquad (18)$$

Because $\frac{AB}{A_1B_1}=\frac{AC}{A_1C_1}\neq 1$ and $A_1, O, O_1$ are collinear, we have that

$$\frac{OA}{OA_1}=\frac{O_1A}{O_1A_1} \qquad (19)$$

This relation shows that $O = O_1$, which is contradictory with the assumption that we made.

If $O = O_1$ then from (17) and (18) we find that

$$\frac{OB}{OB_1}=\frac{OC}{OC_1}$$

which shows that $BC \parallel B_1C_1$ and that $\{O\} = CC_1 \cap BB_1$, therefore the triangles $ABC$ and $A_1B_1C_1$ are homological.

**Observation 5**

The Desargues' theorem is also called the theorem of homological triangles.

**Remark 1**

In the U.S.A. the homological triangles are called perspective triangles. One explanation of this will be presented later.

**Definition 4**

Given a fixed plane $(\alpha)$ and a fixed point $O$ external to the plane $OM$ $(\alpha)$, we name the perspective of a point from space on the plane $(\alpha)$ in rapport to the point $O$, a point $M_1$ of intersection of the line $OM$ with the plane $OM$

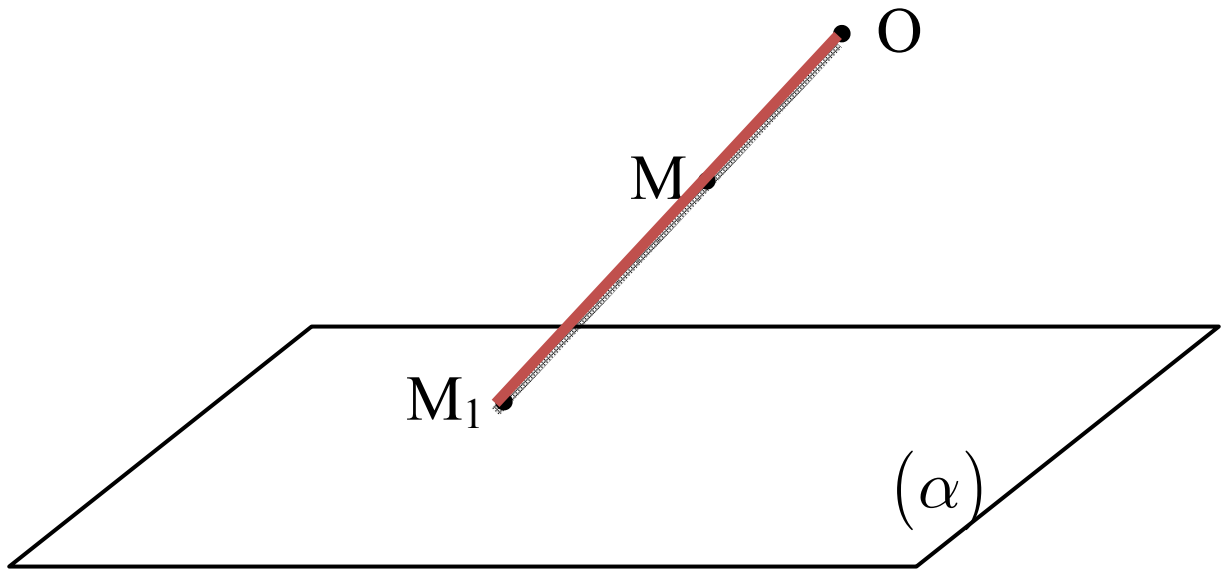

Fig. 4

**Remark 2.**

In the context of these names the theorem of the homological triangles and its reciprocal can be formulated as follows:

If two triangles have a center of perspective then the triangles have a perspective axis.

If two triangles have a perspective axis then the triangles have a perspective center.

**Remark 3.**

Interpreting the plane Desargues' theorem in space or considering that the configuration is obtained through sectioning spatial figures, the proof of the theorem and its reciprocal becomes simple.

We'll illustrate bellow such a proof, precisely we'll prove that if the triangles $ABC$ and $A_1B_1C_1$ are homological and $AB \cap A_1B_1 = \{N\}$, $BC \cap B_1C_1 = \{M\}$, $CA \cap C_1A_1 = \{P\}$, then the points $N$, $M$, $P$ are collinear.

We'll use figure 2 in which $OA_1B_1C_1$ is a triangular pyramidal surface sectioned by the plane $(ABC)$. Because the planes $(ABC)$, $(A_1B_1C_1)$ have a non-null intersection, there is common line and the points $M, N, P$ belong to this line (the homology axis) and the theorem is proved.

In plane: the two triangles from Desargues' theorem can be inscribed one in the other, and in this case we obtain:

**Theorem 3.**

If $A_1, B_1, C_1$ are the Cevians' intersections $AO$, $BO$, $CO$ with the lines $BC$, $CA$, $AB$ and the pairs of lines $AB, A_1B_1$; $BC, B_1C_1$; $CA, C_1A_1$ intersect respectively in the points $N, M, P$, then the points $N, M, P$ are collinear.

**Remark 4.**

The line determined by the points $N, M, P$ (the homology axis of triangles $ABC$ and $A_1B_1C_1$) is called the tri-linear polar of the point $O$ (or the associated harmonic line) in rapport to triangle $ABC$ , and the point $O$ is called the tri-linear pole (or the associated harmonic point) of the line $N, M, P$ .

**Observation 6.**

The above naming is justified by the following definitions and theorems.

**Definition 5.**

Four points $A, B, C, D$ form a harmonic division if

1) The points $A, B, C, D$ are collinear

2) $$\frac{\overline{CA}}{\overline{CB}} = -\frac{\overline{DA}}{\overline{DB}} \qquad (17)$$

If these conditions are simultaneously satisfied we say that the points $C$ and $D$ are harmonic conjugated in rapport to $A, B$ .

Because the relation (17) can be written in the equivalent form

$$\frac{\overline{AC}}{\overline{AD}} = -\frac{\overline{BC}}{\overline{BD}} \qquad (18)$$

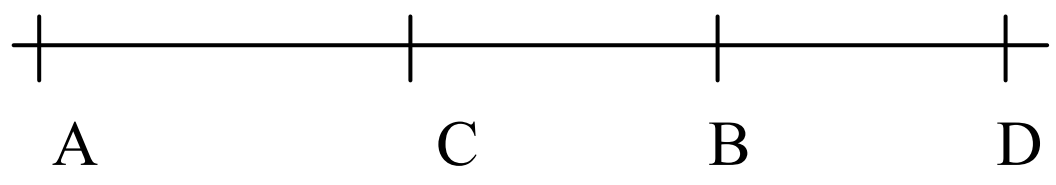

Fig. 5

we can affirm that the points $A, B$ are harmonically conjugated in rapport to $C, D$ . We can state also that $A$ is the harmonic conjugate of $B$ in rapport to $C$ and $D$, or that $D$ is the harmonic conjugate of $C$ in rapport to $A$ and $B$ .

**Observation 7**

It can be proved, through reduction ad absurdum, the unicity of the harmonic conjugate of a point in rapport with other two given points.

**Note 1.**

The harmonic deviation was known in Pythagoras' school ($5^{th}$ century B.C.)

**Theorem 4.**

If $AA_1$, $BB_1$, $CC_1$ are three concurrent Cevians in triangle $ABC$ and $N, M, P$ are the harmonic conjugate of the points $C_1$, $A_1$, $B_1$ in rapport respectively with $A, B$; $B, C$; $C, A$, then the points $N, M, P$ are collinear.

**Proof**

We note $\{M'\} = C_1B_1 \cap CB$

Ceva's theorem gives us:

$$\frac{\overline{A_1B}}{\overline{A_1C}} \cdot \frac{\overline{B_1C}}{\overline{B_1A}} \cdot \frac{\overline{C_1A}}{\overline{C_1B}} = -1 \qquad (19)$$

Menelaus' theorem applied in triangle $ABC$ for transversal $M, B_1, C_1$ leads to:

$$\frac{\overline{M'B}}{\overline{M'C}} \cdot \frac{\overline{B_1C}}{\overline{B_1A}} \cdot \frac{\overline{C_1A}}{\overline{C_1B}} = 1 \quad (20)$$

From relations (19) and (20) we obtain:

$$\frac{\overline{A_1B}}{\overline{A_1C}} = -\frac{\overline{M'B}}{\overline{M'C}} \quad (21)$$

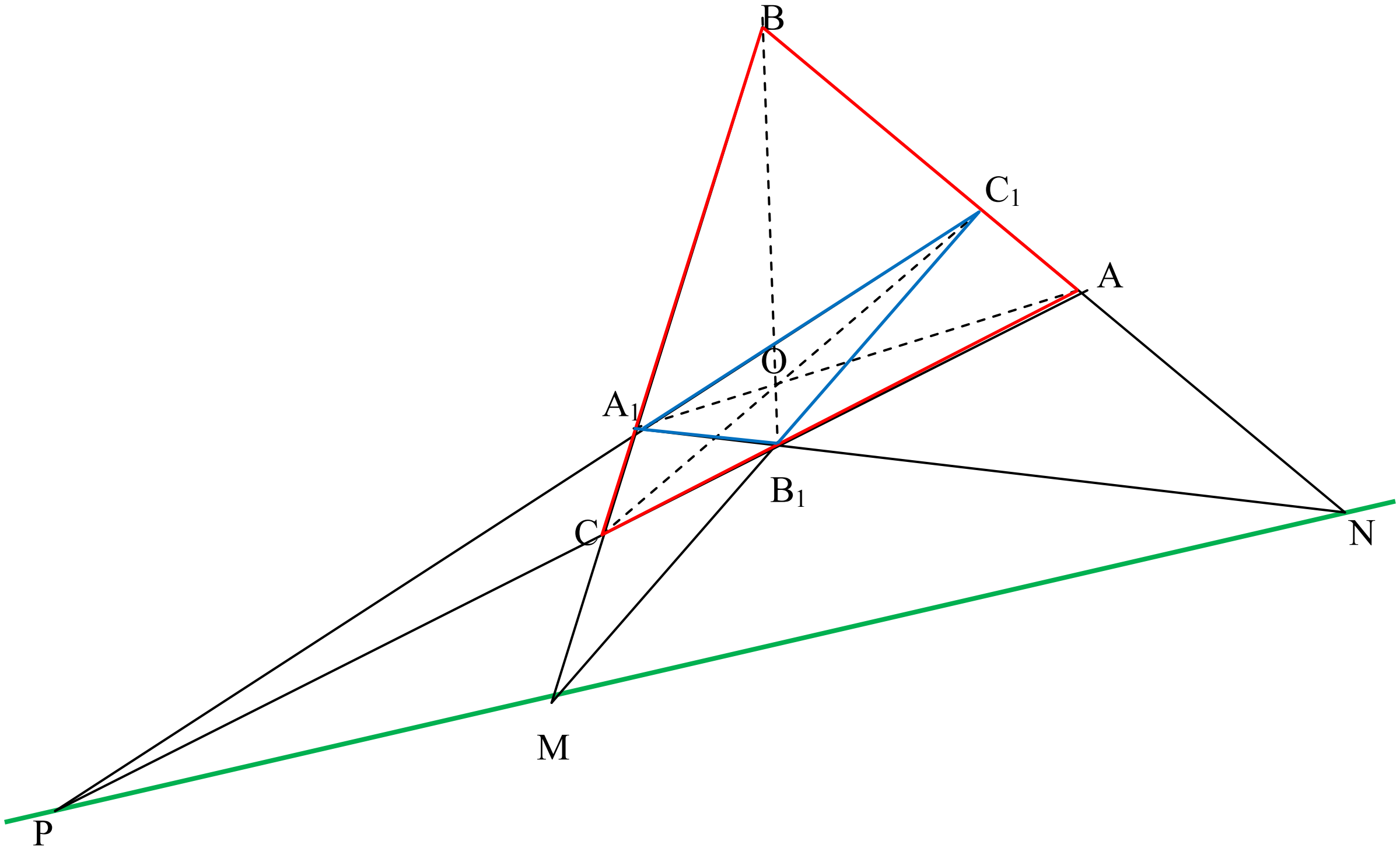

Fig. 6

This relation shows that $M'$ is the harmonic conjugate of point $A_1$ on rapport to $B$ and $C$. But also $M$ is the harmonic conjugate of point $A_1$ on rapport to $B$ and $C$. From the unicity property of a harmonic conjugate of a point, it results that $M = M'$. In a similar way we can prove that $N$ is the intersection of the lines $A_1B_1$ and $AB$ and $P$ is the intersection of lines $A_1C_1$ and $AC$. Theorem 3is, in fact, a particularization of Desargues' theorem, shows that the points $N, M, P$ are collinear.

**Remark 5.**

The theorems 3 and 4 show that we can construct the harmonic conjugate of a given point in rapport with other two given points only with the help of a ruler. If we have to construct the conjugate of a point $A_1$ in rapport to the given points $B, C$ we can construct a similar configuration similar to that in figure 6, and point $M$, the conjugate of $A_1$ it will be the intersection of the lines $BC$ and $B_1C_1$

## 1.2. Some remarkable homological triangles

In this paragraph we will visit several important pairs of homological triangles with emphasis on their homological center and homological axis.

### A. The orthic triangle

**Definition 6**

The orthic triangle of a given triangle is the triangle determined by the given triangle's altitudes' feet.

Theorem 3 leads us to the following

**Proposition 1**

A given triangle and its orthic triangle are homological triangles.

**Observation 8**

In some works the Cevian triangle of a point $O$ from the triangle's $ABC$ plane is defined as being the triangle $A_1B_1C_1$ determined by the Cevians' intersection $AO$, $BO$, $CO$ with $BC$, $CA$, $AB$ respectively. In this context the orthic triangle is the Cevian triangle of the orthocenter of a given triangle.

**Definition 7**

The orthic axis of a triangle is defined as being the orthological axis of that triangle and its orthic triangle.

**Observation 9**

The orthological center of a triangle and of its orthic triangle is the triangle's orthocenter.

**Definition 8**

We call two lines $c,d$ anti-parallel in rapport to lines $a,b$ if $\sphericalangle(a,b)=\sphericalangle(d,a)$.

In figure 7 we represented the anti-parallel lines $c,d$ in rapport to $a,b$.

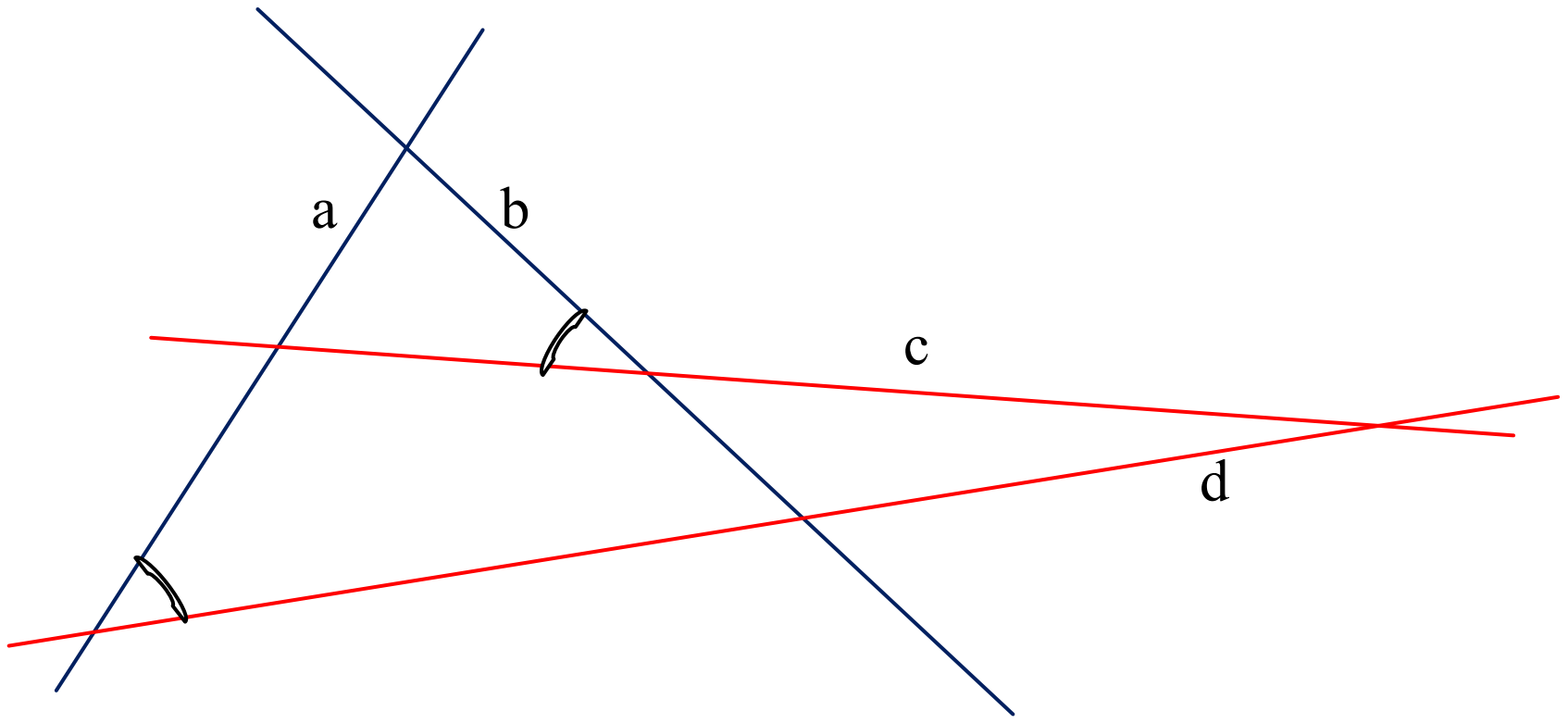

Fig. 7

**Observation 10**

The lines $c, d$ are anti-parallel in rapport to $a, b$ if the quadrilateral formed by these lines is inscribable, If $c, d$ are anti-parallel with the concurrent lines $a, b$ then the lines $a, b$ are also anti-parallel in rapport to $c, d$.

**Proposition 2**

The orthic triangle of a given triangle has the anti-parallel sides with the sides of the given triangle.

**Proof**

In figure 8, $A_1B_1C_1$ is the orthic triangle of the orthological triangle $ABC$.

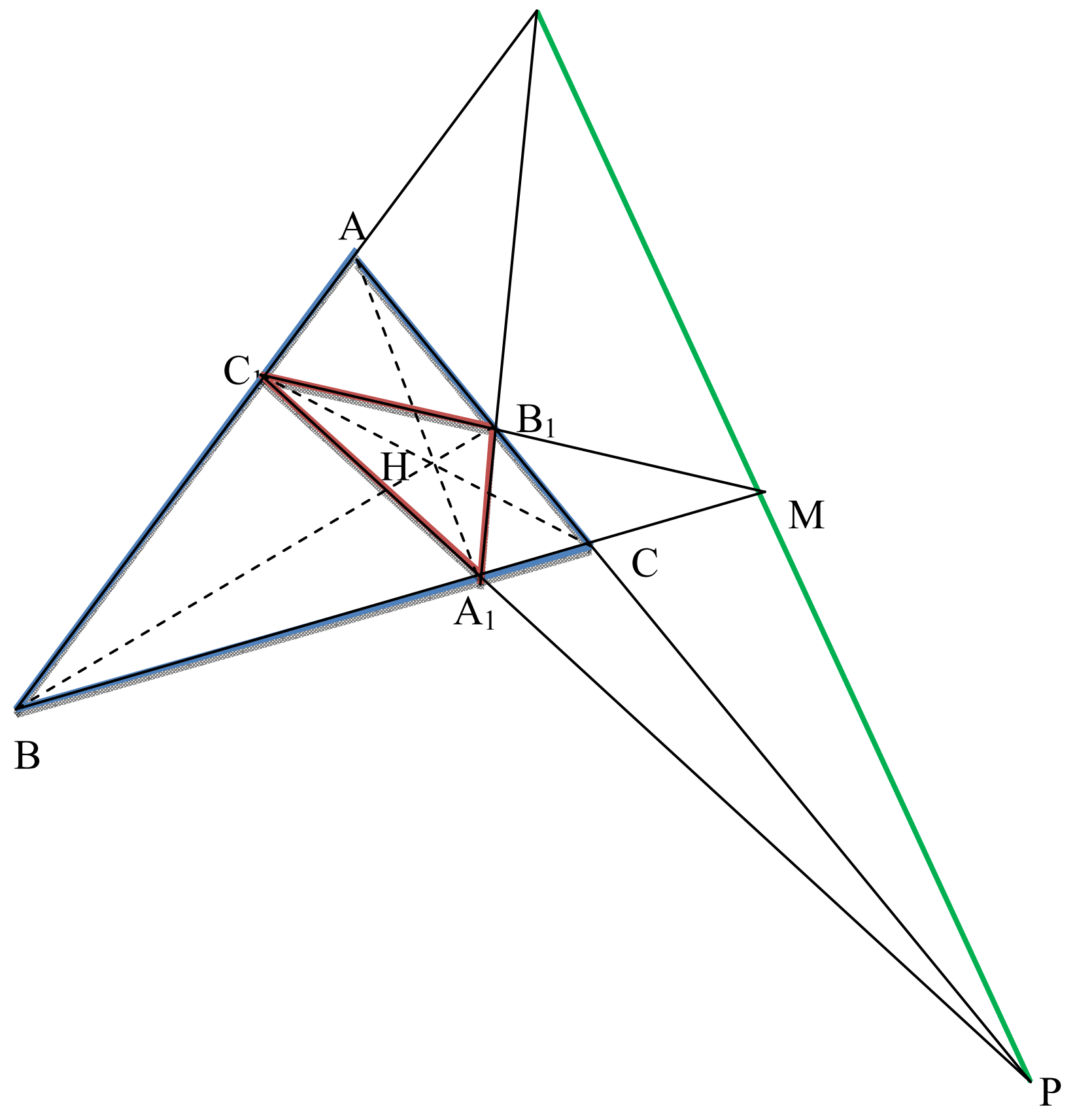

Fig 8

We prove that $B_1C_1$ is anti-parallel to $BC$ in rapport to $AB$ and $AC$. Indeed, the quadrilateral $HB_1AC_1$ is inscribable ($m\sphericalangle(HB_1A)+m\sphericalangle(HC_1A)=90°+90°=180°$), therefore $\sphericalangle AB_1C_1 \equiv \sphericalangle AHC_1$. On the other side $\sphericalangle AHC_1 \equiv \sphericalangle ABC$ (as angles with the sides perpendicular respectively), consequently $\sphericalangle AB_1C_1 \equiv \sphericalangle ABC$. Analogue we can prove that the rest of the sides are anti-parallel.

**B. The Cevian triangle**

**Proposition 3**

The Cevian triangle of the center of a inscribe circle of a triangle and the triangle are homological. The homology axis contains the exterior bisectors' legs of the triangle.

The proof of this proposition results from the theorems 3 and 4. Indeed, if $I$ is the center of the inscribed circle then also $AI$ intersects $BC$ in $A_1$, we have:

$$\frac{A_1B}{A_1C}=\frac{AB}{AC} \tag{22}$$

(the interior bisectors' theorem). If $B_1$ and $C_1$ are the feet of the bisectors' $BI$ and $CI$ and $BC \cap B_1C_1 = \{M\}$ then

$$\frac{MB}{MC}=\frac{AB}{AC} \qquad BC \text{ in} \tag{23}$$

if the exterior bisector of angle $A$ intersects $BC$ in $M'$ then from the theorem of the external bisector we have

$$\frac{M'B}{M'C}=\frac{AB}{AC} \tag{24}$$

From the relations (23) and (24) we find that $M'=M$, therefore the leg of the external bisector of angle $A$ belongs to the homological axis; similarly it can be proved the property for the legs of the external bisectors constructed from $B$ and $C$.

**Proposition 4**

In a triangle the external bisectors of two angles and the internal bisector of the third triangle are concurrent.

**Proof.**

We note by $I_a$ the point of intersection of the external bisector from $B$ and $C$, and with $D$, $E$, $F$ the projections of this point on $BC$, $AB$, $AC$ respectively.

Because $I_a$ belongs to the bisector of angle $B$, we have:

$$I_aD = I_aE \tag{25}$$

Because $I_a$ belongs to the bisector of angle $C$, we have:

$$I_aD = I_aF \tag{26}$$

From (25) and (26) it results that:

$$I_aE = I_aF \tag{27}$$

This relation shows that $I_aE$ belongs to the interior bisector of angle $A$

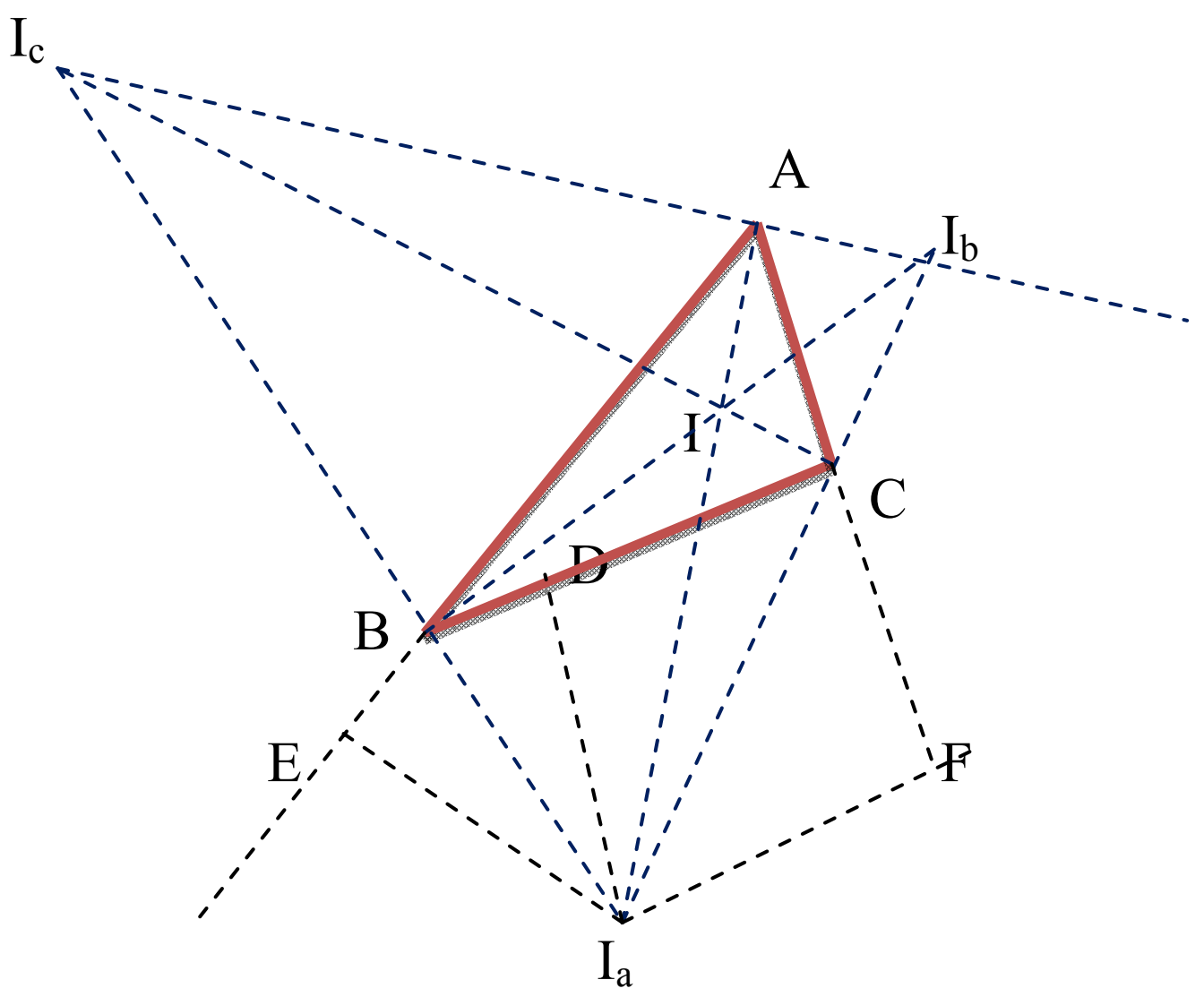

Fig. 9

**Remark 6**

The point $I_a$ is the center of a circle tangent to the side $BC$ and to the extensions of the sides $AB$, $AC$. This circle is called the ex-inscribed circle of the triangle. For a triangle we have three ex-inscribed circles.

**C. The anti-supplemental triangle**

**Definition 9**

The triangle determined by the external bisectors of a given triangle is called the anti-supplemental triangle of the given triangle.

**Observation 11**

The anti-supplemental triangle of triangle $ABC$ is determined the centers of the ex-inscribed circles of the triangle, that is the triangle $I_aI_bI_c$

**Proposition 5**

A given triangle and its anti-supplemental triangle are homological. The homology center is the center of the circle inscribed in triangle and the homological axis is the tri-linear polar of the inscribed circle's center.

**Proof.**

The proof of this property results from the propositions 3 and 3.

**Remark 7**

We observe, without difficulty that for the anti-supplemental triangle $I_aI_bI_c$. The triangle $ABC$ is an orthic triangle (the orthocenter of $I_aI_bI_c$ is $I$ ), therefore the homological axis of triangles $ABC$ and $I_aI_bI_c$ is the orthic axis of triangle $I_aI_bI_c$ that is the line determined by the external bisectors' feet of triangle $ABC$.

### D. The K-symmedian triangle

**Definition 10**

In a triangle the symmetrical of the median in rapport to the interior bisector of the same vertex is called the symmedian.

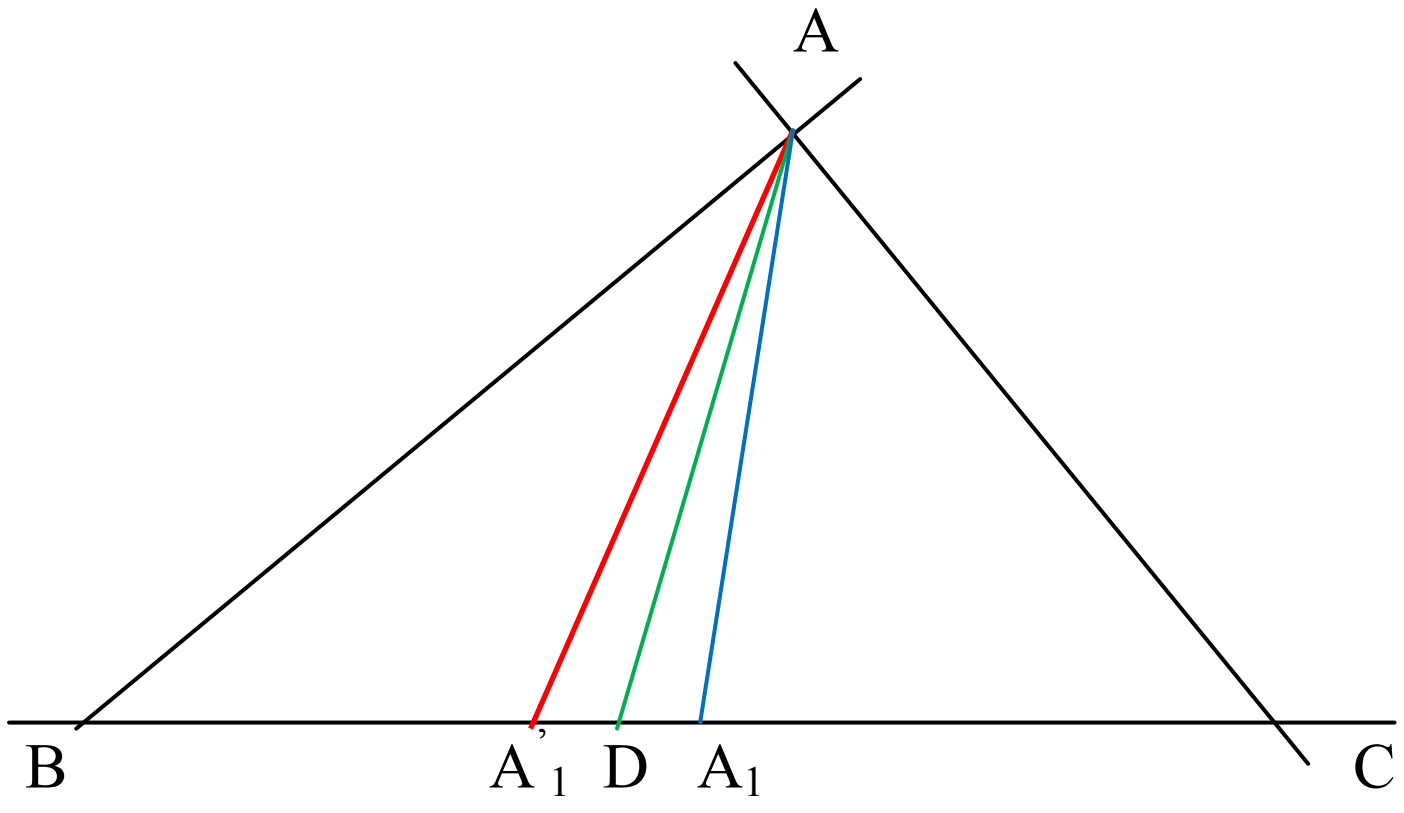

Fig. 10

**Remark 8**

In figure 10, $AA_1$ is the median, $AD$ is the bisector and $A_1'$ is the symmedian.

We have the relation:

$$\sphericalangle BAA_1' \equiv \sphericalangle CAA_1 \qquad (28)$$

Two Cevians in a triangle which satisfy the condition (28) are called isogonal Cevians. Therefore the symmedian is isogonal to the median.

**Theorem 5**

If in a triangle $ABC$ the Cevians $AA_1$ and $AA_1'$ are isogonal, then:

$$\frac{BA_1}{CA_1}\cdot\frac{BB_1'}{CA_1'}=\left(\frac{AB}{AC}\right)^2 \quad \text{(Steiner relation)} \qquad (29)$$

**Proof**

We have

$$\frac{BA_1}{CA_1}=\frac{aria\Delta BAA_1}{aria\Delta CAA_1} \qquad (30)$$

$$aria\vartriangle BAA_1 \frac{1}{2} AB\cdot AA_1\cdot\sin\left(\sphericalangle BAA_1\right)$$

$$aria\vartriangle CAA_1 \frac{1}{2} AC\cdot AA_1\cdot\sin\left(\sphericalangle CAA_1\right)$$

Therefore

$$\frac{BA_1}{CA_1}=\frac{aria\triangle BAA_1}{aria\triangle CAA_1}=\frac{AB\cdot\sin(\sphericalangle BAA_1)}{AC\cdot\sin(\sphericalangle CAA_1)} \quad (31)$$

Similarly

$$\frac{BA_1'}{CA_1'}=\frac{aria\triangle BAA_1'}{aria\triangle CAA_1'}=\frac{AB\cdot\sin(\sphericalangle BAA_1')}{AC\cdot\sin(\sphericalangle CAA_1')} \quad (32)$$

Because

$$\sphericalangle BAA_1 \equiv \sphericalangle CAA_1'$$

$$\sphericalangle CAA_1 \equiv \sphericalangle BAA_1'$$

from the relations (31) and (32) it results the Steiner relation.

**Remark 9**

The reciprocal of theorem 5 is true.

If $AA_1'$ is symmedian in triangle $ABC$, then

$$\frac{BA_1'}{CA_1'}=\left(\frac{AB}{AC}\right)^2 \quad (33)$$

Reciprocal, if $AA_1'\in(BC)$ and relation (33) is true, then $AA_1'$ is symmedian.

**Theorem 6**

The isogonal of the concurrent Cevians in a triangle are concurrent Cevians.

**Proof**.

The proof of this theorem results from the theorem 5 and the Ceva's reciprocal theorem.

**Definition 11**

We call the concurrence point of the Cevians in a triangle and the concurrence point of the their isogonals the isogonal conjugate points.

**Remark 10.**

We can show without difficulty that in a triangle the orthocenter and the center of its circumscribed circle are isogonal conjugated points.

**Definition 12**

The Lemoine's point of a triangle is the intersection of its three symmedians.

**Observation 12**

The gravity center of a triangle and its symmedian center are isogonal conjugated points.

**Propositions 6**

If in a triangle $ABC$ the points $D, E$ belong to the sided $AC$ and $AB$ respectively such that $DE$ and $BC$ to be anti-parallel, and the point $S$ is in the middle of the anti-parallel $(DE)$ then $AS$ is the symmedian in the triangle $ABC$.

**Proof**

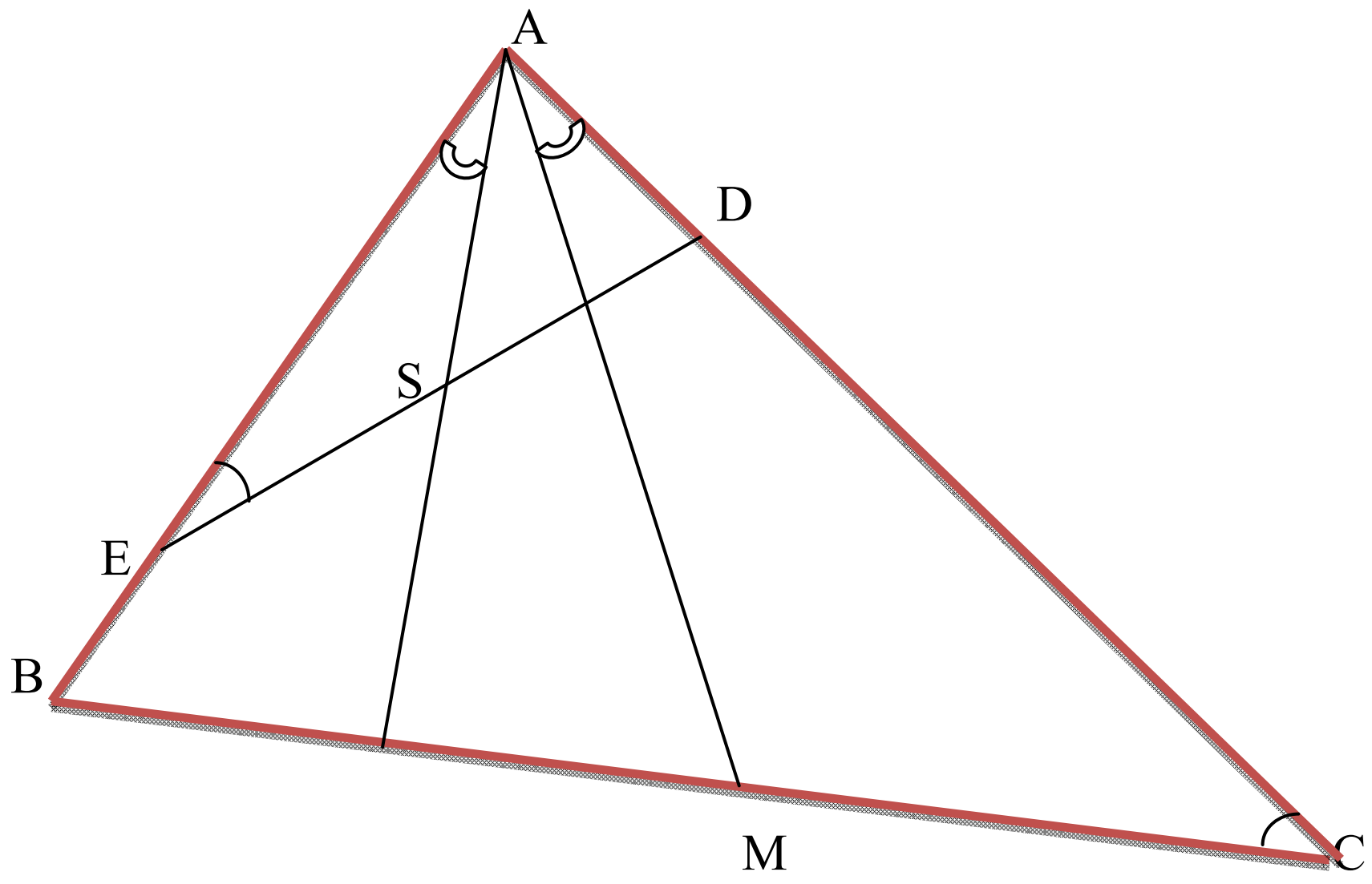

Fig. 11

Let $M$ the middle of the side $(BC)$ (see figure 11). We have $\sphericalangle AED \equiv \sphericalangle ACB$. The triangles $AED$ and $ACB$ are similar. It results that:

$$\frac{AE}{AC} = \frac{ED}{CB} \tag{34}$$

We have also:

$$\frac{ES}{CM} = \frac{ED}{CB} \tag{35}$$

therefore

$$\frac{AE}{AC} = \frac{ES}{CM} \tag{36}$$

Relation (36) along with the fact that $\sphericalangle AES \equiv \sphericalangle ACM$ leads to $\Delta AES \sim \Delta ACM$ with the consequence $\sphericalangle EAS \equiv \sphericalangle ACAM$ which shows that $AS$ is the isogonal of the median $AM$, therefore $AS$ is symmedian.

**Remark 11**

We can prove also the reciprocal of proposition 6 and then we can state that the symmedian of a triangle is the geometrical locus of the centers of the anti-parallels to the opposite side.

**Definition 3**

If in a triangle $ABC$ we note $A_1'$ the leg of the symmedian from $A$ and $A_1''$ is the harmonic conjugate of $A_1'$ in rapport to $B$, $C$, we say that the Cevian $AA_1''$ is the exterior symmedian of the triangle.

**Proposition 7**

The external symmedians of a triangle are the tangents constructed in the triangle's vertices to the triangle's circumscribed circle.

**Proof**

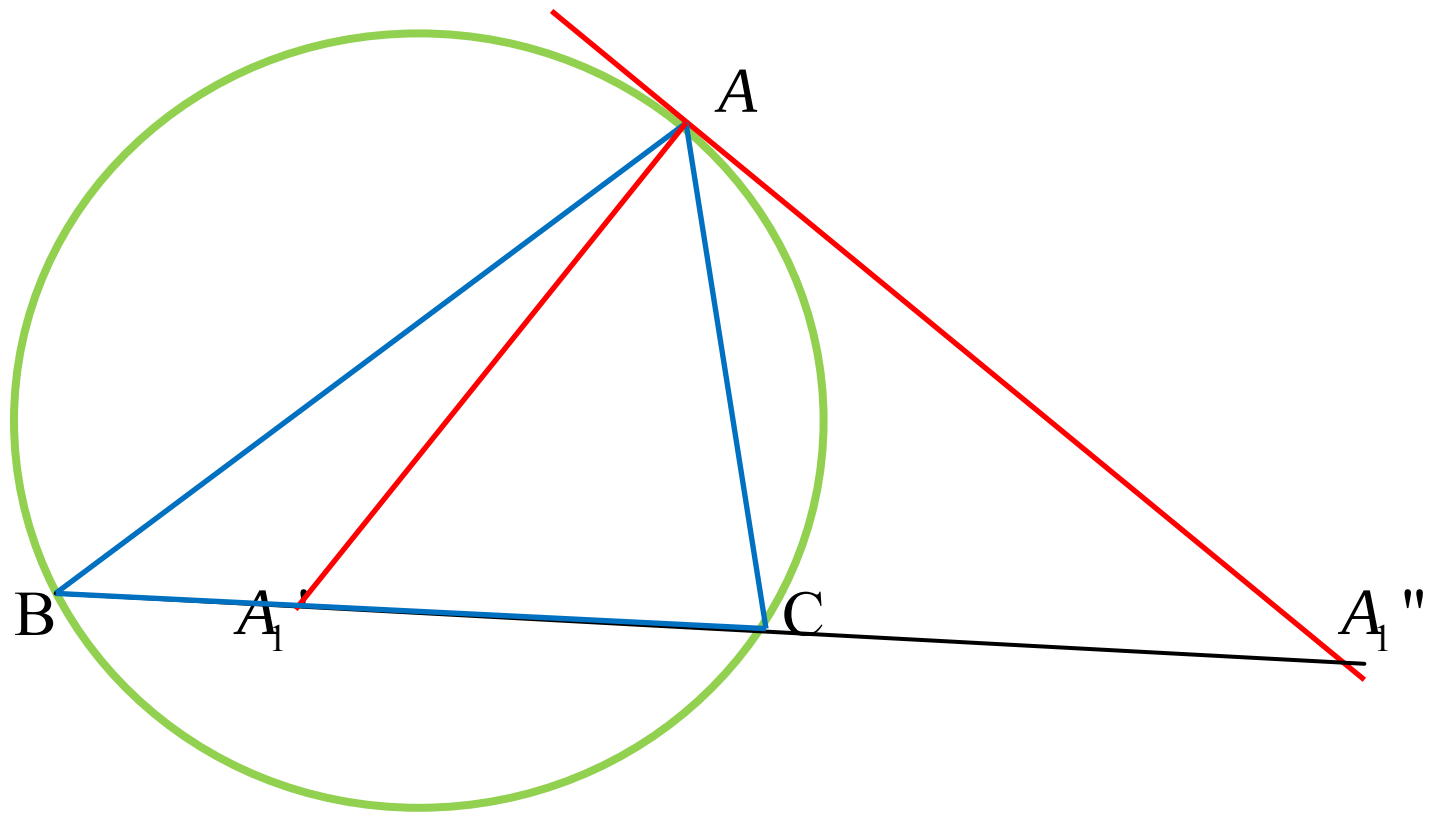

Fig. 12

The triangles $BAA_1"$ and $CAA_1"$ are similar because $\sphericalangle ABA_1" \equiv \sphericalangle CAA_1"$ and $\sphericalangle AA_1"C$ is common (see figure 12).

It results:

$$\frac{A_1"B}{A_1"C} = \frac{AB}{AC} = \frac{A_1"A}{A_1"C} \tag{37}$$

From relation (37) we find that

$$\frac{A_1"B}{A_1"C} = \left(\frac{AB}{AC}\right)^2 \tag{38}$$

This relation along with the relation (33) show that $A_1"$ is the harmonic conjugate of $A_1'$ in rapport to $B, C$, therefore $AA_1"$ is an external symmedian.

**Theorem 7 (Carnot – 1803)**

The tangents constructed on the vertices points of a non-isosceles triangle to its circumscribed circle intersect the opposite sides of the triangle in three collinear points.

The proof of this theorem results as a particular case of theorem 4 or it can be done using the anterior proved results.

**Definition 14**

The line determined by the legs of the exterior medians of a non-isosceles triangle is called the Lemoine's line of the triangle; it is the tri-linear polar of the symmedian center.

From the results anterior obtained , it results the following

**Proposition 8**

The Cevian triangles of the symmedian center of a given non-isosceles triangle are homological. The homology axis is the Lemoine's line of the given triangle.

**Proposition 9**

In a triangle the external symmedians of two vertices and the symmedian of the third vertex are concurrent.

**Proof**

Let $S$ be the intersection of the external symmedians constructed through the vertices $B$ and $C$ of triangle $ABC$ (see figure 13).

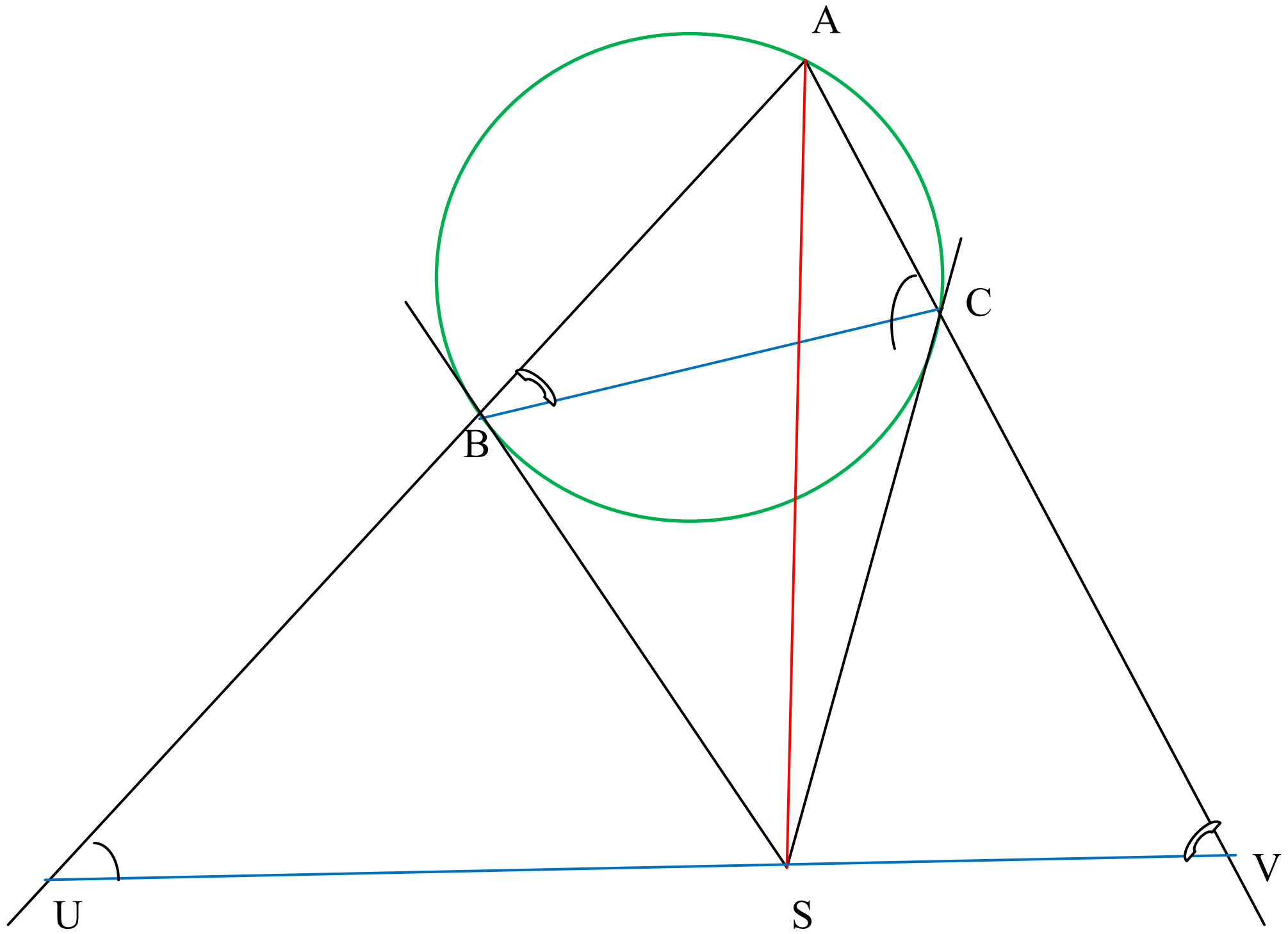

Fig. 13

We construct through $S$ the anti-polar $UV$ to $BC$ ($U \in AB$, $V \in AC$). We have

$$\sphericalangle AUV \equiv \sphericalangle C \text{ and } \sphericalangle AVU \equiv \sphericalangle B$$

On the other side because $BS$, $CS$ are tangent to the circumscribed circle we have

$$SB = SC \text{ and } \sphericalangle CBS \equiv \sphericalangle BCS \equiv \sphericalangle A$$

It results that:

$$\sphericalangle UBS \equiv \sphericalangle C \text{ and } \sphericalangle UCS \equiv \sphericalangle B$$

Consequently, the triangles $SBU$ and $SCV$ are isosceles $SB = SU$; $SC = SV$.

We obtain that $SU = SV$, which based on proposition 6 proves that $AS$ is symmedian.

## E. The tangential triangle

**Definition 15.**

The tangential triangle of a triangle $ABC$ is the triangle formed by the tangents constructed on the vertices $A, B, C$ to the circumscribed circle of the given triangle.

**Observation 13**

In figure 14 we note $T_aT_bT_c$ the tangential triangle of $ABC$. The center of the circumscribed circle of triangle $ABC$ is the center of the inscribed circle in the tangential circle.

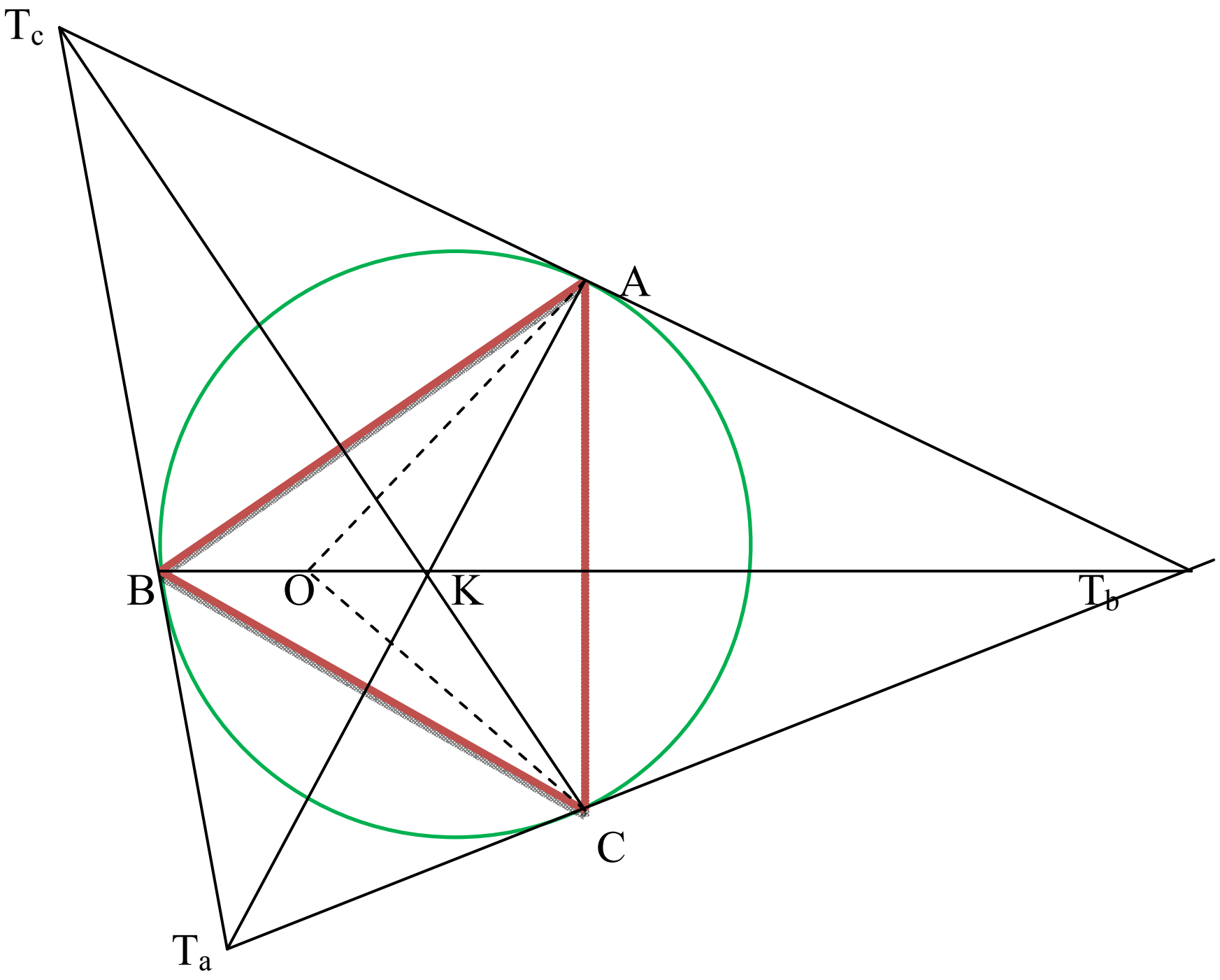

Fug. 14

**Proposition 10**

A non-isosceles given triangle and its tangential triangle are homological. The homological center is the symmedian center of the triangle, and the homology axis is the triangles' Lemoine's line.

**Proof**

We apply proposition 9 it results that the lines $AT_a, BT_b, CT_c$ are symmedians in the triangle $ABC$, therefore are concurrent in the symmedian center $K$, and therefore triangles $ABC$ and $T_aT_bT_c$ are homological. On the other side $AT_b, BT_c, CT_a$ are external symmedians of the triangle $ABC$ and then we apply proposition 8.

**F. The contact triangle**

**Definition 16**

The contact triangle of a given triangle is the triangle formed by the tangential vertices of the inscribed circle in triangle with its sides.

**Observation 14**

In figure 15 the contact triangle of the triangle $ABC$ is noted $C_aC_bC_c$.

**Definition 17**

A pedal triangle of a point from a triangle's plane the triangle determined by the orthogonal projections of the point on the triangle's sides.

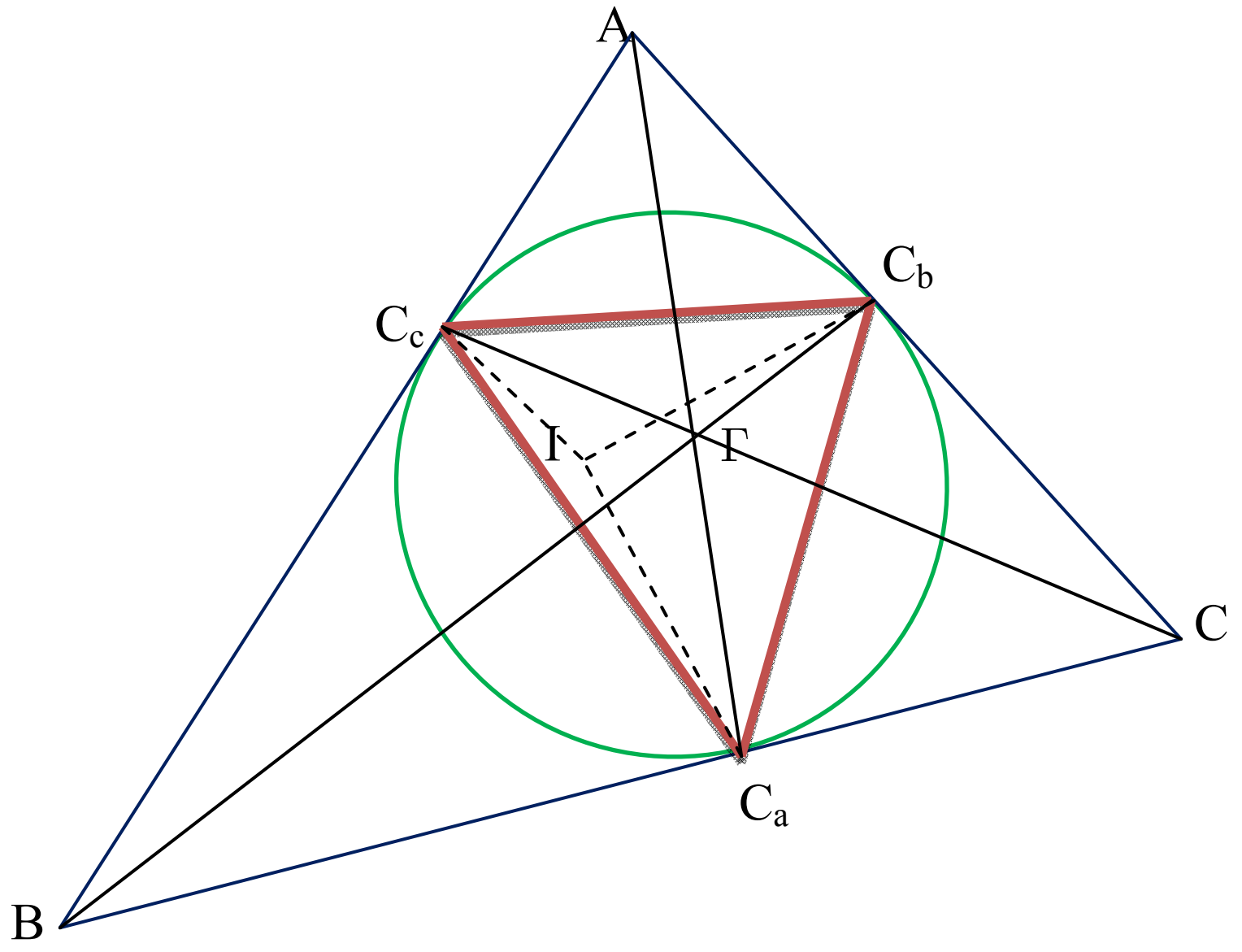

Fig 15

**Observation 15**

The contact triangle of the triangle $ABC$ is the pedal triangle of the center $I$ of the inscribed circle in triangle.

**Proposition 11**

In a triangle the Cevians determined by the vertices and the contact points of the inscribed circle with the sides are concurrent.

The proof results without difficulties from the Ceva's reciprocal theorem and from the fact that the tangents constructed from the triangle's vertices to the circumscribed circle are concurrent.

**Definition 18**

The concurrence point of the Cevians of the contact points with the sides of the inscribed circle in a triangle is called the Gergonne's 's point.

**Observation 16**

In figure 15 we noted the Gergonne's point with $\Gamma$

**Proposition 12**

A non-isosceles triangle and its contact triangle are homological triangles. The homology center is the Gergonne's point, and the homology axis is the Lemoine's line of the contact triangle.

**Proof**

From proposition 11 it results that the Gergonne's point is the homology center of the triangles $ABC$ and $C_aC_bC_c$. The homology axis of these triangles contains the intersections of the opposite sides of the given triangle and of the contact triangle, because, for example, $BC$ is tangent to the inscribed circle, it is external symmedian in triangle $C_aC_bC_c$ and therefore intersects the $C_bC_c$ in a point that belongs to the homology axis of these triangles, that is to the Lemoine's line of the contact triangle.

**Proposition 13**

The contact triangle $C_aC_bC_c$ of triangle $ABC$ and triangle $A_1B_1C_1$ formed by the projections of the centers of the ex-inscribed circles $I_a, I_b, I_c$ on the perpendicular bisectors of $BC$, $CA$ respectively $AB$ are homological. The homology center is the Gergonne's point $\Gamma$ of triangle $ABC$.

**Proof**

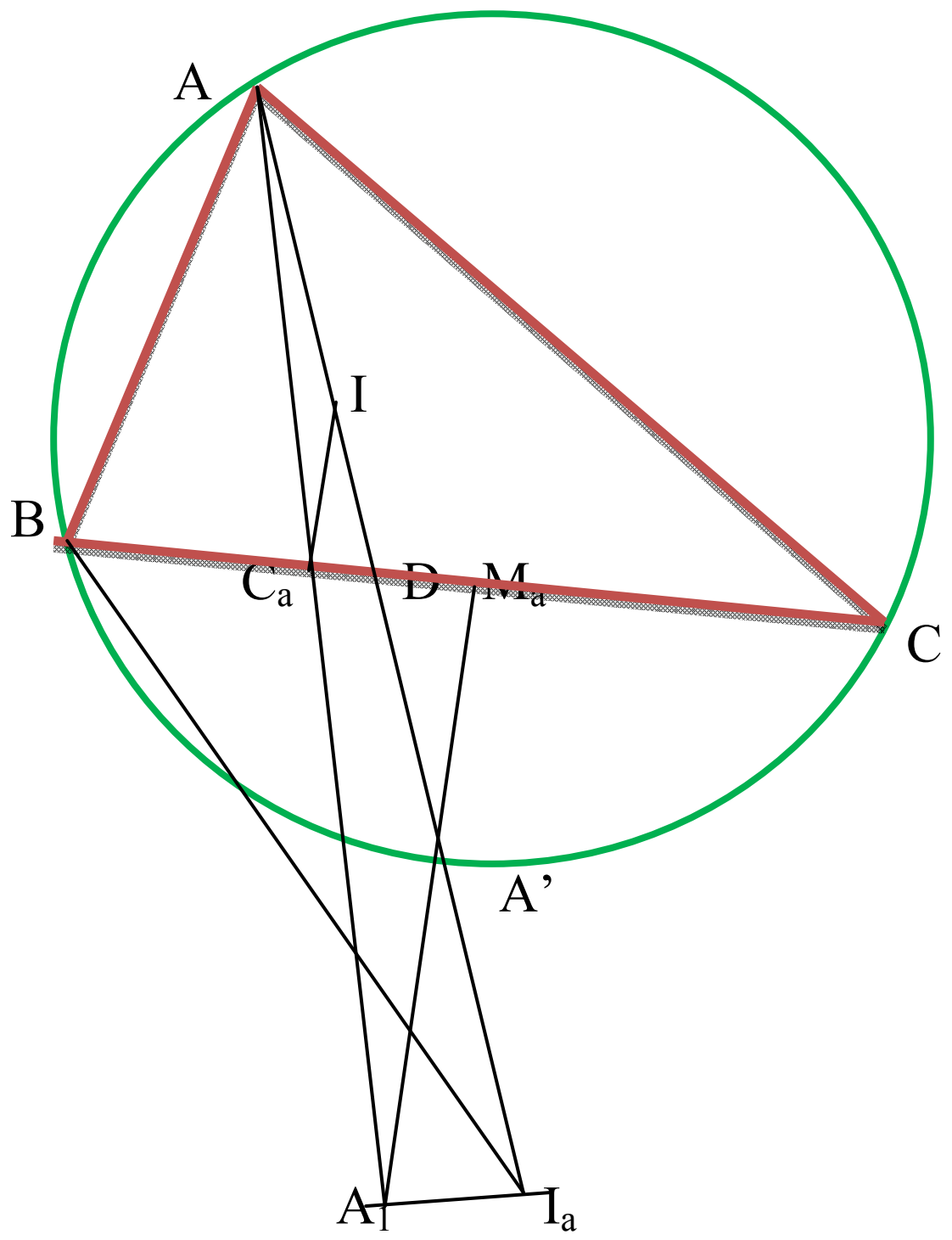

Fig. 16

Let $D$ and $A'$ the intersection of the bisectrix $AI$ with $BC$ and with the circumscribed circle to the triangle $ABC$ and $M_a$ the middle point of $BC$ (see the figure above).

Considering the power of $I_a$ in rapport to the circumscribed circle of triangle $ABC$, we have $I_aA' \cdot I_aA = OI_a^2 - R^2$. It is known that $OI_a^2 = R^2 + 2Rr_a$, where $r_a$ is the A-ex-inscribed circle's radius .

Then

$$I_aA' \cdot I_aA = 2Rr_a \tag{39}$$

The power of $D$ in rapport to the circumscribed circle of triangle $ABC$ leads to:

$$DB \cdot DC = AD \cdot DA'$$

Therefore

$$AD \cdot DA' = \frac{a^2bc}{(b+c)^2} \quad (40)$$

From (39) and (40) we obtain:

$$\frac{I_aA'}{DA'} \cdot \frac{I_aA}{DA} = \frac{2Rr_a(b+c)^2}{a^2bc} \quad (41)$$

Using the external bisectrix' theorem we obtain

$$\frac{I_aA}{DA} = \frac{b+c}{a}$$

and from here

$$\frac{I_aA}{DA} = \frac{b+c}{b+c-a}$$

Substituting in (41) this relation and taking into account that $abc = 4RS$ and $S = r_a(p-a)$ we find that

$$\frac{I_aA'}{DA'} = \frac{b+c}{a} \quad (42)$$

From $DM_a = BM - BD$, we find:

$$DM_a = \frac{a(b-c)}{2(b+c)} \quad (43)$$

The similarity of triangles $A'A_1I_a$ and $A'M_aD$ leads to

$$I_aA' = \frac{b-c}{2}.$$

We note $AA_1 \cap BC = C_a^{'}$, we have that the triangles $AC_a^{'}D$ and $AA_1I_a$ are similar. From:

$$\frac{I_aA}{DA} = \frac{I_aA_1}{DC_a^{'}}$$

we find

$$DC_a^{'} = \frac{(b-c)(p-a)}{(b+c)}$$

Computing $BC_a^{'} = BD - DC_a^{'}$ we find that $BC_a^{'} = p-b$, but we saw that $BC_a = p-b$, therefore $C_a = C_a^{'}$, and it results that $A_1$, $C_a$, $A$ are collinear points and $A_1C_a$ contains the Gergonne's point $(\Gamma)$. Similarly, it can be shown that $I_bC_b$ and $I_cC_c$ pass through $\Gamma$.

## G. The medial triangle

**Definition 19**

A medial triangle (or complementally triangle) of given triangle is the triangle determined by the middle points of the sides of the given triangle.

**Definition 20**

We call a complete quadrilateral the figure $ABCDEF$ where $ABCD$ is convex quadrilateral and $\{E\} = AB \cap CD$, $\{F\} = BC \cap AD$. A complete quadrilateral has three diagonals, and these are $AC$, $BD$, $EF$

**Theorem 8** (Newton-Gauss)

The middle points of the diagonals of a complete quadrilateral are collinear (the Newton-Gauss line).

**Proof**

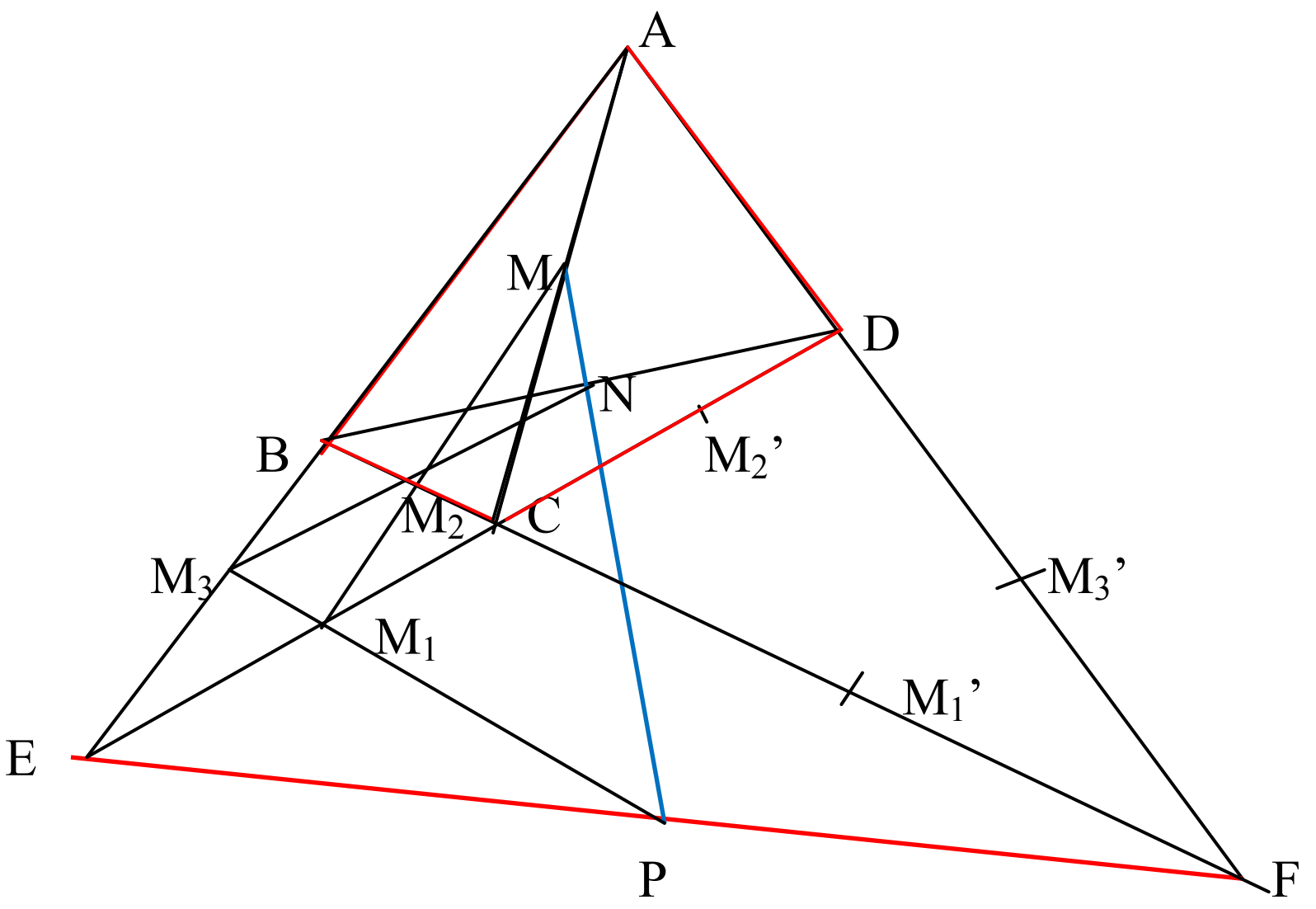

Fig. 17

Let $M_1M_2M_3$ the medial triangle of triangle $BCD$ (see figure 17). We note

$$M_1M_2 \cap AC = \{M\}, \quad M_2M_3 \cap BD = \{N\}, \quad M_1M_3 \cap EF = \{P\}$$

taking into account of the middle lines which come up, the points $M$, $N$, $P$ are respectively the middle points of the diagonals $(AC)$, $(BD)$, $(EF)$ of the complete quadrilateral $ABCDEF$.

We have:

$$MM_1 = \frac{1}{2}AE \quad MM_2 = \frac{1}{2}AB \qquad NM_2 = \frac{1}{2}CD \quad NM_3 = \frac{1}{2}DE$$

$$PM_3 = \frac{1}{2}BF \quad PM_1 = \frac{1}{2}CF$$

Let's evaluate following relation:

$$\frac{MM_1}{MM_2} \cdot \frac{NM_2}{NM_3} \cdot \frac{PM_3}{PM_1} = \frac{AE}{AB} \cdot \frac{DC}{DE} \cdot \frac{FB}{CF}$$

Considering $A, D, F$ transversal in triangle $BCE$ we have, in conformity to Menelaus' theorem, that $\frac{AE}{AB} \cdot \frac{DC}{DE} \cdot \frac{FB}{CF} = 1$ and respectively that $\frac{MM_1}{MM_2} \cdot \frac{NM_2}{NM_3} \cdot \frac{PM_3}{PM_1} = 1$.

From the Menelaus' reciprocal theorem for triangle $M_1M_2M_3$ it results that the points $M$, $N$, $P$ are collinear and the Newton-Gauss theorem is proved.

**Remark 12**

We can consider the medial triangle $M_1'M_2'M_3'$ of triangle $CFD$ ($M_1'$ the middle of $(CF)$ and $M_3'$ the middle of $(CD)$ and $M_3'$ the middle of $(DF)$) and the theorem can be proved in the same mode. Considering this triangle, it results that triangles $M_1M_2M_3$ and $M_1'M_2'M_3'$ have as homological axis the Newton-Gauss line. Their homological center being the intersection of the lines: $MM_1'$, $M_2M_2'$, $M_3M_3'$.

**Proposition 14**

The medial triangle and the anti-supplemental triangle of a non-isosceles given triangle are homological. The homological center is the symmedian center of the anti-supplemental triangle, and the homological axis is the Newton-Gauss line of the complete quadrilateral which has as sides the sides of the given triangle and the polar of the center of the inscribed circle in that triangle.

**Proof**

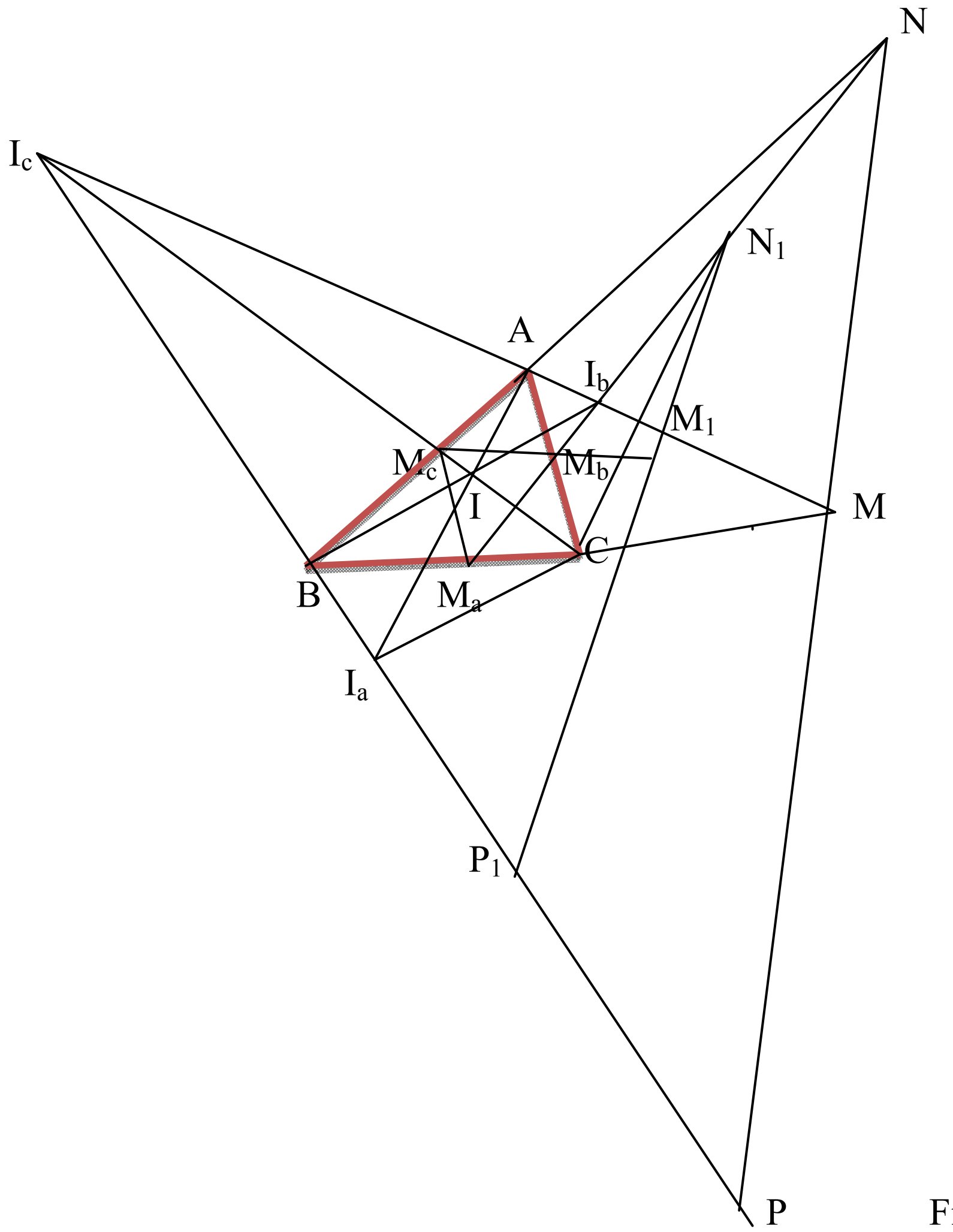

Fig. 18

In a triangle the external bisectrix are perpendicular on the interior bisectrix of the same angles, it results that $I$ is for the anti-supplemental $I_aI_bI_c$ (see figure 18) The orthocenter of the given triangle $ABC$ is the orthic triangle of $I_aI_bI_c$, therefore $BC$ is anti-parallel to $I_bI_c$. In accordance to proposition 6 it results that $I_aM_a$ is symmedian in triangle anti-supplemental, therefore $I_b M_b$ and $I_c M_c$ are symmedians and because these are also concurrent, it result that the triangles $M_a M_b M_c$ and $I_a I_b I_c$ are homological, the homology center being the symmedian center of $I_a I_b I_c$. We note $N, M, P$ the tri-linear polar of $I$, the line determined by the exterior bisectrix feet of the triangle $ABC$.

We note

$$\{N_1\} = M_a M_b \cap I_a I_b,\ \{M_1\} = M_b M_c \cap I_b I_c,\ \{P_1\} = M_a M_c \cap I_a I_c$$

The line $M_1, N_1, P_1$ is the homological axis of triangles $M_a M_b M_c$ and $I_a I_b I_c$, as it can be easily noticed it is the Newton-Gauss line of quadrilateral $NACMBP$ (because $N_1$ is the middle point of $NC$, $M_b$ is the middle of $AC$ and $M_a M_b \parallel AN$, etc.)

**Proposition 15**

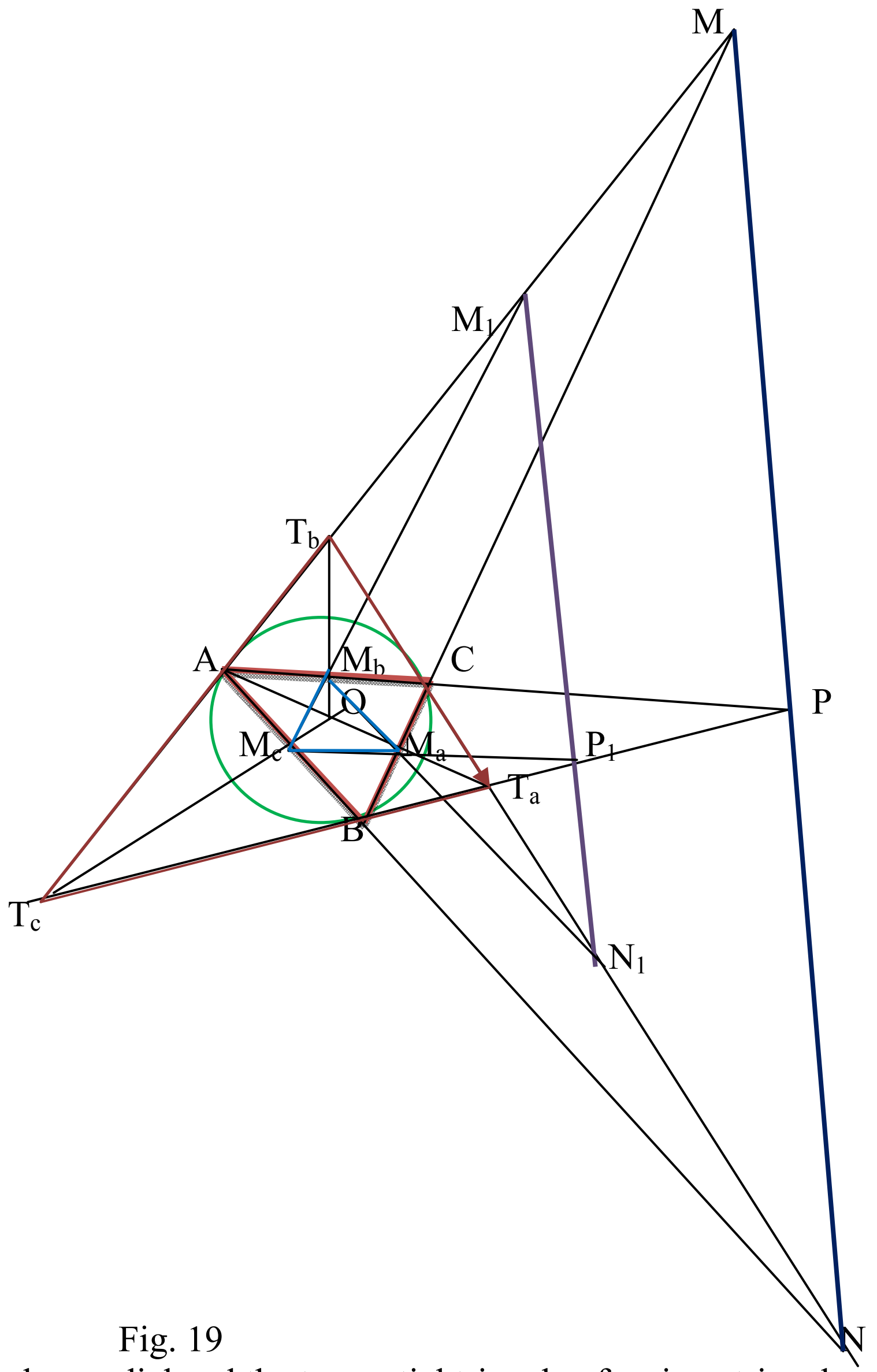

Fig. 19

The triangle non-isosceles medial and the tangential triangle of a given triangle are homological. The homology center is the center of the circumscribed circle of the given triangle, and the homological axis is the Newton-Gauss line of the complete quadrilateral which has as sides the sides of the given triangle and its Lemoine's line.

**Proof**

The lines $M_aT_a, M_bT_b, M_cT_c$ are the perpendicular bisector in the triangle $ABC$ , therefore are concurrent in $O$ (triangle $T_aBC$ is isosceles, etc.) We note $M, N, P$ the homological axis of triangles $ABC$ and $T_aT_bT_c$ (see figure 18). We'll note $\{M_1\} = M_bM_c \cap T_bT_c$. We observe that from the fact that $M_bM_c$ is middle line in triangle $ABC$ , it will pass through $M_1$, the middle point of $(AM)$. $\{M\} = BC \cap T_bT_c$

Similarly, $N_1$ is the middle of $(CN)$ and $P_1$ is the middle of $(BP)$

The triangles $M_aM_bM_c$ and $T_aT_bT_c$ being homological it results that $M_1, N_1, P_1$ are collinear and from the previous affirmations these belong to Newton-Gauss lines of the complete quadrilateral $BNPCAM$ , which has as sides the sides of the triangle $ABC$ and its Lemoine's line $M, N, P$ .

**Remark 13**

If we look at figure 18 without the current notations, with the intention to rename it later, we can formulate the following proposition:

**Proposition 16**

The medial triangle of the contact triangle of a given non-isosceles triangle is homological with the given triangle. The homological center is the center of the inscribed circle in the given triangle, and the homology axis is the Newton-Gauss's line of the complete quadrilateral which has the sides the given triangle's sided and the tri-linear polar is the Gergonne's point of the given triangle.

**H. The cotangent triangle**

**Definition 21**

A cotangent triangle of another given triangle the triangle determined by the tangent points of the ex-inscribed circles with the triangle sides.

**Observation 17**

In figure 19 we note the cotangent triangle of triangle $ABC$ with $J_aJ_bJ_c$.

**Definition 22**

Two points on the side of a triangle are called isometric if the point of the middle of the side is the middle of the segment determined by them.

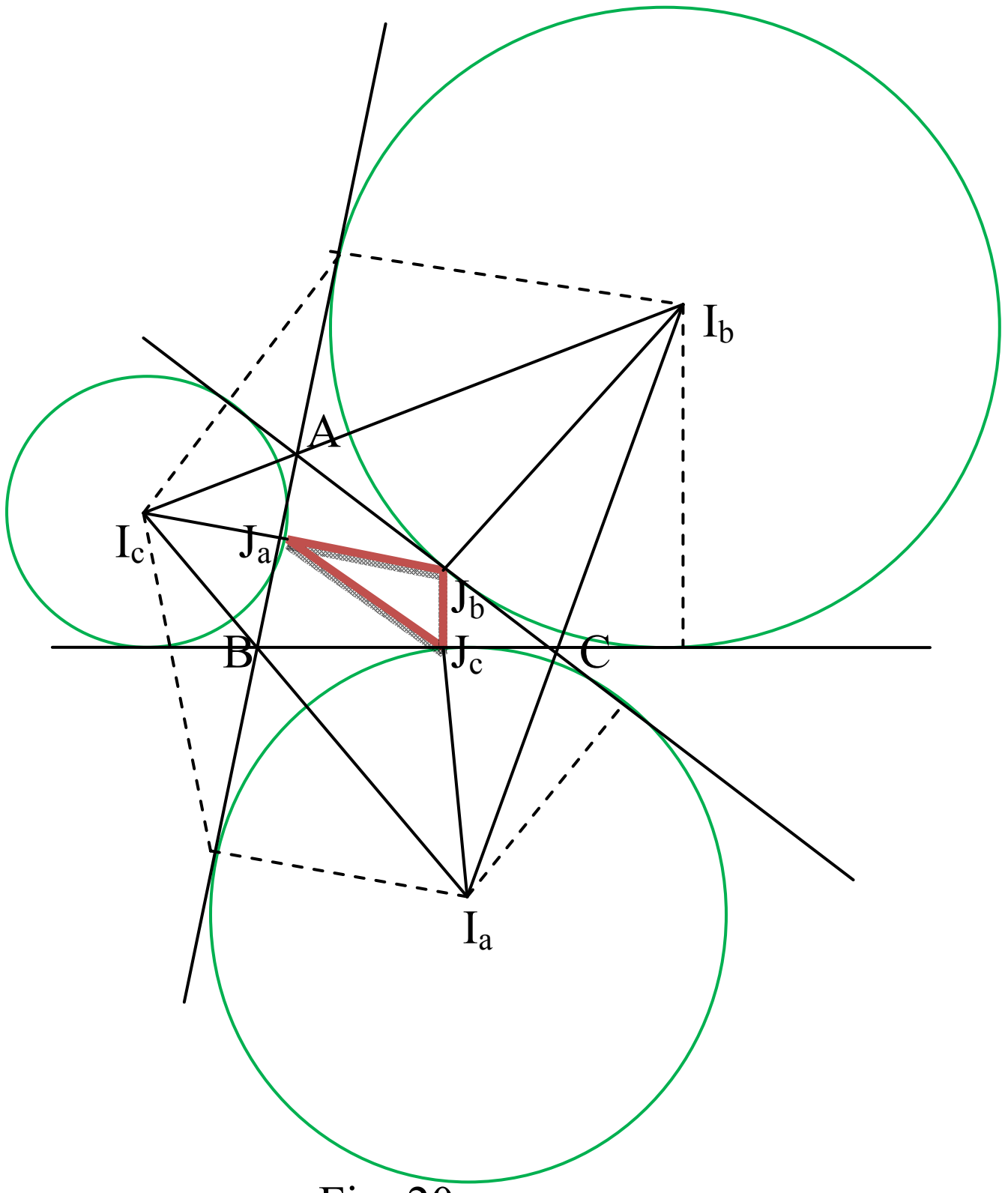

Fig. 20

**Proposition 17**

In a triangle the contact point with a side with the inscribed and ex-inscribed are isometrics

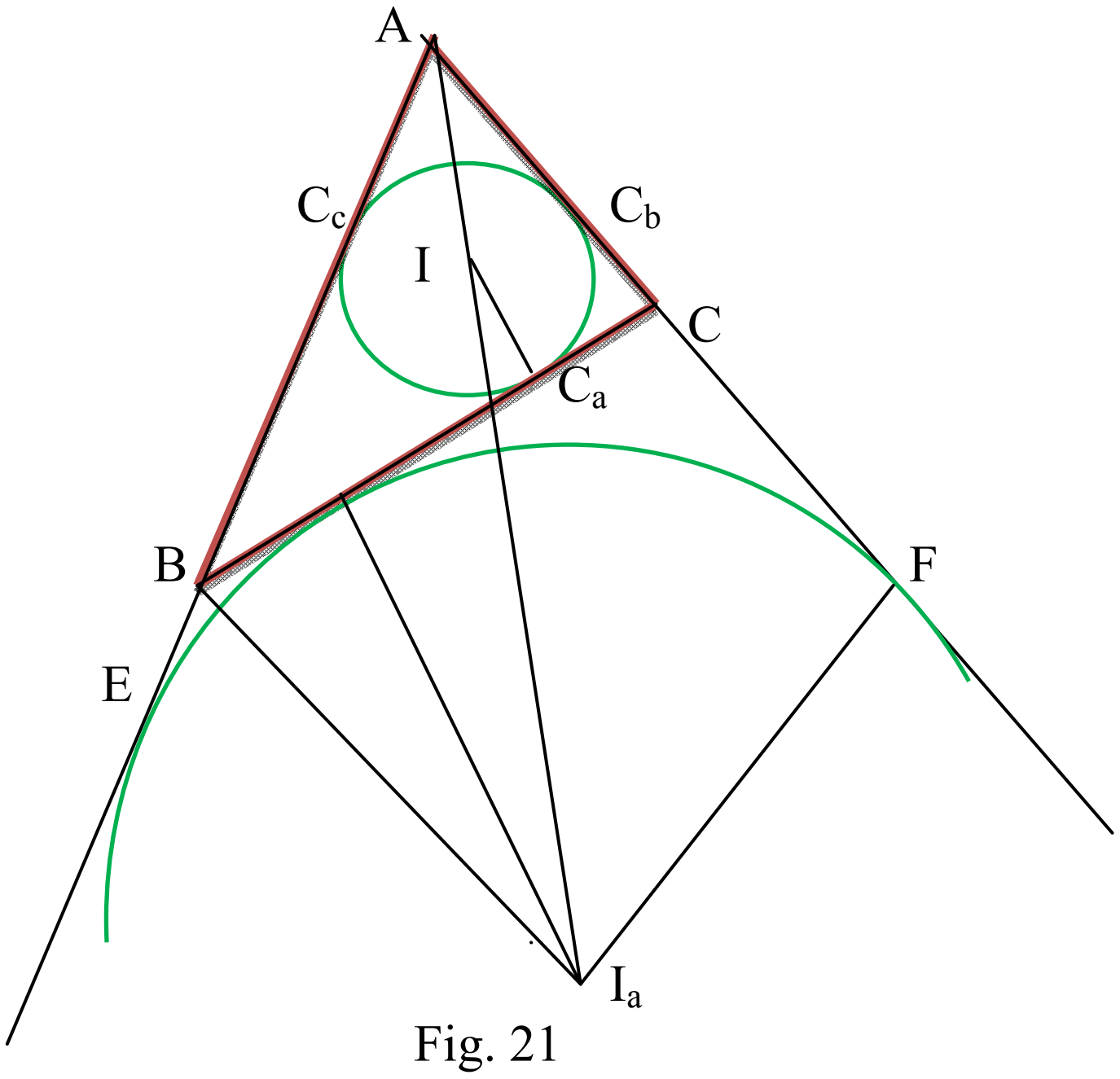

Fig. 21

**Proof**

Consider the triangle $ABC$ and the side $BC$ to have $C_a$ and $I_a$ the contact points of the inscribed and A-ex-inscribed (see figure 21).

To prove that $C_a$ and $I_a$ are isometrics, in other words to prove that are symmetric in rapport to the middle of $(BC)$ is equivalent with showing that $BI_a = CC_a$. Will this through computation, finding the expressions of these segments in function of the lengths $a,b,c$ of the triangle. We'll note:

$$x = AC_b = AC_c$$
$$y = BC_a = BC_c$$
$$z = CC_a = CC_c$$

From the system:

$$\begin{cases} x + y = c \\ y + z = a \\ z + x = b \end{cases}$$

By adding them and taking into account that $a + b + c = 2p$, we find:

$$x = p - a$$
$$y = p - b$$
$$z = p - c$$

Therefore $CC_a = p - c$.

We note with $E, F$ the tangent points of the A ex-inscribed circle with $AB$ and $AC$ and also

$$y' = BI_a = BE$$
$$z' = CI_a = CF$$

We have also $AE = AF$, which gives us:

$$\begin{cases} c + y' = b + z' \\ y' + z' = a \end{cases}$$

From this system we find

$$y' = \frac{1}{2}(a + b - c) = p - c$$

Therefore $BI_a = p - c = CC_c$, which means that the points $I_a, C_a$ are isometric.

**Observation 18**

$AE = AF = p$

**Definition 23**

Two Cevians of the same vertex of the same triangle are called isotomic lines if their base (feet) are isotomic points.

**Theorem 9** (Neuberg's theorem)

The isotomic Cevians of concurrent Cevians are concurrent

**Proof**

In figure 22, let's consider in the triangle $ABC$ the concurrent Cevians in point $P$ noted with $AP_1, BP_2, CP_3$ and the Cevians $AQ_1, BQ_2, CQ_3$ their isotomic lines.

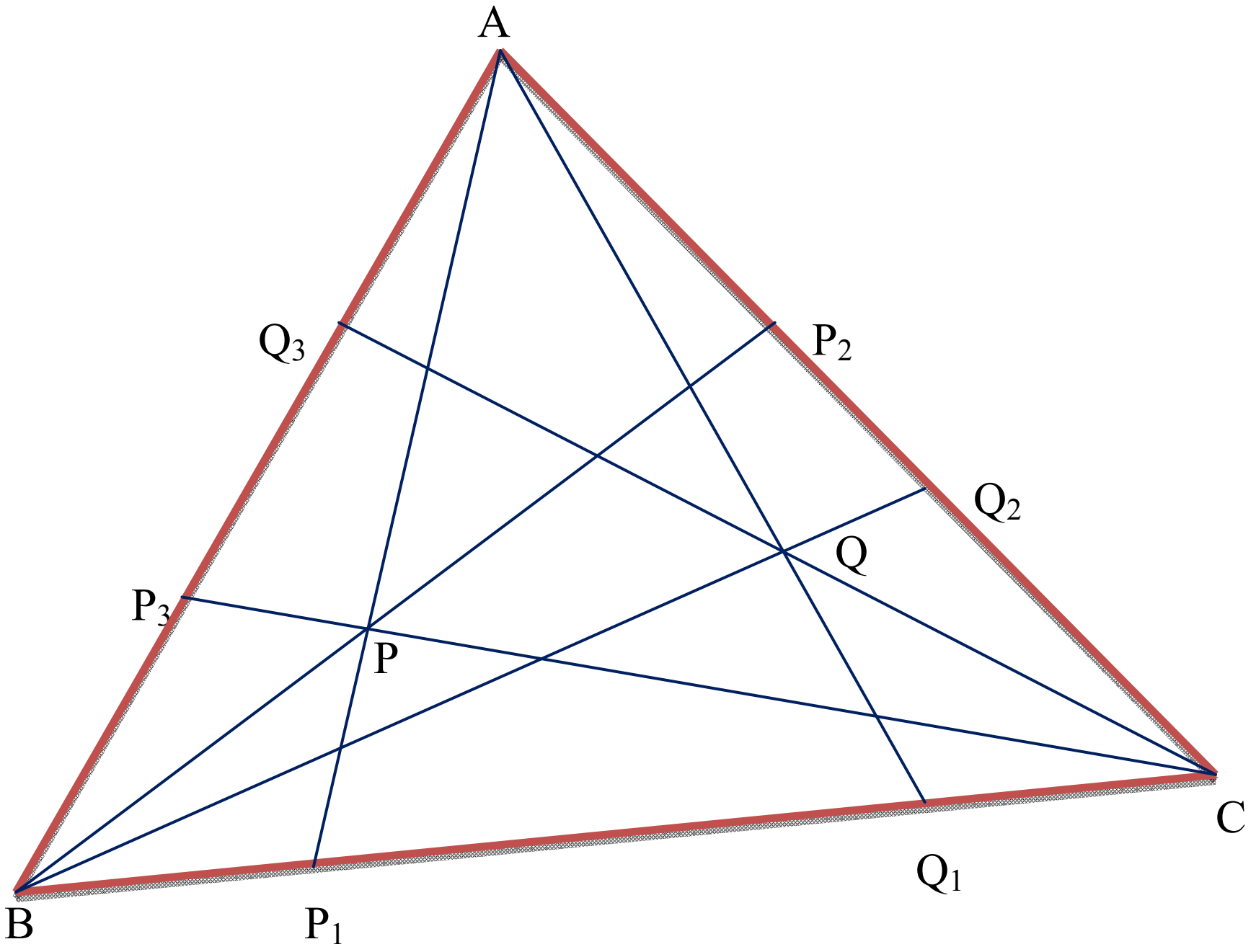

Fig. 22

From Ceva's theorem it results that

$$\frac{P_1B}{P_1C}\cdot\frac{P_2C}{P_2A}\cdot\frac{P_3A}{P_3B}=1 \quad (39)$$

Because

$$P_1B = Q_1C, P_2A = Q_2C, P_1C = Q_1B, P_2C = Q_2A, Q_3A = P_3B, P_3A = Q_3B$$

We can write

$$\frac{Q_1C}{Q_1B}\cdot\frac{Q_2A}{Q_2C}\cdot\frac{Q_3B}{Q_3A}=1 \quad (40)$$

The Ceva's reciprocal theorem implies the concurrence of the Cevians $AQ_1, BQ_2, CQ_3$. We noted with $Q$ their concurrence point.

**Definition 24**

The points of concurrence of the Cevians and their isotomic are called isotomic conjugate.

**Theorem 11.**(Nagel)

In a triangle the Cevians $AI_a, BI_b, CI_c$ are concurrent.

**Proof**

The proof results from theorem 9 and from proposition 15.

**Definition 25**

The conjugate isotomic point of Gergonne's point $(\Gamma)$ is called Nagel's point $(N)$.

**Observation 20**

The concurrence point of the Cevians $AH_a, BH_b, CH_c$ is the Nagel point $(N)$.

**Theorem 10**

Let $A', B', C'$ the intersection points of a line $d$ with the sides $BC, CA, AB$ of a given triangle $ABC$. If $A'', B'', C''$ are isotomic of the points $A', B', C'$ respectively then $A'', B'', C''$ are collinear.

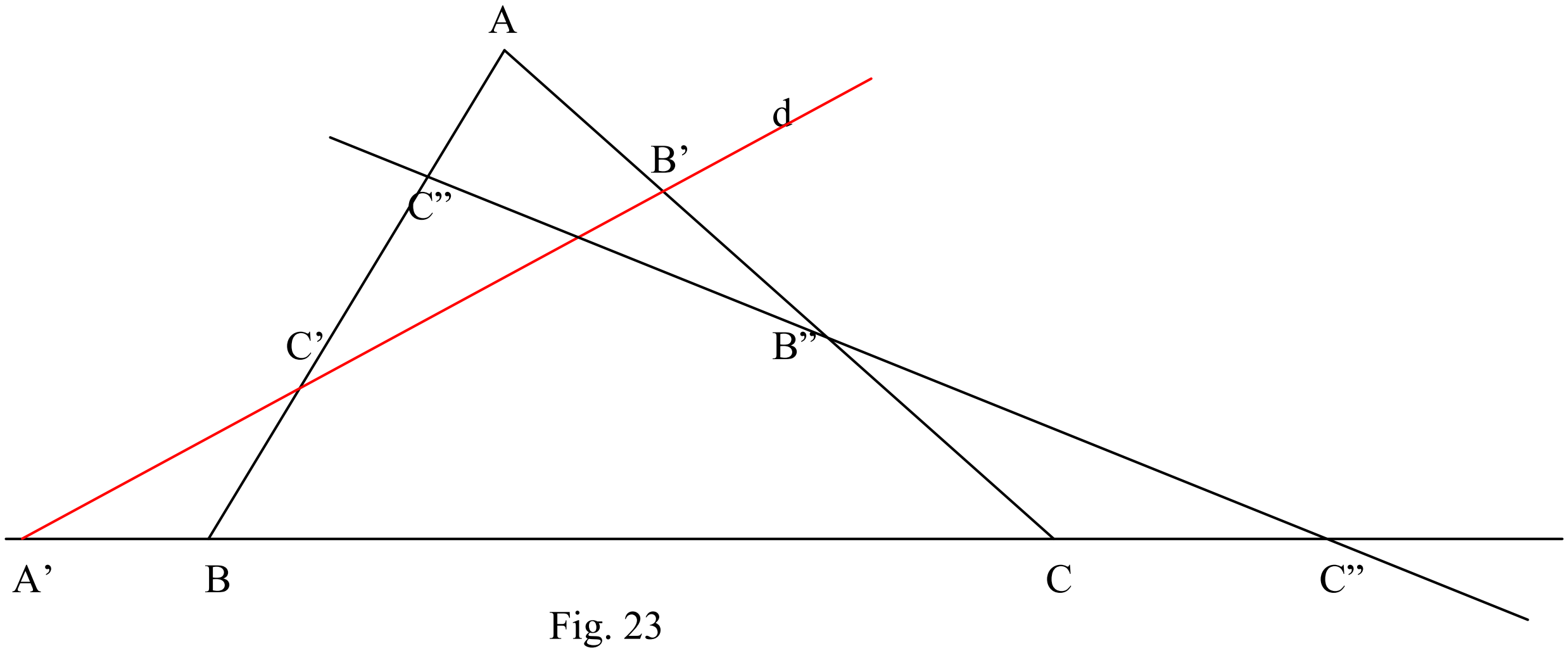

Fig. 23

**Proof**

The $A', B', C'$ being collinear, from the Menelaus' reciprocal theorem it results:

$$\frac{A'B}{A'C}\cdot\frac{B'C}{B'A}\cdot\frac{C'A}{C'B}=1 \qquad (41)$$

Because

$A'B = A''C, A'C = A''B;\ B'A = B''C, B'C = B''A;\ C'A = C''B, C'B = C''A$

And substituting in (41) it result

$$\frac{A''C}{A''B}\cdot\frac{B''A}{B''C}\cdot\frac{C''B}{C''A}=1 \qquad (42)$$

This relation and the Menelaus' s reciprocal theorem shows the collinearity of the points $A'', B'', C''$.

**Remark 14**

The line of the points $A'', B'', C''$ is called the isotomic transversal.

**Proposition 18**

A given non-isosceles triangle and its cotangent triangle are homological triangles. The homology center is the Nagel's point of the given triangle and the homology axis is the isotomic transversal of the Lemoine's line of the contact triangle.

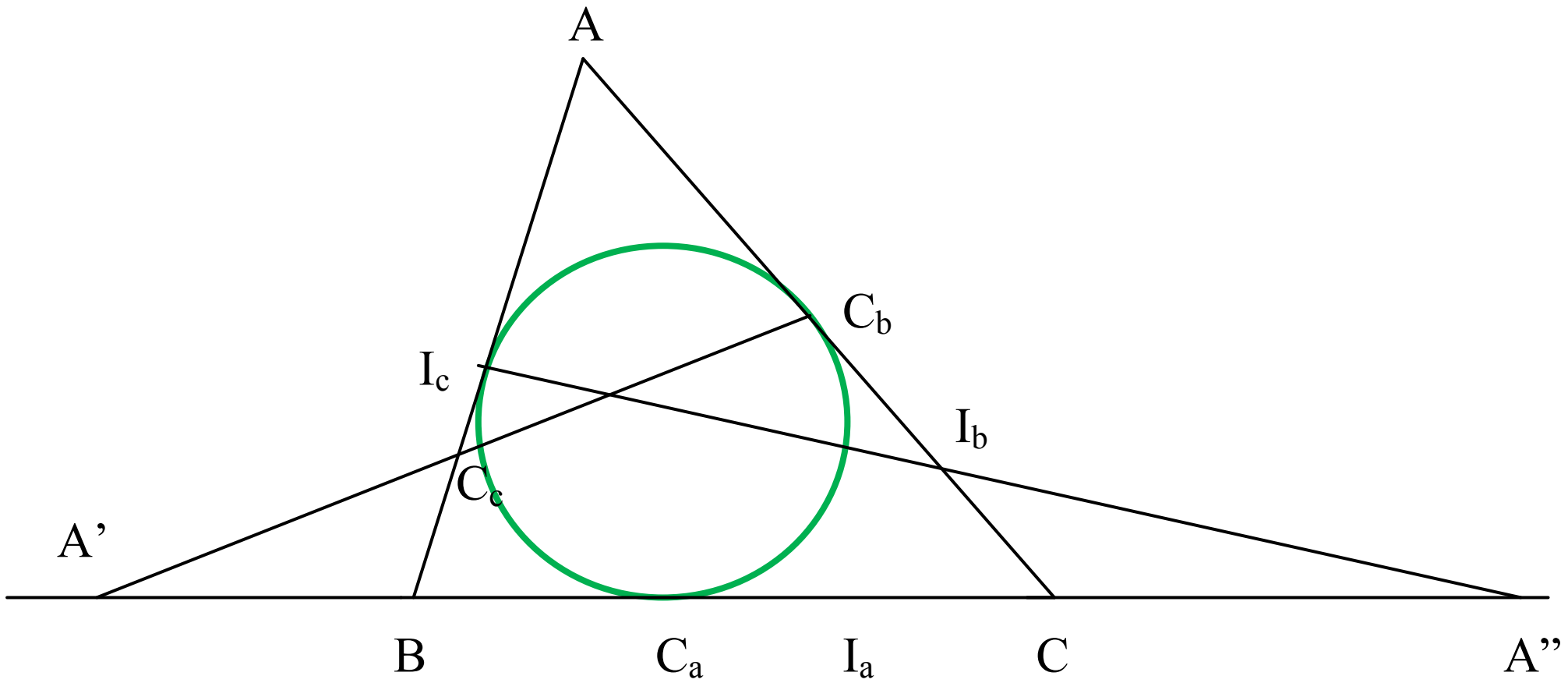

Fig. 24

**Proof**

From theorem 10 it results that the lines $AI_a, BI_b, CI_c$ are concurrent in the Nagel's point, therefore triangle $ABC$ and $I_aI_bI_c$ have the center of homology the point $N$.

We note $\{A'\} = C_cC_b \cap BC$ and $\{A''\} = I_bI_c \cap BC$ (see figure 24).

Will show that $A', A''$ are isotomic points. We know that $A'$ is the harmonic conjugate in rapport to $B$ and $C$.

Also $A''$ is the harmonic conjugate of $I_c$ in rapport to $B, C$ and we also know that $C_a$ and $I_a$ are isotomic points. ($BC_a = CI_a = p-b; CC_a = BI_a$).

We have:

$$\frac{A'B}{A'C} = \frac{C_aB}{C_aC}$$

that is:

$$\frac{A'B}{A'C} = \frac{p-b}{p-c}$$

from which:

$$A'B = \frac{a(p-b)}{p-c} \tag{43}$$

$$\frac{A''C}{A''B} = \frac{I_aC}{I_aB}$$

therefore

$$\frac{A''C}{A''B} = \frac{p-b}{p-c}$$

from which

$$A''C = \frac{a(p-b)}{c-b} \tag{44}$$

Relations (43) and (44) show that $A'$ and $A''$ are conjugated points. Similarly we prove that $B', B''$; C',C'' are conjugated. The homology axes $A', B', C'$ and $A'', B'', C''$ of the contact triangle and the triangle $ABC$ respectively of the co-tangent triangle and of triangle $ABC$ are therefore isotomic transversals.

**Observation 21**

Similarly can be proved the following theorem, which will be the generalization of the above result.

**Theorem 12**

The Cevians triangle of two isotomic conjugated points in given triangle and that triangle are homological. The homology axes are isotomic transversals.

**Proposition 19**

The medial triangle of the cotangent triangle of a non-isosceles triangle $ABC$ is homological to the triangle $ABC$ .

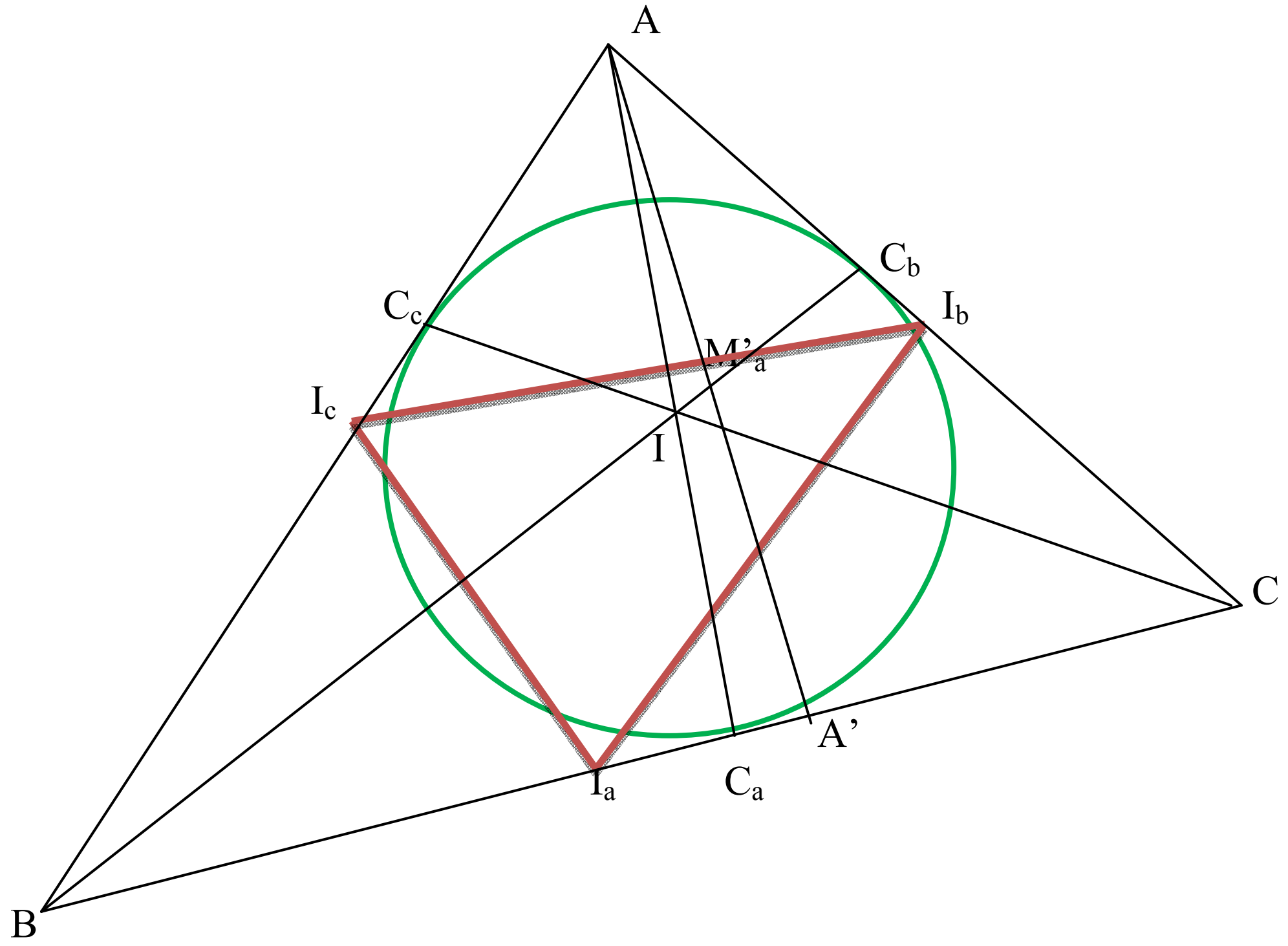

Fig. 25

**Proof**

We note $M'_a, M'_b, M'_c$ the vertexes of the medial triangle of the cotangent circle $I_a I_b I_c$ (see figure 25) and with $\{A'\} = AM'_a \cap BC$. We have

$$\frac{A'B}{A'C} = \frac{aria\Delta ABA'}{aria\Delta ACA'} = \frac{AB\sin BAA'}{AC\sin CAA'} \tag{45}$$

On the other side:

$$aria\Delta AI_aM_a^{'} = aria\Delta AI_bM_a^{'}$$

from here we find that

$$AI_c \sin BAA' = AI_b \sin CAA'$$

therefore

$$\frac{\sin BAA'}{\sin CAA'} = \frac{AI_b}{AI_c} = \frac{CC_b}{BC_b} = \frac{p-c}{p-b} \tag{46}$$

Looking in (45) we find:

$$AD \cdot DA_1 = \frac{a^2bc}{(b+c)^2} \tag{47}$$

With the notations $\{B'\} = BM_b^{'} \cap AC$ and $\{C'\} = CM_c^{'} \cap AB$ we proceed on the same manner and we find

$$\frac{B'C}{B'A} = \frac{a(p-a)}{c(p-c)} \tag{48}$$

$$\frac{C'C}{C'B} = \frac{b(p-b)}{a(p-a)} \tag{49}$$

The last three relations and Ceva's reciprocal theorem, lead us to the concurrency of the lines: $AM'_a, BM'_b, CM'_c$ therefore to the homology of the triangles $ABC$ and $M_a^{'}M_b^{'}M_c^{'}$.

We note

$$\{M\} = I_bI_c \cap BC,$$

$$\{N\} = I_aI_b \cap AB,$$

$$\{P\} = I_aI_c \cap AC$$

and

$$\{M_1\} = M_a^{'}M_b^{'} \cap AB,$$

$$\{N_1\} = M_b^{'}M_c^{'} \cap BC,$$

$$\{P_1\} = M_c^{'}M_a^{'} \cap CA$$

We observe that $N_1$ is the middle of the segment $(I_cN)$, $N_1$ is the middle of the segment $(I_aM)$ and $P_1$ is the middle of the segment $(I_bP)$. In the complete quadrilateral $PMI_bI_aNI_c$ the Newton-Gauss line is $N_1M_1P_1$ therefore the homological axis of triangles $ABC$ and $M_a^{'}M_b^{'}M_c^{'}$.

**Proposition 20**

The cotangent triangle of a given triangle and the triangle $A'B'C'$ formed by the projection of the circle inscribed to triangle $ABC$ on the perpendicular bisectors of the sides $BC, CA, AB$ are homological. The homology center is Nagel's point $N$ of the triangle $ABC$.

**Proof**

Let $A_1$ the intersection of the bisector $A'$ with the circumscribed circle of triangle $ABC$, $M_a$ the middle of $BC$ and $I_a$ the vertex on $BC$ of the cotangent triangle (see figure 26).

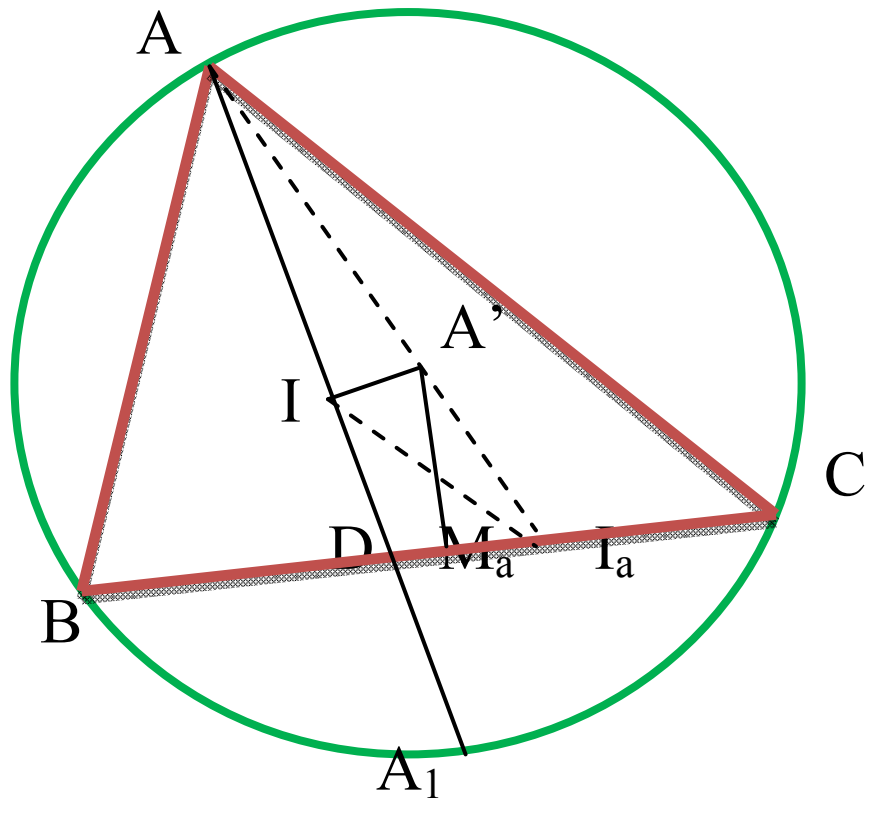

Fig. 26

Consider the power of $I$ in rapport to the circumscribed circle of triangle $ABC$, we have $AI \cdot IA_1 = R^2 - OI^2$. It is known that $DI^2 = R^2 - 2Rr$. Substituting we find

$$AI \cdot IA_1 = 2Rr \qquad (*)$$

Considering the power of $D$ in rapport to the circumscribed circle of the triangle $ABC$ we have $AD \cdot DA_1 = BD \cdot DC$.

We know that

$$BD = \frac{ac}{b+c}$$

$$CD = \frac{ab}{b+c}$$

It results

$$AD \cdot DA_1 = \frac{a^2bc}{(b+c)^2} \qquad (**)$$

On the other side with bisector's theorem we obtain

$$\frac{AI}{ID} = \frac{b+c}{a}$$

$$\frac{AI}{ID} = \frac{b+c}{a+b+c} \qquad (***)$$

Dividing relations (*) and (**) side by side it results

$$\frac{AI}{ID}=\frac{IA_1}{DA_1}=\frac{2Rr\left(b+c\right)^2}{a^2bc} \qquad (****)$$

Taking into account (***) and formula $abc=4RS$ and $S=p\cdot r$, we find

$$\frac{IA_1}{DA_1}=\frac{b+c}{a} \qquad (*****)$$

Because $DM_a=BM_a-BD$, we find that

$$DM_a=\frac{a(p-c)}{2(b+c)}$$

From the fact that triangles $A_1A'I$ and $A_1MD$ are similar we have

$$\frac{IA_1}{DA_1}=\frac{IA'}{DM_a}$$

From here we obtain

$$IA'=\frac{b-c}{2} \qquad (******)$$

We note $I_a^{'}=AA'\cap BC$. From the similarity of the triangles $AIA'$ and $ADI_a^{'}$ it results

$$DI_a'=\frac{(a+b+c)(b-c)}{2(b+c)}$$

$BI_a^{'}=BD+DI_a^{'}=p-c$, this relation shows that $I_a=I_a^{'}$ contact of the A-ex-inscribed circle or $AA'$ is Cevian Nagel of the triangle $ABC$. Similarly it results that $I_bB'$ contains the Nagel's point.

### I. The ex-tangential triangle

**Definition 26**

Let $ABC$ and $I_a, I_b, I_c$ the centers of the A-ex-inscribed circle, B-ex-inscribed circle, C.-ex-inscribed circle. The common external tangents to the ex-inscribed circles (which don't contain the sides of the triangle $ABC$ ) determine a triangle $E_aE_bE_c$ called the ex-tangential triangle of the given triangle $ABC$

**Proposition 21**

The triangles ex-tangential and anti-supplemental of a given non-isosceles triangle are homological.

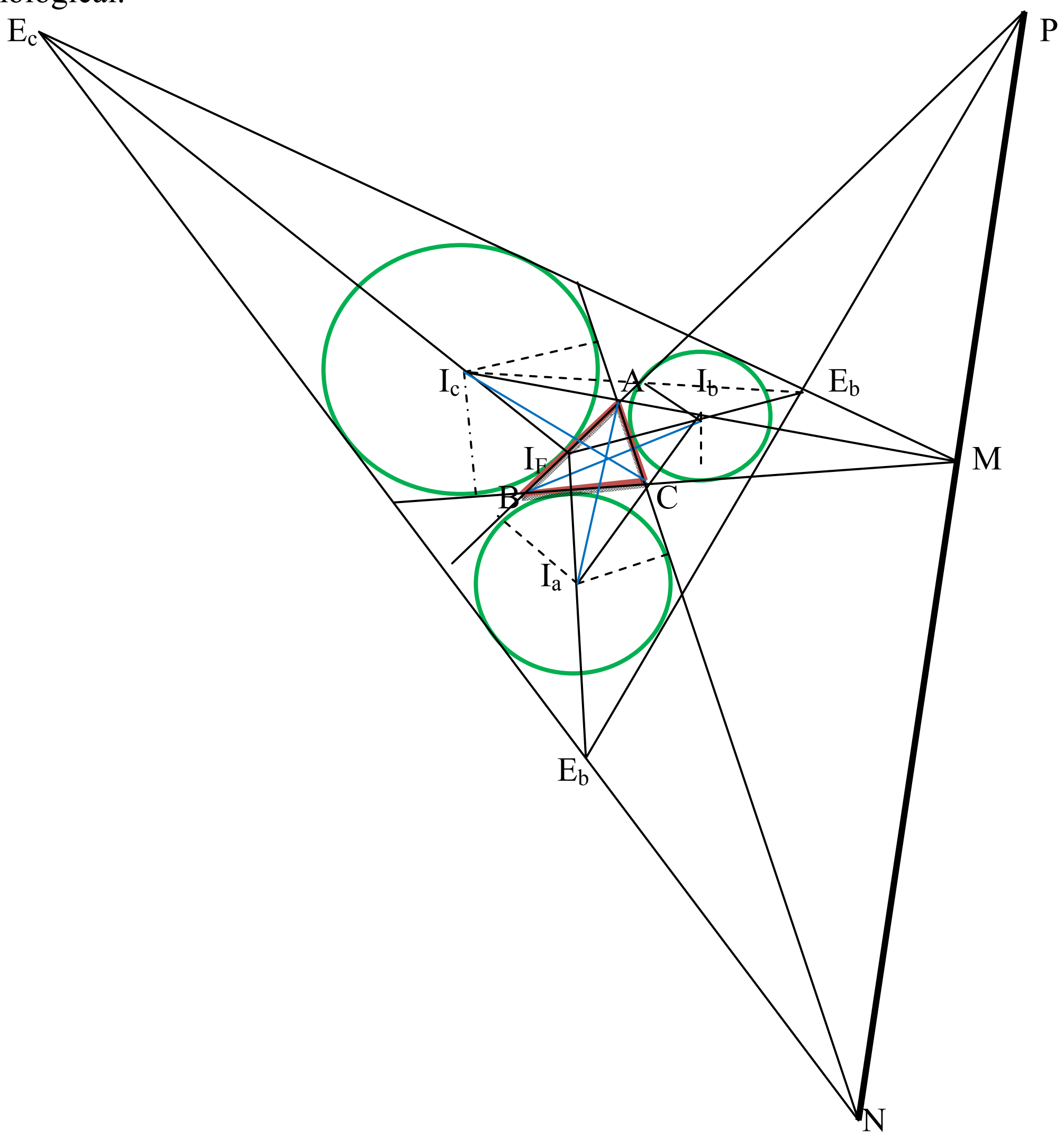

Fig. 27

The homology center is the center of the inscribed circle in the ex-tangential triangle, and the homology axis is the anti-orthic axis of the given triangle.

**Proof**

The lines $E_aI_a$, $E_bI_b$, $E_cI_c$ are the bisectors in the ex-tangential triangle, therefore these are concurrent in the center of the circle inscribed to this circle, noted $I_E$ (see figure 27). Therefore it results that triangles $E_aE_bE_c$ and $I_aI_bI_c$ are homological.

We note

$$I_bI_c \cap E_bE_c = \{M\},$$
$$I_aI_c \cap E_aE_c = \{N\},$$
$$I_aI_b \cap E_aE_b = \{P\}$$

The homology axis of triangles $E_aE_bE_c$ and $I_aI_bI_c$ is, conform to Desargues' theorem the line $M$, $N$, $P$. Because $E_bE_c$ and $BC$ intersect also in $M$, it results that is the feet of the exterior bisector of angle $A$ of triangle $ABC$. Consequently $M$, $N$, $P$ is the anti-orthic axis of triangle $ABC$.

**Remark 15**

From the above affirmation it results that the anti-orthic axis of triangle $ABC$ is the homology axis for triangle $ABC$ as well as for triangle $E_aE_bE_c$. Therefore we can formulate

**Proposition 22**

A given triangle and its ex-tangential triangle are homological. The homology axis is the anti-orthic axis of the given triangle.

**Remark 16**

a) The homology center of triangle $ABC$ and of its ex-tangential $E_aE_bE_c$ is the intersection of the lines $AE_b$, $BE_b$, $CE_c$.

b) From the proved theorem it results (in a particular situation) the following theorem.

**Theorem 13** (D'Alembert 1717-1783)

The direct homothetic centers of three circles considered two by two are collinear, and two centers of inverse homothetic are collinear with the direct homothetic center which correspond to the third center of inverse homothetic.

Indeed, the direct homothetic centers of the ex-inscribed circles are the points $M, N, P$, and the inverse homothetic centers are the points $A, B, C$. More so, we found that the lines determined by the inverse homothetic centers and the vertexes of the ex-tangential triangle are concurrent.

**Observation 22**

Considering a given isosceles triangle, its anti-supplemental triangle $I_aI_bI_c$ and its ex-tangential triangle $E_aE_bE_c$ it has been determined that any two are homological and the homology axis is the anti-orthic axis of the triangle $ABC$. We will see in the next paragraph what relation does exist between the homological centers of these triangles.

**J. The circum-pedal triangle (or meta-harmonic)**

**Definition 27**

We define a circum-pedal triangle (or meta-harmonic) of a point $D$, from the plane of triangle $ABC$, in rapport with the triangle $ABC$ - the triangle whose vertexes are the intersections of the Cevians $AD, BD, CD$ with the circumscribed circle of the triangle $ABC$.

**Remark 17**

Any circum-pedal triangle of any triangle and the given triangle are homological.

**Proposition 23**

The circum-pedal triangle of the orthocenter $H$ of any triangle $ABC$ is the homothetic of the orthic triangle of that triangle through the homothety of center $H$ and of rapport 2.

**Proof.**

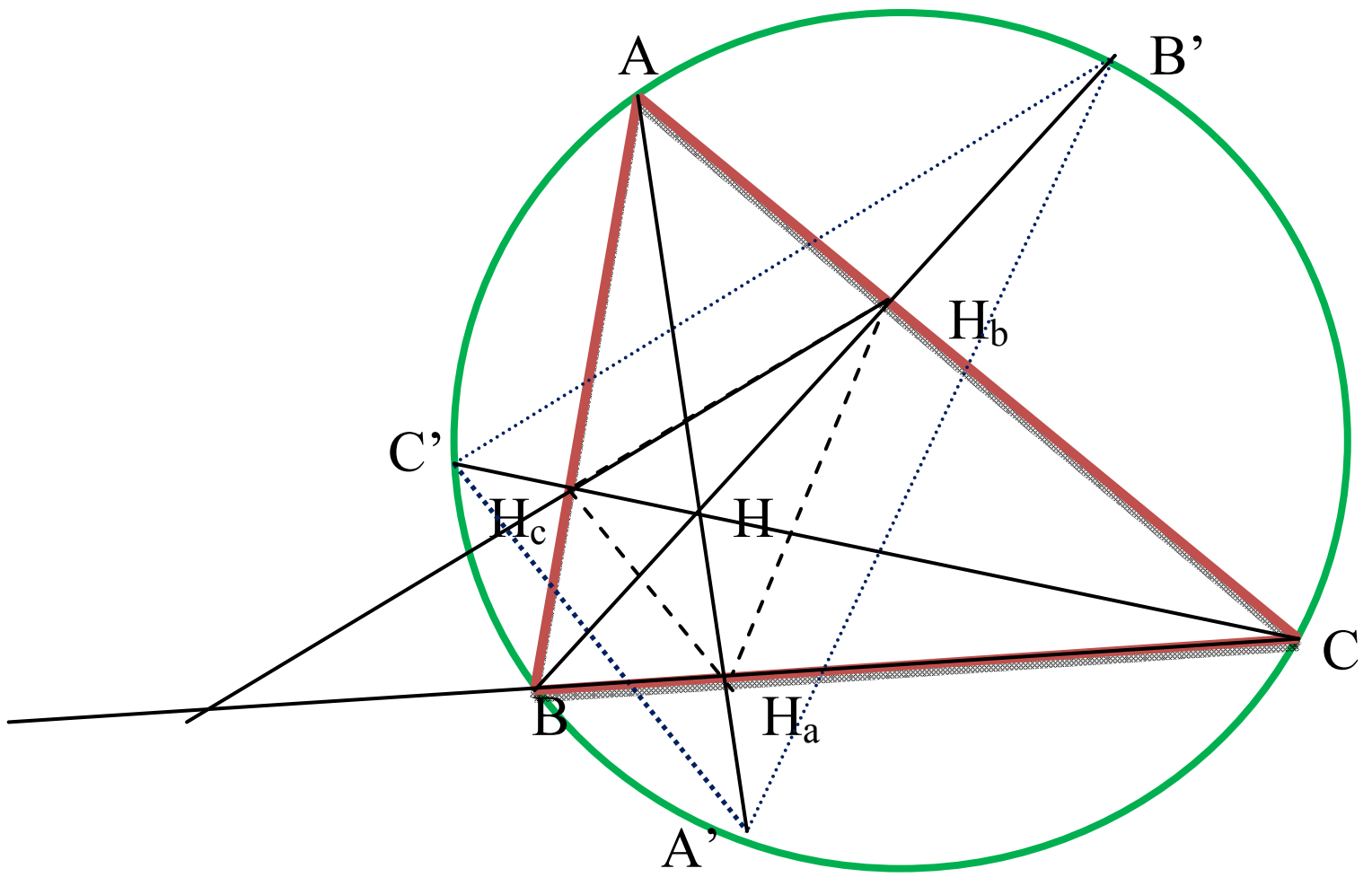

Fig.28

Let $ABC$ a scalene triangle. $H$ its orthocenter, $H_aH_bH_c$ its orthic triangle and $A'B'C'$ its circum-pedal triangle (see figure 28).

Because $\sphericalangle BA'H \equiv \sphericalangle BHA'$ (are inscribed in circle and have as measure $\frac{1}{2}m(AB)$) and $\sphericalangle BCA = \sphericalangle BHA'$ *angles with sides respectively perpendicular, we obtain that* $\sphericalangle BA'H \equiv \sphericalangle BHA'$, therefore the triangle $BHA'$ is isosceles. $BH_a$ being the altitude, it is the median, therefore $HH_a \equiv H_aA'$ or $H_aA' = 2HH_a$ which shows that $A'$ is homothetic to $H_a$. through the homothety $\mathcal{H}(H;2)$. The property is proved similarly for the vertexes $B'$ and $C'$ of the circum-pedal triangle as well as in the case of the rectangle triangle.

**Remark 18**

We will use the proposition 22 under the equivalent form: The symmetric of the orthocenter of a triangle in rapport with its sides belong to the circumscribed circle.

**Proposition 24**

The circum-pedal triangle of the symmedians center of a given triangle has the same symmedians as the given triangle.

**Proof (Efremov)**

Let $K$ the symmedian center of triangle $ABC$ and $DEF$ the pedal triangle of $K$ (see figure 29)

The quadrilateral $KDBF$ is inscribable; it results that

$$\sphericalangle KDF \equiv \sphericalangle KBF \equiv \sphericalangle B'BA \equiv \sphericalangle AA'B' \quad (50)$$

The quadrilateral $KDCE$ is inscribable, it results:

$$\sphericalangle KDE \equiv \sphericalangle KCE \equiv \sphericalangle AA'C' \quad (51)$$

From (50) and (51) we retain that $\sphericalangle EDF \equiv \sphericalangle B'A'C'$.

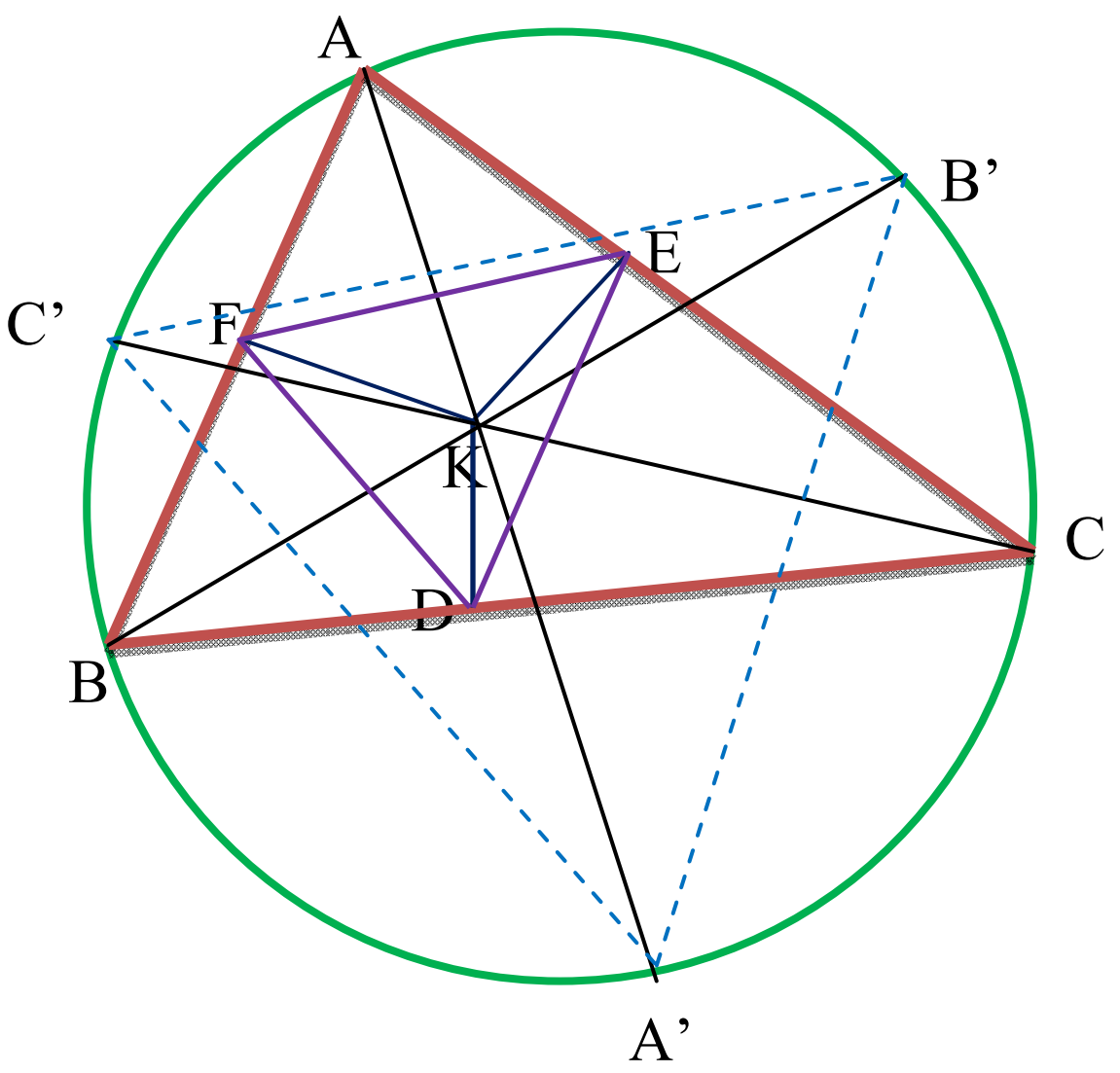

Fig. 29

Similarly we find that $\sphericalangle EDF \equiv \sphericalangle A'B'C'$, therefore $\Delta DEF \sim \Delta A'B'C'$.

Because $K$ is the gravity center of the triangle $DEF$, from $\sphericalangle EDK \equiv \sphericalangle BB'A$ and $\sphericalangle B'A'C' \equiv \sphericalangle EDF$ it results that $AA'$ is a symmedian in the triangle $A'B'C'$. Similarly we show that $BB'$ is a median in the same triangle, and the theorem is proved.

**Remark 18**

The triangles $ABC$ and $A'B'C'$ are called co-symmedians triangle.

**Proposition 24**

The homology axis of two co-symmedians triangles is the Lemoyne's line of one of them.

**Proof**

In the triangle $ABC$ we consider the symmedian center and $A'B'C'$ the circum-pedal triangle of $K$ (see figure 30).

It is known that Lemoine's line of triangle $ABC$ passes through $A_1$ which is the intersection constructed in $A$ to the circumscribed to triangle $ABC$ with $BC$.

We'll note $AA' \cap BC = \{S\}$

We have

$$\Delta ASB \sim \Delta CSA'$$

From where

$$\frac{AB}{A'C} = \frac{AS}{CS} \tag{52}$$

Also

$$\Delta ASC \sim \Delta BSA'$$

We obtain

$$\frac{AC}{A'C} = \frac{AS}{CS} \tag{53}$$

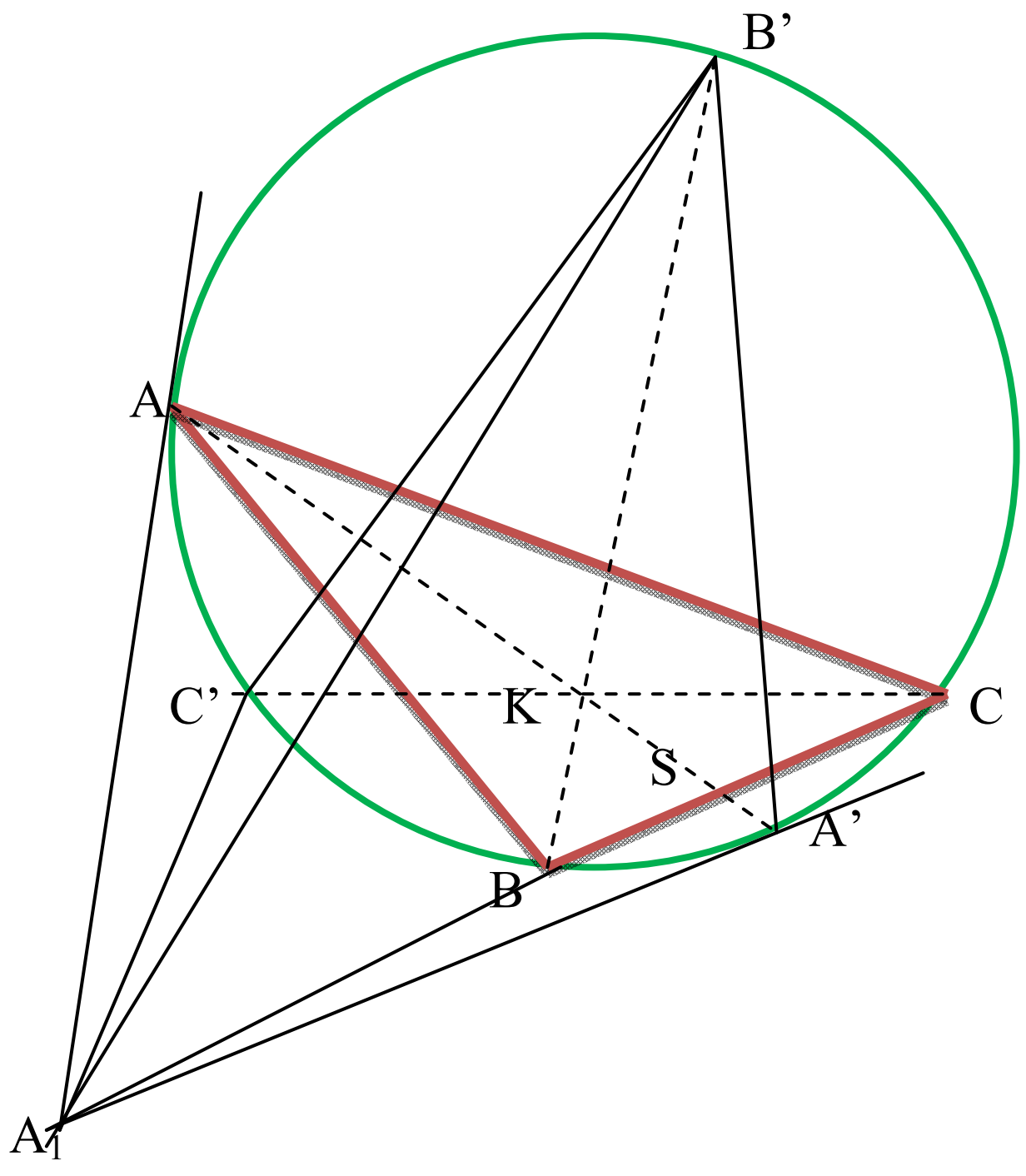

Fig. 30

From (52) and (53) it results

$$\frac{AB}{AC} \cdot \frac{BA'}{A'C} = \frac{BS}{CS} \tag{54}$$

On the other side $AS$ being symmedian we have

$$\frac{BS}{CS} = \frac{AB^2}{AC^2} \tag{55}$$

Therefore it results

$$\frac{BA'}{CA'} = \frac{AB}{AC} \tag{56}$$

The line $A_1A$ is ex-symmedian in triangle $ABC$, we have

$$\frac{A_1B}{A_1C} = \frac{AB^2}{AC^2} \tag{57}$$

We note with $A'_1$ the intersection of tangent in $A'$ with $BC$, because $A'_1A'$ is ex-symmedian in triangle $BA'C$, we have

$$\frac{A_1'B}{A_1'C} = \frac{A'B^2}{A'C^2} \tag{58}$$

Taking into account relation (56) it results

$$\frac{A_1'B}{A_1'C} = \frac{AB^2}{AC^2} \tag{59}$$

From (57) and (59) we find

$$\frac{A_1B}{A_1C} = \frac{A_1'B}{A_1'C}$$

Which shows that $A'_1 \equiv A_1$.

Applying the same reasoning for triangle $A'B'C'$ we find that the tangent in $A'$ intersects $B'C'$ and the tangent in $A$ at the circumscribed circle in the same point on $B'C'$, the tangent from $A$ and the tangent from $A'$ intersects in the unique point $A_1$. Consequently, $B'C'$ and $BC$ intersect in $A_1$.

**K. The Coşniţă triangle**

**Definition 28**

Given a triangle $ABC$, we define Coşniţă triangle relative to the given triangle, the triangle determined by the centers of the circumscribed circles to the following triangles $BOC, COA, AOB$, where $O$ is the center of the circumscribed circle of the given triangle $ABC$

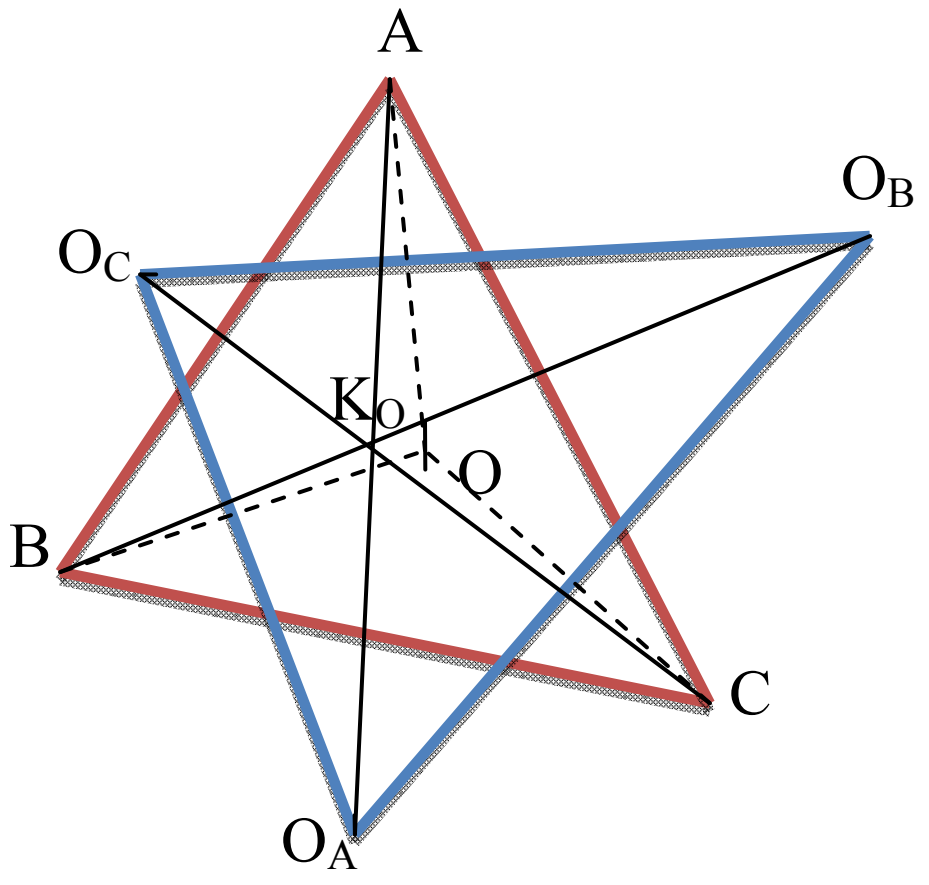

Fig. 31

**Theorem 14 (C. Coşniţă)**

The Coşniţă triangle of triangle $ABC$ and triangle $ABC$ are homological.

**Proof**

In figure 31 we construct a scalene triangle $ABC$ and we note $O_A, O_B, O_C$ the centers of the circumscribed triangles of triangles $BOC$, $COA$, $AOB$.

We have $m\sphericalangle(BOC) = 2A$, $m\sphericalangle(BOO_A) = A$, $m\sphericalangle(ABO_A) = 90° - (\sphericalangle C - \sphericalangle A)$

$m\sphericalangle(ACO_A) = 90° - (\sphericalangle B - \sphericalangle A)$.

The sinuses' theorem applied in triangles $ABO_A$ and $ACO_A$ leads to:

$$\frac{BO_A}{\sin(BAO_A)} = \frac{AO_A}{\sin(90^{\circ} - (\sphericalangle C - \sphericalangle A))} \quad (*)$$

$$\frac{CO_A}{\sin(CAO_A)} = \frac{AO_A}{\sin(90^{\circ} - (\sphericalangle B - \sphericalangle A))} \quad (**)$$

From (*) and (**) we find

$$\frac{\sin(BOO_A)}{\sin(CAO_A)} = \frac{\cos(\sphericalangle C - \sphericalangle A)}{\cos(\sphericalangle B - \sphericalangle A)}$$

(***)

Similarly we obtain

$$\frac{\sin(ABO_B)}{\sin(CBO_B)} = \frac{\cos(\sphericalangle A - \sphericalangle B)}{\cos(\sphericalangle C - \sphericalangle B)}$$

(****)

$$\frac{\sin(BCO_C)}{\sin(ACO_C)} = \frac{\cos(\sphericalangle B - \sphericalangle C)}{\cos(\sphericalangle A - \sphericalangle C)}$$

(*****)

From relations (***), (****), (*****) and the from the trigonometrically variation of Ceva's theorem it results that $AO_A$, $BO_B$, $CO_C$ are concurrent. The concurrence point is noted $K_O$, and it is called the Coşniţă point. Therefore the Coşniţă point is the homology center of triangle and of $ABC$ Coşniţă triangle.

Note

The name of Coşniţă point has been introduced by Rigby in 1997.

**Observation 23**

The theorem can be similarly proved in the case of an obtuse triangle.

**Remark 19**

Triangle $ABC$ and Coşniţă triangle $O_AO_BO_C$ being homological, have the homology axis the Coşniţă line.

**Theorem 15**

The Coşniţă point of triangle $ABC$ is the isogonal conjugate of the center of the circle of nine points associated to triangle $ABC$.

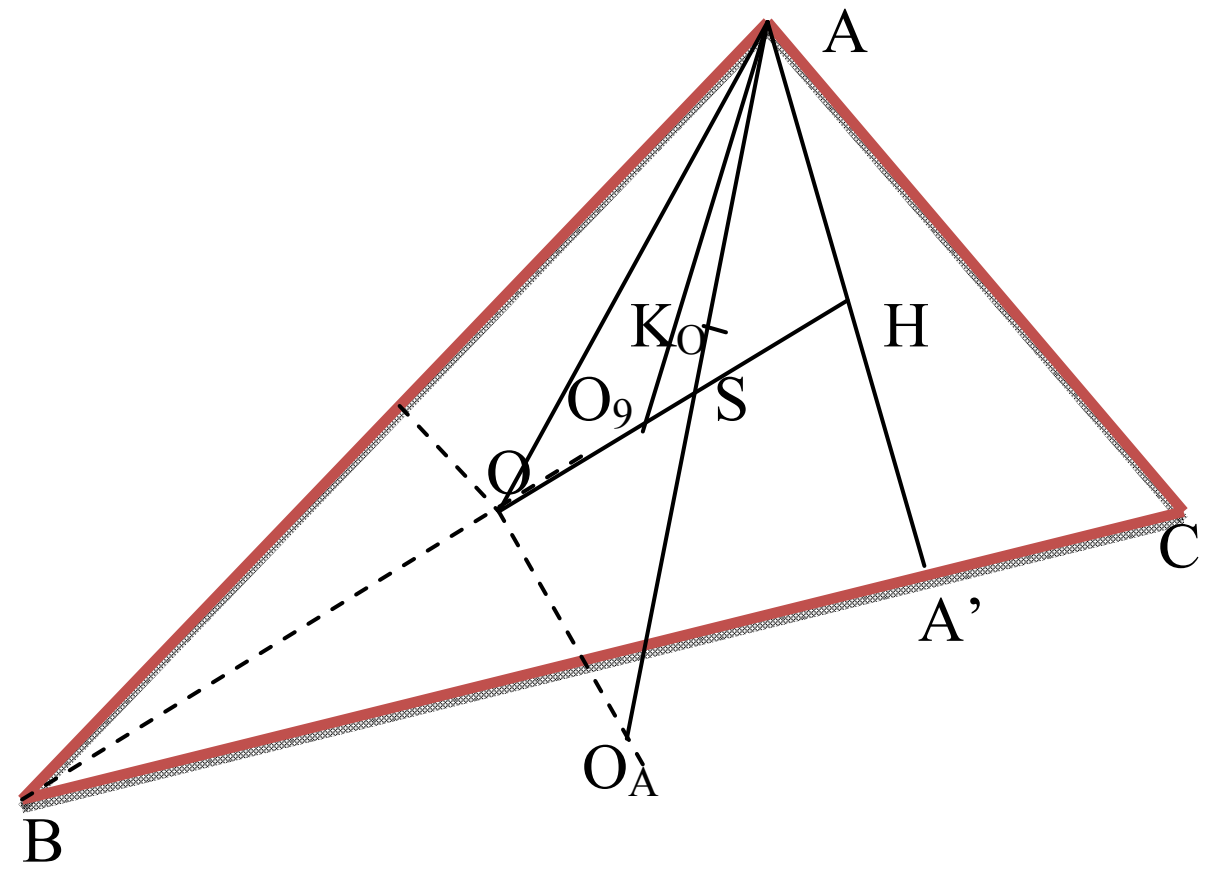

Fig. 32

**Proof**

It is known that the center of the circle of nine points, noted $O_9$ in figure 32 is the middle of segment $OH$, $H$ being the orthocenter of triangle $ABC$.

We note $\{S\}=OO_A \cap OH$.

We'll prove that $\sphericalangle OAK_O \equiv \sphericalangle HAO_9$ which is equivalent with proving that $AS$ is symmedian in triangle $OAH$. This is reduced to prove that $\frac{OS}{SH}=\frac{OA^2}{AH^2}$.

Because $OO_A \parallel AH$ we have that $\triangle OAS \sim \triangle HSA$, it results that $\frac{OS}{SM}=\frac{OO_A}{AH}$.

From the sinuses' theorem applied in triangle $BOC$ we find that $OO_A=\frac{R}{2(\cos A)}$. It is known that in a triangle $AH=2R(\cos A)$, it results that $\frac{OO_A}{AH}=\frac{R^2}{AH^2}$; therefore $AS$ is the isogonal of the median $AO_9$; similarly we prove that $BK_O$ is a symmedian in triangle $BOH$, consequently $K_O$ is the isogonal conjugate of the center of the circle of the nine points.

**Theorem 16** (generalization of the Coşniţă's theorem)

Let $P$ a point in the plan of triangle $ABC$, not on the circumscribed circle or on the triangle's sides; $A'B'C'$ the pedal triangle of $P$ and the points $A_1, B_1, C_1$ such that

$$\overrightarrow{PA'}\cdot\overrightarrow{PA_1}=\overrightarrow{PB'}\cdot\overrightarrow{PB_1}=\overrightarrow{PC'}\cdot\overrightarrow{PC_1}=k\ ,.\ k\in Q^*$$

Then the triangles $ABC$ and $A_1B_1C_1$ are homological.

**Proof**

Let $\alpha, \beta, \gamma$ the barycentric coordinates of $P$. From $\alpha = aria(\Delta PBC)$ we have $PA' = \frac{2\alpha}{a}$, and from $PA' \cdot PA_1 = k$ we find $PA_1 = \frac{ak}{2\alpha}$ (we considered $P$ in the interior of triangle $ABC$, see figure 33). We note $D$ and respectively $P_1$ the orthogonal projections of $A_1$ on $AD$ and $P$ on $A_1D$.

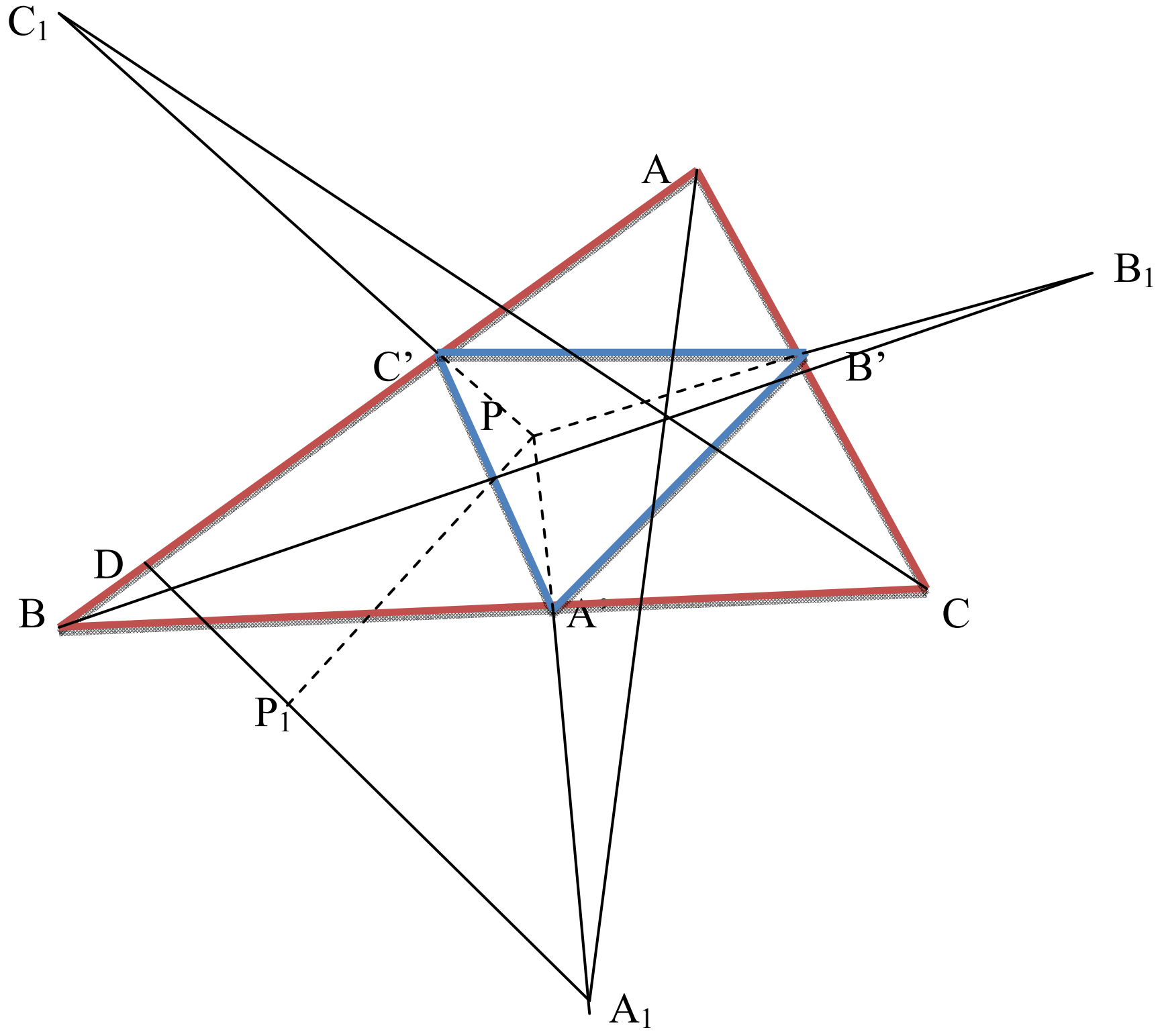

Fig.33

**Proof**

Because $\sphericalangle PA_1D \equiv \sphericalangle ABC$ (angles with perpendicular sides), we have

$$A_1P_1 = PA_1 \cdot \cos B = \frac{aK}{2\alpha}\cos B$$

From $\gamma = aria(\Delta PAB)$ it results $PC' = \frac{2\gamma}{c}$ and $A_1D = A_1P_1 + P_1D = \frac{aK\cos B}{2\alpha} + \frac{2\gamma}{c}$.

We note $A_1 = \alpha_1\beta_1\gamma_1$ and we have

$$\gamma_1 = \frac{kac\cos B + 4\gamma\alpha}{4\alpha},$$

$$\alpha_1 = aria(\Delta A_1 BC) = \frac{1}{2} BC \cdot A_1 A' = \frac{a^2 k - 4\alpha^2}{4\alpha}$$

And similarly as in the $\gamma_1$ computation we obtain

$$\beta_1 = \frac{abk \cos C + 4\alpha\beta}{4\alpha}$$

Or

$$A_1\left(\frac{a^2k - 4\alpha^2}{4\alpha}, \frac{abk\cos C + 4\alpha\beta}{4\alpha}, \frac{ack\cos B + 4\alpha\gamma}{4\alpha}\right)$$

Similarly,

$$B_1\left(\frac{abk\cos C + 4\alpha\beta}{4\beta}, \frac{b^2k - 4\beta^2}{4\beta}, \frac{cbk\cos A + 4\beta\gamma}{4\beta}\right),$$

$$C_1\left(\frac{abk\cos B + 4\alpha\gamma}{4\gamma}, \frac{cbk\cos A + 4\beta\gamma}{4\gamma}, c\frac{b^2k - 4\gamma^2}{4\gamma}\right)$$

Conform with [10] the Cevians $AA_1, BB_1, CC_1$ are concurrent if and only if $\alpha_2\beta_3\gamma_1 = \alpha_3\beta_1\gamma_2$. Because in our case this relation is verified, it results that the lines $AA_1$, $BB_1$, $CC_1$ are concurrent.

The theorem is proved in the same manner also in the case when $P$ is in the exterior of the triangle.

**Remark 20**

The center of homology of the triangles $ABC$ and $A_1B_1C_1$ has been named the generalized point of Coşniţă.

**Observation 24.1**

The conditions from the above theorem have the following geometrical interpretation: the points $A_1, B_1, C_1$ are situated on the perpendiculars from point $P$ on the triangle's sides and are the inverses of the points $A', B', C'$ in rapport to the circle in point $P$ and of radius $|k|$.

**Observation 24.2**

The Coşniţă's theorem is obtained in the particular case $P = O$ and $k = \frac{R^2}{2}$ ($O$ and $R$ are the center and respectively the radius of the circumscribed circle). Indeed, $OA' = R\cos A$, and with the sinuses' theorem applied in triangle $BOC$ we have $\frac{R}{\sin 2A} = 2OA_1$ ($A_1$ being the center of the circumscribed circle of triangle $BOC$ ). From where $OA' \cdot OA_1 = \frac{R^2}{2}$, similarly we find

$$OB' \cdot OB_1 = OC' \cdot OC_1 \frac{R^2}{2}.$$

**Observation 24.3**

It is easy to verify that if we consider $P = I$ (the center of the inscribed circle) and $k = r(r + a)$, $a > 0$ given, and $r$ is the radius of the inscribed circle in the Kariya point. Therefore the above theorem constitutes a generalization of the Kariya's theorem.

# Chapter 2

## Triplets of homological triangles

This chapter we prove of several theorems relative to the homology axes and to the homological centers of triplets of homological triangles.

The proved theorems will be applied to some of the mentioned triangles in the precedent sections, and also to other remarkable triangles which will be defined in this chapter.

### 2.1. Theorems relative to the triplets of homological triangles

**Definition 29**

The triplet $(T_1, T_2, T_3)$ is a triplet of homological triangles if the triangles $(T_1, T_2)$ are homological, the triangles $(T_2, T_3)$ are homological and the triangles $(T_1, T_3)$ are homological.

**Theorem 17**

Given the triplet of triangles $(T_1, T_2, T_3)$ such that $(T_1, T_2)$ are homological, $(T_1, T_3)$ are homological and their homological centers coincide, then

(i) $(T_1, T_2, T_3)$ is a triplet of homological triangles. The homological center of $(T_1, T_3)$ coincides with the center of the previous homologies.

(ii) The homological axes of the pairs of triangles from the triplet $(T_1, T_2, T_3)$ are concurrent, parallel, or coincide.

**Proof**

Let's consider $T_1$ triangle $A_1B_1C_1$, $T_2$ triangle $A_2B_2C_2$, and $T_3$ triangle $A_3B_3C_3$. (See figure 34).

We note $O$ the common homological center of triangles $(T_1, T_2)$ and $(T_2, T_3)$ such that

$$A_1A_2 \cap B_1B_2 \cap C_1C_2 = \{O\}$$

$$A_2A_3 \cap B_2B_3 \cap C_2C_3 = \{O\}$$

From these relations results without difficulty that

$$A_1A_3 \cap B_1B_3 \cap C_1C_3 = \{O\}$$

Consequently, $(T_1, T_3)$ are homological and the homology center is also $O$.

We consider the triangle formed by the intersection of the lines $A_1B_1$, $A_2B_2$, $A_3B_3$, noted in figure $PQR$ and the triangle formed by the intersections of the lines $B_1C_1$, $B_2C_2$, $B_3C_3$, noted $KLM$. We observe that

$$PR \cap KM = \{B_1\}$$

Also

$$PQ \cap KL = \{B_2\}$$

$$RQ \cap ML = \{B_3\}$$

Because $B_1, B_2, B_3$ are collinear from the Desargues' theorem we obtain that the triangles $PQR$ and $KLM$ are homological, therefore $PK, RM, QL$ are concurrent lines.

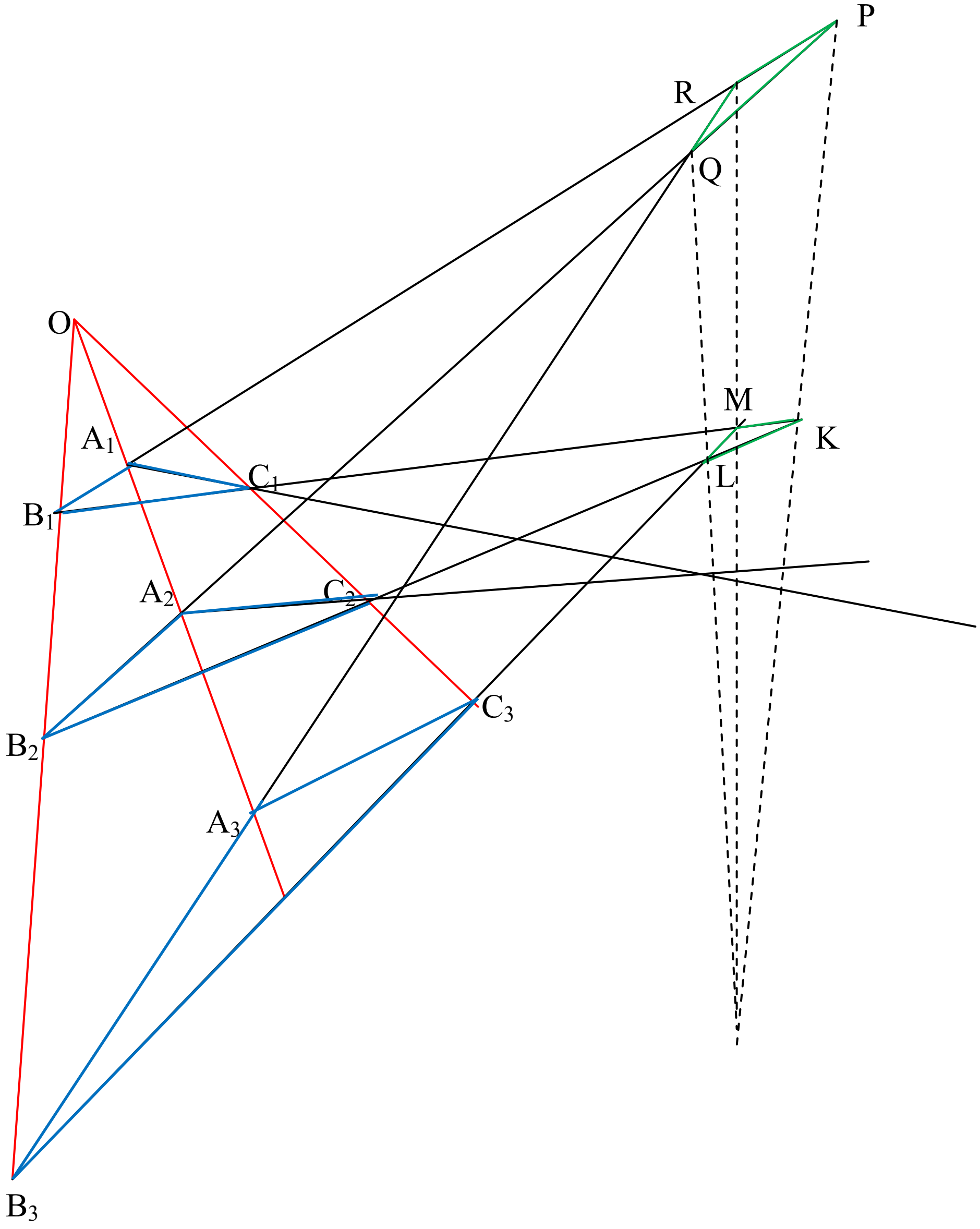

Fig. 34

The line $PK$ is the homology axis of triangles $(T_1, T_2)$, line $QL$ is the homology axis of triangles $(T_2, T_3)$ and $RM$ is the homology axes of triangles $(T_1, T_3)$. Because these lines are concurrent, we conclude that the theorem is proved.

**Remark 21**

a) We can prove this theorem using the space role: if we look at the figure as being a space figure, we can see that the planes $(T_1)$, $(T_2)$ share the line $PK$ and the planes $(T_1)$, $(T_3)$ share the line $QL$. If we note $\{O'\} = PK \cap LQ$ it results that $O'$ belongs to the planes $(T_1)$ and $(T_3)$, because these planes intersect by line $RM$, we find that $O'$ belongs to this line.

The lines $PK, RM, QL$ are homological axes of the considered pairs of triangles, therefore we conclude that these are concurrent in a point $O'$.

b) The theorem's proof is not valid when the triangles $PQR$ and $KLM$ don't exist.

A situation of this type can be when the triangles $(T_1, T_2)$ are homological, the triangles $(T_2, T_3)$ are homological, and the triangles $(T_1, T_3)$ are homothetic. In this case considering the figure s in space we have the planes $(T_1)$ and $(T_3)$ are parallel and the plane $(T_2)$ will intersect them by two parallel lines (the homology axes of the triangles $(T_1, T_2)$ and $(T_1, T_3)$).

c) Another situation when the proof needs to be adjusted is when it is obtained that the given triangles have two by two the same homological axis.

The following is a way to justify this hypothesis. We'll consider in space three lines $d_1, d_2, d_3$ concurrent in the point $O$ and another line $d$ which does not pass through O. Through $d$ we draw three planes $\alpha, \beta, \gamma$ which will intersect $d_1, d_2, d_3$ respectively in the points $A_1, B_1, C_1$; $A_2, B_2, C_2$; $A_3, B_3, C_3$. The three triangles $A_1B_1C_1$, $A_2B_2C_2$, $A_3B_3C_3$ are homological two by two, their homology center is $O$ and the common homology axis is $d$.

**Theorem 18**

Given the triplet of triangles $(T_1, T_2, T_3)$ such that $(T_1, T_2)$ are homological, $(T_2, T_3)$ are homological, and the two homology having the same homological axis, then:

i) The triplet $(T_1, T_2, T_3)$ is homological. The homology axis of triangles $(T_1, T_3)$ coincide with the previous homological axis.

ii) The homological centers of the triangles $(T_1, T_2)$, $(T_2, T_3)$ and $(T_1, T_3)$ are collinear or coincide.

**Proof.**

If

$$T_1 = A_1B_1C_1$$
$$T_2 = A_2B_2C_2$$
$$T_3 = A_3B_3C_3$$

and $M, N, P$ is the common homological axis of triangles $(T_1, T_2)$ and $(T_2, T_3)$ it results that

$$\{M\} = B_1C_1 \cap B_2C_2$$
$$\{N\} = A_1C_1 \cap A_2C_2$$
$$\{P\} = A_1B_1 \cap A_2B_2$$

and

$$\{M\} = B_2C_2 \cap B_3C_3$$
$$\{N\} = A_2C_2 \cap A_3C_3$$
$$\{P\} = A_2B_2 \cap A_3B_3$$

From these relations we find that

$$\{M\} = B_1C_1 \cap B_3C_3$$
$$\{N\} = A_1C_1 \cap A_3C_3$$
$$\{P\} = A_1B_1 \cap A_3B_3$$

which shows that the line $M, N, P$ is homological axis also for triangles $(T_1, T_3)$

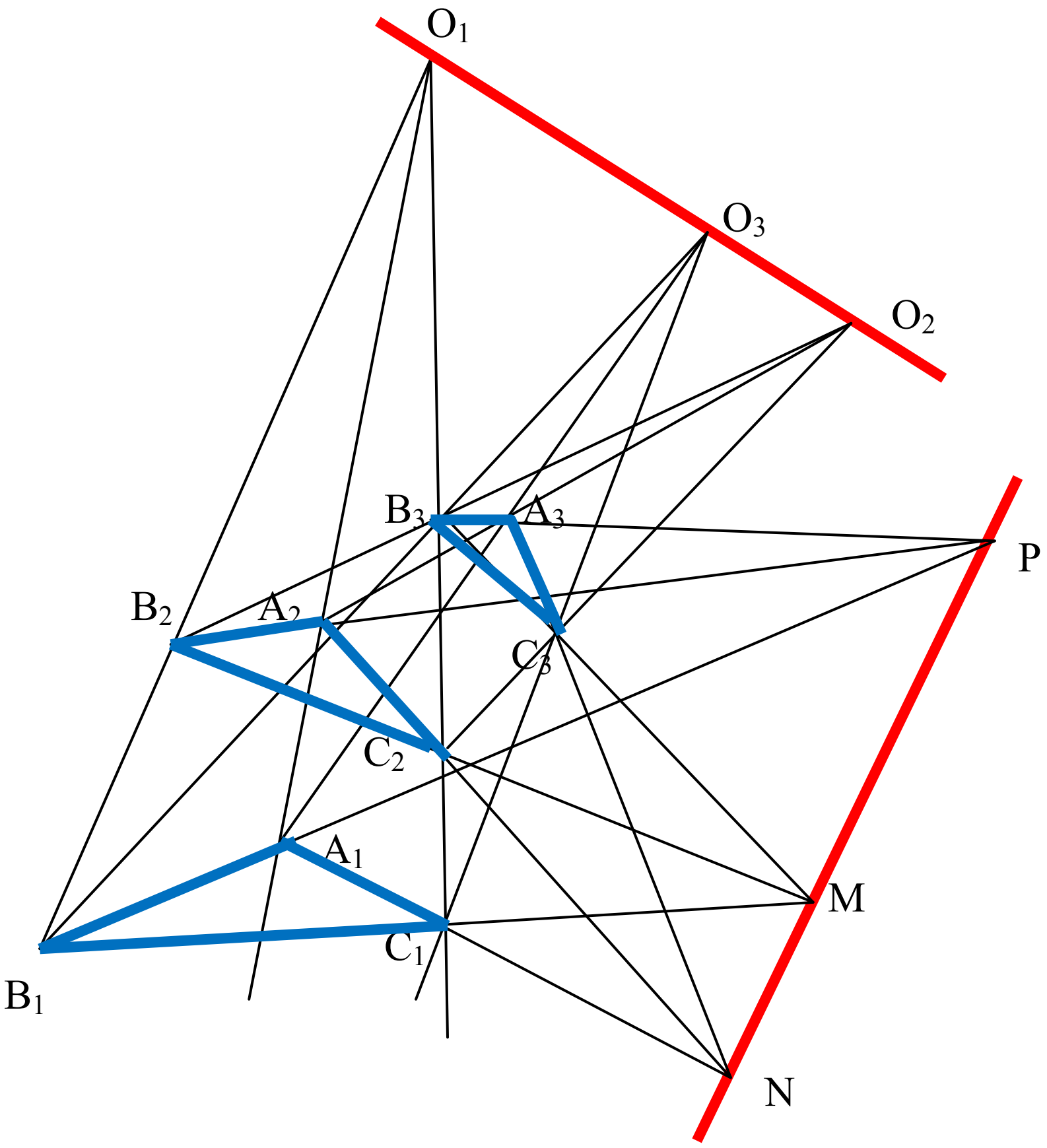

Fig. 35

In figure 35 we noted with $O_1$ the homology center of triangles $(T_1, T_2)$, with $O_2$ the homology center of triangles $(T_2, T_3)$ and with $O_3$ the homology center of triangles $(T_1, T_3)$. Considering triangles $A_1A_2A_3$ and $B_1B_2B_3$ we see that $A_1B_1 \cap A_2B_2 \cap A_3B_3 = \{P\}$, therefore these triangles are homological, their homology center being the point $P$. Their homological axis is determined by the points

Consequently the points $O_1, O_2, O_3$ are collinear.
Verify formula

$$\{O_1\} = A_1A_2 \cap B_1B_2$$
$$\{O_2\} = A_2A_3 \cap B_1B_3$$
$$\{O_3\} = A_1A_3 \cap B_1B_2$$

**Theorem 19** (the reciprocal of theorem 18)

If the triplet of triangles $(T_1, T_2, T_3)$ is homological and the centers of the homologies $(T_1, T_2)$, $(T_2, T_3)$, and $(T_1, T_3)$ are collinear, then these homologies have the same homology axes.

**Proof**

Let $O_1, O_2, O_3$ the collinear homology centers (see fig 35). We'll consider triangles $B_1B_2B_3$, $C_1C_2C_3$ and we observe that these have as homology axes the line that contains the points $O_1, O_2, O_3$. Indeed, $\{O_1\} = B_1B_2 \cap C$, $\{O_2\} = B_2B_3 \cap C_2C_3$, and $\{O_3\} = B_1B_3 \cap C_1C_3$, it results that these triangles have as homology center the point $M$ ($B_1C_1 \cap B_2C_2 \cap B_3C_3 = \{M\}$). Similarly, the triangles $A_1A_2A_3$, $C_1C_2C_3$ have as homology axes the line $O_1O_2O_3$, therefore as homology center the point $M$; the triangles $A_1A_2A_3$, $B_1B_2B_3$ are homological with the homology axis $O_1O_2O_3$ and of center $P$. Theorem14 implies the collinearity of the points $M, N, P$, therefore the theorem is proved.

**Theorem 20** (Véronèse)

Two triangles $A_1B_1C_1$, $A_2B_2C_2$ are homological and

$$\{A_3\} = B_1C_2 \cap B_2C_1$$
$$\{B_3\} = A_1C_2 \cap C_1A_2$$
$$\{C_3\} = A_1B_2 \cap B_1A_2$$

then the triplet $(A_1B_1C_1, A_2B_2C_2, A_3B_3C_3)$ is homological and the three homologies have the homological centers collinear.

**Proof**

Let $O_1$ the homology center of triangles $(T_1, T_2)$, where

$$T_1 = A_1B_1C_1$$
$$T_2 = A_2B_2C_2$$

(see fig. 36) and $A'B'C'$ their homology axis.

We observe that $O_1$ is a homological center for triangles $A_1B_1C_2$, $A_2B_2C_1$, and their homological axis is $C'A_3B_3$.

Also $O_1$ is the homological center for triangles $B_1C_1A_2$, $B_2C_2A_1$; these triangles have as homological axis the line $A', B_3, C_3$.

Similarly, we obtain that the points $B', A_3, C_3$ are collinear, being on the homological axis of triangles $C_1A_1B_2$, $C_2A_2B_1$.

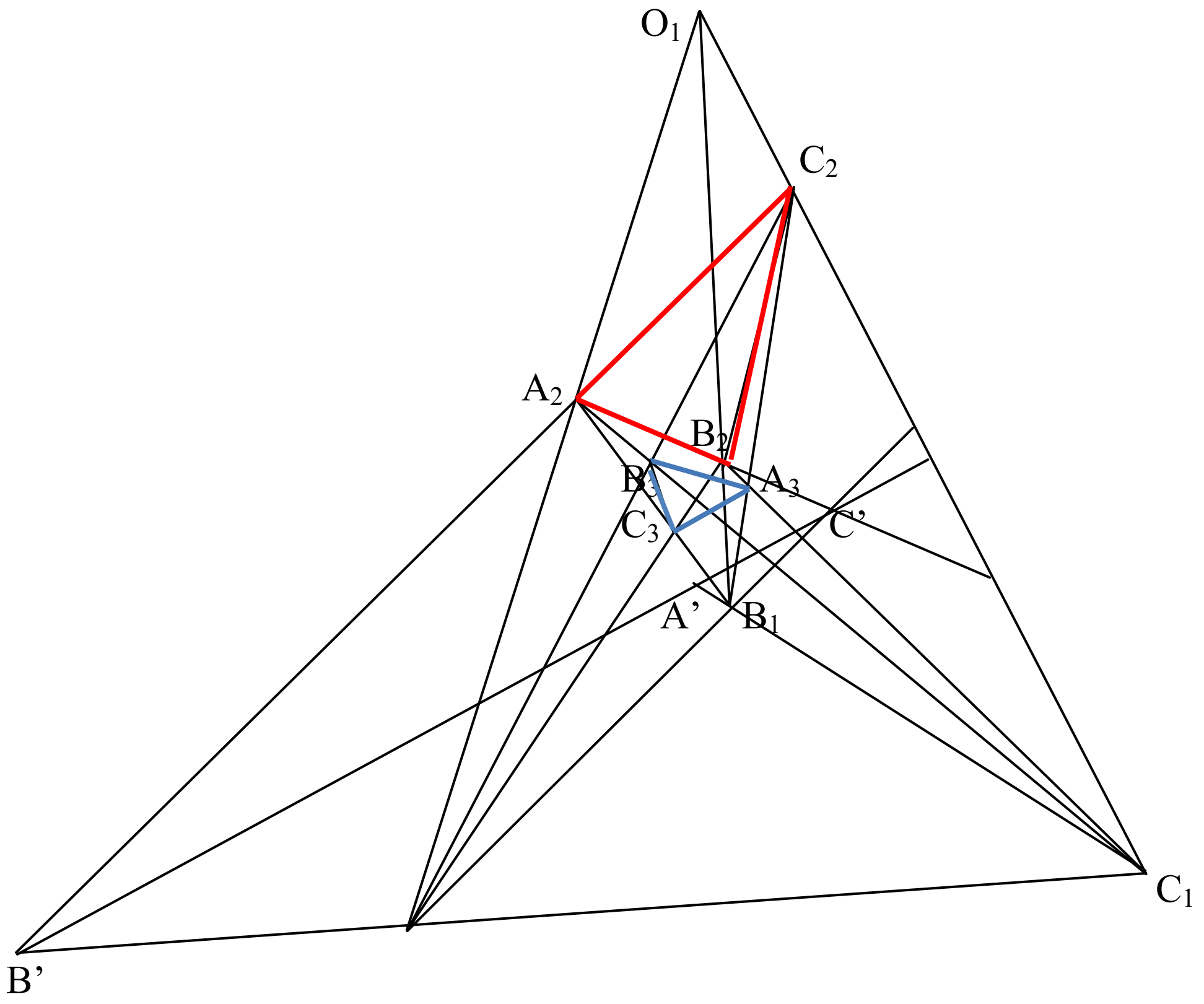

Fig. 36

The triplets of collinear points $(C', A_3, B_3)$, $(B', A_3, C_3)$ and $(A', B_3, C_3)$ show that the triangle $T_3 = A_3B_3C_3$ is homological with $T$ and $T_2$.

The triplet $(T_1, T_2, T_3)$ is homological and their common homological axis is $A', B', C'$. In conformity to theorem 15, it result that their homological centers are collinear.

## 2.2. A remarkable triplet of homological triangles

**Definition 30**

A first Brocard triangle of given triangle is the triangle determined by the projections of the symmedian center of the given triangle on its medians.

**Observation 25**

In figure 37 the first Brocard's triangle of triangle $ABC$ has been noted $A_1B_1C_1$.

**Definition 31**

In a given triangle $ABC$ there exist the points $\Omega$ and $\Omega'$ and an angle $\omega$ such that

$$m(\sphericalangle \Omega AB) = m(\sphericalangle \Omega BC) = m(\sphericalangle \Omega CA) = \omega; \quad m(\sphericalangle \Omega' BA) = m(\sphericalangle \Omega' CA) = m(\sphericalangle \Omega' AB) = \omega$$

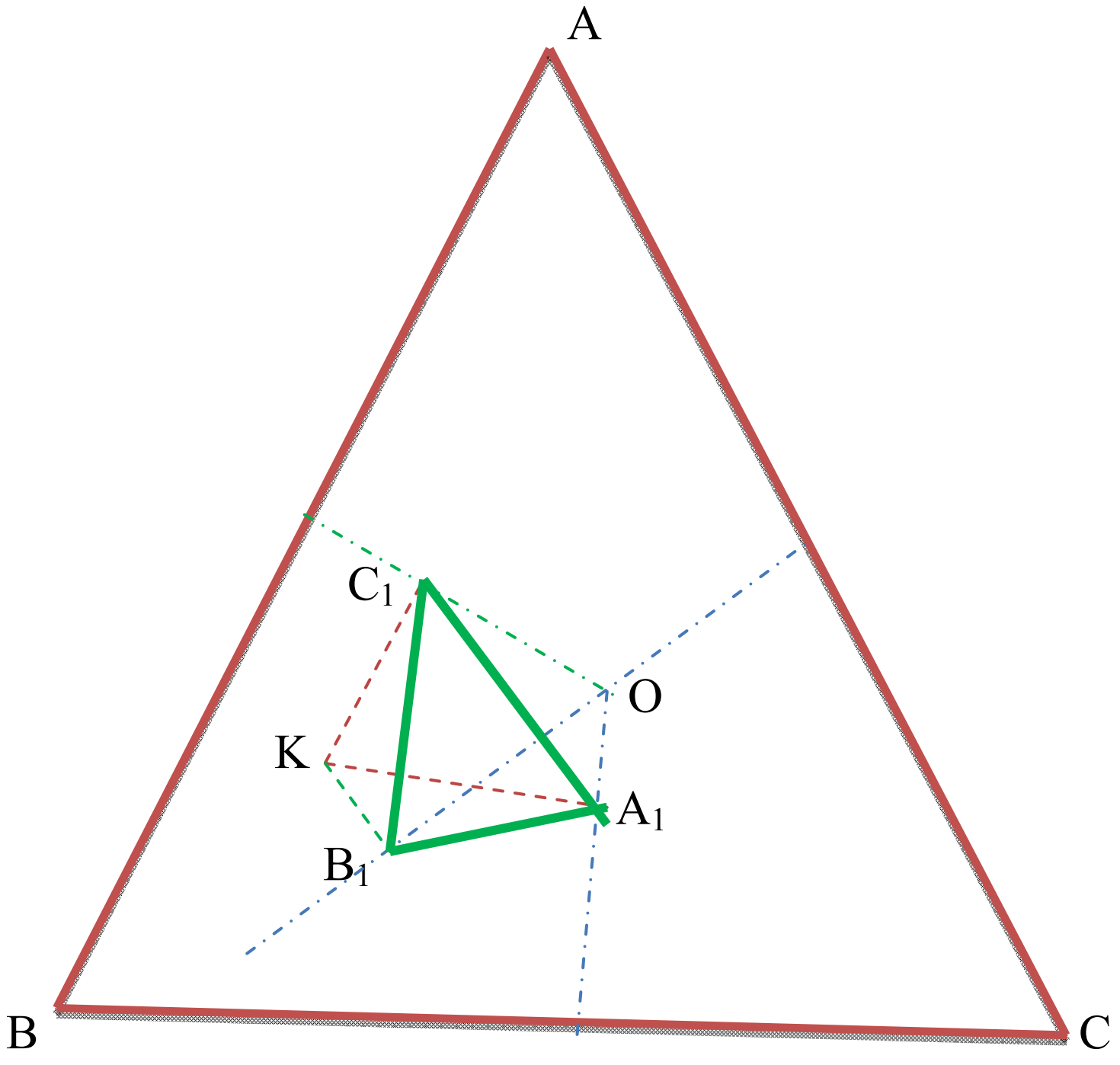

Fig. 37

The points $\Omega$ and $\Omega'$ are called the first, respectively the second Brocard's point, and $\omega$ is called the Brocard angle.

**Definition 32**

An adjunct circle of a triangle is the circle which passes by two vertexes of the given triangle and it is tangent in one of these vertexes to the side that contains it.

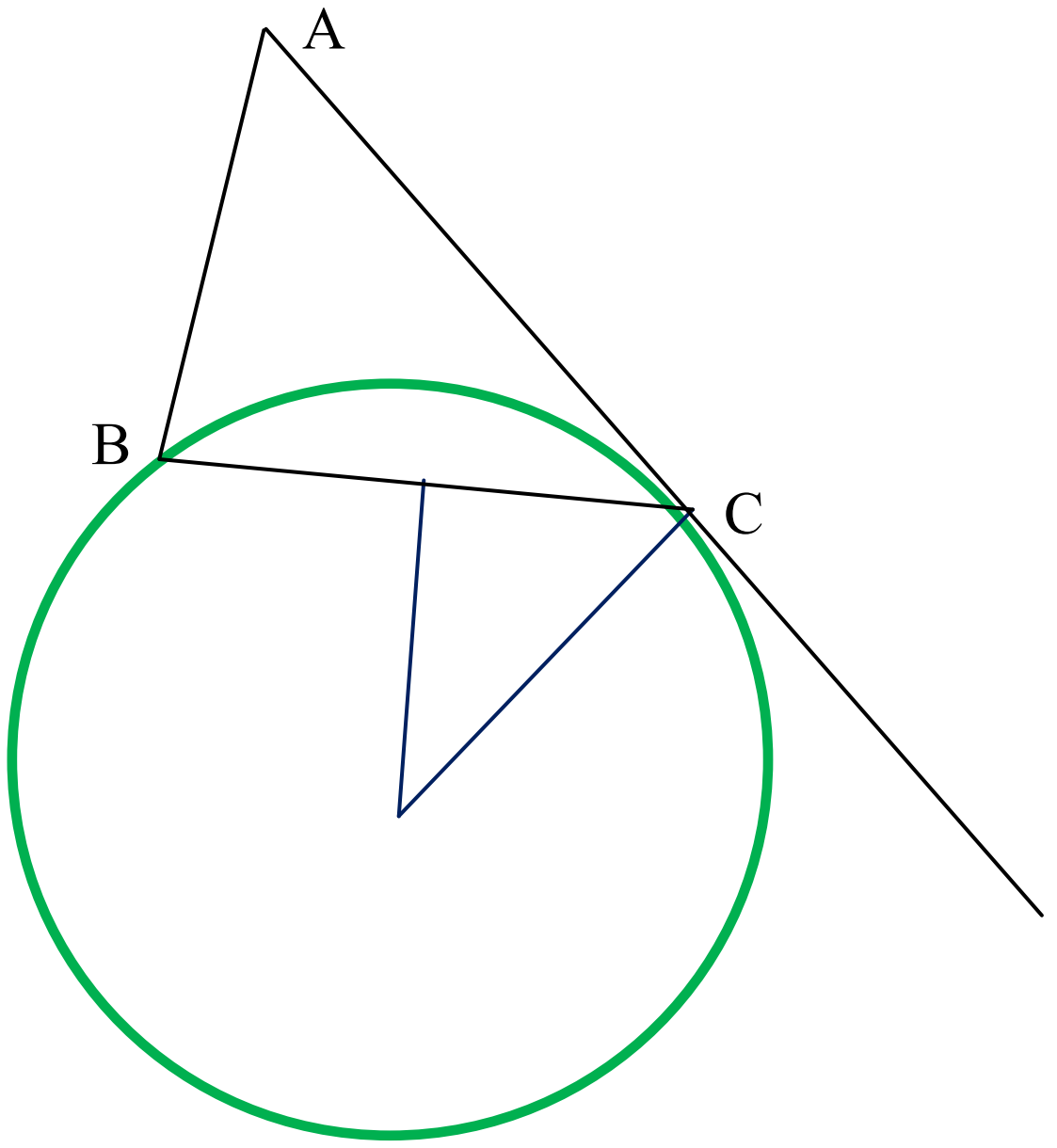

Fig 38

**Observation 26**

In figure 38 it is represented the adjunct circle which passes through $B$ and $C$ and is tangent in $C$ to the side $AC$. Will note this circle $\widehat{BC}$. To the given triangle $ABC$ corresponds six adjunct circles.

**Proposition 26**

The adjoin circles $\widehat{AB},\widehat{BC},\widehat{CA}$ of triangle $ABC$ intersect in Brocard's point $\Omega$.

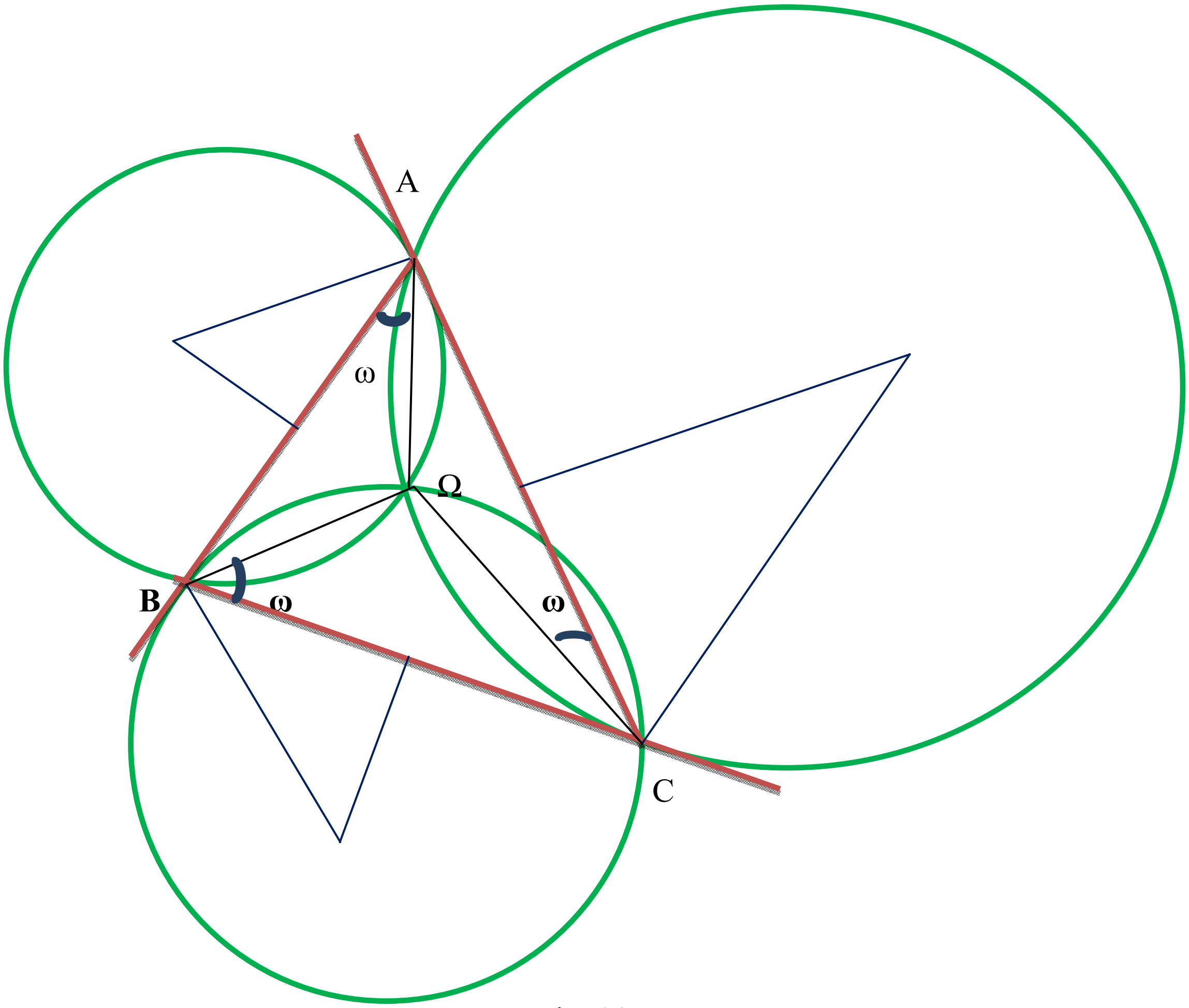

Fig. 39

**Proof**

Let $\Omega$ the second point of intersection of the circles $\widehat{AB},\widehat{BC}$. Then we have that

$$\measuredangle \Omega AB \equiv \measuredangle \Omega BC.$$

Because $B\widehat{C}$ is tangent in $B$ to the side $AC$ we have also the relation

$$\measuredangle \Omega BC \equiv \measuredangle \Omega CA$$

These imply that

$$\measuredangle \Omega CA \equiv \measuredangle \Omega AB$$

And this relation shows that the circumscribed circle to triangle $\Omega CA$ is tangent in the point $A$ to the side $AB$, which means that the adjunct circle $\widehat{CA}$ passes through the Bacard's point $\Omega$.

**Remark 22**

Similarly we prove that the adjunct circles $\overset{\frown}{BA}, \overset{\frown}{AC}, \overset{\frown}{CB}$ intersect in the second Bacard's point $\Omega'$. Bacard's points $\Omega$ and $\Omega'$ are isogonal concurrent.

**Proposition 27**

If $ABC$ is a triangle and $\omega$ is Bacard's angle, then

$$ctg\,\omega = ctgA + ctgB + ctgC \tag{60}$$

Proof

Applying the sinus' theorem in triangles $A\Omega C$, $B\Omega C$ we obtain:

$$\frac{C\Omega}{\sin(A-\omega)} = \frac{AC}{\sin A\Omega C} \tag{61}$$

$$\frac{C\Omega}{\sin\omega} = \frac{BC}{\sin B\Omega C} \tag{62}$$

Because

$$m(\sphericalangle A\Omega C) = 180^\circ - A$$

and

$$m(\sphericalangle B\Omega C) = 180^\circ - C$$

From relations (61) and (62) we find

$$\frac{\sin\omega}{\sin(A-\omega)} = \frac{AC}{BC}\cdot\frac{\sin C}{\sin A} \tag{63}$$

And from here

$$\sin(A-\omega) = \frac{a}{b}\cdot\frac{\sin A}{\sin C}\cdot\sin\omega \tag{64}$$

From the sinus' theorem in triangle $ABC$ we have that $\frac{a}{b} = \frac{\sin A}{\sin B}$ and re-writing (64) we have

$$\sin(A-\omega) = \frac{\sin^2 A\cdot\sin\omega}{\sin B\cdot\sin C} \tag{65}$$

Furthermore

$$\sin(A-\omega) = \sin A\cdot\cos\omega - \sin\omega\cos A$$

And

$$\sin A\cdot\cos\omega - \sin\omega\cos A = \frac{\sin^2 A}{\sin B}\cdot\frac{\sin\omega}{\sin C} \tag{66}$$

Dividing relation (66) by $\sin A\cdot\sin\omega$ and taking into account that $\sin A = \sin(B+C)$ and $\sin(B+C) = \sin B\sin C + \cos B\cdot\cos C$ we'll obtain the relation (60).

**Proposition 28**

In a triangle $ABC$ takes place the following relation:

$$ctg\,\omega = \frac{a^2+b^2+c^2}{4s} \tag{67}$$

**Proof**

If $H_2$ is the projection of the vertex $B$ on the side $AC$, we have:

$$ctgA = \frac{AH_2}{BH_2} = \frac{b \cdot c \cdot \cos A}{2s} \quad (68)$$

From the cosin's theorem in the triangle $ABC$ we retain

$$2bc \cos A = b^2 + c^2 - a^2 \quad (69)$$

Replacing in (68) we obtain

$$ctgA = \frac{b^2 + c^2 - a^2}{4s} \quad (70)$$

Considering the relation (60) and those similar to relation (70) we obtain relation (67)

**Definition 33**

It is called the Brocard's circle of a triangle the circle of who's diameter is determined by the symmedian center and the center of the circumscribed circle of the given triangle.

**Proposition 29**

The first Brocard's triangle of a triangle is similar with the given triangle.

**Proof**

Because $KA_1 \parallel BC$ and $OA' \perp BC$ it results that $m(\sphericalangle KA_1O) = 90^\circ$ (see figure 37). Similarly $m(\sphericalangle KB_1O) = m(\sphericalangle KC_1O) = 90^\circ$ and therefore the first Brocard's triangle is inscribed in the Brocar's circle.

Because $m(\sphericalangle A_1OC_1) = 180^\circ - B$ and the points $A_1, B_1, C_1, O$ are con-cyclic, it results that $\sphericalangle A_1B_1C_1, O = \sphericalangle B$, similarly from $m\sphericalangle(B'OC') = 180^\circ - A$ it results that $m(\sphericalangle B_1OC_1) = m(A')$ but $\sphericalangle B_1OC_1 \equiv \sphericalangle B_1A_1C_1$, therefore $\sphericalangle B_1A_1C_1 \equiv \sphericalangle A$, therefore triangle $A_1B_1C_1$ is similar to the given triangle $ABC$.

**Proposition 30**

If in a triangle $ABC$ we note $K_1, K_2, K_3$ the orthogonal projections of the symmedian center $K$ on the triangle's sides, then the following relation takes place:

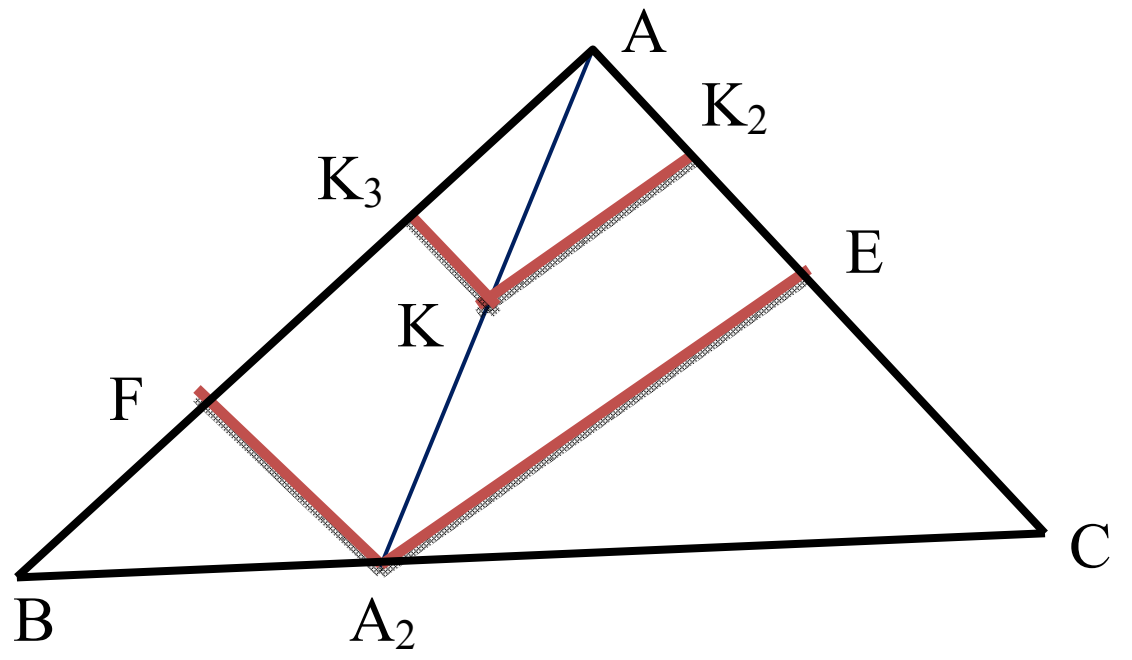

Fig. 40

$$\frac{KK_1}{a} = \frac{KK_2}{b} = \frac{KK_3}{c} = \frac{1}{2} tg\omega \qquad (71)$$

**Proof**

In triangle $ABC$ we construct $AA_2$ the symmedian from the vertex $A$ (see figure 20).

We have $\frac{BA_2}{CA_2} = \frac{c^2}{b^2}$, on the other side

$$\frac{BA_2}{CA_2} = \frac{Aria\Delta BAA_2}{Aria\Delta CAA_2} \qquad (72)$$

Also

$$\frac{Aria\Delta BAA_2}{Aria\Delta CAA_2} = \frac{AB \cdot A_2F}{AC \cdot A_2E} \qquad (73)$$

Where $E$ and $F$ are the projections of $A_2$ on $AC$ and $AB$.

It results that

$$\frac{A_2F}{A_2E} = \frac{c}{b} \qquad (74)$$

From the similarity of triangles $AKK_3$, $AA_2F$ and $AKK_2$, $AA_2E$ we have

$$\frac{KK_3}{KK_2} = \frac{A_2F}{A_2E} \qquad (75)$$

Taking into account (74), we find that $\frac{KK_2}{a} = \frac{KK_2}{c}$, and similarly $\frac{KK_1}{a} = \frac{KK_2}{b}$; consequently we obtain:

(76)

The relation (76) is equivalent to

$$\frac{aKK_1}{a^2} = \frac{bKK_2}{b^2} = \frac{cKK_3}{c^2} = \frac{aKK_1 + bKK_2 + cKK_3}{a^2 + b^2 + c^2} \qquad (77)$$

Because $aKK_1 + bKK_2 + cKK_3 = 2Aria\Delta ABC = 2s$, it result that

$$\frac{KK_1}{a} = \frac{KK_2}{b} = \frac{KK_3}{c} = \frac{2s}{a^2 + b^2 + c^2} \qquad (78)$$

We proved the relation (67) $ctgA = \frac{b^2 + c^2 - a^2}{4s}$ and from this and (78) we'll obtain the relation (71).

**Remark 23**

Because $KK_1 = A_1A'$ ($A'$ being the projection of the vertex $A_1$ of the first Brocard's triangle on $BC$), we find that

$$m(\sphericalangle A_1BC) = m(\sphericalangle A_1CB) = \omega \qquad (79)$$

also

$$m(\sphericalangle C_1AC) = m(\sphericalangle B_1CA) = \omega \qquad (80)$$

and

$$m(\sphericalangle C_1AB) = m(\sphericalangle B_1BA) = \omega \quad (81)$$

**Theorem 21** (J. Neuberg – 1886)

The triangle $ABC$ and its Bacard's first triangle are homological. The homological center is the isotonic conjugate point of the symmedian center.

**Proof**

Let $H_1$ be the projection of vertex $A$ on $BC$ and $A_1^{'}$ the intersection of $AA_1$ with $BC$, and $A'$ the middle of the side $(BC)$.

We found that $m(\sphericalangle A_1BC) = \omega$, therefore $A_1A' = \frac{a}{2}tg\omega$, we have also $AH_1 = \frac{2s}{a}$.

From the similarity of triangles $A_1A'A_2^{'}$ and $AH_1A_2^{'}$ we find $\frac{A_2^{'}A'}{A_2^{'}H_1} = \frac{A_1A'}{AH_1}$, that is:

$$\frac{A_2^{'}A'}{A_2^{'}H_1} = \frac{a^2tg\omega}{4s} \quad (82)$$

We observe that $A_2^{'}A' = \frac{a}{2} - BA_2^{'}$, $A_2^{'}H_1 = BH_1 - BA_2^{'} = c\cos B - BA_2^{'}$.

Getting back to (82), we obtain $\frac{\frac{a}{2} - BA_2^{'}}{c\cos B - BA_2^{'}} = \frac{a^2tg\omega}{4s}$ and from here

$$\frac{a - 2BA_2^{'}}{2c\cos B - a} = \frac{a^2tg\omega}{4s - a^2tg\omega} \quad (83)$$

From (83) taking into account that $tg\omega = \frac{4s}{a^2+b^2+c^2}$ and $2ac\cos B = a^2 + c^2 - b^2$ we find that

$$BA_2^{'} = \frac{a \cdot b^2}{b^2+c^2} \quad (84)$$

then we obtain:

$$BA_2^{'} = \frac{a \cdot b^2}{b^2+c^2} \quad (85)$$

therefore

$$\frac{BA_2^{'}}{CA_2^{'}} = \frac{b^2}{c^2} \quad (86)$$

We note $BA_1 \cap AC = \{B_2^{'}\}$, $CA_1 \cap AB = \{C_2^{'}\}$

then

$$\frac{CB_2^{'}}{AB_2^{'}} = \frac{c^2}{a^2} \quad (87)$$

And

$$\frac{AC_2^{'}}{BC_2^{'}} = \frac{a^2}{b^2} \quad (88)$$

The last three relations along with the reciprocal of Ceva's theorem show that the Cevians
$AA_1, BB_1, CC_1$ are concurrent.

We showed that if $AA_1'$ is symmedian in triangle $ABC$ then $\frac{BA_1'}{CA_1'} = \frac{c^2}{b^2}$ (33). This relation and the relation:

$$\frac{CA_2'}{BA_2'} = \frac{c^2}{b^2} \quad (89)$$

lead us to the equality: $CA_2' = BA_2'$, which shows that the Cevian $AA_1$ is the isotonic of the symmedian from the vertex $A$ of triangle $ABC$; the property is true also for the Cevian $BB_1, CC_1$ and therefore the concurrence point of the Cevian $AA_1, BB_1, CC_1$ is the isotonic conjugate of the symmedian center $K$ of the triangle $ABC$.

In some publications this point is noted by Ω" and is called the third Brocard's point. We will use also this naming convention further on.

**Theorem 22**

The perpendicular constructed from the vertexes $A, B, C$ respectively on the sides $A_1C_1, C_1A_1, A_1B_1$ of the Brocard's first triangle of a given triangle $ABC$ intersect in a point $T$ which belongs to the triangle's circumscribed circle.

Proof

We'll note with $T$ the intersection of the perpendicular constructed from $B$ on $A_1C_1$

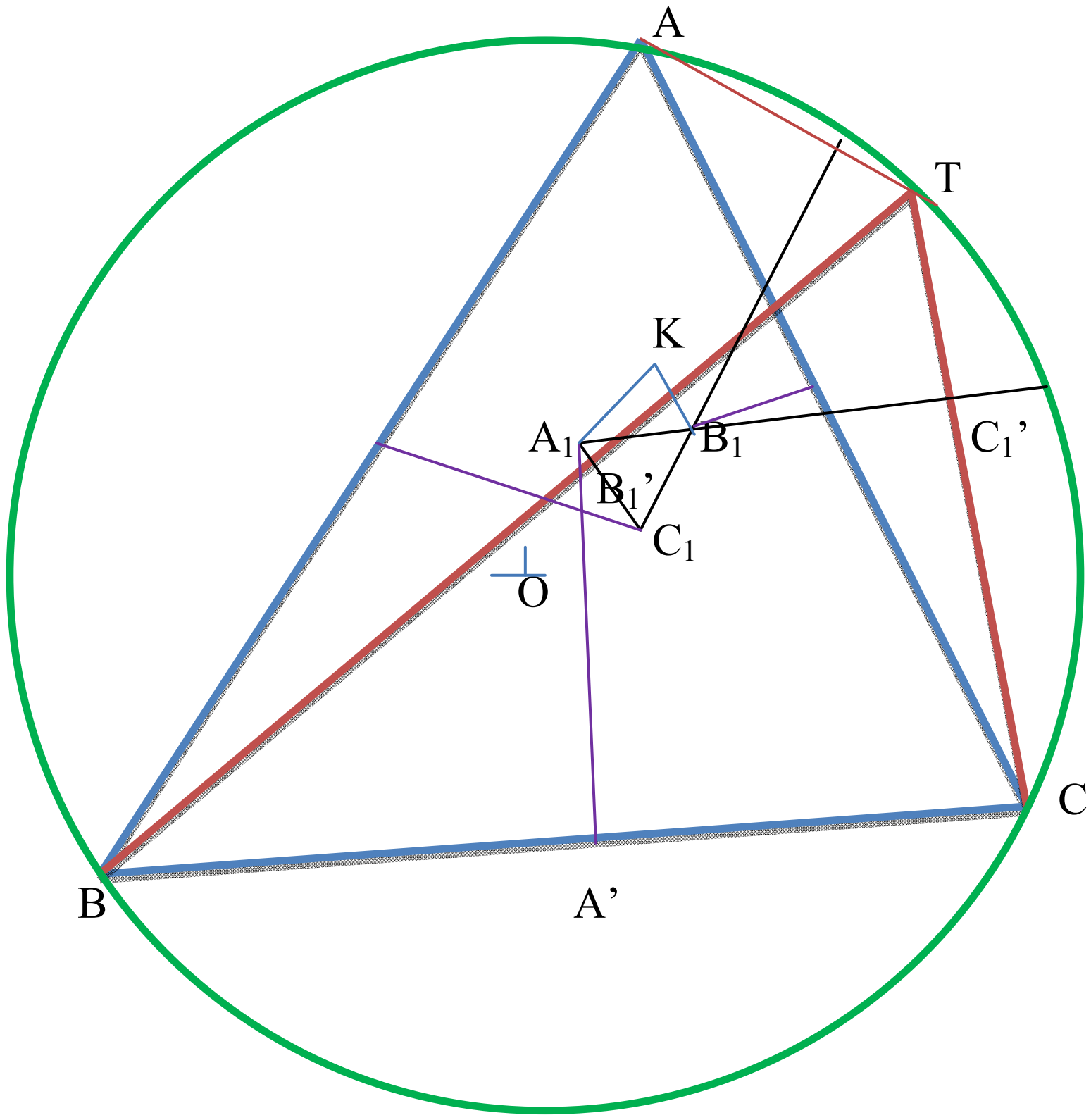

Fig 41

with the perpendicular constructed from $C$ on $A_1B_1$ and let $\{B_1'\} = BT \cap A_1C_1$ and $\{C_1'\} = A_1B_1 \cap CT$ (see figure 41).

We have $m(\sphericalangle B_1'TC_1') = m(\sphericalangle C_1A_1B_1)$. But conform to the proposition $\sphericalangle C_1A_1B_1 \equiv \sphericalangle A$, it results $m(\sphericalangle BTC) = \sphericalangle A$, therefore $T$ belongs to the circumscribed circle of triangle $ABC$.

If $\{A_1'\} = B_1C_1 \cap AT$ let note $T'$ the intersection of the perpendicular constructed from $A$ on $B_1C$ with the perpendicular constructed from $B$ on $A_1C_1$ ; we observe that $m(\sphericalangle B_1'TA_1') = m(\sphericalangle A_1C_1B_1)$ , therefore $T'$ belongs to the circumscribed circle of the triangle $ABC$ The points $T, T'$ belong to the line $BB_1'$ and to the circumscribed circle of triangle $ABC$. It result that $T = T'$ and the proof is complete.

**Remark 24**

The point $T$ from the precedent theorem is called Tarry's point of triangle $ABC$.

Similar can be proved the following theorem:

**Theorem 23**

If through the vertexes $A, B, C$ of a triangle $ABC$ we construct parallels to the sides $B_1C_1, C_1A_1, A_1B_1$ of the first Brocard's triangle of the given triangle, then these parallels are concurrent in a point $S$ on the circumscribed circle of the given triangle.

**Remark 25**

The point $S$ from the previous theorem is called Steiner's point of the triangle $ABC$. It can be easily shown that the Stern's point and Tarry's point are diametric opposed.

**Definition 34**

Two triangles are called equi-brocardian if they have the same Brocard's angle.

**Proposition 31**

Two similar triangles are equi-brocardian

**Proof**

The proof of this proposition results from the relation

$$ctg\,\omega = ctgA + ctgB + ctgC$$

**Remark 26**

A given triangle and its firs Brocard's triangle are equi-brocardian triangle.

**Proposition 32**

If $ABC$ is a triangle and $\overset{\frown}{AB}, \overset{\frown}{AC}$, are its adjoin circles which intersect the second side $AC$ in $E$ respective $F$, then the triangles $BEC, BFC$ are equi-brocardian with the given triangle.

**Proof**

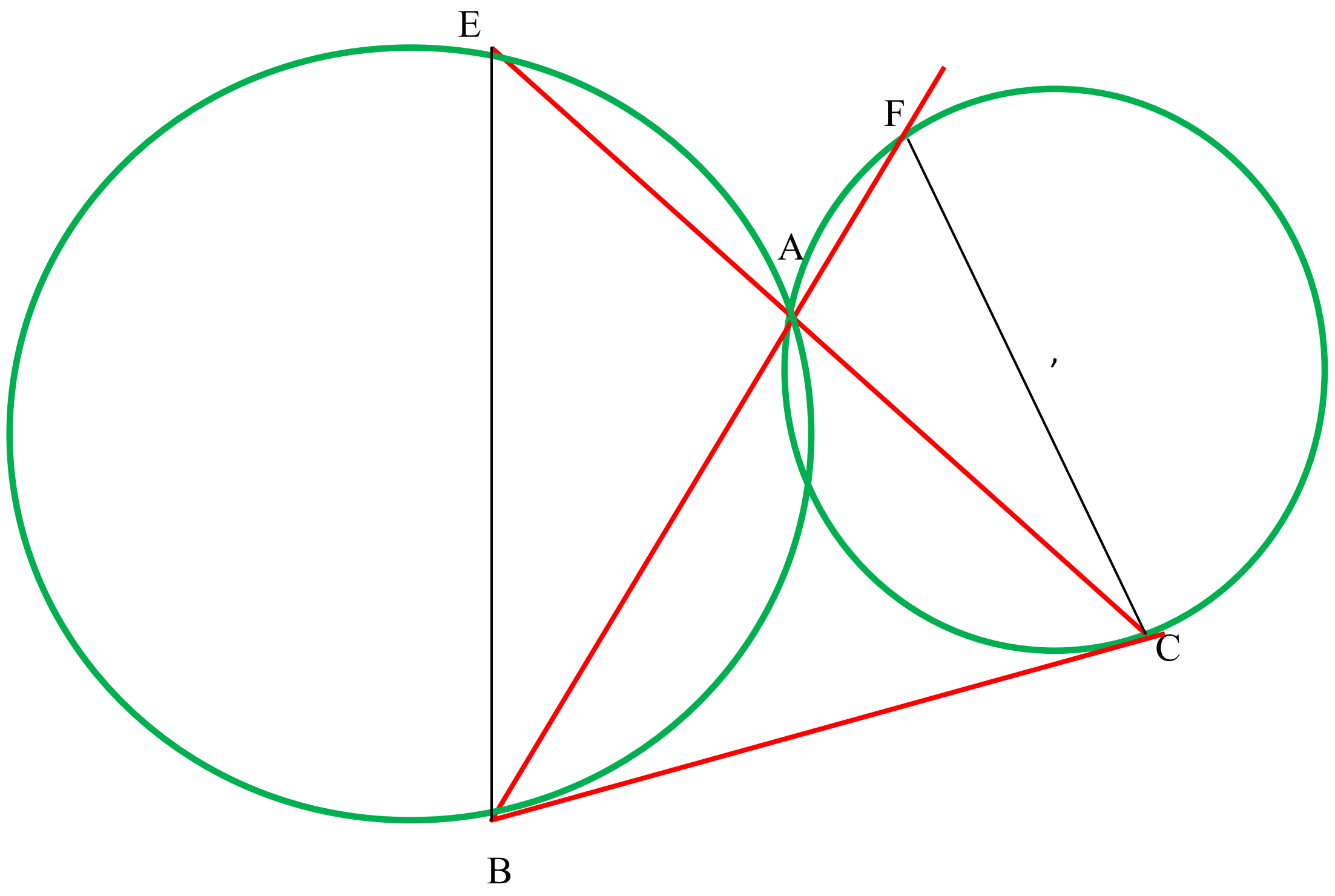

Fig. 42

The triangles $BEC, ABC$ have the angle $C$ in common and $\sphericalangle BEC \equiv \sphericalangle ABC$ because they extend the same arc $\overset{\frown}{AB}$ in the adjoined circle $\overset{\frown}{AB}$.

Therefore these triangles are similar therefore equi-brocardian.

Similarly we can show that triangle $ABC$ is similar to triangle $CBF$, therefore these are also equi-brocardian.

**Theorem 24**

The geometric locus of the points $M$ in a plane that are placed on the same side of the line $BC$ as the point $A$, which form with the vertexes $B$ and $C$ of triangle $ABC$ equi-brocard triangles with it is a circle of center $H_1$ such that $m(\sphericalangle BH_1C) = 2\omega$ and of a radius

$$\eta_1 = \frac{a}{2}\sqrt{ctg^2\omega - 3} \qquad (90)$$

**Proof**

From proposition 32 we find that the points $E, F$ belong evidently to the geometric locus that we seek. $A$ belongs to the geometric locus (see fig. 42).

We suppose that the geometric locus will be a circumscribed circle to triangle $AEF$ . . We'll compute the radius of this circle. We observe that $CA \cdot CE = CB^2$ and $BF \cdot BA = BC^2$, therefore the points $C, B$ have equal power in rapport to this circle, and the power is equal to $a^2$ .

From the precedent relations we find that

$$AE = \frac{|a^2 - b^2|}{b} \text{ and } AF = \frac{|c^2 - a^2|}{c}$$

Applying the sinus' theorem in the triangle $AEF$ we find

$$EF^2 = \frac{c^2\left(a^2-b^2\right)^2 + b^2\left(c^2-a^2\right)^2 + \left(b^2+c^2-a^2\right)\left(a^2-b^2\right)\left(c^2-a^2\right)}{b^2c^2}$$

The sinus' theorem in triangle $AEF$ gives

$$\eta_1 = \frac{EF}{2\sin A}$$

where $\eta_1$ is the radius of the circumscribed circle to triangle $AEF$ .

Because $2bc\sin A = 4s$ and $16s^2 = 2a^2b^2 + 2b^2c^2 + 2c^2a^2 - a^4 - b^4 - c^4$ and taking into account also the relation $ctg\,\omega = \frac{a^2+b^2+c^2}{4s}$ we find

$$\eta_1 = \frac{a}{2}\sqrt{ctg^2\omega - 3}$$

If we note $H_1$ the radius of the circumscribed circle to triangle $AEF$ , taking into account that the points $C, B$ powers in rapport to this circle are equal to $a^2$ and that this value is the square of the tangent constructed from $B$ respectively $C$ to this circle and further more equal to $H_1B^2 - \eta_1^2$ and to $H_1C^2 - \eta_1^2$ , we find that $H_1B = H_1C$ , which means that $H_1A' = \frac{a}{2}ctg\,\omega$ and from here we find that $m(\sphericalangle BH_1C) = 2\omega$.

Let's prove now that if $M$ is a point on the circle $\mathcal{C}(H_1, \eta_1)$ then triangle $MBC$ has the same Brocad angle $\omega$ as $ABC$. We'll note the Brocad's angle of triangle $MBC$ with $\omega'$, then

$$ctg\,\omega' = \frac{MB^2 + MC^2 + BC^2}{4s_{\Delta MBC}} \quad (91)$$

From the median theorem applied in triangle $MBC$ we find that

$$MB^2 + MC^2 = 2MA'^2 + \frac{a^2}{2} \quad (92)$$

$$AriaMBC = \frac{a \cdot MA' \cdot \cos MA'H_1}{2} \quad (93)$$

The cosin's theorem applied in triangle $MA'H_1$ gives

$$\eta_1^2 = MA'^2 + H_1A'^2 - 2MA' \cdot \cos\sphericalangle MA'H_1 \cdot H_1A' \quad (94)$$

We'll compute $4s_{\Delta MBC}$ taking into account relation (93) and substituting in this $2MA' \cdot \cos\sphericalangle MA'H_1$ from relation (94) in which also we'll substitute $H_1A' = \frac{a}{2}ctg\,\omega$ we'll obtain

$$4s_{\Delta MBC} = \frac{2MA'^2 + \frac{3a^2}{2}}{ctg\omega} \quad (95)$$

Now we'll consider the relations (92) and (95) which substituted in (91) give $ctg\omega' = ctg\omega$; therefore $\omega' = \omega$ and the theorem is proved.

**Remark 27**

The geometric locus circle from the previous proved theorem is called the Henberg's circle. If we eliminate the restriction from the theorem's hypothesis the geometric locus would be made by two symmetric Henberg's circles in rapport to $BC$. To a triangle we can associate, in general, six Henberg's circles.

**Definition 35**

We call a Henberg's triangle of a given triangle $ABC$ the triangle $H_1H_2H_3$ formed by the centers of the Henberg's circles.

**Proposition 33**

The triangle $ABC$ and its triangle Henberg $H_1H_2H_3$ are homological. The homology center is the triangle's Tarry's point, $T$.

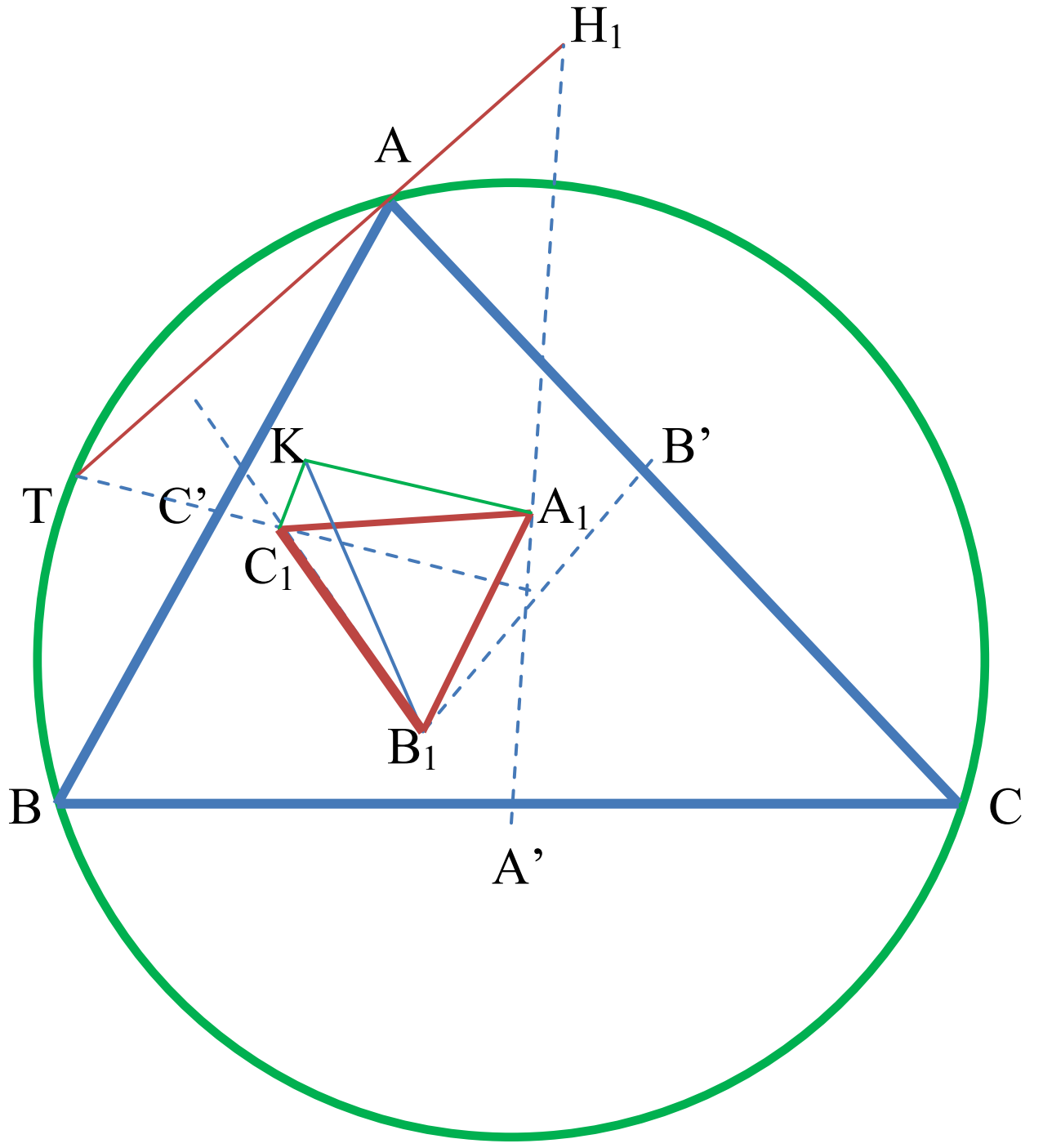

Fig. 43

**Proof**

We proved in theorem 22 that the perpendiculars from $A, B, C$ on the sides of the first Brocad triangle are concurrent in T. We'll prove now that the points $H_1, A, T$ are collinear. It is sufficient to prove that

$$\overrightarrow{H_1A}\cdot\overrightarrow{B_1C_1}=0 \tag{96}$$

We have $\overrightarrow{H_1A}=\overrightarrow{AA'}+\overrightarrow{A'H_1}$ and $\overrightarrow{B_1C_1}=\overrightarrow{C_1C'}+\overrightarrow{C'B'}+\overrightarrow{B'B_1}$ (see figure 43). Then

$$\overrightarrow{H_1A}=\cdot\overrightarrow{B_1C_1}=\left(\overrightarrow{AA'}+\overrightarrow{A'H_1}\right)\left(\overrightarrow{C_1C'}+\overrightarrow{C'B'}+\overrightarrow{B'B_1}\right) \tag{97}$$

Because $\overrightarrow{AA'}=\frac{1}{2}\left(\overrightarrow{AB}+\overrightarrow{AC}\right)$, $\overrightarrow{C'B'}=\frac{1}{2}\overrightarrow{BC}$, we have

$$\overrightarrow{H_1A}\cdot\overrightarrow{B_1C_1}=\frac{1}{2}\overrightarrow{AB}\cdot\overrightarrow{C_1C'}+\frac{1}{4}\overrightarrow{AB}\cdot\overrightarrow{BC}+\frac{1}{2}\overrightarrow{AB}\cdot\overrightarrow{B'B_1}+\frac{1}{2}\overrightarrow{AC}\cdot\overrightarrow{C_1C'}+$$

$$+\frac{1}{4}\overrightarrow{AC}\cdot\overrightarrow{BC}+\frac{1}{2}\overrightarrow{AC}\cdot\overrightarrow{B'B_1}+\overrightarrow{A'H_1}\cdot\overrightarrow{C_1C'}+\frac{1}{2}\overrightarrow{A'H_1}\cdot\overrightarrow{BC}+\overrightarrow{A'H_1}\cdot\overrightarrow{B'B_1} \tag{98}$$

Evidently, $\overrightarrow{AB}\cdot\overrightarrow{C_1C'}=0$, $\overrightarrow{AC}\cdot\overrightarrow{B'B_1}=0$, $\overrightarrow{A'H_1}\cdot\overrightarrow{BC}=0$

On the other side

$$\frac{1}{4}\overrightarrow{AB}\cdot\overrightarrow{BC}=-\frac{1}{4}ac\cdot\cos B,$$

$$\frac{1}{2}\overrightarrow{AC}\cdot\overrightarrow{C_1C'}=\frac{1}{4}bc\cdot tg\omega\sin A,$$

$$\frac{1}{4}\overrightarrow{AC}\cdot\overrightarrow{BC}=-\frac{1}{4}ab\cdot\cos C,$$

$$\overrightarrow{A'H_1}\cdot\overrightarrow{C_1C_1}=\frac{ca}{4}ctg\,\omega\cdot tg\,\omega\cdot\cos B,$$

$$\overrightarrow{A'H_1}\cdot\overrightarrow{B_1B_1}=-\frac{ab}{4}ctg\,\omega\cdot tg\,\omega\cdot\cos C$$

Considering these relations in the relation (98) it will result the relation (96).
Similarly it can be proved that $H_2B$ passes through $T$ and $H_3C$ also contains the point $C$.

**Theorem 25**

Let $ABC$ a triangle and $A_1B_1C_1$ its first Brocard triangle and $H_1H_2H_3$ is its Neuberg's triangle, then these triangles are two by two homological and have the same homological axis.

**Proof**

We roved that the triangle $ABC$ and its first Brocard's triangle are homological and their homological center is $\Omega''$ the third Brocard's point. Also we proved that the triangle $ABC$ and the Neuberg's triangle are homological and the homology center is Tarry's point $T$.

Taking into consideration that $H_1A_1$ is mediator in triangle $ABC$ for $BC$, and that $H_2B_1$ and $H_3C_1$ are also mediator in the same triangle, it results that the Neuberg's triangle and triangle $A_1B_1C_1$ (the first Brocard's triangle) are homological, the homology center being $O$, the center of the circumscribed circle of triangle $ABC$.

Consequently, the triangle $ABC$, its first Brocard's triangle $A_1B_1C_1$ and Hegel's triangle $H_1H_2H_3$ is a triplet of triangles two by two homological having the homology centers $\Omega'',T,O$.

We'll prove now that these homology points are collinear.

The idea for this proof is to compute the $\frac{PO}{PT}$, where we note $\{P\} = AA_1 \cap OT$, and to show that this rapport is constant. It will result then that $P$ is located on $BB_1$ and $CC_1$, and therefore $P = \Omega''$.

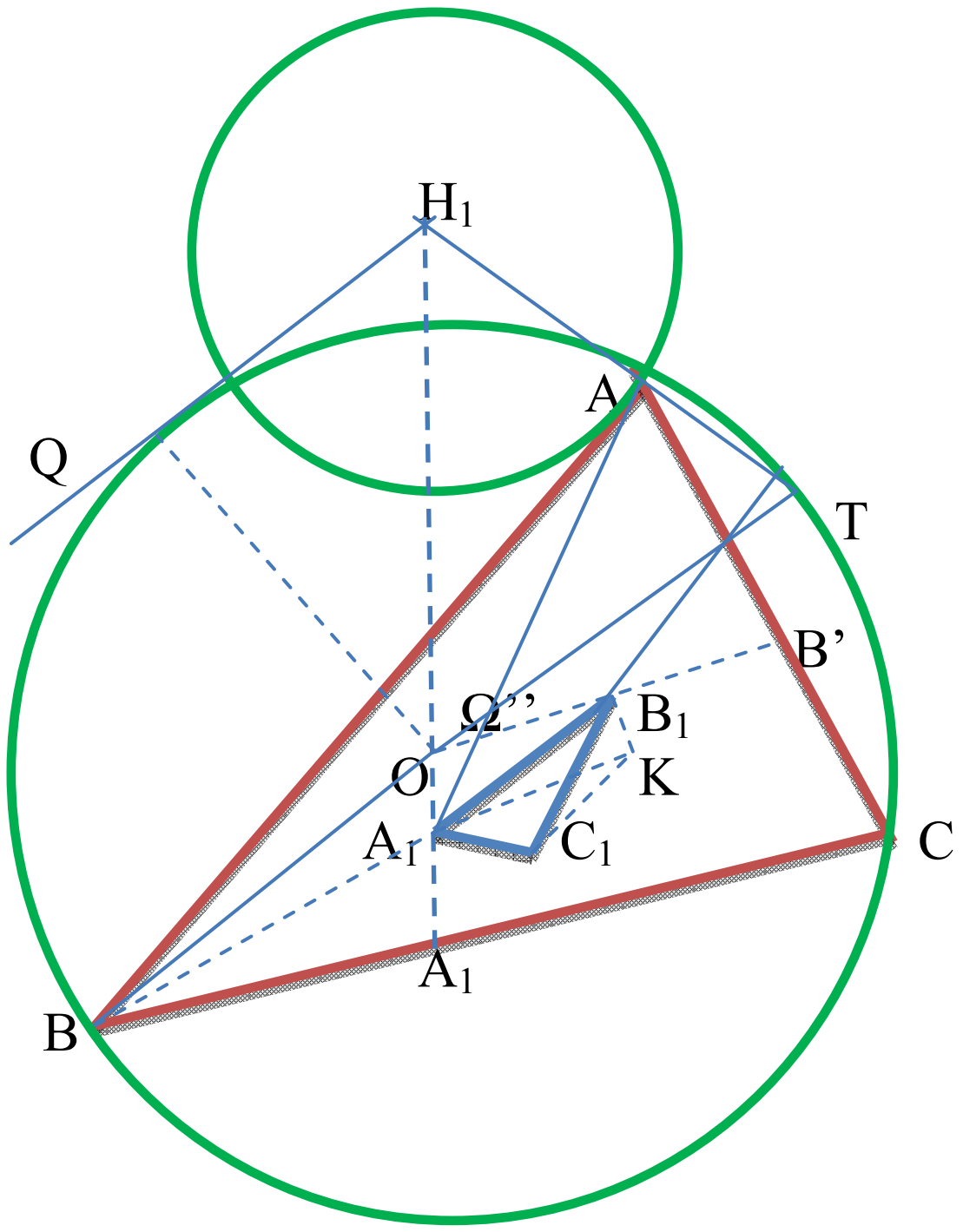

Fig. 44

Let then $\{P\} = AA_1 \cap OT$. Menelaus' theorem applied in the triangle $H_1OT$ for the transversal $A_1, P, A$ gives

$$\frac{A_1O}{A_1H_1} \cdot \frac{AH_1}{AT} \cdot \frac{PT}{PO} = 1 \qquad (99)$$

We have

$$A_1O = OA' - A_1A' = \frac{a}{2}(ctgA - tg\omega)$$

$$A_1H_1 = H_1A' - A_1A' = \frac{a}{2}(ctg\omega - tg\omega)$$

Therefore

$$\frac{A_1O}{A_1H_1} = \frac{ctgA - tg\omega}{ctg\omega - tg\omega} \qquad (100)$$

Considering the power of $H_1$ in rapport to the circumscribed circle of triangle $ABC$ it results

$$H_1A \cdot H_1T = H_1O^2 - OQ^2 = \frac{a}{2}\left(ctg\omega - \frac{a}{2}ctgA\right)^2 - R^2 \qquad (101)$$

We noted $Q$ he tangent point of the tangent from $H_1$ to the circumscribed circle of triangle $ABC$.

We know that

$$H_1A = \eta_1 = \frac{a}{2}\sqrt{ctg^2\omega - 3}$$

$$AT = H_1T - H_1A = \frac{\frac{a^2}{4}\left(ctg\omega - \frac{a}{2}ctgA\right)^2 - R^2 - \eta_1^2}{\eta_1}$$

It results that

$$\frac{A_1H_1}{AT} = \frac{\frac{a^2}{4}\left(ctg^2\omega - 3\right)}{\frac{a^2}{4}ctg^2\omega - \frac{a^2}{2}ctg\omega ctgA + \frac{a^2}{4}ctg^2A - \frac{a^2}{4\sin^2 A} - \frac{a^2}{4}ctg^2\omega + \frac{3a^2}{4}}$$

Because $\frac{1}{\sin^2 A} = 1 + ctg^2 A$ we find

$$\frac{A_1H_1}{AT} = \frac{ctg^2\omega - 3}{2\left(1 - ctg\omega ctgA\right)} \qquad (101)$$

Substituting the relation (100) and (99), we obtain

$$\frac{PT}{PO} = \frac{2\left(1 - ctg^2\omega\right)}{ctg^2\omega - 3} = \mathcal{A} \qquad (102)$$

If we note $P'$ the intersection of $OT$ with $BB_1$, we'll find, similarly $\frac{P'T}{P'O} = \frac{2\left(1 - ctg^2\omega\right)}{ctg^2\omega - 3}$.

It result that $P = P'$, therefore the intersection of the lines $AA_1, BB_1, CC_1$ noted $\Omega''$ coincides with $P$. The triangles from the considered triplet have therefore their homology centers collinear. Applying theorem 19 it results that these have the same homological axis and the theorem is proved.

## 2.3. Other theorems on homological triangles

**Proposition 34**

Let $ABC$ a triangle. We note $D_a, E_a, F_a$ the contact point of the A-ex-inscribed circle with the lines $BC, CA, AB$ respectively. The lines $AD_a, BE_a, CF_a$ are concurrent.

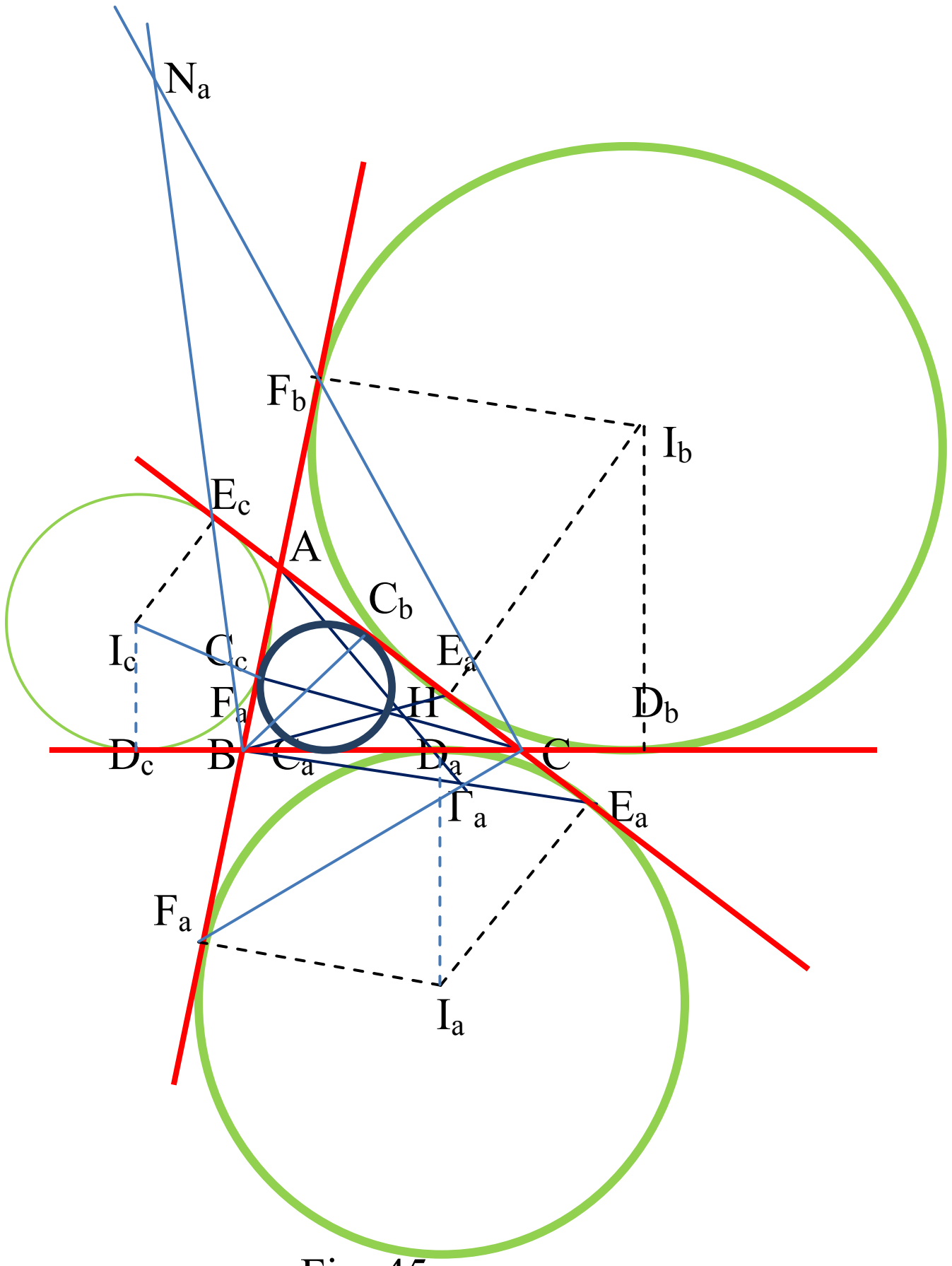

Fig. 45

**Proof**

We have that $AE_a = BF_a$, $BD_a = BF_a$ and $CD_a = CE_a$, then $\frac{D_aB}{D_aC} \cdot \frac{E_aA}{E_aC} \cdot \frac{F_aC}{F_aB} = 1$ and from the Ceva's reciprocal's theorem it results that the lines $AD_a, BE_a, CF_a$ are concurrent.

**Remark 28**

The concurrence point of the lines $AD_a, BE_a, CF_a$ is called the adjoin point of Gergonne's point ($\Gamma$) of triangle $ABC$, and has been noted it $\Gamma_a$. Similarly we define the adjoin points $\Gamma_b, \Gamma_c$. Because $AD_a$ is a Nagel Cevian of the triangle we have the following proposition.

**Proposition 35**

Triangle $ABC$ and triangle $\Gamma_a, \Gamma_b, \Gamma_c$ are homological. Their homology center is the Nagel's point (N) of triangle $ABC$.

**Proposition 36**

In triangle $ABC$ let $C_a$ be the contact point of the inscribed circle with $BC$, $F_b$ the contact point of the B-ex-inscribed circle with $AB$ and $E_c$ the contact point of the C-ex-inscribed circle with $AC$. The lines $AC_a, BE_c, CF_b$ are concurrent.

Proof

$$\frac{C_aB}{C_aC}\cdot\frac{E_cC}{E_cA}\cdot\frac{F_bA}{F_bB}=1 \qquad (1^*)$$

Indeed

$$C_aB = BC_c = AF_c = AE_c \qquad (2^*)$$

$$F_bA = AE_b = CC_b = CC_a \qquad (3^*)$$

$$E_cC = E_cA + AE_b + AF_b$$

But

$$E_cA = AF_c;\ AE_b = AF_b \text{ and } E_bC = AC_b = AC_c = BF_c$$

It result

$$E_cC = AF_c + AF_b + BF_c = BF_c \qquad (4^*)$$

Taking into account (2*), (3*), and (4*) we verified (1*), which shows that the Cevians $AC_a, BE_c, CF_b$ are concurrent.

**Remark 29**

The concurrence point of the Cevians $AC_a, BE_c, CF_b$ is called the adjoin point of Nagel, and we note it $H_a$. Similarly we define the adjoin points $H_b, H_c$ of the Nagel's point N.

Because the $AC_a, BC_b, CC_c$ are concurrent in the Gergonne's point $(\Gamma)$ of the triangle we can formulate the following proposition.

**Proposition 37**

The triangle $ABC$ and the triangle $H_aH_bH_c$ of the adjoin points of Nagel are homological. The homology point is Gergonne's point $(\Gamma)$.

**Theorem 26**

The triangle $\Gamma_a\Gamma_b\Gamma_c$ (having the vertexes in the adjoin Gergonne's points) and the triangle $H_aH_bH_c$ (having the vertexes in the adjoin Nagel's points) are homological. The center of homology belongs to the line $\Gamma H$ determined by the Geronne's and Nagel's points.

**Proof**

The triangle $ABC$ and triangle $H_aH_bH_c$ are homological, their homology center being $H$. We have

$$\{\Gamma_a\} = BH_c \cap CH_b,\ \{\Gamma_b\} = AH_c \cap CH_a,\ \{\Gamma_c\} = AH_b \cap BH_a$$

Applying the Veronese' theorem, it results that triangle $\Gamma_a\Gamma_b\Gamma_c$ is homological with the triangles $H_aH_bH_c$ and $ABC$. Furthermore, this theorem states that the homology centers of triangles $(ABC,\ H_aH_bH_c)$, $(ABC,\ \Gamma_a\Gamma_b\Gamma_c)$, $(\Gamma_a\Gamma_b\Gamma_c,\ H_aH_bH_c)$ are collinear.

We note $S$ the homology center of triangles $(\Gamma_a\Gamma_b\Gamma_c,\ H_aH_bH_c)$. It results that $S$ belongs to line $H\Gamma$.

**Remark 30**

Triangle $ABC$ and triangles $\Gamma_a\Gamma_b\Gamma_c$, $H_aH_bH_c$ have the same homological axis. This conclusion results from the precedent theorem and theorem 19.

**Theorem 27**

If two triangles one inscribed and the other circumscribed to the same triangle are homological with this triangle

Proof

Let triangle $A_1B_1C_1$ circumscribed to triangle $ABC$ and triangle $A_2B_2C_2$ inscribed in triangle $ABC$.

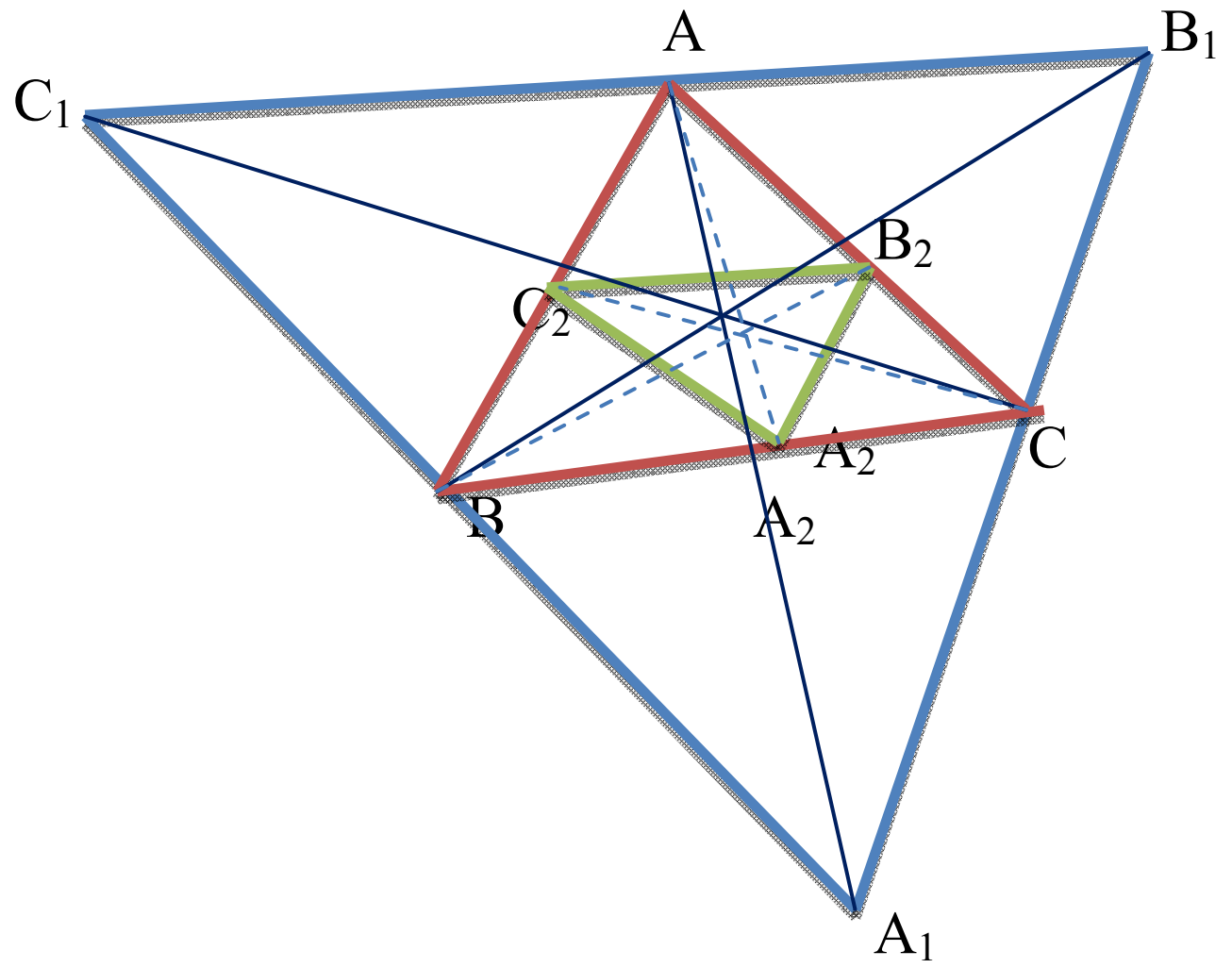

Fig. 46

Because $A_1B_1C_1$ and $ABC$ are homological the lines $A_1A,\ B_1B,\ C_1C$ are concurrent and therefore

$$\frac{AB_1}{AC_1}\cdot\frac{BC_1}{BA_1}\cdot\frac{CA_1}{CB_1}=1 \qquad (1^*)$$

Also $ABC$ and $A_2B_2C_2$ are homological triangles and consequently:

$$\frac{A_2B}{A_2C}\cdot\frac{B_2C}{B_2A}\cdot\frac{C_2A}{C_2B}=1 \qquad (2^*)$$

We have

$$\frac{Aria\Delta A_1CA_2}{Aria\Delta A_1BA_2}=\frac{A_2C}{A_2B}=\frac{A_1C\cdot A_1A_2\cdot\sin\left(CA_1A_2\right)}{A_1B\cdot A_1A_2\cdot\sin\left(BA_1A_2\right)}$$

From here

$$\frac{\sin\left(CA_1A_2\right)}{\sin\left(BA_1A_2\right)}=\frac{A_2C}{A_2B}\cdot\frac{A_1B}{A_1C} \qquad (3^*)$$

Similarly we find

$$\frac{\sin(AB_1B_2)}{\sin(CB_1B_2)} = \frac{B_2A}{B_2C} \cdot \frac{B_1C}{B_1A} \tag{4*}$$

$$\frac{\sin(BC_1C_2)}{\sin(AC_1C_2)} = \frac{C_2B}{C_2A} \cdot \frac{C_1A}{C_1B} \tag{5*}$$

Multiplying the relations (3*), (4*), (5*) and taking into account (1*) and (2*) gives

$$\frac{\sin(CA_1A_2)}{\sin(BA_1A_2)} \cdot \frac{\sin(AB_1B_2)}{\sin(CB_1B_2)} \cdot \frac{\sin(BC_1C_2)}{\sin(AC_1C_2)} = 1$$

This relation and Ceva's theorem (the trigonometric form shows the concurrence of the lines $A_1A_2$, $B_1B_2$, $C_1C_2$ and implicitly the homology of triangles $A_1B_1C_1$ and $A_2B_2C_2$.

**Remark 31**

Theorem 27 has a series of interesting consequences that provide us other connections regarding the triangles homology associate to a given triangle.

**Proposition 38**

The anti-supplementary triangle and the contact triangle of a given triangle $ABC$ are homological.

**Proof**

The anti-supplementary triangle $I_aI_bI_c$ and the contact triangle $C_aC_bC_c$ of triangle $ABC$ are respectively circumscribed and inscribed to $ABC$. Also these triangles are homological with $ABC$ (see proposition 5 and proposition 12), in conformity with theorem 27 these are also homological.

**Proposition 39**

The tangential triangle and the contact triangle of a given triangle $ABC$ are homological.

**Proof**

The proof results from the precedent theorem and from propositions 10 and 12.

**Proposition 40**

The tangential triangle and the orthic triangle of a given triangle $ABC$ are homological. The proof results from theorem 27 and from proposition 10.

**Remark 32**

The homology center of the tangential triangle and of the orthic triangle is called the Gob's point, and traditionally is noted $\Phi$.

**Proposition 41**

The anti-supplementary triangle of a given triangle $ABC$ and its cotangent triangle are homological.

**Proof**

From proposition 5 and from proposition 17 it results that the triangle $ABC$ is homological with triangle $I_aI_bI_c$ and with its cotangent triangle. From theorem 27 it results the conclusion.

**Remark 33**

The homology center of triangle $I_aI_bI_c$ and of the cotangent triangle is noted $V$ and it is called Beven's point . It can be proved that in a triangle the points $I, O$ and $V$ are collinear.

**Definition 36**

An anti-complementary triangle of a given triangle $ABC$ is the triangle formed by the parallel lines constructed in the triangle's vertexes $A, B, C$ to the opposite sides of the given triangle.

**Proposition 36**

The anti-complementary triangle and the orthic triangle of a given triangle are homological.

The proof of this proposition results from the theorem 27 and from the observation that the anti-complementary triangle of triangle $ABC$ is homological with it, the homology center being the weight center of the triangle $ABC$.

**Theorem 28**

Let $ABC$ and $A_1B_1C_1$ two homological triangles having their homology center in the point $O$. The lines $A_1B, A_1C$ intersect respectively the lines $AC, AB$ in the points $M, N$ which determine a line $d_1$. Similarly we obtain the lines $d_2, d_3$. We note $\{A_2\} = d_2 \cap d_3$, $\{B_2\} = d_3 \cap d_1$, $\{C_2\} = d_1 \cap d_2$ , then the triangle $ABC$ and triangle $A_2B_2C_2$ are homological having as homological axis the tri-linear polar of $O$ in rapport to triangle $ABC$.

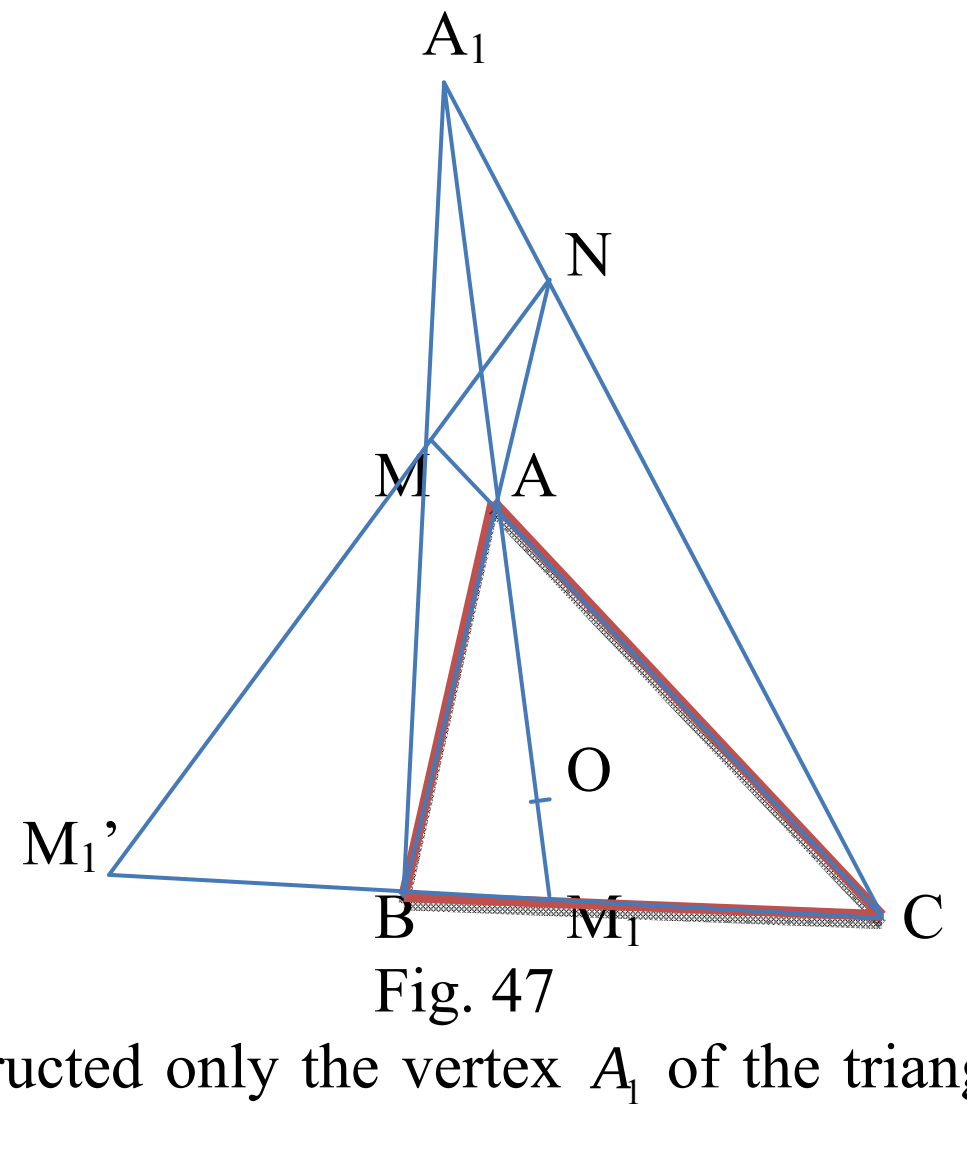

Fig. 47

In figure 47 we constructed only the vertex $A_1$ of the triangle $A_1B_1C_1$ in order to follow the rational easier.

We note $\{M_1\} = A_1A \cap BC$ and $\{M_1^{'}\} = MN \cap BC$. We noticed that $M_1^{'}, B, M, C$ form a harmonic division. Considering in triangle $ABC$ the Cevian $AO$ and $M_1$ its base, it results that the tri-linear polar of $O$ in rapport with triangle $ABC$ intersects $BC$ in $M_1^{'}$. Similarly it can be shown that the lines $d_2, d_3$ intersect the lines $CA, AB$ in points that belong to the tri-linear polar of $O$. Conform to the Desargues' theorem we have that because the triangles $A_2B_2C_2$ and $ABC$ are homological, their homology axis being the tri-linear polar of $O$ in rapport with triangle $ABC$.

**Remark 34**

An analogue property is obtained if we change the role of triangles $ABC$ and $A_1B_1C_1$. We'll find a triangle $A_3B_3C_3$ homological to $A_1B_1C_1$, their homology axis being the tri-linear polar of $O$ in rapport with triangle $A_1B_1C_1$.

An interested particular case of the precedent theorem is the following:

**Theorem 29**

If $ABC$ is a given triangle, $A_1B_1C_1$ is its Cevian triangle and $I_aI_bI_c$ is its anti-supplementary triangle, and if we note $M, N$ the intersection points of the lines $A_1I_b, A_1I_c$ respectively with $I_aI_c, I_aI_b$, $d_1$ the line $MN$, similarly we obtain the lines $d_2, d_3$; let $\{A_2\} = d_2 \cap d_3$, $\{B_2\} = d_3 \cap d_1$, $\{C_2\} = d_1 \cap d_2$, also we note $M_1, N_1$ the intersection points between the lines $B_1I_a, C_1I_a$ respectively with the lines $A_1C_1$, $A_1B_1$. Let $d_4$ the line $M_1N_1$, similarly we obtain $d_5, d_6$, $\{A_3\} = d_5 \cap d_6$, $\{B_3\} = d_4 \cap d_6$, $\{C_3\} = d_4 \cap d_5$. Then

i. $A_1B_1C_1$ and $I_aI_bI_c$ are homological
ii. $A_2B_2C_2$ and $I_aI_bI_c$ are homological
iii. $A_3B_3C_3$ and $A_1B_1C_1$ are homological
iv. The pairs of triangles from above have as homology axis the tri-linear polar of $I$ in rapport to triangle $ABC$ (the anti-orthic axis of triangle $ABC$)

**Proof**

i. Because $AI_a, BI_b, CI_c$ are concurrent in $I$, the center of the inscribed circle, we have that $A_1I_a, B_1I_b, C_1I_c$ are concurrent in $I$, therefore the triangles $A_1B_1C_1$ (the Cevian triangle) and $I_aI_bI_c$ (the triangle anti-supplementary) are homological.

We note $\{I_1\} = C_1B_1 \cap BC$ and $\{A_1^{'}\} = AA_1 \cap B_1C_1$; we noticed that $I_1$ and $A_1^{'}$ are harmonic conjugate in rapport to $C_1$ and $B_1$, because $AA_1^{'}$ is an interior bisector in triangle $AC_1B_1$ it result that that $AI_1$ could be the external bisector of the angle $A$, therefore $I_1$ belongs to the tri-linear polar of $I$ in rapport to $ABC$. On the other side $I_1$ is the intersection between $I_bI_c$ and $B_1C_1$ therefore it belongs to the homological axis of triangles $A_1B_1C_1$ and $I_aI_bI_c$. We note $\{I_2\} = A_1C_1 \cap AC$ and $\{I_3\} = A_1B_1 \cap AB$ and it results that $I_1, I_2, I_3$ are the feet of the external bisectors of triangle $ABC$. Furthermore, it can be shown that $I_1 - I_2 - I_3$ is the tri-linear polar of $I$ in rapport with triangle $A_1B_1C_1$

The proof for ii), iii), and iv) result from theorem 27.

**Remark 35**

a) The pairs of triangles that belong to the triplet $(I_aI_bI_c, ABC, A_1B_1C_1)$ have the same homology center – the point $I$ and the same homological axis anti-orthic of triangle $ABC$.

b) In accordance with theorem 18 it result that $I$ and the homological centers of the triangles $A_2B_2C_2$, $I_aI_bI_c$ and $A_1B_1C_1$ are collinear.

**Theorem 30**

Let $ABC$ a given triangle, $O$ its circumscribed circle and $H$ its orthocenter and $A_1B_1C_1$ its orthic triangle. We note $M, H, P$ the middle points of the segments $(AH),(BH),(CH)$. The perpendiculars constructed in $A, B, C$ on $OM, OH, OP$ from triangle $A_2B_2C_2$. The triplet of triangles ($ABC, A_1B_1C_1, A_2B_2C_2$) is formed by triangles two by two homological having the same homological axis, which is the orthic axis of triangle $ABC$.

**Proof**

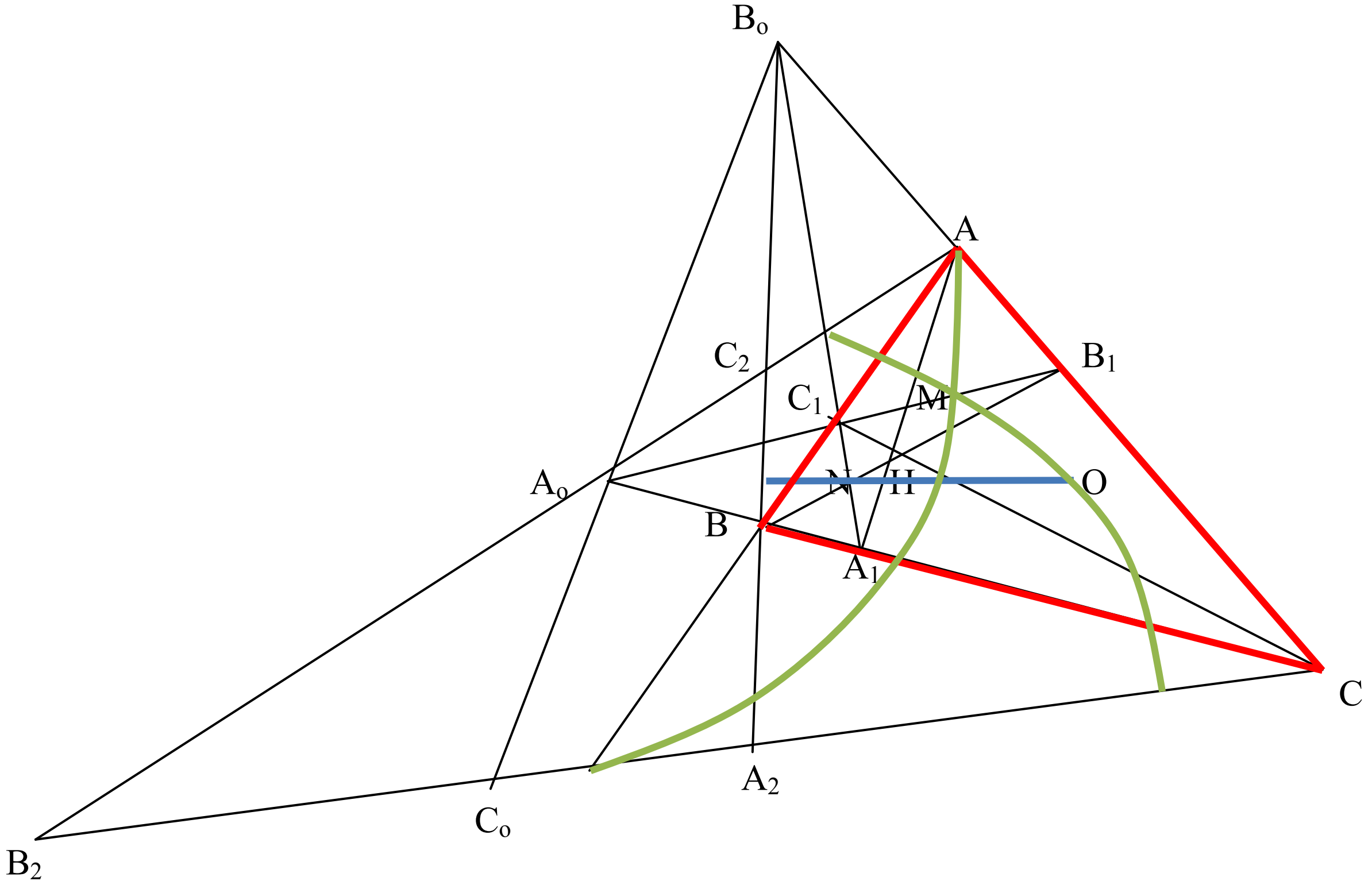

Fig. 48

The center of the circumscribed circle of triangle $A_1B_1C_1$ is the point $M$ (the quadrilateral $AB_1HC_1$ is inscribable, the circle's center being $M$ ).The perpendicular from $A$ on $OM$ , which we note $d_1$ is the radical axis of the arches $(AB_1C_1)$ and $(ABC)$. The line $B_1C_1$ is a radical axis of the circle $AB_1C_1$ and of Euler's circle $(A_1B_1C_1)$. Let $A_o$ the intersection of the lines $d_1$ and $B_1C_1$ (see figure 48). This point is the radical center of the mentioned arcs and $A_o$ belongs also to the tri-linear polar of $H$ in rapport with $ABC$, that is of orthic axis of triangle $ABC$.

Similarly, we note $d_2$ the perpendicular from $B$ on $OH$ and $d_3$ the perpendicular from $C$ on $OP$, $B_o$ the intersection between $d_2$ and $A_1C_1$ and $\{C_o\} = d_3 \cap A_1B_1$. We find that $B_o, C_o$ are on the orthic axis of triangle $ABC$ The homological sides of triangles $ABC$, $A_1B_1C_1$ and $A_2B_2C_2$ intersect in the collinear points $A_o, B_o, C_o$ which belong to the orthic axis of triangle $ABC$.

**Remark 36**

According to theorem 18, the homology centers of the triangle triplet mentioned are collinear.

**Theorem 31**

Let $ABC$ a given triangle, let $O$ the center of the circumscribed triangle, $I$ the center of the inscribed circle and $C_aC_bC_c$ its contact triangle. We note with $M, N, P$ the middle points of the segments $IA, IB, IC$ respectively and the perpendiculars constructed from $A, B, C$ respectively on $OM, ON, OP$ form a triangle $A_1B_1C_1$. The triplet $(ABC, C_aC_bC_c, A_1B_1C_1)$ contains triangles two by two homological having a common homological axis, which is the radical axis of the circumscribed and inscribed circles to triangle $ABC$.

**Proof**

The circumscribed circle of triangle $AC_bC_c$ is the point $M$ , it result that the perpendicular $d_1$ constructed from $A$ on $OM$ is the radical axis of circles $(AC_bC_c)$ and $(ABC)$ On the other side $C_bC_c$ is he radical axis of circles $(AC_bC_c)$ and inscribed to triangle $ABC$. The intersection point $A_o$ of lines $d_1$ and $C_bC_c$ is therefore, the radical center of the mentioned circles; it is situated on the radical axis $d$ of the arcs inscribed and circumscribed to triangle $ABC$.

Similarly are defined the points $B_o, C_o$ , and it results that these belong to line $d$ . Because the corresponding sides of triangles $ABC$, $C_aC_bC_c$, and $A_1B_1C_1$ intersect in the collinear points $A_o, B_o, C_o$. It results, in conformity with the Desargues' theorem that these triangles are two by two homological and their common homological axis is the radical axis of the inscribed circle of triangle $ABC$.

**Remark 37**

We saw that the triangle $ABC$ and its contact triangle $C_aC_bC_c$ are homological, the homology center being $\Gamma$, the Geronne's point, and the homological axis is the Lemoine's line of the contact triangle (Proposition 12). Taking into account the precedent theorem and this result we can make the following statement: The radical axis of the circumscribed and inscribed circles of triangle $ABC$ is the Lemoine's line of the contact triangle of the triangle $ABC$.

**Definition 37**

We call the anti-pedal triangle of point $M$ relative to triangle $ABC$, the triangle formed by the perpendiculars constructed in $A, B, C$ on $MA, MB, MC$ respectively.

**Theorem 32**

Let $M_1, M_2$ two points in the plane of the triangle $ABC$ symmetric in rapport to $O$, which is the center of the circumscribed circle. If $A_1B_1C_1$ is the pedal triangle of $M_1$ and $A_2B_2C_2$ is the anti-pedal triangle of the point $M_2$, then these triangles are homological. The homology axis is the radical axis of the circumscribed circles to triangles $ABC$ and $A_1B_1C_1$, and the homological center is the point $M_1$.

**Proof**

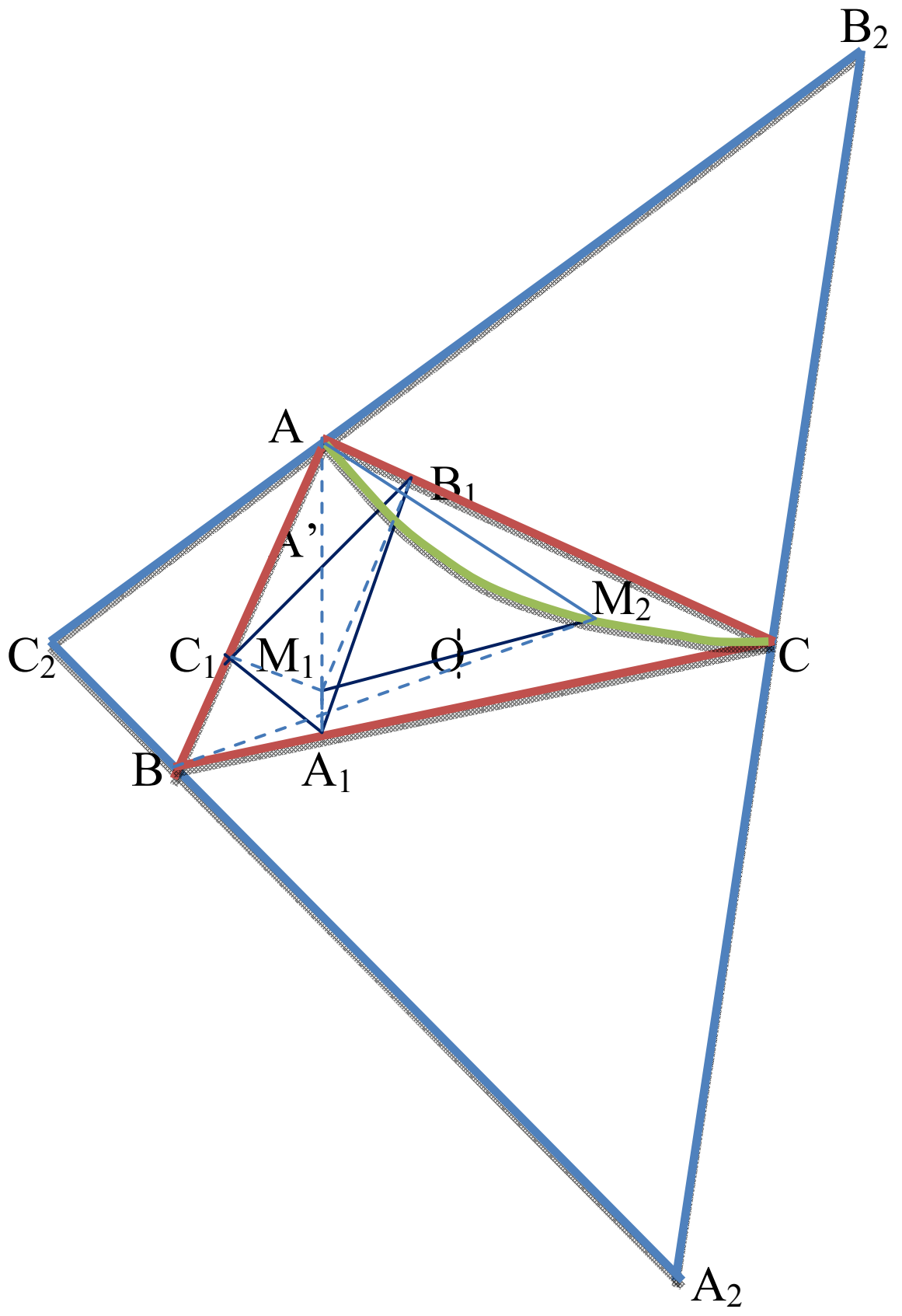

Fig. 49

Let $A', B', C'$ the middle points of the segments $(AM_1), (BM_1), (CM_1)$ (see Fig. 49).
The line $B_2C_2$ is perpendicular on $AM_2$ and because $A'$ is the center of the circumscribed circle.

$B_2C_2$ is the radical axis of the circles $(ABC)$ and $(AB_1C_1)$. On the other side the lane $B_1C_1$ is the radical axis of the circles $(AB_1C_1)$ and $(A_1B_1C_1)$, it result that the point $A_o$, the intersection of the lines $B_1C_1$ and $B_2C_2$ is the radical center of the three mentioned circles, circumscribed to triangles $ABC$ and $A_1B_1C_1$. Similarly it can be proved that the points $B_o, C_o$ in which the lines $A_1C_1$ and $A_2C_2$ respectively $A_1B_1$ and $A_2B_2$ intersect belong to line $d$. Therefore the triangles $A_1B_1C_1$ and $A_2B_2C_2$ are homological having as homological axis line $d$, which is the radical axis of the circumscribed circles to triangles $ABC$ and $A_1B_1C_1$.

Because the line $B_2C_2$ is the radical axis of the circles $(AB_1C_1)$ and $(ABC)$, and $A_2C_2$ is the radical axis of the circles $(BC_1A_1)$ and $(ABC)$, it result that the point $C_2$ is the radical center of these circles, therefore it belongs to the line $M_1C_1$ which is the radical axis of circles $(AB_1C_1)$ and $(BC_1A_1)$, consequently the line $C_1C_2$ passes through $M_1$. Similarly it can be shown that the lines $B_2B_1$ and $A_2A_1$ pass through point $M_1$. Therefore this point is the homological center of the triangles $A_1B_1C_1$ and $A_2B_2C_2$.

**Remark 38**

Proposition 14 can be considered a particular case of this theorem. Therefore we obtain that the homological axis of the medial triangle and of the tangential triangle of a given triangle $ABC$ is the radical axis of the circumscribed circle of the triangle $ABC$ and of Euler's circle of triangle $ABC$.

# Chapter 3

## Bi-homological and tri-homological triangles

In this chapter we'll prove a theorem that expresses the necessary and sufficient condition that characterizes the homology of two triangles.

This theorem will allow us to prove another theorem that states that two triangles are bi-homological then these are tri-homological.

### 3.1. The necessary and sufficient condition of homology

**Theorem 33**

The triangles $ABC$ and $A_1B_1C_1$ are homological if and only if

$$\frac{Aria\Delta(A_1AB)}{Aria\Delta(A_1AC)}\cdot\frac{Aria\Delta(B_1BC)}{Aria\Delta(B_1BA)}\cdot\frac{Aria\Delta(C_1CA)}{Aria\Delta(C_1CB)}=-1$$

**Proof**

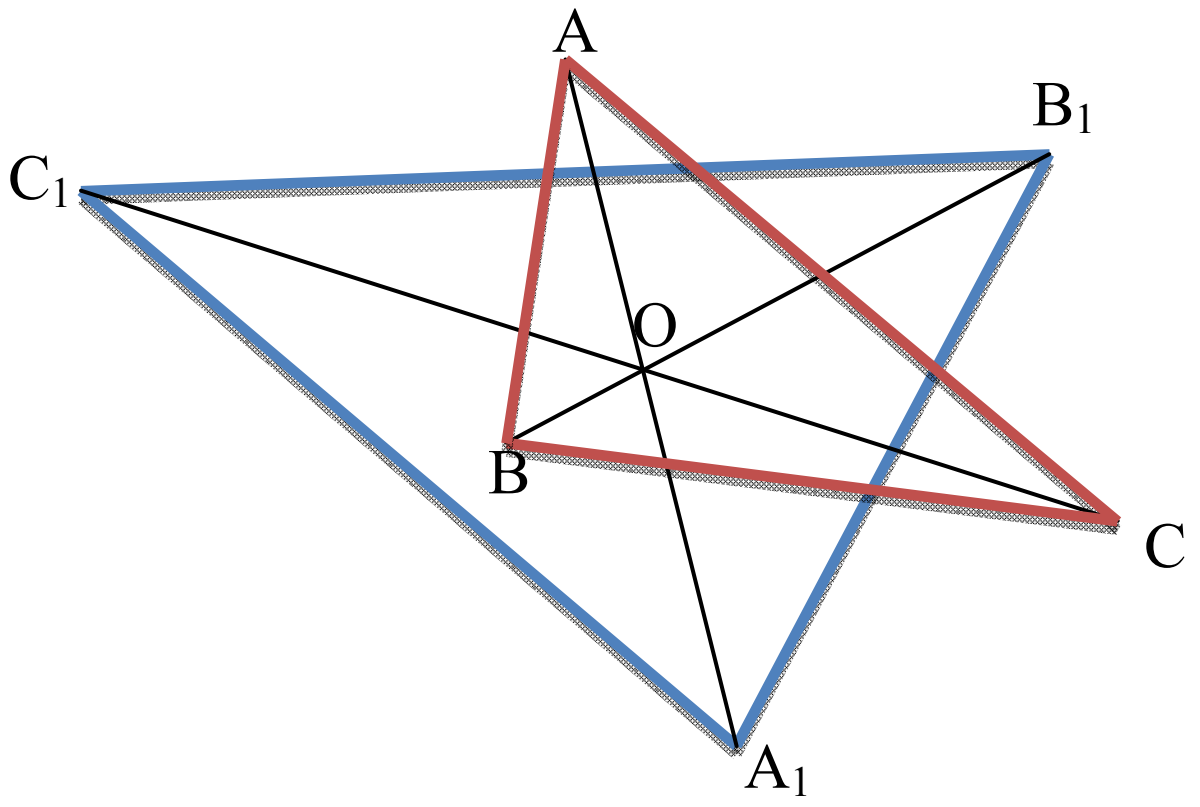

Fig 50

*The condition is necessary.*

The triangles $ABC$ and $A_1B_1C_1$ being homological the lines $AA_1, BB_1, CC_1$ are concurrent in a point $O$ (see Fig. 50).

We have

$$\frac{Aria\Delta(A_1AB)}{Aria\Delta(A_1AC)}=\frac{AB\cdot AA_1\cdot\sin(A_1AB)}{AC\cdot AA_1\cdot\sin(A_1AC)} \quad (1)$$

$$\frac{Aria\Delta(B_1BC)}{Aria\Delta(B_1BA)}=\frac{BC\cdot BB_1\cdot\sin(B_1BC)}{BA\cdot BB_1\cdot\sin(B_1BA)} \quad (2)$$

$$\frac{Aria\Delta(C_1CA)}{Aria\Delta(C_1CB)}=\frac{CA\cdot CC_1\cdot\sin(C_1CA)}{CB\cdot CC_1\cdot\sin(C_1CB)} \quad (3)$$

Multiplying these relations side by side it results

$$\frac{Aria\Delta(A_1AB)}{Aria\Delta(A_1AC)}\cdot\frac{Aria\Delta(B_1BC)}{Aria\Delta(B_1BA)}\cdot\frac{Aria\Delta(C_1CA)}{Aria\Delta(C_1CB)}=\frac{\sin(A_1AB)}{\sin(A_1AC)}\cdot\frac{\sin(B_1BC)}{\sin(B_1BA)}\cdot\frac{\sin(C_1CA)}{\sin(C_1CB)}$$

From Ceva's theorem (the trigonometric form) it results that the relation from the hypothesis is true.

*The condition is sufficient*

If the given relation is satisfied, it results that

$$\frac{\sin(A_1AB)}{\sin(A_1AC)}\cdot\frac{\sin(B_1BC)}{\sin(B_1BA)}\cdot\frac{\sin(C_1CA)}{\sin(C_1CB)}=-1$$

The Ceva's reciprocal theorem gives us the concurrence of the Cevians $AA_1, BB_1, CC_1$, therefore the homology of triangles $ABC$ and $A_1B_1C_1$.

## 3.2. Bi-homological and tri-homological triangles

**Definition 38**

The triangle $ABC$ is direct bi-homological with triangle $A_1B_1C_1$ if triangle $ABC$ is homological with $A_1B_1C_1$ and with $B_1C_1A_1$.

The triangle $ABC$ is direct tri-homological with triangle $A_1B_1C_1$ if is homological with $B_1C_1A_1$ and $B_1A_1C_1$, and $ABC$ is invers tri-homological with triangle $A_1B_1C_1$ if $ABC$ is homological with $A_1C_1B_1$ with $B_1A_1C_1$ and with $C_1B_1A_1$.

**Theorem 34** (Rosanes – 1870)

If two triangles $ABC$ and $A_1B_1C_1$ are direct bi-homological then these are direct tri-homological.

**Proof**

If triangles $ABC$ and $A_1B_1C_1$ are homological then

$$\frac{Aria\Delta(A_1AB)}{Aria\Delta(A_1AC)}\cdot\frac{Aria\Delta(B_1BC)}{Aria\Delta(B_1BA)}\cdot\frac{Aria\Delta(C_1CA)}{Aria\Delta(C_1CB)}=-1 \quad (1)$$

If triangles $ABC$ and $B_1C_1A_1$ are homological, then

$$\frac{Aria\Delta(B_1AB)}{Aria\Delta(B_1AC)}\cdot\frac{Aria\Delta(C_1BC)}{Aria\Delta(C_1BA)}\cdot\frac{Aria\Delta(A_1CA)}{Aria\Delta(A_1CB)}=-1 \quad (2)$$

Taking into consideration that $Aria\Delta(A_1AC)=-aria\Delta(A_1CA)$, $Aria\Delta(B_1BA)=-Aria\Delta(B_1AB)$, $Aria\Delta(C_1CB)=-Aria\Delta(A_1CA)$. By multiplying side by side the relations (1) and (2) we obtain

$$\frac{Aria\Delta(C_1AB)}{Aria\Delta(C_1AC)}\cdot\frac{Aria\Delta(A_1BC)}{Aria\Delta(A_1BA)}\cdot\frac{Aria\Delta(B_1CA)}{Aria\Delta(B_1CB)}=-1 \quad (3)$$

The relation (3) shows that the triangles $ABC$, $C_1A_1B_1$ are homological, therefore the triangles $ABC$ and $A_1B_1C_1$ are direct tri-homological.

**Remark 39**

Similarly it can be proved the theorem: If triangles $ABC$ and $A_1B_1C_1$ are triangles inverse bi-homological, then the triangles are inverse tri-homological.

**Proposition 43**

If triangle $ABC$ is homological with the triangles $A_1B_1C_1$ and $A_1C_1B_1$ then the centers of the two homologies are collinear with the vertex $A$.

The proof of this theorem is immediate.

The Rosanes' theorem leads to a method of construction of a tri-homological triangle with a given triangle, as well as of a triplet of triangles two by two tri-homological as it results from the following theorem.

**Theorem 35**

(i) Let $ABC$ a given triangle and $\Gamma, Q$ two points in its plane. We note $\{A_1\} = BP \cap CQ$, $\{B_1\} = CP \cap AQ$, and $\{C_1\} = AP \cap BQ$.

Triangles $ABC$, $A_1B_1C_1$ are tri-homological.

(ii) If $\{A_2\} = BQ \cap CP$, $\{B_2\} = CQ \cap AP$ and $\{C_2\} = AQ \cap BP$, then the triangles $ABC$, $A_1B_1C_1$, $A_2B_2C_2$ are two by two tri-homological, and their homological centers are collinear.

(iii) We note $\{R\} = AA_1 \cap BB_1$, if the points $P, Q, R$ are not collinear then the triangle $RPQ$ is direct tri-homological with $ABC$ and the invers triangle $A_1B_1C_1$.

**Proof**

(i) From the hypothesis it results that the triangles $ABC$ and $A_1B_1C_1$ are bi-homological $AB_1 \cap BC_1 \cap CA_1 = \{Q\}$, $AC_1 \cap BA_1 \cap CB_1 = \{P\}$ (see Fig. 51). In accordance to Rosanes' theorem it result that $AA_1 \cap BB_1 \cap CC_1 = \{R\}$. Therefore, the triangles $ABC$ and $A_1B_1C_1$ are tri-homological.

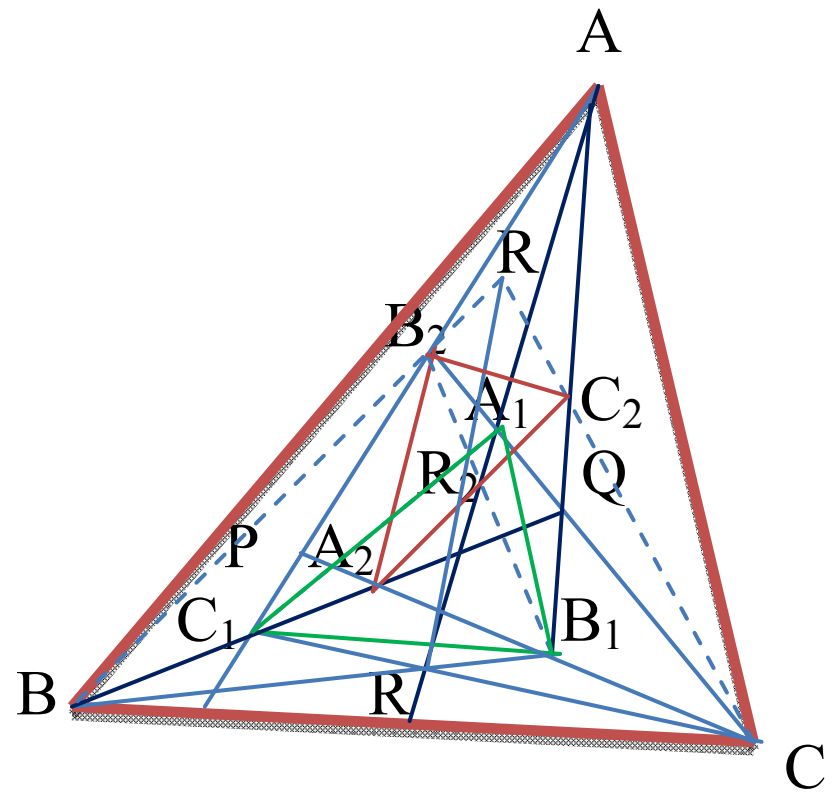

Fig. 51

(ii) We observe that $AB_2 \cap BC_2 \cap CA_2 = \{P\}$ and $AC_2 \cap BA_2 \cap CB_2 = \{Q\}$, therefore the triangles $ABC$ and $A_2B_2C_2$ are bi-homological. In conformity with Rosanes' theorem we have that triangles $ABC$ and $A_2B_2C_2$ are tri-homological, therefore $AA_2, BB_2, CC_2$ are concurrent in a point $R_1$. Also the triangles $A_1B_1C_1$ and $A_2B_2C_2$ are bi-homological having as homological centers the points $Q, P$. In conformity with the same theorem of Rosanes we'll have that $A_1A_2, B_1B_2, C_1C_2$ are concurrent in a point $R_2$. Consequently, the triangles $A_1B_1C_1$ and $A_2B_2C_2$ are tri-homological.

The Veronese's theorem leads us to the collinearity of the homological centers, therefore the points $R, R_1, R_2$ are collinear.

(iii) Indeed triangle $RQP$ is direct tri-homological with $ABC$ because

$RA \cap PB \cap QC = \{A_1\}$, $RB \cap PC \cap QA = \{B_1\}$, $RC \cap PA \cap QB = \{C_1\}$.

Also $RQP$ is invers tri-homological with $A_1B_1C_1$ because

$RA_1 \cap PC_1 \cap QB_1 = \{A\}$, $RB_1 \cap PA_1 \cap QC_1 = \{B\}$ and $RC_1 \cap PB_1 \cap QA_1 = \{C\}$.

**Remark 40**

a) Considering the points $P, R$ and making the same constructions as in the previous theorem we obtain the triangle $A_3B_3C_3$ which along with the triangles $ABC$ and $A_1B_1C_1$ form a triplet of triangles two by two tri-homological.

b) Another triplet of tri-homological triangles is obtained considering the points $R, Q$ and making similar constructions.

c) Theorem 35 shows that given a triangle $ABC$ and two points $P, Q$ in its plane, we can construct an unique triangle $RPQ$ directly tri-homological with the given triangle $ABC$.

d) Considering the triangle $ABC$ and as given points in its plane the Brocard's points $\Omega$ and $\Omega'$, the triangle $A_1B_1C_1$ constructed as in the previous theorem, is the first Brocard's triangle. We find again the J. Neuberg's result that tells us that the triangle $ABC$ and its first Brocard's triangle are tri-homological (see theorem 17). More so we saw that the homological center of triangles $ABC$ and $A_1B_1C_1$ is the isotomic conjugate of the symmedian center, noted $\Omega''$. From the latest resultants obtain lately, it results that the triangle $\Omega\Omega'\Omega''$ is tri-homological with $ABC$.

**Proposition 44**

In the triangle $ABC$ let's consider $A', B', C'$ the feet of its heights and $A_1, B_1, C_1$ the symmetric points of the vertexes $A, B, C$ in rapport to $C', A', B'$, and $A_2, B_2, C_2$ the symmetric points of the vertexes $A, B, C$ in rapport to $B', C', A'$.

If $M_1, N_1, P_1$ are the centers of the circles $BCC_1, CAA_1, ABB_1$ and $M_2, N_2, P_2$ are the centers of the circles $CBB_2, ACC_2, BAA_2$, then the triangles $M_1N_1P_1$ and $M_2N_2P_2$ are tri-homological.

Proof

Let $H$ the orthocenter of the triangle $ABC$, $O$ the center of the circumscribed circle, $O_9$ the center of the Euler's circle of the given triangle and $A'', B'', C''$ the middle points of the sides $BC, CA, AB$ (see Fig. 52).

The points $M_1, N_1, P_1$ are the intersections of the pairs of lines $(OA'', BH)$, $(OB'', CH)$, $(OC'', AH)$ and the points $M_2$ ,$N_2$ , $P_2$ are the intersections of the pairs of lines $(OA'', CH)$, $(OB'', AH)$ and $(OC'', BH)$.

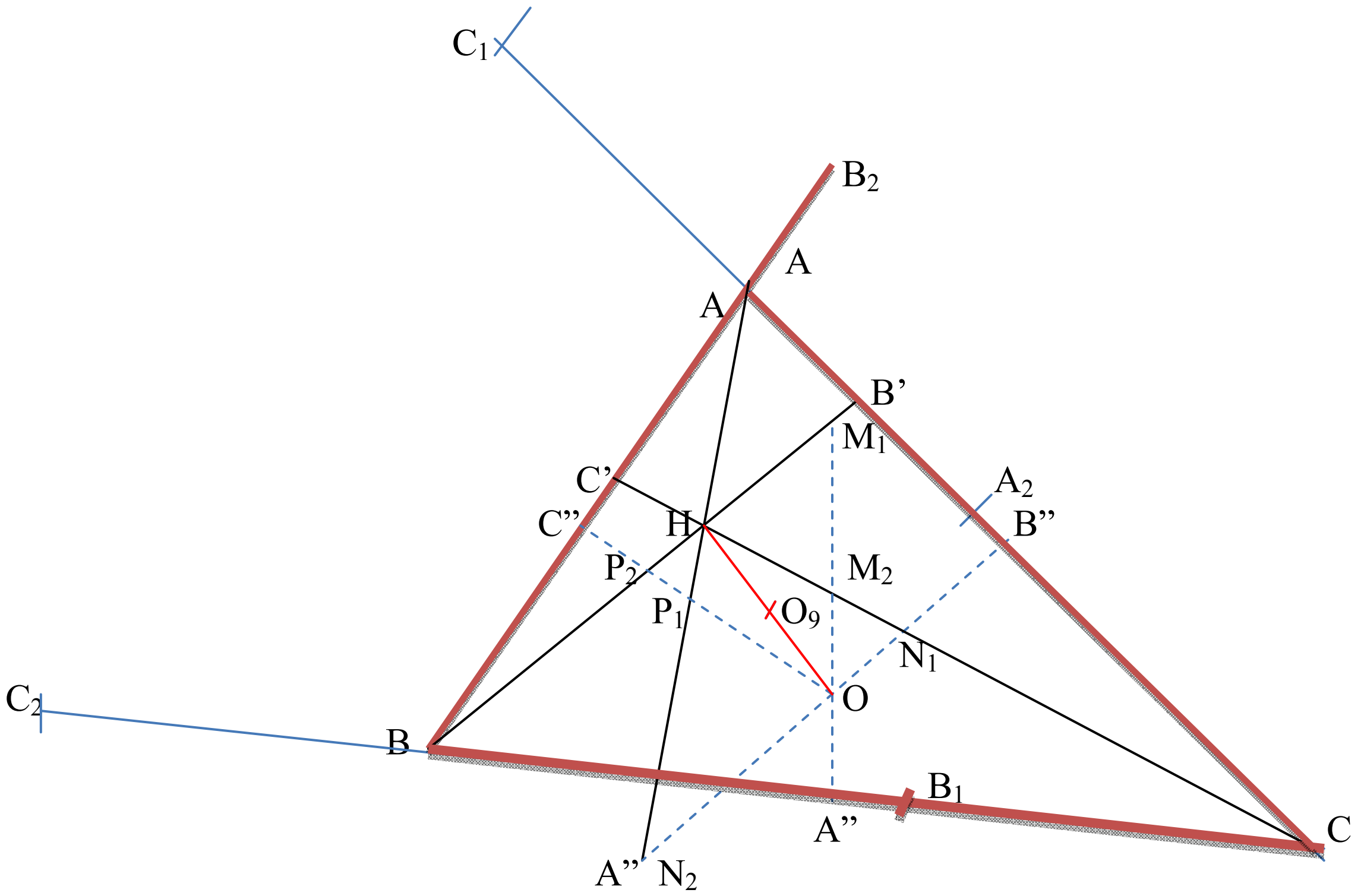

Fig. 52

The triangles $M_1N_1P_1$ and $M_2N_2P_2$ are homological because the lines $M_1M_2, N_1N_2, P_1P_2$ are concurrent in the point $O$. The triangle $M_1N_1P_1$ is homological to $N_2P_2A_2$ because the lines $M_1N_2, N_1P_2, P_1A_2$ are concurrent in the point $O_9$ (indeed $M_1N_2$ is a diagonal in the parallelogram $M_1HN_2O$).

The triangles $M_1N_1P_1$ and $M_2N_2P_2$ being bi-homological, it results that are also tri-homological, and the proposition is proved.

**Observation 27**

The homology of triangles $M_1N_1P_1$ and $M_2N_2P_2$ results also directly by observing that the lines $M_1P_2, N_1N_2, P_1N_2$ are concurrent in $H$ .

The homological centers of triangles tri-homological $M_1N_1P_1$ and $M_2N_2P_2$ are collinear (these belong to Euler's line of triangle $ABC$).

**Definition 39**

We say that the triplet of triangles $T_1T_2T_3$ is tri-homological if any two triangle from the triplet are tri-homological.

**Theorem 36 (**Gh. D. Simionescu)

If the triangles $T_1, T_2$ are tri-homological and $T_3$ is the triangle formed by the homology axes of triangles $(T_1, T_2)$, then

i) The triplet $(T_1, T_2, T_3)$ is tri-homological

ii) The homological axes of any pairs of triangles from the triplet are the sides of the other triangle.

**Proof**

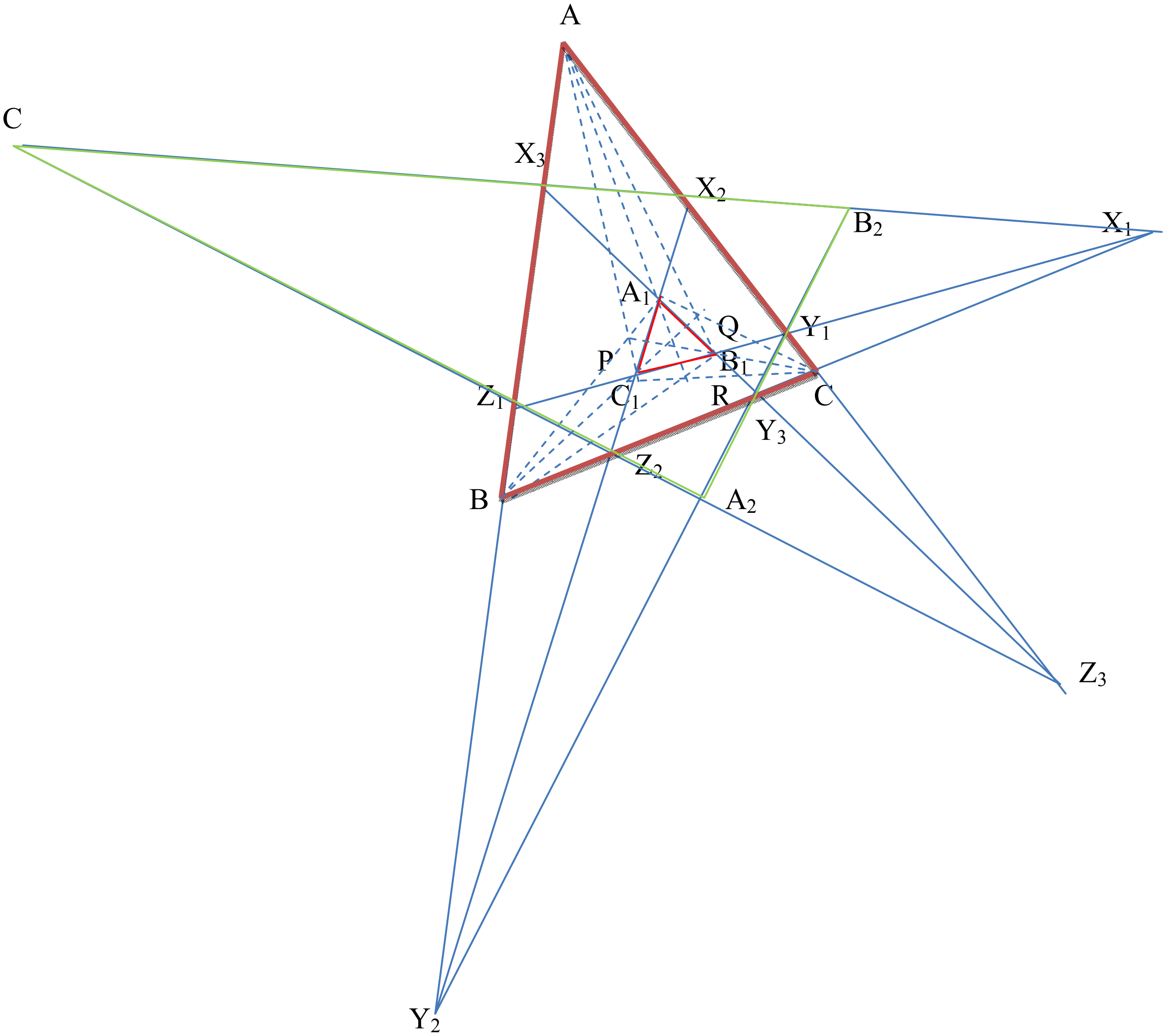

Fig. 53

Let $T_1, T_2$ the triangles $ABC$ and $A_1B_1C_1$ (see figure 53) tri-homological. We noted $X_1, X_2, X_3$ the homological axis of these triangles that corresponds to the homological center $R$, $\{R\} = AA_1 \cap BB_1 \cap CC_1$, $Y_1, Y_2, Y_3$ corresponding to the homological axis of the homological center $P$, $\{P\} = AC_1 \cap BA_1 \cap CB_1$, and $Z_1, Z_2, Z_3$ the homological axis of triangles $T_1, T_2$ corresponding to the homological axis of the homological center $Q$, $\{Q\} = AB_1 \cap BC_1 \cap CA_1$.
Also, we note $\{A_2\} = Z_1Z_2 \cap Y_1Y_2$, $\{B_2\} = X_1X_2 \cap Y_1Y_2$, $\{C_2\} = X_1X_2 \cap Z_1Z_2$, and let $T_3$ the triangle $A_2B_2C_2$.

If we consider triangles $T_2, T_3$ we observe that $B_1C_1 \cap B_2C_2 = \{X_1\}$, $A_1B_1 \cap A_2B_2 = \{Y_3\}$, $A_1C_1 \cap A_2C_2 = \{Z_2\}$.

The points $X_1, Z_2, Y_3$ belong to line $BC$ therefore are collinear and consequently the triangles $T_2, T_3$ are homological. Analyzing the same triangles we observe that $B_1C_1 \cap A_2B_2 = \{Y_1\}$, $A_1C_1 \cap B_2C_2 = \{X_2\}$, $A_1B_1 \cap A_2C_2 = \{Z_3\}$. The points $X_2, Y_1, Z_3$ are collinear being on the line $AC$, Therefore the triangles $T_2, T_3$ are double homological. From Rosanes' theorem or directly, it results that $(T_2, T_3)$ are tri-homological, the third homological axis being $AB$.

Similarly, if we consider the triangles $(T_1, T_3)$ will find that these are tri-homological.

**Lemma 1**

In triangle $ABC$, $AA'$ and $AA''$ are isotonic Cevian. Let $M \in (AA')$ and $N \in (AA'')$ such that $MN$ is parallel with $BC$. We note $\{P\} = CN \cap AB$ and $\{Q\} = BM \cap AC_2$. Prove that $PQ \parallel BC$

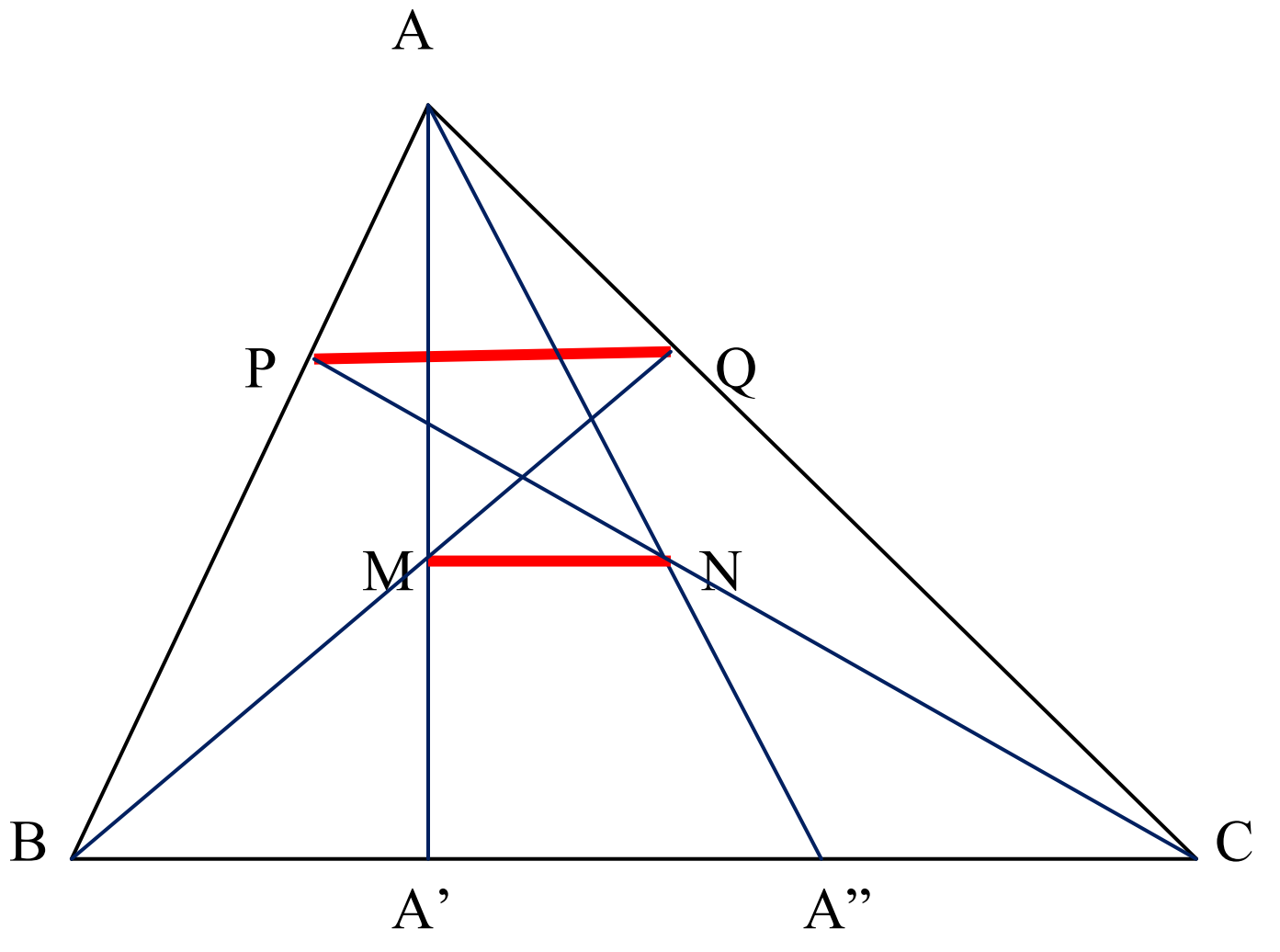

Fig. 54

**Proof**

We'll apply the Menelaus' theorem in triangles $AA'C$ and $AA''B$ for the transversals $B-M-Q$ respectively $C-N-P$ we have

$$\frac{BA'}{BC} \cdot \frac{MA}{MA'} \cdot \frac{QC}{QA} = 1$$

$$\frac{CA''}{CB}\cdot\frac{NA}{NA''}\cdot\frac{PB}{PA}=1$$

Therefore

$$\frac{BA'}{BC}\cdot\frac{MA}{MA'}\cdot\frac{QC}{QA}=\frac{CA''}{CB}\cdot\frac{NA}{NA''}\cdot\frac{PB}{PA}$$

From here and taking into account that $BA'=CA''$ and $\frac{MA}{MA'}=\frac{NA}{NA''}$ it results $\frac{QC}{QA}=\frac{PB}{PA}$ with Menelaus' theorem we obtain $PQ \parallel BC$.

**Theorem 37** (Caspary)

If $X, Y$ are points isotomic conjugate in a triangle $ABC$ and the parallels constructed through $X$ to $BC$, $CA$ respective $AB$ intersect $AY, BY, CY$ respectively in $A_1, B_1$ and $C_1$ then the triangles $ABC$ and $A_1B_1C_1$ are tri-homological triangles.

**Proof**

We note $AA', BB', CC'$ the Cevians concurrent in $X$ and $AA'', BB'', CC''$ their isotonic. See figure 55

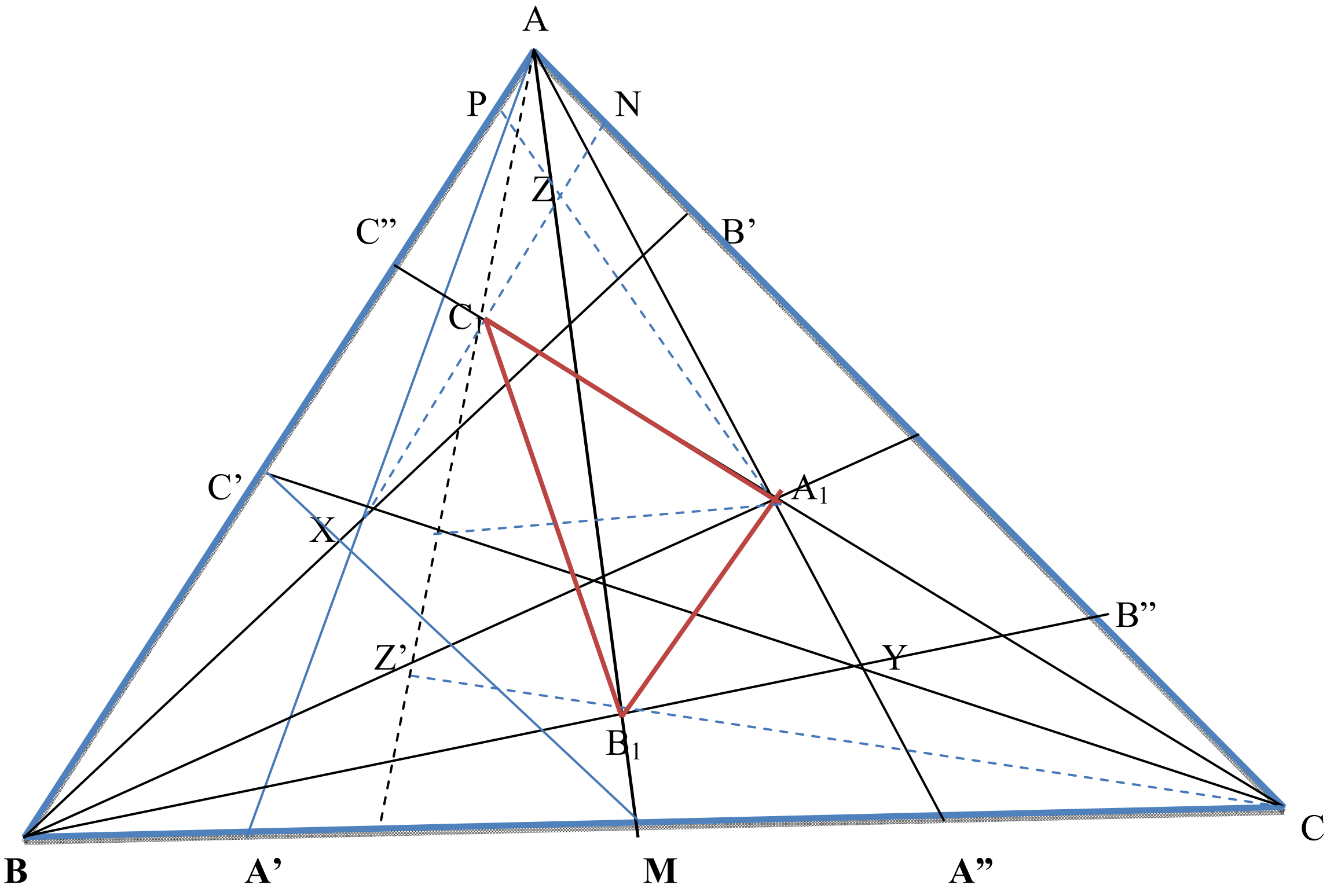

Fig. 55

The vertexes of triangle $A_1B_1C_1$ are by construction on the $AA'', BB'', CC''$, therefore $Y$ is the homological center of triangles $ABC$ and $A_1B_1C_1$.

We note with $M, N, P$ the intersections of the lines $AB_1, BC_1, CA_1$ respectively with $BC, CA, AB$, using Lemma 1 we have that $MC' \parallel AC,\ NA' \parallel AB,\ PB' \parallel BC$, therefore

$$\frac{C'B}{C'A}=\frac{MB}{MC},\ \frac{A'C}{A'B}=\frac{NC}{NA},\ \frac{B'A}{B'C}=\frac{PA}{PB}$$

Because $AA', BB', CC'$ are concurrent, from the Ceva's theorem it results

$$\frac{C'B}{C'A}\cdot\frac{A'C}{A'B}\cdot\frac{B'A}{B'C}=-1$$

also

$$\frac{MB}{MC}\cdot\frac{NC}{NA}\cdot\frac{PA}{PB}=-1$$

which shows that the Cevians $AB_1, BC_1, CA_1$ are concurrent in a point $Z$ and consequently the triangles $ABC$ and $A_1B_1C_1$ are homological.

Similarly it can be proved that the Cevians $AC_1, BA_1, CB_1$ are concurrent in a point $Z'$. Therefore the triangles $ABC$ and $A_1B_1C_1$ are tri-homological, the homology centers being the points $Y, Z, Z'$.

**Remark 41**

The triangle $A_1B_1C_1$ from the Caspary's theorem is called the first Caspary triangle. The triangle $A_2B_2C_2$ analog constructed to $A_1B_1C_1$ drawing parallels to the sides of the triangle $ABC$ through the point $Y$ is called the second triangle of Caspary.

## 1.3. Tri-homological equilateral triangles which have the same center

In this section will enounce a lemma regarding the tri-homology of equilateral triangles inscribed in another equilateral triangle, and then using this lemma we'll prove a theorem accredited to Dan Barbilian, a Romanian mathematician (1895-1961)

**Lemma 2**

Let $A_1B_1C_1$ an equilateral triangle with a center $O$ and $A_2B_2C_2$, $A_3B_3C_3$ equilateral triangles inscribed in $A_1B_1C_1$ ($A_2, A_3 \in (B_1C_1), B_2, B_3 \in (A_1C_1), C_2, C_3 \in (A_1B_1)$)

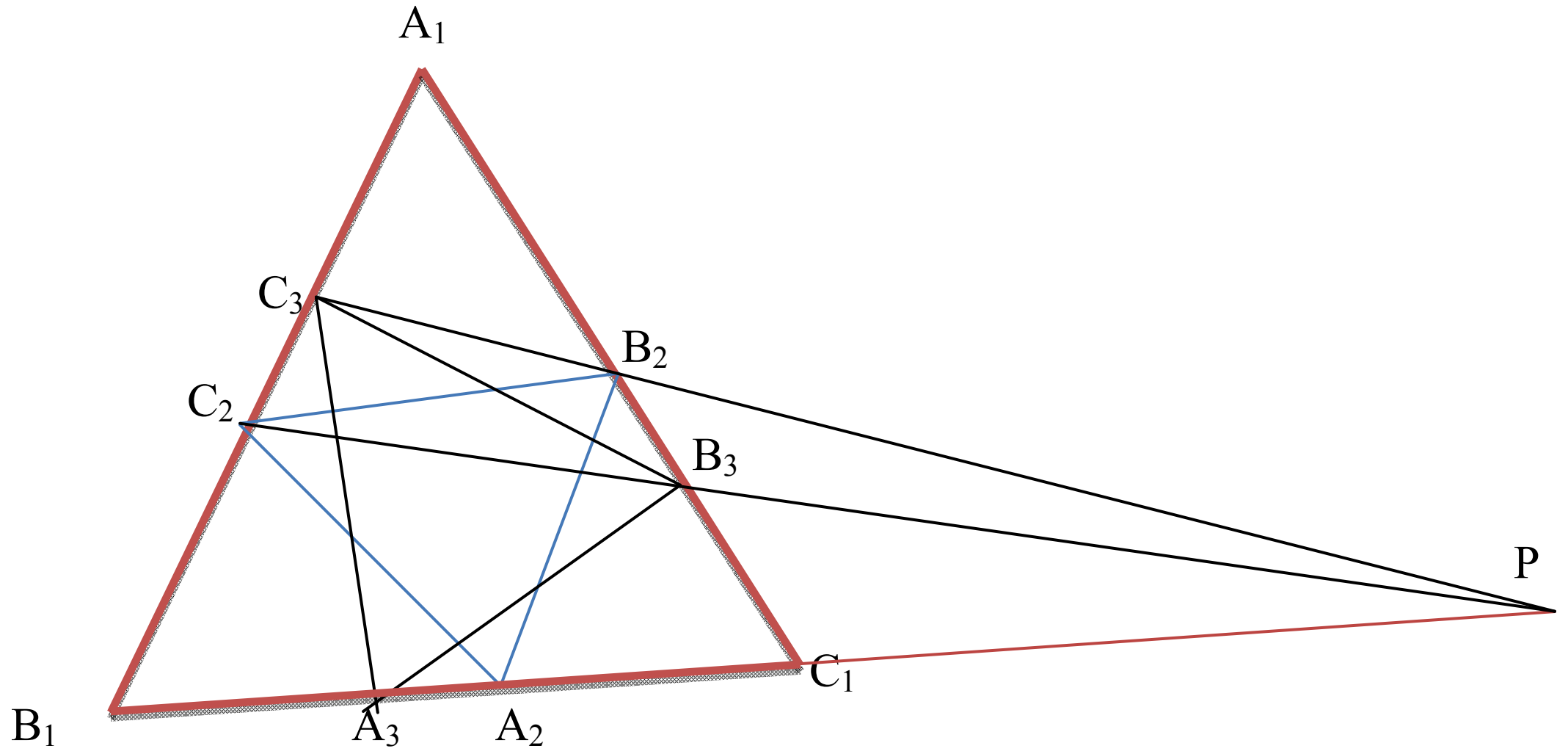

Fig. 56

Then

(i) The triangles $A_1B_1C_1$ and $A_2B_2C_2$ are tri-homological

(ii) The triangles $A_1B_1C_1$ and $A_3B_3C_3$ are tri-homological

(iii) The triangles $A_2B_2C_2$ and $A_3B_3C_3$ are tri-homological.

**Proof**

i. We observe that $A_1B_2, B_1A_2, C_1C_2$ are concurrent in $C_1$, therefore $A_1B_1C_1$ and $B_2A_2C_2$ are homological with the homological center in $C_1$. Similarly results that $B_1$ is the homological center of triangles $A_1B_1C_1$ and $C_2B_2A_2$ and $A_1$ is the homological center of triangles $A_1B_1C_1$ and $A_2C_2B_1$.

ii. Similar to (i).

iii. We'll prove that the triangles $A_1B_2C_2$, $B_1C_2A_2$, and $C_1A_2B_2$ are congruent. Indeed if $m(\sphericalangle A_1B_2C_2) = \alpha$ then $m(\sphericalangle A_1C_2B_2) = 120° - \alpha$ and $m(\sphericalangle B_2C_2A_2) = \alpha$. Therefore $\sphericalangle A_1B_2C_2 \equiv \sphericalangle B_1C_2A_2$ and $\sphericalangle A_1C_2B_2 \equiv \sphericalangle B_1A_2C_2$; we know that $B_2C_2 = C_2A_2$, it results that $\Delta A_1B_2C_2 \equiv \Delta B_1C_2A_2$. Similarly we find that $\Delta B_1C_2A_2 \equiv \Delta C_1A_2B_2$. From these congruencies we retain that

$$A_1B_2 = B_1C_2 = C_1A_2 \qquad (104)$$

In the same way we establish that $\Delta A_1C_3B_3 \equiv \Delta B_1A_3C_3 \equiv \Delta C_1B_3A_3$ with the consequence

$$A_1C_3 = B_1A_3 = C_1B_3 \qquad (105)$$

We will prove that the lines $C_3B_2, B_3C_2, A_3A_2$ are concurrent.

We note $C_3B_2 \cap A_3A_2 = \{P\}$ and $C_2C_3 \cap A_2A_3 = \{P'\}$.

The Menelaus' theorem applied in the triangle $A_1B_1C_1$ for the transversals $P, B_2, C_3$ and $P', B_3, C_2$ provides us with the relations:

$$\frac{PC_1}{PB_1} \cdot \frac{C_3B_1}{C_3A_1} \cdot \frac{B_2A_1}{B_2C_1} = 1 \quad (106)$$

$$\frac{P'C_1}{P'B_1} \cdot \frac{C_2B_1}{C_2A_1} \cdot \frac{B_3A_1}{B_3C_1} = 1 \quad (107)$$

From (106) and (105) we find that $C_3B_1 = A_1B_3, B_2A_1 = C_2B_1, C_3A_1 = B_3C_1, B_2C_1 = C_2A_1$, we come back to the relations (106) and (107) and we find that $\frac{PC_1}{PB_1} = \frac{P'C_1}{P'B_1} \cdot$ and from here we see that $P = P'$ with the consequence that $B_2C_3 \cap B_3C_2 \cap A_2A_3 = \{P\}$.

Similarly we prove that the lines $A_2B_3, B_2A_3, C_2C_3$ are concurrent in $Q$, therefore the triangles $A_2B_2C_2$ and $B_3A_3C_3$ are homological and the lines $A_2C_3, B_2B_3, C_2A_3$ are concurrent in a point $R$, therefore the triangles $A_2B_2C_2$ and $A_3B_3C_3$ are homological.

**Remark 42**

It can be proved that the triangles $A_1B_1C_1$, $A_2B_2C_2$ and $A_3B_3C_3$ have the same center $O$. If we note $\{A_4\} = B_2C_3 \cap C_2A_3; \{B_4\} = A_2B_3 \cap A_3C_2; \{C_4\} = A_2B_3 \cap C_3B_2$ then triangle $A_4B_4C_4$ is equilateral with the same center $O$ and from the Lemma it results it is homological with each of the triangles $A_1B_1C_1$, $A_2B_2C_2$ and $A_3B_3C_3$.

**Theorem 38** (D. Barbilian – 1930)

If $A_1B_1C_1$ and $A_2B_2C_2$ are two equilateral triangles having the same center $O$ and the vertexes notation is in the same rotation sense, then the triangles are three times homological as follows:

$(A_1B_1C_1), (C_2B_2A_2), (A_1B_1C_1, B_2A_2C_2), (A_1B_1C_1, A_2C_2B_2)$

**Proof**

We note

$$\{A_3\} = B_1B_2 \cap C_1C_2,$$
$$\{B_3\} = A_1A_2 \cap C_1C_2,$$
$$\{C_3\} = A_1A_2 \cap B_1B_2.$$

See figure 57.

We notice that

$$\Delta OB_1C_2 \equiv \Delta OC_1C_2 \equiv \Delta OA_1B_2 \text{ (SAS)}$$

it results

$$B_1C_2 = C_2A_2 = A_1B_2 \quad (108)$$

also

$\Delta OC_1C_2 \equiv \Delta OB_1B_2 \equiv \Delta OA_1A_2$ (SAS)

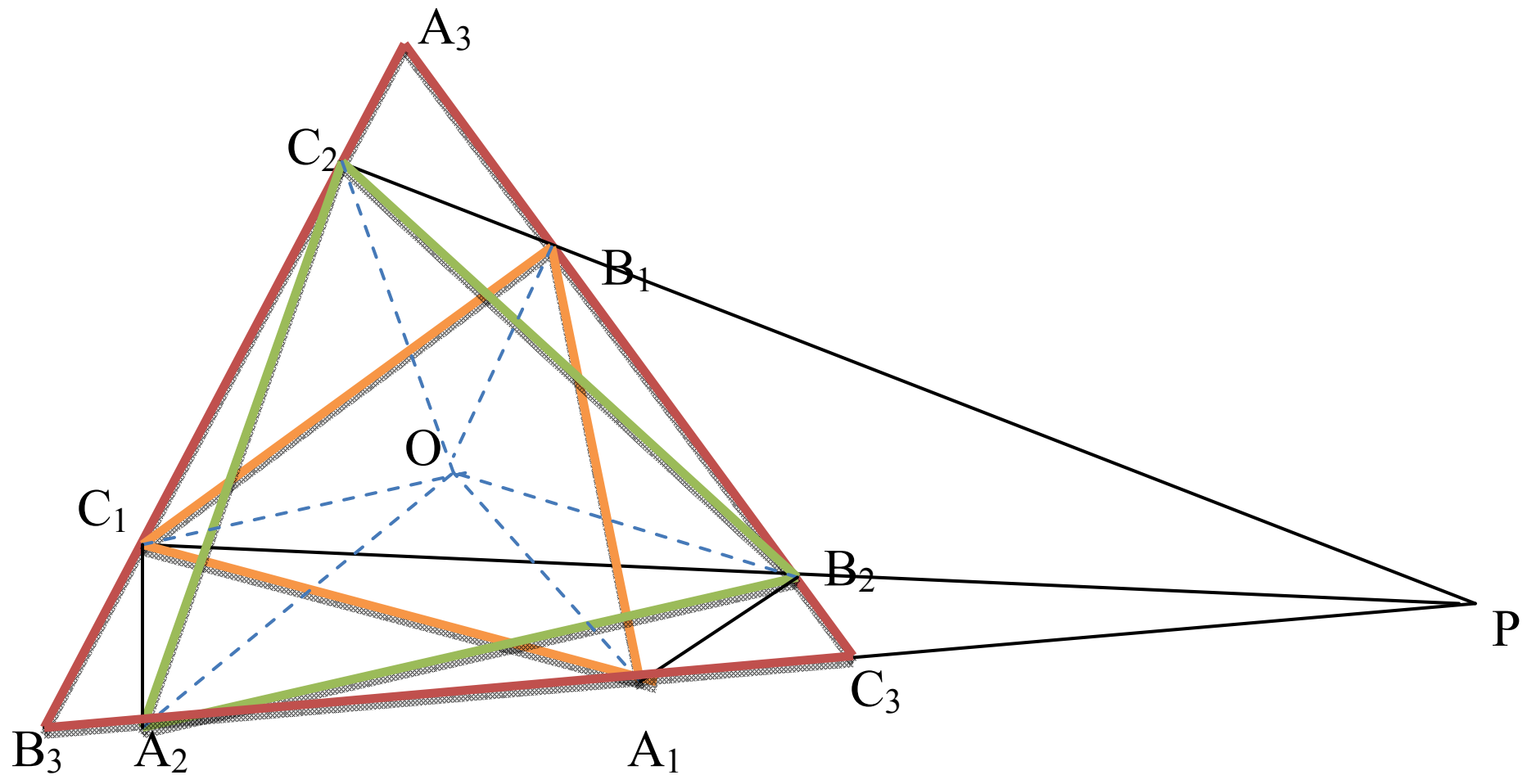

Fig. 57

it results

$$C_1C_2 = B_1A_2 = A_1A_2 \qquad (109)$$

We have also $\Delta B_1C_1C_2 \equiv \Delta C_1A_1A_2 \equiv \Delta B_1A_1A_2$ (SSS)

We obtain that $\sphericalangle B_1C_2C_1 \equiv \sphericalangle C_1A_2A_1 \equiv \sphericalangle B_1B_2A_1$.

From what we proved so far it result that $\Delta A_3B_1C_2 \equiv \Delta B_3C_1A_2 \equiv \Delta C_3A_1B_2$ with the consequence $\sphericalangle A_3 \equiv \sphericalangle B_3 \equiv \sphericalangle C_3$, which shows that the triangle $A_3B_3C_3$ is equilateral.

Applying lemma for the equilateral triangles $A_1B_1C_1$ and $A_2B_2C_2$ inscribed in the equilateral triangle $A_3B_3C_3$ it result that the triangles $A_1B_1C_1$ and $A_2B_2C_2$ are tri-homological.

## 1.4. The Pappus' Theorem

**Theorem 39** (Pappus – 3rd century)

If the vertexes of a hexagon are successively on two given lines, then the intersections of the opposite sides are collinear.

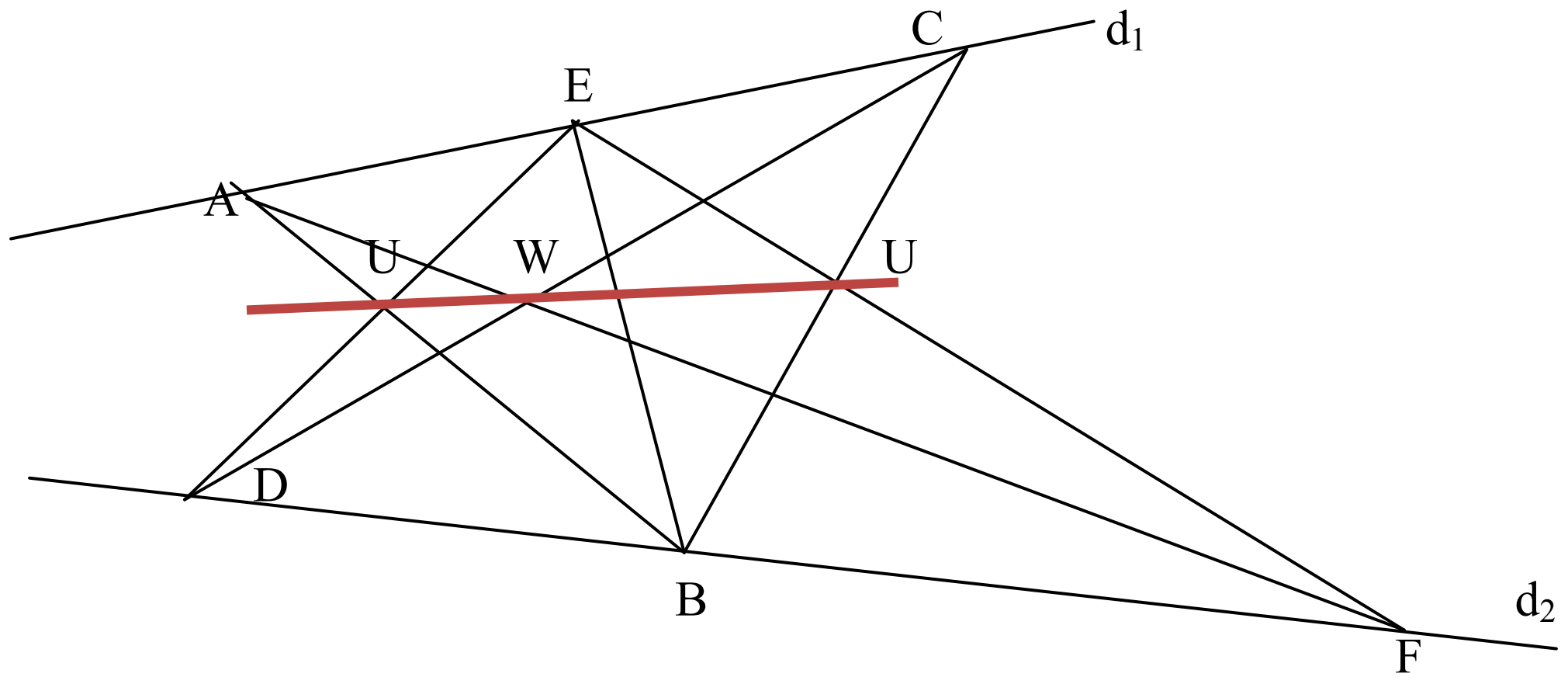

Fig. 58

**Proof**

Let $ABCDEF$ a hexagon with the vertexes $A, C, E$ on line $d_1$ and vertexes $B, D, F$ on the line $d_2$ (see figure 58).

We note $\{U\} = AB \cap DE$; $\{V\} = BC \cap EF$ and $\{W\} = CD \cap FA$.

The triangle determined by the intersections of the lines $AB, EF, CD$ and he triangle determined by the intersections of the lines $BC, DE, FA$ are twice homological having as homological axes the lines $d_1, d_2$.

In accordance with theorem 24 these triangles are tri-homological, the third homological axis is the line to which belong the points $U, V, W$.

**Remark 43**

The Pappus' theorem can be directly proved using multiple time the Menelaus' theorem.

## 1.5. The duality principle

A line and a point are called incidental if the point belongs to the line or the line passes through the point.

**Definition 40**

A duality is a transformation which associates bi-univoc to a point a line. It is admitted that this correspondence preserves the incidental notion; in this mode to collinear points correspond concurrent lines and reciprocal.

If it I considered a theorem T whose hypothesis implicitly or explicitly appear points, lines, incident and it is supposed that its proof is completed, then if we change the roles of the points with the lines reversing the incidence, it is obtained theorem T' whose proof is not necessary.

**Theorem 40** (The dual of Pappus' theorem)

If we consider two bundles each of three concurrent lines $S(a,b,c)$, $S'(a',b',c')$ such that the lines $a, b'$ and $b, a'$ intersect in the points $C_1, C_2$; $a, c'$ and $c, a'$ intersect in the points $B_1, B_2$ and the lines $b, c'$ and $c, b'$ intersect in $A_1, A_2$, then the lines $A_1A_2, B_1B_2, C_1C_2$ are concurrent.

**Proof**

Analyzing the figure 59 we observe that it is obtain by applying the duality principle to Pappus' theorem.

Indeed, the two bundles $S(a,b,c)$, $S'(a',b',c')$ correspond to the two triplets of vertexes of a hexagon situated on the lines $d_1, d_2$ to which correspond the points $S$ and $S'$.

The Pappus' theorem proves the collinearity of $U, V, W$ which correspond to the concurrent lines $A_1A_2, B_1B_2, C_1C_2$.

**Observation 28**

The dual of Pappus' theorem can be formulated in an important particular case.

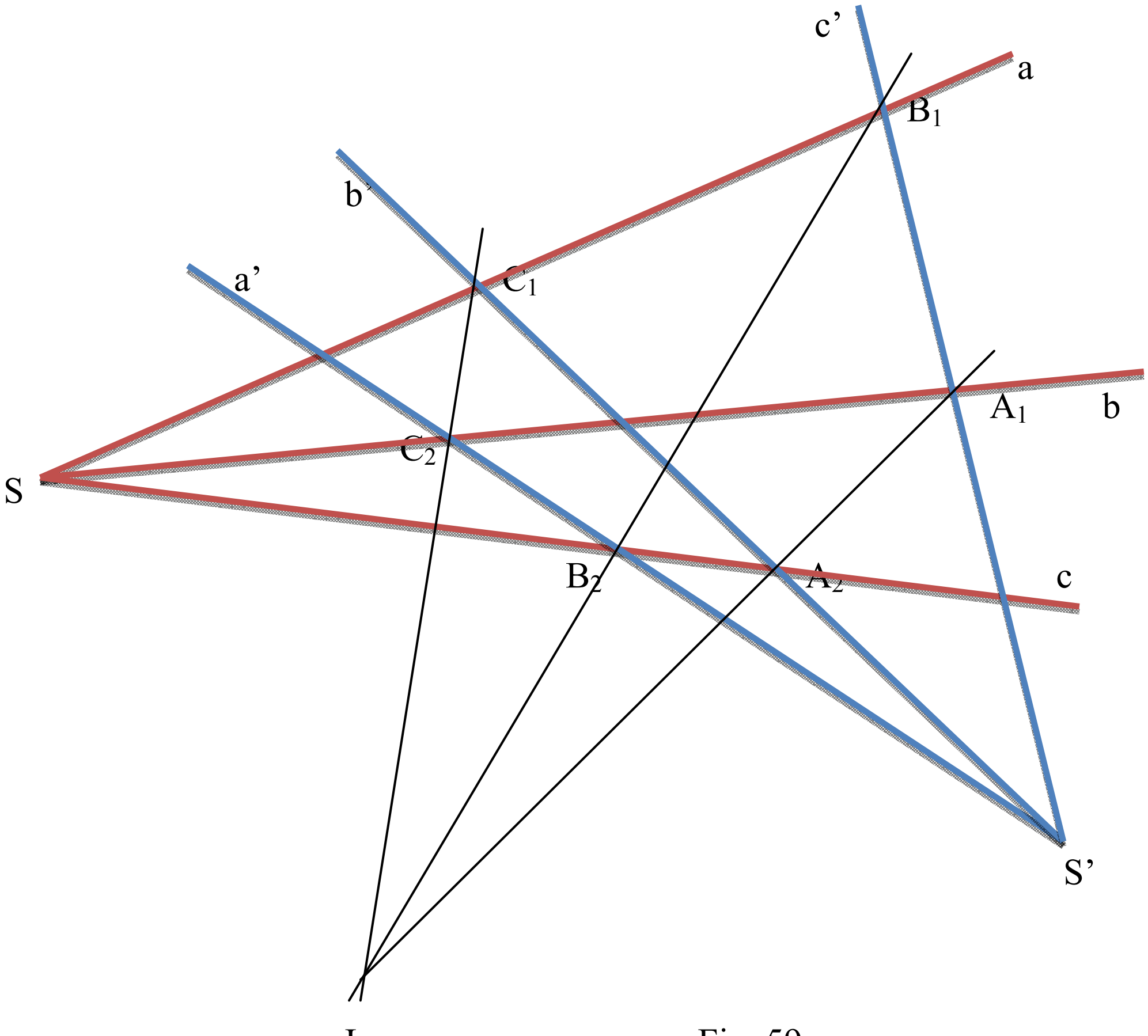

Fig. 59

**Theorem 41**

We consider a complete quadrilateral and through the vertexes $E, F$ we construct two secants, which intersect $AD, BC$ in the points $E_1, E_2$ and $AB, CD$ in the points $F_1, F_2$.

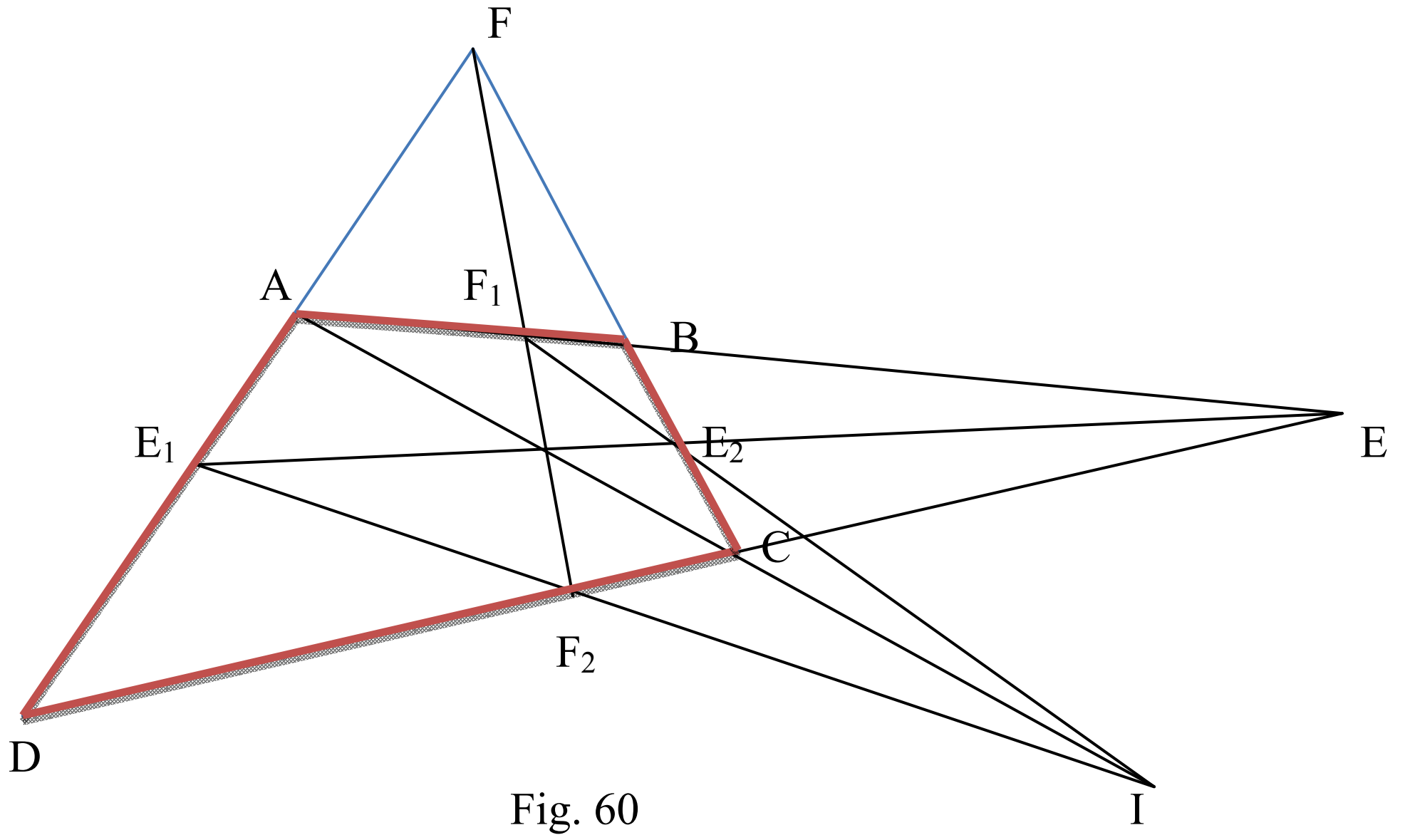

Fig. 60

Then the lines $E_1F_2, F_1E_2$ intersect on the diagonal $AC$ and the lines $E_1E_2, F_1F_2$ intersect on the diagonal $BD$.

Indeed, this theorem is a particular case of the precedent theorem. It is sufficient to consider the bundles of vertexes $E, F$ and of lines $(CD, E_1E_2, AB)$ respectively $(AD, F_1F_2, BC)$ see figure 60, and to apply theorem 28.

From what we proved so far, it result that the triangles $BF_1E_2$ and $DF_2E_1$ are homological, therefore $BD, F_1F_2, E_2E_1$ are concurrent.

**Remark 44**

The dual of Pappus' theorem leads us to another proof for theorem 34 (Rosanes).

We prove therefore that two homological triangles are tri-homological.

We consider the triangles $ABC, A'B'C'$ bi-homological. Let $S$ and $S'$ the homology centers: $S$ the intersection of the lines $AA', BB', CC'$ and $\{S'\} = AB' \cap BC' \cap CA'$( see figure 61). We'll apply theorem 41 for the bundles $S(AA', BB', CC')$ and $S'(C'B, A'C, B'A)$.

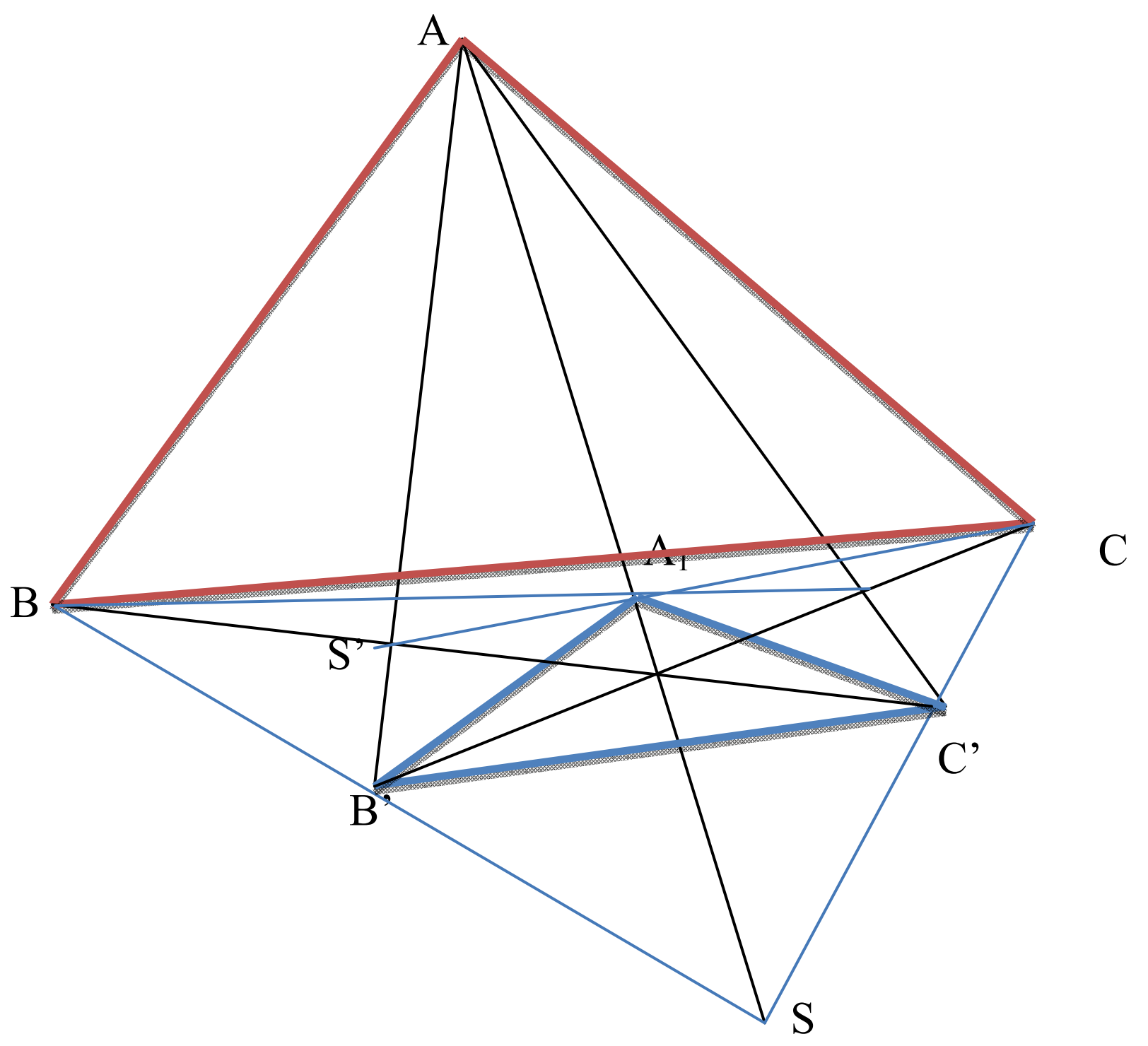

Fig. 61

We observe that

$AA' \cap A'C = \{A'\}$,

$C'B \cap BB' = \{C'\}$,

$BB' \cap B'A = \{B'\}$,

$CC' \cap A'C = \{C\}$,

$$AA' \cap B'A = \{A\},$$

$$CC' \cap C'B = \{C'\}$$

Therefore the lines $BA', CB', AC'$ are concurrent which shows that the triangles $ABC$ and $C'A'B'$ are homological, thus the triangles $ABC$ and $A'B'C'$ are tri-homological.

**Theorem 42**

In triangle $ABC$ let's consider the Cevians $AA_1, BB_1, CC_1$ in $M_1$ and $AA_2, BB_2, CC_2$ concurrent in the point $M_2$. We note $A_3, B_3, C_3$ the intersection points of the lines $(CC_1, BB_2)$, $(AA_1, CC_2)$ respectively $(BB_1, AA_2)$, and $A_4, B_4, C_4$ the intersection points of the lines $(CC_2, BB_1)$, $(AA_2, CC_1)$ respectively $(BB_2, AA_1)$, then

(i) The triangles $A_3B_3C_3$ and $A_4B_4C_4$ are homological, and we note their homological center with $P$.

(ii) The triangles $ABC$ and $A_3B_3C_3$ are homological, their homological center being noted $Q$.

(iii) The triangles $ABC$ and $A_4B_4C_4$ are homological, their homological center being noted with $R$

(iv) The points $P, Q, R$ are collinear.

**Proof**

(i)

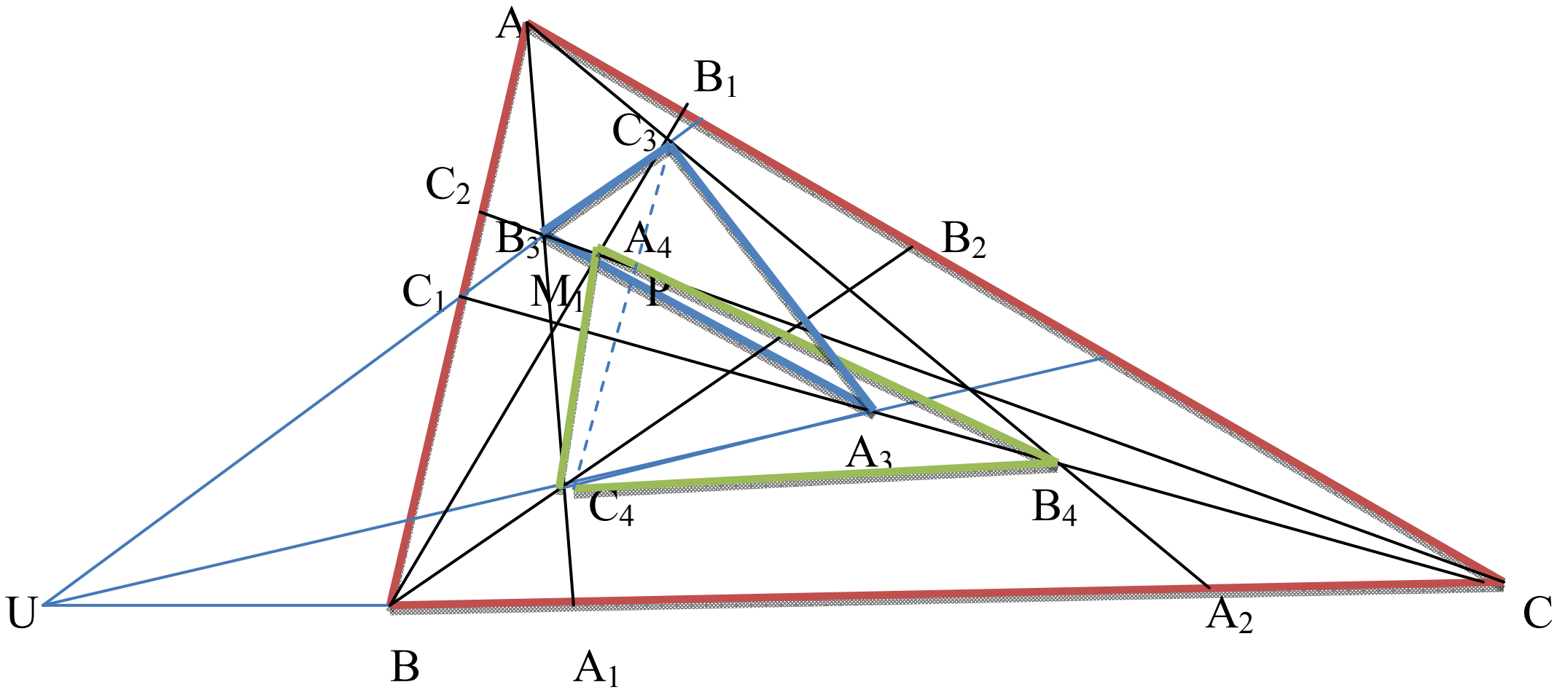

Fig. 62

Let consider the point $P$ the intersection of $A_3A_4$ and $C_3C_4$ with the sides of the hexagon $C_4M_1A_3A_4M_2C_3$, which has each three vertexes on the lines $BB_1$, $BB_2$. In conformity with Pappus' theorem the opposite lines $C_4M_1$, $A_4M_2$; $M_1A_3$, $M_2C_3$ ; $A_3A_4$ , $C_3C_4$ intersect in collinear points .

These points are $B_3, B_4$ and $P$; therefore the line $B_3B_4$ passes through $P$, and thus the triangles $A_3B_3C_3$, $A_4B_4C_4$ are homological. We note $U, V, W$ their homological axis, therefore

$$\{U\} = B_3C_3 \cap B_4C_4$$

$\{V\} = A_3C_3 \cap A_4C_4$

$\{W\} = A_3B_3 \cap A_4B_4$

(ii)

We consider the hexagon $C_4B_4M_1C_3B_3M_2$ , which each of its vertexes on $AA_1$ respectively $AA_2$.

The opposite sides $(B_3C_3, B_4C_4)$, $(C_3M_1, C_4M_2)$, $(M_1B_4, M_2B_3)$ intersect in the collinear points $U, B, C$. It results that the point $U$ is on the side $BC$ and similarly the points $V, W$ are on the sides $AC, AB$.

Consequently the triangle $ABC$ is homological with $A_3B_3C_3$.

(iii)

From the fact that $U, V, W$ are respectively on $BC, AC, AB$, from their collinearity and from the fact that $U$ is on $B_4C_4$, $V$ belongs to line $A_4C_4$, and $W$ belongs to the line $A_4B_4$, it results that the triangles $ABC$ and $A_4B_4C_4$ are homological.

(iv)

The lines $BC, B_3C_3, B_4C_4$ have $U$ as common point, we deduct that the triangles $BB_3B_4$ and $CC_3C_4$ are homological. Consequently their opposite sides intersect in three collinear points, and these points are $P, Q, R$.

**Remark 45**

The point (iv) of the precedent theorem could be proved also by applying theorem 18.

Indeed, the triangles $(ABC, A_3B_3C_3, A_4B_4C_4)$ constitute a homological triplet, and two by two have the same homological axis, the line $U, V, W$. It results that their homological centers i.e. *P, Q, R* are collinear.

# Chapter 4

## Homological triangles inscribed in circle

This chapter contains important theorems regarding circles, and certain connexions between them and homological triangles.

### 4.1. Theorems related to circles.

**Theorem 43** (L. Carnot-1803)

If a circle intersects the sides $BC, AC, AB$ of a given triangle $ABC$ in the points $A_1A_2; B_1B_2; C_1C_2$ respectively, then the following relation takes place

$$\frac{A_1B}{A_1C}\cdot\frac{A_2B}{A_2C}\cdot\frac{B_1C}{B_1A}\cdot\frac{C_2C}{B_2A}\cdot\frac{C_1A}{C_1B}\cdot\frac{C_2A}{C_2B}=1 \qquad (110)$$

**Proof**

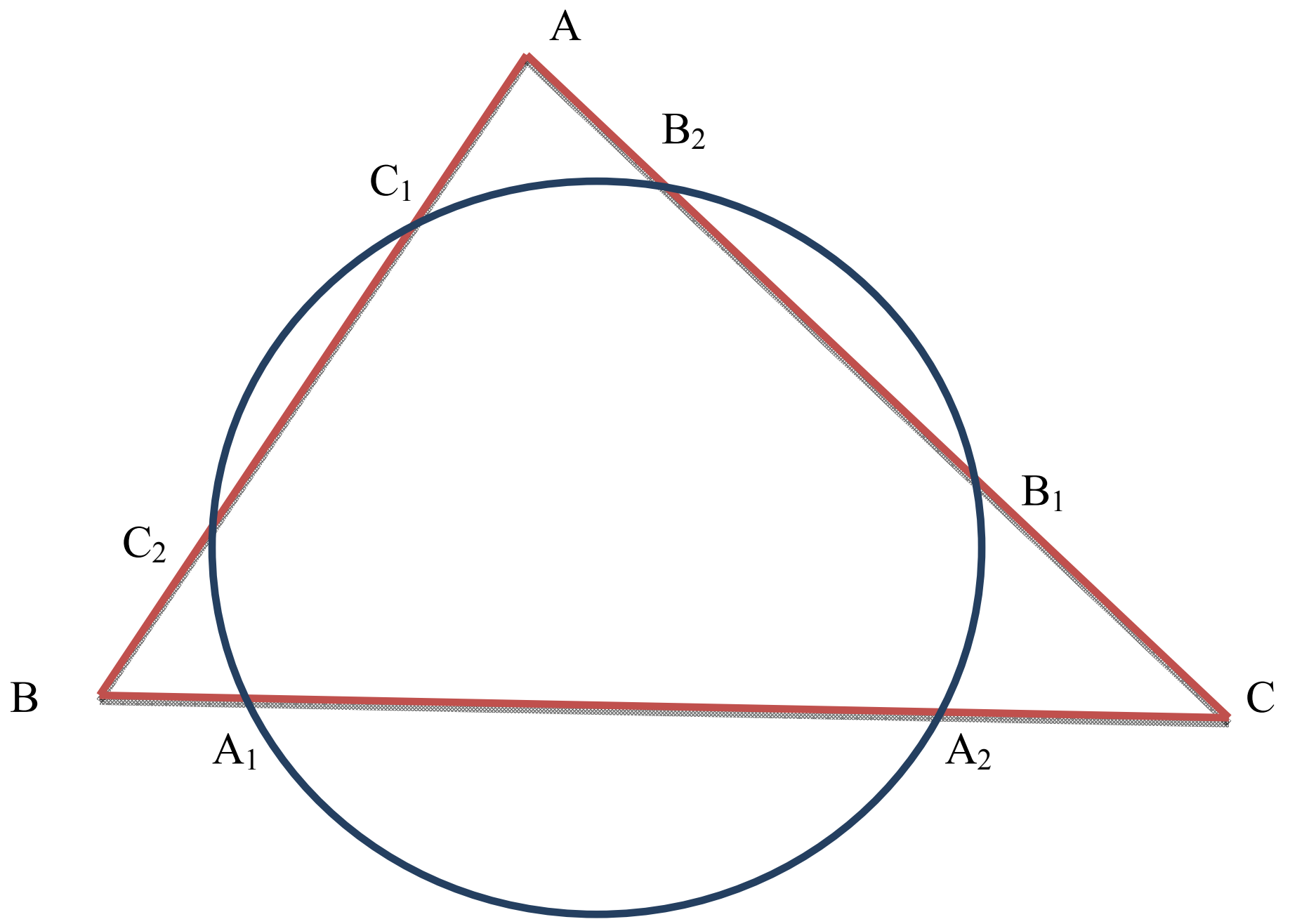

Fig. 63

We consider the power of the points $A, B, C$ in rapport to the given circle, see figure 63. We obtain:

$$AC_1\cdot AC_2 = AB_1\cdot AB_2 \qquad (111)$$

$$BA_1\cdot BA_2 = BC_1\cdot BC_2 \qquad (112)$$

$$CA_1\cdot CA_2 = CB_1\cdot CB_2 \qquad (113)$$

From these relations it results relation (110)

**Remark 46**

From Carnot's relation we observe that if

$$\frac{A_1B}{A_1C}\cdot\frac{B_1C}{B_1A}\cdot\frac{C_1A}{C_1B}=-1$$

then

$$\frac{A_2B}{A_2C}\cdot\frac{B_2C}{B_2A}\cdot\frac{C_2A}{C_2B}=-1$$

and this relation proves the following theorem:

**Theorem 44** (Terquem)

If we construct a circle through the legs of three Cevians concurrent in a triangle, then it will intersect the legs of other concurrent Cevians.

**Theorem 45** (Pascal – 1640)

The opposite sides of a hexagon inscribed in a circle intersect in collinear points.

**Proof**

Let $Q_1R_1Q_2R_2Q_3R_3$ the inscribed hexagon in a circle (see figure 64).

We note

$$\{A_1\}=R_3Q_3\cap Q_2R_2$$
$$\{B_1\}=Q_1R_1\cap R_3Q_3$$
$$\{C_1\}=R_1Q_1\cap R_2Q_2$$
$$\{A_2\}=Q_1R_3\cap R_1Q_2$$
$$\{B_2\}=R_1Q_2\cap R_2Q_3$$
$$\{C_2\}=R_2Q_3\cap Q_1R_3$$

See figure 64.

Applying the Carnot's theorem we have:

$$\frac{Q_1B_1}{Q_1C_1}\cdot\frac{R_1B_1}{R_1C_1}\cdot\frac{Q_2C_1}{Q_2A_1}\cdot\frac{R_2C_1}{R_2A_1}\cdot\frac{Q_3A_1}{Q_3B_1}\cdot\frac{R_3A_1}{R_3B_1}=1$$

This relation can be written:

$$\frac{Q_1B_1\cdot Q_2C_1\cdot Q_3A_1}{Q_1C_1\cdot Q_2A_1\cdot Q_3B_1}=\frac{R_1C_1\cdot R_2A_1\cdot R_3B_1}{R_1B_1\cdot R_2C_1\cdot R_3A_1}$$

Taking into consideration this relation, it results that the triangles $A_1B_1C_1$ and $A_2B_2C_2$ are homological (the lines $A_1A_2, B_1B_2, C_1C_2$ are concurrent), therefore the opposite sides of the hexagon $Q_1R_1Q_2R_2Q_3R_3$ intersect in collinear points.

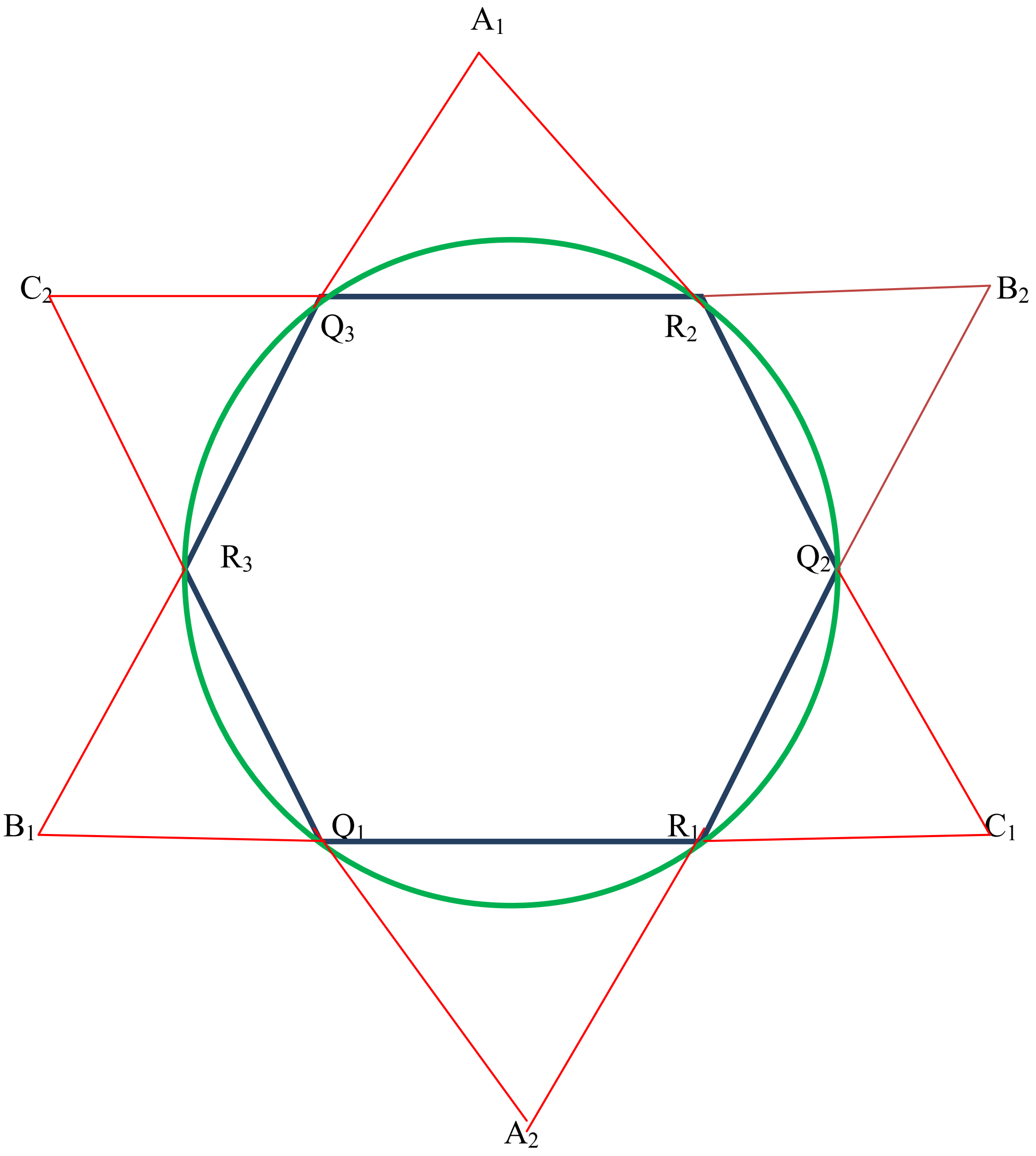

Fig. 64

**Remark 47**

The Pascal's theorem is true also when the inscribed hexagon is non-convex. Also, Pascal's theorem remains true when two or more of the hexagon's vertexes coincide. For example two of the vertexes $C$ and $C'$ of the inscribed hexagon $ABCC'DE$ coincide, then we will substitute the side $CD$ with the tangent in $C$ to the circumscribed circle.\

**Theorem 46**

If $ABCDE$ is an inscribed pentagon in a circle, $M, N$ are the intersection points of the sides $AB$ and $CD$ respectively $BC$ and $DE$ and $P$ is the intersection point of the tangent constructed in $C$ to the pentagon's circumscribed circle with the side $DE$, then the points $M, N, P$ are collinear.

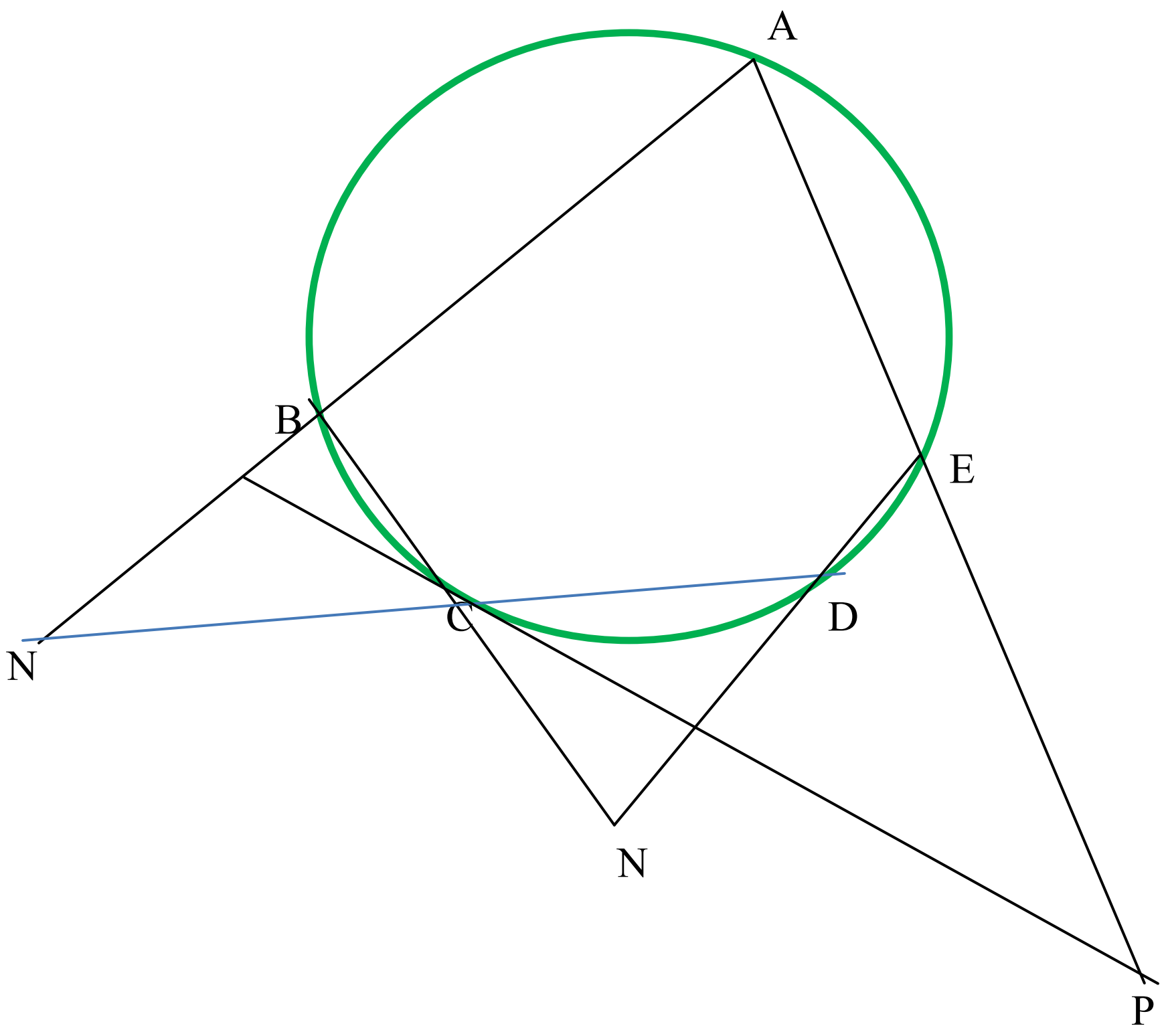

Fig. 65

**Observation 29**

If in an inscribed hexagon $AA'BCC'D$ we suppose that two pairs of vertexes coincide, the figure becomes an inscribed quadrilateral, which we can consider as a degenerated hexagon. The sides being $AB$, $BC$, $CC'$ - tangent in $C$, $C'D \to CD$, $DA' \to DA$, $AA' \to$ tangent.

**Theorem 47**

In a quadrilateral inscribed in a circle the opposite sides and the tangents in the opposite vertexes intersect in four collinear points.

**Remark 48**

This theorem can be formulated also as follows.

**Theorem 48**

If $ABCDEF$ is a complete quadrilateral in which $ABCD$ is inscribed in a circle, then the tangents in $A$ and $C$ and the tangents in $B$ and $D$ to the circumscribed circle intersect on the quadrilateral's diagonal $EF$.

**Observation 30**

The figure corresponds to theorem

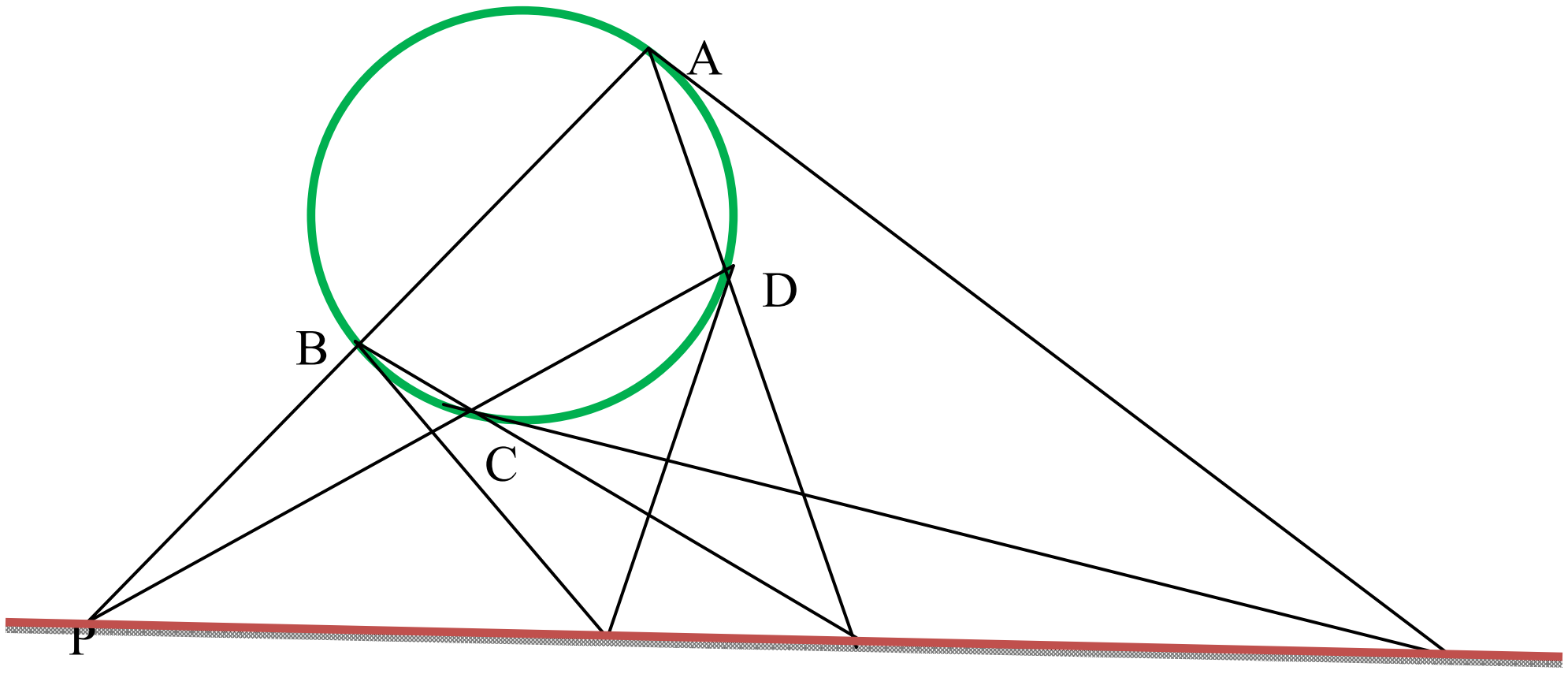

Fig. 66

**Remark 49**

If we apply the theorem of Pascal in the degenerated hexagon AA'BB'CC' where the points A,A'; B,B'; C, C' coincide and the sides $AA', BB', CC'$ are substituted with the tangents constructed in $A, B, C$ to the circumscribed to triangle $ABC$, we obtain theorem 7 (Carnot)

**Theorem 49** (Chasles - 1828)

Two triangles reciprocal polar with a circle are homological.

**Proof**

Let $ABC$ and $A_1B_1C_1$ two reciprocal polar triangles in rapport with the circle of radius $r$ (see figure 67).

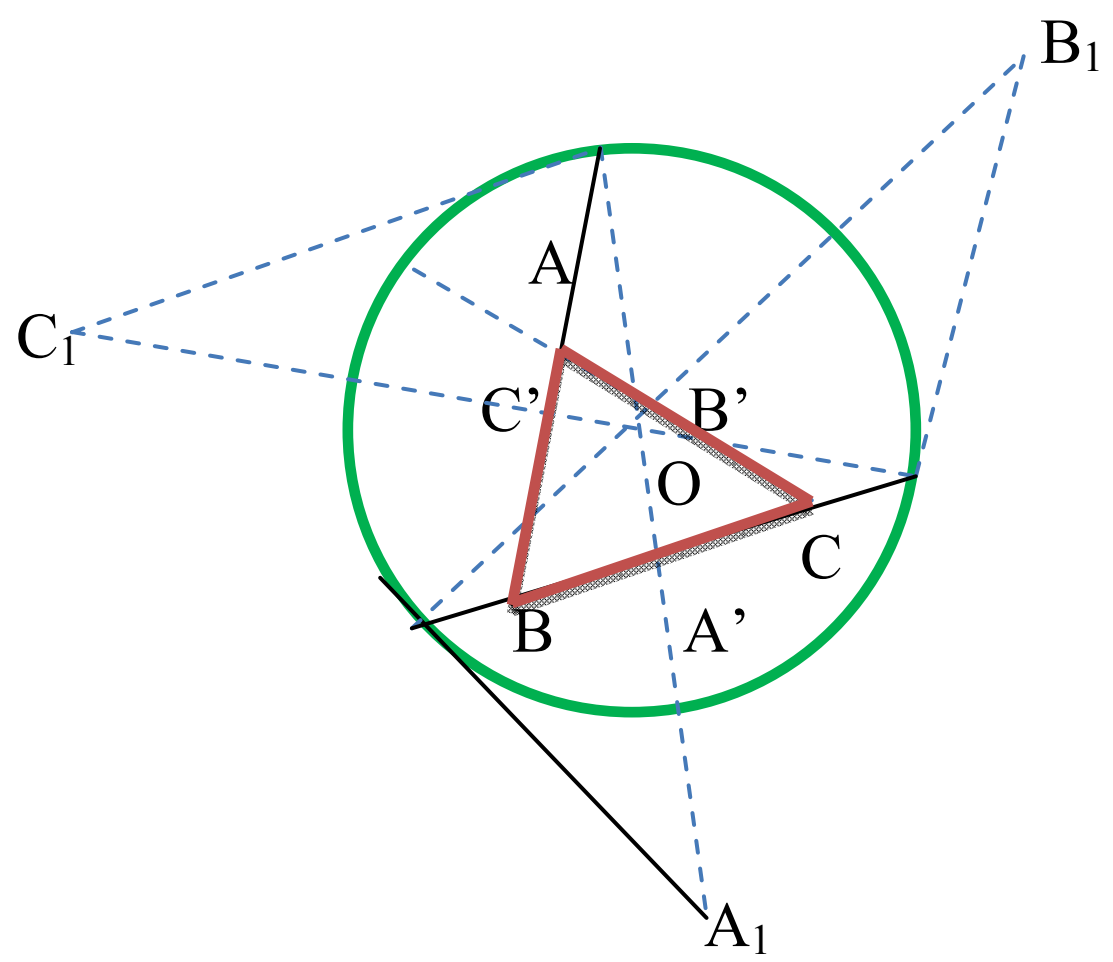

Fig. 67

We consider that $BC$ is the polar of $A_1$, $CA$ is the polar of $B_1$ and $AB$ is the polar of $C_1$ Therefore, $OA' \cdot OA_1 = r^2$ . Also $OB' \cdot OB_1 = OC' \cdot OC_1 = r^2$ . We noted with $A', B', C'$ the orthogonal projections of the point $O$ on $BC, CA, AB$ respectively.

Applying the Coşniţă theorem (its generalization), it results that the lines $AA_1, BB_1, CC_1$ are concurrent, consequently the triangles $ABC$ and $A_1B_1C_1$ are homological.

**Remark 50**

If one considers the points $A, B, C$ on the circle of center $O$, then the sides of the triangle $A_1B_1C_1$ will be tangents in $A, B, C$ to the circumscribed circle to triangle $ABC$, and the homological center of the triangles is the Gergonne's point of the triangle $A_1B_1C_1$.

**Theorem 50** (Brianchon -1806)

If a hexagon $ABCDEF$ is circumscribed to a circle then the diagonals $AD, BE, CF$ are concurrent.

**Proof**

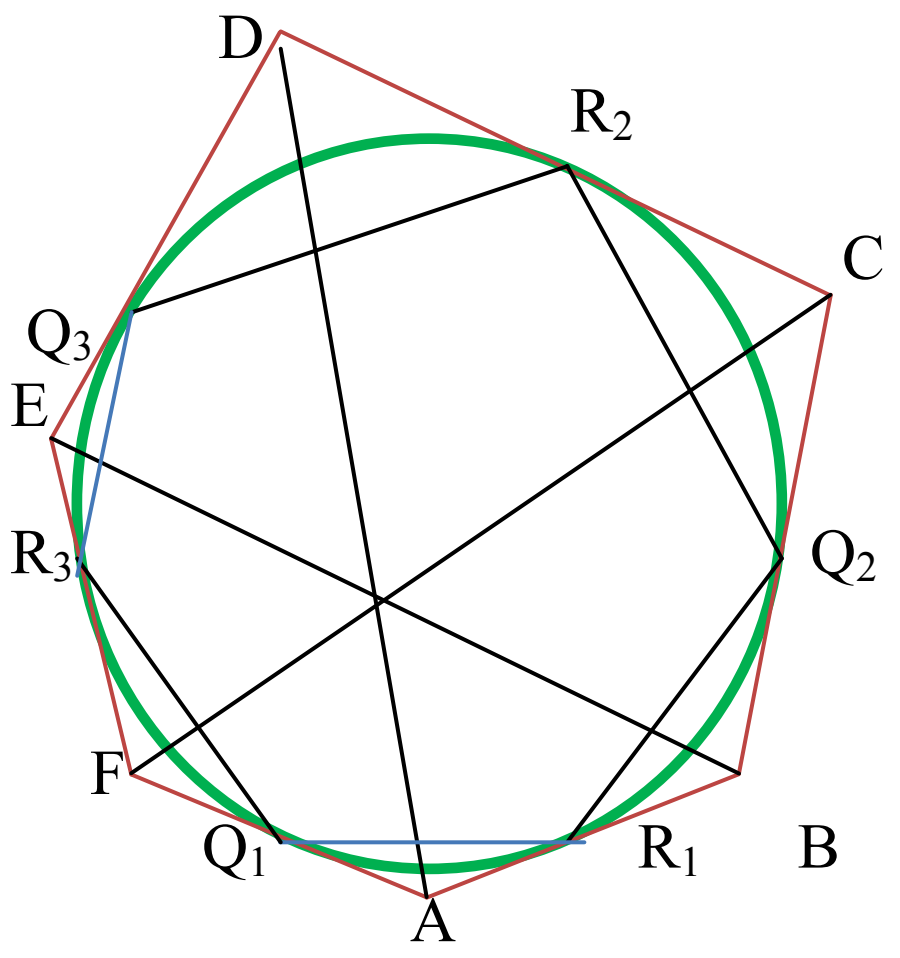

Fig 68

We will transform by duality Pascal's theorem 45 in in relation with the inscribed hexagon $Q_1R_1Q_2R_2Q_3R_3$ in rapport to the circle in which the hexagon is inscribed. Therefore to the line $Q_1R_1$ corresponds the point $A$ of intersection of tangents constructed in the points $R_1, Q_1$ on the circle (the polar of the points $R_1, Q_1$). Similarly, we obtain the vertexes $B, C, D, E, F$ of the hexagon $ABCDEF$ circumscribed to the given circle.

To the intersection point of the opposite sides $Q_1R_1$ and $Q_3R_3$ corresponds the line determined by the pols of these lines that is the diagonal $AD$.

Similarly we find that the diagonals and $BE$ correspond to the intersection point of the other two pairs of opposite sides.

Because the intersection points of the opposite sides of the inscribed hexagon are collinear, it will result that the polar, that is $AD, BE, CF$ are concurrent and Brianchon's theorem is proved.

**Remark 51**

The Brianchon's theorem remains true also if the hexagon is degenerate in the sense that two sides are prolonged.

In this case we can formulate the following theorem

**Theorem 51**

In a pentagon circumscribable the diagonals and the lines determined by the opposite points of tangency are concurrent

**Remark 53**

The Newton's theorem is obtained by duality transformation.

If the hexagon $ABCDEF$ from the Brianchon's theorem is degenerated, in the sense that the three pairs of sides are in prolongation, we obtain as a particular case the Gergonne's theorem.

## 4.1. Homological triangles inscribed in a circle.

Theorem 53 (Aubert – 1899)

Let $ABC$ and $A'B'C'$ two homological triangles inscribed in the same circle, $P$ their homological center and $I$ an arbitrary point on the circumscribed circle. The line $IA'$ intersects the side $BC$ in $U$, similarly are obtained the points $V,W$. The points $U,V,W$ are on the line that passes through the point $P$.

**Proof**

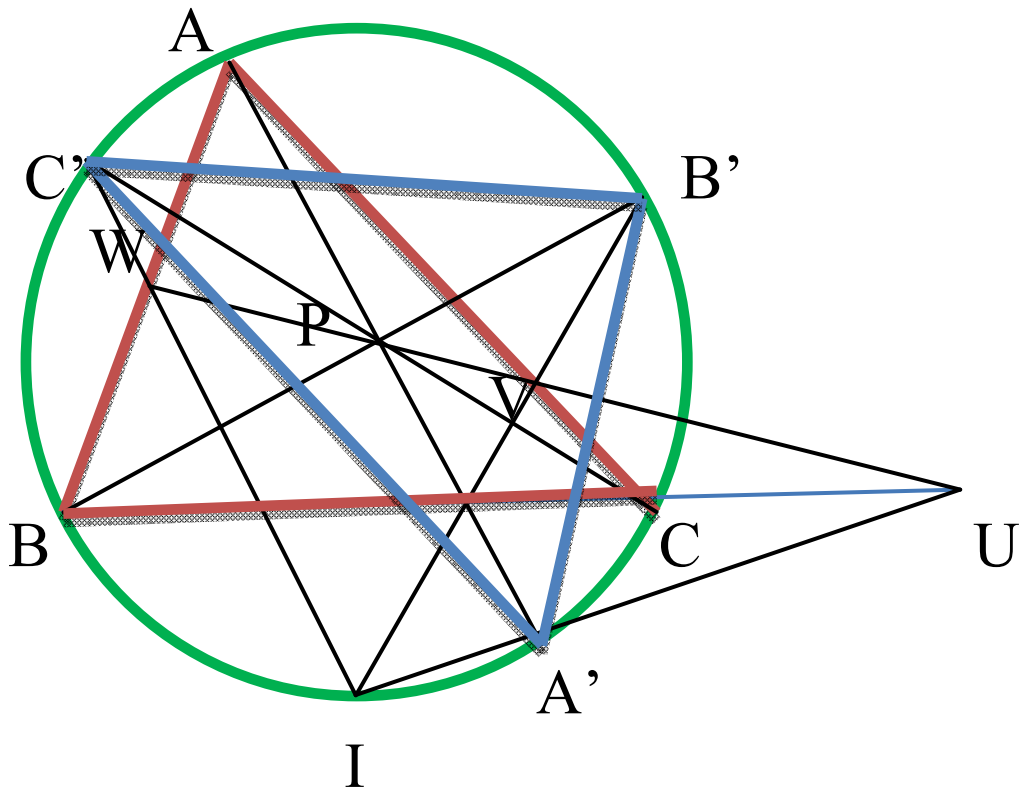

Fig 69

Consider the inscribed hexagon $IB'BACC'$ (see figure 69) and apply the Pascal's theorem.

It is obtain that the intersection points $V,P,W$ of the opposite sides $IB'$ with $AC$ of $B'B$ with $CC'$ and $BA$ with $C'I$ are collinear.

We consider the inscribed hexagon $IA'CCBB'$, and applying the Pascal's theorem, we find that the points $U,P,V$ are collinear.

From these two triplets of collinear points found, it results the collinearity of the points $U,V,W$ and $P$.

Transforming by duality the Aubert's theorem we obtain:

**Theorem 54** (the dual theorem of Aubert)

Let $ABC$ and $A'B'C'$ two homological triangles of axis $d_1$ circumscribed to a given circle and $t$ an arbitrary tangent to the circle that intersects the sides of the triangle $A'B'C'$ in $A"B"C"$. Then the triangles $ABC$ and $A"B"C"$ are homological, their homological center belonging to the line $d_1$.

**Theorem 55**

If $P,Q$ are isogonal conjugated points in the triangle $ABC$ and $P_1P_2P_3$ and $Q_1Q_2Q_3$ are their pedal triangles, we note with $X_1$ the intersection point between $P_2Q_3$ and $P_3Q_2$; similarly we define the points $X_2,X_3$. Then $X_1,X_2,X_3$ belong to the line $PQ$.

**Proof**

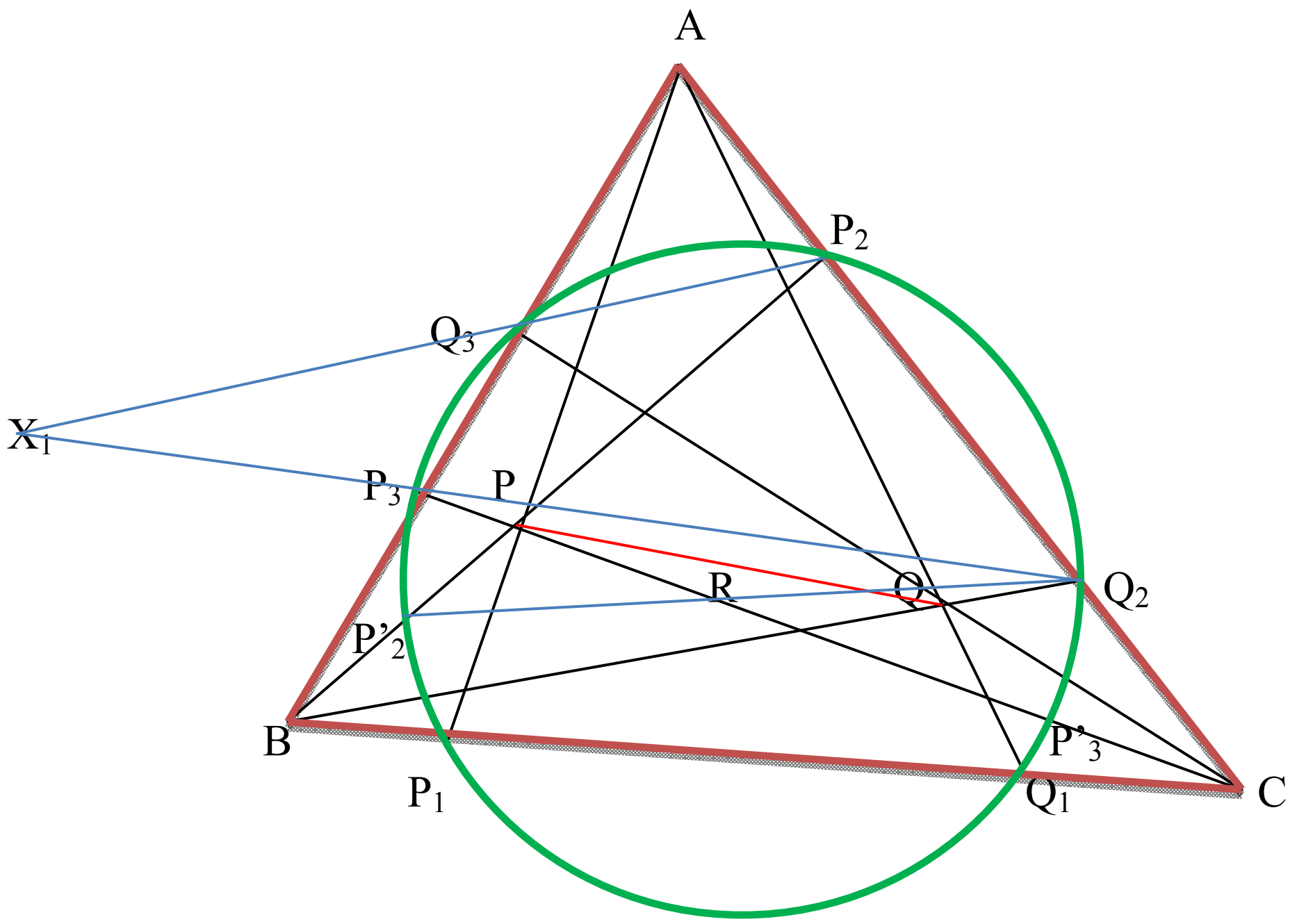

Fig70

It is known that the points $P_1Q_1Q_2P_2Q_3P_3$ are on a circle with the center $R$ which is the middle of $PQ$ (the circle of the 6 points) see figure 70.

We note with $P_2^{'}$ the intersection of the lines $Q_2R$ and $BP_2$ (the point $P_2^{'}$ belongs to the circle. Similarly $P_3^{'}$ is the intersection of the lines $Q_3R$ and $CP_3$ (the point $P_3^{'}$ is on the circle of the six points).

Applying the Pascal's theorem in the inscribed hexagon $P_2Q_3P_3^{'}P_3Q_2P_2^{'}$ it results that the points $X_1, R$ and $P$ are collinear. Similarly it can be shown that $X_2$ and $X_3$ belong to the line $PQ$.

**Theorem 56** (Alasia's theorem)

A circle intersects the sides $AB, BC, CA$ of a triangle $ABC$ in the points $A, D'$; $E, E'$ respectively $F, F'$. The lines $DE', EF', FD'$ determine a triangle $A'B'C'$ homological with triangle $ABC$.

**Proof**

We note

$$\{A'\} = DE' \cap EF'$$

$$\{B'\} = FD' \cap EF'$$

$$\{C'\} = FD' \cap DE'$$

$$\{B''\} = A'C' \cap AC$$

$$\{A''\} = B'C' \cap BC$$

$$\{C''\} = A'B' \cap AB$$

See figure 71

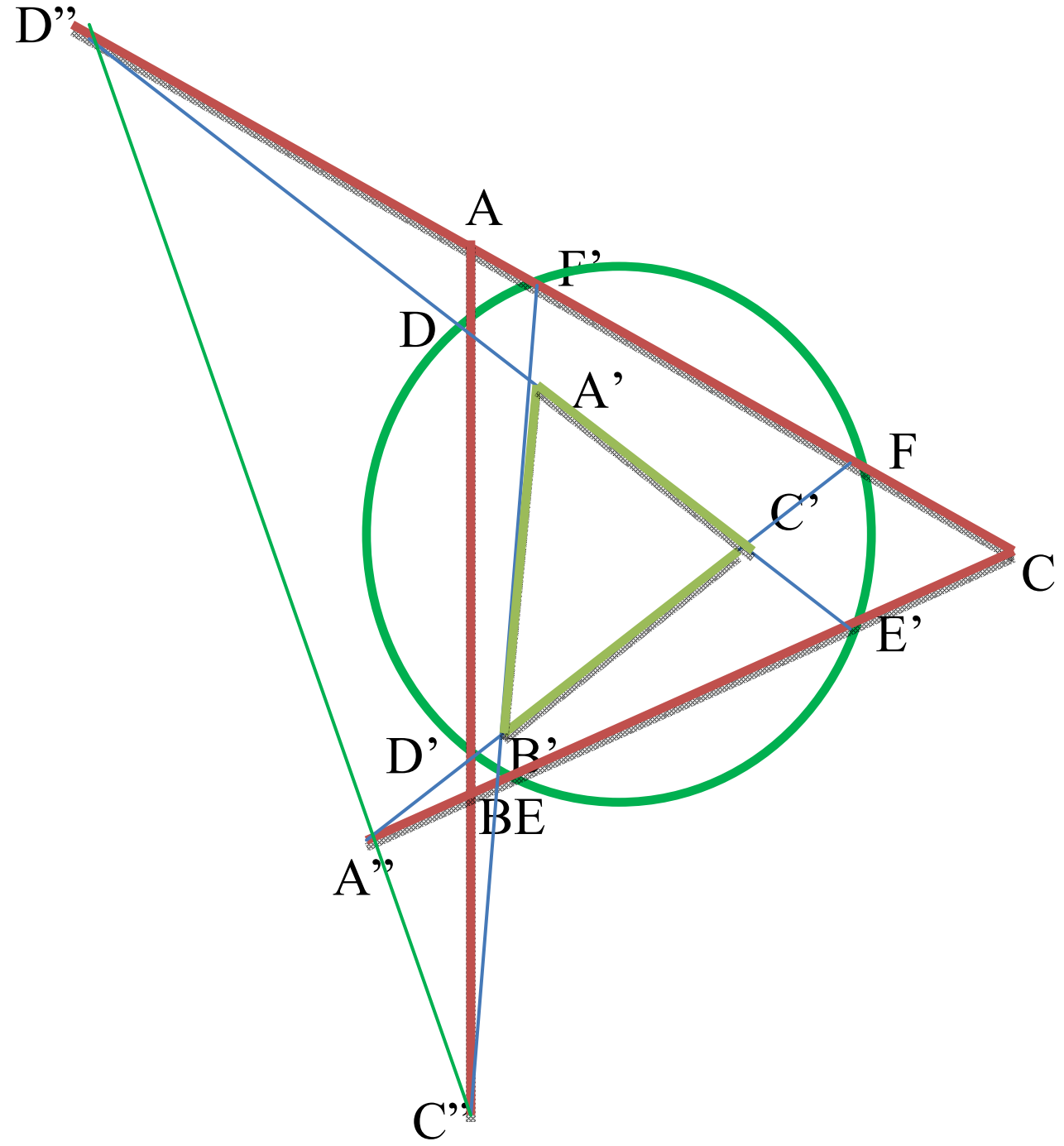

Fig.71

We apply Menelaus' theorem in the triangle $ABC$ for the transversals $A'', D', F$ ; $B'', D, E'$; $C'', E, F'$ obtaining

$$\frac{A''B}{A''C} \cdot \frac{D'A}{D'B} \cdot \frac{FC}{FA} = 1$$

$$\frac{B''C}{B''A} \cdot \frac{E'B}{E'C} \cdot \frac{DA}{DB} = 1$$

$$\frac{C''A}{C''B} \cdot \frac{F'C}{F'A} \cdot \frac{EB}{EC} = 1$$

From these relations and taking into account the Carnot's theorem it results

$$\frac{A''B}{A''C} \cdot \frac{B''C}{B''A} \cdot \frac{C''A}{C''B} = \frac{D'B}{D'A} \cdot \frac{FA}{FC} \cdot \frac{E'C}{E'B} \cdot \frac{DB}{DA} \cdot \frac{F'A}{F'C} \cdot \frac{EC}{EB} = 1$$

From the Menelaus' theorem, it results that $A'', B'', C''$ are collinear and from the reciprocal of the Desargues' theorem we obtain that the triangles $ABC$ and $A'B'C'$ are homological.

**Theorem 57** (the dual theorem of Alasia)

If $A_1B_1C_1D_1E_1F_1$ is a circumscribable hexagon and we note $X, Y, Z$ the intersection of the opposite sides $(B_1C_1, E_1F_1), (A_1B_1, D_1E_1)$ respectively $(C_1D_1, A_1F_1)$, then the triangles $A_1B_1C_1$ and $XYZ$ are homological.

**Proof**

See figure 72

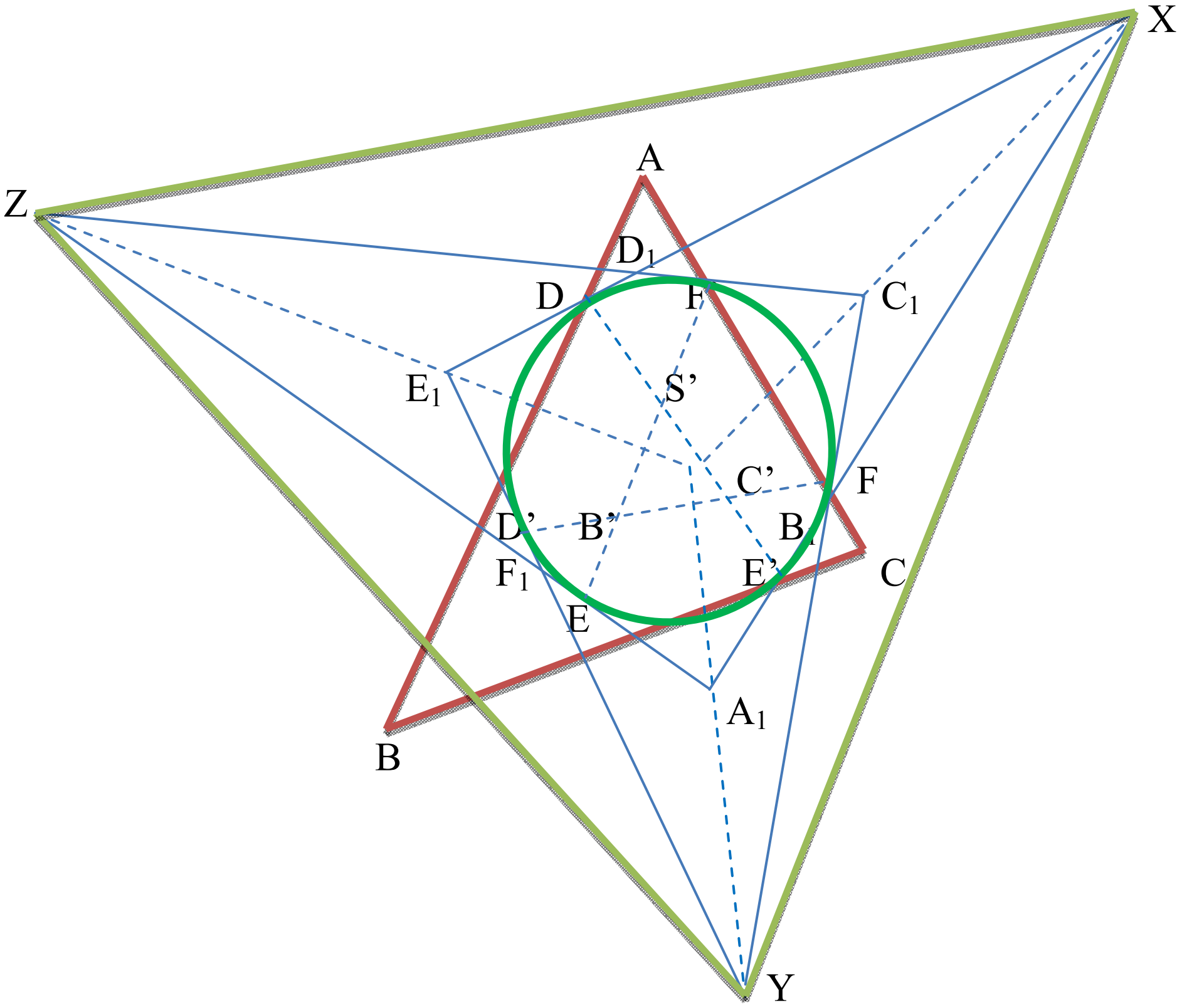

Fig. 72

We will transform by duality in rapport to the circle the Alasia's theorem, see figure 72.

To the points $D, D', E, E', F, F'$ will correspond to the tangents constructed in these points to the circle (their polar). To the line $EE'$ will correspond the point $A_1$ which is the intersection of the tangents on $E$ and $F'$, and to the line $D'F$ corresponds the intersection $X$ of the tangents constructed in $D', F$ (see figure 52), therefore to the intersection point $A''$ between $BC, D'F$ corresponds its polar, that is the line $XA_1$, similarly to the intersection point $B''$ between $DE'$ and $AC$ corresponds the line $C_1Y$, and to the point $A''$ corresponds line $E_1Z$. The points $A'', B'', C''$ are collinear (Alasia's theorem). It results that the their polar are concurrent, consequently the lines $A_1X, C_1Y, E_1Z$ are concurrent and the triangles $A_1B_1C_1$ and $XYZ$ are homological.

**Theorem 58**

Let $ABC$ and $A_1B_1C_1$ two homological triangles inscribed in the circle $\mathcal{C}(O, R)$ having the homology center $P$ and the axis $(d)$. If $A'B'C'$ and $A_1^{'}B_1^{'}C_1^{'}$ are their tangential triangles, then these are homological having the same center $P$ and axis $(d)$.

**Proof**

$$\{U\} = BC \cap B_1C_1$$

$$\{V\} = CA \cap A_1C_1$$

$\{W\} = AB \cap A_1B_1$

See figure 73

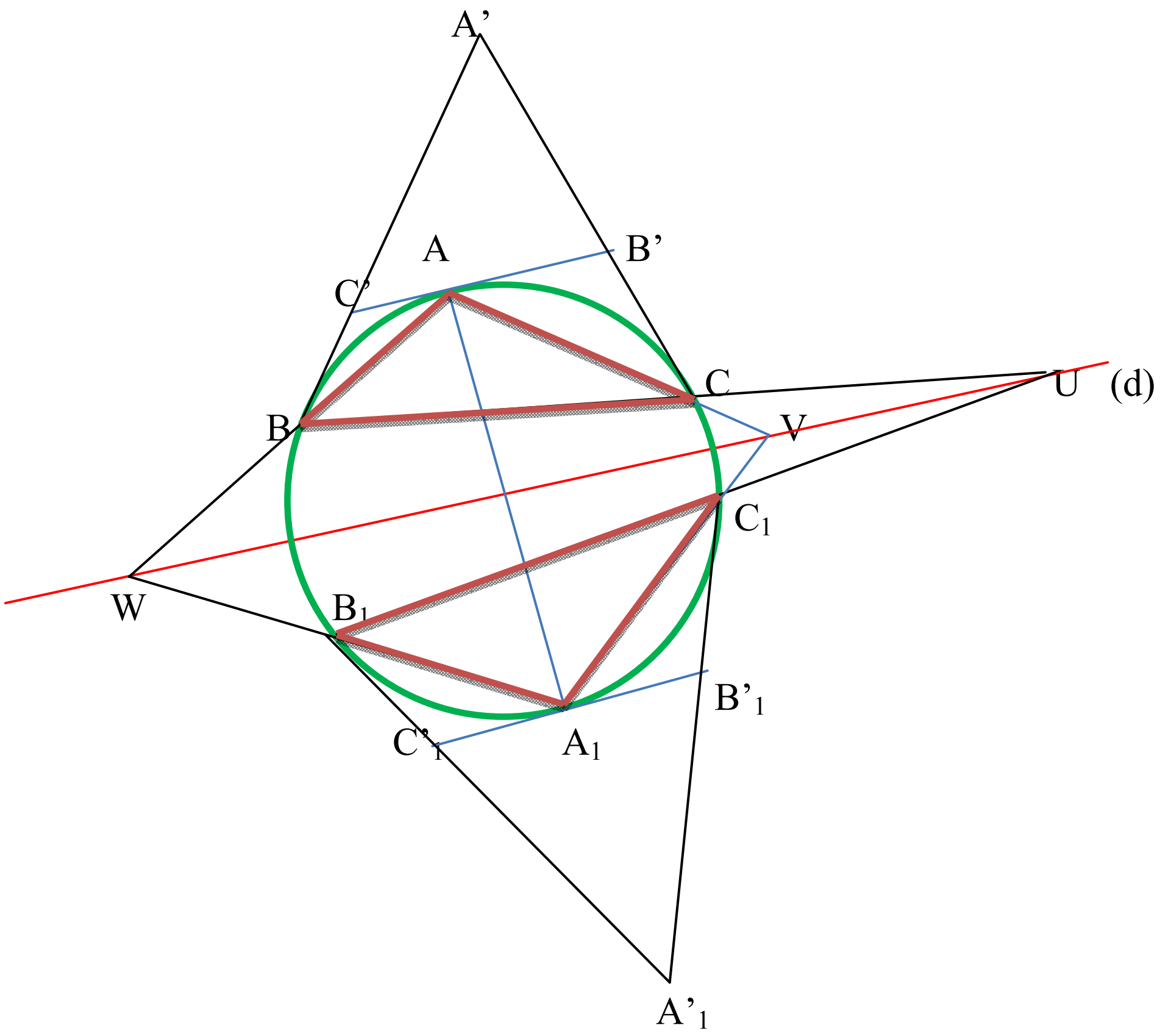

Fig. 73

The points $U, V, W$ belong to the homological axis $(d)$ .From the theorem (47) applied to the quadrilaterals $ABB_1A_1$, $ACC_1A$ and $BCC_1B_1$ it results that the polar of the point $P$ in rapport with the circle $(O)$ is the line $(d)$.

Because the polar of $A$ is $B'C'$ and the polar of $A_1$ is $B_1^{'}C_1^{'}$ it results that $\{U'\} = B'C' \cap B_1^{'}C_1^{'}$ is the pole of the line $AA_1$, but $AA_1$ passes through $P$, therefore the pole $U'$ of $AA_1$ belongs to the polar of $P$, that is to the line $(d)$. Similarly it results that $\{W'\} = A'B' \cap A_1^{'}B_1^{'}$ , therefore $W'$ belongs to $(d)$ and the triangles $A'B'C'$ and $A_1^{'}B_1^{'}C_1^{'}$ have as homological axis the line $(d)$.

The line $A'A_1^{'}$ is the polar of the point $U$ because $A'$ is the pole of $BC$ and $A_1^{'}$ is the pole of $B_1^{'}C_1^{'}$, therefore $A'A_1^{'}$ is the polar of a point on the polar of $P$, therefore $A'A_1^{'}$ passes through $P$. Similarly results that $B'B_1^{'}$ and $C'C_1^{'}$ pass through $P$.

**Remark 54**

If the triangles $ABC$ and $A_1B_1C_1$ are quasi-median then taking into consideration the precedent theorem and the proposition (25) we obtain that the triangles $ABC$ and $A_1B_1C_1$ and their tangential $A_1B_1C_1$ and $A_1^{'}B_1^{'}C_1^{'}$ form a quartet of triangles two b two homological having the same homological center and the same homological axis, which is the symmedian center, respectively the Lemoire's line of triangle $ABC$.

**Theorem 59** (Jerabeck)

If the lines which connect the vertexes of the triangle $ABC$ with two points $M', M''$ intersect the second time the triangle's circumscribed circle in the points $A', B', C'$ and $A'', B'', C''$, then the triangle determined by the lines $A'A'', B'B'', C'C''$ is homological with the triangle $ABC$.

**Proof**

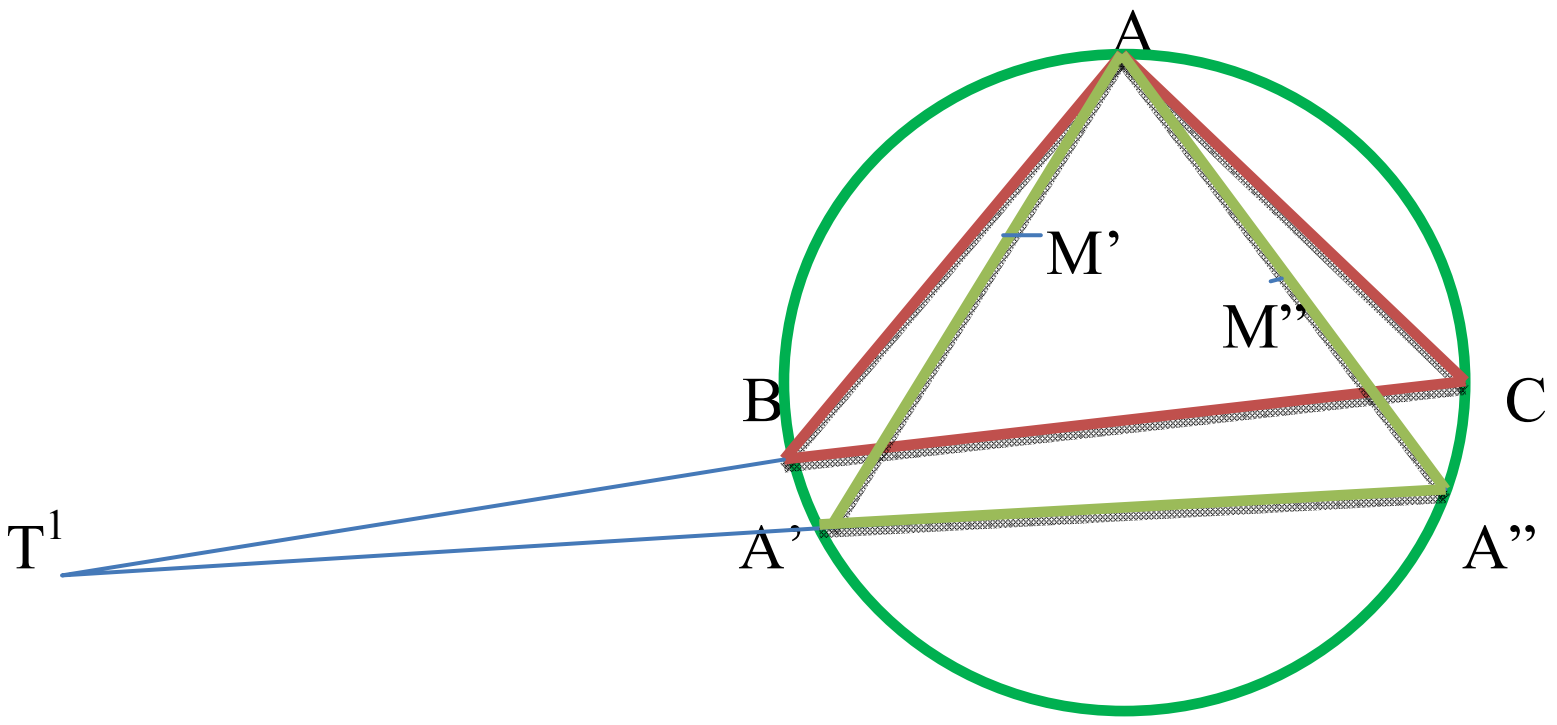

Fig. 74

We'll consider the $M', M''$ in the interior of the triangle $ABC$ and we note

$$m(\sphericalangle BAA') = \alpha,\ m(\sphericalangle CAA'') = \alpha'$$

$$m(\sphericalangle CBB') = \beta,\ m(\sphericalangle ABB'') = \beta'$$

$$m(\sphericalangle ACC') = \gamma,\ m(\sphericalangle BCC'') = \gamma'$$

$$\{T_1\} = A'A'' \cap BC$$

$$\{T_2\} = C'C'' \cap AC$$

$$\{T_3\} = B'B'' \cap AB$$

From the similarity of the triangles $T_1BA'$ and $T_1A''C$ we have

$$\frac{T_1B}{T_1A''} = \frac{T_1A'}{T_1C} = \frac{BA'}{A''C}$$

From the sinus' theorem in the triangles $BAA', CAA''$ we find

$$BA' = 2R\sin\alpha \text{ and } CA'' = 2R\sin\alpha'$$

consequently,

$$\frac{T_1B}{T_1A''} = \frac{T_1A'}{T_1C} = \frac{\sin\alpha}{\sin\alpha'} \qquad (1)$$

Also from the sinus's theorem applied in the triangles $T_1BA'$ and $T_1CA'$

$$\frac{T_1B}{\sin(A-\alpha')} = \frac{T_1A'}{\sin(A-\alpha)} \qquad (2)$$

$$\frac{T_1C}{\sin(A-\alpha)}=\frac{T_1A''}{\sin(A-\alpha')} \quad (3)$$

From the relations (1), (2) and (3) we obtain

$$\frac{T_1B}{T_1C}=\frac{\sin\alpha}{\sin\alpha'}\cdot\frac{\sin(A-\alpha')}{\sin(A-\alpha)} \quad (4)$$

Similarly we find

$$\frac{T_2B}{T_2A}=\frac{\sin\beta}{\sin\beta'}\cdot\frac{\sin(B-\beta')}{\sin(B-\beta')} \quad (5)$$

and

$$\frac{T_3A}{T_3B}=\frac{\sin\gamma}{\sin\gamma'}\cdot\frac{\sin(C-\gamma')}{\sin(C-\gamma')} \quad (6)$$

The relations (4), (5), (6) along with Ceva's theorem (the trigonometric variant) lead to the collinearity of the points $T_1, T_2, T_3$ and implicitly to the homology of the triangles $ABC$ and $A'''B'''C'''$, where we noted $\{A'''\}=B'B''\cap C'C''; \{B'''\}=A'A''\cap C'C''; \{C'''\}=B'B''\cap A'A''$.

**Remark 55**

If $M'=M''=M$ we'll obtain the following theorem:

1) The tangential triangle of the circumpedal triangle of the point $M\neq I$ from the interior of triangle $ABC$ and the triangle $ABC$ are homological.
2) If $M=I$ the triangle $ABC$ and the tangential triangle of the circumpedal triangle of $I$ (the center of the inscribed circle) are homothetic.
3) The triangles $ABC$ and $A'''B'''C'''$ are homothetic in the hypothesis that $M', M''$ are isogonal conjugate in the triangle $ABC$.

Bellow will formulate the dual theorem of the precedent theorem and of the Jerabeck's theorem.

**Theorem 60**

Let $ABC$ be a given triangle, $C_aC_bC_c$ its contact triangle and $T_1-T_2-T_3$ an external transversal of the inscribed circle $T_1\in BC, T_2\in CA, T_3\in AB$. If $A'$ is the second tangential point with the inscribed circle of the tangent constructed from $T_2$ ($A'\neq C_a$), $B''$ is the tangency point with the inscribed circle of the tangent constructed from $T_2$ ($B'\neq C_b$), and C' is the tangency point with the inscribed circle of the tangent constructed from $T_3$ ($C'\neq C_c$), then the triangles ABC and A'B'C' are homological.

**Theorem 61** (the dual theorem of Jerabeck's theorem)

Let $ABC$ an arbitrary given triangle and its contact triangle; we consider two transversals $T_1T_2T_3$ and $T_1'T_2'T_3'$ exterior to the inscribed circle ($T_1T_1'\in BC$, etc.) and we note $A', A''$ the tangent points with the inscribed circle of the tangents constructed from $T_1$ respectively $T_1'$ ($A', A''$ different of $C_a$), also we note with $A'''$ the intersection point of these tangents. Similarly are obtained the points $B', B'', B'''$ and $C', C'', C'''$. Then the triangle $ABC$ is homological with each of the triangles $A'B'C'$, $A''B''C''$ and $A'''B'''C'''$.

**Proof**

We'll consider the configuration from the Jerabeck's theorem. (see figure 75) and we'll transform it through reciprocal polar.

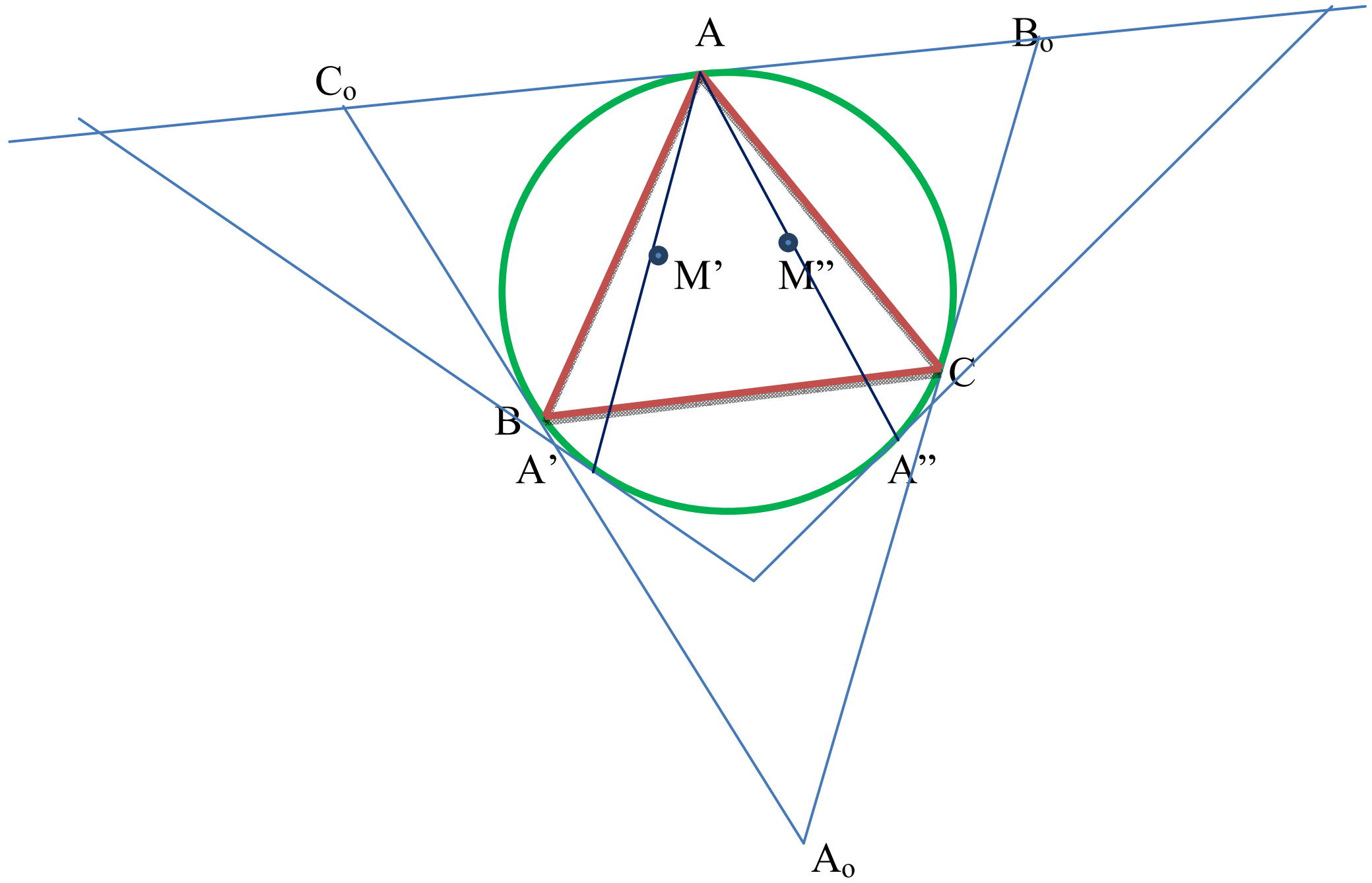

Fig. 75

Therefore if $ABC$ is the inscribed triangle in the circle $O$ and $AA, AA''$ are the Cevians from the hypothesis, we observe that to the point $A$ corresponds the tangent in $A$ to the circumscribe circle of the triangle $ABC$ and similarly to $B_2 - C$ we note the triangle formed by these tangents $A_oB_oC_o$. To point $A'$ corresponds the tangent in $A'$ constructed to the circle and in the same manner to the point $A''$ corresponds the tangent to the circle. We'll note $A_1$ the intersection point of the tangents.

Through the considered duality, to the line $BC$ corresponds the point $A_o$ (its pole), and to the line $A'A''$, its pole noted with $A_1$.

Because $BC$ and $A'A''$ intersect un a point, it results that that point is the pole of the line $A_oA_1$. Because the intersection points of the lines $BC$ and $A'A''$; $AC$ and $B'B''$; $AB$ and $C'C''$ are collinear, it will result that the lines $A_oA_1, B_oB_1, C_oC_1$ are concurrent. We note $A_o$ with $A$, $A_1$ with $A_3$, $A$ with $C_o$, etc. we obtain the dual theorem of Jerabeck's theorem.

**Theorem 62**

Let $ABC$, $A_1B_1C_1$ two homological triangles inscribed in a given circle. The tangents in $A_1, B_1, C_1$ to the circle intersect the lines $BC, CA, AB$ in three collinear points.

**Proof**

If we consider he triangles $ABC$ and $A_1B_1C_1$ homological with the center $O$, we note

$m(\sphericalangle BAA_1) = \alpha$, $m(\sphericalangle CBB_1) = \beta$ and $m(\sphericalangle ACC_1) = \gamma$

From the sinus's theorem in $ABA_1$ and we find

$$A_1B = 2R\sin\alpha,\ A_1C = 2R\sin(A-\alpha)$$

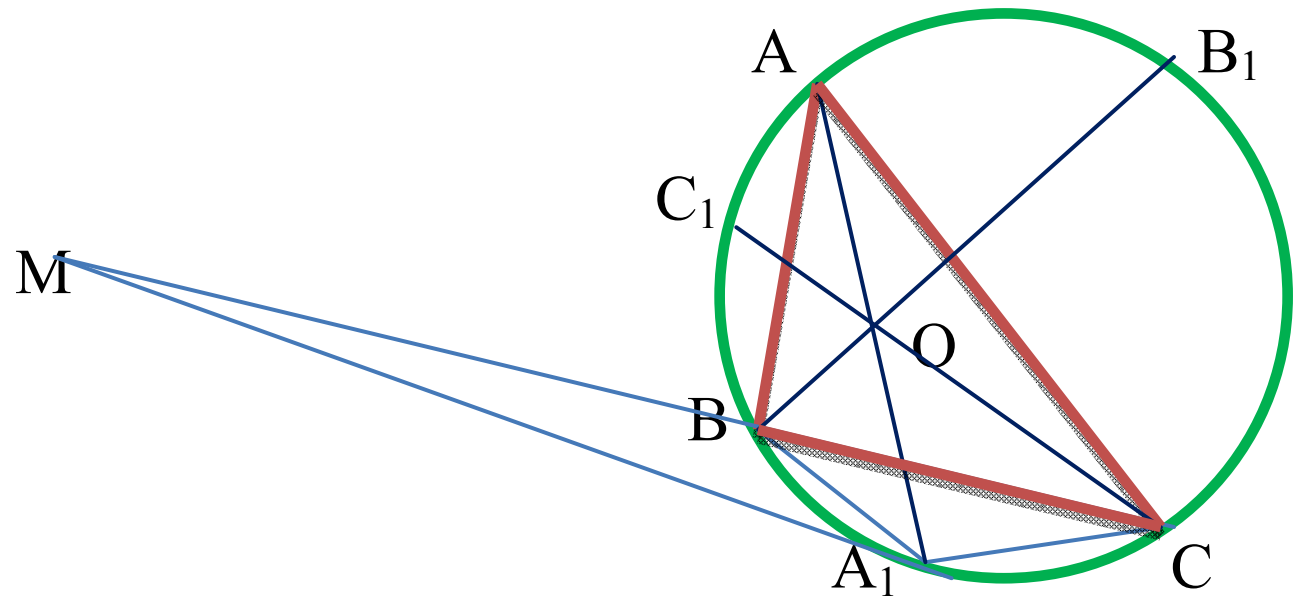

Fig. 76

Therefore

$$\frac{A_1B}{A_1C} = \left(\frac{\sin\alpha}{\sin(A-\alpha)}\right)^2$$

Similarly $\frac{NC}{NA} = \left(\frac{\sin\beta}{\sin(B-\beta)}\right)^2$, $\frac{PA}{PB} = \left(\frac{\sin\gamma}{\sin(C-\gamma)}\right)^2$.

Using the reciprocal of Menelaus' theorem immediately results the collinearity of the points $M, N, P$

**Observation 31**

a) Similarly, it result that the tangents constructed in $A, B, C$ to the circumscribed triangle intersect the sides $B_1C_1, A_1C_1, C_1B_1$ in collinear points.

b) The theorem can be formulated also as follows: The tangential; triangle of homological triangle with a given triangle (both inscribed in the same circle) is homological with the given triangle.

**Theorem** (I. Pătraşcu)

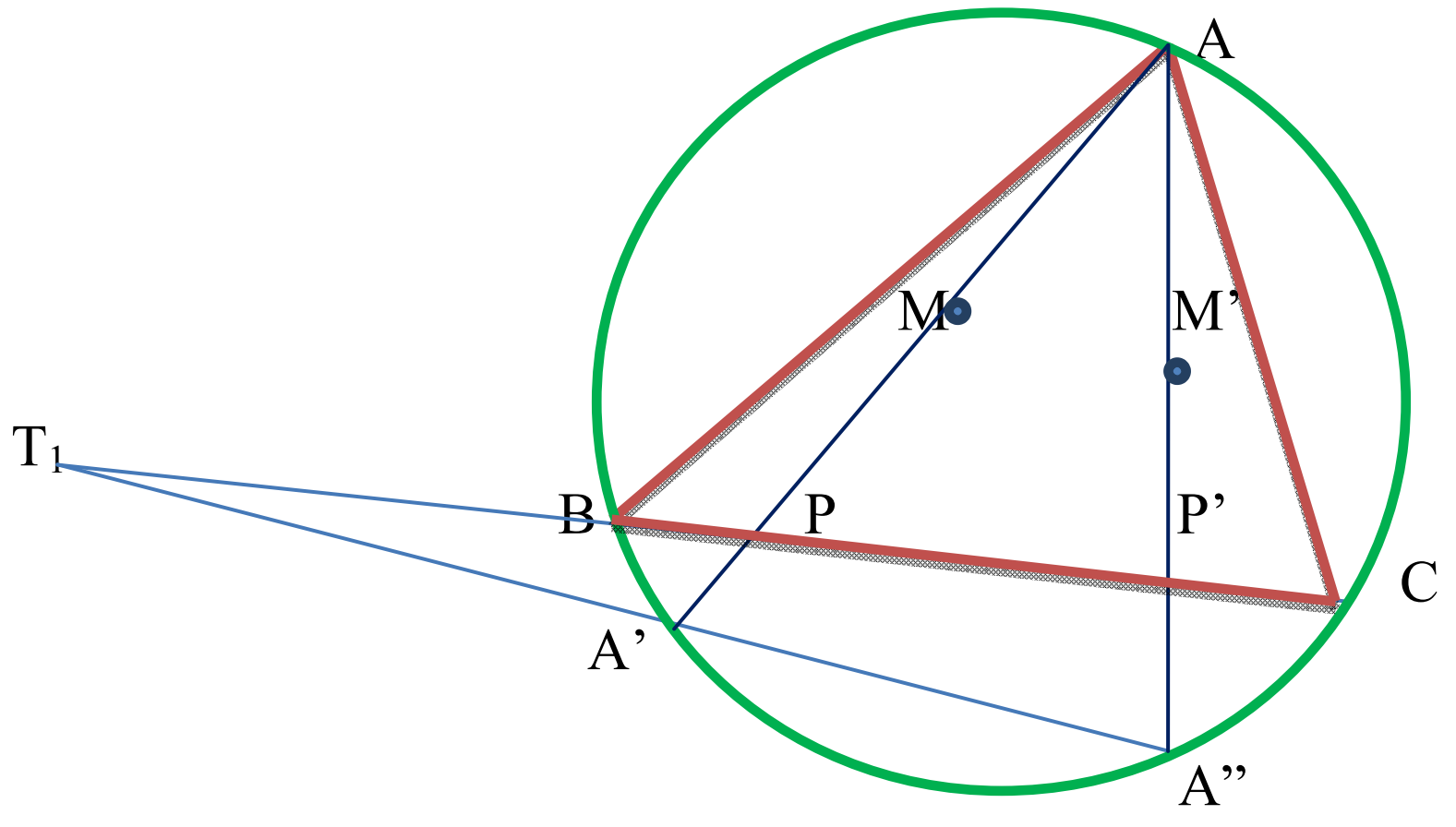

Fig. 77

Let $M, M'$ be two isotomic points conjugate inside of triangle $ABC$, and $A'B'C'$, $A"B"C"$ their circumpedal triangles. The triangle determined by the lines $A'A", B'B", C'C"$ is homological with the triangle $ABC$. The homology axis of these triangles is the isotomic transversal of the Lemoine's line to triangle $ABC$.

**Proof**

Let $\{T_1\} = BC \cap A'A"$ , $\{P\} = AA' \cap BC$ , $\{P'\} = AA" \cap BC$ . We note $\alpha = m\left(\widehat{BAA'}\right)$, $\alpha' = m\left(\widehat{CAA"}\right)$, following the same process as in Jerabeck's theorem we obtain:

$$\frac{T_1B}{T_1C} = \frac{\sin\alpha}{\sin\alpha'} \cdot \frac{\sin(A-\alpha')}{\sin(A-\alpha)}$$

From Jerabeck's theorem results that the triangle $ABC$ and the triangle formed by the lines $A'A", B'B", C'C"$ are homological.

From the sinus' theorem we have

$$\frac{\sin\alpha}{BP} = \frac{\sin\sphericalangle APB}{AB}$$

Also

$$\frac{\sin(A-\alpha)}{PC} = \frac{\sin\sphericalangle APC}{AC}$$

Because $\sin\sphericalangle APB = \sin\sphericalangle APC$, from the precedent relations we retain that

$$\frac{\sin\alpha}{\sin(A-\alpha)} = \frac{PC}{BP} \cdot \frac{AC}{AB}$$

Similarly

$$\frac{\sin\alpha'}{\sin(A-\alpha')} = \frac{P'B}{P'C} \cdot \frac{AB}{AC}.$$

The Cevians $AP, AP'$ being isometric we have $BP = P'C$ and $BP' = CP$.

We find that $\frac{T_1B}{T_1C} = \left(\frac{AC}{AB}\right)^2$, it is known that the exterior symmedian of the vertex $A$ in the triangle $ABC$ is tangent in $A$ to the circumscribed circle and if $T'$ is its intersection with $BC$ then $\frac{T_1'C}{T_1'B} = \left(\frac{AC}{AB}\right)^2$. This relation and the precedent show that $T_1$ and $T_1'$ are isotomic points. Similarly, if $T_2$ is the intersection of the line $B'B"$ with $AC$, we can show that $T_2$ is the isotomic of the external symmedian leg from the vertex $B$, and if $T_3$ is the intersection of the line $C'C"$with $AC$, $T_3$ is the isotomic of the external symmedian leg from the vertex $C$ of the triangle $ABC$. The homology axis of the triangle from the hypothesis is $T_1T_2T_3$ and it is the isotomic transversal of the line determined by the legs of the external symmedian of triangle $ABC$, that is the Lemoine's line of triangle $ABC$.

# Chapter 5

# Proposed problems. Open problems

## 5.1. Proposed problems

**1.** If $ABCD$ parallelogram, $A_1 \in (AB), B_1 \in (BC), C_1 \in (CD), D_1 \in (DA)$ such that the lines $A_1D_1, BD, B_1C_1$ are concurrent, then

(i) The lines $AC, A_1C_1, B_1D_1$ are concurrent;

(ii) The lines $A_1B_1, C_1D_1, AC$ are concurrent

(Florentin Smarandache, Ion Pătraşcu)

**2**. Let $ABCD$ a quadrilateral convex such that

$$\{E\} = AB \cap CD$$
$$\{F\} = BC \cap AD$$
$$\{P\} = BD \cap EF$$
$$\{R\} = AC \cap EF$$
$$\{O\} = AC \cap BD$$
$$(AB), (BF), (CA)$$

We note $G, H, I, J, K, L, P, O, R, M, N, Q, U, V, T$ respectively the middle point of the segments $(AB), (BF), (CA), (AD), (AE), (DE), (CE), (BE), (BC), (CF), (DF), (DC)$. Prove that:

(i) Triangle $POR$ is homological with each of the triangles $GHI, JKL, MNQ, UVT$;

(ii) The triangles $GHI, JKL$ are homological;

(iii) The triangles $MNQ, UVT$ are homological;

(iv) The homology centers of the triangles $GHI, JKL, POR$ are collinear;

(v) The homology centers of the triangles $MNQ, UVT, POR$ are collinear

(Florentin Smarandache, Ion Pătraşcu)

**3.** Let $ABC$ a triangle and $A_1B_1C_1$, $A_2B_2C_2$ isotomic triangles inscribed in $ABC$. Prove that if the triangles $ABC$ and $A_2B_2C_2$ are homological then:

(i) The triangle $ABC$ and the triangle $A_2B_2C_2$ are homological and their homology center is the isotomic conjugate of the homology center of the triangles $ABC$ and $A_1B_1C_1$

(ii) The triangle $ABC$ and the medial triangle of the triangle $A_1B_1C_1$ are homological.

**4**. Let $A_1B_1C_1$ and $A_2B_2C_2$ equilateral triangles having the same center $O$. We note:

$$\{A_3\} = B_1B_2 \cap C_1C_2$$
$$\{B_3\} = A_1A_2 \cap C_1C_2$$
$$\{C_3\} = A_1A_2 \cap B_1B_2$$

Prove that:

(i) $(A_2B_2) \equiv (B_1C_2) \equiv (C_1A_2)$;

(ii) $(A_2B_2) \equiv (B_1B_2) \equiv (C_1C_2)$;

(iii) $(A_2B_2) \equiv (B_1A_2) \equiv (C_1B_2)$;

(iv) The triangle $A_3B_3C_3$ is equilateral and has its center in the point $O$;

(v) The triangle $A_2B_2C_2$ and the triangle $A_3B_3C_3$ are tri-homological.

(Ion Pătraşcu)

**5**. If a circle passes through the vertexes $B, C$ of the triangle $ABC$ and intersect the second time $AB$ in $E$ and $AC$ in $D$, and we note $F$ the intersection of the tangent in $D$ to the circle with $BC$ and with $G$ the intersection of the tangent to the circle constructed in $C$ with the line $DE$ then the points $A, F, G$ are collinear.

**6**. Prove that in circumscribed octagon the four cords determined by the contact points with the circle of the opposite sides are concurrent

(Ion Pătraşcu, Florentin Smarandache)

**7.** Let two external circles in a plane. It is known the center of a circle, construct the center of the other circle only with the help of a unmarked ruler.

**8**. Let $A_1B_1C_1$ an inscribed triangle in the triangle $ABC$ such that the triangles $ABC$ and $A_1B_1C_1$ are homological. Prove that if $\overrightarrow{AA_1} + \overrightarrow{BB_1} + \overrightarrow{CC_1} = \vec{0}$ then the homology center of the triangles $ABC$ and $A_1B_1C_1$ is the weight center of the triangle $ABC$.

**9.** Let $A'B'C'$ the pedal triangle of a point in rapport with the triangle $ABC$. A transversal intersects the sides $BC, CA, AB$ in the points $U, V, W$. The lines $AU, BV, CW$ intersect $B'C', C'A', A'B'$ respectively in the points $U', V', W'$. Prove that $U', V', W'$ are collinear.

**10**. Let $A_1B_1C_1$ and $A_2B_2C_2$ the pedal triangles of the points $M_1, M_2$ in rapport with the triangle $ABC$. We note $A_3, B_3, C_3$ the intersection points of the lines $B_1C_1$ and $B_2C_2$; $C_1A_1$ and $C_2A_2$; $A_1B_1$ and $A_2B_2$. We note $A_4, B_4, C_4$ the intersection points of the lines $B_1C_2$ and $B_2C_1$; $C_2A_1$ and $C_1A_2$; $A_1B_2$ and $A_2B_1$. Prove that:

(i) The sides of the triangle $A_3B_3C_3$ pass through the vertexes of the triangle $ABC$;

(ii) The points $A_4, B_4, C_4$ belong to the line $M_1M_2$;

(iii) The sides of the triangle $A_3B_3C_3$ are the polar of the opposite vertexes in rapport to the sides of the triangle $ABC$ taken two by two and pass through the points $A_4, B_4, C_4$;

(iv) The lines $AA_3, BB_3, CC_3$ are concurrent ;

(v) The triangle $A_3B_3C_3$ is homological with the triangles $A_1B_1C_1$ and $A_2B_2C_2$, their homological centers being on the homological axis of the triangles $ABC$ and $A_3B_3C_3$.

` (G.M. 1903-1904)

**11**. Let $ABCDEF$ a complete quadrilateral and $M$ point in side of the triangle $BCE$. We note $PQR$ the pedal triangle of $M$ in rapport with $BCE$ ($P \in (BC)$, $Q \in (BE)$). We also note $PR \cap AC = \{U\}$, $PQ \cap BD = \{V\}$, $RQ \cap EF = \{W\}$. Prove that the points $U, V, W$ re collinear.

(Ion Pătraşcu)

**12**. Prove that the tangential triangle of the triangle $ABC$ and the circumpedal triangle of the weight center of the triangle $ABC$ are homological.
Note: The homology center is called the Exeter point of the triangle $ABC$.

**13**. In the triangle $ABC$ let $U-V-W$ a transversal $U \in BC$, $V \in AC$ .
The points $U' \in AU$, $V' \in BV$, $W' \in CW$ are collinear and $A' \in BC$ such that the points $U', C', B'$ are collinear, where $\{B'\} = A'W' \cap AC$ and $\{C'\} = AB \cap A'V'$.
Prove that the triangles $ABC$ and $A'B'C'$ are homological.

(Ion Pătraşcu)

**14**. Let $ABC$ a given random triangle, $I$ the center of the inscribed circle and $C_aC_bC_c$ its contact triangle. The perpendiculars constructed from $I$ on $IA, BI, CI$ intersect the sides $BC, CA, AB$ respectively in the points $A', B', C'$. Prove that the triangle formed by the lines $AA'', BB', CC'$ homological with the triangle $C_aC_bC_c$.

**15**. If $A'B'C'$, $A''B''C''$ are inscribed triangles in the triangle $ABC$ and homological with it and if we note

$$\{A'''\} = C'B' \cap C''B'', \{B'''\} = A'C' \cap A''C'', \{C'''\} A'B' \cap A''B'',$$

then the triangle $A'''B'''C'''$ is homological with each of the triangles
$ABC$, $A'B'C'$, $A''B''C''$.

**16**. The complementary (medial) of the first Brocard triangle associated to the triangle $ABC$ and the triangle $ABC$ are homological

(Stall)

**17**. The tri-linear poles of the Longchamps's line of a triangle coincide with one of the homology center of this triangle and of the first triangle of Brocard .

(Longchamps)

**18**. If $BCA'$, $CAB'$, $ABC'$ are similar isosceles triangles constructed on the sides of the triangle $ABC$ in its interior or exterior.. Then the homological axis of the triangles $ABC$ and $A'B'C'$ is perpendicular on the line that connects their homological center with the center of the circle circumscribed to triangle $ABC$ . The perpendiculars constructed from the vertexes $A, B, C$ on the sides $B'C', C'A', A'B'$ are concurrent in a point of the same line.

**19**. Let $ABC$ a triangle and $A'B'C'$ the pedal triangle of the center of the circumscribed circle to triangle $ABC$ . We'll note $A'', B'', C''$ middle points of the segments $AA', BB', CC'$. The homology axis of the triangles $ABC$ and $A''B''C''$ is the tri-linear polar of the orthocenter of the triangle $ABC$.

**20**. Let $A_1B_1C_1, A_2B_2C_2$ two triangles circumscribed to triangle $ABC$ and homological with it, their homological centers being $M_1, M_2$. The lines $A_1M_2, B_1M_2, C_1M_2$ intersect $BC, CA, AB$ in $A', B', C'$. The lines $A_2M_1, B_2M_1, C_2M_1$ intersect $BC, CA, AB$ in the same points $A', B', C'$. The lines $AA', BB', CC'$ are concurrent.

**21.** In a triangle the lines determined by the feet of the height constructed from $B, C$ , the lines determined by the feet of the bisectors of the angles $B, C$ and the line of the contact points of the inscribed circle with the sides $AC, AB$ are concurrent in a point $U$ . Similarly, are defined the points $V, W$ . Prove that the lines $AU, BV, CW$ are concurrent.

**22**. Let $ABC$ a triangle and $A_1B_1C_1$ the circumpedal triangle of the circumscribed circle $I$ and $C_aC_bC_c$ the contact triangle of the triangle $ABC$.
Prove that

(i) $A_1C_a, B_1C_b, C_1C_cB$ are concurrent in the isogonal of the Nagel's point $N$ of the triangle $ABC$ ;

(ii) The isogonal of the Nagel's point, the center of the inscribed circle and the center of the circumscribed circle are collinear.

(Droz)

**23**. The perpendicular bisectors of an arbitrary given triangle intersects the opposite sides of the triangle in three collinear points.

(De Longchamps)

**24**. The homology axis of the triangle $ABC$ and of its orthic triangle is the radical axis of the circumscribed arcs and of the Euler's circle.

**25**. $A'B'C'$ are the projections of the center of the inscribed circle in the triangle $ABC$ on the perpendicular bisectors of the triangle, then the lines $AA', BB', CC'$ are concurrent in a point which is the isotonic conjugate of the Gergonne's point of the triangle.

(De Longchamps)

**26**. Prove that in an arbitrary triangle, the Lemoine's point, the triangle orthocenter and Lemoine's point of the orthic triangle are collinear

(Vigarie)

**27**. Prove that a triangle is isosceles if and only if the intersection of a median with a symmedian is a Brocard's point of the triangle.

(Ion Pătraşcu)

**28**. In an arbitrary triangle $ABC$ let $A'B'C'$ the circumpedal triangle of the center of the inscribed circle $I$. We'll note $A_1, B_1, C_1$ the intersections of the following pairs of lines $(BC, B'C'), (AC, A'C'), (AB, A'B')$. If $O$ is the center of the circumscribed circle of the triangle $ABC$, prove that the lines $OI, A_1B_1$ are perpendicular

(Ion Pătraşcu)

**29**. In the random triangle $ABC$ let $C_aC_bC_c$ its contact triangle and $A'B'C'$ the. pedal triangle of the center of the inscribed circle $I$.
We'll note $\{U\} = B'C' \cap C_bC_a$, $\{V\} = A'C' \cap C_aC_c$, $\{W\} = A'B' \cap C_aC_b$.
Prove that the perpendiculars constructed from $A, B, C$ respectively on $IU, IV, IW$ intersect the lines $C_bC_c, C_aC_c, C_aC_b$ in three collinear points.

(Ion Pătraşcu)

**30**. Let $ABCD$ a trapeze and $M, N$ the middle points of the bases $AB$ and $CD$, and $E \in (AD)$ different of the middle point of $(AD)$. The parallel through $E$ to the base intersects $(BC)$ in $F$.
Prove that the triangles $BMF, DNE$ homological.

**31**. Consider the triangle $ABC$ and the transversal $A_1, B_1, C_1$ ( $A_1 \in BC$ , $B_1 \in CA$ , $C_1 \in AB$). The lines $BB_1, CC_1$ intersect in $A_2$ ; the lines $CC_1, AA_1$ intersect in $B_2$ and the lines $AA_1, BB_1$ intersect in $CC_2$. Prove that the lines $A_1A_2, B_1B_2, C_1C_2$ are concurrent.

(Gh. Ţiţeica)

**32**. Let $ABCD$ an inscribed quadrilateral to a circle and $A', B', C', D'$ the tangency points of the circle with the sides $AB, BC, CD, DA$.

We'll notes $\{A_1\} = A'B' \cap C'D', \{B_1\} = A'D' \cap B'C'$.
Prove that
(i) The lines $AA_1, BD', DA'$ are concurrent;

(ii) The lines $BB_1, AB', CA'$ are concurrent.

(Ion Pătraşcu)

**33**. Let $ABC$ an arbitrary triangle , we note $D, E, F$ the contact points of the inscribed circle with the sides $BC, CA, AB$ and with $M, N, P$ the middle of the arches $BC, CA, AB$ of the circumscribed circle.
Prove that:

(i) The triangles $MNP$ and $DEF$ are homothetic, the homothety center being the point $L$;

(ii) We note $A_1, B_1, C_1$ the intersections of the segments $(LA), (LB), (LC)$ with the inscribed circle in the triangle $ABC$ and $A_2, B_2, C_2$ the intersection points of the $B_1C_1$ with $LM$ ; $A_1C_1$ with $LN$ and $A_1B_1$ with $LP$. Prove that the triangles $A_1B_1C_1$, $A_2B_2C_2$ are homothetic.

(Ion Pătraşcu)

**34**. Let $ABC$ an arbitrary triangle and $A'B'C'$ it contact triangle. We'll note $A'', B'', C''$ the diametric opposite points of $A', B', C'$ in the inscribed circle in the triangle $ABC$.
Prove that the triangles $ABC$, $A''B''C''$ are homological.

**35**. Let $ABCDEF$ a complete quadrilateral and $M \in (AC), N \in (BD), P \in (EF)$ such that

$$\frac{AM}{MC} = \frac{BN}{ND} = \frac{EP}{PF} = k.$$

Determine the value of $k$ the points $M, N, P$ are collinear.

(Ion Pătraşcu)

**36**. If $I$ is the center of the inscribed circle in the triangle $ABC$ and $D, E, F$ are the centers of the inscribed circles in the triangles $BIC, CIA, AIB$ respectively, then the lines $AD, BE, CF$ are concurrent. (The first point of Velliers)

**37**. Let $I_aI_bI_c$ the centers of the ex-inscribed circles corresponding to triangle $ABC$ and $I_1I_2I_3$ the centers of the inscribed circles in the right triangles $BI_aC$, $CI_bA$, $AI_cB$ respectively. Prove that the lines $AI_1$, $BI_2$, $CI_3$ are concurrent. (The second point of Velliers).

**38**. Let $I$ the center of the inscribed circle in the triangle $ABC$ and $I_1, I_2, I_3$ the centers of the ex-inscribed circles of triangles $BIC, CIA, AIB$ (tangents respectively to the sides $(BC), (CA), (AB)$).
Prove that the lines $AI_1$, $BI_2$, $CI_3$ are concurrent.

(Ion Pătraşcu)

**39**. Let $I_a, I_b, I_c$ the centers of the ex-inscribed circles to the triangle $ABC$ and $I_1, I_2, I_3$ the centers of the ex-inscribed circles of triangles $BI_aC$, $CI_bA$, $AI_cB$ respectively (tangent respectively to the sides $(BC), (CA), (AB)$).
Prove that the lines $AI_1$, $BI_2$, $CI_3$ are concurrent.

(Ion Pătraşcu)

**40**. Let $ABC$ a triangle inscribed in the circle $\mathcal{C}(O, R)$, $P$ a point in the interior of triangle and $A_1B_1C_1$ the circumpedal triangle of $P$. Prove that triangle $ABC$ and the tangential triangle of the triangle $A_1B_1C_1$ are homological.

**41**. Let $ABC$ a random triangle. $I$ the center of its inscribed circle and $I_aI_bI_c$ its anti-supplementary triangle. We'll note $O_1$ the center of the circle circumscribed to triangle $I_aI_bI_c$ and $M, N, P$ the middle of the small arches $\overset{\frown}{BC}, \overset{\frown}{CA}, \overset{\frown}{AB}$ from the circumscribed circle to triangle $ABC$. The perpendiculars from $I_a, I_b, I_c$ constructed respectively on $O_1M_1, O_1N_1, O_1P_1$ determine a triangle $A_1B_1C_1$.
Prove that the triangles $A_1B_1C_1$, $I_aI_bI_c$ are homological, the homology axis being the tri-linear polar of $I$ in rapport with the triangle $ABC$.

**42**. Let $\Omega$ the Brocard's point of triangle $ABC$ and $A_1B_1C_1$ the circumpedal triangle of $ABC$. We note $A_2B_2C_2$ the triangle, which has as vertexes the diametric vertexes $A_1, B_1, C_1$. We'll note $\{A_2\} = BC \cap A_1B_1$, $\{B_2\} = AC \cap B_1C_1$, $\{C_2\} = AB \cap C_1A_1$, $\{A_3\} = BA_1 \cap CB_1$, $\{B_3\} = AC_1 \cap CB_1$, $\{C_3\} = AC_1 \cap BA_1$.
Prove that the triangles $A_2B_2C_2$, $A_3B_3C_3$ are homological and their homological axis is the perpendicular on the line $OO_1$, where $O$ is the center of the circumscribed circle to triangle $A_2B_2C_2$.

(Ion Pătraşcu)

**43**. Let $I$ and $O$ the centers of the inscribed and circumscribed circle to the triangle $ABC$, and $J$ the symmetric point of $I$ in rapport with $O$. The perpendiculars in $A, B, C$ respectively on the lines $AJ, BJ, CJ$ form an homological triangle with the contact triangle of $ABC$. The homology axis of the two triangles is the radical axis of the inscribed and circumscribed circles and the homological center is $I$.

(C. Ionescu-Bujor)

**44**. Are given the circles $(C_a),(C_b),(C_c)$ that intersect two by two. Let $A_1,A_2$ the common points of the pair of circles $(C_b)$, $(C_a)$; $B_1,B_2$ the common points of the pair of circles $(C_c),(C_a)$ and $C_1,C_2$ the common points of the pair of circles $(C_a),(C_b)$. Prove that the triangles $A_1B_1C_1$, $A_2B_2C_2$ are homological, having the homology center in the radical center of the given circles and as homological axis the radical axis of the circles $A_1B_1C_1$, $A_2B_2C_2$.

**45**. Let $ABC$ a triangle, $H$ its orthocenter and $L$ the symmetric of $H$ in rapport with the center of the circumscribed to triangle $ABC$. He parallels constructed through $A,B,C$ to the sides of triangle $ABC$ form an homological triangle with the pedal triangle $A_1B_1C_1$ of the point $L$. The homology axis of the two triangles is the radical axis of the circles $(ABC),(A_1B_1C_1)$, and the homology center is the point $L$.

(C. Ionescu-Bujor)

**46**. Let $ABC$ and $A'B'C'$ two conjugate triangles inscribed in the same circle $(O)$. We'll note $A_1B_1C_1$ the triangle formed by the intersections of the lines $(BC,B'C')$, $(AC,A'C')$, $(AB,A'B')$ and $A_2B_2C_2$ the triangle formed by the lines $AA'',BB',CC'$. Prove that the triangles $A_1B_1C_1$, $A_2B_2C_2$ are homological, having as homological axis the perpendicular on the line determined by the center of the circumscribed circle to the triangle $A_1B_1C_1$ and $O$.

**47**. Let $I$ the center of the inscribed circle in the triangle $ABC$; the lines which pass through $I$ intersect the circumscribed circles of triangles $BIC,CIA,AIB$ respectively in the points $A_1,B_1,C_1$. Prove that:

(i) The projections $A_1',B_1',C_1'$ of the points $A_1,B_1,C_1$ on the sides $BC,CA,AB$ are collinear.

(ii) The tangents in $A_1,B_1,C_1$ respectively to the circumscribed circles to triangles $BIC,CIA,AIB$ form an homological triangle with the triangle $ABC$.

. (Gh. Ţiţeica)

**48**. Let $O$ the center of the circumscribed circle to a random triangle $ABC$, we'll note $P,Q$ the intersection of the radius $OB$ with the height and the median from $A$ and $R,S$ the intersections of the radius $OA$ with the height and median from $B$ of triangle $ABC$.

Prove that the lines $PR, QS$ and $AB$ are concurrent.

(Ion Pătraşcu)

**49**. Let $ABC$ an inscribed triangle in a circle of center $O$. The tangent in $B, C$ to the circumscribed circle to the triangle intersect in a point $P$. We'll note $Q$ the intersection point of the median $AM$ with the circumscribed circle and $R$ the intersection point of the polar to $BC$ constructed through $A$ with the circumscribed circle.
Prove that the points $P, Q, R$ are collinear.

**50**. Two lines constructed through the center $I$ of the inscribed circle of the triangle $ABC$ intersects the circles $(BIC), (CIA), (AIB)$ respectively in the points $A_1, A_2$; $B_1, B_2$; $C_1, C_2$.
Prove that the lines $A_1A_2, B_1B_2, C_1C_2$ form a homological triangle with $ABC$.

**51**. In the scalene triangle $ABC$ let consider $AA_1, BB_1, CC_1$ the concurrent Cevians in $P$ ($A_1 \in (BC), B_1 \in (CA), C_1 \in (AB)$), and $AA_2, BB_2, CC_2$ the isogonal Cevians to the anterior Cevians, and $O$ their intersection point. We'll note $P_1, P_2, P_3$ the orthogonal projections of the point $P$ on the $BC, AC, AB$; $R_1, R_2, R_3$ the middle point of the segments $(AP), (BP), (CP)$.
Prove that if the points $P_1, P_2, P_3, Q_1, Q_2, Q_3, R_1, R_2, R_3$ are concyclic then these belong to the circle of nine points of the triangle $ABC$.

(Ion Pătraşcu)

**52**. Prove that the perpendiculars constructed from the orthocenter of a triangle on the three concurrent Cevians of the triangle intersect the opposite sides of the triangle in three collinear points.

**53**. If $X, Y, Z$ are the tangency points with the circumscribed circle to the triangle $ABC$ of the mix-linear circumscribed circle corresponding to the angles $A, B, C$ respectively, then the lines $AX, BY, CZ$ are concurrent.

(P. Yiu)

**54**. Let $ABC$ a given triangle and $A'B'C'$, $A"B"C"$ two circumscribed triangles to $ABC$ and homological with $ABC$. Prove that in the triangle formed by the lines $A'A", B'B", C'C"$ is homological with the triangle $ABC$.

**55.** Let $ABC$ a triangle which is not rectangular and $H$ is its orthocenter, and $P$
A point on $(AH)$. The perpendiculars constructed from $H$ on $BP, CP$ intersect $AC, AB$ respectively in $B_1, C_1$. Prove that the lines $B_1C_1, BC$ are parallel.

(Ion Pătraşcu)

**56**. Let $MNPQ$ a quadrilateral inscribed in the circle $(O)$. We note $U$ the intersection point of its diagonals and let $[AB]$ a cord that passes through $U$ such that $AU = BU$. We note $\{U\} = AU \cap MN$, $\{W\} = PQ \cap BU$.

Prove that $UV = UW$

(The butterfly problem)

**57**. Let $A_1, A_2, A_3, A_4$ four points non-concyclic in a plane. We note $p_1$ the power of the point $A_1$ in rapport with the circle $(A_2A_3A_4)$, $p_2$ the power of the point $A_2$ in rapport with the circle $(A_1A_3A_4)$, with $p_3$ the power the point $A_3$ in rapport to the circle $(A_1A_2A_4)$ and $p_4$ the power of the point $A_4$ in rapport with the circle $(A_1A_2A_3)$. Show that it takes place the following relation $\frac{1}{p_1}+\frac{1}{p_2}+\frac{1}{p_3}+\frac{1}{p_4}=0$.

(Şerban Gheorghiu, 1945)

**58**.The quadrilaterals $ABCD, A'B'C'D'$ are conjugated. Prove that triangle $BCD$, $B'C'D'$ are homological. (The quadrilaterals follow the same sense in notations.)

**59**. Prove that an arbitrary given triangle $ABC$ is homological with the triangle where $A_1, B_1, C_1$ are the vertexes of an equilateral triangle constructed in the exterior of the triangle $ABC$ on it sides. (The homology center of these triangles is called the Toricelli-Fermat point).

**60**. Consider the points $A', B', C'$ on the sides $(BC)$, $(CA)$, $(AB)$ of the triangle $ABC$ which satisfy simultaneously the following conditions:

(i) $A'B^2 + B'C^2 + C'A^2 = A'C^2 + B'A^2 + C'B^2$;

(ii) The lines $AA', BB', CC'$ are concurrent

Prove that

a) The perpendiculars constructed from $A'$ on $BC$, from $B'$ on $AC$, from $C'$ on $AB$ are concurrent in appoint $P$;

b) The perpendiculars constructed from $A'$ on $B'C'$, from $B'$ on $A'C'$, from $C'$ on $A'B'$ are concurrent in a point $P'$;

c) The points $P$ and $P'$ are isogonal conjugated;

d) If $A'', B'', C''$ are the projections of the point $P'$ on $BC$, $AC$, $AB$, then the points $A', A'', B', B'' C', C''$ are concyclic;

e) The lines $AA'', BB'' CC''$ are concurrent.

(Florentin Smarandache, Ion Pătraşcu)

**61**. Consider a triangle $ABC$. In its exterior are constructed on the sides $(BC)$, $(CA)$, $(AB)$ squares. If $A_1, B_1, C_1$ are the centers of the three squares, prove that the triangles $ABC$, $A_1B_1C_1$ are homological. (The homology center is called the Vecten's point.)

**62**. A semi-circle has the diameter $(EF)$ situated on the side $(BC)$ of the triangle $ABC$ and it is tangent in the points $P, Q$ to the sides $AB$, $AC$.
Prove that the point $K$ common to the lines $EQ$, $FP$ belong the height from $A$ of the triangle $ABC$

**63**. In the triangle $ABC$ we know that $BC^2 = AB \cdot AC$; let $D, E$ the legs of the bisectrices of the angles $C, B$ ($D \in (AB)$, $E \in (AC)$). If $M$ is the middle of $(AB)$, $N$ is the middle of $(AC)$ and $P$ is the middle of $(DE)$
Prove that $M, N, P$ are collinear

(Ion Pătraşcu)

**64**. Let $ABC$ a right triangle in $A$. We'll construct the circles
$\mathcal{C}(A;BC), \mathcal{C}(B;AC), \mathcal{C}(C;AB)$:
Prove that:

(i) The circles $\mathcal{C}(A;BC), \mathcal{C}(B;AC), \mathcal{C}(C;AB)$ pass through the same point $L$;

(ii) If $Q$ is the second point of intersection of the circles $\mathcal{C}(B;AC), \mathcal{C}(C;AB)$, then the points $A, B, L, C, Q$ are concyclic;

(iii) If $P$ is the second point of intersection of the circles $\mathcal{C}(C;AB)$, $\mathcal{C}(A;BC)$ and $R$ is the second point of intersection of the circles $\mathcal{C}(B;AC)$, $\mathcal{C}(A;BC)$, then the points $P, Q, A, R$ are collinear

(Ion Pătraşcu)

**66**. If triangle $ABC$ is a scalene triangle and $a', b', c'$ are the sides of its orthic triangle then $4(a'b' + b'c' + c'a') \le a^2 + b^2 + c^2$

(Florentin Smarandache)

**67**. In an arbitrary triangle $ABC$ let $D$ the foot of the height from $A$, $G$ its weight center and $P$ the intersection of the line $(DG$ with the circumscribed circle to triangle $ABC$. Prove that $GP = 2GD$.

(Ion Pătraşcu)

**68**. Let $(AB)$ a cord in given circle. Through its middle we construct another cord $(CD)$. The lines $AC, BD$ intersect in a point $E$, and the lines $AD, BC$ intersect in a point $F$. Prove that the lines $EF, AB$ are parallel.

**69**. Let $ABCDEF$ complete quadrilateral ($\{E\} = AB \cap CD, \{F\} = ABC \cap AD$). A line intersects $(CD)$ and $(AB)$ in $C_1, A_1$.
Prove that:

(i) The lines $A_1B_1, C_1D_1, AC$ are concurrent;

(ii) The lines $B_1C_1, A_1D_1, BD$ are concurrent.

**70**. In a triangle $ABC$ let $AA'AA''$ two isotomic Cevians and $PQ \| BC$, $P \in (AB)$, $Q \in (AC)$. We'll note $\{M\} = BQ \cap AA'$, $\{N\} = CP \cap AA'$. Prove that $MN \| BC$.

**71**. In the triangle $ABC$ let consider $AA', BB', CC'$ the concurrent Cevians in the point $P$. Determine the minimum values of the expressions:

$$E(P) = \frac{AP}{PA'} + \frac{BP}{PB'} + \frac{CP}{PC'}$$

$$F(P) = \frac{AP}{PA'} \cdot \frac{BP}{PB'} \cdot \frac{CP}{PC'}$$

where $A' \in [BC], B' \in [CA], C' \in [AB]$

(Florentin Smarandache)

**72**. Let triangles $A_1B_1C_1$, $A_2B_2C_2$ such that

$$B_1C_1 \cap B_2C_2 = \{P_1\}, B_1C_1 \cap A_2C_2 = \{Q_1\}, B_1C_1 \cap A_2B_2 = \{R_1\}$$

$$A_1C_1 \cap A_2C_2 = \{P_2\}, A_1C_1 \cap A_2C_2 = \{Q_2\}, A_1C_1 \cap C_2B_2 = \{R_2\}$$

$$A_1B_1 \cap A_2B_2 = \{P_3\}, A_1B_1 \cap B_2C_2 = \{Q_3\}, A_1B_1 \cap C_2A_2 = \{R_3\}$$

Prove that

(i) $$\frac{P_1B_1}{P_1C_1} \cdot \frac{P_2C_1}{P_2A_1} \cdot \frac{P_3A_1}{P_3B_1} \cdot \frac{Q_1B_1}{Q_1C_1} \cdot \frac{Q_2C_1}{Q_2A_1} \cdot \frac{Q_3A_1}{P_3B_1} \cdot \frac{R_1B_1}{R_1C_1} \cdot \frac{R_2C_1}{R_2A_1} \cdot \frac{R_3A_1}{R_3B_1} = 1$$

(ii)Triangles $A_1B_1C_1$, $A_2B_2C_2$ are homological (the lines $A_1A_2, B_1B_2, C_1C_2$ are concurrent) if and only if $\frac{Q_1B_1}{Q_1C_1} \cdot \frac{Q_2C_1}{Q_2A_1} \cdot \frac{Q_3A_1}{P_3B_1} = \frac{R_1B_1}{R_1C_1} \cdot \frac{R_2C_1}{R_2A_1} \cdot \frac{R_3A_1}{R_3B_1}$

**73**. On a line we consider three fixed points $A, B, C$. Through the points $A, B$ we construct a variable arc, and from $C$ we construct the tangents to the circle $CT_1, CT_2$. Show that the line $T_1T_2$ passes through a fixed point.

**74**. In the triangle $ABC$ we construct the concurrent Cevians $AA_1, BB_1, CC_1$ such that $A_1B^2 + B_1C^2 + C_1A^2 = AB_1^2 + BC_1^2 + CA_1^2$ and one of them is a median. Show that the other two Cevians are medians or that the triangle $ABC$ is isosceles.

(Florentin Smarandache)

**75**. Let $ABC$ a triangle and $A_1, B_1, C_1$ points on its exterior such that

$\sphericalangle A_1BC \equiv \sphericalangle C_1BA, \sphericalangle C_1AB \equiv \sphericalangle B_1AC, \sphericalangle B_1CA \equiv \sphericalangle A_1CB$.

Prove that the triangles $ABC$, $A_1B_1C_1$ are homological.

**76**. Show that if in a triangle we can inscribe three conjugated squares, then the triangle is equilateral.

**77**. Let a mobile point $M$ on the circumscribed circle to the triangle $ABC$. The lines $BM, CM$ intersect the sides $AC, AB$ in the points $N, P$.
Show that the line $NP$ passes through a fixed point.

**78**. Given two fixed points $A, B$ on the same side of a fixed line $d$ ($AB \nparallel d$). A variable circle which passes through $A, B$ intersects $d$ in $C, D$. Let $\{M\} = AC \cap BD$, $\{N\} = AD \cap BC$. Prove that the line $MN$ passes through a fixed point

**79**. In the triangle $ABC$ we have $A = 96°, B = 24°$. Prove that $OH = a - b$, ($O$ is the center of the circumscribed circle and $H$ is the orthocenter of the triangle $ABC$)

(Ion Pătraşcu), G.M 2010

**80**. Let $ABCD$ a quadrilateral convex inscribed in a circle with the center in $O$. We will note $P, Q, R, S$ the middle points of the sides $AB, BC, CD, DA$. If $M$ is a point such that $2\overrightarrow{OM} = \overrightarrow{OA} + \overrightarrow{OB} + \overrightarrow{OC} + \overrightarrow{OD}$ and $T$ is the intersection of the lines $PR, QS$,
Prove that:

(i) $2\overrightarrow{OT} = \overrightarrow{OM}$

(ii) If we note $P', Q', R', S'$ the orthogonal projection of the points $P, Q, R, S$ respectively on $CD, DA, AB, BC$, the lines $PP', QQ', RR', SS'$ are concurrent. (The intersection point is called the Mathot's point of the quadrilateral.)

**81**. Three conjugated circles are situated in the interior of a triangle and each of them is tangent to two of the sides of the triangle. All three circles pass through the same point. Prove that the circles' common point and the centers of the inscribed and circumscribed to the triangle are collinear points.

(O.I.M. 1981)

**82**. Let $O$ the center of the circumscribed circle of the triangle $ABC$ and $AA', BB', CC'$ the heights of the triangle. The lines $AO, BO, CO$ intersect respectively the lines $B'C', A'C', A'B'$ in $A_1, B_1, C_1$.

Prove that the centers of the circumscribed circles to triangles $AA'A_1, BB'B_1, CC'C_1$ are collinear.

(A. Angelescu, G.M.)

**83**. Let $ABC$ an equilateral triangle and $AA_1, BB_1, CC_1$ three Cevians concurrent in this triangles.

Prove that the symmetric of each of the Cevians in rapport to the opposite side of the triangle $ABC$ are three concurrent lines.

**84**. Given three circles $\mathcal{C}(O_1, R), \mathcal{C}(O_2, R), \mathcal{C}(O_3, R)$ that have a common point $O$. Let $A', B', C'$ the diametric points of $O$ in the three circles. The circumscribed circle to the triangles $B'OC', C'OA', A'OB'$ intersect the second time the circles $\mathcal{C}(O_1, R), \mathcal{C}(O_2, R), \mathcal{C}(O_3, R)$ respectively in the points $A_1', B_1', C_1'$.

Prove that the points $A_1', B_1', C_1'$ are collinear.

(Ion Pătraşcu)

85 Let $d_1, d_2, d_3$ three parallel lines. The triangles $ABC$, $A'B'C'$ have $A, A'$ on $d_1$, $B, B'$ on $d_2$, $C, C'$ on $d_3$ and the same weight center $G$.

If $\{U\} = BC \cap B'C'; \{V\} = AC \cap A'C'; \{W\} = AB \cap A'B'$.

Prove that the points $U, V, W, G$ are collinear.

**86**. In the triangle $ABC$ let $AA', BB', CC'$ its interior bisectrices. The triangle determined by the mediators of the segments $(AA'), (BB'), (CC')$ is called the first triangle of Sharygin; prove that the triangle $ABC$ and its first triangle of Sharygin are homological, the homology axis being the Lemain's line of the triangle $ABC$.

**87**. Let $ABCDEF$ a hexagon inscribed in a circle. We note $\{M\} = AC \cap BD, \{N\} = BE \cap CF, \{P\} = AE \cap DF$.

Prove that the points $M, N, P$ are collinear.

**88**. Let $ABC$ a triangle in which $AB < AC$ and let triangle $SDE$, where $S \in (BC)$, $D \in (AB)$, $E \in (AC)$ such that $\Delta SDE \sim \Delta ABC$ ($S$ being different of the middle of $(BC)$).

Prove that $\frac{BS}{CS} = \left(\frac{AB}{AC}\right)^2$

(Ion Pătraşcu G.M. 205)

**89**. Let $ABC$ a random given triangle, $C_aC_bC_c$ its contact triangle and $P$ a point in the interior of triangle $ABC$. We'll note $A_1, B_1, C_1$ the intersections of the inscribed circle with the semi-lines $(C_aP, (C_bP, (C_cP$.
Prove that the triangles $ABC$, $A_1B_1C_1$ are homological.

**90**. In a random triangle $ABC$, $O$ is the center of the circumscribed circle and $O_1$ is the intersection point between the mediator of the segment $(OA)$ with the parallel constructed through $O$ to $BC$. If $A'$ is the projection of $A$ on $BC$ and $D$ the intersection of the semi-line $(OA'$ with the circle $C(O_1, O_1A)$, prove that the points $B, O, C, D$ are concyclic.

(Ion Pătraşcu)

**91**. If in triangle $ABC$, $I$ is the center of the inscribed circle, $A', B', C'$ are the projections of $I$ on $BC, CA, AB$ and $A'', B'', C''$ points such that $\overrightarrow{IA''} = \overrightarrow{KIA'}$, $\overrightarrow{IB''} = \overrightarrow{KIB'}$, $\overrightarrow{IC''} = \overrightarrow{KIC'}$, $K \in R^*$, then the triangle $ABC$, $A''B''C''$ are homological. The intersection point of the lines $AA'', BB'', CC''$ is called the Kariya point of the triangle $ABC$.

**92.** Let $ABC$ a random triangle and $C_aC_bC_c$ its contact triangle. The perpendiculars in the center of the inscribed circle $I$ of the triangle $ABC$ on $AI, BI, CI$ intersect $BC, CA, AB$ respectively in points $A_1, B_1, C_1$ and the tangents to the inscribed circle constructed in these points intersect $BC, CA, AB$ respectively in the points $A_2, B_2, C_2$. Prove that

i. The points $A_1, B_1, C_1$ are collinear;
ii. The points $A_2, B_2, C_2$ are collinear;
iii. The lines $A_1B_1$, $A_2B_2$ are parallel.

(Ion Pătraşcu)

**93**. Show that the parallels constructed through orthocenter of a triangle to the external bisectrices of the triangle intersect the corresponding sides of the triangle in three collinear points.

**94**. Let a circle with the center in $O$ and the points $A,B,C,D,E,F$ on this circle. The circumscribed circles to triangles $(AOB),(DOE);(BOC),(EOF);(COD),(FOA)$ intersect the second time respectively in the points $A_1,B_1,C_1$.
Prove that the points $O,A_1,B_1,C_1$ are concyclic.

**95**. Let $ABC$ an isosceles triangle $AB=AC$ and $D$ a point diametric opposed to $A$ in the triangle's circumscribed circle. Let $E\in(AB)$ and $\{P\}=DE\cap BC$ and $F$ is the intersection of the perpendicular in $P$ on $DE$ with $AC$. Prove that $EF=BE=CF$.

(Ion Pătraşcu)

**96**. In the triangle $ABC$ we'll note $A'B'C'$ the circumpedal triangle of the center of the circumscribed circle in the triangle $ABC$ and $O_1,O_2,O_3$ the symmetric point to the center $O$ of circumscribed circle to triangle $ABC$ in rapport to $B'C',A'B',A'C'$ respectively.
Prove that the triangles $ABC$, $O_1O_2O_3$ are homological, the homology center being the Kariya's point of the triangle $ABC$.

**97**. Let $ABCDEF$ a complete quadrilateral in which $BE=DF$. We note $\{G\}=BD\cap EF$.
Prove that the Newton-Gauss lines of the quadrilaterals $ABCDEF$, $EFDBAG$ are perpendicular

(Ion Pătraşcu)

**98**. Let in a random triangle $ABC$, the point $O$ the center of the circumscribed circle and $M_a$, $M_b$, $M_c$ the middle points of the sides $(BC),(CA),(AB)$. If $K\in R^*$ and $A',B',C'$ three points such that $\overrightarrow{OA'}=K\cdot\overrightarrow{OM_a},\overrightarrow{OB'}=K\cdot\overrightarrow{OM_b},\overrightarrow{OC'}=K\cdot\overrightarrow{OM_c}$. Prove that the triangles $ABC$, $A'B'C'$ are homological. (The homology center is called the Franke's point of the triangle $ABC$.)

**99**. Let $ABC$ a scalene triangle and $M,N,P$ the middle points of the sides $(BC),(CA),(AB)$. We construct three circles with the centers in $M,N,P$, and which intersect the $(BC),(CA),(AB)$ respectively in $A_1,A_2$; $B_1,B_2$; $C_1,C_2$ such that these six points are concyclic.
Prove that the three circle of centers $M,N,P$ have as radical center the orthocenter of the triangle $ABC$.

(In connection with the problem 1 –O.I.M-2008, Ion Pătraşcu)

**100**. Let $M$ an arbitrary point in the plane of triangle $ABC$.
Prove that the tangents constructed in $M$ to the circles $(BMC),(CMA),(AMB)$ intersect $(BC),(CA),(AB)$ in collinear points.

(Cezar Coşniţă)

## 5.2. Open problems

In this section we selected and proposed a couple of problems for which we didn't find a solution or for which there is not a known elementary solution.

**1.** A diameter of the circle $C(R,O)$ circumscribed to triangle $ABC$ intersects the sides of the triangle in the points $A_1, B_1, C_1$. We'll consider $A', B', C'$ the symmetric points in rapport with the center $O$. Then the lines $AA', BB', CC'$ are concurrent in a point situated on the circle.

(Papelier)

**2.** A transversal intersects the sides $BC, CA, AB$ of a triangle $ABC$ in the points $A', B', C'$. The perpendiculars constructed on $A', B', C'$ on the sides $BC, CA, AB$ form a triangle $A"B"C"$.
Prove that the triangles $ABC$, $A"B"C"$ are homological, the homology center belonging to the circumscribed circle of the triangle $ABC$.

(Cezar Coşniţă)

**3.** Let $ABC$ a triangle, $A_1B_1C_1$ the Caspary's first triangle and the triangle $YZZ'$ formed by the homological centers of the triangles $ABC$, $A_1B_1C_1$.
Show that these triangles have the same weight center.

(Caspary)

**4.** In triangle $ABC$ we'll note with $M, N, P$ the projections of the weight center $G$ on the sides $BC, CA, AB$ respectively.
Show that if $AM, BM, CP$ are concurrent then the triangle $ABC$ is isosceles.

(Temistocle Bîrsan)

**5.** Let triangle $ABC$, the Cevians $AA_1, BB_1, CC_1$ concurrent in $Q_1$ and $AA_2, BB_2, CC_2$ concurrent in $Q_2$. We note $\{X\} = B_1C_1 \cap B_2C_2, \{Y\} = C_1A_1 \cap C_2A_2, \{Z\} = A_1B_1 \cap A_2B_2$.
Show that

a. $AX, BY, CZ$ are concurrent;

b. The points $A, Y, Z$ ; $B, Y, Z$ and $C, Y, Z$ are collinear.

(Cezar Coşniţă)

**6.** Through the point $M$ of a circumscribed circler to a triangle $ABC$ we'll construct the circles tangent in the points $B, C$ to $AB$ and $AC$. These circle intersect for the second time in a point $A'$ situated on the side $BC$. If $B', C'$ are the points obtained in a similar mode as $A'$, prove that the triangles $ABC$, $A'B'C'$ are homological.

(Cezar Coşniţă)

**7.** Prove that the only convex quadrilateral $ABCD$ with the property that the inscribed circles in the triangles $AOB, BOC, COD, DOA$ are congruent is a rhomb ($\{O\} = AC \cap BD$).

(Ion Pătraşcu)

**8.** Let $AA_1$, $AA_2$; $BB_1, BB_2$; $CC_1, CC_2$ three pairs of isogonal Cevians in rapport with the angles $A, B, C$ of the triangle $ABC$ ( $A_i \in BC, B_i \in CA, C_i \in AB$ ). The points $X, Y, Z$ and $X', Y', Z'$ being defined as follows $\{X\} = A_1B_1C_1 \cap A_2C_2, \{Y\} = A_1B_2 \cap A_2C_1$, etc. Prove the following:

a) The lines $AX, BY, CZ$ are concurrent in a point $P$;
b) The lines The lines $AX', BY', CZ'$ are concurrent in a point $P'$;
c) The points $P$, $P'$ are isogonal conjugated.

(Temistocle Bîrsan)

**9.** If $T_1, T_2, T_3$ are triangles in plane such that $(T_1, T_2)$ are tri-homological, $(T_2, T_3)$ are tri-homological, $(T_3, T_1)$ are tri-homological and these pairs of tri-homological triangles have each of them in common two homological centers, then the three homological centers left non-common are collinear?

**10.** If $ABCDEF$ is a complete quadrilateral, $U, V, W$ are collinear points situated respectively on $AC, BD, EF$ and there exists $Q$ on $(BE)$ such that $\{P\} = VQ \cap BC; \{R\} = WQ \cap CE$ and $U, P, R$ are collinear, then the triangles $BCE, QRP$ are homological.

**11.** It is known that a triangle $ABC$ and its anti-supplementary triangle are homological. Are homological also the triangle $ABC$ and the pedal triangle of a point $M$ from the triangle's $ABC$ plane. It is also known that the anti-supplementary of triangle $ABC$ and the pedal of $M$ are homological.
What property has $M$ if the three homological centers of the pairs of triangles mentioned above are collinear?

**12.** If $P$ is a point in the plane of triangle $ABC$, which is not situated on the triangle's circumscribed circle or on its sides; $A'B'C'$ is the pedal triangle of $P$ and $A_1, B_1, C_1$ three points such that $\overrightarrow{PA'} \cdot \overrightarrow{PA_1} = \overrightarrow{PB'} \cdot \overrightarrow{PB_1} = \overrightarrow{PC'} \cdot \overrightarrow{PC_1} = K$, $K \in R^*$. Prove that the triangles $ABC$, $A_1B_1C_1$ are homological.

(The generalization of the Cezar Coşniță's theorem).

# Chapter 6

# Notes

## 6.1. Menelaus' theorem and Ceva's theorem

### I. Menelaus' theorem

**Definition 1**

A line, which intersects the three sides of a triangle is called the triangle transversal. If the intersection points with $BC, CA, AB$ are respectively $A_1, B_1, C_1$, we'll note the transversal $A_1 - B_1 - C_1$.

**Theorem 1** (Chapter 1, 2)

If in the triangle $ABC$, $A_1 - B_1 - C_1$ is a transversal, then

$$\frac{A_1B}{A_1C} \cdot \frac{B_1C}{B_1A} \cdot \frac{C_1A}{C_1B} = 1 \qquad (1)$$

**Proof**

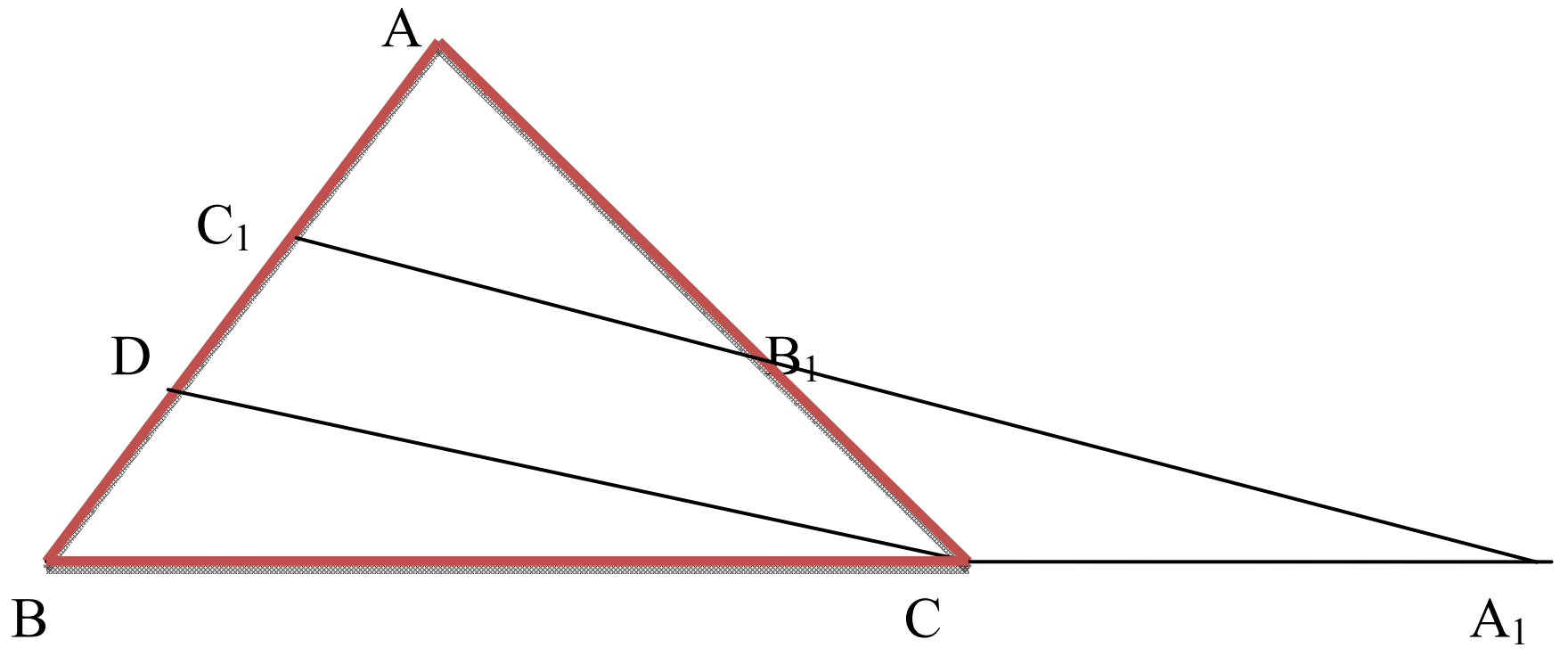

Fig.1

We construct $CD \parallel A_1B_1$, $D \in AB$ (see Fig.1). Using the Thales' theorem awe have

$$\frac{A_1B}{A_1C} = \frac{C_1B}{C_1D}, \quad \frac{B_1C}{B_1A} = \frac{C_1D}{C_1A}$$

Multiplying side by side these relations, after simplification we obtain the relation (1).

**Theorem 2** (The reciprocal of Menelaus' theorem)

If the points $A_1, B_1, C_1$ are respectively on the sides $BC, CA, AB$ of the triangle $ABC$ and it takes place the relation (1), then the points $A_1, B_1, C_1$ are collinear.

**Proof**

We suppose by absurd that $A_1 - B_1 - C_1$ is not a transversal. Let then $\{C_1^{'}\} = A_1B_1 \cap AB$, $C_1^{'} \neq C_1$. We'll apply the Menelaus' theorem for the transversal $A_1 - B_1 - C_1^{'}$, and we have:

$$\frac{A_1B}{A_1C}\cdot\frac{B_1C}{B_1A}\cdot\frac{C_1^{'}A}{C_1^{'}B}=1 \qquad (2)$$

From the relations (1) and (2) we find:

$$\frac{C_1A}{C_1B}=\frac{C_1^{'}A}{C_1^{'}B}$$

And from this it results that $C_1 = C_1^{'}$, which contradicts the initial supposition. Consequently, $A_1 - B_1 - C_1$ is a transversal in the triangle $ABC$.

## II. Ceva's theorem

**Definition 2**

A line determined by a vertex of a triangle and a point on the opposite side is called the triangle's Cevian.

**Note**

The name of this line comes from the Italian mathematician Giovanni Ceva (1647-1734)

**Theorem 3** (G. Ceva – 1678)

If in a triangle $ABC$ the Cevians $AA_1, BB_1, CC_1$ are concurrent then

$$\frac{\overline{A_1B}}{\overline{A_1C}}\cdot\frac{\overline{B_1C}}{\overline{B_1A}}\cdot\frac{\overline{C_1A}}{\overline{C_1B}}\frac{C_1^{'}A}{C_1^{'}B}=-1 \qquad (3)$$

**Proof**

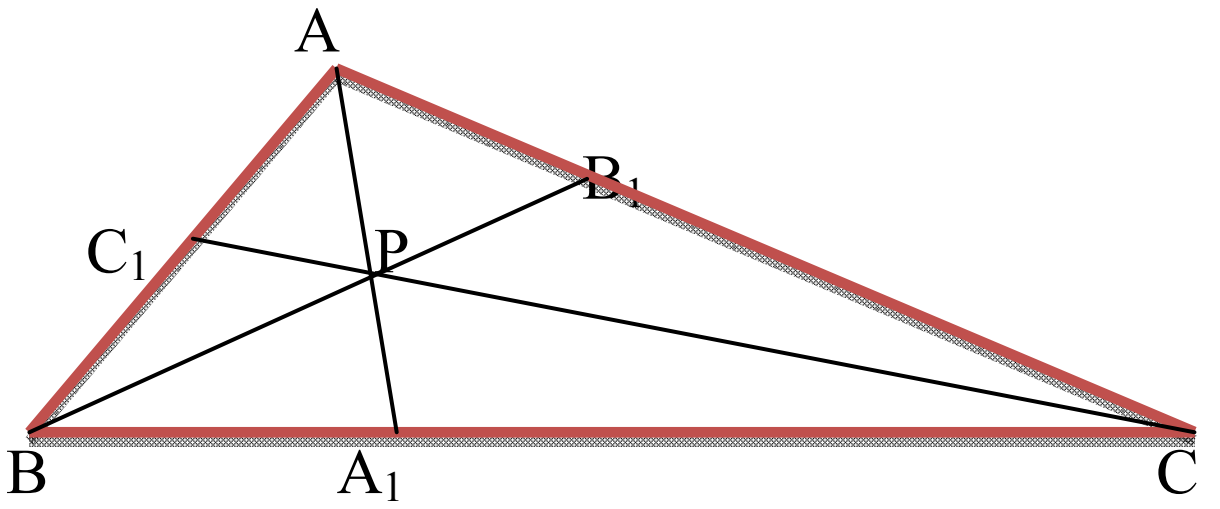

Fig. 2

Let $\{O\} = AA_1 \cap BB_1 \cap CC_1$. We'll apply the Menelaus' theorem in the triangles $AA_1C, AA_1B$ for the transversals $B - P - B_1, C - P - C_1$, we obtain

$$\frac{\overline{BA_1}}{\overline{BC}}\cdot\frac{\overline{B_1C}}{\overline{B_1A}}\cdot\frac{\overline{PA}}{\overline{PA_1}}=1 \qquad (4)$$

$$\frac{\overline{CB_1}}{\overline{CA_1}} \cdot \frac{\overline{PB_1}}{\overline{PA_1}} \cdot \frac{\overline{C_1A}}{\overline{C_1B}} = 1 \qquad (5)$$

By multiplying side by side the relations (4) and (5) and taking into account that

$$\frac{\overline{BC}}{\overline{CB}} = -1 \text{ and } \frac{\overline{BA_1}}{\overline{CA_1}} = \frac{\overline{A_1B}}{\overline{A_1C}}$$

We'll obtain he relation (3).

**Theorem 4** (The reciprocal of Ceva's theorem)

If $AA_1, BB_1, CC_1$ are Cevians in the triangle $ABC$ such that the relation (3) is true, then these Cevians are concurrent.

**Proof**

The proof is done using the method of reduction ad absurdum.

**Lemma 1**

If $A_1$ is a point on the side $BC$ of the triangle $ABC$, then

$$\frac{\overline{AC}}{\overline{A_1C}} = \frac{\sin \sphericalangle AA_1C}{\sin \sphericalangle A_1AC}.$$

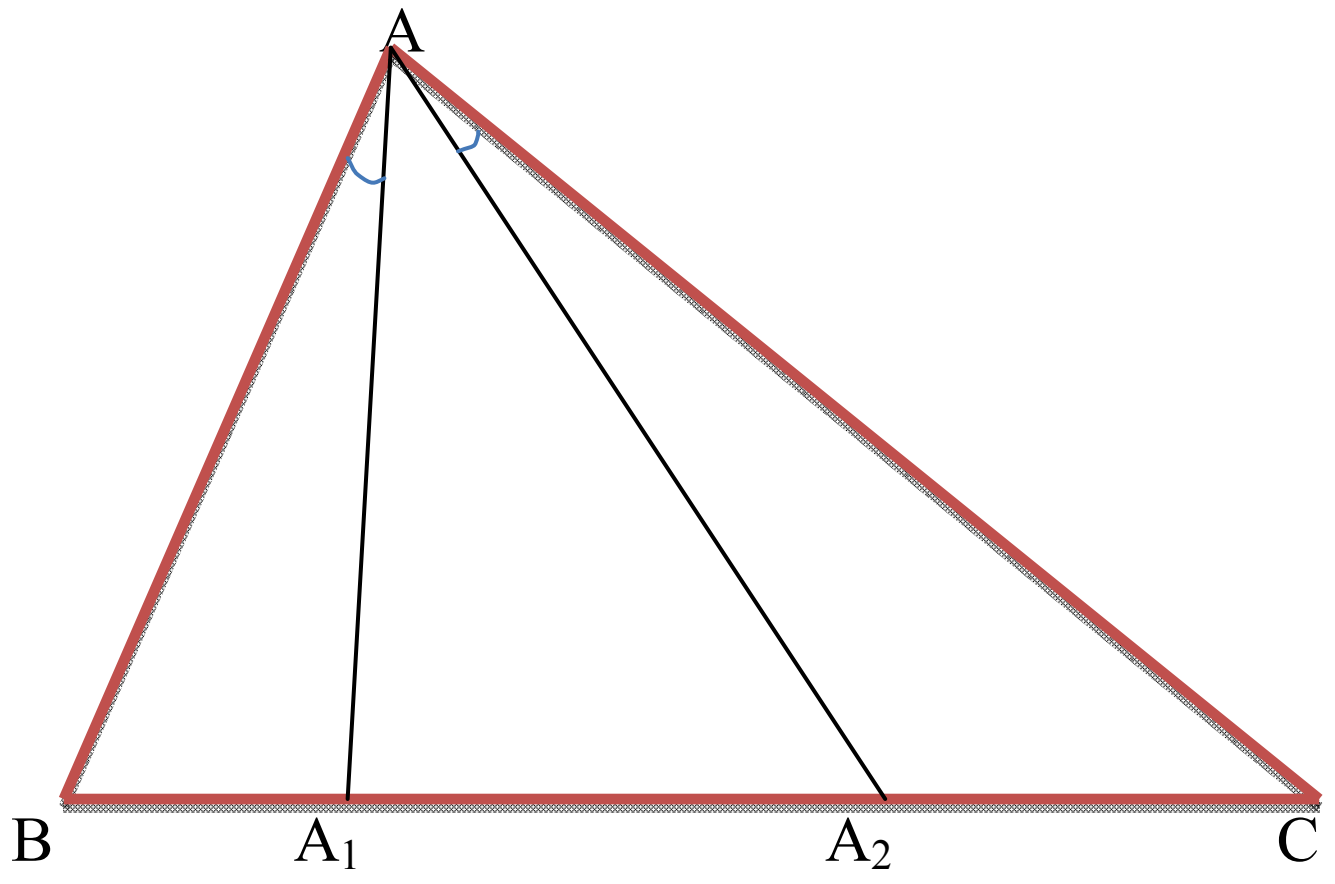

Fig. 3

Because $\sin \sphericalangle AA_1B = \sin \sphericalangle AA_1C$ and by multiplying the precedent relations side by side we obtain the relation from the hypothesis.

**Observation 1**

The relation from the hypothesis can be obtained also from

$$\frac{\overline{A_1B}}{\overline{A_1C}} = \frac{S_{ABA_1}}{S_{ACA_1}}$$

**Corollary**

If $AA_1$ and $AA_2$ are isogonal Cevians in the triangle $ABC$ ($\sphericalangle A_1AB = \sphericalangle A_2AC$), (see Fig. 3), then from Lemma 1 it results:

$$\frac{\overline{A_1B}}{\overline{A_1C}} \cdot \frac{\overline{A_2C}}{\overline{A_2C}} = \left(\frac{AB}{AC}\right)^2$$

**Theorem 5** (The trigonometric form of Ceva's theorem)
In the triangle $ABC$, the Cevians $AA_1, BB_1, CC_1$ are concurrent if and only if

$$\frac{\sin \sphericalangle A_1AB}{\sin \sphericalangle A_1AC} \cdot \frac{\sin \sphericalangle B_1BC}{\sin \sphericalangle B_1BA} \cdot \frac{\sin \sphericalangle C_1CA}{\sin \sphericalangle C_1CB} = -1$$

To prove this it can be used the Lemma 1 and theorems 3,4.

**Theorem 6** (The trigonometric form of Menelaus' theorem)
Three points $A_1, B_1, C_1$ situated respectively on the opposite sides of the triangle $ABC$ are collinear if and only if

$$\frac{\sin \sphericalangle A_1AB}{\sin \sphericalangle A_1AC} \cdot \frac{\sin \sphericalangle B_1BC}{\sin \sphericalangle B_1BA} \cdot \frac{\sin \sphericalangle C_1CA}{\sin \sphericalangle C_1CB} = -1$$

## III. Applications

1. If $AA_1, BB_1, CC_1$ are three Cevians in the triangle $ABC$ concurrent in the point $P$ and $B_1C_1$ intersects $BC$ in $A_2$, $A_1B_1$ intersects $AB$ in $C_2$ and $C_1A_1$ intersects $AC$ in $B_2$, then

(i) The points $A_2, B_2, C_2$ are harmonic conjugates of the points $A_1, B_1, C_1$ in rapport to $B, C; C, A$ respectively $A, B$;

(ii) The points $A_2, B_2, C_2$ are collinear

**Proof**

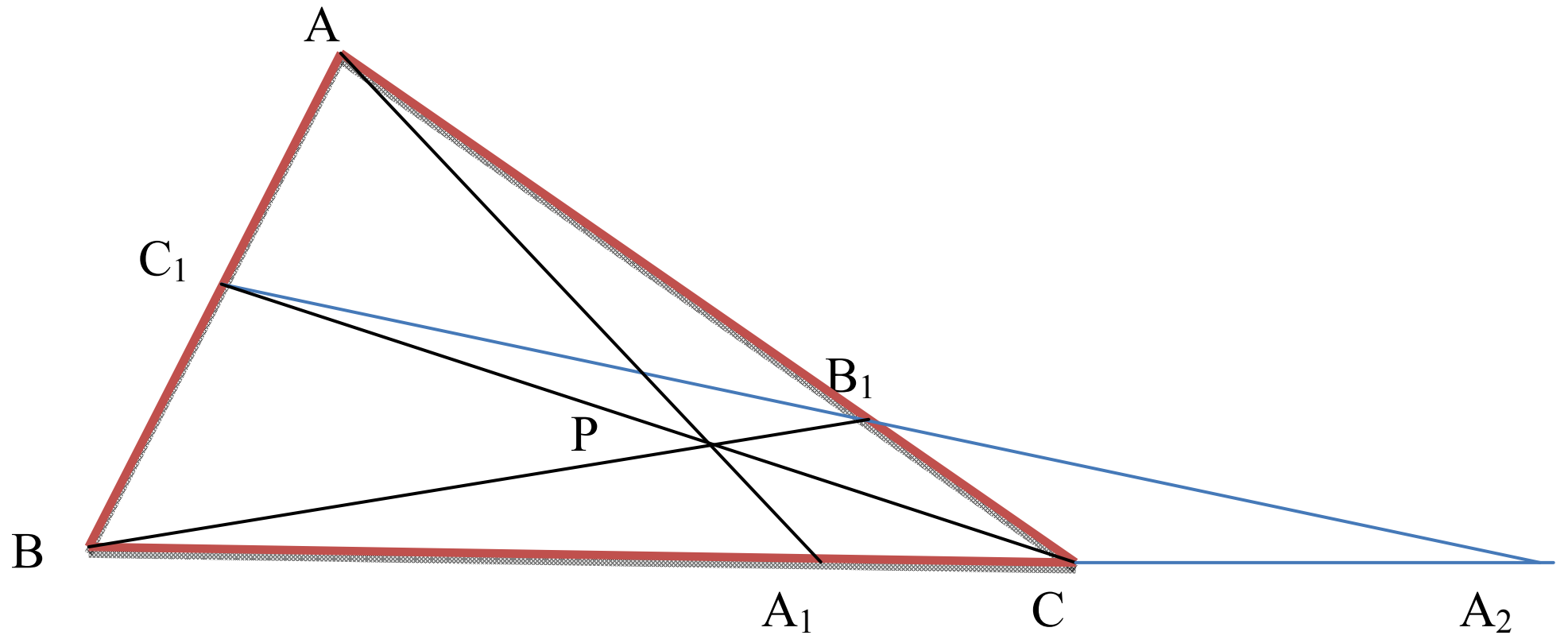

Fig. 4

(i) From the Ceva's theorem we have

$$\frac{\overline{A_1B}}{\overline{A_1C}} \cdot \frac{\overline{B_1C}}{\overline{B_1A}} \cdot \frac{\overline{C_1A}}{\overline{C_1B}} \frac{C_1'A}{C_1'B} = -1$$

From the Menelaus' theorem applied in the triangle $ABC$ for the transversal $A_2 - B_1 - C_1$ it results:

$$\frac{\overline{A_2B}}{\overline{A_2C}} \cdot \frac{\overline{B_1C}}{\overline{B_1A}} \cdot \frac{\overline{C_1A}}{\overline{C_1B}} = 1$$

Dividing side by side the precedent relations it results $\frac{\overline{A_1B}}{\overline{A_1C}} = -\frac{\overline{A_2B}}{\overline{A_2C}}$ which shows that the point $A_2$ is the harmonic conjugate of the point $A_1$ in rapport with the points $B, C$.

(ii) This results from (i) and from the reciprocal of Menelaus' theorem.

**Remark 1**

The line $A_2, B_2, C_2$ is called the harmonic associated to the point $P$ or the tri-linear polar of the point $P$ in rapport to the triangle $ABC$. Inversely, the point $P$ is called the tri-linear pole or the harmonic pole associated to the line $A_2B_2C_2$.

2. If $P_1$ is a point halfway around the triangle from $A$, that is $\overline{AB} + \overline{BP_1} = \overline{VC} + \overline{CA}$ an $P_2, P_3$ are similarly defined, then the lines $AP_1, BP_2, CP_3$ are concurrent (The concurrence point is called the Nagel's point of the triangle)

**Proof**

We find that $A_2 = p - b, P_1A_3 = c - p$ and the analogues ad then it is applied the reciprocal Ceva's theorem.

# Annex 1

## Important formulae in the triangle geometry

The cosine theorem

$$a^2 = b^2 + c^2 - 2bc\cos A$$

$$\sin\frac{A}{2} = \sqrt{\frac{(p-b)(p-c)}{bc}};\ \cos\frac{A}{2} = \sqrt{\frac{p(p-a)}{bc}};\ tg\frac{A}{2} = \sqrt{\frac{(p-b)(p-c)}{p(p-c)}}$$

The median relation

$$4ma^2 = 2\left(b^2 + c^2\right) - a^2$$

The tangents' theorem

$$\frac{tg\frac{A-B}{2}}{tg\frac{A+B}{2}} = \frac{a-b}{a+b}$$

$$r = \frac{s}{p},\ \ R = \frac{abc}{4S};\ r_a = \frac{s}{p-a}$$

$$a^2 + b^2 + c^2 = 2\left(p^2 - r^2 - 4Rr\right)$$

$$a^3 + b^3 + c^3 = 2p\left(p^2 - 3r^2 - 6Rr\right)$$

$$a^4 + b^4 + c^4 = 2\left[p^4 - 2p^2r(4R+3r) + r^2(4R+r)^2\right]$$

$$ab + bc + ca = p^2 + r^2 + 4Rr$$

$$abc = 4Rpr$$

$$16S^2 = 2\left(a^2b^2 + b^2c^2 + c^2a^2\right) - \left(a^4 + b^4 + c^4\right)$$

We prove:

$$ab + bc + ca = p^2 + r^2 + 4Rr$$

From $S^2 = p(p-a)(p-b)(p-c)$ and $S^2 = p^2 - r^2$ we find $(p-a)(p-b)(p-c) = p\cdot r^2$.
It results

$$p^3 - (a+b+c)p^2 + (ab+bc+ca)p - abc = pr^2$$

$$\Leftrightarrow p^3 - 2p^3 + (ab+bc+ca)p - 4Rpr = pr^2$$

From here we retain $ab + bc + ca = p^2 + r^2 + 4Rr$; $a^2 + b^2 + c^2 = 2\left(p^2 - r^2 - 4Rr\right)$
$a^3 + b^3 + c^3 = 2p\left(p^2 - 3r^2 - 6Rr\right)$.

We use the identity $a^3 + b^3 + c^3 - 3abc = (a+b+c)\left(a^2 + b^2 + c^2 - ab - bc - ca\right)$.

It results

$$a^3+b^3+c^3=2p\left(2p^2-2r^2-8Rr-4Rr\right)+12Rp^2=2p\left(p^2-3r^2-6Rr\right).$$

## C. Distances between remarkable points in triangle geometry

**1.** $CG^2=\frac{1}{9}\left(9R^2-2p^2+2r^2+8Rr\right)$

**2.** $OH^2=9R^2-2p^2+2r^2+8Rr$

**3.** $OI^2=R^2-2Rr$

**4.** $OI_a^2=R^2+2Rr_a$

**5.** $OK^2=R^2-\frac{3a^2b^2c^2}{\left(a^2+b^2+c^2\right)^2}$

**6.** $O\Gamma^2=R^2-\frac{4p^2r(R+r)}{(4R+r)^2}$

**7.** $ON^2=(R-2r)^2$

**8.** $GH^2=\frac{4}{9}\left(9R^2-2p^2+2r^2+8Rr\right)$

**9.** $GI^2=\frac{1}{9}\left(p^2+5r^2-16Rr\right)$

**10.** $GI_a^2=\frac{2}{9}\left(a^2+b^2+c^2\right)-\frac{1}{6(p-a)}\cdot\left(-a^2+b^2+c^2\right)+2Rra$

**11.** $GK^2=\frac{2\left(a^2+b^2+c^2\right)^3-3\left(a^2+b^2+c^2\right)\left(a^4+b^4+c^4\right)-27a^2b^2c^2}{9\left(a^2+b^2+c^2\right)^2}$

**12.** $G\Gamma^2=\frac{4}{9(4R+r)^2}p^2\left[\left(4R^2+8Rr-5r^2\right)-r(4R+r)^3\right]$

**13.** $GN^2=\frac{4}{9}\left(p^2+5r^2-16Rr\right)$

**14.** $HI^2=4R^2-p^2+3r^2-4Rr$

**15.** $HI_a^2=4R^2+2r_a^2+r^2+4Rr-p^2$

**16.** $HK^2=4R^2-\left(a^2+b^2+c^2\right)+\frac{2\left(a^2+b^2+c^2\right)^3\left(a^2b^2+b^2c^2+c^2a^2\right)-3a^2b^2c^2}{\left(a^2+b^2+c^2\right)^2}$

**17.** $H\Gamma^2=4R^2\left[1-\frac{2p^2(2R-r)}{R(4R+r)^2}\right]$

**18.** $HN^2=4R(R-2r)$

**19.** $II_a^2 = 4R(r_a - r)$

**20.** $I\Gamma^2 = r^2\left[1 - \frac{3p^2}{(4R+r)^2}\right]$

**21.** $IN^2 = p^2 + 5r^2 - 16Rr$

**22.** $IK^2 = \frac{4r^2R\left[p^2(R+r) - r(r+4R)^2\right]}{(p^2 - r^2 - 4Rr)^2}$

We'll prove the formulae 14 and 15

The position's vector of the inscribed circle in the triangle $ABC$ is

$$\overrightarrow{PI} = \frac{1}{2p}\left(a\overrightarrow{PA} + b\overrightarrow{PB} + c\overrightarrow{PC}\right)$$

If $H$ is the orthocenter of the triangle $ABC$ then $\overrightarrow{HI} = \frac{1}{2p}\left(a\overrightarrow{HA} + b\overrightarrow{HB} + c\overrightarrow{HC}\right)$

Let's evaluate $\overrightarrow{HI}\cdot\overrightarrow{HI}$

$$HI^2 = \frac{1}{4p^2}\left(a^2HA^2 + b^2HB^2 + c^2HC^2 + 2ab\overrightarrow{HA}\overrightarrow{HB} + 2bc\overrightarrow{HB}\overrightarrow{HC} + 2ac\overrightarrow{HA}\overrightarrow{HC}\right)$$

If $A_1$ is the middle of $BC$, we have $AH = 2OA_1$ then $AH^2 = 4R^2 - a^2$, similarly $BH^2 = 4R^2 - b^2$ and $CH^2 = 4R^2 - c^2$. Also, $\overrightarrow{HA}\overrightarrow{HB} = \left(\overrightarrow{OA} + \overrightarrow{OB}\right)\left(\overrightarrow{OC} + \overrightarrow{OA}\right) = 4R^2 - \frac{1}{2}\left(a^2 + b^2 + c^2\right)$.

Taking into consideration that $a^2 + b^2 + c^2 = 2\left(p^2 - r^2 - 4Rr\right)$, we find

$$\overrightarrow{HA}\overrightarrow{HB} = \overrightarrow{HB}\overrightarrow{HC} = \overrightarrow{HA}\overrightarrow{HC} = 4R^2 + r^2 + 4Rr - p^2$$

Coming back to $HI^2$ we have

$$HI^2 = \frac{1}{4p^2}4R^2\left(a^2 + b^2 + c^2\right) - \left(a^4 + b^4 + c^4\right)\left(4R^2 + r^2 + 4Rr - p^2\right)\left(2ab + 2bc + 2ac\right)$$

But $ab + bc + ac = r^2 + p^2 + 4Rr$ and $16S^2 = 2a^2b^2 + 2b^2c^2 + 2c^2a^2 - a^4 - b^4 - c^4$

After some computation it results

$$IH^2 = 4R^2 + 4Rr + 3r^2 - p^2$$

The position vector of the center of the A-ex-inscribed circle is

$$\overrightarrow{PI_a} = \frac{1}{2(p-a)}\left(-a\overrightarrow{PA} + b\overrightarrow{PB} + c\overrightarrow{PC}\right)$$

We have $\overrightarrow{HI_a} = \frac{1}{2(p-a)}\left(-a\overrightarrow{HA} + b\overrightarrow{HB} + c\overrightarrow{HC}\right)$

We evaluate $\overrightarrow{HI_a}\overrightarrow{HI_a}$ and we have

$$HI_a^2 = \frac{1}{4(p-a)^2}\left(a^2HA^2 + b^2HB^2 + c^2HC^2 - 2ab\overrightarrow{HA}\overrightarrow{HB} - 2ac\overrightarrow{HA}\overrightarrow{HC} + 2bc\overrightarrow{HB}\overrightarrow{HC}\right)$$

We have

$$\overrightarrow{HA}\overrightarrow{HB} = \overrightarrow{HB}\overrightarrow{HC} = \overrightarrow{HA}\overrightarrow{HC} = 4R^2 - p^2 + r^2 - 4Rr$$

$$HI_a^2 = \frac{1}{4(p-a)^2}4R^2\left(a^2+b^2+c^2\right) - \left(a^4+b^4+c^4\right) + \left(4R^2 - p^2 + r^2 + 4Rr\right)(2bc - 2ab - 2ac)$$

From $2(p-a) = b+c-a$ it result that $4(p-a)^2 = a^2+b^2+c^2+2bc-2ab-2ac$, consequently

$$2bc - 2ab - 2ac = 4(p-a)^2 - \left(a^2+b^2+c^2\right)$$

Using

$$16S^2 = 2a^2b^2 + 2b^2c^2 + 2c^2a^2 - a^4 - b^4 - c^4$$

$$a^2b^2 + b^2c^2 + c^2a^2 = (ab+bc+ca)^2 - 2abc(a+b+c) = \left(r^2 + p^2 + 4Rr\right)^2 - 4pabc$$

After few computations we have

$$HI_a^2 = 4R^2 + 2ra^2 + r^2 + 4Rr - p^2$$

**Note**

In general, the formulae from this section can be deducted using the barycentric coordinates. See [10].

## D. Application

**Theorem** (Feuerbach)

In a triangle the circle of the nine points is tangent to the inscribed circle of the triangle and to the triangle's ex-inscribed circles.

**Proof**

We'll apply the median's theorem in the triangle $OIH$ and we have

$$4IO_9^2 = 2\left(OI^2 + IH^2\right) - OH^2$$

Because $OI^2 = R^2 - 2Rr$, $OH^2 = 9R^2 + 2r^2 + 8Rr - 2p^2$ and $IH^2 = 4R^2 + Rr + 3r^2 - p^2$

We obtain $IO_9 = \frac{R}{2} - r$, relation which shows that the circle of the none points (which has the radius $\frac{R}{2}$) is tangent to the inscribed circle.

We'll apply the median's theorem in the triangle $OI_aH$:

$$4I_aO_9^2 = 2\left(OI_a^2 + I_aH^2\right) - OH^2$$

Because $OI_a^2 = R^2 + 2Rr_a$ (Feuerbach) and $I_aH^2 = 4R^2 + 2r_a^2 + r^2 + 4Rr - p^2$ we obtain

$I_aO_9 = \frac{R}{2} + r_a$. This relation shows that the circle of the nine points and the S-ex-inscribed circle are exterior tangent. Similarly it can be show that the circle of the nine points and the B-ex-inscribed circles and C-ex-inscribed circles are tangent.

**Note**

In an article published in the G.M 4/1982, the Romanian professor Laurentiu Panaitopol proposed to find the strongest inequality of the type

$$R^2 + hr^2 \geq a^2 + b^2 + c^2$$

And proves that it is:

$$8R^2 + 4r^2 \geq a^2 + b^2 + c^2.$$

Taking into account $IH^2 = 4R^2 + 2r^2 - \frac{a^2 + b^2 + c^2}{2}$ and that $IH^2 \geq 0$ we find the inequality and its geometrical interpretation.

## 6.3. The point's power in rapport to a circle

**Theorem 1**

Let $P$ a point and $C(O,r)$ in a plane. Two lines constructed through $P$ intersect the circle in the points $A, B$ respectively $A', B'$; then takes place the equality:

$$\overrightarrow{PA}\cdot\overrightarrow{PB}=\overrightarrow{PA'}\cdot\overrightarrow{PB'}$$

**Proof**

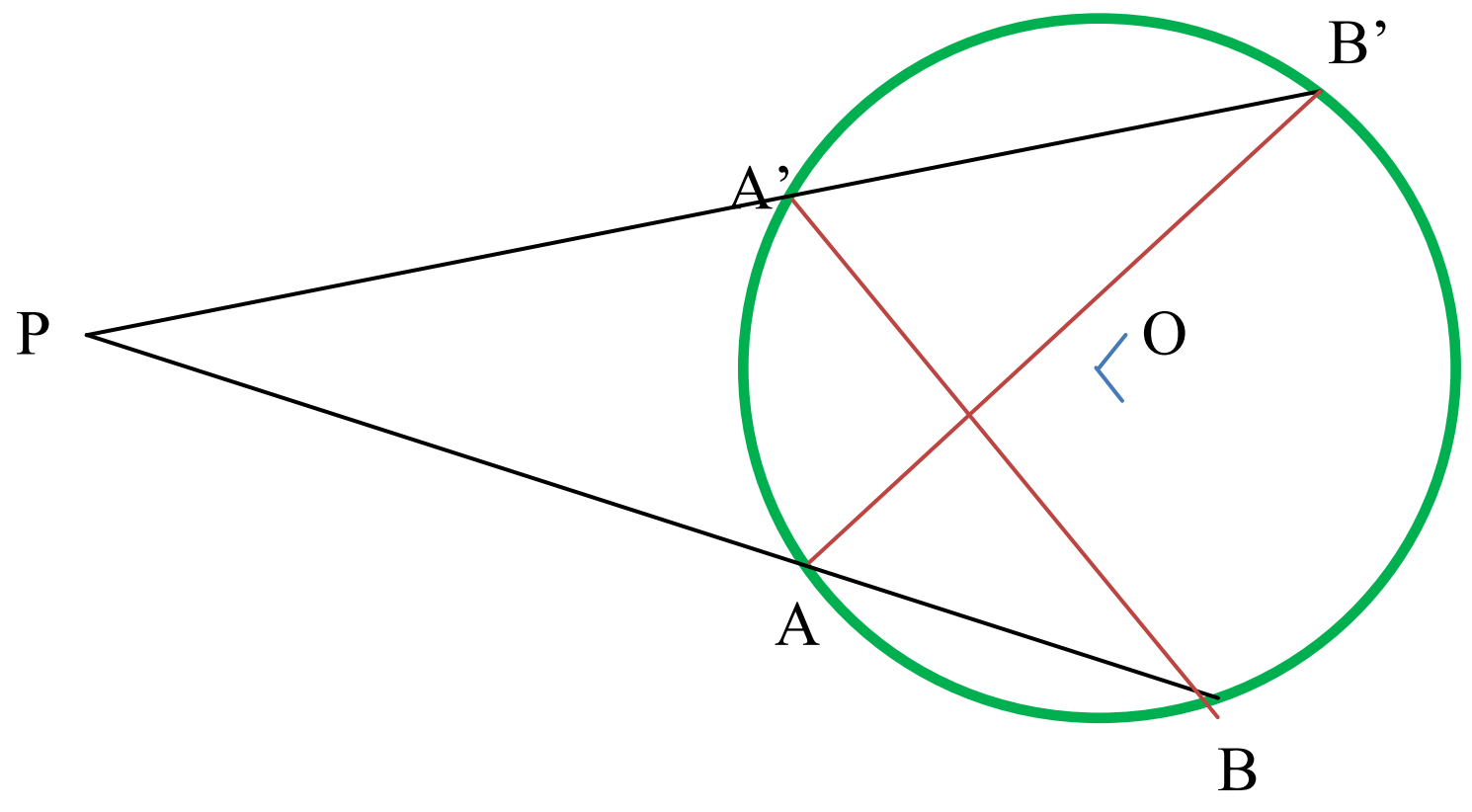

Fig.1

From the similarity of the triangles $PAB'$ and $PA'B$ (see Fig.1) it results

$$\frac{PA}{PA'}=\frac{PB'}{PB}$$

therefore

$$\overrightarrow{PA}\cdot\overrightarrow{PB}=\overrightarrow{PA'}\cdot\overrightarrow{PB'}.$$

Similarly it can be proved the relation from the hypothesis if the point $P$ is on the circle's interior . See Fig.2

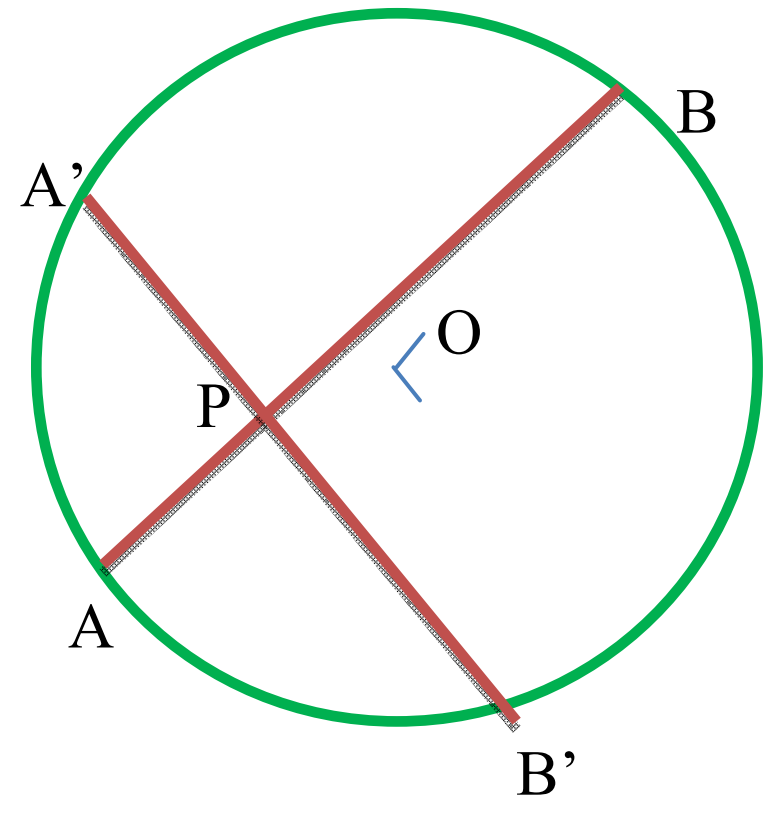

Fig. 2

**Corollary 1**

If a variable secant passes through a fix point $P$ and intersects the circle $C(O,r)$ in the points $A, B$, then the scalar product $\overrightarrow{PA}\cdot\overrightarrow{PB}$ is constant

**Definition 1**

The point's $P$ power in rapport to the circle $C(O,r)$ is the number $p(P)=\overrightarrow{PA}\overrightarrow{PB}$ where $A, B$ are the intersections with the circle of a secant that passes through $P$.

**Theorem 2**

The power of a point $P$ (such that $OP=d$ ) in rapport to the circle $C(O,r)$ is $p(P)=d^2-r^2$

**Proof**

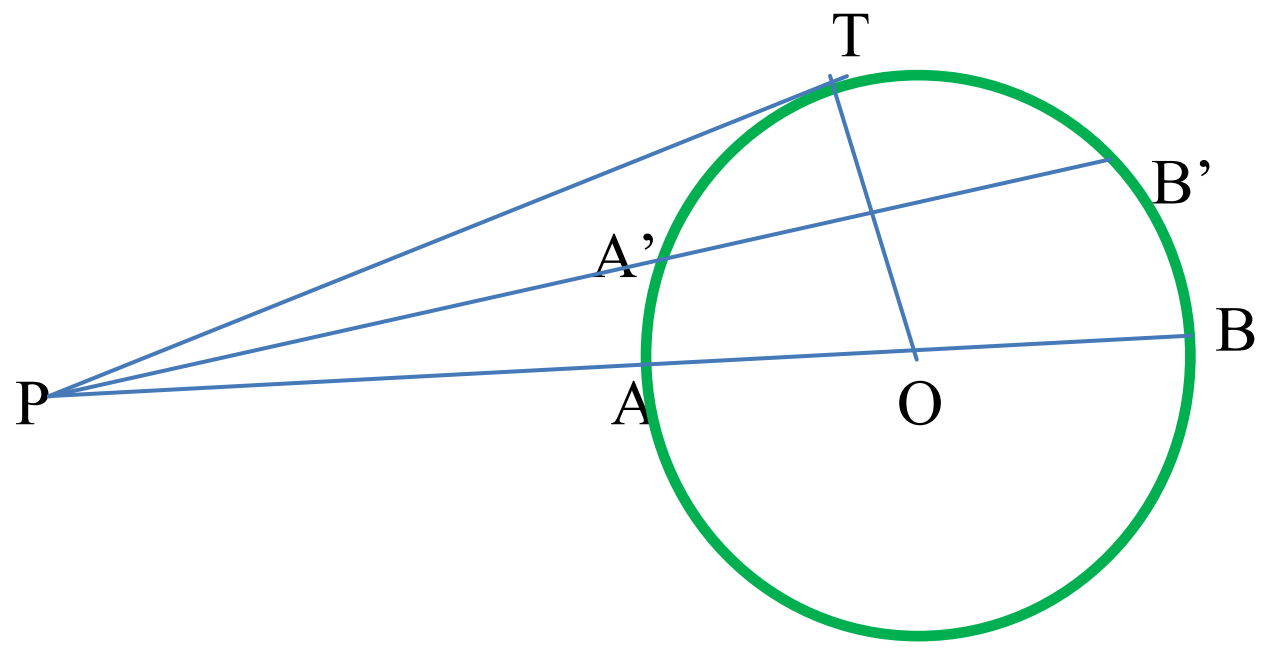

Fig. 3

If $OP=d>r$ (The point $P$ is external to the circle). See fig. 3. Then constructing the secant $PO$ we have $\overrightarrow{PA}\overrightarrow{PB}=PA\cdot PB=(d-r)(d+r)$, and we find $p(P)=d^2-r^2$

If $PO=d<r$ (see Fig. 4)

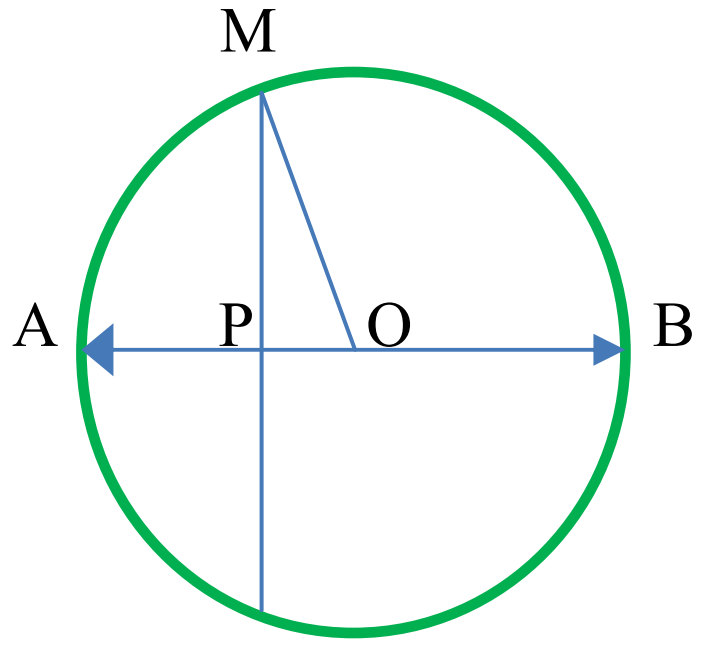

Fig. 4

Then $p(P)=\overrightarrow{PA}\overrightarrow{PB}=-PA\cdot PB=-(r-d)(r+d)=d^2-r^2$

**Remark 1**

a) If the point $P$ is in the exterior of the circle then $p(P)$ is positive number.

b) If the point $P$ is interior to the circle the point's $P$ power is a negative number (the vectors $\overrightarrow{PA}, \overrightarrow{PB}$ have opposite sense). If $p=0$ then $P(O)=-r^2$.

c) If the point $P$ is on the circle, its power is null (because one of the vectors $\overrightarrow{PA}$ or $\overrightarrow{PB}$ is null)
d) If the point $P$ is exterior to the circle $p(P) = \overrightarrow{PA}\overrightarrow{PB} = PT^2$ where $T$ is the tangency point with the circle of a tangent constructed from $P$. Indeed,

$$PT^2 = PO^2 - OT^2$$

e) If $P$ is interior to the circle $p(P) = \overrightarrow{PA}\overrightarrow{PB} = -PM^2$, where $M \in \mathcal{C}(O,r)$ such that $m\widehat{MPA} = 90°$. Indeed, $MP^2 = OM^2 - OP^2 = r^2 - d^2$.

**Note**

The name of the power of point in rapport to a circle was given by the mathematician Jacob Steiner (1796-1863), in 1832 and it is explained by the fact that in the definition of the power of a point appears the square of the length of a segment, and in Antiquity, the mathematician Hippocrates (sec V BC) used the expression "power of a segment" to define the square of a segment.

**Applications**

1. Determine the geometrical locus $d$ of the points from the plane of a given circle, which have in rapport with this circle a constant power.

**Solution**

Let $M$ a point with the property from hypothesis $p(M)=k$ (constant). But $p(M)=d^2-r^2$, we noted $p(M)=d^2-r^2$, $r$ the radius of the circle with the center in $O$. It results that $d^2=r^2+k$, therefore $=\sqrt{r^2+k}$ = const., therefore, the points with the given property are placed on a circle $C(O,d)$. Because it can be easily shown that any point on the circle $C(O,d)$ has the given property, it results that the geometrical locus is this circle. Depending of $k>0$ or $k<0$ the circle geometrical locus is on the exterior or in the interior of the given circle.

2. Determine the geometrical locus of the points on a plane, which have equal powers in rapport to two given circles. (the radical axis of two circles)

**Solution**

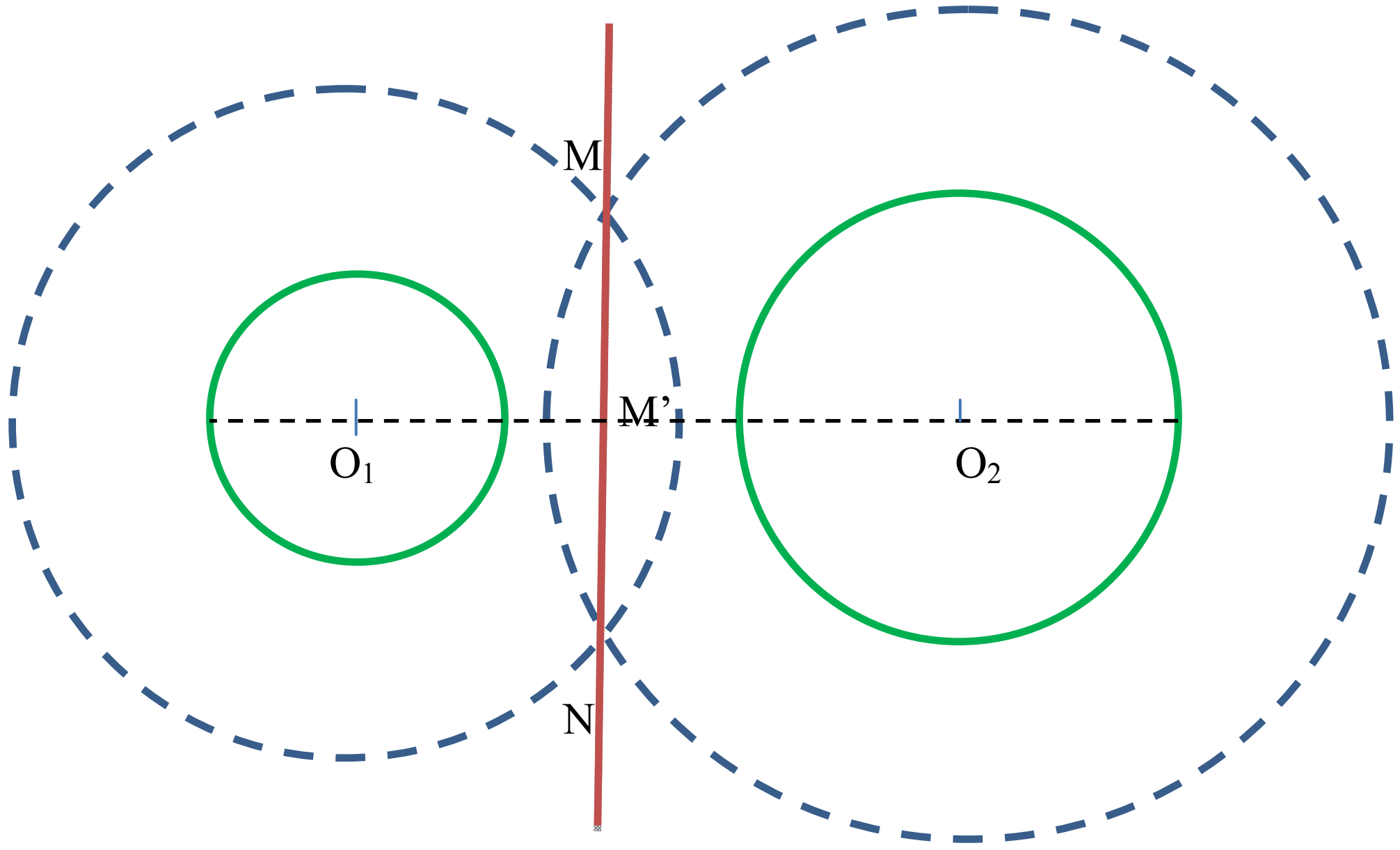

Fig. 5

Let's consider the circles $C(O_1,r_1),C(O_2,r_2)$, $O_1O_2>r_1+r_2$, see Fig. 5. If $M$ is a point such that its power in rapport to the two circles, $k_1,k_2$ are equal, then taking into consideration the result from the previous application, it result that $M$ is on the exterior of the given circles $C(O_1,d_1),C(O_2,d_2)$, where $d_1=\sqrt{k_1+r_1^2}$, $d_2=\sqrt{k_2+r_2^2}$. The point $M$ is on the exterior of the given circles (if we would suppose the contrary we would reach a contradiction with $k_1=k_2$).

In general the circles $C(O_1,d_1),C(O_2,d_2)$ have yet another common point $N$ and $O_1O_2$ is the perpendicular mediator of the segment $(MN)$, therefore $M,N$ belong to a perpendicular constructed on $O_1O_2$.

Let $M'$ the projection of $M$ on $O_1O_2$, $\{M'\}=MN\cap O_1O_2$, applying the Pythagoras' theorem in the triangles $MO_1M'$, $MO_2M'$ we'll find

$$d_1^2-d_2^2=O_1M'^2-O_2M'^2=r_1^2-r_2^2=const$$

From here it result that $M'$ is a fixed point, therefore the perpendicular from $M$ on $O_1O_2$ is a fixed line on which are placed the points with the property from the hypothesis

**Reciprocal**

Considering a point $P$ which belongs to the fixed perpendicular from above, we'll construct the tangents $PT_1$, $PT_2$ to the given circles; the fact that $P$ belongs to the respective fixed perpendicular is equivalent with the relation $PO_1^2-PO_2^2=r_1^2-r_2^2$. From the right triangles $PT_1O_1$, $PT_2O_2$, it results that $PO_1^2=PT_1^2+O_1T_1^2$; $PO_2^2=PT_2^2+O_2T_2^2$. We obtain that $PT_1=PT_2$, which shows that $P$ has equal powers in rapport to the given circles. Therefore the geometrical locus is the fixed perpendicular line on the centers' line of the given circles.

This line, geometrical locus, is called the radical axis of the two circles.

**Remark 2**

a) If the two given circles are conjugated ( $O_1\neq O_2$ ) then the radical axis is the perpendicular mediator of the segment $(O_1O_2)$.
b) If the circles from the hypothesis are interior, the radical axis is placed on the exterior of the given circles.
c) If the given circles are interior tangent or exterior tangent, the radical axis is the common tangent in the contact point.
d) If the circles are secant, the radical axis is the common secant.
e) If the circles are concentric, the geometrical locus is the null set.

3. Determine the geometrical locus of the points from plane, which have equal powers in rapport with three given circles. (The radical center of three circles).

**Solution**

Let $C(O_1,r_1),C(O_2,r_2)$, $C(O_3,r_3)$ three given circles and we'll consider them two by two external and such that their centers $O_1,O_2,O_3$ being non-collinear points (see Fig. 6). If $Q_3$ is the radical axis of the circles $C(O_1,r_1),C(O_2,r_2)$, and $Q_1$ is the radical axis of the circles $C(O_2,r_2),C(O_3,r_3)$, we note $\{R\}=Q_1\cap Q_3$ (this point exists, its non-existence would contradict the hypothesis that $O_1,O_2,O_3$ are non-collinear).

It is clear that the point $R$ having equal powers in rapport to the circles $C(O_1,r_1),C(O_2,r_2)$, it will have equal powers in rapport to the circles $C(O_2,r_2),C(O_3,r_3)$ and in rapport to $C(O_1, r_1), C(O_3, r_3)$, therefore it belongs to the radical axis $Q_2$ of these circles.

Consequently, $R$ is a point of the looked for geometrical locus.

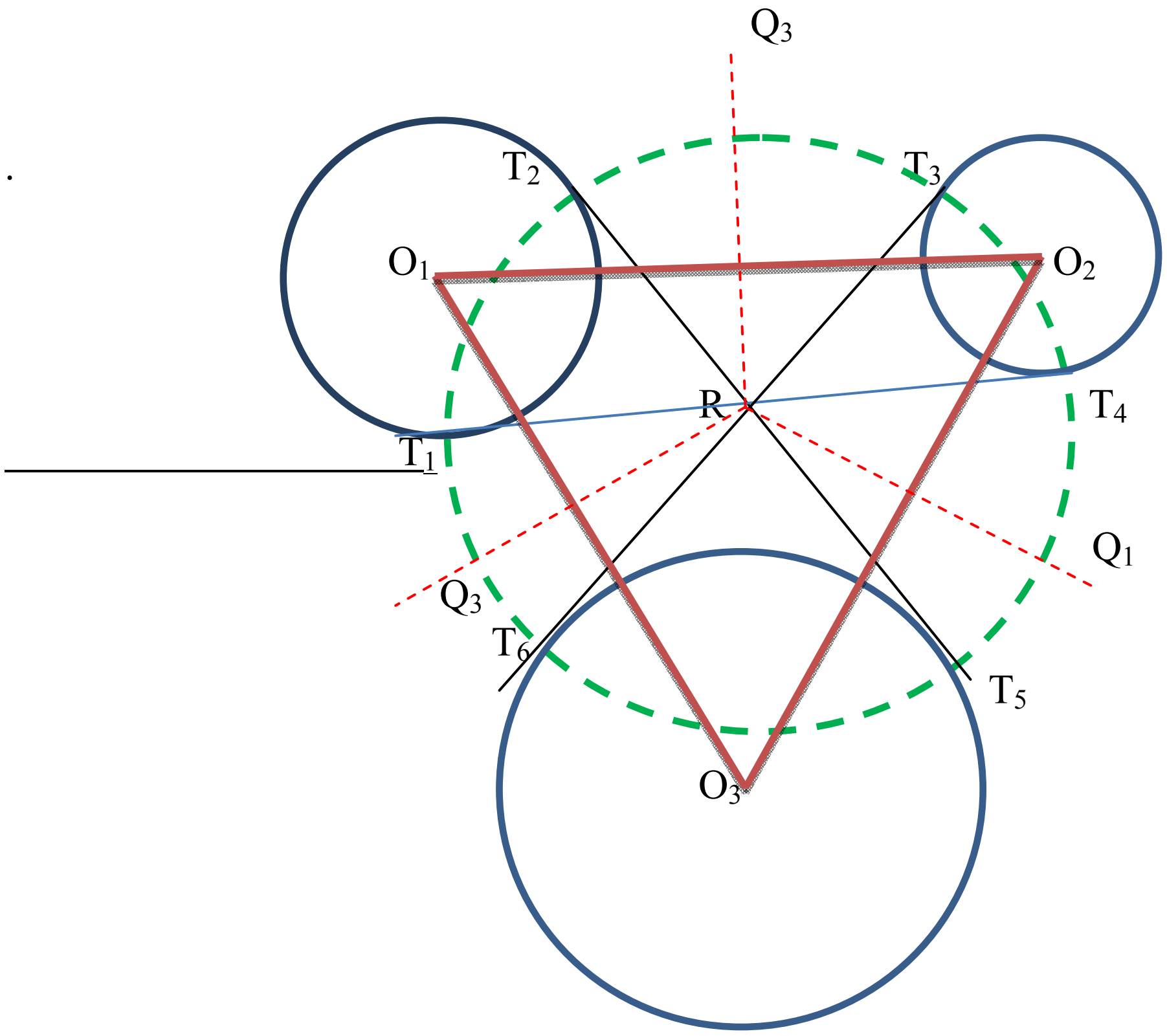

Fig. 6

It can be easily shown that $R$ is the unique point of the geometrical locus. This point – geometrical locus – is called the radical center of the given circles.

From the point $R$ we can construct the equal tangents $RT_1, RT_2$; $RT_3, RT_4$; $RT_5, RT_6$ to the given circles.

The circle with the center in $R$ is called the radical circle of the given circles.

**Remark 3**

a) If the centers of the given circles are collinear, then the geometrical locus, in general, is the null set, with the exception of the cases when the three circles pass through the same two common points. In this case, the geometrical locus is the common secant. When the circles are tangent two by two in the same point, the geometrical locus is the common tangent of the three circles.

b) If the circles' centers are distinct and non-collinear, then the geometrical locus is formed by one point, the radical center.

## 6.4. The non-harmonic rapport

**Definition 1**

If $A, B, C, D$ are four distinct points, in this order, on a line, we call their harmonic rapport, the result of the rapport in which the points $B$ and $D$ divide the segment $(AC)$.

We note $r = (ABCD) = \frac{BA}{BC} : \frac{DA}{DC}$

**Remark 1**

a) From the above definition it results that if $(ABCD) = (ABCD')$, then the points $D, D'$ coincide
b) If $r = -1$, the rapport is called harmonic.

**Definition 2**

If we consider a fascicle of four lines $a, b, c, d$ concurrent in a point $O$ and which determine on a given line $e$ the points $A, B, C, D$ such that, we can say about the fascicle formed of the four lines that it is a harmonic fascicle of vertex $O$, which we'll note it $O(ABCD)$. The lines $OA, OB, OC, OD$ are called the fascicle's rays.

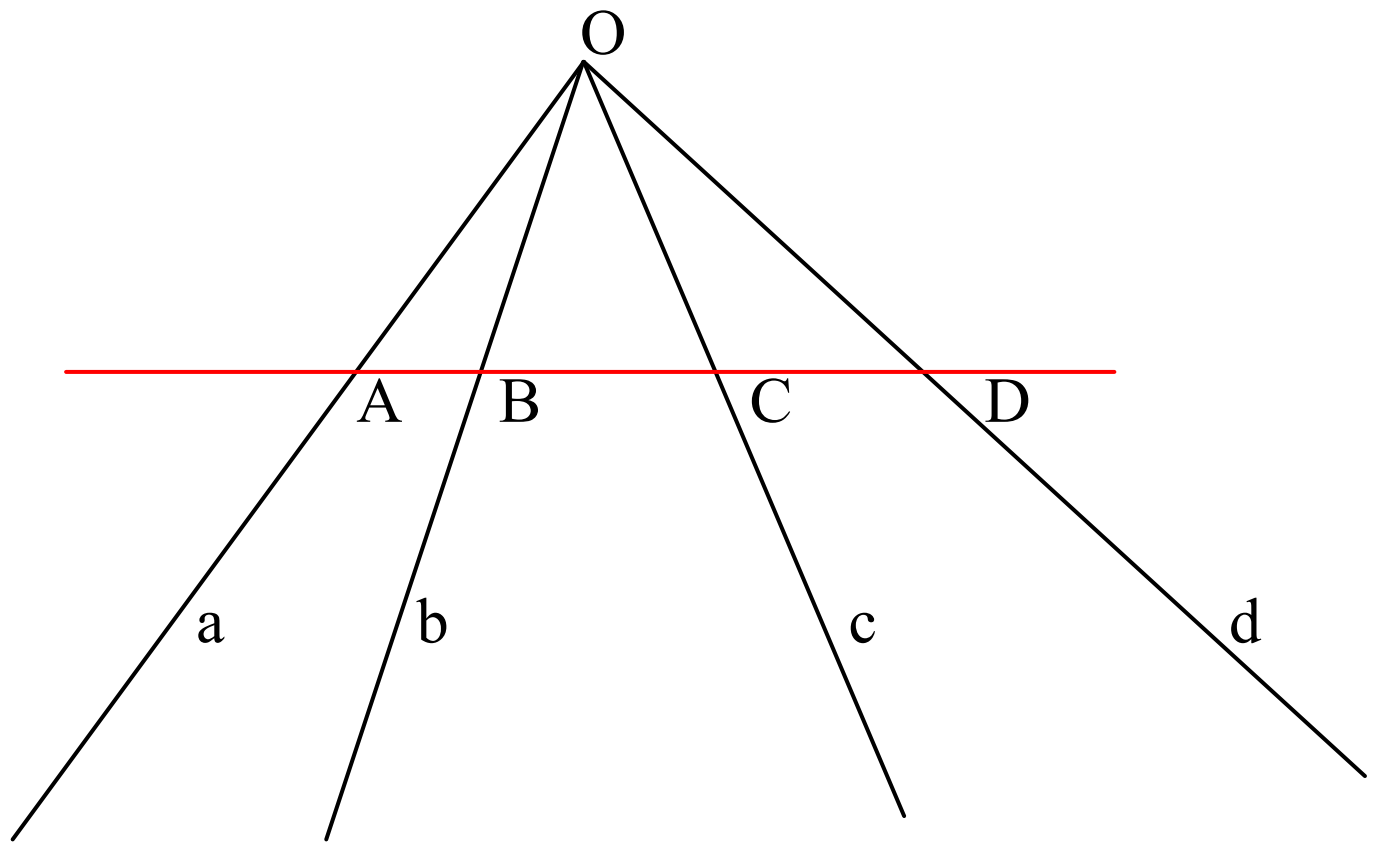

Fig. 1

**Property 1** (Pappus)

A fascicle of four rays determine on any secant a harmonic rapport constant

**Proof**

We'll construct through the point $C$ of the fascicle $O(ABCD)$ the parallel to $OA$ and we'll note $U$ and $V$ its intersections with the rays $OB, OD$ (see Fig. 2)

We have $(ABCD) = \frac{BA}{BC} : \frac{DA}{DC} = \frac{OA}{CU} : \frac{OA}{CV}$

Then $(ABCD) = \frac{CV}{CU}$

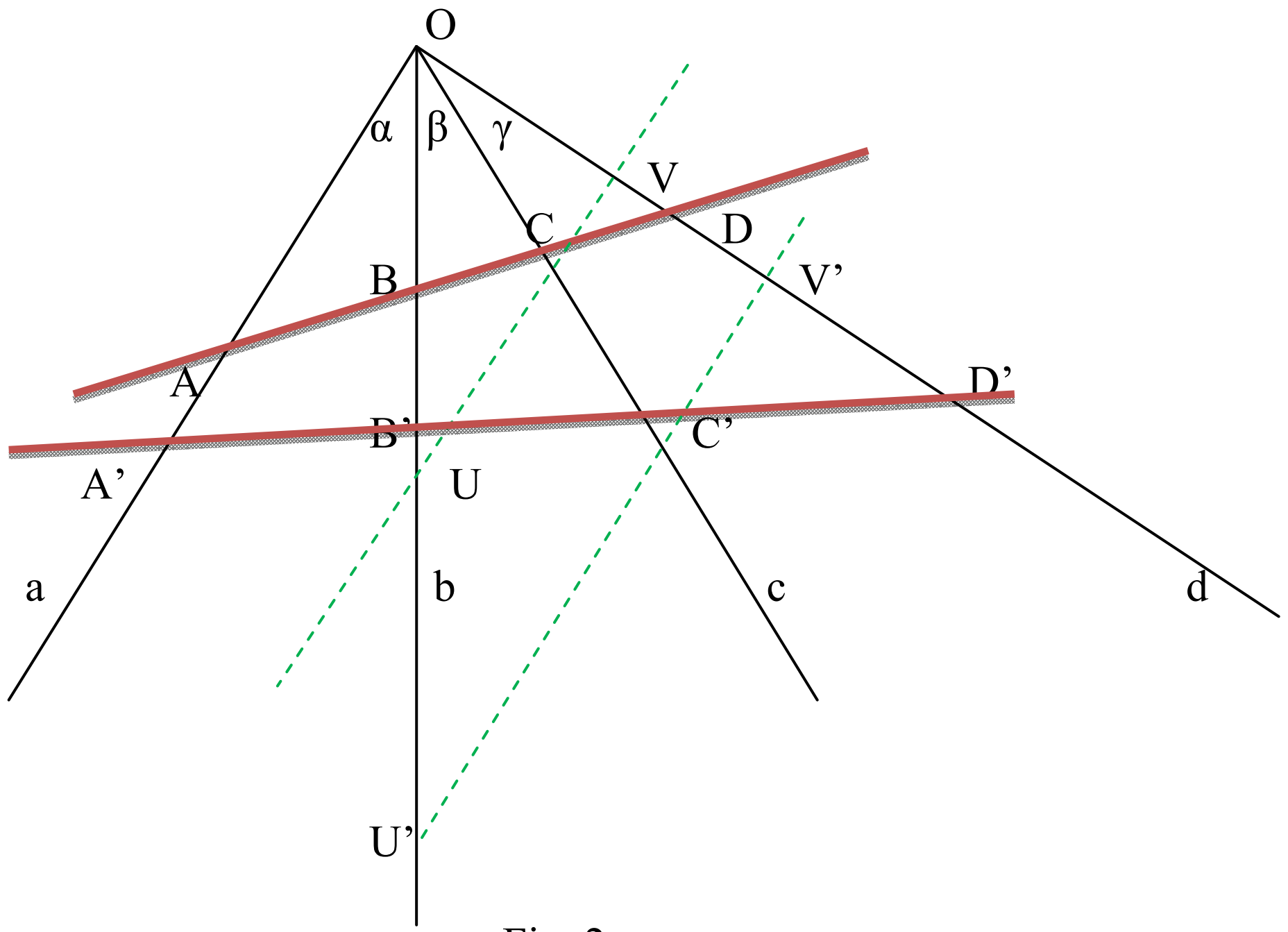

Fig. 2

Cutting the fascicle with another secant, we obtain the points $A'B'C'D'$, and similarly we obtain $(A'B'C'D') = \frac{C'V'}{C'U'} = \frac{CV}{CU} = (ABCD)$

**Property 2**

To fascicles whose rays make angles respectively equal have the same harmonic rapport

**Proof**

Let $O(ABCD)$ a fascicle in which we note

$$\alpha = m(\sphericalangle AOB),\ \beta = m(\sphericalangle BOC),\ \gamma = m(\sphericalangle COD)$$

(see Fig. 2)

We have

$$(ABCD) = \frac{BA}{BC} : \frac{DA}{DC} = \frac{S_{AOB}}{S_{BOC}} : \frac{S_{AOD}}{S_{COD}}$$

$$(ABCD) = \frac{AO \cdot BO \cdot \sin\alpha}{BO \cdot CO \cdot \sin\beta} : \frac{AO \cdot DO \cdot \sin(\alpha+\beta+\gamma)}{CO \cdot DO \cdot \sin\gamma} = \frac{\sin\alpha}{\sin\beta} : \frac{\sin(\alpha+\beta+\gamma)}{\sin\gamma}$$

Because the harmonic rapport is in function only of the rays' angles, it results the hypothesis' statement.

**Property 3**

If the rays of a fascicle $O(ABCD)$ are intersected by a circle that passes through $O$, respectively in the points $A_1, B_1, C_1, D_1$, then

$$(ABCD) = \frac{B_1A_1}{B_1C_1} : \frac{D_1A_1}{D_1C_1}$$

**Proof**

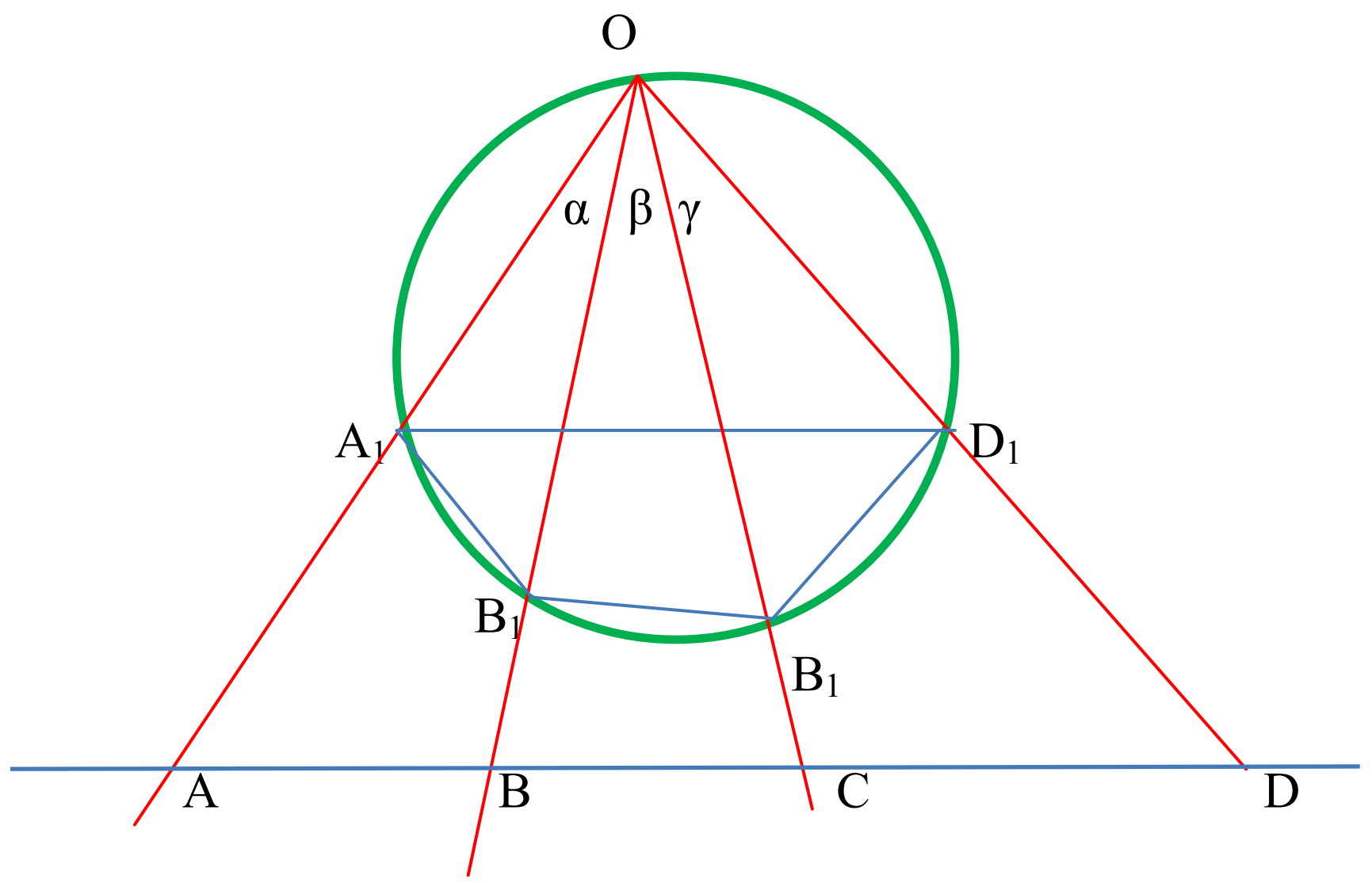

Fig. 3

From the sinus' theorem we have $A_1B_1 = 2R\sin\alpha$, $B_1C_1 = 2R\sin\beta$, $C_1D_1 = 2R\sin\gamma$, $D_1A_1 = 2R\sin(\alpha+\beta+\gamma)$. On the other side $(ABCD) = \frac{\sin\alpha}{\sin\beta} : \frac{\sin(\alpha+\beta+\gamma)}{\sin\gamma}$.

We'll note $(ABCD) = (A_1B_1C_1D_1)$.

**Remark 2**

From what we proved before it results:

If $A, B, C, D$ are four point on a circle and $M$ is a mobile point on the circle, then the harmonic rapport of the fascicle $M(ABCD)$ is constant.

**Theorem 1**

If two harmonic division have a homology point common, then the lines determined by the rest of the homological points are concurrent.

**Proof**

Let the harmonic divisions $(ABCD)$, $(AB'C'D')$ with the homological point A common.

We'll note $BB' \cap CC' = \{O\}$. We'll consider the fascicle $O(ABCD)$, and we'll note $OD \cap d' = \{D_1\}$, where $d'$ is the line of the points $(AB'C'D')$ (see Fig. 4)).

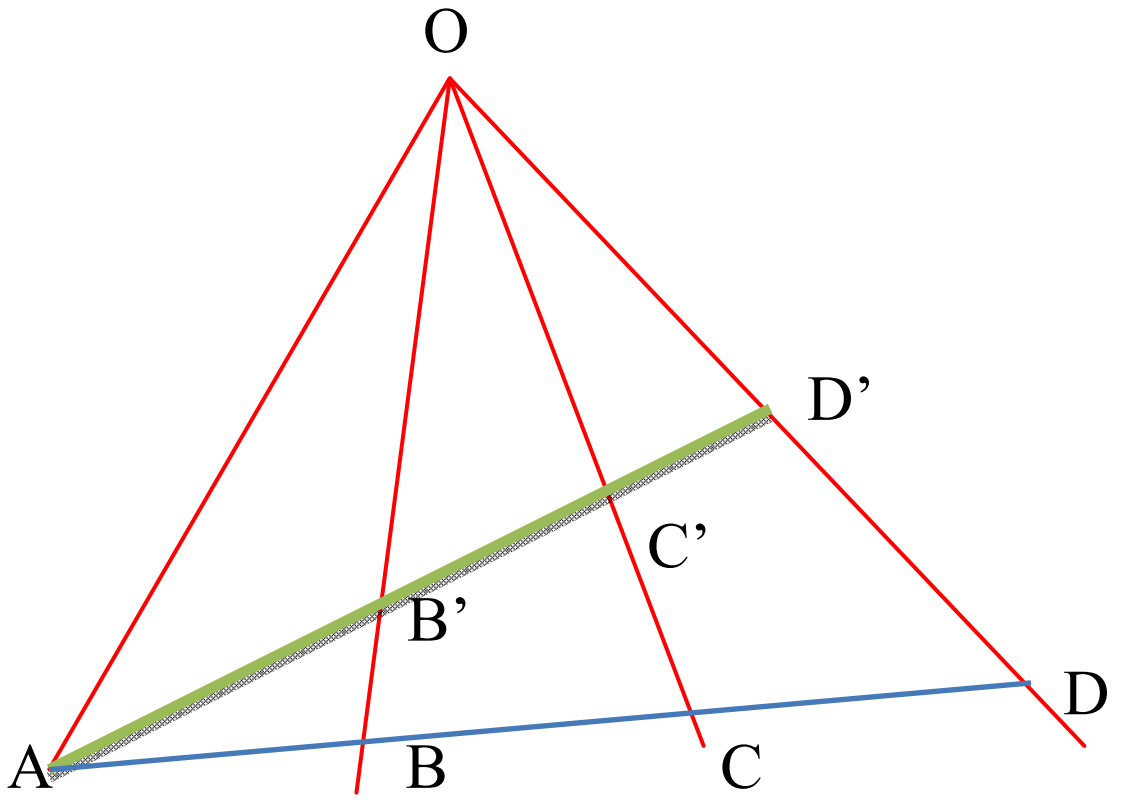

Fig. 4

The fascicle $O(ABCD)$ being intersected by $d'$, we have $(ABCD) = (AB'C'D_1)$.

On the other side, we have $(AB'C'D_1) = (AB'C'D')$, therefore $D_1 = D'$ thus the lines $BB', CC', DD'$ are concurrent in the point $O$.

**Theorem 2**

If two fascicles have the same non-harmonic rapport, and a common homological ray, then the rest of the rays intersect in collinear points.

**Proof**

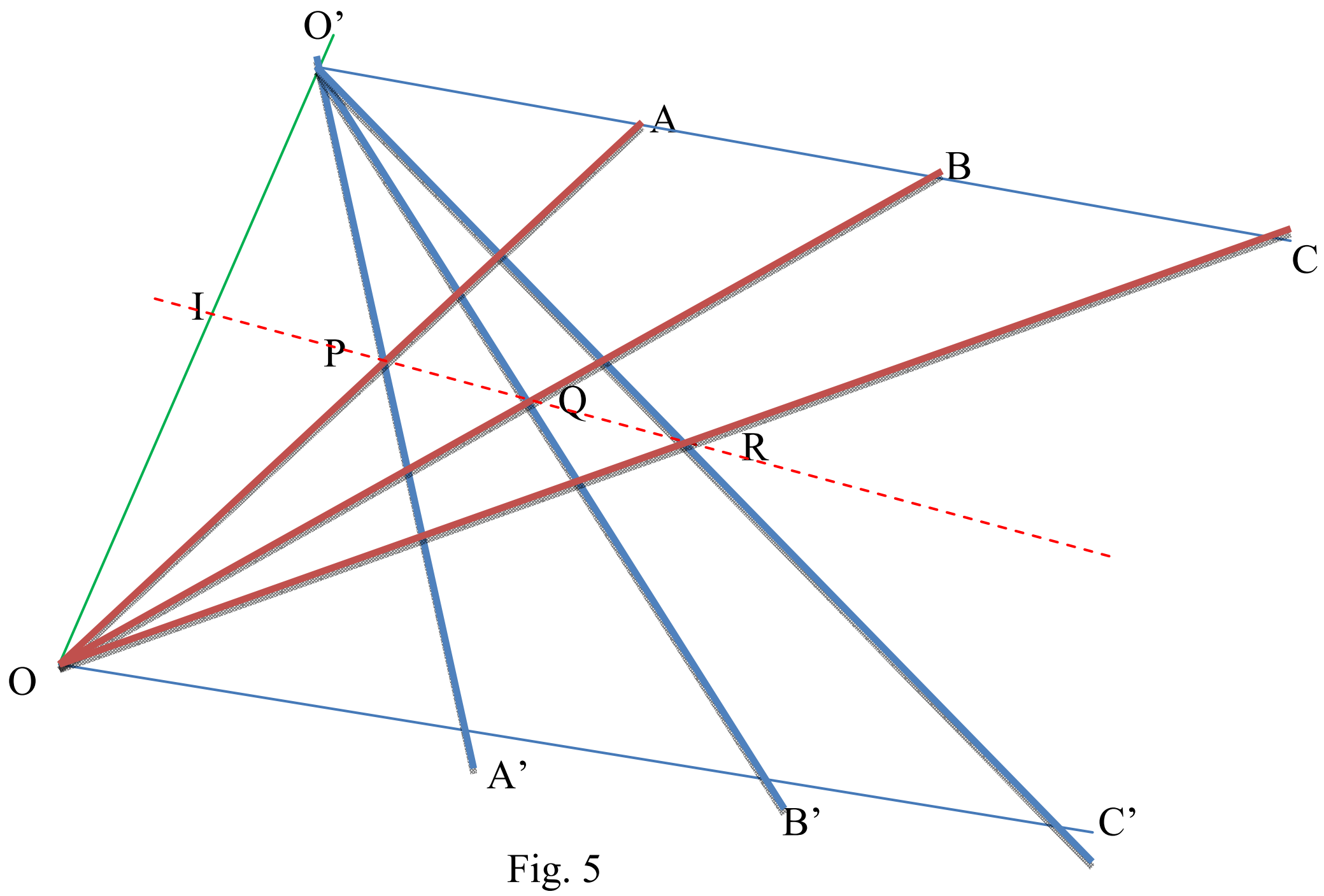

Fig. 5

Le the fascicles $O(O'ABC), O'(OAB'C')$ that have the common homological ray $OO'$ and $(O'ABC) = (OAB'C')$ (see Fig. 5).

Let

$$\{P\} = OA \cap O'A'; \{Q\} = OB \cap O'B'; \{I\} = OO' \cap PQ; \{R'\} = PQ \cap O'C'$$

We have

$$(IPQR) = (IPQR') = (O'ABC) = (OAB'C'),$$

therefore $R' = R$.

We obtain that $\{R\} = OC \cap O'C'$, therefore the homological rays intersect in the points $P, Q, R$.

**Applications**

1.

If a hexagon $ABCDEF$ is inscribed in a circle, its opposing sites intersect in collinear points (B. Pascal – 1639)

**Proof**

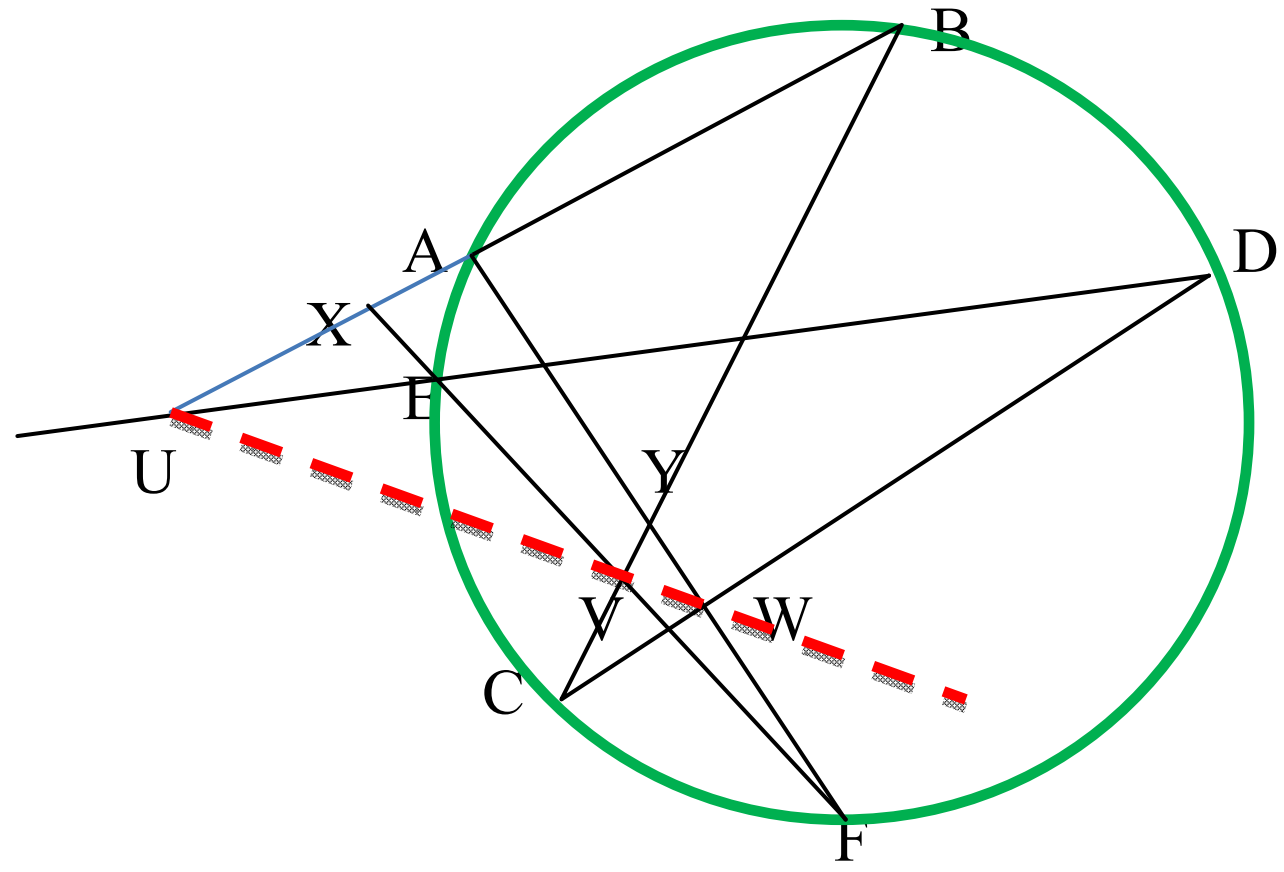

Fig. 6

Let

$$\{U\} = AB \cap DE; \{V\} = BC \cap EF; \{W\} = CD \cap FA$$ (see Fig. 6.

We saw that the harmonic rapport of four points on a circle is constant when the fascicle's vertex is mobile on a circle, therefore we have: $E(ABCDF) = C(ABDF)$.

We'll cut these equal fascicles, respectively with the secants $AB, AF$.

It results $(ABUX) = (AYWF)$.

These two non-harmonic rapports have the homological point $A$ in common. Thus that the lines $BY, UW, XF$ are concurrent, therefore $\{V\} = BY \cap XF$ belongs to the line $UW$, consequently $U, V, W$ are concurrent and the theorem is proved.

2. If the triangles $ABC$, are placed such that the lines $AA', BB', CC'$ are concurrent, then their homological sides intersect in three collinear points. (G. Desargues – 1636)

**Proof**

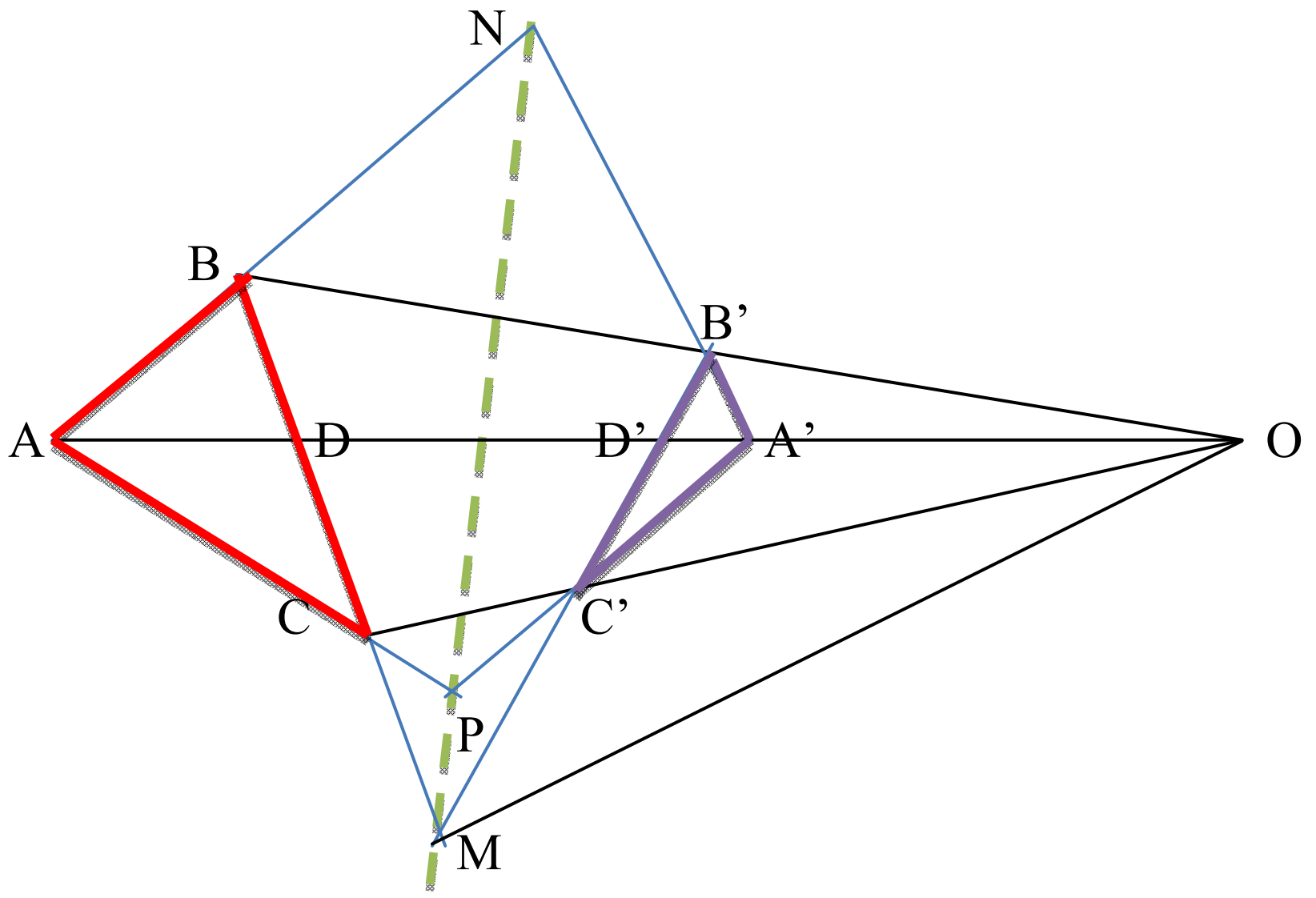

Fig. 7

Let $O$ be the intersection of the lines $AA', BB', CC'$ and

$$\{M\} = BC \cap BC'; \{N\} = AB \cap A'B'; \{P\} = AC \cap A'C'$$

We'll consider the fascicle $(OBACM)$ and we'll cut with the secants $BC, B'C'$, then we have

$$(BDCM) = (B'D'C'M)$$

We noted $\{D\} = OA \cap BC$ and $\{D'\} = OA' \cap B'C'$.

Considering the fascicles $A(BDCM) = A'(B'D'C'M)$ we observe that have the homologue ray $AA'$ in common, it results that the homological rays $AB, A'B'$; $AC, A'B'$; $AM, A'M$ intersect in the collinear points $N, P, M$.

3.

In a circle we consider a point $M$ in its interior.

Through $M$ we construct the cords $(PQ)$, $(RS)$, $(KL)$ such that $(MK) \equiv (ML)$.

If $MK \cap PS = \{U\}$; $ML \cap RQ = \{V\}$, then $MU = MV$ (the butterfly problem).

**Proof**

Considering the points $K, S, Q, L$, we have

$$P(KSQL) = R(KSQL)$$

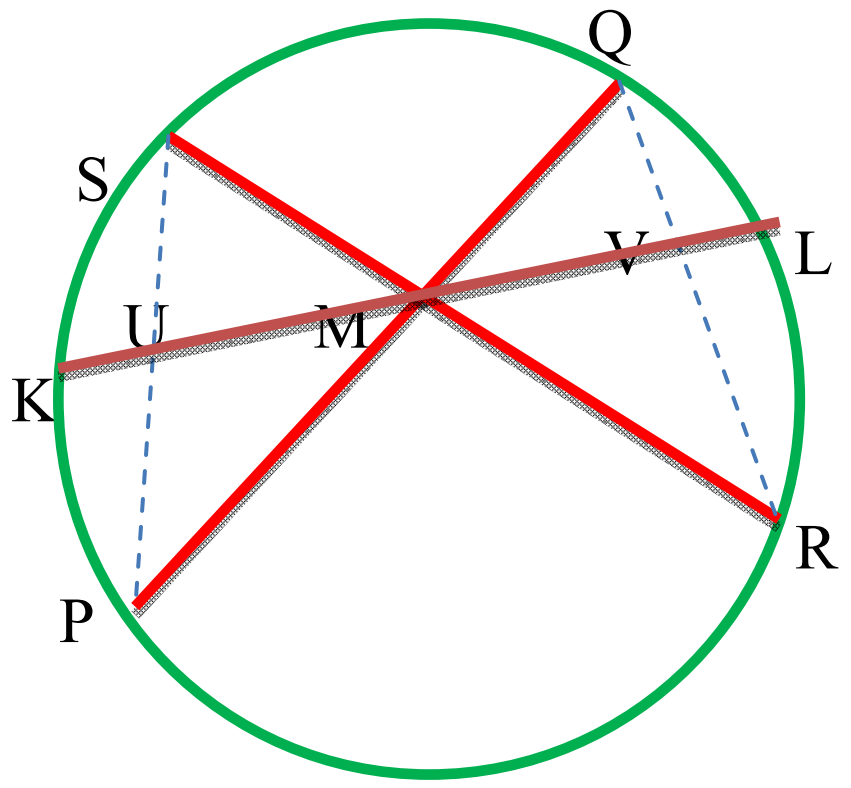

Fig.8

(The harmonic rapport of the fascicles being $\frac{SK}{SQ} : \frac{LK}{LQ}$).

Intersecting these fascicles with the line $KL$, we have $(KULM) = (KMVL)$.

Therefore

$$\frac{UK}{UM} : \frac{LK}{LM} = \frac{MK}{MV} : \frac{LK}{LM}.$$

Taking into account hat $(MK) \equiv (ML)$ it results $\frac{UK}{UM} = \frac{LV}{MV}$

and from here

$$MU = MV$$

### 6.5. Pole and polar in rapport with a circle

**Definition 1**

Let circle $C(O,R)$, a point $P$ in its plane, $P \neq O$, and a point $P'$ such that

$$\overrightarrow{OP} \cdot \overrightarrow{OP'} = R^2 \qquad (1)$$

We call the perpendicular $p$ constructed in the point $P'$ on the line $OP$, the polar of the point $P$ in rapport with the circle, and about the point $P$ we say that it is the pole of the line $p$

**Remark 1**

a) If $P$ belongs to the circle, its polar is the tangent in $P$ to $C(O,R)$.

Indeed, the relation (1) leads to $P' \neq P$

b) If $P$ is in the interior of the circle, its polar $p$ is an external line of the circle.

c) If $P$ and $Q$ are two points such that $m(\sphericalangle POQ) = 90°$, and $p$, $q$ are their polar, then from the definition 1 it results that $p \perp q$.

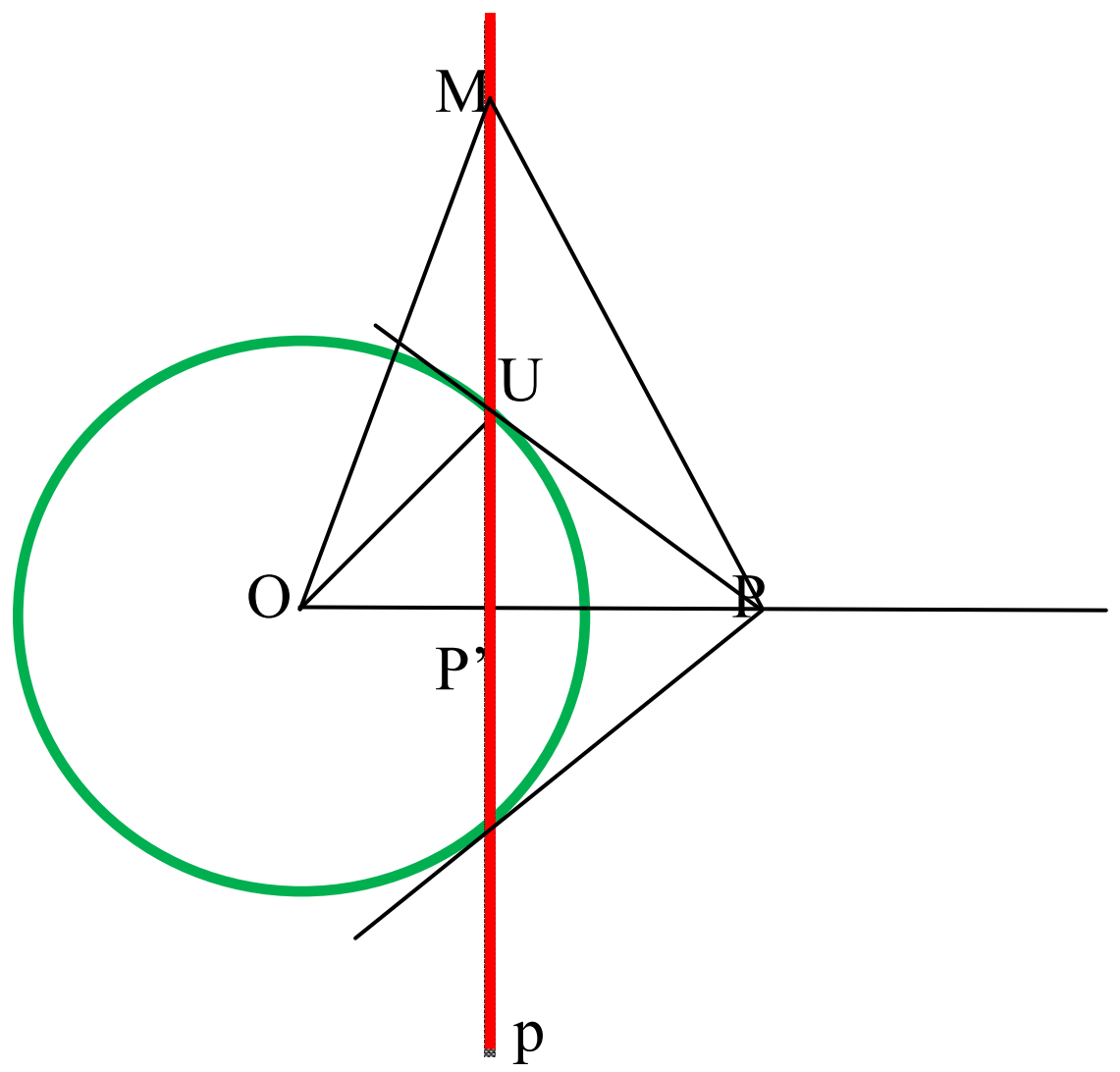

Fig. 1

**Proposition 1**

If the point $P$ is in the exterior of a circle, its polar is determined by the contact points with the circle of the tangents constructed from $P$ to the circle. (see Fig. 1)

**Proof**

Let $U,V$ the contact points of the tangents from the point $P$ to $C(O,R)$. In the right triangle $OUP$, if we note $P''$ the orthogonal projection of $U$ on $OP$, we have $OU^2 = OP'' \cdot OP$ (the right triangle's side theorem ). But from (1) we have that $OO \cdot OP' = R^2$, it results that $P'' = P'$, and therefore $U$ belongs to the polar of the point $P$. Similarly, $V$ belongs to the polar, therefore $UV$ is the polar of the point $P$.

**Theorem 1** (the polar characterization)

The point $M$ belongs to the polar of the point $P$ in rapport to the circle $C(O,R)$, if and only if

$$MO^2 - MP^2 = 2R^2 - OP^2 \qquad (2)$$

**Proof**

If $M$ is an arbitrary point on the polar of the point $P$ in rapport to circle $C(O,R)$, then $MP' \perp OP$ (see Fig. 1) and

$$MO^2 - MP^2 = \left(P'O^2 + P'M^2\right) - \left(P'P^2 + P'M^2\right) = P'O^2 - P'P^2 =$$

$$= OU^2 - P'U^2 + P'U^2 - PU^2 = R^2 - \left(OP^2 - R^2\right) = 2R^2 - OP^2.$$

Reciprocally, if $M$ is in the plane of the circle, such that the relation (2) is true, we'll note $M'$ the projection of $M$ on $OP$; then we have

$$M'O^2 - M'P^2 = \left(MO^2 - M'M^2\right) - \left(MP^2 - M'M^2\right) = MO^2 - MP^2 = 2R^2 - OP^2 \qquad (3)$$

On the other side

$$P'O^2 - P'P^2 = 2R^2 - OP^2 \qquad (4)$$

From (3) and (4) it result that $M' = P'$, therefore $M$ belongs to the polar of the point $P$.

**Theorem 2** (Philippe de la Hire)

If $P, Q, R$ are points which don't belong to a circle, then

1) $P \in q$ if and only if $R \in p$

(If a point belongs to the polar of another point, then the second point also belongs to the polar of the first point in rapport with the same circle.)

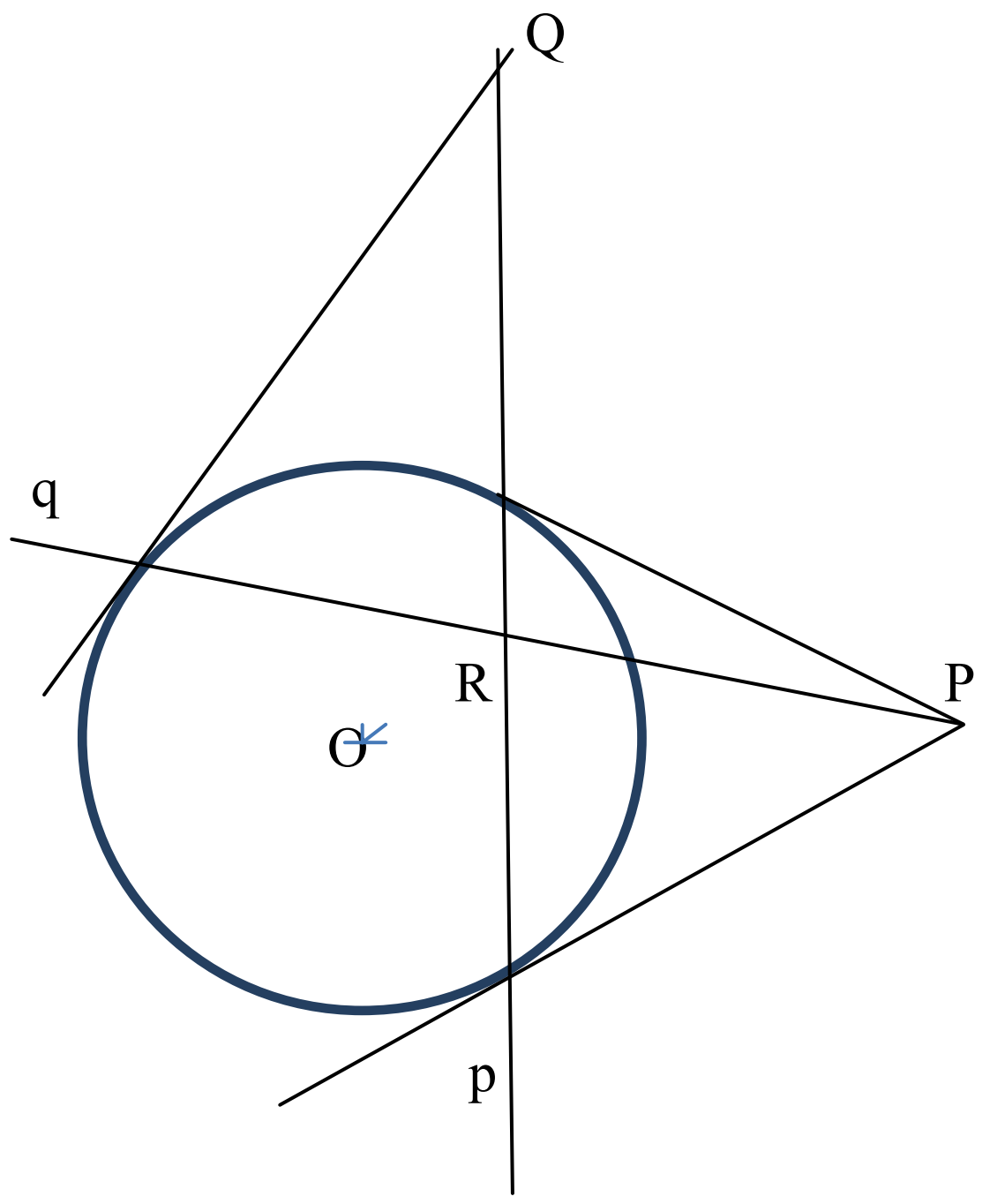

Fig. 2

2) $r = PQ \Leftrightarrow R \in p \cap q$

(The pole of a line that passes through two points is the intersection of the polar of the two points.)

**Proof**

1) From theorem 1 we have $P \in q \Leftrightarrow PO^2 - PQ^2 = 2R^2 - OQ^2$. Therefore

$$QO^2 - OP^2 = 2R^2 - OP^2 \Leftrightarrow Q \in p$$

2) Let $R \in p \cap q$; from 1) it results $P \in r$ and $Q \in r$, therefore $r = PQ$.

**Observation 1**

From theorem 2 we retain:

a) The polar of a point, which is the intersection of two given lines in rapport with a circle, is the line determined by the poles of those lines.
b) The poles of concurrent lines are collinear points and reciprocally, the polar of collinear points are concurrent lines.

**Transformation by duality in rapport with a circle**

The duality in rapport with circle $C(O,R)$ is a geometric transformation which associates to any point $P \neq O$ from plane its polar, and which associates to a line from the plane it pole.

By duality we, practically, swap the lines' and points' role; Therefore to figure $\mathcal{F}$ formed of points and lines, through duality corresponds a new figure $\mathcal{F}'$ formed by the lines(the polar of the points from figure $\mathcal{F}$) and from points (the poles of the lines of figure $\mathcal{F}$) in rapport with a given circle..

The duality has been introduced in 1823 by the French mathematician Victor Poncelet.

When the figure $\mathcal{F}$ is formed from points, lines and, eventually a circle, and if these belong to a theorem T, transforming it through duality in rapport with the circle, we will still maintain the elementary geometry environment, and we obtain a new figure $\mathcal{F}$, to which is associated a new theorem T', which does not need to be proved.

From the proved theorems we retain:

- If a point is situated on a line, through duality to it corresponds its polar, which passes through the line's pole in rapport with the circle.
- To the line determined by two points correspond, by duality in rapport with a circle, the intersection point of the polar of the two points.
- To the intersection point of two lines correspond, by duality in rapport with a circle, the line determined by the poles of these lines.

**Observation 2**

The transformation by duality in rapport with a circle is also called the transformation by reciprocal polar.

**Definition 2**

Two points that belong each to the polar of the other one in rapport with a given circle, are called conjugated points in rapport with the circle, and two lines each passing through the other's pole are called conjugated lines in rapport with the circle.

**Definition 3**

If through a point $P$ exterior to the circle $C(O,R)$ we construct a secant which intersects the circle in the points $M,N$, and the point $Q$ is on this secant such that $\frac{PM}{PN}=\frac{QM}{QN}$, we say about the points $P,Q$ that they are harmonically conjugated in rapport to the circle $(O)$ (see Fig. 3)

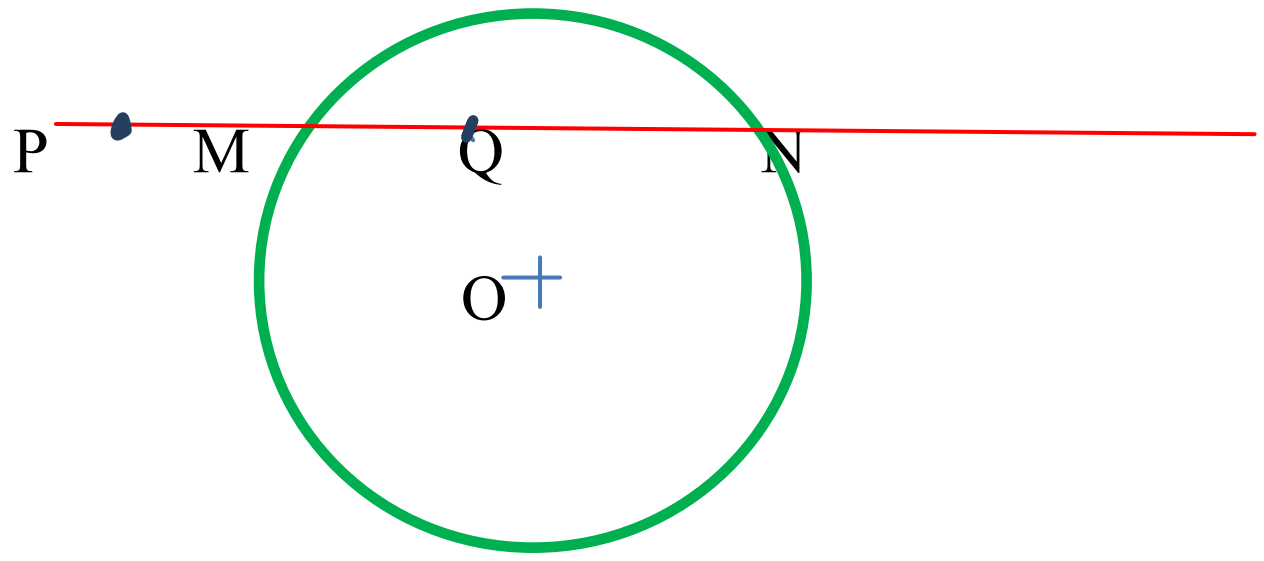

Fig. 3

We say that the points $P,M,Q,N$ form a harmonic division.

**Theorem 3**

If a line that passes through two conjugated points in rapport with a circle is the secant of the circle, then the points are harmonic conjugated in rapport with the circle.

**Proof**

Let $P,Q$ The conjugated points in rapport to the circle $C(O,R)$ and $M,N$ the intersections of the secant $PQ$ with the circle (see Fig. 4), and the circle circumscribed to triangle $OMN$. The triangles $OMP$ and $OP'M$ are similar ( $\sphericalangle MOP \equiv \sphericalangle P'OM$ and $\sphericalangle MOP \equiv \sphericalangle OP'M$ having equal supplements ).

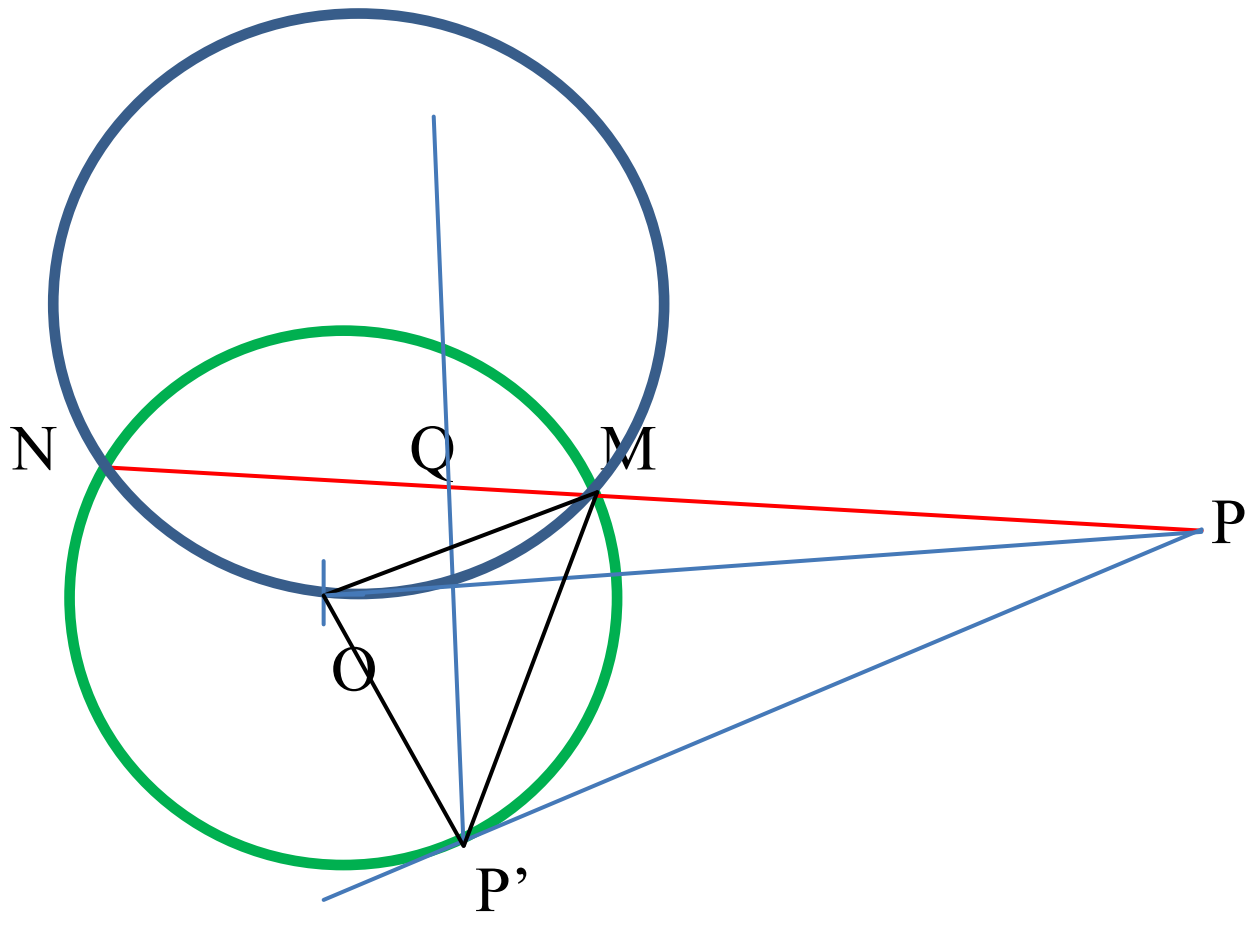

Fig. 4

It results $\frac{OM}{OP'} = \frac{OP}{OM}$, therefore $OP \cdot OP' = R^2$ which shows that $P'$ belongs to the polar of $P$, in rapport to the circle $C(O,R)$.

Because the point $Q$ also belongs to the polar, it results that $QP'$ is the polar of $P$, therefore $QP' \perp OP$. Because $\sphericalangle OMN \equiv \sphericalangle OP'M$ and $\sphericalangle ONM \equiv \sphericalangle MP'P$ it result that $PP'$ is the exterior bisector of the triangle $MP'N$. Having $QP' \perp P'P$, we obtain that $P'Q$ is the interior bisector in the triangle $MP'N$. The bisector's theorem (interior and exterior) leads to $\frac{QM}{QN} = \frac{PM}{PN}$, therefore the points $P, Q$ are harmonically conjugate in rapport with the circle ($O$)

**Observation 3**

a) It can be proved that also the reciprocal of the previous theorem, if two points are harmonically conjugate in rapport with a circle, then any of these points belongs to the polar of the other in rapport with the circle.

b) A corollary of the previous theorem is: the geometrical locus of the harmonically conjugate of a point in rapport with a given circle is included in the polar of the point in rapport with the given circle.

**Theorem 4**

If $ABCD$ is a quadrilateral inscribed in the circle $(O)$ and $\{P\} = AB \cap CD$, $\{Q\} = BC \cap AD$, $\{R\} = AC \cap BD$ then

a) The polar of $P$ in rapport with the circle $(O)$ is $QR$;

b) The polar of $Q$ in rapport with the circle $(O)$ is $PR$;

c) The polar of $R$ in rapport with the circle $(O)$ is $PQ$.

**Proof**

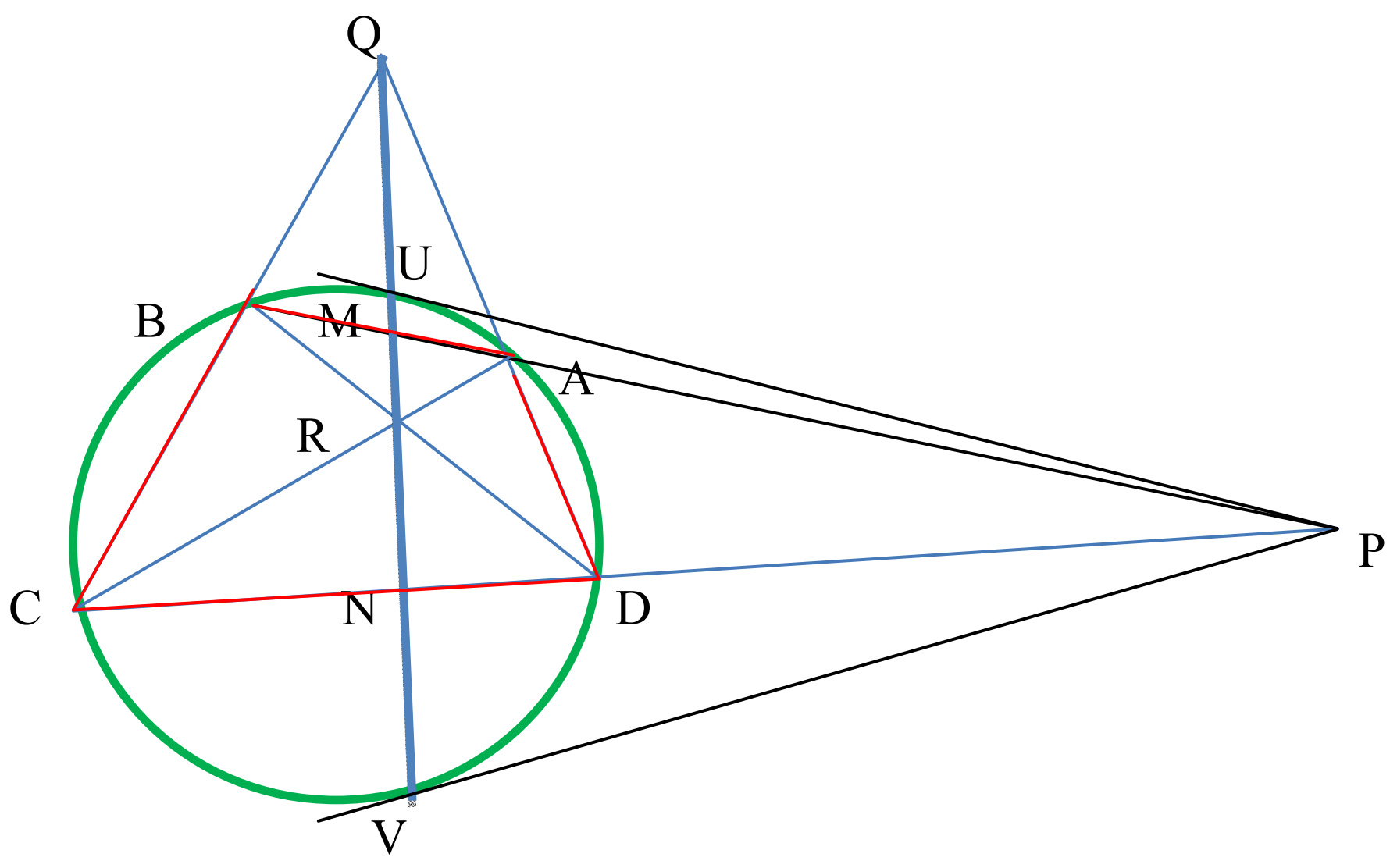

Fig. 5

It is sufficient to prove that $\frac{MA}{MB}=\frac{PA}{PB}$ and $\frac{ND}{NC}=\frac{PD}{PC}$ where $M,N$ are the intersections of the line $QR$ with $(AB)$, $(CD)$ respectively. (See Fig. 5).

We have

$$\frac{MA}{MB}=\frac{QA}{QB}\cdot\frac{\sin MQA}{\sin MQB} \quad (1)$$

$$\frac{RB}{RD}=\frac{QB}{QD}\cdot\frac{\sin BQR}{\sin DQR} \quad (2)$$

$$\frac{QA}{QD}=\frac{AC}{CD}\cdot\frac{\sin QCA}{\sin QCD} \quad (3)$$

$$\frac{RD}{RB}=\frac{CD}{CB}\cdot\frac{\sin ACD}{\sin RCB} \quad (4)$$

Multiplying side by side the relations (1), (2), (3), (4) and simplifications, we obtain

$$\frac{MA}{MB}=\frac{AC}{BC}\cdot\frac{\sin ACD}{\sin QCD} \quad (5)$$

On the other side

$$\frac{PA}{PB}=\frac{AC}{BC}\cdot\frac{\sin ACD}{\sin QCD} \quad (6)$$

From (5) and (6) we obtain

$$\frac{MA}{MB}=\frac{PA}{PB}$$

therefore, $M$ is the harmonic conjugate of the point $P$ in rapport with the circle. Similarly, we prove that $N$ is the harmonic conjugate of the point $P$ in rapport with the circle, therefore the polar of $P$ is $QR$.

Similarly we prove b) and c).

**Definition 4**

Two triangles are called reciprocal polar in rapport with a circle if the sides of one of the triangle are the polar of the vertexes of the other triangle in rapport with the circle. A triangle is called auto conjugate in rapport with a circle if its vertexes are the poles of the opposite sides.

**Observation 4**

Theorem 4 shows that the triangle $PQR$ is auto conjugate.

**Applications**

1) If $ABCD$ is a quadrilateral inscribed in a circle of center $O$ and $AB \cap CD = \{P\}$, $AC \cap BD = \{R\}$, then the orthocenter of the triangle $PQR$ is the center $O$ of the quadrilateral circumscribed circle.

**Proof**

From the precedent theorem and from the fact that the polar of a point is perpendicular on the line determined by the center of the circle and that point, we have that $OP \perp QR$, $OR \perp PQ$, which shows that $O$ is the triangle's $PQR$ orthocenter.

**Observation 5**

This theorem can be formulated also as follows: The orthocenter of a auto conjugate triangle in rapport with a circle is the center of the circle.

2) **Theorem** (Bobillier)

If $O$ is a point in the plane of the triangle $ABC$ and the perpendiculars constructed in $O$ on $AO$, $BO$, $CO$ intersect respectively $BC$, $CA$, $AB$ in the points $A_1, B_1, C_1$ are collinear.

**Proof**

Let's consider a circle with the center in $O$, the triangle $ABC$ and we execute a duality transformation of Fig. 6 in rapport with the circle. We have $p(BC)=A'$ (the pole of the line $BC$ is the point $A'$), $p(CA)=B'$, $p(CB)=C'$. The polar of $A$ will be $B'C'$, the polar of $B$ will be $A'C'$ and the polar of $C=AC\bigcap BC$ is $A'B'$.

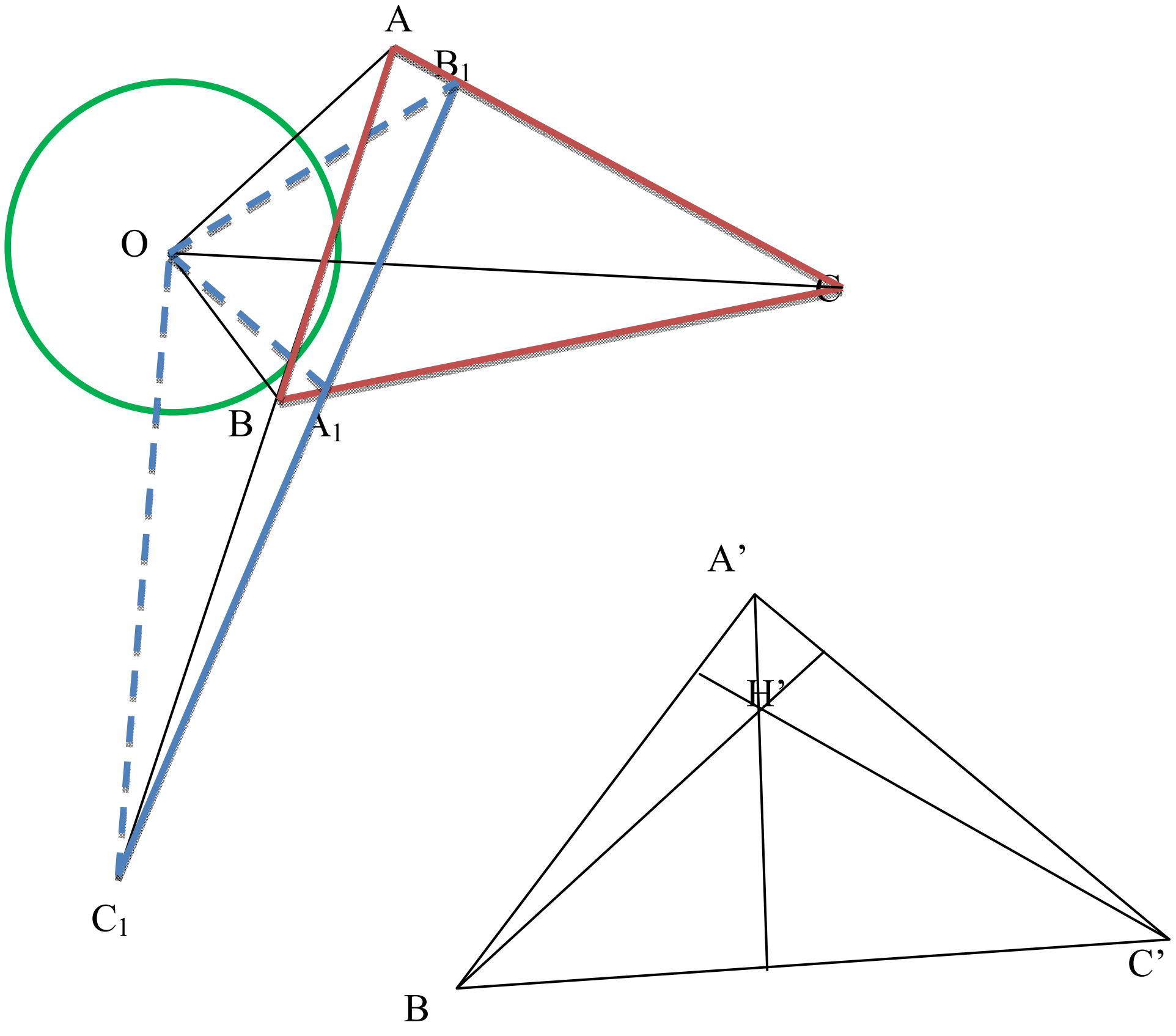

Fig. 6

Because $OA\perp OA_1$, it result that the polar of $A$ will be the perpendicular on the polar of $A_1$, therefore the polar of $A_1$ will be the perpendicular from $A'$ on $BC$. Similarly, the polar of $B_1$ will be the height from $B'$ of the triangle $A'B'C'$. And the polar of $C_1$ will be the height from $C'$ of the triangle $A'B'C'$. Because the heights of the triangle $A'B'C'$ are concurrent, it results that the orthocenter $H'$ of these triangle is the pole of the line determined by the points $A_1, B_1, C_1$.

3) **Theorem** (Pappus)

If on the side $(OX$ of the angle $XOA$ we consider the points $A, B, C$, and on the side $(OY$ the points $A_1, B_1, C_1$ such that $AB_1$ intersects $BA_1$ in the point $K$, $BC_1$ and $CB_1$ intersects in the point $L$, and $AC_1$ and $CA_1$ intersect in the point $M$, then the points $K, L, M$ are collinear.

**Proof**

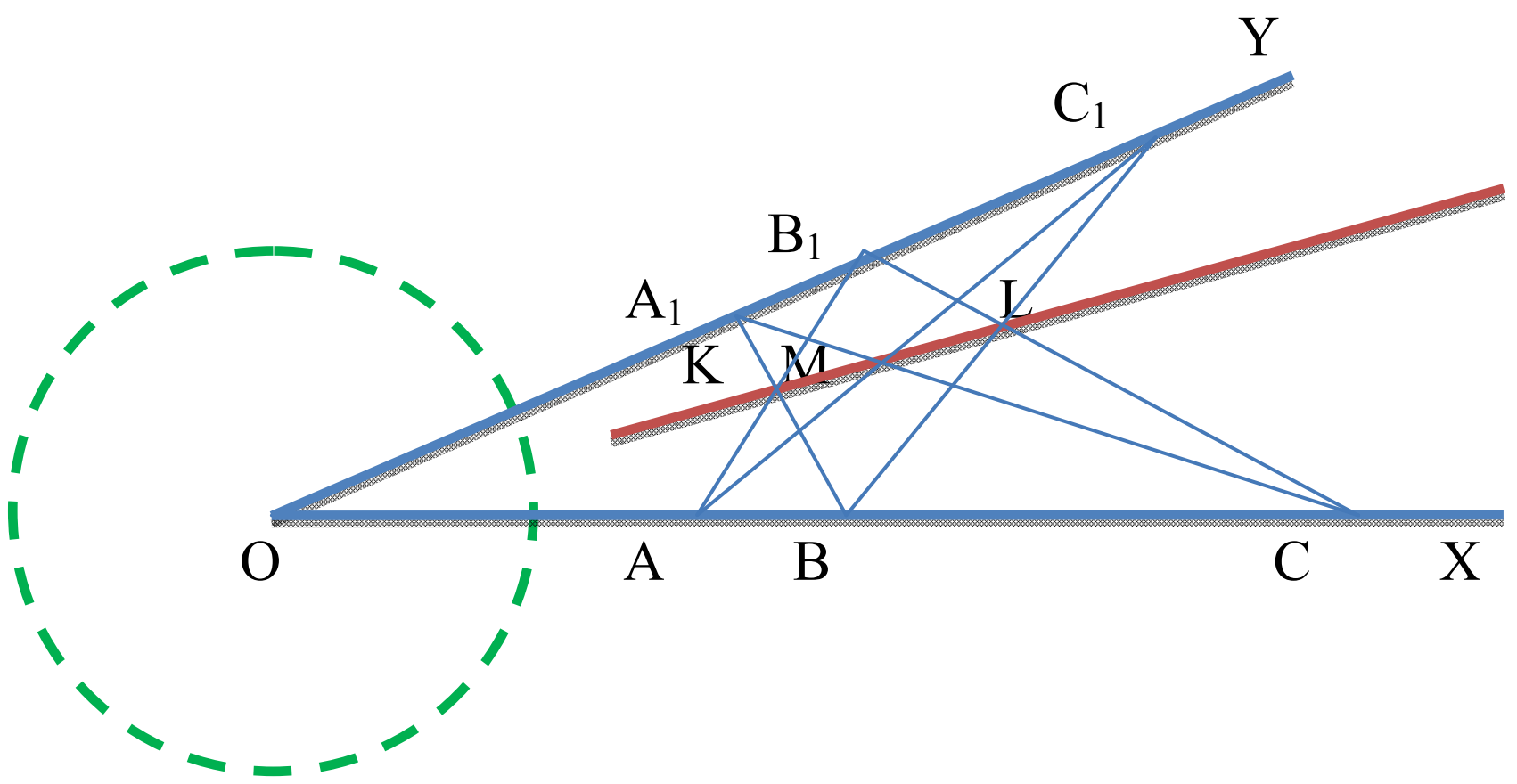

Fig. 7

We consider a circle with the center in $O$ and we'll transform by duality in rapport with this circle figure 7.

Because the points $A, B, C$ are collinear with the center $O$ of the circle in rapport with which we perform the transformation, it results that the polar a, b, c of these points will be parallels lines.

Similarly the polar of the points $A_1, B_1, C_1$ are the parallel lines $a_1$, $b_1$, $c_1$ (see Fig. 8).

The polar of point $K$ will be the line determined by poles of the lines $AB_1$ and $BA_1$ ($\{K\}=BA_1\bigcap AB_1$), therefore it will be the line $B_{c1}C_{b1}$.

Similarly the polar of $M$ will the line $A_{c1}C_{a1}$.

It can be proved without difficulty that the lines $A_{b1}B_{a1}, B_{c1}C_{b1}, A_{c1}C_{a1}$ are concurrent in a point $T$ . The polars being concurrent it means that their poles, i.e. the points $K, L, M$ are collinear, and the theorem is proved.

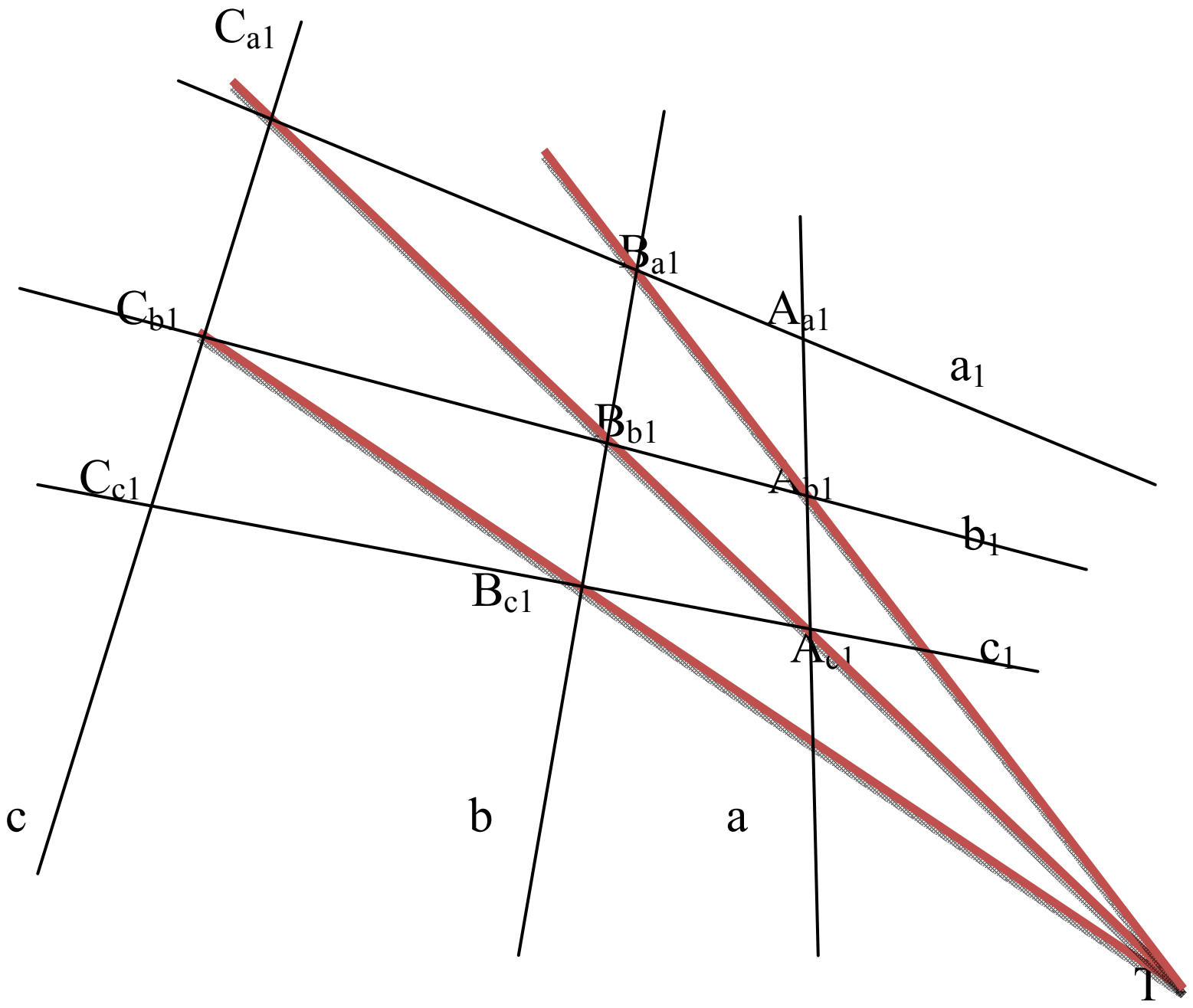

Fig. 8

## 6.6. Homothety

**Homothety definition**

Let $\pi$ a plane and $O$ a fixed point in this plan, and $k \in \mathbb{R}$, $k \neq 0$.

**Definition 1**

The homothety of center $O$ and of rapport $k$ is the transformation of the plane $\pi$ through which to any point $M$ from the plane we associate the point $M'$ such that $\overrightarrow{OM'} = k\overrightarrow{OM}$

We'll note $h_{(O,k)}$ the homothety of center $O$ and rapport $k$. The point $M' = h_{(O,k)}(M)$ is called the homothetic of point $M$ .

**Remark 1**

If $k > 0$, the points $M, M'$ are on the same side of the center $O$ on the line that contains them. In this case the homothety is called direct homothety.

If $k < 0$, the points $M, M'$ are placed on both sides of $O$ on the line that contains them. In this case the homothety is called inverse homothety.

In both situation described above the points $M, M'$ are direct homothetic respectively inverse homothetic.

If $k = -1$ the homothety $h_{(O,-1)}$ is the symmetry of center $O$

**Properties**

Give the homothety $h_{(O,k)}$ and a pair of points $M, N$, ($O \notin MN$) then $[MN] \| [M'N']$, where $M', N'$ are the homothetic of the points $M, N$ through the considered homothety and $\frac{M'N'}{MN} = |k|$.

**Proof**

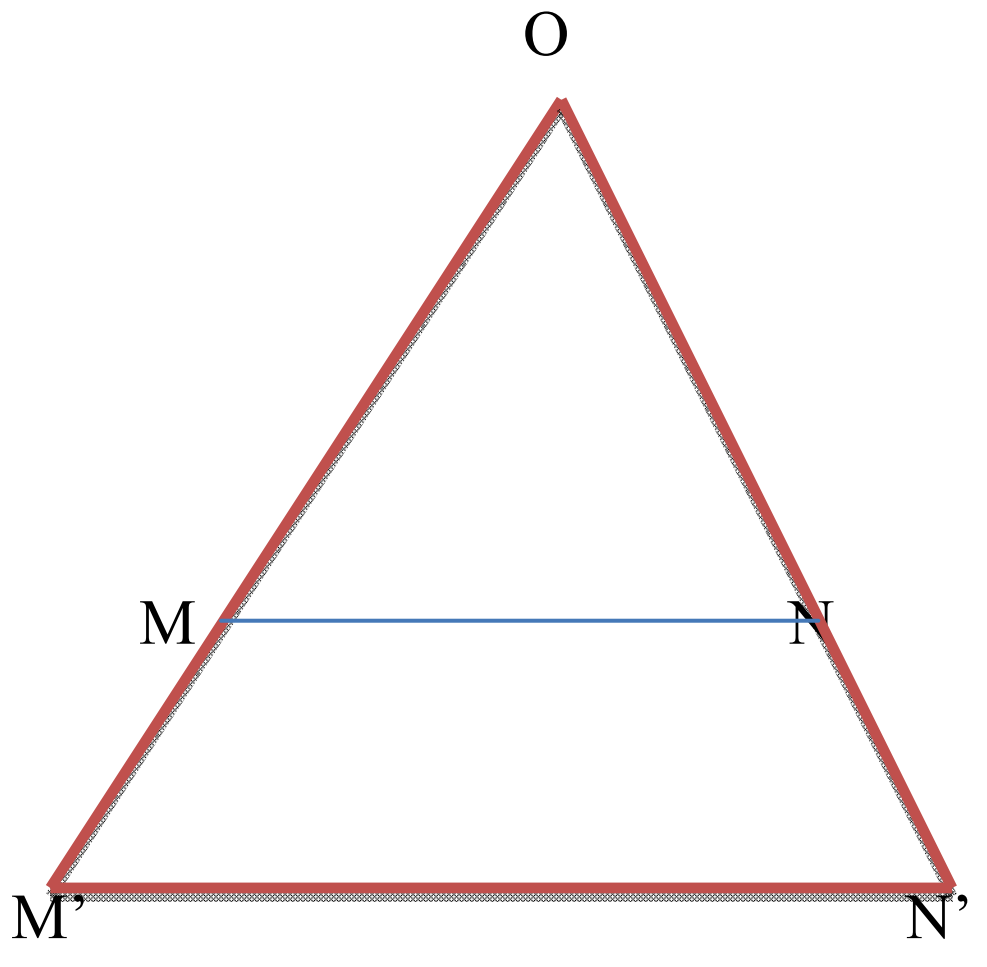

Fig. 1

From $\frac{OM'}{OM}=|k|,\ \frac{ON'}{ON}=|k|$ we have that $\Delta OMN \equiv \Delta OM'N'$, and therefore $\frac{M'N'}{MN}=|k|$.

We have that $\sphericalangle OMN \equiv \sphericalangle OM'N'$. It results that $M',N' \parallel MN$.

**Remark 2**

1. If we consider three collinear points $M, N, P$, then their homothetic points $M',N',P'$ are also collinear. Therefore, the homothety transforms a line (which does not contain the homothety center) in a parallel line with the given line.

The image, through a homothety, which passes through the homothety center, is that line.

2. If we consider a triangle $ABC$ and a homothety $h_{(O,k)}$ the image of the triangle through the considered is a triangle $A'B'C'$ similar with $ABC$. The similarity rapport being $\frac{A'B'}{AB}=|k|$. Furthermore, the sides of the two triangles are parallel two by two.

This result can be extended using the following: the homothety transforms a figure $\mathcal{F}$ in another figure $\mathcal{F}'$ parallel with it.

The reciprocal of this statement is also true, and we'll prove it for the case when the figure is a triangle.

**Proposition 2**

Let $ABC$ and $A'B'C'$ two triangles, where $AB \parallel A'B',\ AC \parallel A'C',\ BC \parallel B'C'$ ; $AB \neq A'B'$, then there exists a homothety $h_{(O,k)}$ such that $h_{(O,k)}(\Delta ABC)=\Delta A'B'C'$

**Proof**

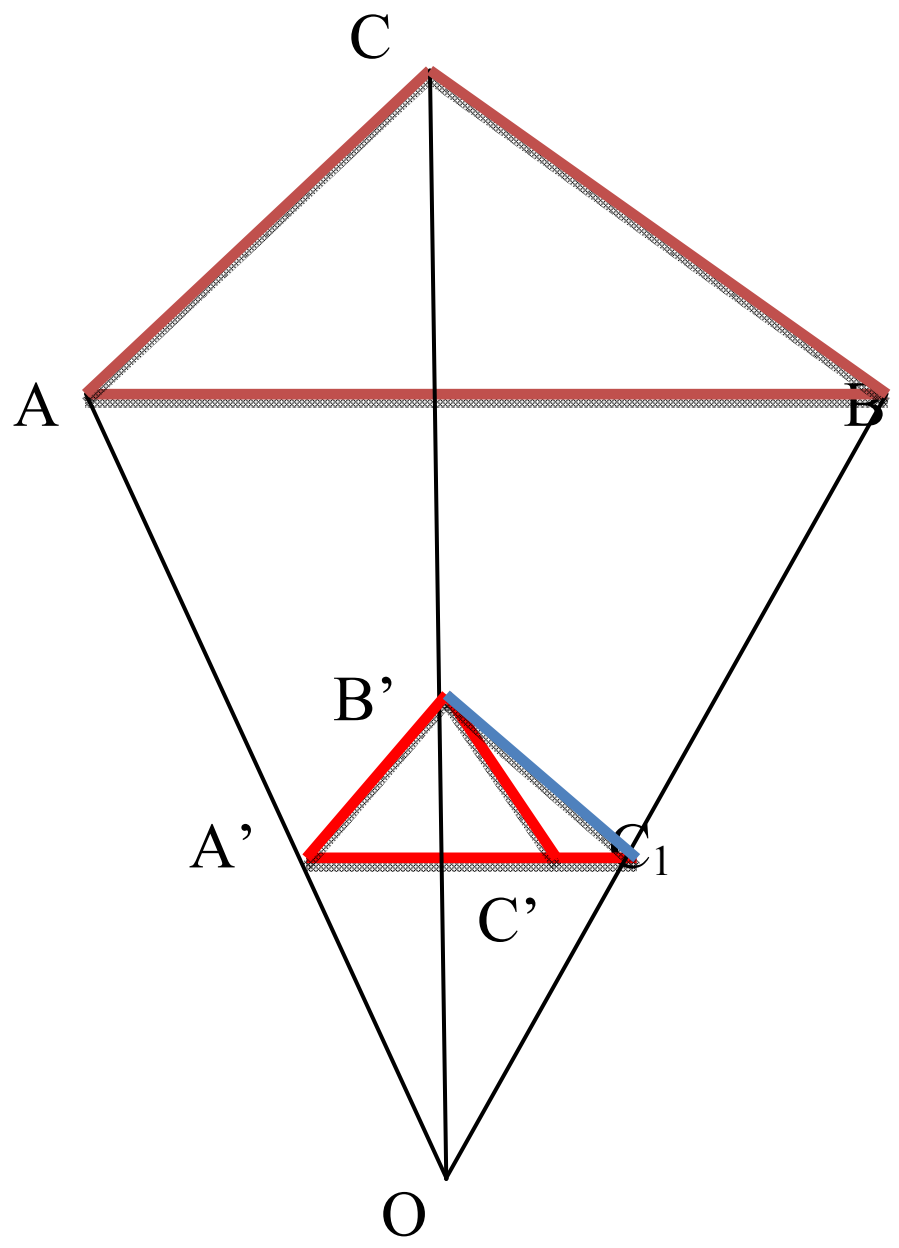

Fig. 2

Let $\{O\} = AA' \cap BB'$ and $\{C_1\} = OC \cap A'C'$ (see Fig. 2)

We have $\frac{OA'}{OA} = \frac{OB'}{OB} = \frac{A'B'}{AB} = \frac{A'C'}{AC} = \frac{A'C_1}{AC}$.

Therefore, $C_1 = C_o^{'}$, the triangles $ABC$, $A'B'C'$ are homothetic through the homothety of center $O$ and of rapport $\frac{A'B'}{AB} = |k|$.

3. **The product (composition) of two homothety**

a) **The product of homothety that have the same center**

Let $h_{(O,k_1)}$ and $h_{(O,k_2)}$ two homothety of the same center and $M$ a point in plane.

We have $h_{(O,k_1)}(M) = M'$, where $\overrightarrow{OM'} = k_1\overrightarrow{OM}$, similarly $h_{(O,k_2)}(M) = M''$, where $\overrightarrow{OM''} = k_2\overrightarrow{OM}$.

If we consider $\left(h_{(O,k_1)} \circ h_{(O,k_2)}\right)(M)$, we have

$$\left(h_{(O,k_1)} \circ h_{(O,k_2)}\right)(M) = h_{(O,k_1)}\left(h_{(O,k_2)}(M)\right) = h_{(O,k_1)}(M'') = M'''$$

where

$$\overrightarrow{OM'''} = k_1\overrightarrow{OM''} = k_1k_2\overrightarrow{OM} = (k_1k_2)\overrightarrow{OM}$$

Therefore $h_{(O,k_1)} \circ h_{(O,k_2)}(M) = h_{(O,k_1k_2)}(M)$.

It result that the product of two homothety of the same center is a homothety of the same center and of a rapport equal with the product of the rapports of the given homothety.

If we consider $\left(h_{(O,k_1)} \circ h_{(O,k_2)}\right)(M)$ we obtain $h_{(O,k_2)}\left(h_{(O,k_1)}(M)\right) = h_{(O,k_2)}(M') = M_1$, where $\overrightarrow{OM_1} = k_2\overrightarrow{OM'} = k_1k_2\overrightarrow{OM} = (k_1k_2)\overrightarrow{OM}$. But we noted $(k_1k_2)\overrightarrow{OM} = \overrightarrow{OM'''}$, therefore $M_1 = M'''$ and $h_{(O,k_1)} \circ h_{(O,k_2)} = h_{(O,k_2)} \circ h_{(O,k_1)}$, in other words the product of two homothety is commutative.

The following proposition is true.

**Proposition 3**

The homothety of a plane having the same center form in rapport with the composition an Abel group.

**Remark 3**

The inverse of the homothety $h_{(O,k)} : \pi \to \pi$ is the homothety $h_{(O,\frac{1}{k})} : \pi \to \pi$

b) **The product of homotheties of different centers**

Let $h_{(O,k_1)}$ and $h_{(O,k_2)}$ two homothety in the plane $\pi$ and $\mathcal{F}$ a figure in this plane. Transforming the figure through the homothety $h_{(O,k_1)}$, we obtain a figure $\mathcal{F}_1$ parallel with $\mathcal{F}$. If we transform the figure $\mathcal{F}_1$ through the homothety $h_{(O,k_2)}$ we'll obtain figure $\mathcal{F}_2$ parallel with $\mathcal{F}_1$. Because the parallelism relation is a transitive relation, it results that figure $\mathcal{F}$ is parallel with $\mathcal{F}_2$, therefore $\mathcal{F}_2$ can be obtained from $\mathcal{F}$ through a homothety. Let's see which is the center and the rapport of this homothety.

The line $O_1O_2$ passes through the center of the first homothety, therefore it is invariant through it, also it contains the center of the second homothety. Therefore it is also invariant through this homothety. It results that $O_1O_2$ in invariant through the product of the given homothety which will have the center on the line $O_1O_2$.

**Proposition 4**

The product of two homothety of centers, different points, and of rapport $k_1, k_2$, such that $k_1 \cdot k_2 \neq 1$ is a homothety with the center on the line if the given centers of homotheties and of equal rapport with the product of the rapports of the given homotheties.

**Proof**

Let $h_{(O_1,k_1)}, h_{(O_2,k_2)}$ the homotheties, $O_1 \neq O_2$ and $\mathcal{F}$ a given figure.

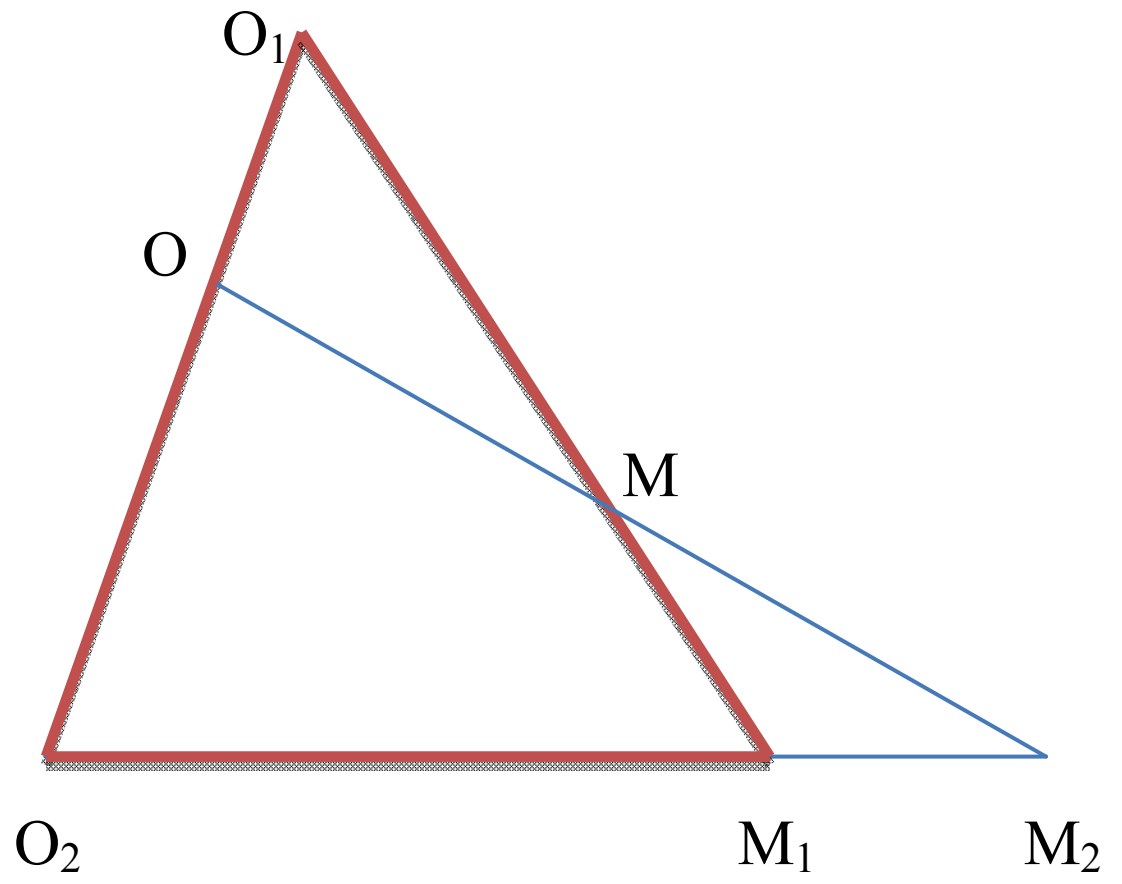

Fig. 3

We note $\mathcal{F}_1 = h_{(O_1,k_1)}(\mathcal{F})$, $\mathcal{F}_2 = h_{(O_2,k_2)}(\mathcal{F}_1)$, If $M \in \mathcal{F}$, let $M_1 = h_{(O_1,k_1)}(M)$, then $M_2 \in \mathcal{F}_1$.
$\overrightarrow{O_1M_1} = k_1\overrightarrow{O_1M}$ ; $M_2 = h_{(O_2,k_2)}$,
therefore

$$\overrightarrow{O_2M_1} = k_2\overrightarrow{O_2M_1}.$$

We note $\{O\} = MM_2 \cap O_1O_2$ (see figure 3).

Applying the Menelaus' theorem in the triangle $MM_1M_2$ for the transversal $O_1 - O - O_2$, we obtain

$$\frac{\overline{O_1M}}{\overline{O_1M_1}} \cdot \frac{\overline{OM_2}}{\overline{OM}} \cdot \frac{\overline{O_2M_1}}{\overline{O_2M_2}} = 1.$$

Taking into account that

$$\frac{\overline{O_1M_1}}{\overline{O_1M}} = k_1 \text{ and } \frac{\overline{O_2M_2}}{\overline{O_2M_1}} = k_2,$$

we obtain

$$\frac{\overline{OM_2}}{\overline{OM}} = k_1k_2 .$$

Therefore the point $M_2$ is the homothetic of the point $M$ through the homothety $h_{(O,k_1k_2)}$.

In conclusion $h_{(O_1,k_1)} \circ h_{(O_2,k_2)} = h_{(O,k_1k_2)}$, where $O_1, O_2, O$ are collinear and $k_1 \cdot k_2 \neq 1$.

**Remark 4.**

The product of two homotheties of different centers and of rapports of whose product is equal to 1 is a translation of vectors of the same direction as the homotheties centers line.

**Applications**

**1.**

Given two circles non-congruent and of different centers, there exist two homotheties (one direct and the other inverse) which transform one of the circles in the other one. The centers of the two homotheties and the centers of the given circles forma harmonic division.

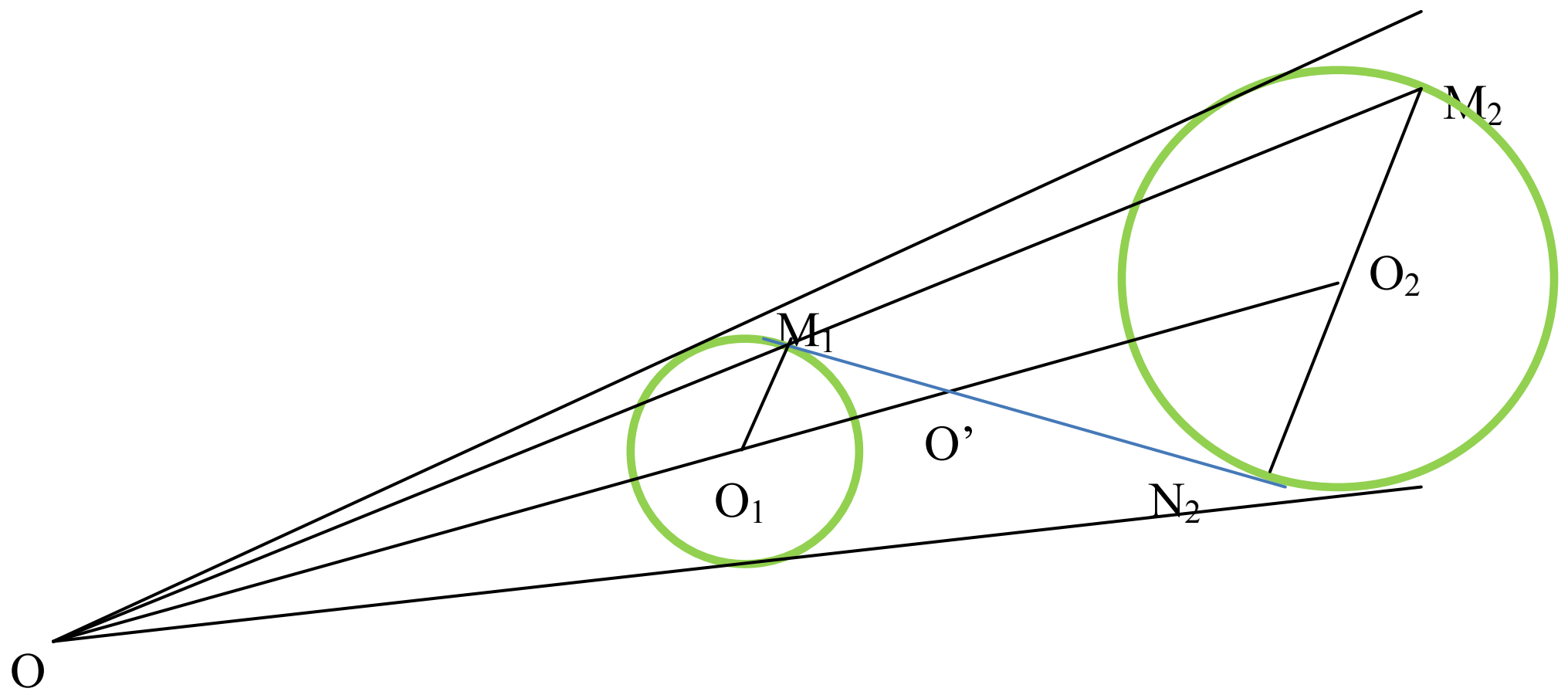

Fig.4

**Proof**

Let $\mathcal{C}(O_1, r_1)$ and $\mathcal{C}(O_2, r_2)$ the given circles, $r_1 < r_2$ (see figure 4)

We construct two parallel radiuses in the same sense: $O_1M_1, O_2M_2$ in the given circles. We note with $O$ the intersection of the lines $O_1O_2$ and $M_1M_2$. From the similarity of the triangles $OM_1O_1$ and $OM_2O_2$, it results

$$\frac{\overrightarrow{OO_1}}{\overrightarrow{OO_2}} = \frac{\overline{O_1M_1}}{\overline{O_1M}} = \frac{\overline{OM_1}}{\overline{OM_2}} = \frac{r_1}{r_2}$$

It results that the point $O$ is fix and considering the point $M_2$ mobile on $\mathcal{C}(O_2, r_2)$, there exists the homothety $h_{(O,\frac{r_1}{r_2})}$, which makes to the point $M_2$, the point $M_1 \in \mathcal{C}(O_1, r_1)$. $\overrightarrow{OM_1} = \frac{r_1}{r_2}\overrightarrow{OM_2}$. Through the cited homothety the circle $\mathcal{C}(O_2, r_2)$ has as image the circle $\mathcal{C}(O_1, r_1)$. If the point $N_2$ is the diametric opposed to the point $M_2$ in $\mathcal{C}(O_2, r_2)$ and $\{O'\} = O_1O_2 \cap M_1N_2$, we find

$$\frac{\overrightarrow{OO_1}}{\overrightarrow{O'O_2}} = \frac{\overline{O'M_1}}{\overline{O'N_2}} = -\frac{r_1}{r_2}$$

Therefore the circle $\mathcal{C}(O_1, r_1)$ is obtained from $\mathcal{C}(O_2, r_2)$ through the homothety $h_{(O',-\frac{r_1}{r_2})}$.

The relation $\frac{\overrightarrow{OO_1}}{\overrightarrow{OO_2}} = -\frac{\overline{O'O_1}}{\overline{O'O_2}}$ shows that the points $O, O_1, O'O_2$ form a harmonic division.

**Remark 5**

The theorem can be proved similarly and for the case when the circles: interior, exterior tangent and interior tangent.

In the case of tangent circles one of the homothety centers is the point of tangency. If the circles are concentric, then there exists just one homothety which transforms the circle $\mathcal{C}(O_1, r_1)$ in the circle $\mathcal{C}(O_2, r_2)$, its center being $O = O_1 = O_2$ and the rapport $\frac{r_1}{r_2}$

**2. G. Monge Theorem**

If three circles are non-congruent two by two and don't have their centers collinear, then the six homothety centers are situated in triplets on four lines.

**Proof**

Let $S_1, S_1'$ the direct and inverse homothety centers of the circle $\mathcal{C}(O_2, r_2)$, $\mathcal{C}(O_3, r_3)$, similarly $S_2, S_2'$, $S_3, S_3'$. In figure 5 we considered $r_1 < r_2 < r_3$.

We prove the collinearity of the centers $S_1, S_2, S_3$.

Through homothety the circle $\mathcal{C}(O_1, r_1)$ gets transformed in the circle $\mathcal{C}(O_2, r_2)$, and through the homothety $h_{(S_1,\frac{r_3}{r_2})}$ the circle $\mathcal{C}(O_2, r_2)$ gets transformed I n the circle $\mathcal{C}(O_3, r_3)$. By composing these two homotheties we obtain a homothety of a rapport $\frac{r_3}{r_1}$ and of a center which is collinear with $S_1$, $S_3$ and placed on the line $O_1O_3$.

This homothety transforms the circle $\mathscr{C}(O_1, r_1)$ in the circle $\mathscr{C}(O_3, r_3)$, therefore its center is the point $S_2$, and therefore $S_1, S_2, S_3$ are collinear.

Similarly we prove the theorem for the rest of the cases.

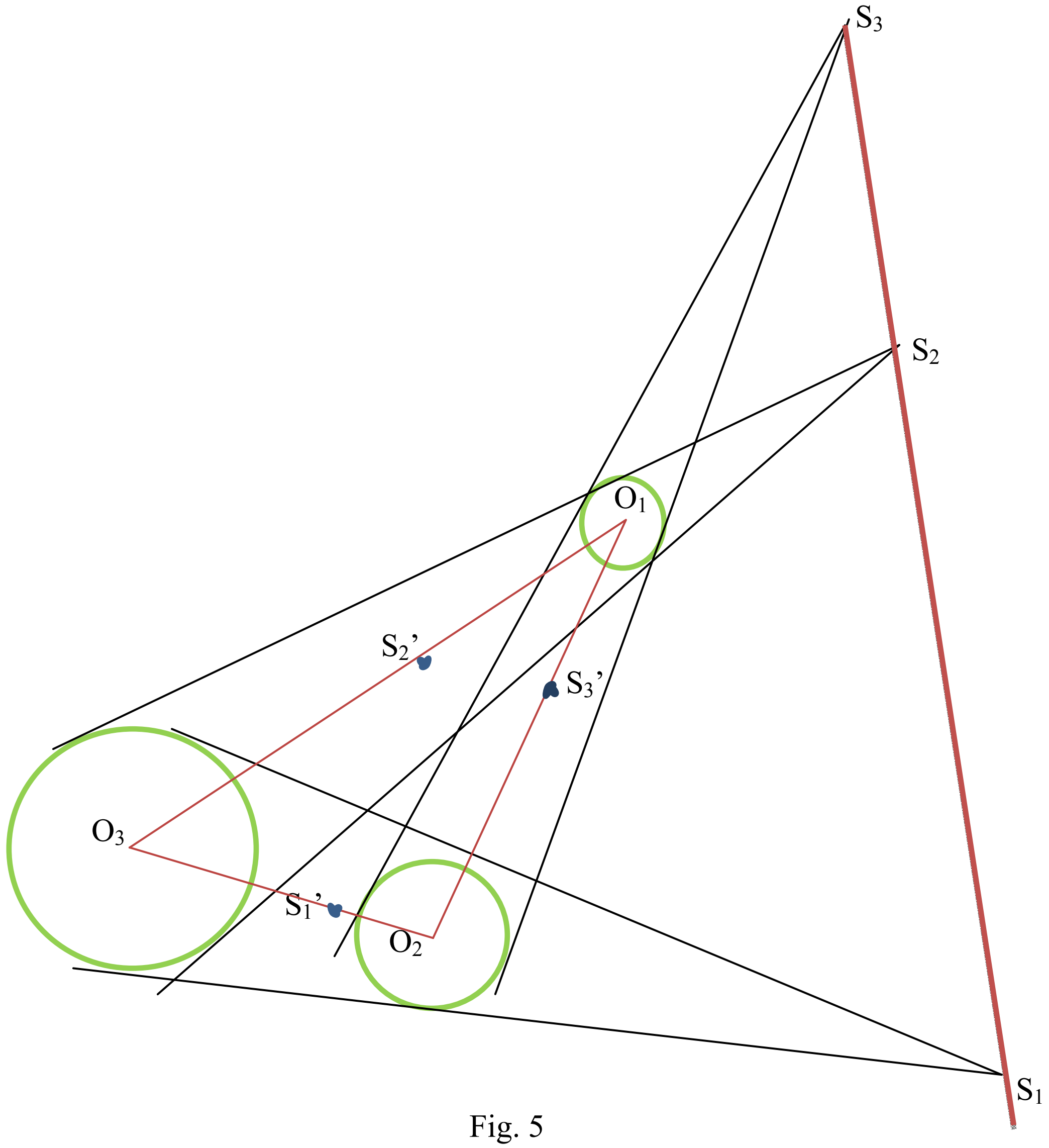

Fig. 5

**3.**

In a triangle the heights' feet, the middle points of the sides and the middle segments determined by the triangle orthocenter with its vertexes are nine concyclic points (the circle of the nine points).

**Proof**

It is known that the symmetric points $H_1, H_2, H_3$ of the orthocenter $H$ in rapport to the triangle's $ABC$ sides are on the triangle's circumscribed circle. Then, considering the homothety $h_{(H,\frac{1}{2})}$, we obtain that the circumscribed circle gets transformed through this in the circumscribed circle to the orthic triangle $A'B'C'$ of the triangle $ABC$ (see figure 6).

The center of this circle will be the middle of the segment $(OH)$, we'll note this point $O_9$

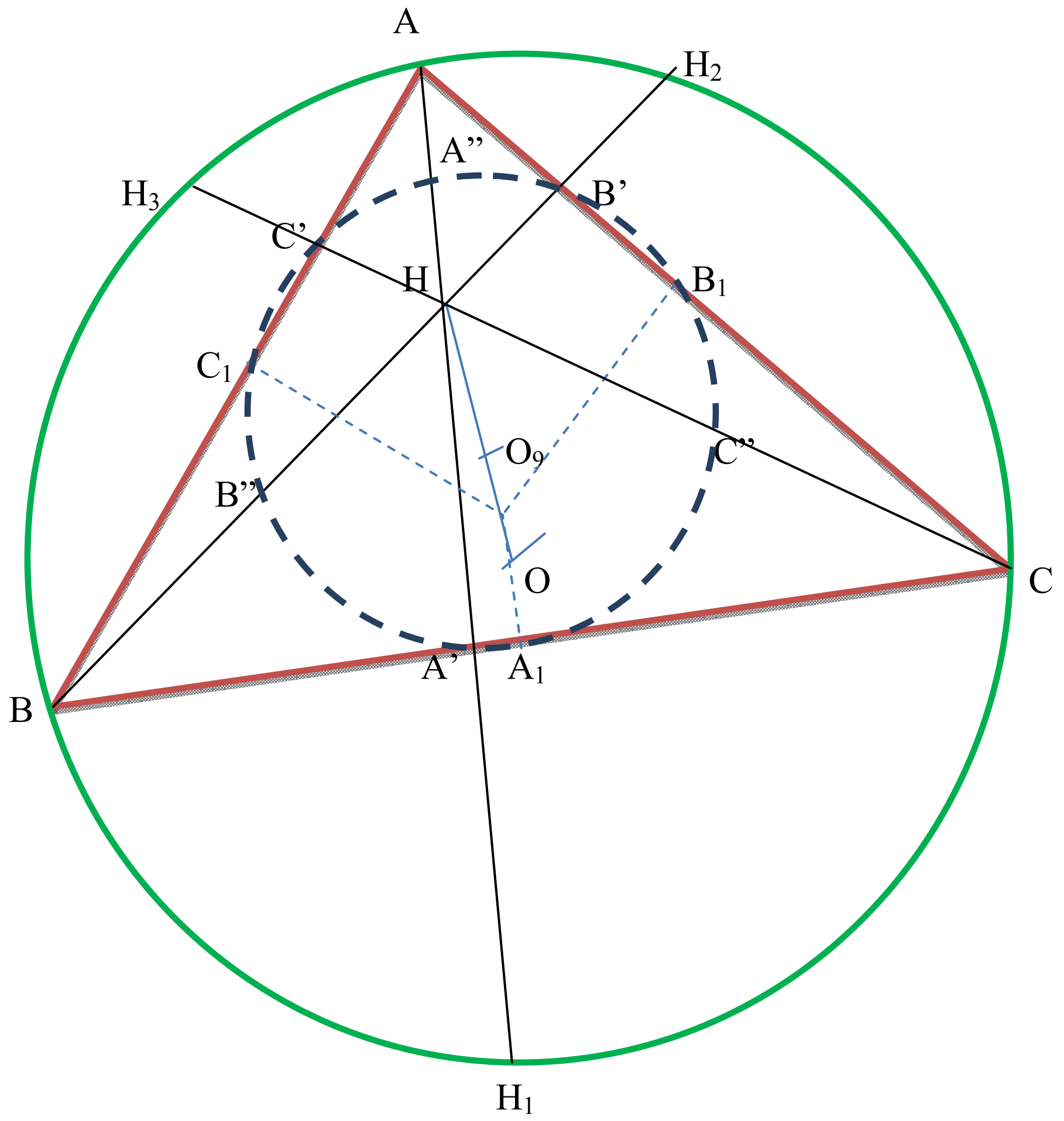

Fig. 6

And the radius of the circumscribed circle to triangle $A'B'C'$ will be $\frac{R}{2}$.

Also, on this circle will be situated the points $A'', B'', C''$ which are the middle of the segments $(AH), (BH), (CH)$ the symmetric of the points $A, B, C$ through the considered homothety.

The medial triangle $A_1B_1C_1$ has its sides parallel with the sides of the triangle $ABC$, therefore these are homothetic. Through the homothety $h_{(G,\frac{1}{2})}$ the circumscribed circle to the triangle $ABC$ gets transformed in the circle $\mathcal{C}\left(O_9, \frac{R}{2}\right)$, which contains the middle points $A_1, B_1, C_1$ of the sides of the triangle $ABC$.

Therefore, the points $A', B', C', A_1, B_1, C_1, A'', B'', C''$ belong to the circle $\mathcal{C}\left(O_9, \frac{R}{2}\right)$, which is called the circle of the nine points.

**Remark**

The circumscribed circle and the circle of the nine points are homothetic and their direct, and inverse homothety centers are the points $H, G$.

In conformity with application 1, it results that the points $O, G, O_9, H$ are collinear and these form a harmonic division. The line of the points $O, G, H$ is called The Euler's line.

## 6.7. Inversion in plane

### A. Definition, Notations

**Definition 1**

Let $O$ a fixed point in the plane $\pi$ and $k$ a real number not null. We call an inversion of pole $O$ and of module (power) $k$ the geometrical transformation, which associates to each point $M \in \pi \setminus \{0\}$ the point $M' \in \pi \setminus \{0\}$ such that :

1. $O, M, M'$ are collinear
2. $\overrightarrow{OM} \cdot \overrightarrow{OM'} = k$.

We'll note $i_0^k$ the inversion of pole $O$ and of module $k$.

The point $M' = i_0^k(M)$ is called the inverse (image) of the point $M$ through the inversion of pole $O$ and power $k$. The points $M$ and $i_0^k(M)$ are called homological points of the inversion $i_0^k$.

**Remark 1.**

a) If $k > 0$, then the inversion $i_0^k$ is called a positive inversion, and if $k < 0$ the inversion is called a negative inversion.
b) From the definition it results that a line $d$ that passes through the inversion's pole, through the inversion $i_0^k$ has as image the line $d \setminus \{0\}$.
c) From the definition of inversion it results that the point $M$ is the inverse of the point $M'$ through the inversion $i_0^k$.

### B. The image of a figure through an inversion

We consider the positive inversion $i_0^k$, we saw that the lines that pass through $O$ are invariant through this inversion. We propose to find the images (inverses) of some remarkable figures such as the circle and the line through this inversion.

**Theorem 1**

If $i_0^k$ is a positive inversion, then the circle $\mathcal{C}(O, \sqrt{k})$ is invariant point by point through this inversion, then through this inversion the interior of the circle $\mathcal{C}(O, \sqrt{k})$ transforms in the exterior of the circle $\mathcal{C}(O, \sqrt{k})$ and the reciprocal being also true.

**Proof**

Let $M \in \mathcal{C}(O, \sqrt{k})$, we have

$$\overrightarrow{OM} \cdot \overrightarrow{OM'} = \left(\sqrt{k}\right)^2 = k$$

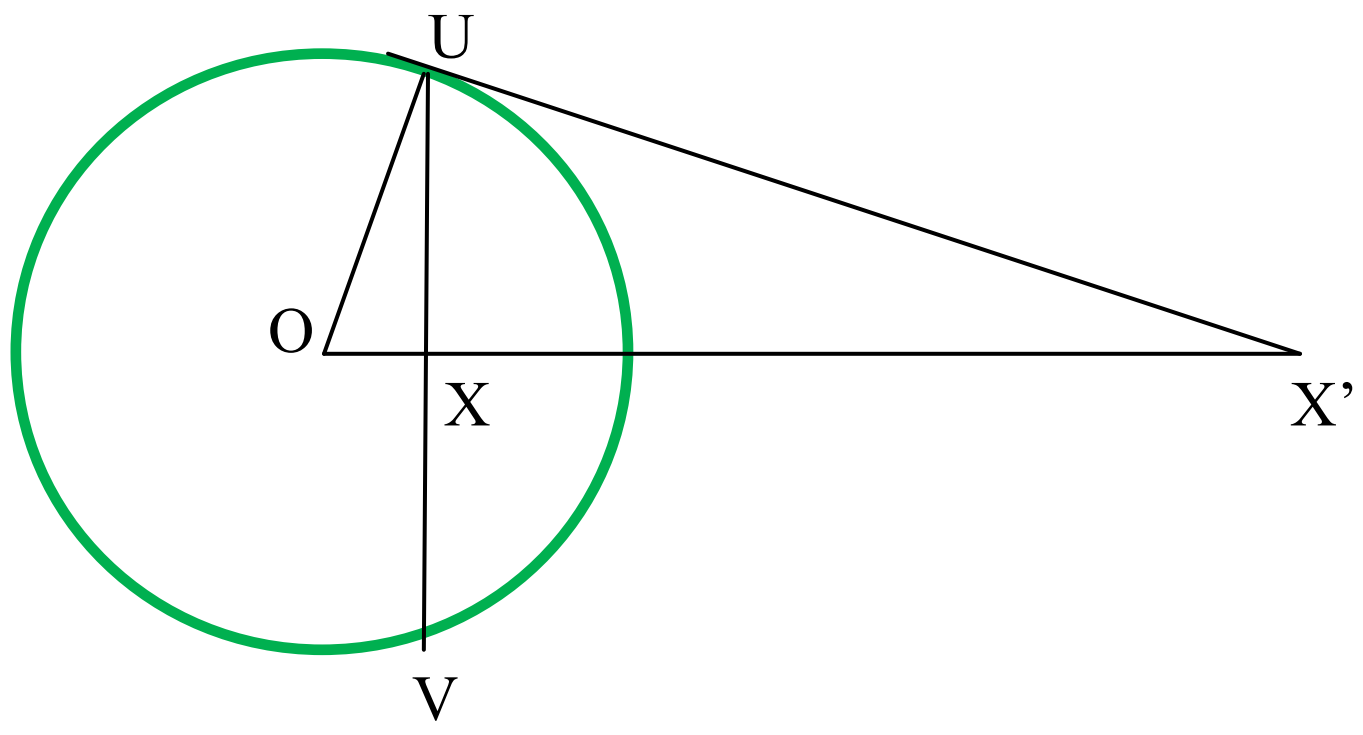

Fog. 1

Therefore $i_0^k(M) = M$ for any $M \in \mathcal{C}\left(O, \sqrt{k}\right)$, therefore $i_0^k\left(\mathcal{C}\left(O, \sqrt{k}\right)\right) = \mathcal{C}\left(O, \sqrt{k}\right)$.

Let now $X \in Int\mathcal{C}\left(O, \sqrt{k}\right)$, we construct a perpendicular cord in $X$ on $OX$, let it be $(UV)$ (see figure 1). The tangent in $U$ to $\mathcal{C}\left(O, \sqrt{k}\right)$ intersects $OX$ in $X'$. From the legs' theorem applied in the right triangle $OUX'$ it results that $\overrightarrow{OX} \cdot \overrightarrow{OX'} = OU^2$, therefore $\overrightarrow{OX} \cdot \overrightarrow{OX'} = k$, which shows that ${}_0^k(X) = X'$, evidently $X' \in Ext\mathcal{C}\left(O, \sqrt{k}\right)$.

If we consider $X \in Ext\mathcal{C}\left(O, \sqrt{k}\right)$, constructing the tangent $XT$ to the circle $\mathcal{C}\left(O, \sqrt{k}\right)$ and the projection $X'$ of the point $T$ on $OX$, we find that

$$\overrightarrow{OX} \cdot \overrightarrow{OX'} = k,$$

therefore

$$i_0^k(X) = X' \text{ and } X' \in Int\mathcal{C}\left(O, \sqrt{k}\right).$$

**Remark 2**

From this theorem it results a practical method to construct the image of a point $X$ through a positive inversion $i_0^k$.

**Definition 2**

If $i_0^k$ is a positive conversion, we say that the circle $\mathcal{C}\left(O, \sqrt{k}\right)$ is the fundamental circle of the inversion $i_0^k$ or the inversion circle.

**Theorem 2**

The image of a line $d$ that does not contain the pole $O$ of the inversion $i_0^k$, through this inversion, is a circle which contains the pole $O$, but from which we exclude $O$, and which has the diameter which passes through $O$ perpendicular on the line $d$.

**Proof**

We'll try to find the geometrical locus of the points $M'$, from plane, with the property that $\overrightarrow{OM} \cdot \overrightarrow{OM'} = k$, where $M$ is a mobile point on $d$.

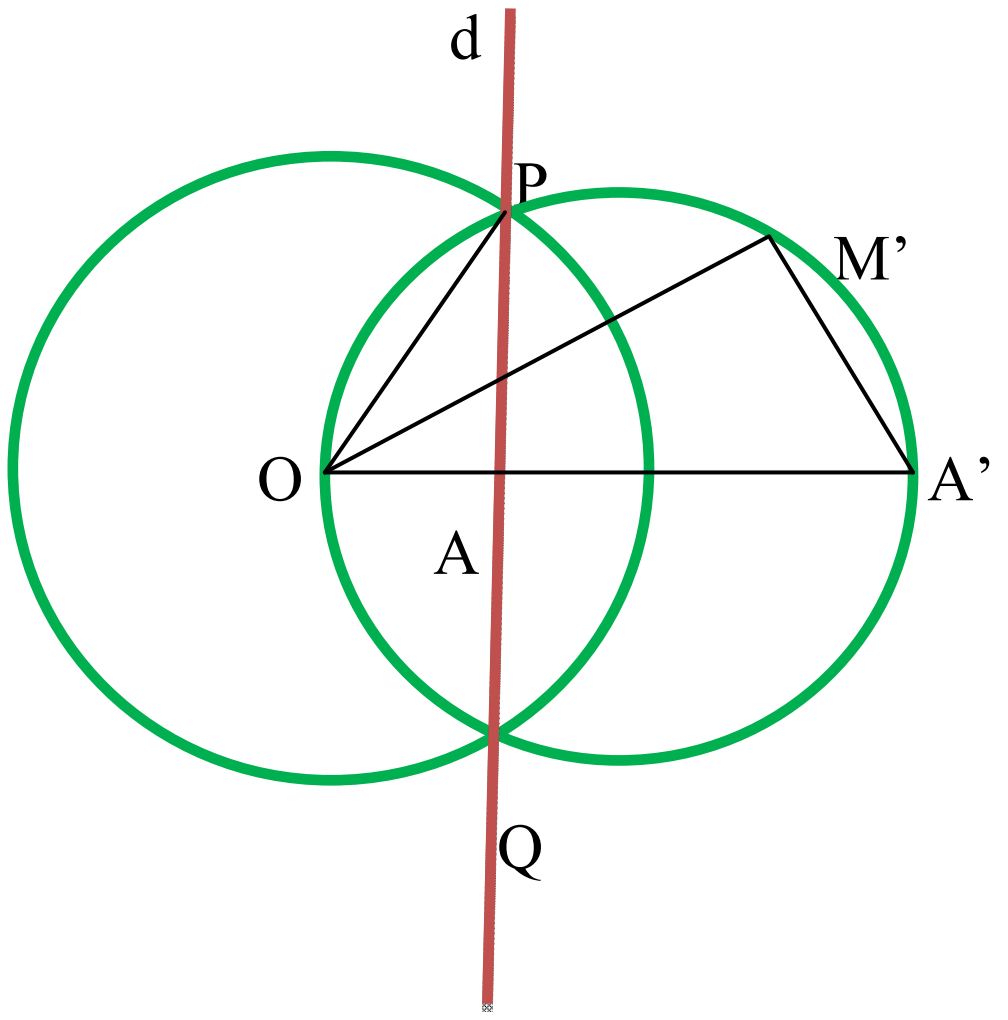

Fig. 2

We'll consider the point $A$ the orthogonal projection of the point $O$ on the line $d$; let $A'$ be the inverse of the point $A$ through the inversion $i_0^k$. We have $\overrightarrow{OA} \cdot \overrightarrow{OA'} = k$. If $M$ is a random point on $d$ and $M' = i_0^k(M)$, we have $\overrightarrow{OM} \cdot \overrightarrow{OM'} = k$. The relation $\overrightarrow{OA} \cdot \overrightarrow{OA'} = \overrightarrow{OM} \cdot \overrightarrow{OM'}$ shows that the quadrilateral $AA'M'M$ is inscribable (see figure 2). Because $m(\sphericalangle MAA') = 90°$ it results that also $m(\sphericalangle MM'A') = 90°$, consequently, taking into consideration that $OA$ is constant, therefore $OA'$ is constant, it results that form $M'$ the segment $(OA')$ is seen under a right triangle, which means that the geometric locus of the point $M'$ is the circle whose diameter is $OA'$. The center of this arc is the middle of the segment $(OA')$. If we'll consider $N'$ a random point on this circle and considering $\{N\} = ON' \cap d$, then $NN'A'A$ is an inscribable quadrilateral and $\overrightarrow{OA} \cdot \overrightarrow{OA'} = \overrightarrow{ON} \cdot \overrightarrow{ON'} = k$, therefore $N'$ is the inverse of $N$ through $i_0^k$.

**Remark 3**

a) If $P, Q$ are the intersections of the line $d$ with the fundamental circle of the inversion (these points do not always exist) we observe that the inverse of these points are the points themselves, therefore these are also located on the circle image through $i_0^k$ of the line $d$. The line $PQ$ is the radical axis of the circle $\mathcal{C}(O, \sqrt{k})$ and of the inverse circle to the line $d$. In general, the line $d$ is the radical axis of the circle $OA^2 + O_1A^2 = OO_1^2$ and of the image circle of the line $d$ through the image $i_0^k$.

b) The radius of the circle's image of the line $d$ through the positive inversion $i_0^k$ is equal to $\frac{k}{2a}$, where $a$ is the distance from $O$ to the line $d$.

c) The points of quartet constructed from two pairs of homological points through an inversion are concyclic if neither of them is the inversion pole.

Because of the symmetry of the relation through which are defined the inverse points it is true the following theorem.

**Theorem 3.**

The image through the positive inversion $i_0^k$ of a circle which passes through $O$ (from which we exclude the point $O$) is a line (the radical axis of the given circle and of the fundamental circle of the inversion $\mathcal{C}\left(O, \sqrt{k}\right)$

**Theorem 4.**

The image through the positive inversion $i_0^k$ of a circle, which does not contain its center is a circle which does not contain he pole of the inversion $O$.

**Proof**

Let the given circle $\mathcal{C}(O_1, r_1)$ and the positive inversion $i_0^k$, $O \notin \mathcal{C}(O_1, r_1)$.

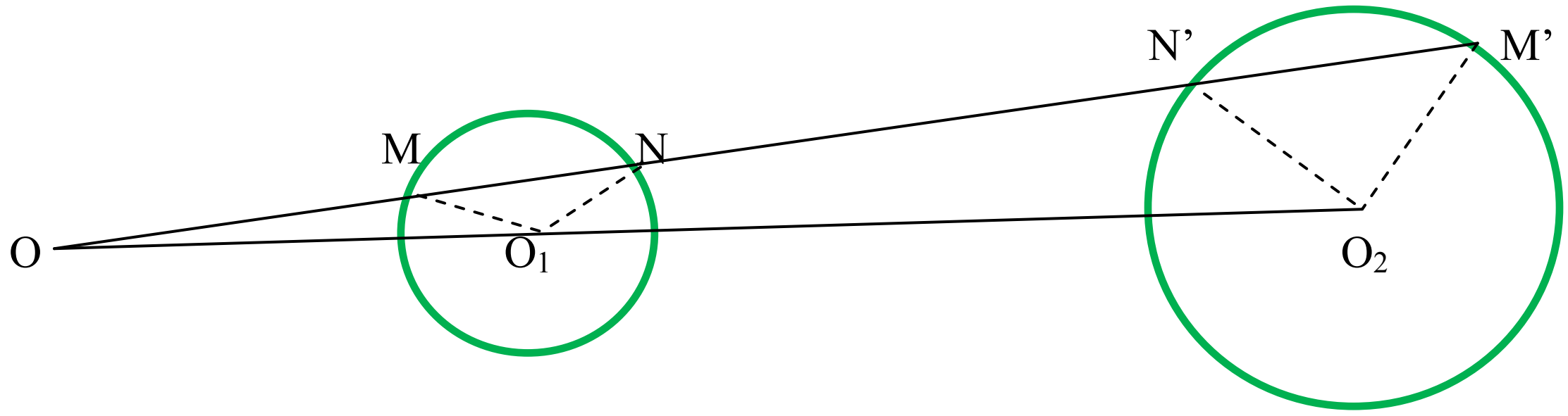

Fig. 3

We’ll consider the secant $O, M, N$ for the given circle and let $M' = i_0^k(M)$, $N' = i_0^k(N)$, see figure 3.

We have

$$\overrightarrow{OM} \cdot \overrightarrow{OM'} = k \qquad (1)$$

$$\overrightarrow{ON} \cdot \overrightarrow{ON'} = k \qquad (2)$$

It is known that $\overrightarrow{OM} \cdot \overrightarrow{ON} = const$ (the power of the point $O$ in rapport to the circle $\mathcal{C}(O_1, r_1)$.

We note

$$\overrightarrow{OM} \cdot \overrightarrow{ON} = p \qquad (3)$$

From the definition of inversion we have that the points $O, M, N, M', M'$ are collinear.

The relations (1) and (3) lead to $\frac{\overrightarrow{OM'}}{\overrightarrow{OM}} = \frac{k}{r}$; from (2) and (3) we obtain $\frac{\overrightarrow{ON'}}{\overrightarrow{OM}} = \frac{k}{p}$.

These relations show that the point $M'$ is the homothetic of the point $N$ through the homothety $h_o^{\frac{k}{p}}$ (also the point $N'$ is the homothetic of the point $M$ through the same homothety), consequently the geometric locus of the point $M'$ is a circle which is the homothetic of the circle $\mathcal{C}(O_1, r_1)$ through the homothety $h_o^{\frac{k}{p}}$. We will note this circle $\mathcal{C}(O_2, r_2)$, where $r_2 = \frac{k}{p} r_1$ and $p = \left| OO_1^2 - r_1^2 \right|$.

**Remark 4**

If the power of the pole $O$ of the inversion $i_0^k$ in rapport to the given circle $\mathcal{C}(O_1, r_1)$ is equal with $k$, then the circle $\mathcal{C}(O_1, r_1)$ is invariant through $i_0^k$.

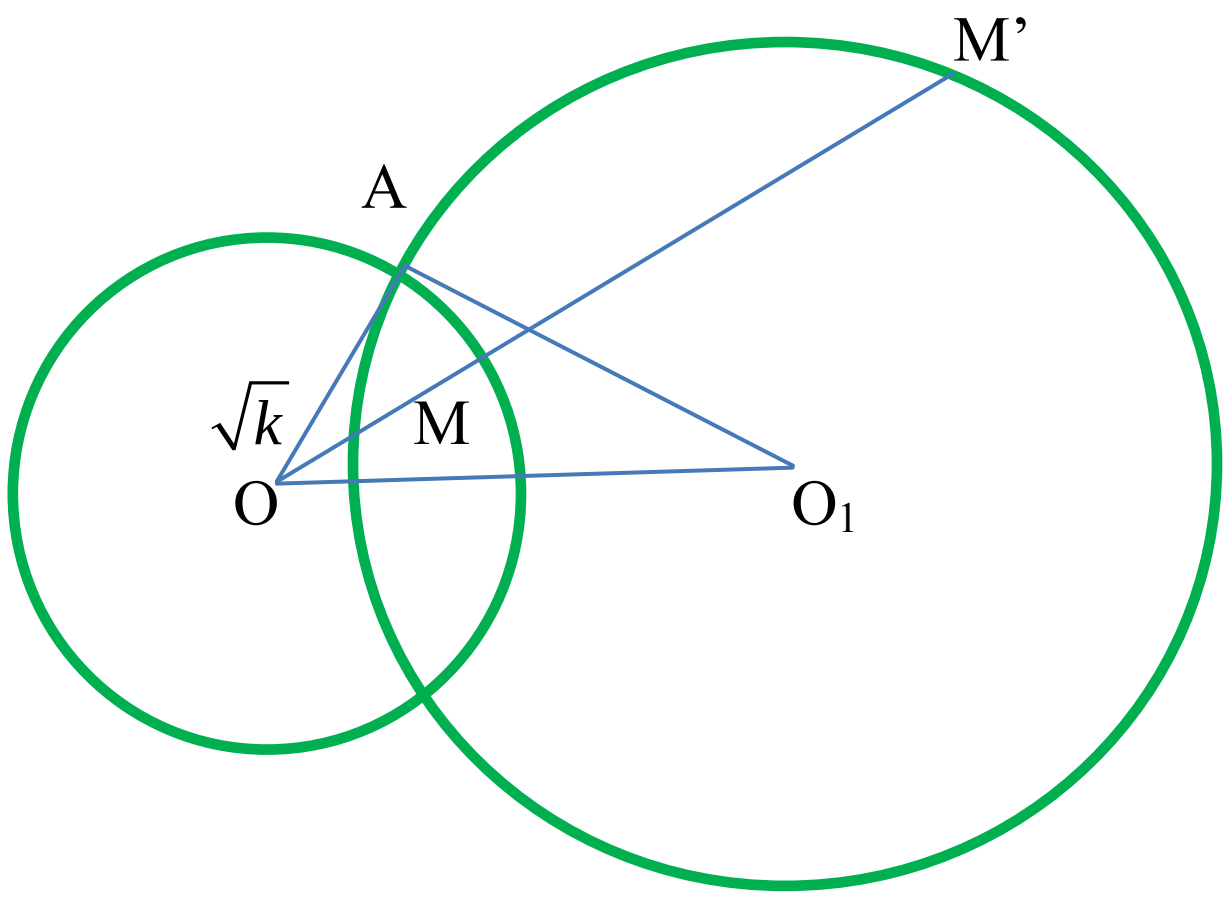

Fig. 4

Indeed, if $M$ belongs to the circle $\mathcal{C}(O_1, r_1)$ and $M'$ is the second intersection of the line $OM$ with the circle, we have:

$$\overrightarrow{OM} \cdot \overrightarrow{OM'} = k = OA^2.$$

This shows that $OA^2 + O_1A^2 = OO_1^2$, therefore the circles $\mathcal{C}\left(O, \sqrt{k}\right)$, $\mathcal{C}(O_1, r_1)$ are orthogonal (see figure 4).

## C. The construction with the compass and the ruler of the inverses of a line and of a circle

### 1. The construction of the inverse of a line

Let $\mathcal{C}\left(O, \sqrt{k}\right)$, the inversion $i_0^k$ and the line $d$.

If the line $d$ is external to the circle $\mathcal{C}\left(O, \sqrt{k}\right)$ we construct the orthogonal projection of $O$ on $d$, then the tangent $AT$ to the circle $\mathcal{C}\left(O, \sqrt{k}\right)$. We construct the projection $A'$ of $T$ on $OA$, we have $A' = i_0^k(A)$.

We construct the circle of diameter $[OA']$, this without the point $O$ represents $i_0^k(d)$ (see figure 5)

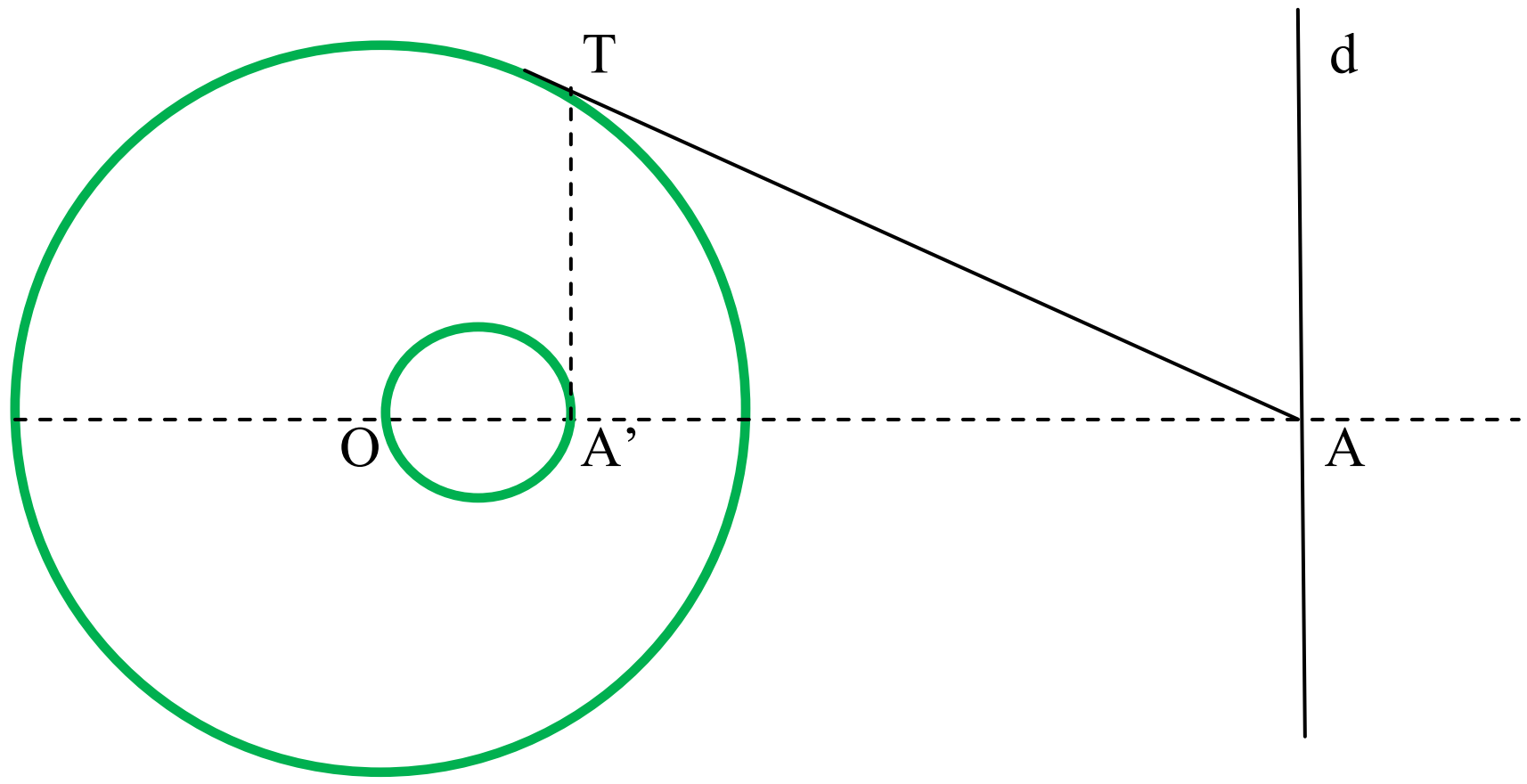

Fig.5

If the line $d$ is tangent to the circle $\mathcal{C}\left(O,\sqrt{k}\right)$, we know that the points of the circle $\mathcal{C}\left(O,\sqrt{k}\right)$ are invariant through the inversion $i_0^k$, therefore if the line $d$ is tangent in $A$ to the circle $\mathcal{C}\left(O,\sqrt{k}\right)$, the point $A$ has as inverse through $i_0^k$ the point $A$.

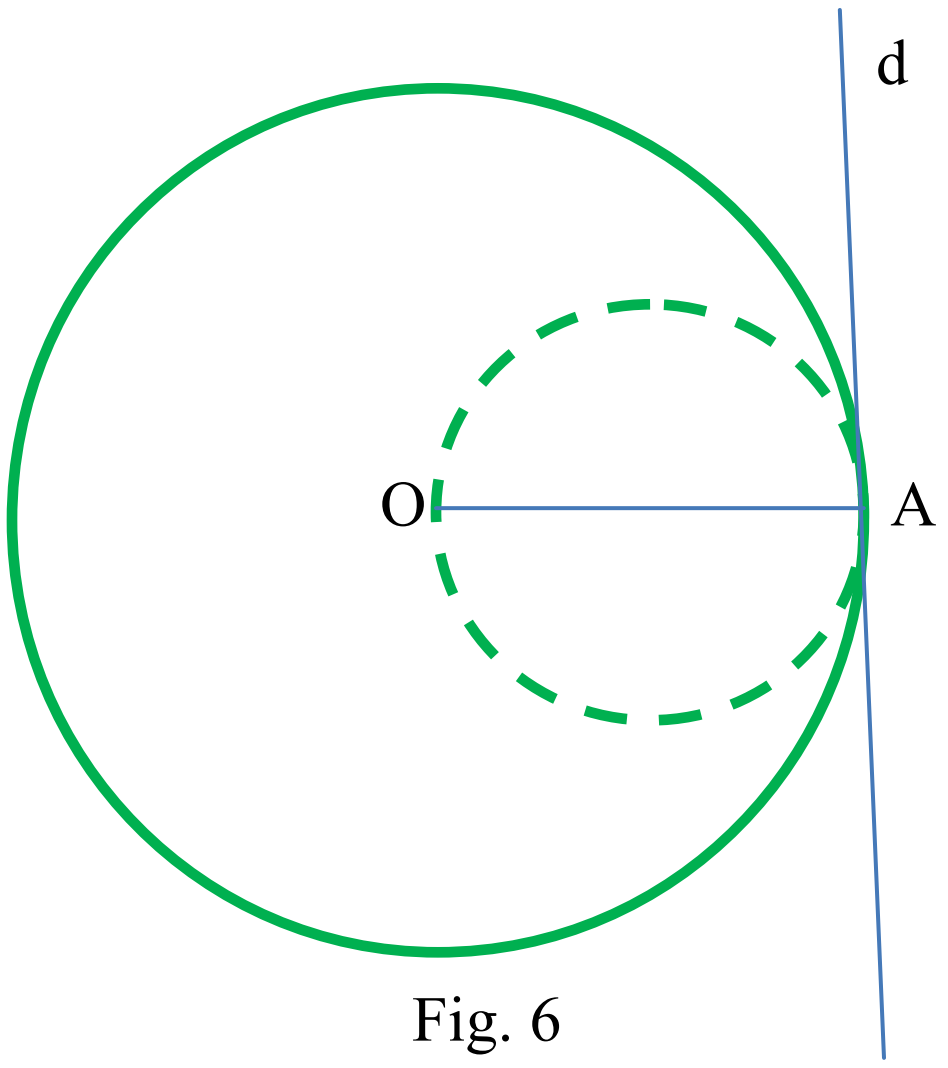

Fig. 6

The image will be the circle of diameter $[OA]$ from which we exclude the point $O$; this circle is tangent interior to the fundamental circle $\mathcal{C}\left(O,\sqrt{k}\right)$.

If the line $d$ is secant to the circle $\mathcal{C}\left(O,\sqrt{k}\right)$ and $O \notin d$, then the image through $i_0^k$ of the line $d$ will be the circumscribed circle to the triangle $OAB$ from which we exclude the point $O$.

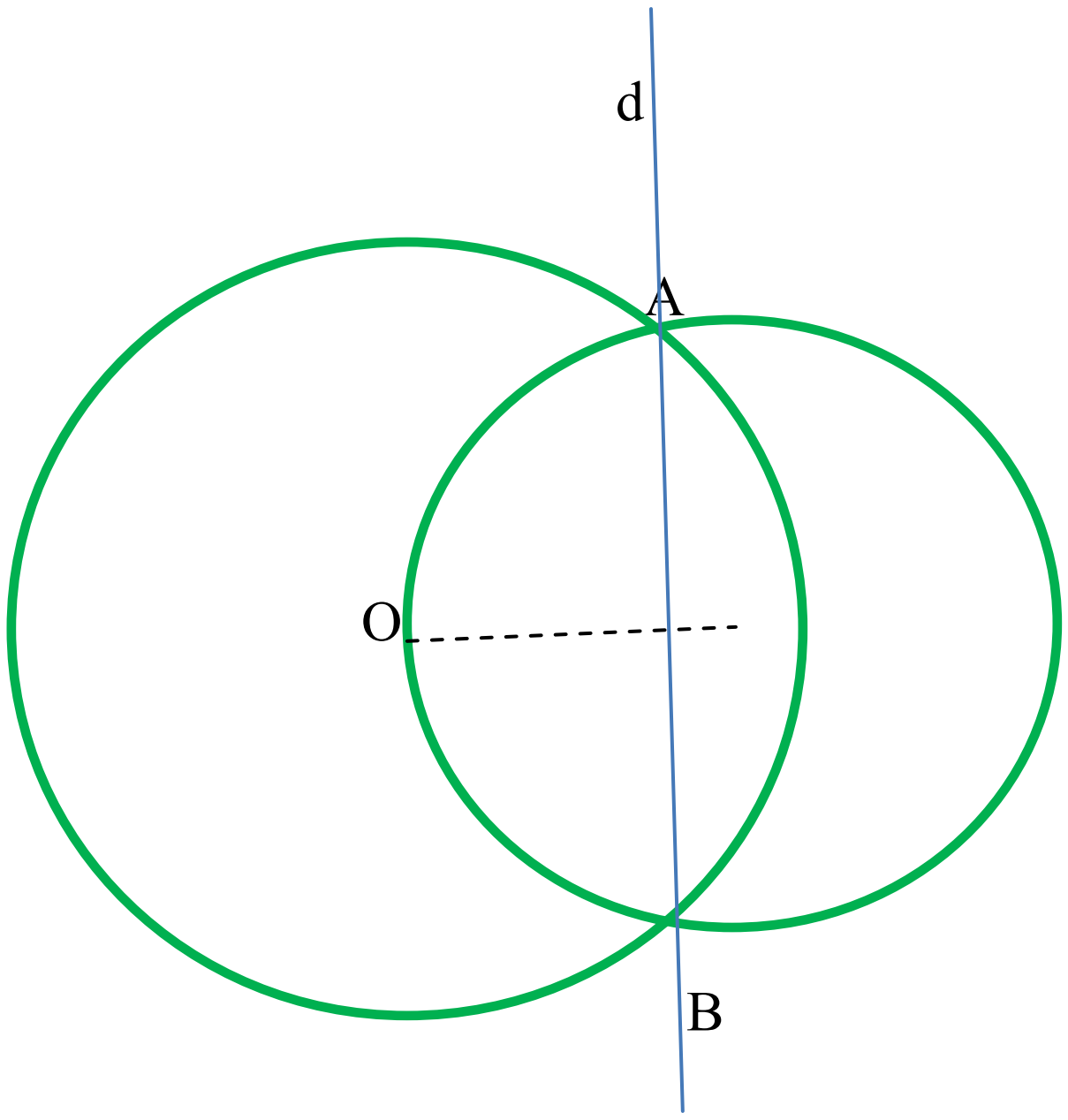

Fig. 7

### 2. The construction of the inverse of a circle

If the circle $\mathcal{C}(O_1, r_1)$ passes through $O$ and it is interior to the fundamental circle $\mathcal{C}(O, \sqrt{k})$, we construct the diametric point $A'$ of the point $O$ in the circle $\mathcal{C}(O_1, r_1)$. We construct the tangent in $A'$ to the circle $\mathcal{C}(O_1, r_1)$ and we note with $T$ one of its points of intersection with $\mathcal{C}(O, \sqrt{k})$. We construct the tangent in $T$ to the circle $\mathcal{C}(O, \sqrt{k})$ and we note with $A$ its intersection with the line $OA'$.
We construct the perpendicular in $A$ on $OA'$; this perpendicular is the image of the circle $\mathcal{C}(O_1, r_1) / \{O\}$ through $i_0^k$ (see figure 5).

If the circle $\mathcal{C}(O_1, r_1)$ passes through $O$ and it is tangent in interior to the circle $\mathcal{C}(O, \sqrt{k})$. The image through $i_0^k$ of the circle $\mathcal{C}(O_1, r_1) / \{O\}$ is the common tangent of the circles $\mathcal{C}(O, \sqrt{k})$ and $\mathcal{C}(O_1, r_1)$.

If the circle $\mathcal{C}(O_1, r_1)$ passes through $O$ and it is secant to the circle $\mathcal{C}(O, \sqrt{k})$, the image through $i_0^k$ of the circle $\mathcal{C}(O_1, r_1)$ from which we exclude the point $O$ is the common secant of the circles $\mathcal{C}(O, \sqrt{k})$ and $\mathcal{C}(O_1, r_1)$.

If the circle $\mathcal{C}(O_1, r_1)$ is secant to the circle $\mathcal{C}(O, \sqrt{k})$ and it does not passes through $O$, we'll note with $A, B$ the common points of the circles $\mathcal{C}(O_1, r_1)$ and $\mathcal{C}(O, \sqrt{k})$.

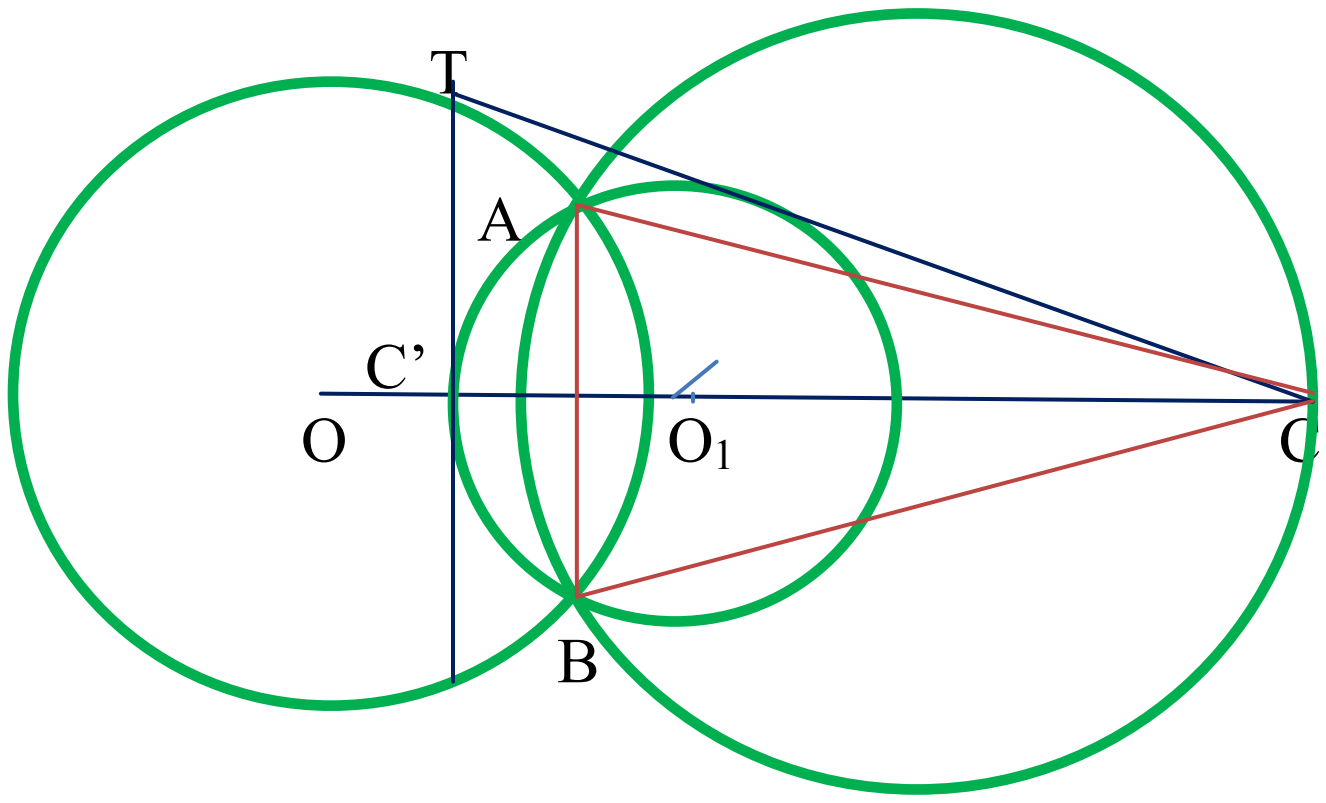

Fig. 8

Let $\{C'\}=(OO_1)\cap\mathcal{C}(O_1,r_1)$. We construct the tangent in $C'$ to the circle $\mathcal{C}(O_1,r_1)$ and we note with $T$ one of its intersection points with $\mathcal{C}(O,\sqrt{k})$. We construct the tangent in $T$ to $\mathcal{C}(O,\sqrt{k})$ and we note with $C$ the intersection of this tangent to the line $OO'$.

We construct the circumscribed circle to the triangle $ABC$. This circle is the image through $i_0^k$ of the circle $\mathcal{C}(O_1,r_1)$ (see figure 8).

If the circle $\mathcal{C}(O_1,r_1)$ is tangent interior to the circle $\mathcal{C}(O,\sqrt{k})$ and it does not passes through $O$. Let $A$ the point of tangency of the circles. We note $\{A'\}=(OA)\cap\mathcal{C}(O_1,r_1)$. We construct the tangent in $A'$ to the circle $\mathcal{C}(O_1,r_1)$ and we note with $T$ one of the intersection points of the tangent with $\mathcal{C}(O,\sqrt{k})$. We construct the tangent in $T$ to the circle $\mathcal{C}(O,\sqrt{k})$ and we note with $A''$ its the intersection with the line $OO_1$. We construct the circle of diameter $[AA'']$, this circle is the inverse of the circle $\mathcal{C}(O_1,r_1)$ through $i_0^k$.

If the circle $\mathcal{C}(O_1,r_1)$ is tangent in the exterior to the circle $\mathcal{C}(O,\sqrt{k})$. Let $A$ the point of tangency of the circles, we construct $A'$ the diametric of $A$ in the circle $\mathcal{C}(O_1,r_1)$, we construct the tangent $A'T$ to the circle $\mathcal{C}(O,\sqrt{k})$ and then we construct the orthogonal projection $A''$ of the point $T$ on $OA$. We construct the circle with the diameter $[AA'']$, this circle is the inverse of the circle $\mathcal{C}(O_1,r_1)$ through $i_0^k$.

If the circle $\mathcal{C}(O_1,r_1)$ is exterior to the circle $\mathcal{C}(O,\sqrt{k})$. Let $\{A,B\}=(OO_1)\cap\mathcal{C}(O_1,r_1)$, we construct the tangents $AT,BP$ to $\mathcal{C}(O,\sqrt{k})$, then we construct the projections $A',B'$ of the point $T$ respectively $P$ on $OO_1$. The circle of diameter $[A'B']$ will be the circle image through $i_0^k$ of the circle $\mathcal{C}(O_1, r_1)$ (see figure 9).

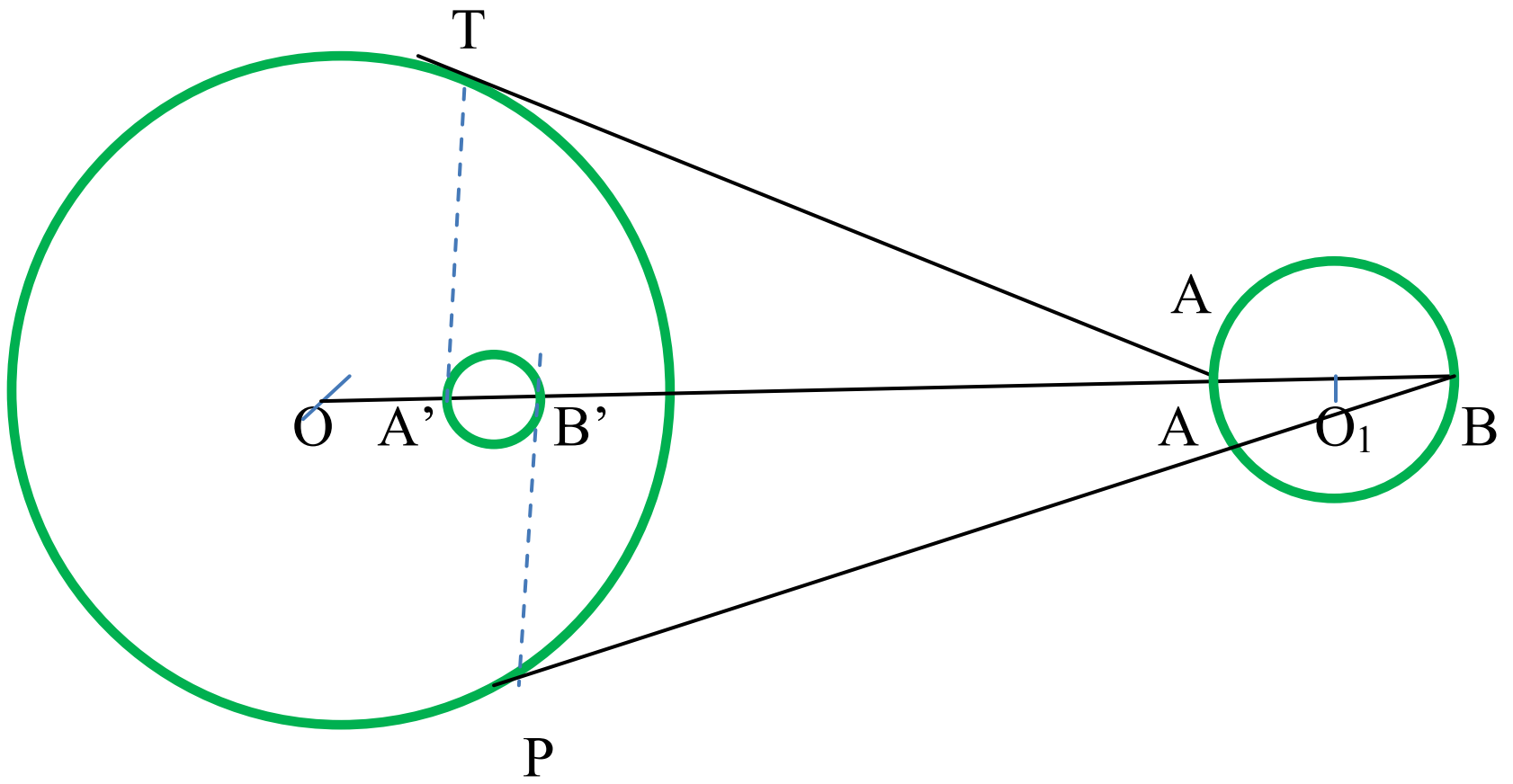

Fig. 9

**Remark 5**

If the circle $\mathcal{C}(O_1, r_1)$ is concentric with $\mathcal{C}(O, \sqrt{k})$, then also its image through the $i_0^k$ will be a concentric circle with $\mathcal{C}(O, \sqrt{k})$.

## D. Other properties of the inversion

**Property 1**

If $M, N$ are two non-collinear points with the pole $O$ of the inversion $i_0^k$ and which are not on the circle $\mathcal{C}(O, \sqrt{k})$, then the points $M, N$, $i_0^k(M)$, $i_0^k(N)$ are concyclic and the circle on which these are situated is orthogonal to the circle $\mathcal{C}(O, \sqrt{k})$

**Proof**

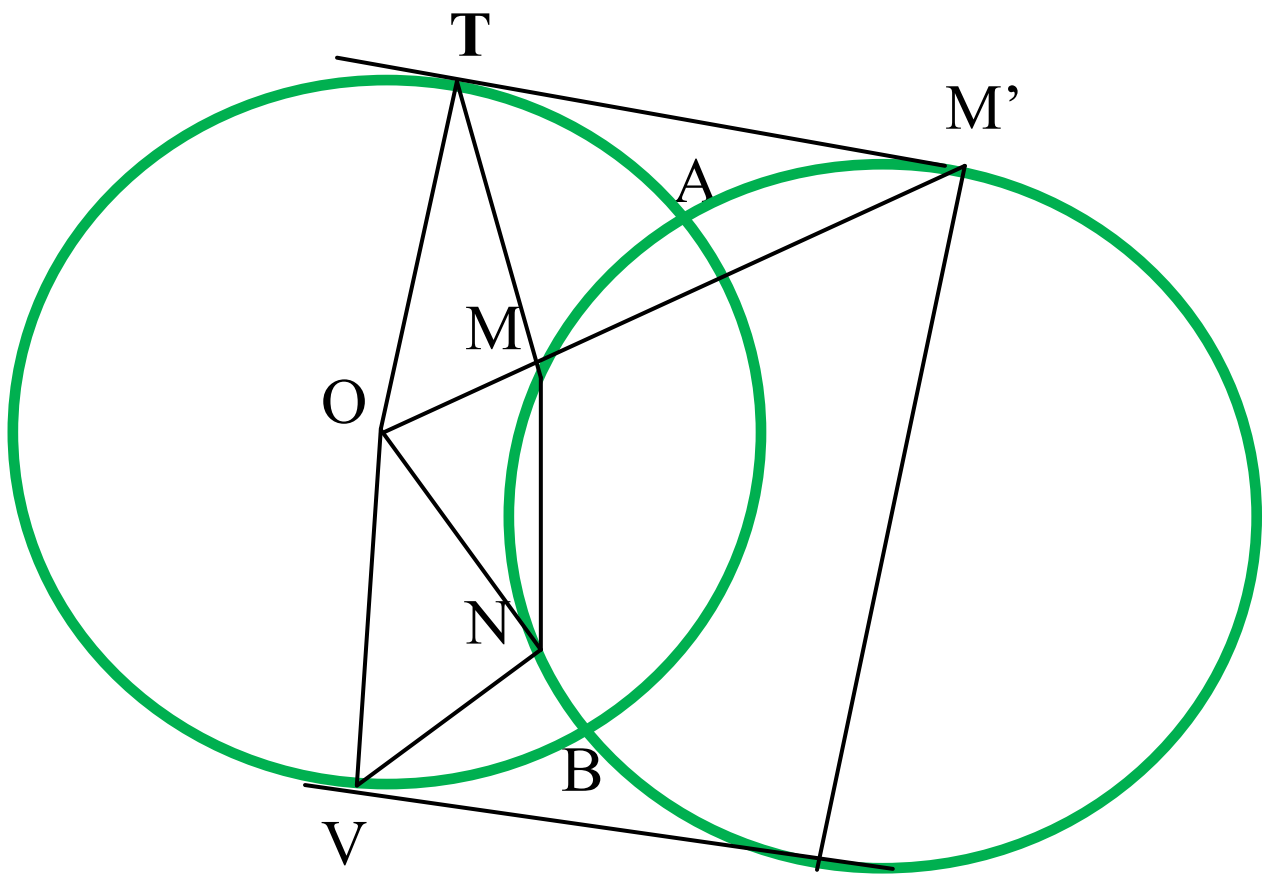

N’

Fig. 10

Let $M, N$ in the interior of the circle $\mathcal{C}(O, \sqrt{k})$ (see figure 10). We construct $M' = i_0^k(M)$ and $N' = i_0^k(N)$. We have $OM \cdot OM' = ON \cdot ON' = k$. It results $\frac{ON}{OM} = \frac{ON'}{OM'}$, which along with $\sphericalangle MON \equiv \sphericalangle N'OM'$ shows that the triangles $OMN, ON'M'$ are similar. From this similarity, it results that $\sphericalangle OMN \equiv \sphericalangle ON'M'$, which show that the points $M, N, N', M'$ are concyclic . If we note with $A, B$ the intersection points of the circles $\mathcal{C}(O, \sqrt{k})$ with that formed by the points $M, N, N', M'$, and because $OM \cdot OM' = k$, it results $OA^2 = OM \cdot OM'$, therefore $OA$ is tangent to the circle of the points $M, N, N', M'$, which shows that this circle is orthogonal to the fundamental circle of the inversion $\mathcal{C}(O, \sqrt{k})$.

**Property 2.**

If $M, N$ are two points in plane and $M', N'$ their inversion through the positive inversion $i_0^k$, then

$$M'N' = k \frac{MN}{OM \cdot ON}$$

**Proof**

We observed that the triangles $OMN$ and $ON'M'$ are similar (see figure 10), therefore

$$\frac{M'N'}{NM} = \frac{OM'}{ON}$$

It results that

$$\frac{M'N'}{NM} = \frac{OM' \cdot OM}{OM \cdot ON} = \frac{k}{OM \cdot ON}$$

And from here

$$M'N' = k \frac{MN}{OM \cdot ON}$$

**Definition 3**

The angle of two secant circles in the points $A, B$ is the angle formed by the tangents to the two circles constructed in the point $A$ or the point $B$.

**Observation 1**

If two circles are orthogonal, then their angle is a right angle. If two circles are tangent, then their angle is null.

**Definition 4**

The angle between a secant line to a circle in $A, B$ and the circle is the angle formed by the line with one of the tangents constructed in $A$ or in $B$ to the circle.

**Observation 2**

If a line contains the center of a circle, its angle with the given circle is right. If a line is tangent to a circle, its angle with the circle is null.

**Theorem**

Through an inversion the angle between two lines is preserved, a line and a circle, two circles.

**Proof**

If the lines pass through the pole $O$ of the inversion, because these are invariant, their angle will remain invariant. We saw that if the lines do not pass through the pole $O$ inversion, their images through $i_0^k$ are two circles which pass through the point $O$ and which have the diameters constructed through $O$ perpendicular on $d_1, d_2$.

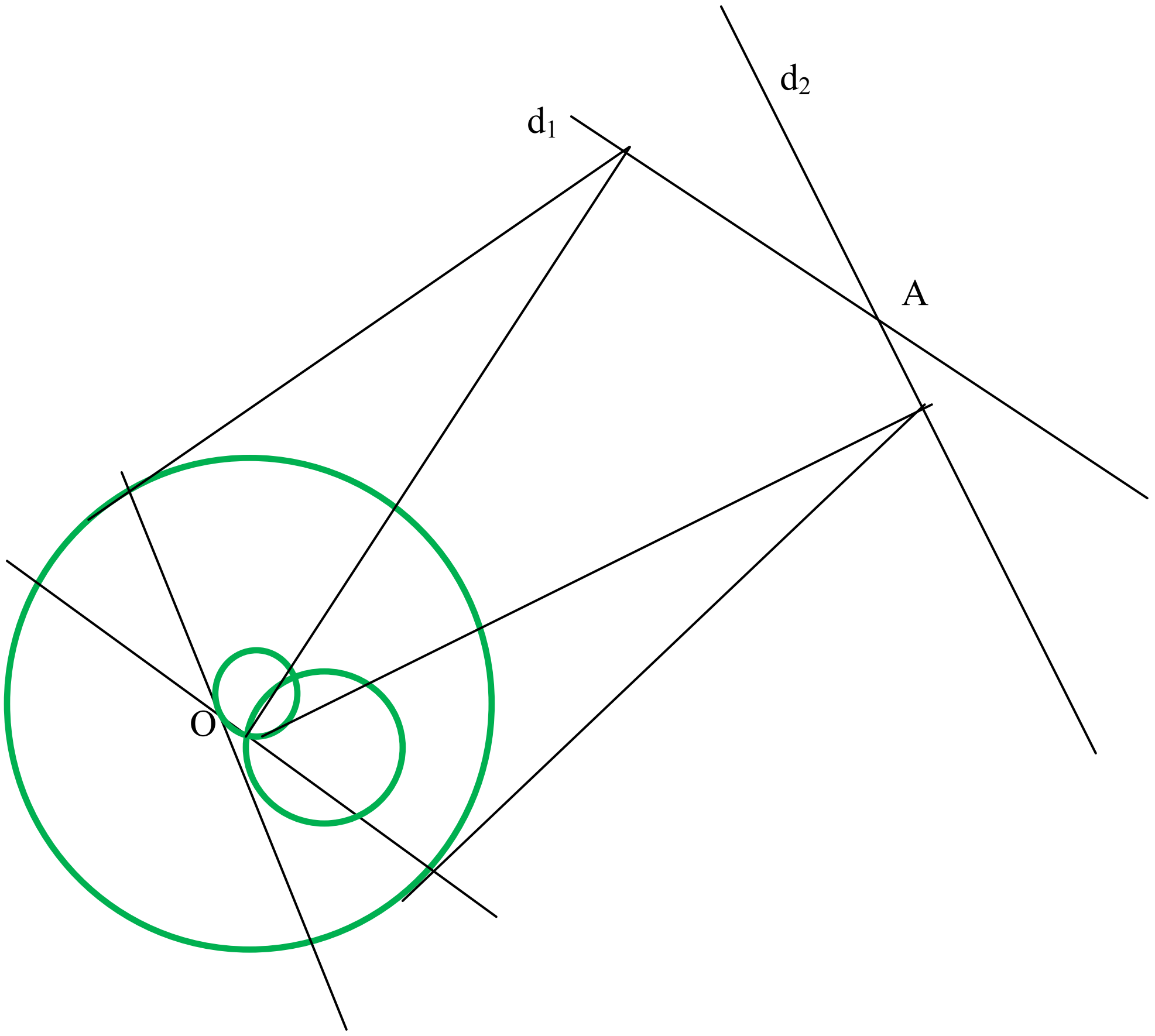

Fig. 11

The angle of the image-circles of the lines is the angle formed by the tangents constructed in $O$ to theses circles; because these tangents are perpendicular on the circles' diameters that pass through $O$, it means that these are parallel with the lines $d_1, d_2$ and therefore their angle is congruent with the angle of the lines $d_1, d_2$.

If a line passes through the inversion's pole and the other does not, the theorem is proved on the same way.

If a line passes through the inversion's pole and the circle secant with the line, then the line's image is that line and the arch's imagine is the given circle through the homothety of center $O$, the inversion's pole.

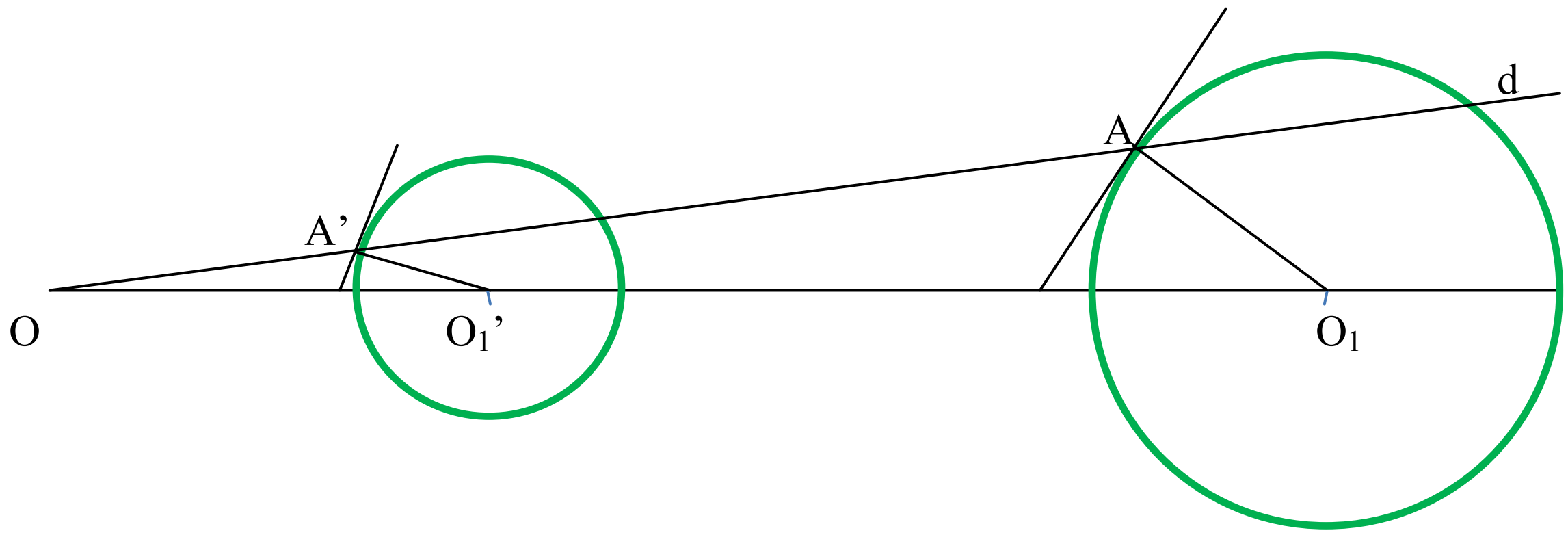

Fig. 12

We'll note with $A$ one of the intersections of the given line $d$ with the given circle $\mathcal{C}(O_1, r)$ and with $\mathcal{C}(O_1', r')$ the image circle of the given circle through the image $i_0^k$, we'll have $O_1'A' \parallel O_1A$ and the angle between $d$ and $\mathcal{C}(O_1, r)$ equal to the angle between $\mathcal{C}(O_1', r')$ and $d$ as angle with their sides parallel (see figure 12) .

The case of the secant circles which do not contain the inversion center is treated similarly as the precedent ones.

**Remark 6**

a) The property of the inversion to preserve the angles in the sense that the angle of two curves is equal with the angle of the inverse curves in the inverse common point suggests that the inversion be called conform transformation.
b) The setoff all homotheties and of the inversions of the plane of the same center form an algebraic group structure. This group of the inversions and homotheties of the plane of the same pole $O$ is called the conform group of center $O$ of the plane . The set of the lines and circles of the plane considered in an ensemble is invariant in rapport with the group's transformations conform in the sense that a line of the group or a circle are transformed also in a line or a circle
c) Two orthogonal circles which don't pass through the pole of the inversion are transformed through that inversion in two orthogonal circles.
d) Two circles tangent in a point will have as inverse two parallel lines through a pole inversion – their point of tangency.

**Applications**

**1.** If $A, B, C, D$ are distinct points in plane, then

$$AC \cdot BD \le AB \cdot CD + AD \cdot BC$$

**(The Ptolomeus Inequality)**

**Proof**

We'll consider the inversion $i_A^k$, $k > 0$, and let $B', C', D'$ the images of the points $B, C, D$ through this inversion.

We have

$$AC \cdot BD \leq AB \cdot AB' = AC \cdot AC' = AD \cdot AD' = k$$

Also

$$B'C' = \frac{k \cdot BC}{AB \cdot AC}, \quad C'D' = \frac{k \cdot CD}{AC \cdot AD}, \quad D'B' = \frac{k \cdot DB}{AD \cdot AB}$$

Because the points $B', C', D'$ determine, in general a triangle, we have

$$D'B' \leq B'C' + C'D'.$$

Taking into consideration this relation and the precedents we find the requested relation.

**Observation**

The equality in the Ptolomeus is achieved if the points $A, B, C, D$ are concyclic. The result obtained in this case is called the I theorem of Ptolomeus.

**2. Feuerbach's theorem** (1872)

Prove that the circle of nine points is tangent to the ex-inscribed and inscribed circles to a given triangle.

**Proof**

The idea we'll use for proving this theorem is to find an inversion that will transform the circle of nine points in a line and the inscribed and ex-inscribed circles tangent to the side $BC$ to be invariant. Then we show that the imagine line of the circle of nine points is tangent to these circles.

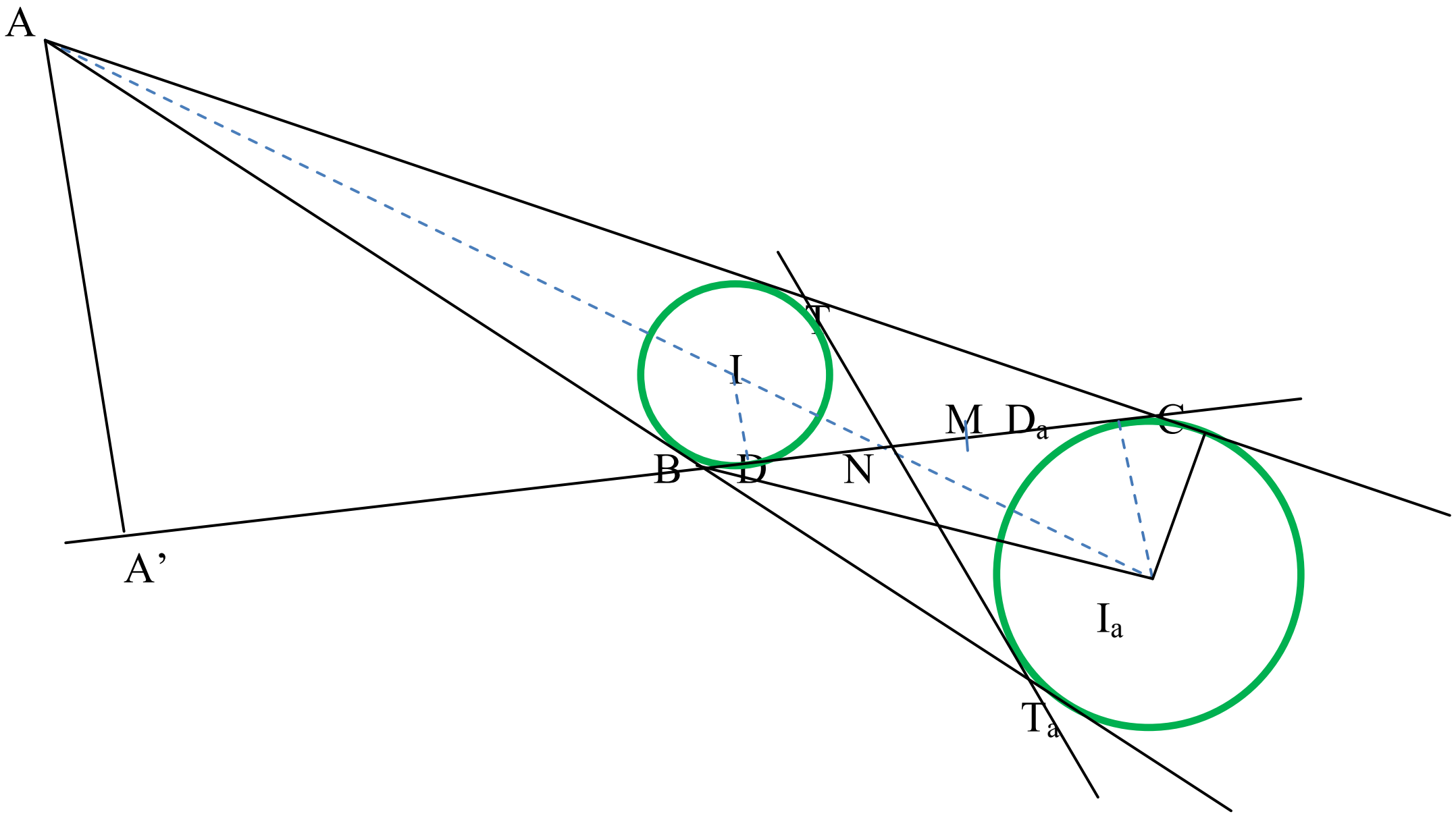

Fig. 13

Let $A'$ the projection of $A$ on $BC$. $D$ and $D_a$ the projections of $I, I_a$ on $BC$, $M$ the middle of $(BC)$ and $N$ the intersection of the bisector $(AI$ with $(BC)$. See figure 12.

It is known that the points $D, D_a$ are isotonic, and we find:

$$MD = MD_a = \frac{(b-c)}{2}$$

Without difficulties we find $MD^2 = MA' \cdot MN$.

Considering $i_M^{\left(\frac{(b-c)}{2}\right)^2}$, we observe that the inverse of $A'$ through this inversion is $N$. Therefore, the circle of nine points, transforms in a line which passes through $N$ and it is perpendicular on $MO_9$, $O_9$ being the center of the circle of nine points, that is the middle of the segment $(ON)$. Because $MO_9$ is parallel to $AO$, means that the perpendicular on $MO_9$ will have the direction of the tangent in $A$ to the circumscribed circle to the triangle $ABC$, in other words of a antiparallel to $BC$ constructed through $N$, and this is exactly the symmetric of the line $BC$ in rapport with the bisector $AN$, which is the tangent $TT_a$ to. The inscribed and ex-inscribed circles tangent to $BC$ remain invariant through the considered inversion, and the inverse of the circle of the nine points is the tangent $TT_a$ of these circles.

This property being true after inversion also, it result that the Euler's circle is tangent to the inscribed circle and ex-inscribed circle $(I_a)$. Similarly, it can be proved that the circle of the nine points is tangent to the other ex-inscribed circles.

**Note**

The tangent points of the circle of the nine points with the ex-inscribed and inscribed circles are called the Feuerbach points.

# Chapter 7

## Solutions and indications to the proposed problems

**1.**

a) The triangles $AA_1D_1$, $CB_1C_1$ are homological because the homological sides intersect in the collinear points $B, D, P$ (we noted $\{P\} = A_1D_1 \cap B_1C_1 \cap BD$).

b) Because $\{P\} = A_1D_1 \cap B_1C_1 \cap BD$, it result that the triangles $DC_1D_1$, $BB_1A_1$ are homological, then the homological sides $DC_1$, $BB_1$; $DD_1$, $BA_1$; $D_1C_1$, $A_1B_1$ intersect in three collinear points $C, A, Q$ where $\{Q\} = D_1C_1 \cap A_1B_1$.

2.

i) Let $O_1, O_2, O_3$ the middle points of the diagonals $(AC),(BD),(EF)$. These points are collinear – the Newton-Gauss line of the complete quadrilateral $ABCDEF$. The triangles $GIN, ORP$ have as intersections as homological sides the collinear points $O_1, O_2, O_3$, therefore these are homological.

ii) $GI \cap JK = \{O_1\}$, $GH \cap JL = \{O_2\}$, $HI \cap KL = \{O_3\}$ and $O_1, O_2, O_3$ collinear show that $GIH, JKL$ are homological

iii) Similar with ii).

iv) We apply the theorem “If three triangles are homological two by two and have the same homological axis, then their homology centers are collinear

v) Similar to iv).

**3.**

i) The Cevians $AA_2, BB_2, CC_2$ are the isotomic of the concurrent Cevians $AA_1, BB_1, CC_1$, therefore are concurrent. Their point of concurrency is the isotomic conjugate of the concurrence point of the Cevians $AA_1, BB_1, CC_1$.

ii) We note $M_aM_bM_c$ the medial triangle of the triangle $A_1B_1C_1$ ($M_a$ - the middle point of $(B_1C_1)$, etc.) and $\{A'\} = AM_a \cap BC$, $\{B'\} = BM_b \cap AC$, $\{C'\} = AM_c \cap AB$ .

We have

$$\frac{BA'}{CA'} = \frac{c \cdot \sin \sphericalangle BAA'}{b \cdot \sin \sphericalangle CAA'}.$$

Because $Aria_{\Delta AC_1M_a} = Aria_{\Delta AB_1M_a}$ we have

$$AC_1 \cdot \sin \sphericalangle BAA' = AB_1 \cdot \sin \sphericalangle CAA',$$

therefore

$$\frac{\sin \sphericalangle BAA'}{\sin \sphericalangle CAA'} = \frac{AB_1}{AC_1} = \frac{z}{c-y}.$$

We noted $x = A_1C_1$, $y = BC_1$, $z = AB_1$,

and therefore

$$\frac{BA'}{CA'}=\frac{c}{b}\cdot\frac{z}{c-y}.$$

Similarly we find

$$\frac{CB'}{AB'}=\frac{a}{c}\cdot\frac{y}{a-x},\ \frac{AC'}{BC'}=\frac{b}{a}\cdot\frac{x}{b-z}.$$

We have

$$\frac{BA'}{CA'}\cdot\frac{CB'}{AB'}\cdot\frac{AC'}{BC'}=\frac{x\cdot y\cdot z}{(a-x)(b-y)(c-z)}.$$

But the triangles $ABC$, $A_1B_1C_1$ are homological,

therefore

$$\frac{x\cdot y\cdot z}{(a-x)(b-y)(c-z)}=1,$$

and therefore

$$\frac{BA'}{CA'}\cdot\frac{CB'}{AB'}\cdot\frac{AC'}{BC'}=1.$$

In conformity with Ceva's theorem it result that $AM_a$, $BM_b$, $CM_c$ are concurrent.

iii) Similar to ii)

**4**

i) $\sphericalangle B_1OC_1=\sphericalangle C_2OA_2=120^\circ$,

it results that $\sphericalangle B_1OC_2\equiv\sphericalangle C_1OA_2$.

Similarly $\sphericalangle C_1OA_2\equiv\sphericalangle A_1OB_2$. It result $\Delta A_1OA_2\equiv\sphericalangle B_1OC_2\equiv\sphericalangle C_1OA_2$.

ii) $\Delta A_1OA_2\equiv\sphericalangle B_1OC_2\equiv\sphericalangle C_1OA_2$ (S.A.S).

iii) $\Delta A_1B_1B_2\equiv\Delta B_1C_1C_2\equiv\Delta C_1A_1A_2$ (S.S.S) from here we retain that $\sphericalangle A_1B_2B_1\equiv\sphericalangle B_1C_2C_1\equiv\sphericalangle A_1A_2C_1$, and we obtain that $\Delta C_1A_2B_2\equiv\Delta B_1C_2A_2\equiv\Delta A_1B_2C_1$ (S.A.S).

iv) $\Delta A_3B_1C_2\equiv\Delta B_3C_1A_2\equiv\Delta C_3A_1B_2$ (A.S.A). It result that $\sphericalangle A_3\equiv\sphericalangle B_3\equiv\sphericalangle C_3$, therefore $A_3B_3C_3$ is equilateral, $\Delta OB_1A_3\equiv\Delta OC_1B_3\equiv\Delta OA_1C_3$ (S.A.S). It result that $OA_3=OB_3=OC_3$, therefore $O$ is the center of the equilateral triangle $A_3B_3C_3$.

v) For this we'll apply the D. Barbilian's theorem.

**5**.

We apply the Pascal's theorem for the degenerate hexagon $BCCDDE$.

**6**.

See http://vixra.org/abs/1103.0035 (Ion Pătraşcu, Florentin Smarandache, "A Property of the Circumscribed Octagon" – to appear in *Research Journal of Pure Algebra*, India)

**7.**

If $O$ is the center of the known circle, we construct firstly the tangents $OU, OV$

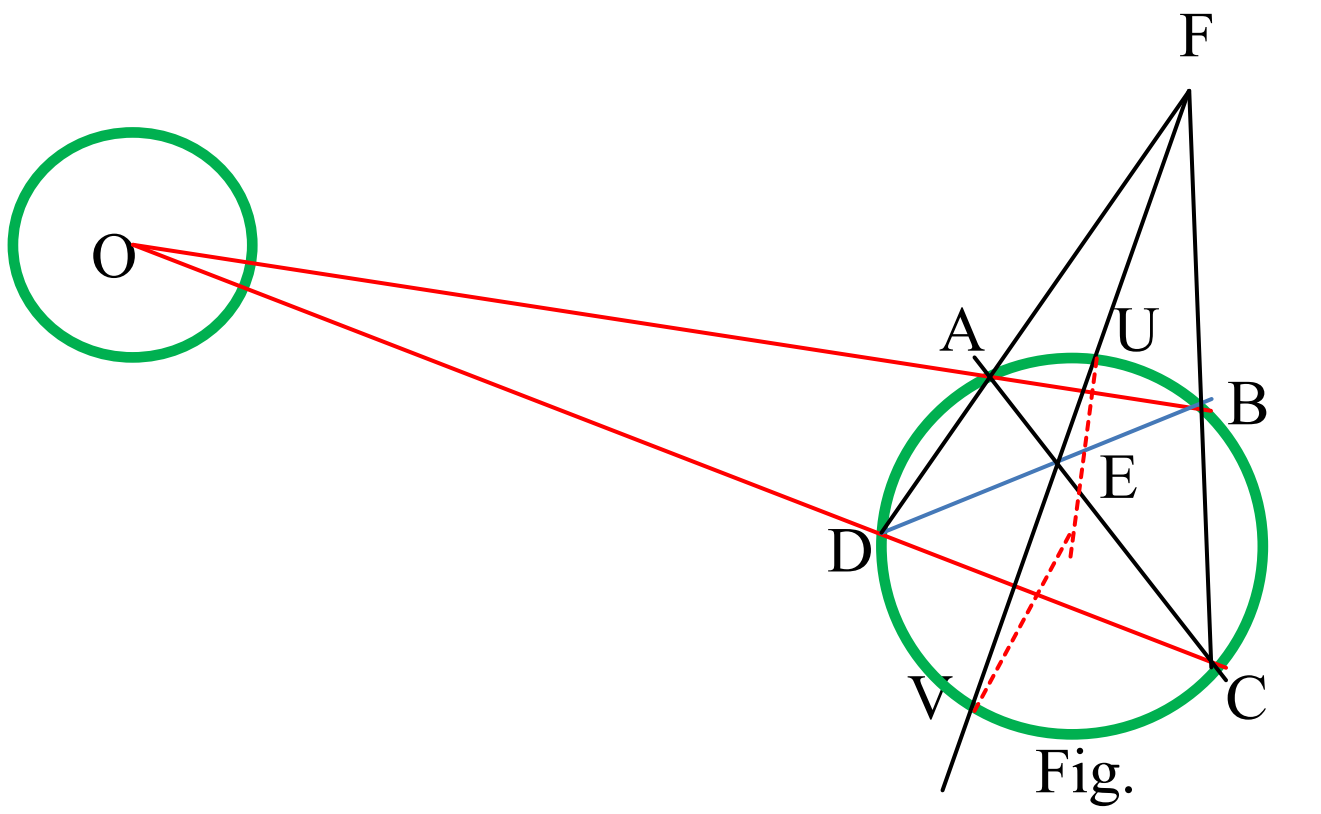

Fig.

to the circle for which we don't know its center in the following manner:

We construct the secant $O, A, B$ and $O, C, D$. We construct $\{E\} = AC \cap BD$, $\{F\} = AD \cap BC$. Construct $EF$ and we note $U, V$ its intersections with the circle (see figure.)

We, practically, constructed the polar of the point $O$ in rapport to the circle whose center $O'$ we do not know. It is $UV$, and as it is known it is perpendicular on $OO'$. To obtain $O'$ we'll construct the perpendiculars in $U, V$ on $OU, OV$. The intersection of these perpendiculars being the point $O'$.

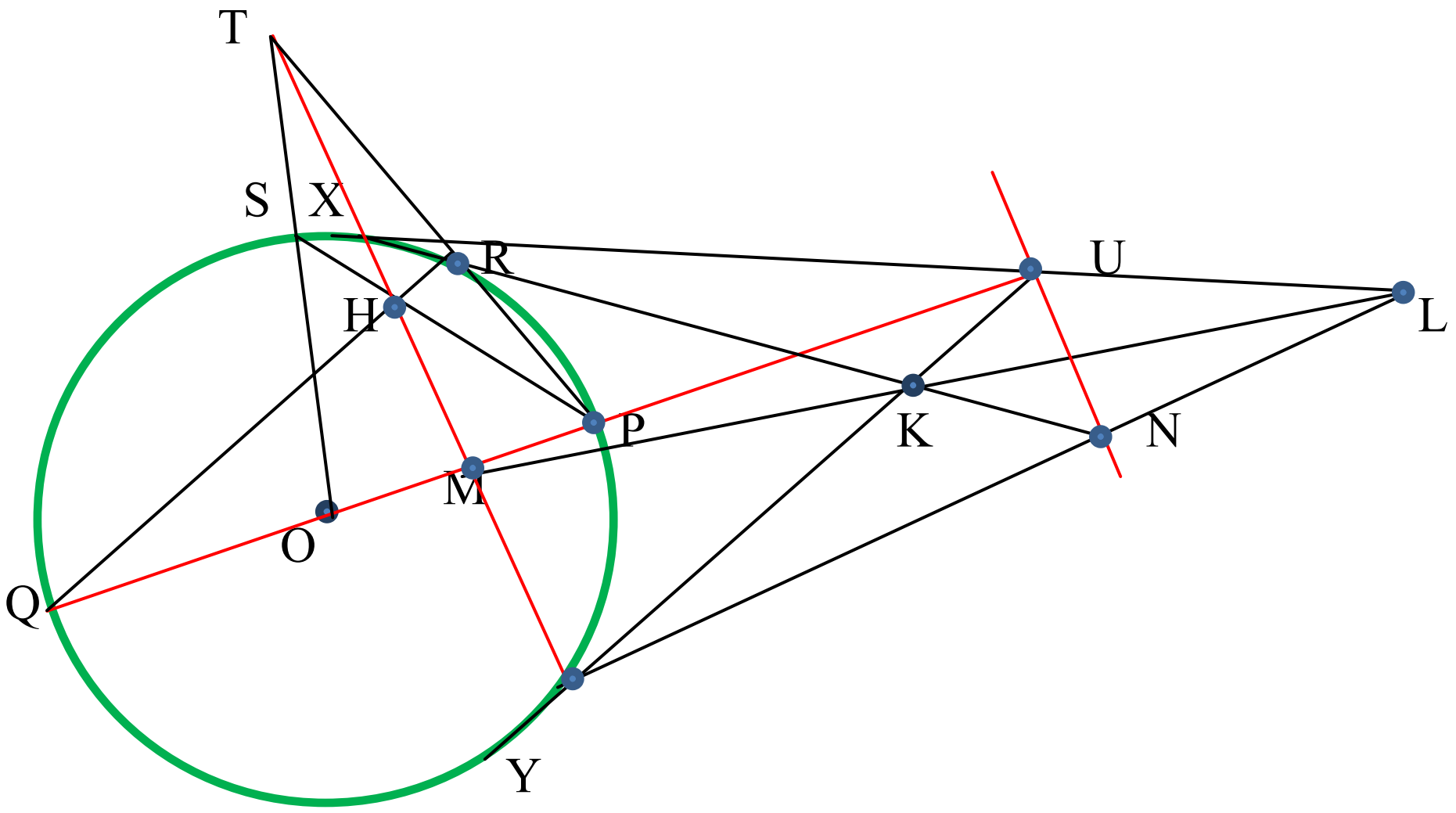

Fig.

We'll construct a perpendicular on $OU$. We note $P, Q$ the intersection of $OU$ with the circle $\mathcal{C}(O)$. We consider a point $R$ on $\mathcal{C}(O)$, we construct $PR, QR$ and consider a point $T$ on $PR$ such that $(OT \cap \mathcal{C}(O) = \{S\}$. We connect $S$ and $P$, note $\{H\} = RQ \cap PS$. Connect $T$ and $H$ and note $X, Y$ the intersection points of $TH$ with the circle $\mathcal{C}(O)$. We have that $TH \perp OU$ because in the triangle $PTQ$, the point $H$ is its orthocenter.

Through $U$ we'll construct a parallel to $XY$. We note $\{M\} = XY \cap OU$, connect $X$ with $U$, consider $L$ on $XU$, connect $L$ with $M$, $Y$ with $U$, note $\{K\} = YU \cap LM$ and $\{N\} = XK \cap LY$. We have $UN \parallel XY$ and therefore $UN$ will contain $O'$.

We will repeat this last construction for a tangent $OV$ and we construct then the perpendicular in $V$ on $OV$, let it be $VG$. The intersection between $UN$ and $VG$ is $O'$, which is the center which needed to be constructed.

**8.**

Let $k = \frac{\overrightarrow{A_1B}}{\overrightarrow{A_1C}}$, $p = \frac{\overrightarrow{B_1C}}{\overrightarrow{C_1A}}$, $q = \frac{\overrightarrow{C_1A}}{\overrightarrow{C_1B}}$. We have

$$\overrightarrow{AA_1} = \frac{1}{1+k}\left(\overrightarrow{AB} + k\overrightarrow{AC}\right);\ \overrightarrow{BB_1} = \frac{1}{1+p}\left(\overrightarrow{BC} + p\overrightarrow{BA}\right);\ \overrightarrow{CC_1} = \frac{1}{1+q}\left(\overrightarrow{CA} + q\overrightarrow{CB}\right)$$

We obtain

$$\overrightarrow{AB}(1+p)(1+q) + \overrightarrow{BC}(1+k)(1+q) + \overrightarrow{CA}(1+k)(1+p) + k\overrightarrow{AC}(1+p)(1+q) +$$
$$+ p\overrightarrow{BA}(1+k)(1+q) + q\overrightarrow{BC}(1+k)(1+p) = \vec{0}$$

After computations we have

$$\overrightarrow{AB}(1+p)(1+q-p-pk) + \overrightarrow{BC}(1+k)(1-qp) + \overrightarrow{AC}(1+p)(kq-1) = \vec{0}$$

But $\overrightarrow{AC} = \overrightarrow{AB} + \overrightarrow{BC}$, it results

$\overrightarrow{AB}(1+p)(1+k)(q-p) + \overrightarrow{BC}(k-p)(1+q) = \vec{0}$. We must have

$$\begin{cases} (1+p)(1+k)(q-p) = 0 \\ (k-p)(1+q) = 0 \end{cases}$$

This will take place every time when $p = q = k = -1$ which show that $AA_1, BB_1, CC_1$ are medians and that the homological center of the given triangle is $G$, which is the weight center of the triangle $ABC$.

**9.**

We note $m\sphericalangle(UAB) = \alpha$, $m\sphericalangle(ABV) = \beta$, $m\sphericalangle(ACW) = \gamma$ (see figure …)

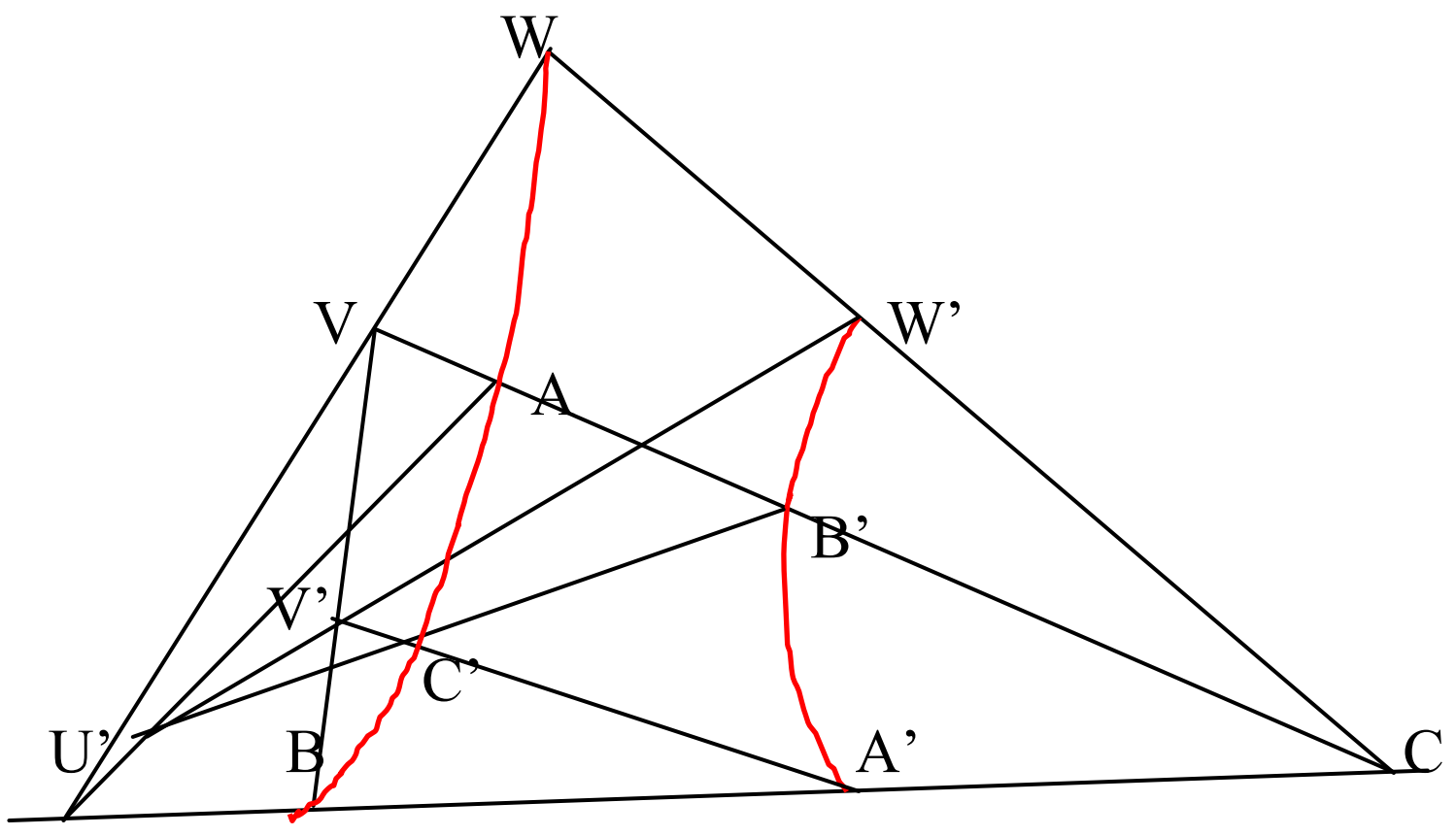

U

Fig.

We have $\frac{UB}{UC}=\frac{Aria\Delta UAB}{Aria\Delta UAC}$, therefore

$$\frac{UB}{UC}=\frac{AB\sin\alpha}{AC\sin(A+\alpha)} \quad (1)$$

$$\frac{U'B'}{U'C'}=\frac{Aria\Delta U'AB'}{Aria\Delta U'AC'}=\frac{AB'\sin(A+\alpha)}{AC'\sin\alpha}$$

Taking into account (1) we find

$$\frac{U'B'}{U'C'}=\frac{AB'}{AC'}\cdot\frac{AB}{AC}\cdot\frac{UC}{UB} \quad (2)$$

Similarly

$$\frac{V'C'}{V'A'}=\frac{BC'}{BA'}\cdot\frac{BC}{AB}\cdot\frac{VA}{VC} \quad (3)$$

$$\frac{W'A'}{W'B'}=\frac{CA'}{CB'}\cdot\frac{CA}{CB}\cdot\frac{WB}{WA} \quad (4)$$

From relations (2), (3), (4) taking into consideration the Menelaus and Ceva's theorems it results that $U',V',W'$ is a transversal in the triangle $A'B'C'$.

**10**.

i) We note $\{A'\}=BC\cap B_1C_1$, $\{A''\}=BC\cap B_2C_2$. In the complete quadrilaterals $CB_1C_1BAA'$, $BC_2B_2CAA''$, the points $B,A_1,C,A'$ respective $C,A_2,B,A''$ are harmonic divisions (see A). The fascicules $(C_1;BA_1M_1B_1)$, $(C_2;B_2M_2A_2B)$ are harmonic and gave the ray $C_1C_2$ in common, then the intersection points $C,A_3,B_3$ of the homological rays $(C_1M_1,C_2M_2),(C_1B_1,C_2B_2),(C_1A_1,C_2A_2)$ are collinear, therefore the side $A_3B_3$ passes through $C$

Similarly, it can be shown that the sides $B_3C_3,C_3A_3$ pass through $A$ respectively $B$.

ii) We'll apply the Pappus' theorem for the non-convex hexagon $C_2A_1AA_2C_1C$ which has three vertexes on the sides $BA,BC$. The opposite sides $(C_2A_1,A_2C_1),(A_1A,C_1C),(AA_2,CC_2)$ intersect in the collinear points $B_4,M_1,M_2$, therefore $B_4$ belongs to the line $M_1M_2$. Similarly we can show that the points $A_4,C_4$ are on the line $M_1M_2$.

iii) Let $\{S_2\}=AC\cap A_3C_3$. Because the fascicle $(B_2;C_2BA_2C)$ is harmonic, it results that $B,A_3,S_2,C_3$ is a harmonic division, therefore $AC_3$ is the polar of the point $A_3$ in rapport with the sides $AB,AC$ of angle $A$. In the complete quadrilateral $C_2B_1B_2C_1AA_4$ the line $AA_4$ is the polar of the point $A_3$ in rapport with the sides $AB,A$ of the angle $A$. It results that the polars $AC_3,AA_4$ coincide, consequently, the point $A_4$ is situated on the line $B_3C_3$ which passes through $A$. The proof is similar for $B_4,C_4$.

iv) We note $\{S_1\} = BC \cap B_3C$ and $\{S_3\} = AB \cap A_3B_3$. The fascicle $(C_1, B_1CA_1B)$ is harmonic, it results that the points $C, A_3, S_3, B_3$ form a harmonic division. This harmonic division has the point $A_3$ in common with the harmonic division $B, A_3, S_2, C_3$, it results that the line $S_2S_3$ passes through the intersection point of the lines $BC$, $B_3C_3$, which is $S_1$, consequently the points $S_1, S_2, S_3$ are collinear.

Considering the triangles $ABC, A_3B_3C_3$ we observe that $S_1, S_2, S_3$ is their homological axis, therefore these are homological, therefore the lines $AA_3, BB_3, CC_3$ are concurrent.

v) Let $\{S_1\} = C_1C_3 \cap B_3C_3$ and $\{A_5\} = C_1B_1 \cap B_3C_3$. $A_5$ is the harmonic conjugate of $A_3$ in rapport with the points $B_1, C_1$. We'll consider the complete quadrilateral $C_1A_1B_1Q_1B_3C_3$; in this quadrilateral the diagonal $A_1Q_1$ intersects the diagonal $C_1B_1$ in the point $A_3'$, which is the harmonic conjugate of $A_5$ in rapport with $C_1B_1$. Therefore, the points $A_3', A_3$ coincide. It results that the triangles $A_3B_3C_3$, $A_1B_1C_1$ are homological, the homological center being the point $Q_1$. Similarly, we can prove that the triangles $A_3B_3C_3$, $A_2B_2C_2$ are homological, their homological center being $Q_2$. The homological centers $Q_1$, $Q_2$, evidently are on the homological line $S_1S_2S_3$ because these are on the polar of the point $A_3$ in rapport with $BC$ and $B_3C_3$, which polar is exactly the line $S_1S_2S_3$, which is the homological axis of the triangles $ABC, A_3B_3C_3$.

**11**.

$A-D-F$ is transversal for the triangle $CBE$. We'll apply, in this case, the result obtained in problem 9.

12.

Let $T_AT_BT_C$ the tangential triangle of triangle $ABC$, and the circumpedal triangle $A'B'C'$ of $G$.

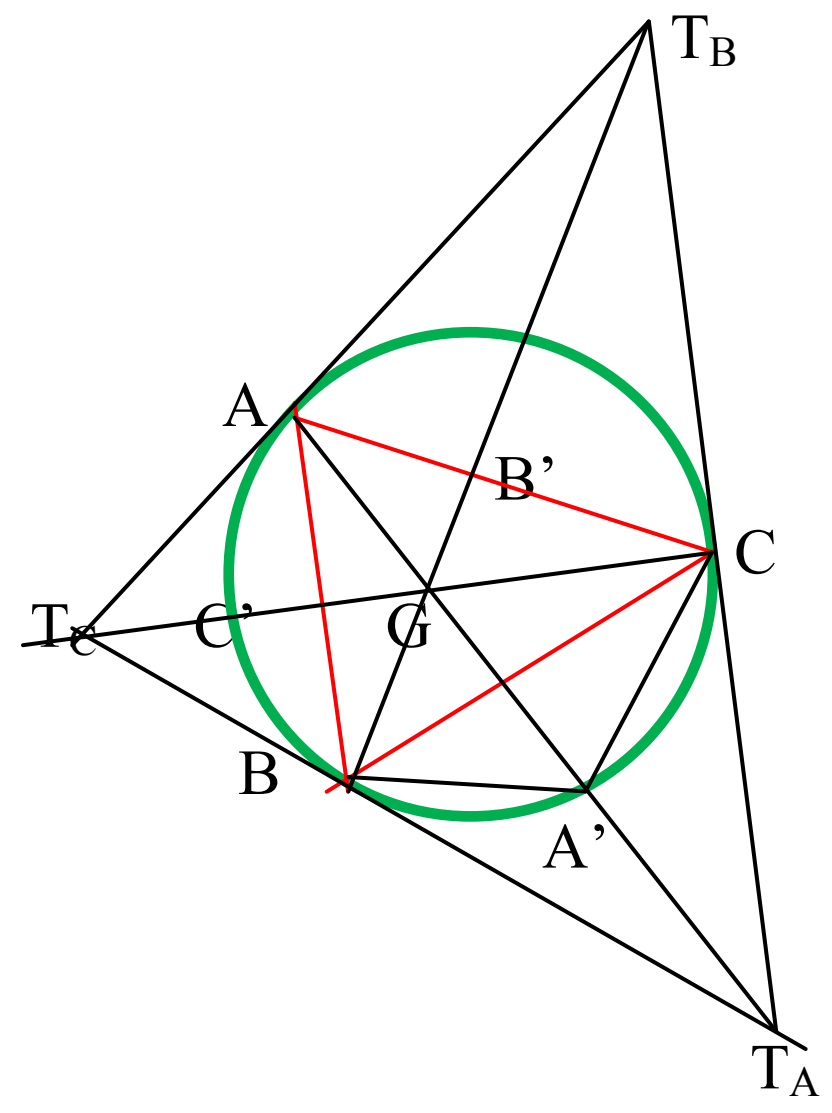

We note $\sphericalangle BAA' = \alpha, \sphericalangle CBB' = \beta, \sphericalangle ACC' = \gamma$, we have

$$\frac{\sin\alpha}{\sin(A-\alpha)}\cdot\frac{\sin\beta}{\sin(B-\beta)}\cdot\frac{\sin\gamma}{\sin(C-\gamma)}=1$$

$$m\sphericalangle T_A BA' = \alpha, m\sphericalangle T_A CA' = A-\alpha$$

The sinus' theorem in the triangles $BA'T_A, CA'T_A$ implies $\frac{T_A A'}{\sin\alpha}=\frac{BA'}{\sin BT_A A'}$, $\frac{T_A A'}{\sin(A-\alpha)}=\frac{CA'}{\sin CT_A A'}$.

We find that

$$\frac{\sin\sphericalangle BT_A A'}{\sin\sphericalangle CT_A A'}\cdot\frac{\sin\alpha}{\sin\sin(A-\alpha)}\cdot\frac{BA'}{CA'}$$

But also the sinus' theorem implies $BA' = 2R\sin\alpha$, $CA' = 2R\sin(A-\alpha)$

Therefore,

$$\frac{\sin\sphericalangle BT_A A'}{\sin\sphericalangle CT_A A'}=\left(\frac{\sin\alpha}{\sin(A-\alpha)}\right)^2$$

Similarly we compute

$$\frac{\sin\sphericalangle CT_B B'}{\sin\sphericalangle AT_B B'}=\left(\frac{\sin\beta}{\sin(B-B)}\right)^2,\ \frac{\sin\sphericalangle AT_C C'}{\sin\sphericalangle BT_C C'}=\left(\frac{\sin\gamma}{\sin(C-\gamma)}\right)^2$$

We'll apply then the Ceva's theorem.

**Observation**

It is possible to prove that the Exeter's point belongs to the Euler's line of the triangle $ABC$.

**13**.

We'll ration the same as for problem 9.

We find

$$\frac{U'B'}{U'C'}=\frac{AB'}{AC'}\cdot\frac{AB}{AC}\cdot\frac{UC}{UB}$$

$$\frac{V'C'}{V'A'}=\frac{BC'}{BA'}\cdot\frac{BC}{AB}\cdot\frac{VA}{VC}$$

$$\frac{W'A'}{W'B'}=\frac{CA'}{CB'}\cdot\frac{CA}{CB}\cdot\frac{WB}{WA}$$

By multiplying these relations side by side and taking into consideration the fact that $U-V-W$ is a transversal in the triangle $A'B'C'$, and using the Ceva's theorem, we find that $AA', BB', CC'$ are concurrent.

**14**.

We prove that $A', B', C'$ are collinear with the help of Menelaus' theorem (or with Bobillier's theorem).

We note $\{P\} = AA' \cap CC', \{Q\} = BB' \cap CC', \{R\} = AA' \cap BB'$.

In the complete quadrilateral $C'ACA'BB'$ the points $C, C'$ are harmonically conjugated in rapport to $P, Q$. The angular polar of the point $C'$ is the line $RC$ and this passes through the intersection of the lines $AQ, BP$, consequently, the triangles $ABC, QPR$ are homological. It is known that the triangles $C_aC_bC_c$, $ABC$ are homological. Because $C_aC_bC_c$ is inscribed in the triangle $ABC$, $PQR$ is circumscribed to triangle $ABC$, both being homological with $ABC$. Applying the theorem , it results that the triangles $C_aC_bC_c$, $PQR$ are homological.

**15.**

Let $m(\sphericalangle BAA''') = \alpha$, $m(\sphericalangle CBB''') = \beta$, $m(\sphericalangle ACC''') = \gamma$, see attached figure

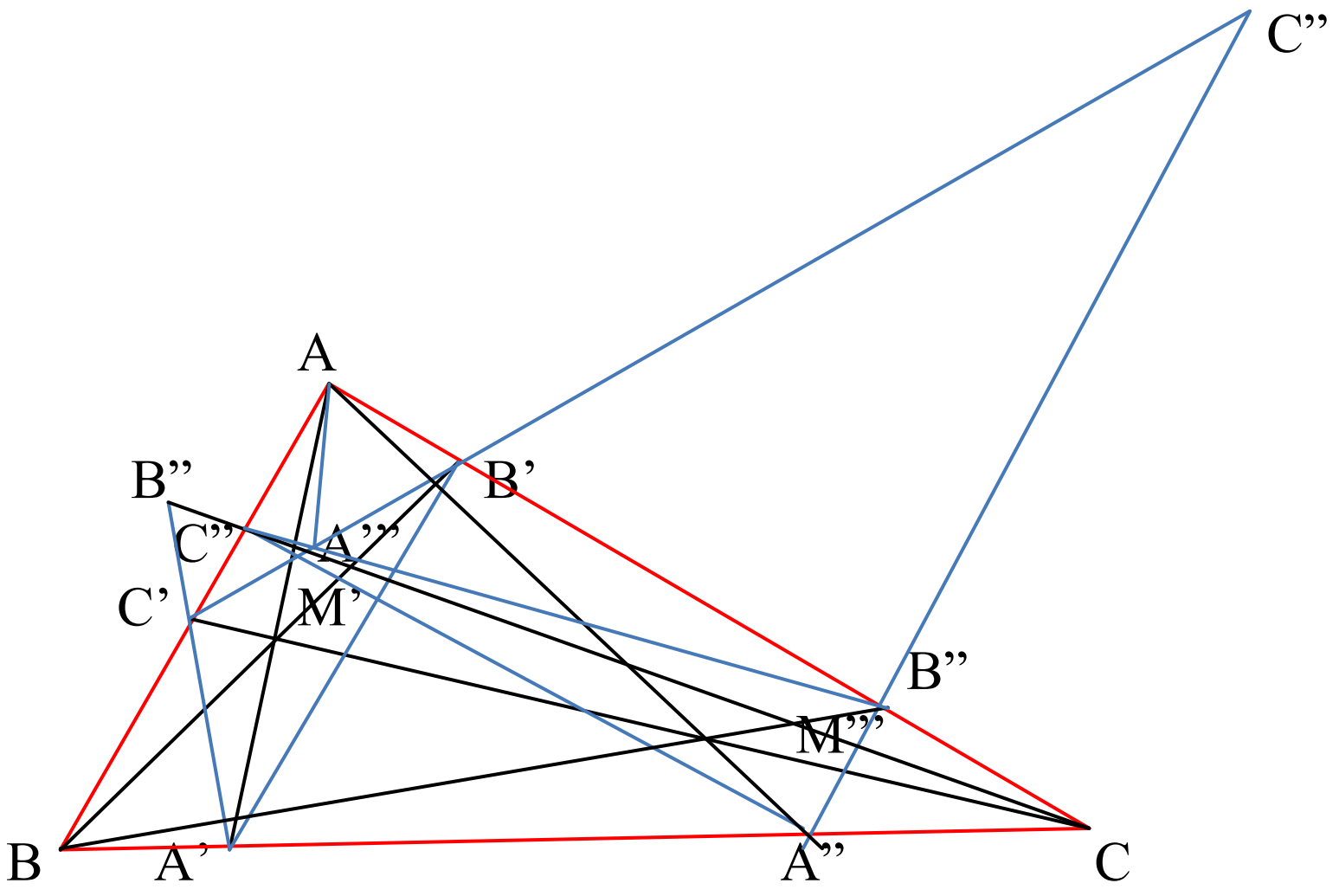

The sinus' theorem in the triangles $AC''A''', AB'A'''$ lead to

$$\frac{\sin\alpha}{A'''C''} = \frac{\sin C''}{AA'''} \qquad (1)$$

$$\frac{\sin(A-\alpha)}{A'''B'} = \frac{\sin B'}{AA'''} \qquad (2)$$

From the relations (1) and (2) it result

$$\frac{\sin\alpha}{\sin(A-\alpha)} = \frac{\sin C''}{\sin B'} \cdot \frac{B'''C''}{A'''B'} \qquad (3)$$

The sinus' theorem applied to the triangles $A'''C'C'', A'''B'B''$ gives:

$$\frac{A'''C''}{\sin C'} = \frac{C'C''}{\sin A'''} \qquad (4)$$

$$\frac{A'''B'}{\sin B''} = \frac{B'B''}{\sin A'''} \qquad (5)$$

From these relations we obtain

$$\frac{A'''C''}{A'''B'} = \frac{\sin C'}{\sin B''} \cdot \frac{C'C''}{B'B''} \qquad (6)$$

From (6) and (3) we retain

$$\frac{\sin\alpha}{\sin(A-\alpha)}=\frac{\sin C'}{\sin B'}\cdot\frac{\sin C''}{\sin B''}\cdot\frac{C'C''}{B'B''} \quad (7)$$

Similarly we obtain:

$$\frac{\sin\beta}{\sin(B-\beta)}=\frac{\sin A'}{\sin C'}\cdot\frac{\sin A''}{\sin C''}\cdot\frac{A'A''}{C'C''} \quad (8)$$

$$\frac{\sin\gamma}{\sin(C-\gamma)}=\frac{\sin B'}{\sin A'}\cdot\frac{\sin B''}{\sin A''}\cdot\frac{B'B''}{A'A''} \quad (9)$$

The relations (7), (8), (9) and the Ceva's reciprocal theorem lead us to

$$\frac{\sin\alpha}{\sin(A-\alpha)}=\frac{\sin\beta}{\sin(B-\beta)}=\frac{\sin\gamma}{\sin(C-\gamma)}=-1$$

Therefore to the concurrency of the lines $AA''', BB''', CC'''$ and implicitly to prove the homology of the triangles $ABC, A'''B'''C'''$.

To prove the homology of the triangles $A'B'C'$, $A'''B'''C'''$, we observe that

$$\frac{C'A'''}{B'A'''}=\frac{Aria_{\Delta AC'A''}}{Aria_{\Delta AB'A'''}}=\frac{AU}{AB'}\cdot\frac{\sin\alpha}{\sin(A-\alpha)}$$

Similar

$$\frac{A'B'''}{C'B'''}=\frac{BA'}{BC'}\cdot\frac{\sin\beta}{\sin(B-\beta)}$$

$$\frac{B'C'''}{A'C'''}=\frac{CB'}{CA'}\cdot\frac{\sin\gamma}{\sin(C-\gamma)}$$

We'll apply the Ceva's theorem. Similarly is proved the homology of the triangles $A''B''C''$, $A'''B'''C'''$.

**Observation**

This theorem could have been proved in the same manner as it was proved theorem 10.

**16**. $A_1B_1C_1$ the first Brocard triangle ($A_1$ is the projection of the symmedian center on the mediator on $(BC)$,etc.).

We've seen that $\sphericalangle A_1BC=\sphericalangle A_1CB=\sphericalangle B_1CA=\sphericalangle B_1CA=\sphericalangle C_1AB=\sphericalangle C_1BA=\omega$ (Brocard's angle).

If $M_a$ is the middle point of the side $(B_1C_1)$ and if we note $M_a(\alpha_1,\beta_1,\gamma_1)$ where $\alpha_1,\beta_1,\gamma_1$ are the barycentric coordinates of $M_a$, that is $\alpha_1 = Aria_{\Delta M_aBC}$, $\beta_1 = Aria_{\Delta M_aCA}$, $\gamma_1 = Aria_{\Delta M_aAB}$.

We find:

$$\alpha_1=\frac{a}{4\cos\omega}\left[c\sin(B-\omega)+b\sin(C-\omega)\right]$$

$$\beta_1=\frac{b}{4\cos\omega}\left[b\sin(B-\omega)+c\sin(A-\omega)\right]$$

$$\gamma_1 = \frac{c}{4\cos\omega}\left[c\sin\omega + b\sin(A-\omega)\right]$$

If $M_b M_c$ are the middle points of the sides $(A_1C_1)$ respectively $(A_1B_1)$ ; $M_b(\alpha_2,\beta_2,\gamma_2), M_c(\alpha_3,\beta_3,\gamma_3)$ we obtain the following relations:

$$\begin{cases} \alpha_2 = \dfrac{a}{4\cos\omega}\left[a\sin\omega + c\sin(B-\omega)\right] \\ \beta_2 = \dfrac{b}{4\cos\omega}\left[a\sin(C-\omega) + c\sin(A-\omega)\right] \\ \gamma_2 = \dfrac{c}{4\cos\omega}\left[c\sin\omega + a\sin(B-\omega)\right] \end{cases}$$

$$\begin{cases} \alpha_3 = \dfrac{a}{4\cos\omega}\left[a\sin\omega + b\sin(C-\omega)\right] \\ \beta_3 = \dfrac{b}{4\cos\omega}\left[b\sin\omega + b\sin(C-\omega)\right] \\ \gamma_3 = \dfrac{c}{4\cos\omega}\left[a\sin(B-\omega) + b\sin(A-\omega)\right] \end{cases}$$

We'll use then the result that $AM_a, BM_b, CM_c$ are concurrent if and only if $\alpha_2\beta_3\gamma_1 = \alpha_3\beta_1\gamma_2$.
Because in a triangle we have the following relations

$$\sin(A-\omega) = \sin\omega\cdot\frac{a^2}{bc}$$

$$\sin(B-\omega) = \sin\omega\cdot\frac{b^2}{ac}$$

$$\sin(C-\omega) = \sin\omega\cdot\frac{c^2}{ab}$$

The precedent relation will be verified.

**17**.

The De Longchamps's line is isotomic transversal to the Lemoine's line (the tri-linear polar of the symmedian $K$ of the triangle $ABC$). We have seen that the isotomic conjugate of the symmedian center is the homological center of the triangle $ABC$ and of its first Brocard triangle. Therefore, this point is the tri-linear point of the De Longchamp' line.

**18.**

Let's suppose that the isosceles similar triangles $BA'C, CB'A, AC'B$ are constructed in the exterior of triangle $ABC$ and that $m(\sphericalangle CBA') = m(\sphericalangle ACB') = m(\sphericalangle ABC') = x$.

We note $\{A_1\} = AA' \cap BC$, we have $\dfrac{BA_1}{CA_1} = \dfrac{Aria_{\Delta ABA'}}{Aria_{\Delta ACA'}}$. We obtain: $\dfrac{BA_1}{CA_1} = \dfrac{AB}{AC}\cdot\dfrac{\sin(B+x)}{\sin(C+x)}$.

Similarly we obtain: $\frac{CB_1}{AB_1}=\frac{BC}{AB}\cdot\frac{\sin(C+x)}{\sin(A+x)}$ and $\frac{AC_1}{BC_1}=\frac{AC}{BC}\cdot\frac{\sin(A+x)}{\sin(B+x)}$.

With the Ceva's theorem it results that $AA', BB', CC'$ are concurrent; therefore, the triangles $ABC, A'B'C'$ are homological and we note the homology center with $P$. It is observed that the triangles $A'B'C'$ and $ABC$ are homological because the perpendiculars from $A', B', C'$ on the sides $BC, CA, AB$ of $ABC$ are concurrent in $O$, which is the center of the circumscribed to the triangle $ABC$. This point is the otology center of the triangles $A'B'C'$ and $ABC$. In other words the triangles $ABC, A'B'C'$ are orthological. Their second orthological center is the intersection point $Q$ of the perpendiculars constructed from $A, B, C$ respectively on $B'C', C'A', A'B'$. In accordance with the Sondat's theorem, it results that $O, P, Q$ are collinear and their line is perpendicular on the homological axis of triangles $ABC, A'B'C'$.

**19**.

If we note $A_1$ the intersection point of the tri-linear polar of the orthocenter $H$ of the triangle $ABC$ with $BC$, and if we use the Menelaus' theorem, we find

$$\frac{A_1C}{A_1B}=\frac{tgB}{tgC} \qquad (1)$$

Let $\{A_1'\}=B''C''\cap BC$. Applying the sinus' theorem in the triangle $BCC'$, we find that $CC'=\frac{a\sin B}{\cos(B-A)}$. Similarly, $BB''=\frac{a\sin C}{\cos(C-A)}$, therefore $CC''=\frac{a\sin B}{2\cos(B-A)}$ and $BB''=\frac{a\sin C}{2\cos(C-A)}$.

We'll apply the Menelaus' theorem for the transversal $A_1'-B''-C''$ in the isosceles triangle $BOC$:

$$\frac{A_1'B}{A_1'C}\cdot\frac{OB''}{BB''}\cdot\frac{CC''}{OC''}=-1 \qquad (1)$$

It results

$$\frac{A_1'C}{A_1'B}=\frac{\dfrac{a\sin B}{2\cos(B-A)}}{R-\dfrac{a\sin B}{2\cos(B-A)}}\cdot\frac{R-\dfrac{a\sin C}{2\cos(C-A)}}{\dfrac{a\sin C}{2\cos(C-A)}}$$

$$\frac{A_1'C}{A_1'B}=\frac{\sin B}{\sin C}\cdot\frac{2R\cos(C-A)-a\sin C}{2R\cos(B-A)-a\sin B}$$

Substituting $2R=\frac{a}{\sin A}$ and after several computations we find

$$\frac{A_1C}{A_1B}=\frac{tgB}{tgC} \qquad (2)$$

From (1) and (2) we find that $A_1'=A_1$ and the problem is resolved.

**20**.

Solution given by Cezar Coşniță

Consider $ABC$ as a reference triangle, and let $(\alpha_1,\beta_1,\gamma_1)$, $(\alpha_2,\beta_2,\gamma_2)$ the barycentric coordinates of the points $M_1, M_2$.

The equations of the lines $A_1C_1$, $A_1B_1$ are

$$\frac{\gamma}{\gamma_1}+\frac{\alpha}{\alpha_1}=0,\quad \frac{\alpha}{\alpha_1}+\frac{\beta}{\beta_1}=0.$$

The line's $A_1M_1$ equation is $\frac{\gamma}{\gamma_1}+\frac{\alpha}{\alpha_1}+k\left(\frac{\alpha}{\alpha_1}+\frac{\beta}{\beta_1}\right)=0$, where $k$ is determined by the condition

$$\frac{\gamma_2}{\gamma_1}+\frac{\alpha_2}{\alpha_1}+k\left(\frac{\alpha_2}{\alpha_1}+\frac{\beta_2}{\beta_1}\right)=0.$$

When $k$ changes, we obtain the following equation

$$\left(\frac{\alpha_2}{\alpha_1}+\frac{\beta_2}{\beta_1}\right)\left(\frac{\gamma}{\gamma_1}+\frac{\alpha}{\alpha_1}\right)-\left(\frac{\gamma_2}{\gamma_1}+\frac{\alpha_2}{\alpha_1}\right)\left(\frac{\alpha}{\alpha_1}+\frac{\beta}{\beta_1}\right)=0$$

Considering $\alpha=0$ we have the equation

$$\left(\frac{\alpha_2}{\alpha_1}+\frac{\beta_2}{\beta_1}\right)\frac{\gamma}{\gamma_1}=\left(\frac{\gamma_2}{\gamma_1}+\frac{\alpha_2}{\alpha_1}\right)\frac{\beta}{\beta_1}$$

which give the coordinates for $A'$.

We observe that the coordinates $\beta,\gamma$ of the point $A'$.

We observe that the line $AA'$ passes through the point whose coordinates are

$$\left(\frac{\alpha_1}{\frac{\beta_2}{\beta_1}+\frac{\gamma_2}{\gamma_1}},\frac{\beta_1}{\frac{\gamma_2}{\gamma_1}+\frac{\alpha_2}{\alpha_1}},\frac{\gamma_1}{\frac{\alpha_2}{\alpha_1}+\frac{\beta_2}{\beta_1}}\right)$$

The symmetry of this expression shows that the similar lines $BB', CC'$ pass also through the same point $M$. If in these expressions we swap the indexes 1 and 2, we obtain the coordinates of the common point of three analogues lines. But the two groups of coordinates coincide with $M$. If $M_2$ is the weight center of the triangle $ABC$, then the coordinates of $M$ are $\frac{1}{\beta_1+\gamma_1}$, $\frac{1}{\gamma_1+\alpha_1}$, $\frac{1}{\alpha_1+\beta_1}$, therefore $M$ is the reciprocal of the complementary of $M_1$. For example, if $M_1$ is the reciprocal of the center of the circumscribed circle.

**21**.

The equations of the given lines are

$$-\left(b^2+c^2-a^2\right)x+\left(c^2+a^2-b^2\right)y+\left(a^2+b^2-c^2\right)z=0$$

$$-bcx+cay+abz=0$$

$$-(b+c-a)x+(c+a-b)y+(a+b-c)z=0.$$

To be concurrent it is necessary that

$$\Delta = \begin{vmatrix} b^2+c^2-a^2 & c^2+a^2-b^2 & a^2+b^2-c^2 \\ b+c-a & c+a-b & a+b-c \\ bc & ca & ab \end{vmatrix}$$

is null.
If we multiply the 3rd and 2nd lines and add the result to the first line, we obtain a determinant with two proportional lines, consequently $\Delta = 0$ and the lines are concurrent in a point $U$, which has the barycentric coordinates

$$U\left(-a(b-c)(b+c-a),\ b(c-a)(c+a-b\ ),\ c(a-b)(a+b-c)\right)$$

Similarly we find the coordinates of the points $V, W$.

The lines $AU, BV, CW$ are associated to the point $R$.
From the results obtained, we have that the tri-linear polar of the orthocenter, the Gergone's point and the center of the circumscribed circle are three concurrent points.

**22.**

i) Let note $\{N'\} = OI \cap A_1C_a$, we have that $\frac{N'I}{N'O} = \frac{\Omega}{R}$, therefore $N'I = OI\frac{R}{R-r}$.

This shows that $N'$ is a fixed point on the line $OI$; similarly it results that $B_1C_b, C_1C_c$ pass through $N'$.

ii) If we note $\{D'\} = AN' \cap BC$, we can show that $AD', AD_a$, where $D_a$ is the contact with $BC$ of the A-ex-inscribed circle are isogonal Cevians by using the Steiner's relation $\frac{BD'}{D'C} \cdot \frac{BD_a}{D_aC} = \frac{c^2}{b^2}$. To compute $BD'$ we apply the Menelaus theorem in the triangle $ADD'$ for the transversal $N' - C_a - A_1$.

**23.**

If the perpendicular bisector $AD$ intersects $BC$ in $A'$, it is observed that $AA'$ is tangent to the circumscribed circle to the triangle $ABC$, therefore the line of the points given in the problem is the Lemoine's line.

**24.**

Let $A'B'C'$ the orthic triangle of the triangle $ABC$. Because the quadrilateral $BCB'C'$ is inscribable, we have $A_1B \cdot A_1C = A_1C' \cdot A_1B'$, where $\{A_1\} = B'C' \cap BC$, therefore $A_1$ has equal power in rapport to the circumscribed circle and the Euler's circle (the circumscribed circle of the triangle $A'B'C'$), therefore $A_1$ is on the radical axis of these circles, similarly $B_1$, $C_1$ belong to this radical axis.

**25.**

We note $\{A''\} = AA' \cap BC$, the point $M$ is the middle of $(BC)$ and $A_1$ is the projection of $A$ on $BC$. From the similarity of the triangles $AA_1A''$, $A'MA''$ it results $\frac{A_1A''}{MA''} = \frac{h_a}{r}$, therefore $\frac{A_1M}{MA''} = \frac{h_a - r}{r}$. But $h_a = \frac{2s}{a}$, $r = \frac{s}{p}$ and $A_1M = \frac{a}{2} - c\cos B$, $MA'' = \frac{a}{2} - CA''$. After computations we find $CA'' = p - b$. It is known that $BC_a = p - b$, $C_a$ is the projection of $I$ on $BC$, therefore $AC_a$ Gergonne Cevian and it result that $AA''$ is its isotonic.

Similarly $BB'$, $CC'$ are the isotomics of the Gergonne's Cevians. The concurrency point is the Nagel's point of the triangle.

**26.**

The barycentric coordinates of the orthocenter $H$ are $H(tgA, tgB, tgC)$, of the symmedian center are $K(a^2, b^2, c^2)$. We note $K'$ the symmedian center of the orthic triangle $A'B'C'$ of the triangle $ABC$. Because $\sphericalangle B'A'C = \sphericalangle C'A'B = \sphericalangle A$ and the radius of the circumscribed circle of the triangle $A'B'C'$. The radius of the Euler's circle is $\frac{R}{2}$, we have that $B'C' = R\sin 2A$, therefore $K'(R^2\sin^2 2A, R^2\sin^2 2B, R^2\sin^2 2C)$. Because

$$\begin{vmatrix} tgA & tgB & tgC \\ a^2 & b^2 & c^2 \\ R^2\sin^2 2A & R^2\sin^2 2B & R^2\sin^2 2C \end{vmatrix}$$

is null $(a^2 = 4R^2\sin^2 A,\ \sin^2 2A = 4\sin^2 A\cos^2 A)$, it results that $H, K, K'$ are collinear.

**27.**

Let $ABC$ an isosceles triangle $AB = AC$, $BB'$ the symmedian from $B$ and $CC'$ the median.

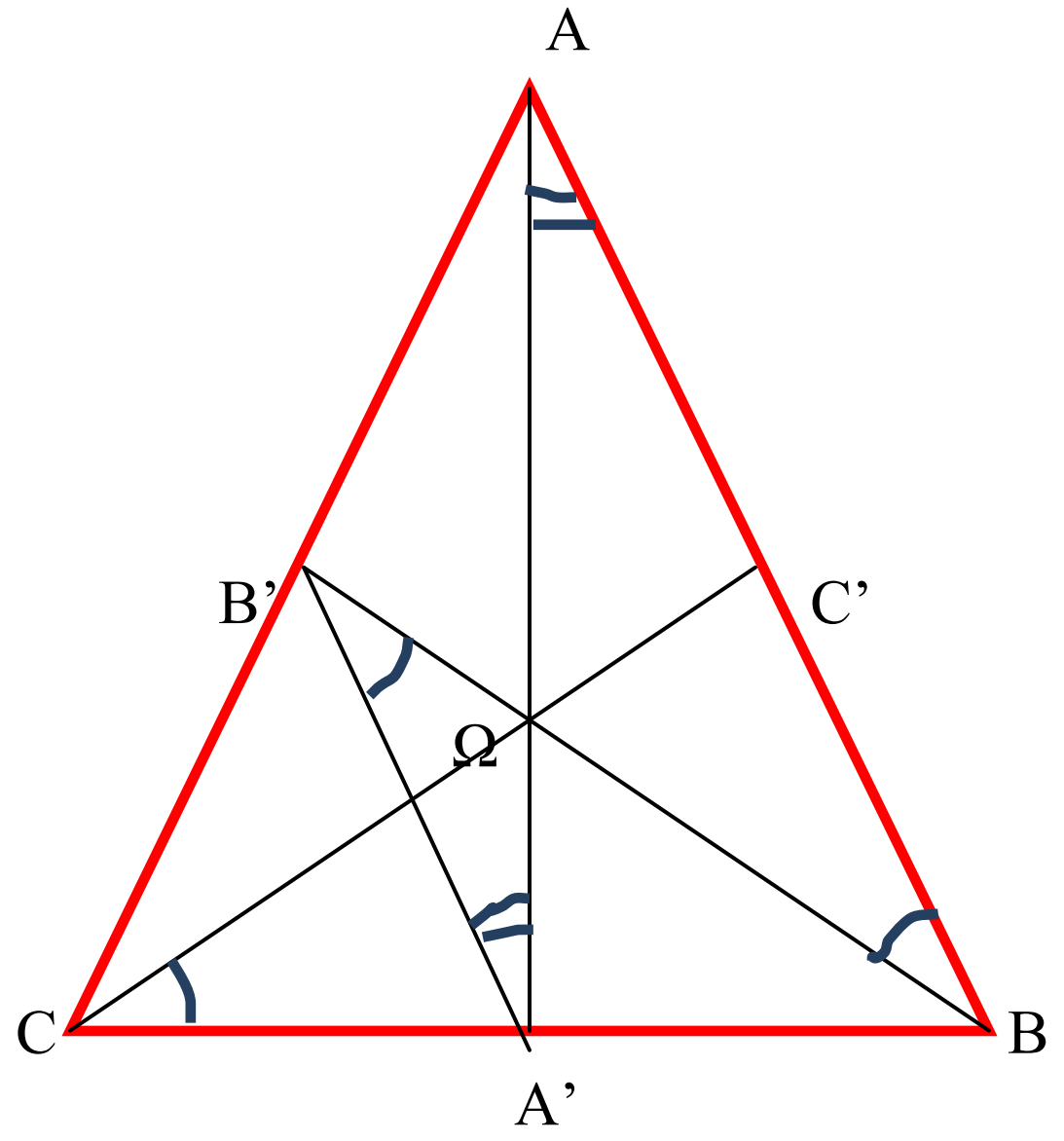

We'll note $\{\Omega\} = BB' \cap CC'$ and $\{A'\} = A\Omega \cap BC$ (see the above figure).

We have $\sphericalangle ABC \equiv \sphericalangle ACB$, it results $\sphericalangle \Omega CA \equiv \sphericalangle \Omega BC$ and $\sphericalangle \Omega BA \equiv \sphericalangle \Omega CB$. From the Ceva's theorem applied in the triangle $ABC$ it results $\frac{A'B}{A'C} \cdot \frac{B'C}{B'A} \cdot \frac{C'A}{C'B} = 1$, from which $\frac{A'B}{A'C} = \frac{B'A}{B'C}$, and with the reciprocal theorem of Thales we retain that $A'B' \parallel AB$. Then $\sphericalangle B'A'A \equiv \sphericalangle BAA'$ and $\sphericalangle ABB' \equiv \sphericalangle BBA'$ therefore $\sphericalangle \Omega B'A' \equiv \sphericalangle \Omega CA$ from which $\sphericalangle \Omega CA \equiv \sphericalangle \Omega AB \equiv \sphericalangle \Omega BC$, which means that $\Omega$ is a Brocard point of the triangle $ABC$.

Reciprocal

Let $\Omega$ a Brocard's point, therefore $\sphericalangle \Omega AB \equiv \sphericalangle \Omega BC \equiv \sphericalangle \Omega CA$, $BB'$ symmedian, $CC'$ the median, and $\{A'\} = A\Omega \cap BC$. From the Ceva's and Thales's theorems we retain that $A'B' \parallel AB$, therefore $\sphericalangle BAA' \equiv \sphericalangle AA'B'$ and $\sphericalangle ABB' \equiv \sphericalangle BB'A'$. Then $\sphericalangle \Omega A'B' \equiv \sphericalangle \Omega CB'$, therefore the quadrilateral $\Omega A'CB'$ is inscribable, from which $\sphericalangle \Omega CA' \equiv \sphericalangle \Omega B'A' \equiv \sphericalangle B'BA$. Therefore $m(\sphericalangle B) = m(\sphericalangle ABB') + m(\sphericalangle B'BC) = m(\sphericalangle C'CB) + m(\sphericalangle C'CA) = m(\sphericalangle C)$. We conclude that the triangle $ABC$ is isosceles.

**28**.

In the inscribable quadrilateral $B'A'BA$. We'll note $\{P\} = AB' \cap A'B$. According to Broard's theorem

$$OI \perp PC_1 \qquad (1)$$

In the inscribable quadrilateral $CAC'A'$ we'll note $\{Q\} = AC' \cap A'C$, the same Brocard's theorem leads to

$$OI \perp QB_1 \qquad (2)$$

In the quadrilateral $CBC'B'$ we note $\{R\} = BC' \cap CB'$, it results that

$$OI \perp RA \qquad (3)$$

From Pascal's theorem applied in the inscribed hexagon AB'CA'BC' we obtain that the points $P, Q, R$ are collinear, on the other side $A_1, B_1, C_1$ are collinear being on the homological axis of the homological triangles $ABC$ and $A'B'C'$.

From the relations (1), (2), (3) we find that $OI \perp A_1B_1$.

**29**.

In problem 21 we saw that $AU, BU, CW$ are concurrent. The polar of the point $A$ in rapport with the inscribed circle in the triangle $ABC$ is $C_bC_c$. The polar of the point $U$ in rapport with the same circle is the perpendicular constructed from $A$ on the line $UI$. The intersection of this perpendicular with $C_bC_c$ is the point $U'$ which is the pole of the line $AU$. Similarly we'll construct the poles $V', W'$ of the lines $BV, CW$. The poles of concurrent lines are collinear points, therefore $U', V', W'$ are collinear.

30.

Let $\{P\} = AC \cap MF$ and $\{P'\} = AC \cap HE$. We'll apply the Menelaus' theorem in the triangle $ABC$ for the transversal $M - F - P$. We'll obtain $\frac{PC}{PA} \cdot \frac{MA}{MB} \cdot \frac{FB}{FC} = 1$, therefore

$$\frac{PC}{PA} = \frac{FC}{FB} \quad (1)$$

The same theorem of Menelaus applied in the triangle $ADC$ for the transversal $E - H - P'$ leads to

$$\frac{P'C}{P'A} = \frac{ED}{EA} \quad (2)$$

Because

$EF \parallel AB$, from (1) and (2) and the fact that $\frac{FC}{FB} = \frac{ED}{EA}$ it results that $P \equiv P'$, therefore the lines $MF, EN$ intersect the line $AC$. But $\{C\} = DH \cap BF$, $\{A\} = BM \cap DE$, therefore the triangles $BMF, DNE$ are homological.

**31**.

It can be observed that the triangles $A_1B_1C$, $A_2B_2C_2$ are homological and their homological axis is $A - B - C_1$; from the reciprocal of the Desargues' theorem it results that $A_1A_2, B_1B_2, CC_2$ are concurrent.

**32.**

i) It can be observed that the triangles $ABC$, $A_1D'A'$ are polar reciprocal in rapport with the circle, then it will be applied the Charles' theorem

ii) Similar as in i).

**33**.

i) Let $\{L\} = OI \cap MD$. We have

$$\frac{LI}{LO} = \frac{r}{R} \quad (1)$$

If $\{L'\} = OI \cap NE$,

we have

$$\frac{L'I}{L'O} = \frac{r}{R} \quad (2)$$

We note $\{L''\} = OI \cap PF$, it results that.

$$\frac{L''I}{L''O} = \frac{r}{R} \quad (3)$$

The relations (1), (2), (3) lead to $L = L' = L''$

The point $L$ is the center of the homothety which transforms the circumscribed circle to the triangle $ABC$ in the inscribed circle in the triangle $ABC$ ($h_L\left(\frac{r}{R}\right)$).

iii) Through $h_L\left(\frac{r}{R}\right)$ the points $B,C,A$ have as images the points $B_1,C_1,A_1$, and the points $D,E,F$ are transformed in $A_2,B_2,C_2$.

Because the lines $AD,BE,CF$ are concurrent in the Gergonne's point $\Gamma$.

**34.**

Let $I_a$ the center of the A-ex-inscribed circle to the triangle $ABC$. Through the homothety $h_I\left(\frac{r}{R}\right)$ the inscribed circle is transformed in the A-ex-inscribed circle. The image of the point $A'$ through this homothety is the point $A_1$, which is situated on the A-ex-inscribed circle, such that $A'I \parallel C_aA_1$. Through the same homothety the transformed of the point $A''$ is the contact of the point $D_a$ with the line $BC$ of the A-ex-inscribed circle, therefore $AA''$ passes through $D_a$ ( $AA''$ is the Nagel Cevian), similarly $BB'',CC''$ are the Nagel Cevians. Consequently the point of concurrency of the lines $AA''$, $BB'',CC''$ is $N$ - the Nagel's point.

**35.**

Let $\overrightarrow{AB}=\vec{u}$, $\overrightarrow{AD}=\vec{v}$, $\overrightarrow{AE}=p\vec{u}$, $\overrightarrow{AF}=q\vec{v}$.

We'll have $\overrightarrow{AN}=\dfrac{\vec{u}+k\vec{v}}{1+k}$, $\overrightarrow{AP}=\dfrac{p\vec{u}+kq\vec{v}}{1+k}$.

We'll write the vector $\overrightarrow{AC}$ in two modes:

$$\overrightarrow{AC}=\frac{p\vec{u}+x\vec{v}}{1+x}=\frac{\vec{u}+yq\vec{v}}{1+y}.$$

By making the coefficients of $\vec{u}$ and $\vec{v}$ we

We'll find $y=\dfrac{x}{p\cdot q}$, $y=\dfrac{q(p-1)}{q-1}$.

$$\overrightarrow{AC}=\frac{p(q-1)\vec{u}+q(p-1)\vec{v}}{pq-1};\ \overrightarrow{AM}=\frac{k}{k+1},\ \overrightarrow{AC}=\frac{k}{k+1}\cdot\frac{p(q-1)\vec{u}+q(p-1)\vec{v}}{pq-1}.$$

$$\overrightarrow{NP}=\overrightarrow{AP}-\overrightarrow{AN}=\frac{(p-1)\vec{u}+k(q-1)\vec{v}}{1+k}$$

$$\overrightarrow{MP}=\overrightarrow{AP}-\overrightarrow{AM}$$

$$\overrightarrow{MP}=\frac{p(pq-1-kq+k)}{(k+1)(pq-1)}\vec{u}+\frac{kpq(q-1)}{(k+1)(pq-1)}\vec{v}$$

The vectors $\overrightarrow{NP}$, $\overrightarrow{MP}$ are collinear if and only if

$$\frac{p-1}{1+k}\cdot\frac{(k+1)(pq-1)}{p(pq-1-kq+k)}=\frac{kq}{1+k}\cdot\frac{(k+1)(pq-1)}{kpq(q-1)}\Leftrightarrow k=1.$$

The line obtained in the case $k=1$ is the line Newton-Gauss of the complete quadrilateral $ABCDEF$.

**36**.

We have $m(\sphericalangle ABD)=\frac{3}{4}\widehat{B}$, $m(\sphericalangle ACD)=\frac{3}{4}\widehat{C}$.

The sinus' theorem implies $\frac{BD}{\sin BAD}=\frac{AD}{\sin\frac{3}{4}B}$; $\frac{CD}{\sin CAD}=\frac{AD}{\sin\frac{3}{4}C}$.

On the other side $\frac{DC}{\sin\frac{1}{4}B}=\frac{BD}{\sin\frac{1}{4}C}$.

We find

$$\frac{\sin\sphericalangle BAD}{\sin\sphericalangle CAD}=\frac{\sin\frac{1}{4}C\sin\frac{3}{4}B}{\sin\frac{1}{4}B\sin\frac{3}{4}C}$$

We continue the rational in the same manner and we use then the trigonometric variant of Ceva's theorem.

**37**.

$m(\sphericalangle ABI_1)=45°+\frac{3}{4}\widehat{B}$, $m(\sphericalangle ACI_1)=45°+\frac{3}{4}\widehat{C}$

We have

$$\frac{BI_1}{\sin BAI_1}=\frac{AI_1}{\sin\left(45°+\frac{3}{4}\widehat{B}\right)};\ \frac{CI_1}{\sin CAI_1}=\frac{AI_1}{\sin\left(45°+\frac{3}{4}\widehat{C}\right)}$$

From the triangle $BI_1C$ we retain that

$$\frac{BI_1}{CI_1}=\frac{\sin\frac{1}{4}(A+B)}{\sin\frac{1}{4}(A+C)};\ \frac{BI_1}{CI_1}=\frac{\sin\left(45°-\frac{1}{4}\widehat{C}\right)}{\sin\left(45°-\frac{1}{4}\widehat{B}\right)}.$$

We obtain

$$\frac{\sin\sphericalangle BAI_1}{\sin\sphericalangle CAI_1}=\frac{\sin\left(45°-\frac{1}{4}\widehat{C}\right)\cdot\sin\left(45°+\frac{3}{4}\widehat{B}\right)}{\sin\left(45°-\frac{1}{4}\widehat{B}\right)\sin\left(45°+\frac{3}{4}\widehat{C}\right)}$$

Similarly we find

$$\frac{\sin\sphericalangle CBI_2}{\sin\sphericalangle ABI_2}=\frac{\sin\left(45°-\frac{1}{4}\widehat{A}\right)\cdot\sin\left(45°+\frac{3}{4}\widehat{C}\right)}{\sin\left(45°-\frac{1}{4}\widehat{C}\right)\sin\left(45°+\frac{3}{4}\widehat{A}\right)}$$

$$\frac{\sin\sphericalangle ACI_3}{\sin\sphericalangle BCI_3}=\frac{\sin\left(45^\circ-\frac{1}{4}\widehat{B}\right)\cdot\sin\left(45^\circ+\frac{3}{4}\widehat{A}\right)}{\sin\left(45^\circ-\frac{1}{4}\widehat{A}\right)\sin\left(45^\circ+\frac{3}{4}\widehat{B}\right)}$$

**38**.

$$m\left(\sphericalangle ABI_1\right)=90^\circ+\frac{3}{4}B;\quad m\left(\sphericalangle ACI_1\right)=90^\circ+\frac{3}{4}C$$

$$\frac{\sin\sphericalangle BAI_1}{BI_1}=\frac{\sin\sphericalangle ABI_1}{AI_1},\ \frac{\sin\sphericalangle CAI_1}{CI_1}=\frac{\sin\sphericalangle ACI_1}{AI_1}$$

We obtain

$$\frac{\sin\sphericalangle BAI_1}{\sin\sphericalangle ACI_1}=\frac{BI_1}{CI_1}\cdot\frac{\cos\frac{3}{4}B}{\cos\frac{3}{4}C}$$

But the sinus' theorem in the triangle $BI_1C$ gives us:

$$\frac{BI_1}{CI_1}=\frac{\sin\sphericalangle BCI_1}{\sin\sphericalangle CBI_1}.$$

Then we obtain

$$\frac{BI_1}{CI_1}=\frac{\cos\frac{1}{4}C}{\cos\frac{1}{4}B},\ \frac{\sin\sphericalangle BI_1C}{\sin\sphericalangle CAI_1}=\frac{\cos\frac{1}{4}C\cos\frac{3}{4}B}{\cos\frac{1}{4}B\cos\frac{3}{4}C}$$

Similarly we find

$$\frac{\sin\sphericalangle CBI_2}{\sin\sphericalangle ABI_2}=\frac{\cos\frac{1}{4}A\cos\frac{3}{4}C}{\cos\frac{1}{4}C\cos\frac{3}{4}A}$$

and

$$\frac{\sin\sphericalangle ACI_3}{\sin\sphericalangle BCI_3}=\frac{\cos\frac{1}{4}B\cos\frac{3}{4}A}{\cos\frac{1}{4}A\cos\frac{3}{4}B}$$

**39.**

$$\sphericalangle CBI_a=\sphericalangle ABI_c=90^\circ-\frac{1}{2}\sphericalangle B,\ \sphericalangle I_cBC=90^\circ-\frac{1}{2}\sphericalangle B,$$

$$\sphericalangle I_1BC=45^\circ-\frac{1}{4}\left(\sphericalangle B\right),\ \sphericalangle ABI_1=\frac{3}{4}\left(\sphericalangle B\right)-45^\circ,\ \sphericalangle BCI_a=\sphericalangle ACI_b=90^\circ-\frac{1}{2}\sphericalangle CI$$

$$\sphericalangle I_aCB=90^\circ+\frac{1}{2}\left(\sphericalangle C\right),\ \sphericalangle I_1CB=45^\circ+\frac{1}{4}\left(\sphericalangle C\right),\ \sphericalangle ACI_1=45^\circ+\frac{3}{4}\left(\sphericalangle C\right).$$

We'll obtain

$$\frac{\sin \sphericalangle I_1AB}{\sin \sphericalangle I_1AC}=\frac{\sin\left(45^\circ+\frac{1}{4}C\right)}{\sin\left(45^\circ+\frac{1}{4}B\right)}\cdot\frac{\sin\left(\frac{3}{4}B-45^\circ\right)}{\sin\left(45^\circ-\frac{3}{4}C\right)};$$

$$\frac{\sin \sphericalangle I_2BC}{\sin \sphericalangle I_2BA}=\frac{\sin\left(45^\circ+\frac{1}{4}A\right)}{\sin\left(45^\circ+\frac{1}{4}C\right)}\cdot\frac{\sin\left(45^\circ-\frac{3}{4}C\right)}{\sin\left(45^\circ-\frac{3}{4}A\right)};$$

$$\frac{\sin \sphericalangle I_3CA}{\sin \sphericalangle I_3CB}=\frac{\sin\left(45^\circ+\frac{1}{4}B\right)}{\sin\left(45^\circ+\frac{1}{4}A\right)}\cdot\frac{\sin\left(45^\circ-\frac{3}{4}A\right)}{\sin\left(\frac{3}{4}B-45^\circ\right)}.$$

**40.**

Let $T_1T_2T_3$ the tangential triangle of the triangle $A_1B_1C_1$.
We have

$$\sphericalangle T_1C_1A \equiv \sphericalangle AA_1C_1 \text{ and } \sphericalangle T_1B_1A \equiv \sphericalangle AA_1B_1;$$

$$\frac{AC_1}{\sin C_1T_1A}=\frac{AT_1}{\sin T_1C_1A};\ \frac{AB_1}{\sin B_1T_1A}=\frac{AT_1}{\sin T_1B_1A}$$

Therefore

$$\frac{\sin C_1T_1A}{\sin B_1T_1A}=\frac{AC_1}{AB_1}\cdot\frac{\sin T_1C_1A}{\sin T_1B_1A}$$

But

$$\frac{\sin T_1C_1A}{\sin T_1B_1A}=\frac{\sin AA_1C_1}{\sin AA_1B_1}$$

On the other side $\sin \sphericalangle AA_1C_1 = 2R - AC_1$ and $\sin\left(\sphericalangle AA_1B_1\right) = 2R - AB_1$
We'll obtain that

$$\frac{\sin C_1T_1A}{\sin B_1T_1A}=\left(\frac{AC_1}{AB_1}\right)^2, \text{ etc.}$$

**41.**

$I$ is the orthocenter of the orthic triangle of the triangle $I_aI_bI_c$ (that is of the triangle $ABC$). The point $M$ is the center of the circumscribed circle o the triangle $IBC$ (this circle passes through the point $I_a$). The perpendicular constructed from $I_a$ on $O_1M$ is the radical axis of the circumscribed circle of the triangle $I_aI_bI_c$ and $IBC$.

On the other side $BC$ is the radical axis of the circumscribed circles $IBC$ and $ABC$, it results that that the intersection between $BC$ and $A_1C_1$ is the radical center of the circumscribed circles of the triangles $I_aI_bI_c$, $ABC$ and $IBC$ - this point belongs to the tri-linear polar of the point $I$ in rapport to the triangle $ABC$.

**42.**

We observe that the triangle $ABC$ is congruent with the triangle $C_1A_1B_1$. Because $AB = C_1A_1$ and $AC_1 = BA_1$, it results that the quadrilateral $AC_1BA_1$ is an isosceles trapeze.

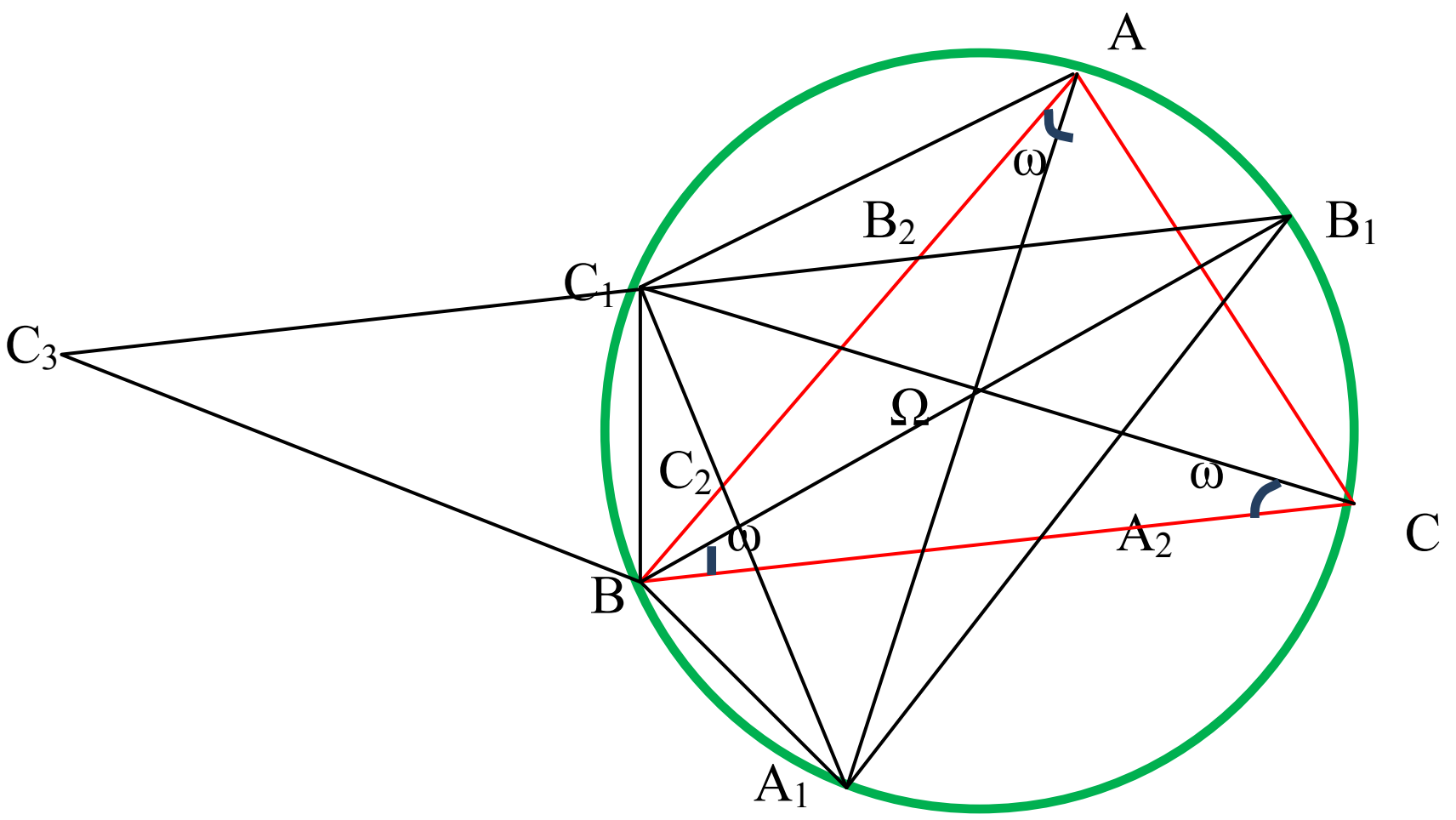

$C_2C_3$ is perpendicular on $BC_1$ and because it passes through its middle it will contain also the center $O$ of the circumscribed circle of the triangle $ABC$. Similarly we show that $B_2C_3$ and $A_2A_3$ pass through the center $O$. The triangles $A_2B_2C_2, A_3B_3C_3$ are homological and their homology center is $O$.

We'll note

$$\{L\} = A_2B_2 \cap A_3B_3,\ \{M\} = B_2C_2 \cap B_3C_3,\ .$$

The line $L - M - N$ is the homology axis of the triangles $A_2B_2C_2, A_3B_3C_3$. From the congruency of the angles $C_2AB_2, B_1C_1C_2$ it results that the quadrilateral $B_2AC_1C_2$ is inscribable, therefore

$$MC_1 \cdot MA = MC_2 \cdot MB_2.$$

This equality shows that the point $M$ has equal powers in rapport to the circle $(O)$ and in rapport to circle $(O_1)$, therefore $M$ belongs to the radical axis of these circles.

Similarly we can show that $L, N$ belong to this radical axes also. This gives us also a new proof of the triangles' homology from the given problem.

**43.**

Let $A_1B_1C_1$ the contact triangle of the triangle $ABC$ (the pedal triangle of $I$). We note $A', B', C'$ the middle points of the segments $(AI),(BI),(CI)$.

We note $A_2B_2C_2$ the anti-pedal triangle of the contact point $I$ with the triangle $ABC$. Because $B_2C_2$ is perpendicular on $AI$ and $A'O \parallel AI$.

It results that $B_2C_2$ is perpendicular on $A'O$ as well, therefore $B_2C_2$ is the radical axis of the circles $(AB_1C_1)$, $(ABC)$. The line $B_1C_1$ is the radical axis of the circles $(AB_1C_1)$ and $(A_1B_1C_1)$, it results that the point $\{A_3\} = B_1C_1 \cap B_2C_2$ is the radical center of the above three circles.
Therefore $A_3$ is on the radical axis of the circumscribed and inscribed circles to the triangle $ABC$.

Similarly, it results that the intersection points $B_3, C_3$ of the pairs of lines $(C_1A_1, C_2A_2)$, $(A_1B_1, A_2B_2)$ are on the radical axis of the inscribed and circumscribed circles.

The point $A_2$ is the radical center of the circles $(BA_1C_1),(CA_1B_1),(ABC)$, it, therefore, belongs to the radical axis of the circles $(BA_1C_1),(CA_1B_1)$; this is the line $A_2A_1$ which passes through .

**44**.

The lines $A_1A_2$, $C_1C_2$ are common cords of the given circles and are concurrent in the radical center of these circles; it results that the triangles $A_1B_1C_1$, $A_2B_2C_2$ are homological.

We'll note

$$\{A_3\} = B_1C_1 \cap B_2C_2 .$$

The quadrilateral $B_1C_1B_2C_2$ is inscribed in the circle $(C_a)$

It results that

$$A_3B_1 \cdot A_3C_1 = A_3B_2 \cdot A_3C_2 .$$

This relation shows that $A_3$, which belongs to the homology axis of triangles $A_1B_1C_1$, $A_2B_2C_2$, belongs to the radical axis of the circumscribed circles to these triangles.

**45**.

We will use the same method as in the problem 43.

**46**.

Because $AB = A'B'$ and $AA' = BB'$, it results that $AB'BA'$ is an isosceles trapeze, $C_1C_2$ id the perpendicular bisector of the cord $AB'$, therefore it passes through $O$. Similarly, it results that $A_1A_2$ and $B_1B_2$ pass through $O$; this point is the homology center of the triangles.

If $L, M, N$ are the intersections points of the opposite sides of the triangles $A_1B_1C_1$, $A_2B_2C_2$ ($\{M\} = B_1C_1 \cap B_2C_2$), we have that $L, M, N$ are collinear, and these belong to the homology axis of the above triangles.

We observe that the quadrilateral $AB'B_1C_1$ is inscribable, therefore $MA \cdot MA' = MB_1 \cdot MC_1$
This equality shows that the point $M$ has equal powers in rapport to the circles $(O)$ and $(O_1)$. We noted with $(O_1)$ the center of the circumscribed circle of the triangle $A_1B_1C_1$, therefore $M$ belongs to the radical axis of these circles.

Similarly, $N$ belongs to this radical axis, therefore $L-M-N$ is perpendicular on $OO_1$.

**47.**

i) We'll note with $\alpha$ the measure of the angle formed by the line constructed from $I$ with the bisector $AI_a$ ($I_a$ being the center of the A-ex-inscribed circle).

$$CA_1^{'} = CA_1 \cdot \cos \widehat{A_1CA_1^{'}} = CA_1 \cdot \sin\left(\alpha - \frac{C}{2}\right)$$

$$CA_1 = II_a \cdot \sin \widehat{A_1BC} = II_a \cos\left(\frac{B}{2} + \alpha\right)$$

$$CA_1^{'} = II_a \cos\left(\frac{B}{2} + \alpha\right) \cdot \sin\left(\alpha - \frac{C}{2}\right)$$

$$BA_1^{'} = II_a \cos\left(\alpha - \frac{C}{2}\right) \cdot \sin\left(\frac{B}{2} + \alpha\right)$$

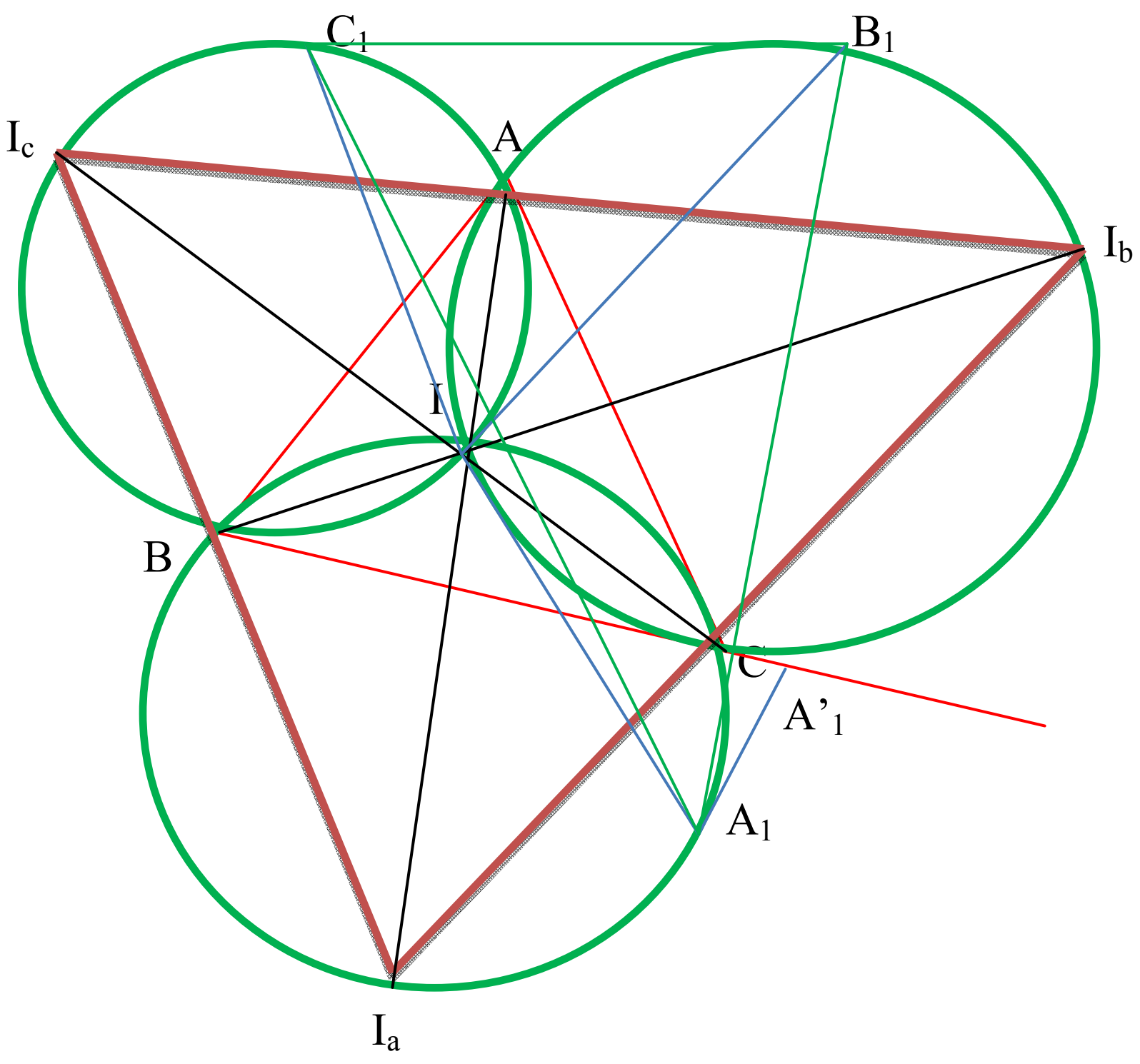

Therefore

$$\frac{A_1^{'}B}{A_1^{'}C} = \frac{tg\left(\frac{B}{2} + \alpha\right)}{tg\left(\alpha - \frac{C}{2}\right)}$$

Similarly

$$\frac{B_1^{'}C}{B_1^{'}A} = \frac{ctg\alpha}{tg\left(\frac{B}{2}+\alpha\right)}$$

$$\frac{C_1^{'}A}{C_1^{'}B} = tg\alpha \cdot tg\left(\alpha - \frac{C}{2}\right)$$

Then is applied the Menelaus' theorem.

ii) We'll note $A_1^{''}$ the intersection of the tangent in $A_1$ with $BC$.
We find that

$$\frac{A_1^{''}C}{A_1^{''}B} = \left(\frac{A_1C}{A_1B}\right)^2 = \frac{\cos^2\left(\frac{B}{2}+\alpha\right)}{\cos^2\left(\frac{C}{2}-\alpha\right)}$$

If $B_1^{''}$ and $C_1^{''}$ are the intersections of the tangents constructed in $B_1$ respectively in $C_1$ to the circles $CIA, BIA$ with $AC$ respectively $AB$, we'll have

$$\frac{B_1^{''}A}{B_1^{''}C} = \left(\frac{B_1A}{B_1C}\right)^2 = \frac{\sin^2\alpha}{\cos^2\left(\frac{B}{2}+\alpha\right)},$$

$$\frac{C_1^{''}B}{C_1^{''}A} = \frac{\cos^2\left(\frac{C}{2}-\alpha\right)}{\sin^2\alpha}$$

Then we apply the Menelaus' theorem

**48**.

The triangles
$APQ, BRS$ are homological because

$AP \cap BR = \{H\}$, $QR \cap SB = \{G\}$, $PQ \cap RS = \{O\}$

and $H, G, O$ - collinear (the Euclid line). According to the reciprocal of the Desargues' theorem, the lines $AB, PR, QS$ are concurrent.

**49**.

It is known that $AP$ (see the figure below) is symmedian in the triangle $ABC$, therefore if we note with $S$ the intersection point of the segment $(AP)$ with the circle, we have that

$$arc\,\widehat{BS} \equiv arc\,\widehat{CQ}$$

Therefore $S,Q$ are symmetric in rapport to $PM$.

Because $PM$ is the mediator of $BC$, it results that $PM$ passes through $O$, and $AR \parallel BC$ leads to the conclusion that $PM$ is also the mediator of the cord $\widehat{AR}$, therefore $R$ is the symmetric of $A$ in rapport with $PM$.

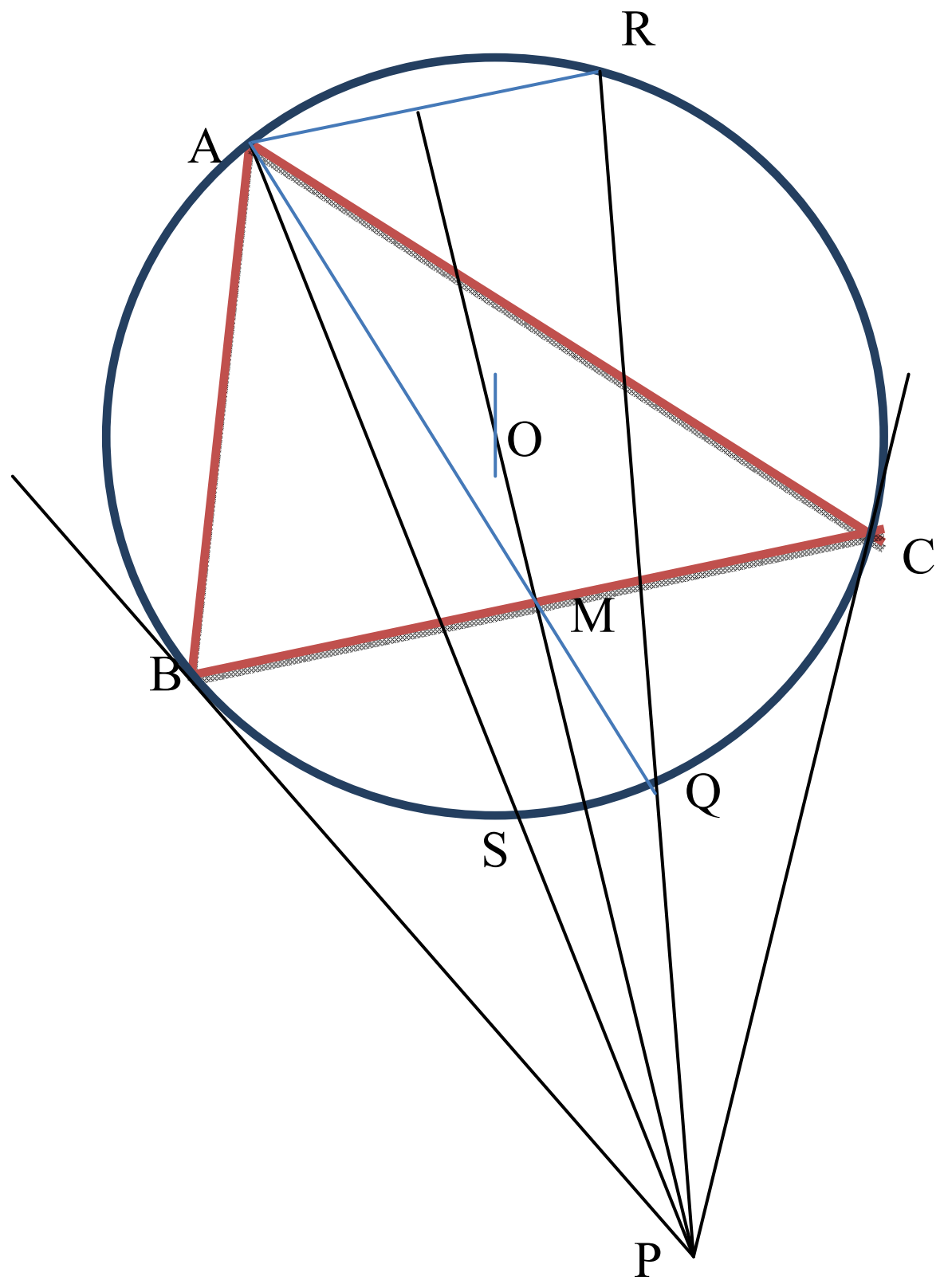

Because $A,S,P$ are collinear and $Q,R$ are the symmetric of $S,A$ in rapport to $PM$, it results that $P,Q,R$ are collinear

**50.**

We will note $\alpha = m\sphericalangle A_1II_a$, $I_a$ is the center of the A-ex-inscribed circle, $\beta = m\sphericalangle A_2II_a$.

We have

$$\Delta A_3BA_2 \sim \Delta A_3A_1C$$

From this similarity we find

$$\frac{A_3B}{A_3C} = \left(\frac{BA_2}{A_1C}\right)^2 \cdot \frac{A_3A_1}{A_3A_2}$$

The sinus' theorem implies

$$\frac{A_3A_1}{\sin BCA_1} = \frac{CA_1}{\sin A_1A_3C}$$

$$\frac{A_3A_2}{\sin CBA_2} = \frac{BA_2}{\sin BA_3A_2}$$

From here

$$\frac{A_3A_1}{A_3A_2} = \frac{CA_1}{BA_2} \cdot \frac{\sin BCA_1}{\sin CBA_2}$$

$$\frac{A_3B}{A_3C} = \frac{BA_2}{CA_1} \cdot \frac{\sin BCA_1}{\sin CBA_2} = \frac{\sin CBA_2}{\sin CBA_1} \cdot \frac{\sin BCA_1}{\sin CBA_2}$$

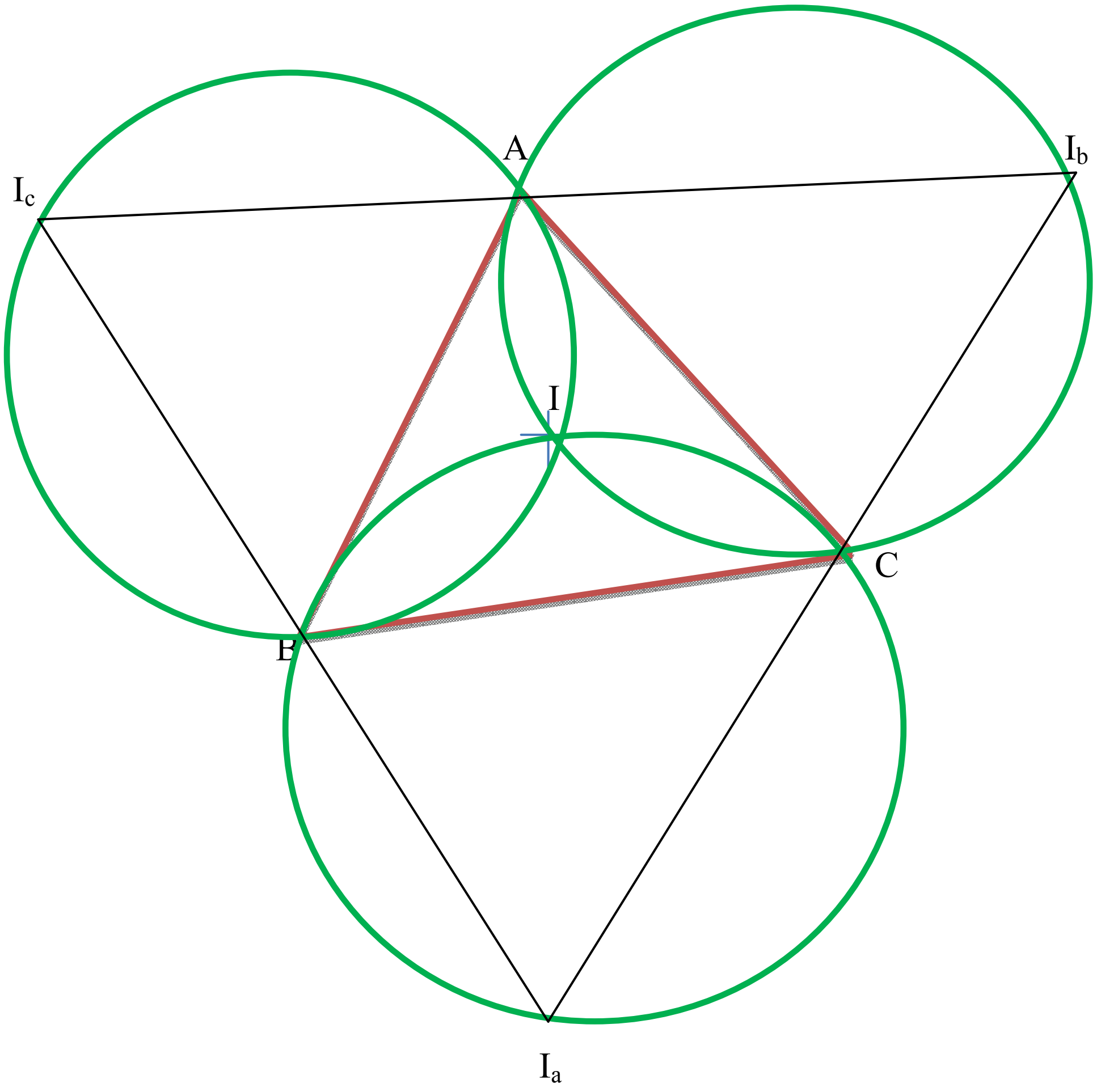

$$BCA_2 = \frac{A}{2} + \frac{B}{2} + \beta = 90° - \frac{C}{2} + \beta = 90° - \left(\frac{C}{2} - \beta\right)$$

$$\sin BCA_2 = \cos\left(\frac{C}{2} - \beta\right)$$

$$\sin BCA_1 = \cos\left(\frac{C}{2} - \alpha\right)$$

$$\sin CBA_1 = \cos\left(\frac{B}{2}+\alpha\right)$$

$$\sin CBA_2 = \cos\left(\frac{B}{2}+\beta\right)$$

Therefore

$$\frac{A_3B}{A_3C} = \frac{\cos\left(\frac{C}{2}-\alpha\right)\cdot\cos\left(\frac{C}{2}-\beta\right)}{\cos\left(\frac{B}{2}+\alpha\right)\cdot\cos\left(\frac{B}{2}+\beta\right)} \quad (1)$$

Similarly

$$\frac{B_3C}{B_3A} = \frac{\cos\left(\frac{B}{2}+\alpha\right)\cdot\cos\left(\frac{B}{2}+\beta\right)}{\sin\alpha\sin\beta} \quad (2)$$

$$\frac{C_3A}{C_3B} = \frac{\sin\alpha\sin\beta}{\cos\left(\frac{C}{2}-\alpha\right)\cdot\cos\left(\frac{C}{2}-\beta\right)} \quad (3)$$

The relations (1), (2), (3) and the Menelaus' theorem lead to the solution of the problem.

**51.**

The mediators of the segments $P_1Q_1, P_2Q_2, P_3Q_3$ pass through the point $O_9$, which is the middle of the segment $PQ$ and which is the center of the circle which contains the points $P_1, P_2, P_3, Q_1, Q_2, Q_3, R_1, R_2, R_3$

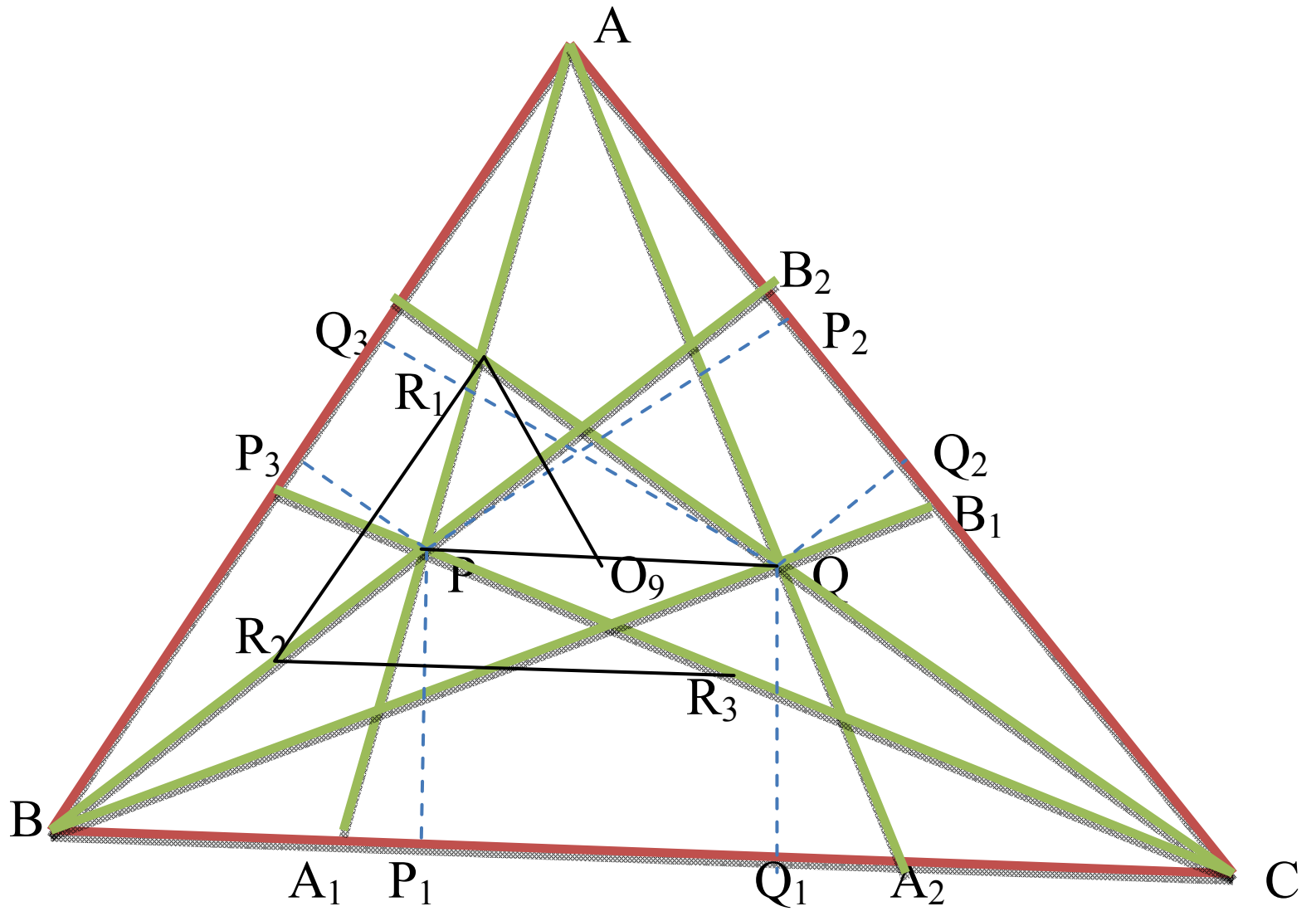

Because $R_1R_2$ is middle line in the triangle $PAB$, it results that

$$\sphericalangle PAB \equiv \sphericalangle PR_1R_2 \quad (1)$$

Also, $R_1R_3$ is middle line in the triangle $PAC$, and $R_1O_9$ is middle line in the triangle $PAQ$, therefore we obtain

$$\sphericalangle O_9R_1R_3 \equiv \sphericalangle QAC \qquad (2)$$

The relations (1), (2) and the fact that $AP$ and $AQ$ are isogonal Cevians lead to

$$\sphericalangle PR_1R_2 \equiv \sphericalangle O_9R_1R_3 \qquad (3) \quad .$$

The point $O_9$ is the center of the circumscribed circle of the triangle $R_1R_2R_3$. Considering relation (3) and a property of the isogonal Cevians, we obtain that in the triangle $R_1R_2R_3$ the line $R_1P$ is a height and from $R_2R_3 \parallel BC$ we obtain that $AP$ is height in the triangle $ABC$.

Applying the same rational one more time, it will result that $BP$ is height in the triangle $ABC$, consequently $P$ will be the orthocenter of the triangle $ABC$, $Q$ will be the center of the circumscribed circle to the triangle $ABC$, and $O_9$ will be the center of the circle of nine points of the triangle $ABC$.

**52**.

Let $H$ the orthocenter of the triangle and $AA_1^{'}, BB_1^{'}, CC_1^{'}$ three concurrent Cevians . We note $m(\sphericalangle BAA_1) = \alpha$, $m(\sphericalangle CBB_1) = \beta$, $m(\sphericalangle ACC_1) = \gamma$. $A_1$ is the intersection of $BC$ with the perpendicular from $H$ on $AA_1^{'}$.

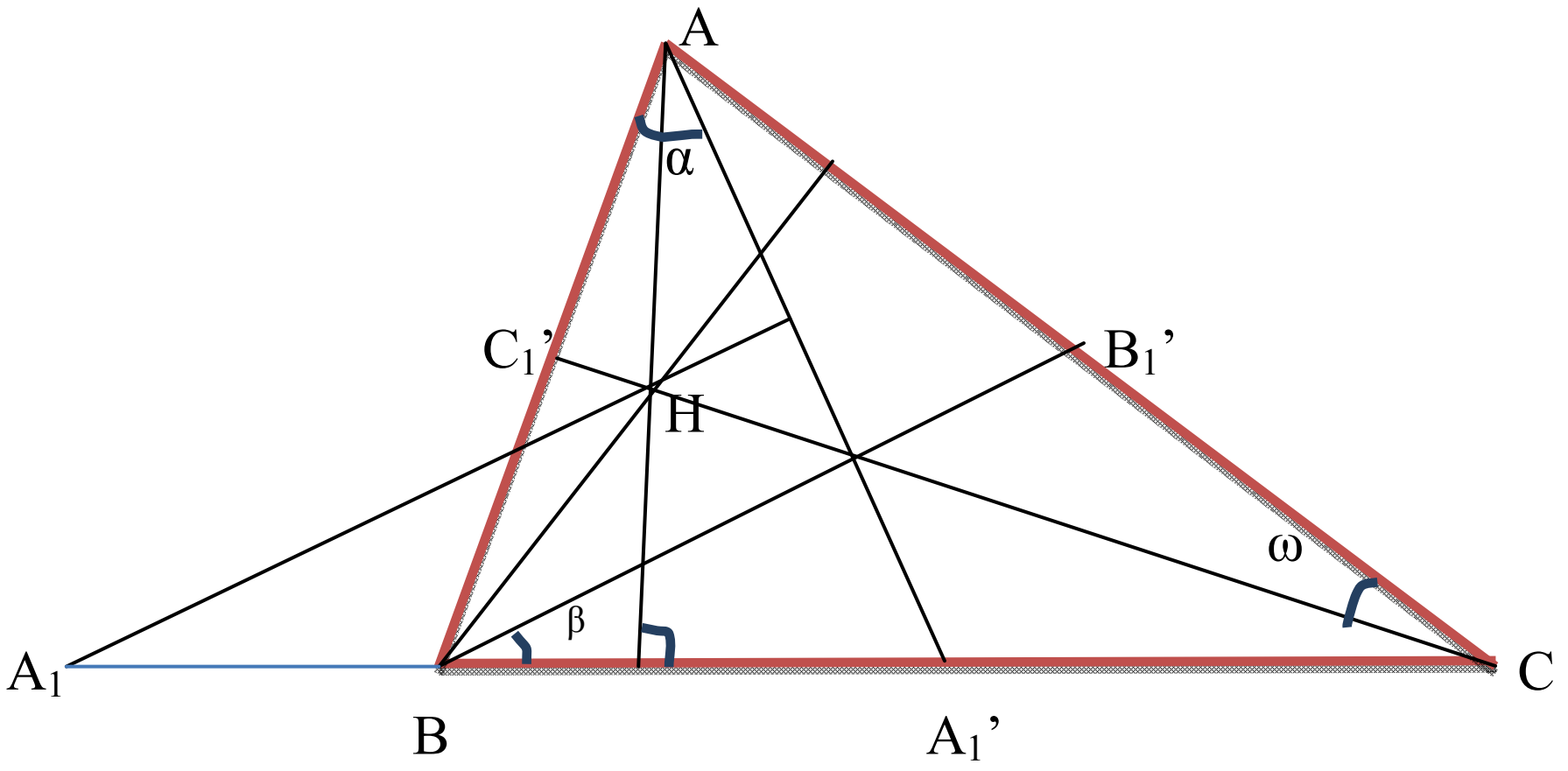

Because $AA_1^{'}, BB_1^{'}, CC_1^{'}$ are concurrent we have

$$\frac{\sin\alpha}{\sin(A-\alpha)} \cdot \frac{\sin\beta}{\sin(B-\beta)} \cdot \frac{\sin\gamma}{\sin(C-\gamma)} = 1$$

On the other side

$$\frac{A_1B}{A_1C} = \frac{AriaA_1HB}{AriaA_1HC} = \frac{A_1H \cdot \sin A_1HB \cdot HB}{A_1H \cdot \sin A_1HC \cdot HC}$$

Because $\widehat{A_1HB} \equiv \widehat{A_1^{'}AC}$ (as angles with the sides perpendicular), it results

$$\left(\widehat{A_1HB}\right) = A - \alpha,\ m\widehat{A_1HC} = m\widehat{A_1HB} + m\widehat{BHC} = 180^\circ - \alpha,$$

It results

$$\frac{A_1B}{A_1C}=\frac{\sin(A-\alpha)}{\sin\alpha}\cdot\frac{HB}{HC}$$

Similarly, we find

$$\frac{B_1C}{B_1A}=\frac{\sin(B-\beta)}{\sin\beta}\cdot\frac{HC}{HA} \text{ and } \frac{C_1A}{C_1B}=\frac{\sin(C-\gamma)}{\sin\gamma}\cdot\frac{HA}{HB}$$

We, then apply the Menelaus' theorem.

The proof is similar for the case when the triangle is obtuse.

**53.**

We'll prove firstly the following lemma:

**Lemma**

If $(BC)$ is a cord in the circle $\mathscr{C}(O',r')$ is tangent in the point $A$ to the circle $(O)$ and in the point $D$ to the cord $(BC)$, then the points $A,D,E$ where $E$ is the middle of the arch $BC$ which does not contain the point $A$ are collinear. More than that we have $ED\cdot EA=EB^2\cdot EC^2$.

**Proof**

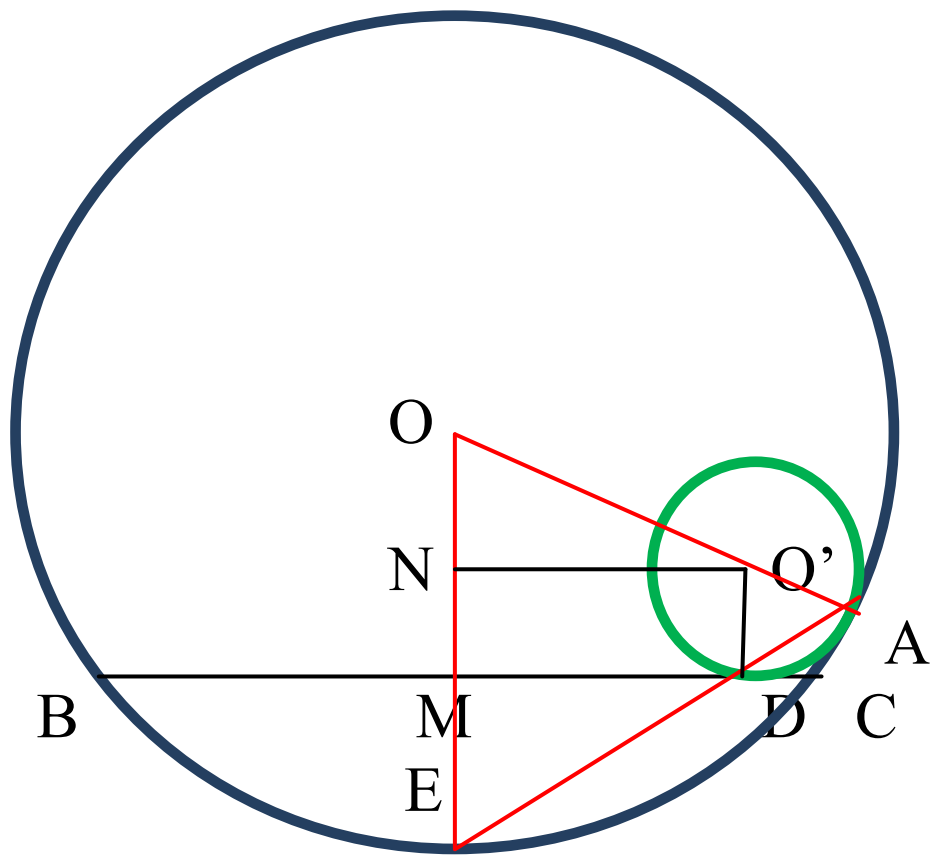

The triangles $AOE, AO'D$ are similar, because the points $A,O',O$ are collinear $OE$ is parallel with $O'D$ (see the figure above), and $\frac{OA}{O'A}=\frac{OE}{O'D}=\frac{R}{r'}$.

It results then that $\sphericalangle OAD\equiv\sphericalangle OAE$, therefore $A,D,E$ are collinear.

We'll note with $N$ the projection of $O'$ on $OE$ and $M$ the intersection point between $BC$ and $OE$, and $x=ON$.

We have $NM=r'$, $ME=R-x-r'$, $NE=R-x$, $OO'=R-r'$.

$ED\cdot EA$ represents the power of the point $E$ in rapport to the circle $\mathscr{C}(O',r')$ and it is equal to $EO'^2-r'^2$.

Because $EO'^2=O'N^2+NE^2=OO'^2-O'N^2+(R-x)$ we obtain

$$EO'^2 - r'^2 = 2R(R - r' - x).$$

On the other side $EC^2 = MC^2 + ME^2 = OC^2 + OM^2 + ME^2$ or $EC^2 = MC^2 = R^2 + (x + r')^2 + (R - r' - x)^2$

We find that $EC^2 = 2R(R - r' - x)$

Therefore, $EB^2 = EC^2 = ED \cdot EA$, and the Lemma is proved.

The proof of the theorem of P. Yiu

We note $A_1$ the intersection n point of the line $BC$ with the line $AX$ and with $X_2, X_3$ the points of tangency of the lines $AC, AB$ respectively with the inscribed circle A-mix-linear (see the figure below).

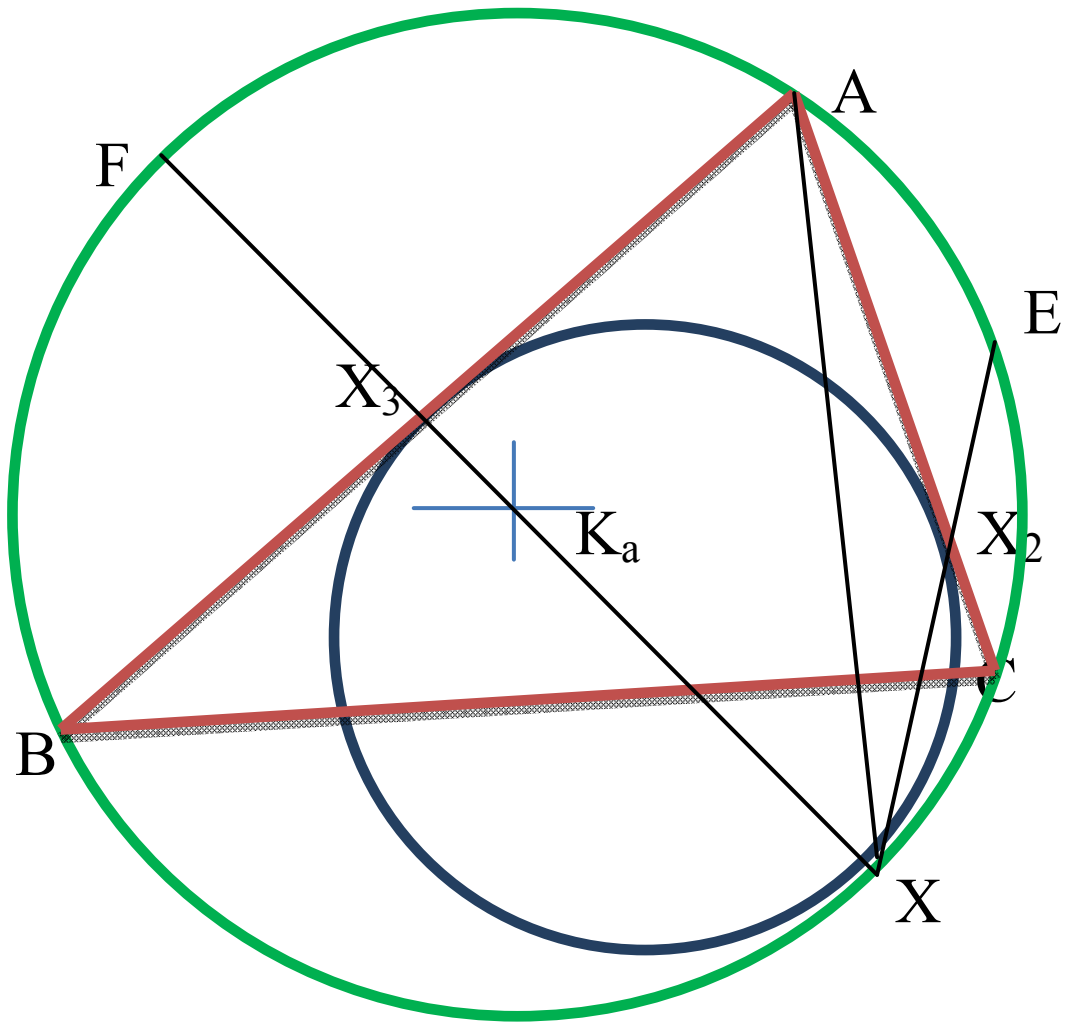

We have

$$\frac{BA_1}{CA_1} = \frac{areaBA_1X}{areaCA_1X} = \frac{BX \cdot XA_1 \cdot \sin BXA_1}{CX \cdot XA_1 \cdot \sin CXA_1} = \frac{BX \cdot \sin BXA_1}{CX \cdot \sin CXA_1} \qquad (1)$$

According to the lemma, the points $X, X_2, E$ (the middle of the arc $\overset{\frown}{AC}$) and the points $X, X_3, F$ (the middle of the arch $\overset{\frown}{AB}$) are collinear, consequently, $XX_2$ and $XX_3$ are bisectors in the triangles $AXC$ respectively $AXB$.

The bisector's theorem applied in these triangles leads to

$$\frac{AX_2}{CX_2} = \frac{AX}{CX}, \quad \frac{AX_3}{BX_3} = \frac{AX}{BX}.$$

From these relations we obtain

$$\frac{BX}{CX} = \frac{BX_3}{CX_2}.$$

We'll substitute in the relation (1), and we find:

$$\frac{BA_1}{CA_1} = \frac{BX_3}{CX_2} \cdot \frac{\sin C}{\sin B} \qquad (2)$$

But $BX_3 = AX - AX_3$, $CX_2 = AC - AX_2$, $AX_2 = \frac{AI}{\cos\frac{A}{2}}$, $AI = \frac{r}{\sin\frac{A}{2}}$.

It results

$$AX_2 = \frac{r}{\cos\frac{A}{2}\sin\frac{A}{2}} = \frac{2r}{\sin A}$$

From the sinus' theorem we retain that $\sin A = \frac{a}{2R}$.

We obtain $AX_2 = \frac{4Rr}{a}$ (it has been taken into consideration that $4RS = abc$ and $S = pr$ ). Therefore $AX_2 = \frac{bc}{p}$.

The relation (2) becomes:

$$\frac{BA_1}{CA_1} = \frac{(p-b)c^2}{(p-c)b^2}$$

Similarly, we find $\frac{CB_1}{AB_1} = \frac{(p-c)a^2}{(p-a)c^2}$, $\frac{AC_1}{BC_1} = \frac{(p-a)b^2}{(p-b)a^2}$

Because $\frac{BA_1}{CA_1}\cdot\frac{CB_1}{AB_1}\cdot\frac{AC_1}{BC_1} = 1$, in conformity to the reciprocal Ceva's theorem, it result that the Cevians $AX, BY, CZ$ are concurrent. The coordinates of their point of concurrency are the barycentric coordinates : $\frac{a^2}{(p-a)};\frac{b^2}{(p-b)};\frac{c^2}{(p-c)}$ and it is $X(56)$ in the Kimberling list.

The point $X(56)$ is the direct homothety center of the inscribed and circumscribed circles of the given triangle.

**54**.

Solution given by Gh. Ţiţeica

Let $A'B'C', A"B"C"$ the circumscribed triangles of the given triangle $ABC$ and homological with $ABC$.

Their homological centers being $M$ and respectively $N$ (the triangles $A'B'C', A"B"C"$ are called antipedal triangles of the points $M$, $N$ in rapport to $ABC$.

We note $A_1, A_2$ the intersections with $BC$ of the lines $AA', AA"$.

We have

$$(AA_1MA') = (AA_2NA') = -1$$

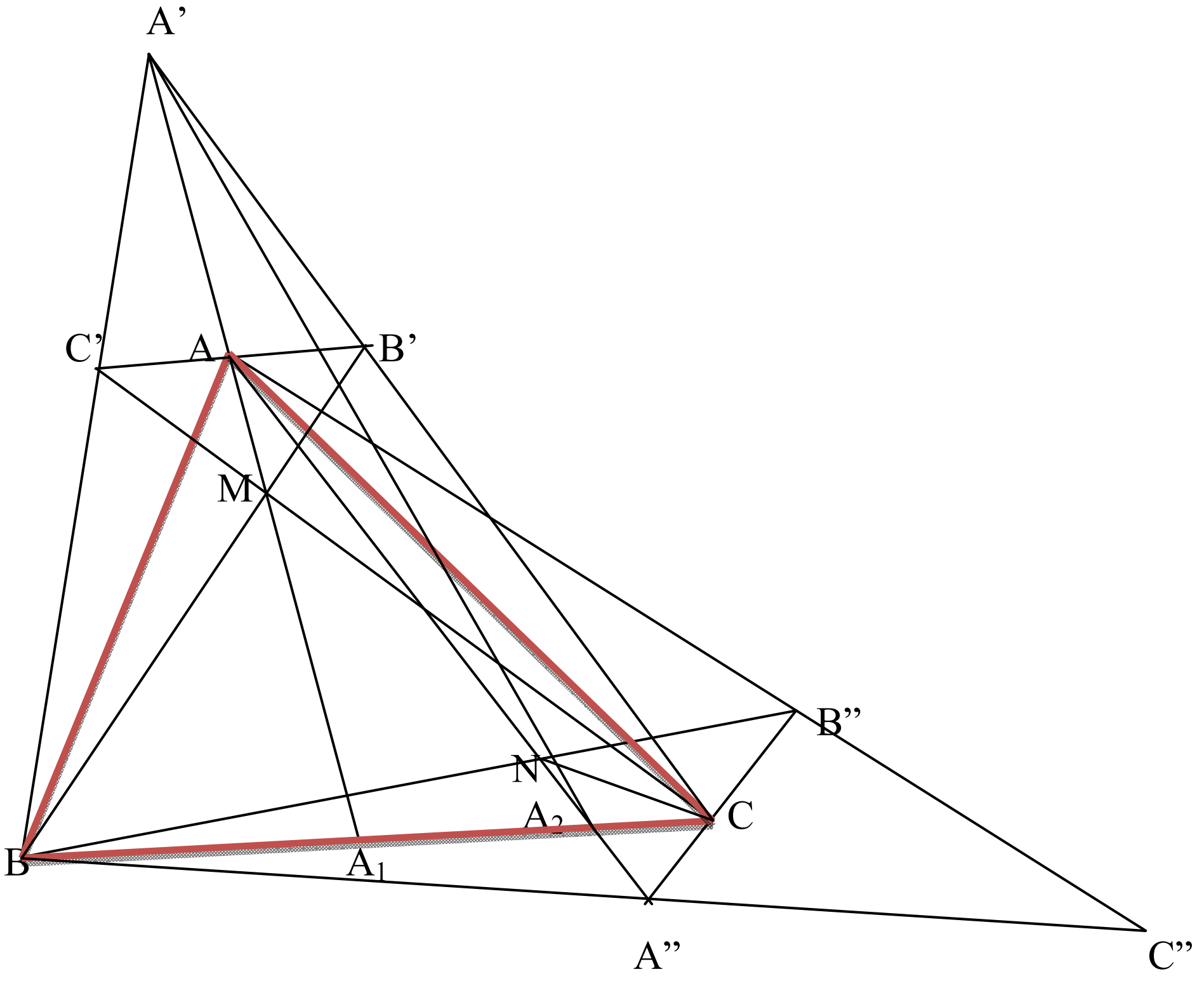

The line $A'A"$ intersects the line $BC$ in a point of the line $MN$.

**55**.

The triangles $B_1C_1A, BPC$ are orthological because the perpendiculars constructed from the vertexes of the first triangle on the sides of the second triangle are concurrent.

Indeed, the perpendicular constructed from $B_1$ on $BP$, the perpendicular constructed from $A$ on $BC$ and the perpendicular constructed from $C_1$ on $CP$ are concurrent in $H$. The reciprocal is also true: the perpendicular constructed from $B$ on $AB_1$, the perpendicular constructed from $AB_1$ on $AC_1$ and the perpendicular constructed from $P$ on $B_1C_1$ are concurrent.

Because the first two perpendiculars are the heights of the triangle $ABC$, it means that also the third perpendicular passes through $H$, therefore $PH$ is perpendicular on $BC$, it results that $BC \parallel B_1C_1$.

**56**.

We'll show firstly that

$$\frac{VM \cdot UN}{UV^2} = \frac{WP \cdot WQ}{UW^2}$$

The sinus' theorem in the triangles $MUV, NUV$ leads to

$$\frac{UV}{\sin\alpha} = \frac{VM}{\sin\alpha_1}, \frac{UV}{\sin\beta} = \frac{VN}{\sin\alpha_2}$$

(see the figure bellow)

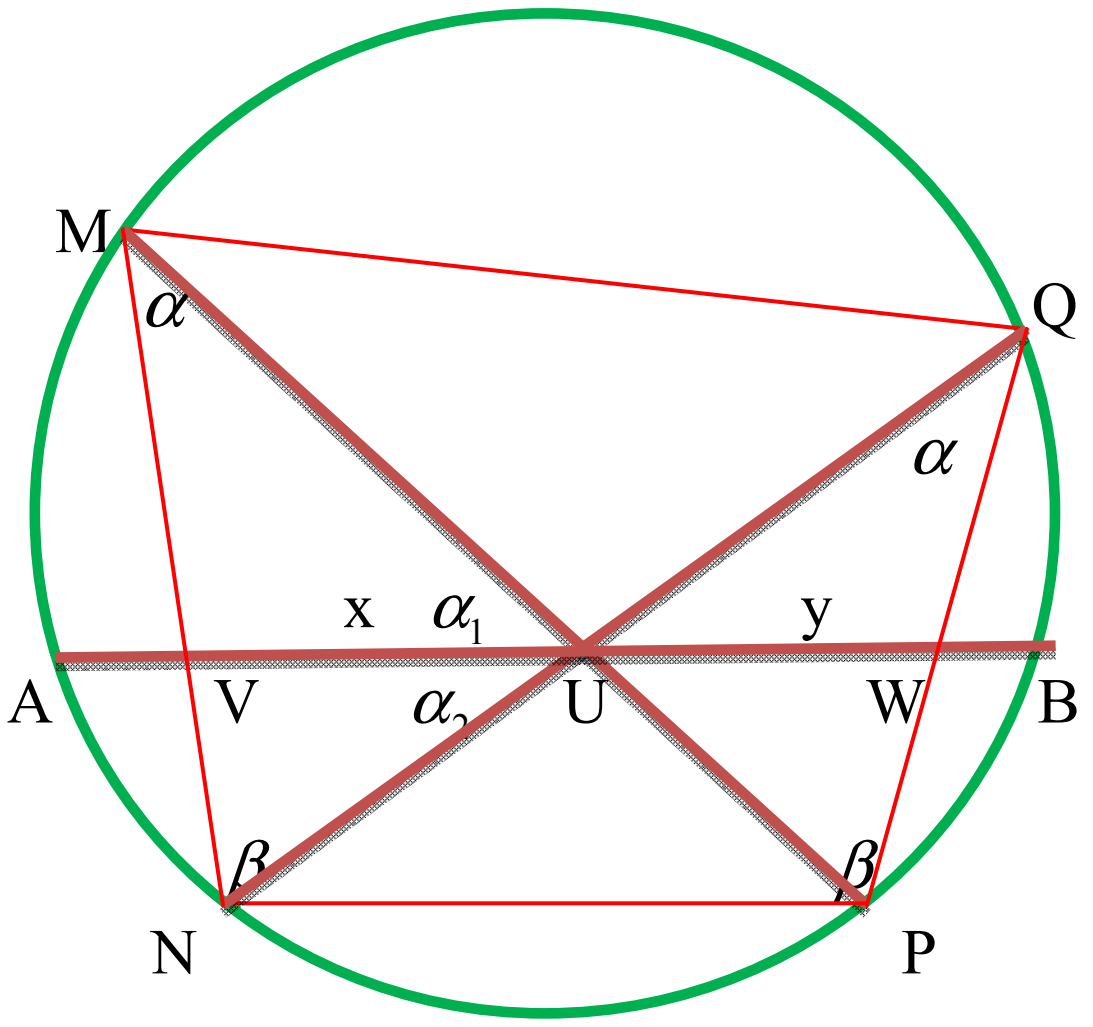

Therefore,

$$\frac{VM \cdot UN}{UV^2} = \frac{\sin\alpha_1 \sin\alpha_2}{\sin\alpha \text{in}\,\beta}.$$

Similarly, from the triangles $UWQ, UWP$

we find

$$\frac{WP \cdot WQ}{UW^2} = \frac{\sin\alpha_1 \sin\alpha_2}{\sin\alpha \text{in}\,\beta}.$$

If we note $UA = UB = a$, $UV = x$, $UW = y$, and considering the power of $U$ in rapport with the circle, we have:

$$VM \cdot UN = (a - x)(a + x) = a^2 - x^2$$

Considering the power of $W$ in rapport with the circle, we have $WP \cdot WQ = a^2 - y^2$

From the relation proved above we have $\frac{a^2 - x^2}{x^2} = \frac{a^2 - y^2}{y^2}$, and from here it results $x = y$.

**57.**

We will note $c_1 = (A_2A_3A_4)$ and $\mu_1$ the power of $M$ in rapport with the circle $c_1$, etc.

The set of the points $M$ whose powers $\mu_i$ in rapport to the circles $c_i$, which satisfy the linear relation

$$\frac{\mu_1}{p_1}+\frac{\mu_2}{p_2}+\frac{\mu_3}{p_3}+\frac{\mu_4}{p_4}=1 \qquad (1)$$

is a circle whose equation is satisfied by the points $A_1, A_2, A_3, A_4$.

For example for the point $A_1$ we have $\mu_1 = p_1$ and $\mu_2 = \mu_3 = \mu_4 = 0$. The points $A_i$ being, hypothetical, arbitrary, it results that (1) is an identity, which is satisfied by any point from the plane. In particular, it takes place when the point $M$ is at infinite; we have $\mu_1 = MO_i^2 - r_i^2$, where $O_i$ and $r_i$ are the center and the radius of the circle $c_i$. We'll divide by $\mu_1$ and because $\mu_i : \mu_1 \to 1$ when $M \to \infty$, the relation (1) is reduced to

**58.**

Let $\{I\} = CC' \cap BB'$; The triangle $ABB'$ is isosceles, therefore, $\sphericalangle ABB' = \sphericalangle AB'B$.

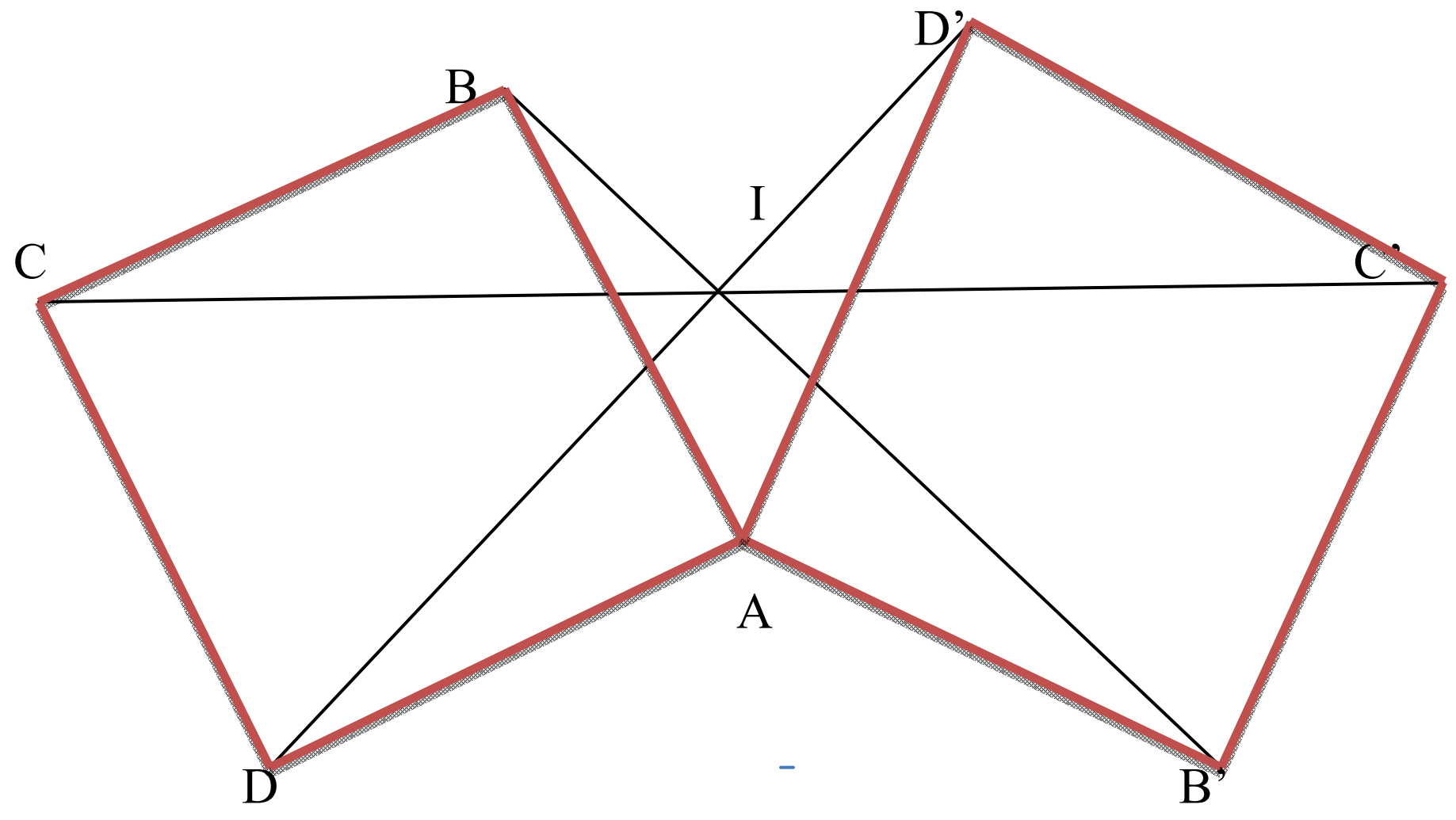

We note $m(\sphericalangle ABB') = \alpha$. Then we have:

$$\frac{c_i}{\sin(90^\circ+\alpha)}=\frac{BC}{\sin BIC}=\frac{B'C'}{\sin B'IC'}=\frac{C'I}{\sin(90^\circ-\alpha)}$$

But

$$\sin(90^\circ+\alpha)=\sin(90^\circ-\alpha)=\cos\alpha$$

Then we deduct that $CI = C'I$.

We'll note $\{I'\} = CC' \cap DD'$ and $m(\sphericalangle ADD') = m(\sphericalangle AD'D) = \beta$.

We have

$$\frac{CI'}{\sin\left(90^\circ-\beta\right)}=\frac{CD}{\sin CI'D}=\frac{O'D'}{\sin C'ID'}=\frac{C'I'}{\sin\left(90^\circ+\beta\right)}$$

It results that $CI'=C'I'$, and we deduct that $I=I'$, which point is the middle of the segment $CC'$, consequently $CC'\cap DD'\cap BB'=\{I\}$.

**Observation**

The problem is true also when the squares are not congruent. The intersection point is the second point of intersection of the circles circumscribed to the squares.

**59.**

We note $m(\sphericalangle A_1AC)=\alpha$, $m(\sphericalangle B_1BC)=\beta$, $m(\sphericalangle C_1CB)=\gamma$.

Because $AB_1, BA_1, CA_1$ are concurrent we can write

$$\frac{\sin\alpha}{\sin(A-\alpha)}\cdot\frac{\sin(B+60^\circ)}{\sin 60^\circ}\cdot\frac{\sin 60^\circ}{\sin(C+60^\circ)}=-1$$

Similarly

$$\frac{\sin\beta}{\sin(B-\beta)}\cdot\frac{\sin(C+60^\circ)}{\sin 60^\circ}\cdot\frac{\sin 60^\circ}{\sin(A+60^\circ)}=-1$$

$$\frac{\sin\gamma}{\sin(B-\gamma)}\cdot\frac{\sin(A+60^\circ)}{\sin 60^\circ}\cdot\frac{\sin 60^\circ}{\sin(B+60^\circ)}=-1$$

We multiply these three relations and we find

$$\frac{\sin\alpha}{\sin(A-\alpha)}\cdot\frac{\sin\beta}{\sin(B-\beta)}\cdot\frac{\sin\gamma}{\sin(C-\gamma)}=-1$$

which shows that the lines $AA_1, BB_1, CC_1$ are concurrent.

**60.**

a) Let $P$ the intersection point of the perpendicular in $A'$ on $BC$ with the perpendicular in $B'$ on $AC$.

We have $PB^2 - PC^2 = A'B^2 - A'C^2$, $PC^2 - PA^2 = B'C^2 - B'A^2$

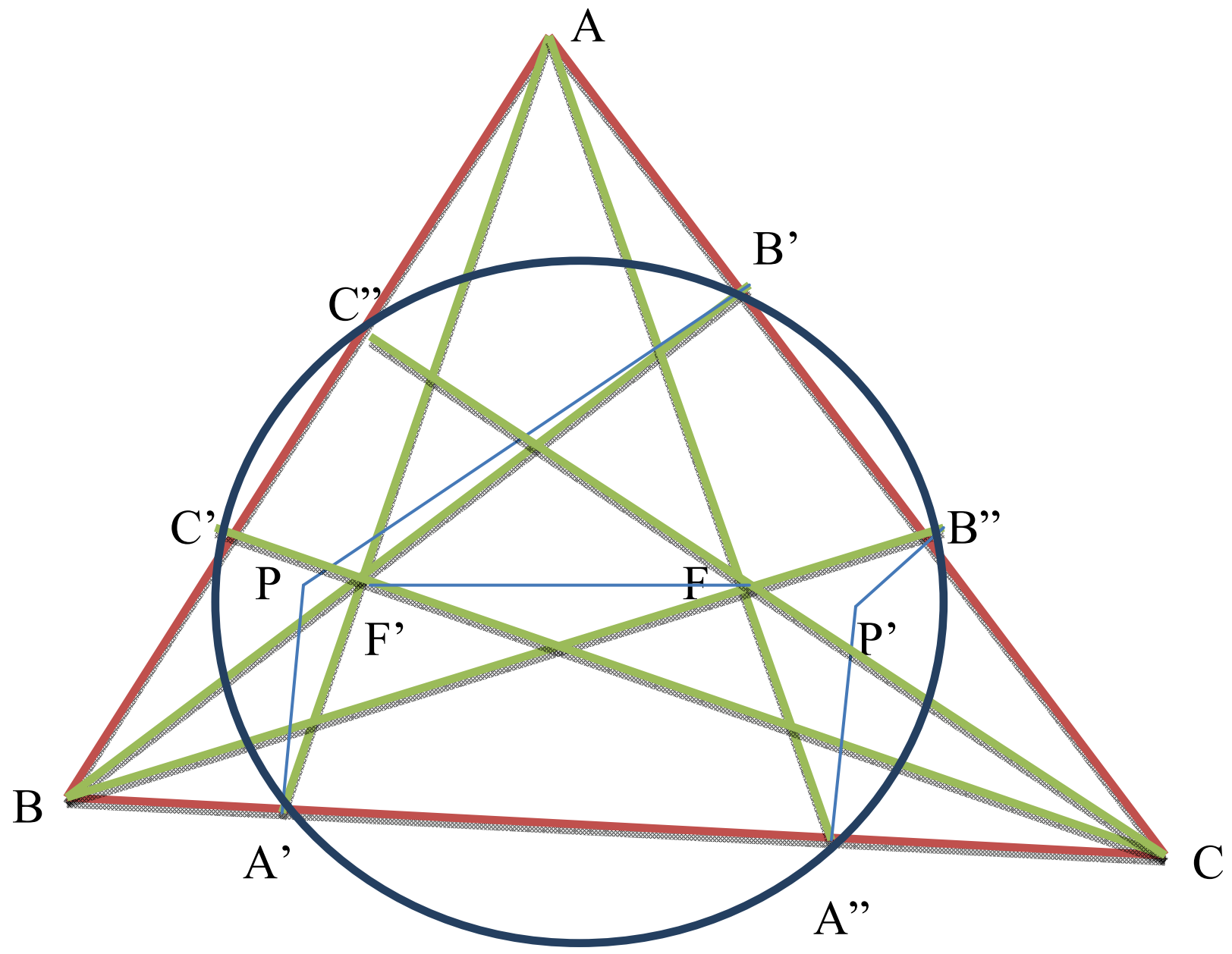

Adding side by side these relations, it results

$$PB^2 - PA^2 = A'B^2 - A'C^2 + B'C^2 + B'A^2 \qquad (1)$$

If $C_1$ is the projection of $P$ on $AB$, then

$$PB^2 - PA^2 = C_1B^2 - C_1A^2 \qquad (2)$$

From (1) and (2) it results $C_1 \equiv C'$

b) Let $A_1, B_1, C_1$ the projections of the points $A, B, C$ on $B'C', C'A', A'B'$ respectively.

We have

$$A_1C'^2 - A_1B'^2 = C'A^2 - B'A^2$$

$$B_1C'^2 - B_1A'^2 = C'B^2 - A'B^2$$

$$C_1A'^2 - C_1B'^2 = A'C^2 - B'C^2$$

From these relations it results $A_1C'^2 + B_1A'^2 + C_1B'^2 = A_1B'^2 + B_1C'^2 + C_1B'^2$ which is a relation of the type from a) for the triangle $A'B'C'$.

Using a similar method, it results that triangle $A_1B_1C_1$ is the pedal triangle of $P'$.

c) The quadrilateral $AB'PC'$ is inscribable, therefore $\sphericalangle APB' \equiv \sphericalangle AC'B'$ and because these angles have as complements the angles $\sphericalangle CAP, \sphericalangle B'AP'$, it results that these angles are also congruent, therefore the Cevians $AP, AP'$ are isogonal.

d) We observe that the mediators of the segments $AP, AP'$ pass through the point $F$, which is the middle of the segment $PP'$.

We show that $F$ is the center of the circle that contains the points from the given statement.

We'll note $m(\sphericalangle P'AC) = m(\sphericalangle PAB) = \alpha$, $AP = x, AP' = x'$.

We'll use the median's theorem in the triangles $C'PP', B'PP'$ to compute $C'F, B'F$.

$$4C'F^2 = 2\left(PC'^2 + P'C'^2\right) - PP'^2$$

$$4B'F^2 = 2\left(PB'^2 + P'B'^2\right) - PP'^2$$

$$PC' = x\sin\alpha,\ P'C'^2 = P'C''^2 + C''C'^2,\ P'C'' = x'\sin(A-\alpha)$$

$$AC'' = x'\cos(A-\alpha),\ AC' = x\cos\alpha$$

$$P'C'^2 = x'^2 + \sin^2(A-\alpha) + \left(x'\cos(A-\alpha) - x\cos\alpha\right)^2 =$$

$$= x'^2 + x^2\cos^2\alpha - 2xx'\cos\alpha\cos(A-\alpha)$$

$$4C'F^2 = 2\left[x'^2 + x^2\cos^2\alpha - 2xx'\cos\alpha\cos(A-\alpha)\right] - PP'^2$$

$$4C'F^2 = 2\left[x'^2 + x^2 - 2xx'\cos\alpha\cos(A-\alpha)\right] - PP'^2$$

Similarly, we find the expression for $B'F^2$, and it will result that $C'F = B'F$ ,therefore $C', C'', B'', B'$ are concyclic.

e) We'll consider the power of the points $A, B, C$ in rapport with the circle of the points $A', A'', B', B'', C', C''$.

We have

$$AB' \cdot AB'' = AC' \cdot AC''$$
$$BA' \cdot BA'' = BC' \cdot BC''$$
$$CA' \cdot CA'' = CB' \cdot CB''$$

Multiplying these relations and using the reciprocal of the Ceva's theorem, it results the concurrency of the lines $AA'', BB'', CC''$.

**61**.

Follow the same method as for problem 59. Instead of 60 we'll use 45 .

**62**.

Let $O$ the center of the circumscribed circle and $H$ the projection of $x$ on $BC$. The quadrilateral $KHFP$ is inscribable, from here it results that

$$m(\sphericalangle KHP) = m(\sphericalangle KFP) = \frac{1}{2}m\left(\overset{\frown}{PQ}\right) = m(\sphericalangle KOP).$$

On the other side $\Delta AOP \equiv \Delta AOQ$, then $m(\sphericalangle AOP) = \frac{1}{2}m(\sphericalangle QOP)$. We obtain that $\sphericalangle KHP \equiv \sphericalangle AOP$ and then their complements are congruent, therefore $m(\sphericalangle AHO) = 90°$. It results that the quadrilateral $AHOP$ is an inscribable quadrilateral, therefore $A, K, H$ are collinear.

**63**.

From the bisectrices' theorem we find $\frac{DB}{DA}=\frac{BC}{AC}$, $\frac{EA}{EB}=\frac{AB}{BC}$ and the given relation implies $\frac{BC}{AC}=\frac{AB}{BC}$, we deduct that

$$\frac{DB}{DA}=\frac{EA}{EC} \qquad (1)$$

We'll use the reciprocal of the Menelaus theorem in the triangle $ADE$ for the transversal $M-P-N$. We have to compute $\frac{MD}{MA}\cdot\frac{NA}{NE}\cdot\frac{PE}{PD}$. The point $P$ is the middle of $(DE)$, therefore $\frac{PE}{PD}=1$. From (1) we find $\frac{DA}{MA}=\frac{EC}{NA}$, therefore $\frac{MA-MD}{MA}=\frac{CN-NE}{NA}$, then $\frac{MD}{MA}=\frac{NE}{NA}$.

Therefore, $\frac{MD}{MA}\cdot\frac{NA}{NE}\cdot\frac{PE}{PD}=1$ and $M,P,N$ are collinear.

**64.**

We'll note $m(\sphericalangle MAB)=\alpha$, $m(\sphericalangle MBC)=\beta$, $m(\sphericalangle MCA)=\gamma$,

$$\{A_1'\}=BC\cap AA_1,\ \{B_1'\}=AC\cap BB_1,\ \{C_1'\}=AB\cap CC_1$$

We have

$$\frac{A_1'B}{A_1'C}=\frac{Aria(ABA_1)}{Aria(ACA_1)}=\frac{AB\cdot\sin(60^\circ+\beta)\cdot A_1B}{AC\cdot\sin(120^\circ-\gamma)\cdot A_1C}$$

But $AB=AC, A_1B=MB, A_1C=MC, \sin(120^\circ-\gamma)=\sin(60^\circ+\gamma)$.

Therefore, $\frac{A_1'B}{A_1'C}=\frac{\sin(60^\circ+\beta)}{\sin(60^\circ+\gamma)}\frac{MB}{MC}$.

Similarly $\frac{B_1'C}{B_1'A}=\frac{\sin(60^\circ+\gamma)}{\sin(60^\circ+\alpha)}\frac{MC}{MA}$ and $\frac{C_1'A}{C_1'B}=\frac{\sin(60^\circ+\alpha)}{\sin(60^\circ+\beta)}\frac{MA}{MB}$

Then we apply the reciprocal of the Ceva's theorem.

**65**.

i) Let $L,Q$ the intersection points of the circles $\mathcal{C}(C,AB),\mathcal{C}(B,AC)$.

We have

$$\Delta ABC\equiv\Delta LCB \text{ (S.S.S)},$$

it results

$$\sphericalangle CLB=90^\circ \text{ and } \sphericalangle LCB\equiv\sphericalangle ABC$$

which leads to $LC \parallel AB$.

Having also $LC=AB$, $m(\sphericalangle CAB)=90^\circ$,

we obtain that the quadrilateral $ABLC$ is a rectangle, therefore $AL=BC$ and $L$ ,belongs to the circle $\mathcal{C}(A,BC)$.

ii) $\Delta CQB \equiv \Delta BAC$ (S.S.S), it results that $\sphericalangle CQB = 90°$.

Because $ABLC$ is a rectangle that means that $A, B, L, C$ are on the circle with the center in $O$, the middle of $BC$.

The triangle $BQC$ is a right triangle, we have

$$QO = \frac{1}{2} BC = OC = OB,$$

therefore the point $Q$ is on the circle of the points $A, B, L, C$.

iii) $\Delta OCA \equiv \Delta ABR$, it results that

$$m(\sphericalangle PAC) + m(\sphericalangle CAR) + m(\sphericalangle BAR) = 180°,$$

consequently, the points $P, A, R$ are collinear.

But

$$m(\sphericalangle PQL) = 90°,\ m(\sphericalangle LQR) = 90° \ (L, B, R \text{ collinear}).$$

It results

$$m(\sphericalangle QPR) = 180°,$$

therefore $P, Q, R$ are collinear.

We saw that $P, A, R$ are collinear, we deduct that $P, Q, A, R$ are collinear.

**66**.

If $H$ is the orthocenter of the triangle $AB'C'$, from the sinus' theorem we have

$$\frac{B'C'}{\sin A} = AH.$$

But $AH = 2R\cos A$

We obtain

$$B'C' = 2R\cos A \sin A \text{ or } B'C' = 2R\sin 2A.$$

Because $2R = \frac{\alpha}{\sin A}$, we find

$$a' = B'C' = R\cos A.$$

From the cosine's theorem we have

$$\cos A = \frac{b^2 + c^2 - a^2}{2bc}$$

$$4(a'b' + b'c' + c'a') =$$

$$= 4\left[\frac{ab(b^2+c^2-a^2)(a^2+c^2-b^2)}{4abc^2} + \frac{bc(a^2+c^2-b^2)(a^2+b^2-c^2)}{4a^2bc} + \frac{ac(b^2+c^2-a^2)(a^2+b^2-c^2)}{4ab^2c}\right]$$

We deduct that

$$4(a'b' + b'c' + c'a') = \frac{c^4 - (b^2 - a^2)^2}{c^2} + \frac{a^4 - (c^2 - b^2)^2}{a^2} + \frac{b^4 - (c^2 - a^2)^2}{b^2} =$$

$$= a^2 + b^2 + c^2 - \left[\left(\frac{b^2 - a^2}{c}\right)^2 + \left(\frac{c^2 - b^2}{a}\right)^2 + \left(\frac{c^2 - a^2}{b}\right)^2\right] \le a^2 + b^2 + c^2.$$

We observe that if the triangle $ABC$ is equilateral, then we have an equality in the proposed inequality.

**67**.

Let $\{M\} = AG \cap BC$, we consider the homothety $h_G^{-\frac{1}{2}}$ the image of the circumscribed circle of the triangle $ABC$ (the circle circumscribed to the medial triangle). On the circle of the nine points is also the point $D$, therefore it is the homothetic of a precise point on the circumscribed circle. This point is exactly the point $P$, therefore

$$GD = \frac{1}{2} G'P .$$

**68**.

If $M$ is the middle of the cord $AB$, $N$ is the middle of the cord $CD$ and $P$ is the middle of the segment $EF$ then $M, N, P$ are collinear (the Newton-Gauss line of the complete quadrilateral $DACBFE$.)

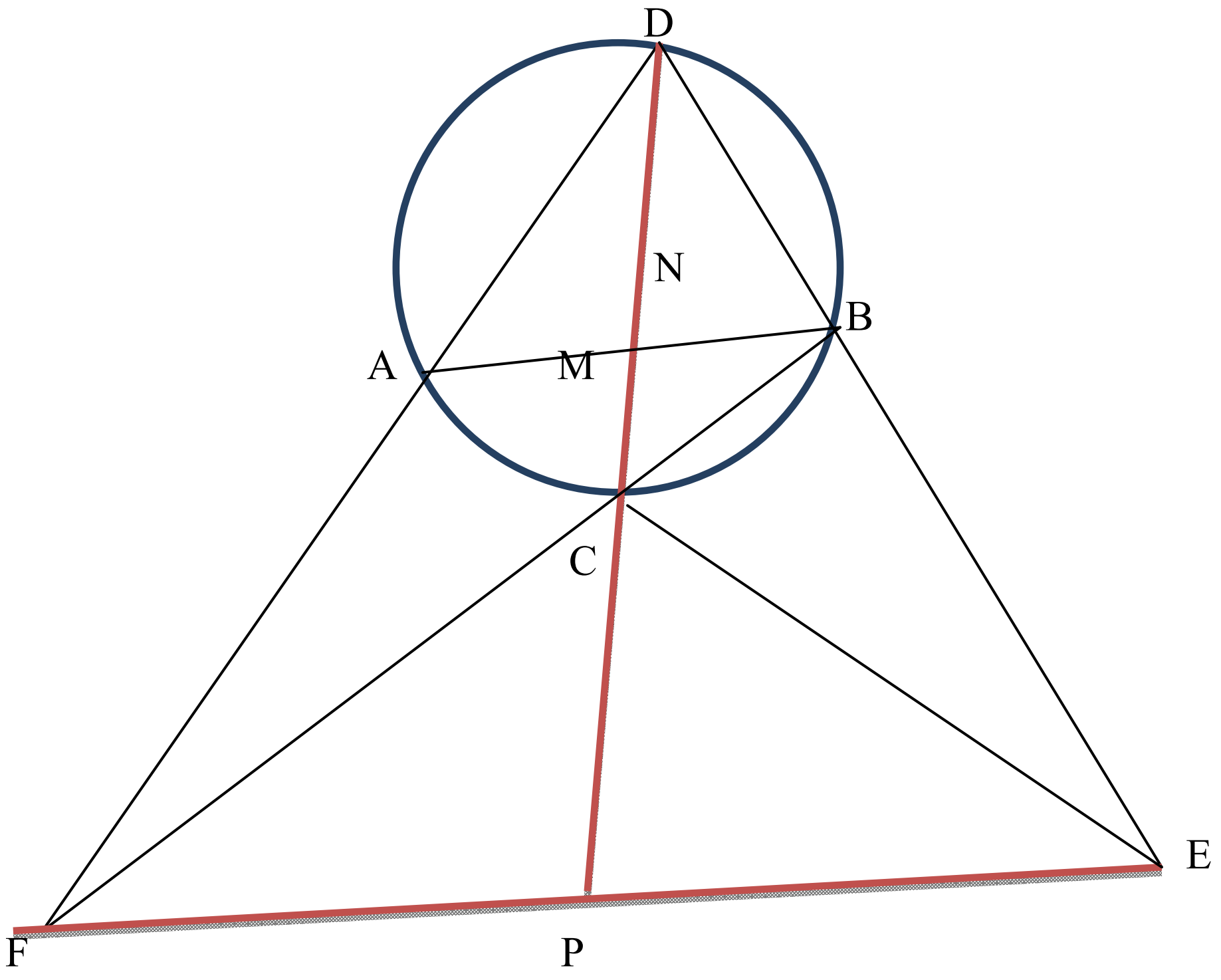

(see the figure above). It results that in the triangle $DEF$, $DP$ is median. But $C$ being on this median, we show then that Ceva's theorem $i_c$ from $DP, FB, EA$ concurrency, it results $AB \parallel EF$.

**69.**

This is the Pappus' theorem correlative.

**70**.

The Menelaus' theorem in the triangles $AA'C, AA''B$ for the transversals $B-M-Q,\ C-N-P$ lead us to

$$\frac{BA'}{BC}\cdot\frac{MA}{MA'}\cdot\frac{QC}{QA}=1 \qquad (1)$$

$$\frac{CA''}{CB}\cdot\frac{PB}{PA}\cdot\frac{NA}{NA''}=1 \qquad (2)$$

Because $\frac{PB}{PA}=\frac{QC}{QA}$ and $BA'=CA''$ from the equality of the relations (1) and (2), it results $\frac{MA}{MA'}=\frac{NA}{NA''}$, which implies $MN \parallel BC$.

**71**.

Applying the Van Aubel theorem for the triangle $ABC$ we have:

$$\frac{AP}{PA'}=\frac{AC'}{CB'}+\frac{AB'}{B'C} \qquad (1)$$

$$\frac{BP}{PB'}=\frac{BA'}{A'C}+\frac{BC'}{C'A} \qquad (2)$$

$$\frac{CP}{PC'}=\frac{CA'}{A'B}+\frac{CB'}{B'A} \qquad (3)$$

We'll note $\frac{AC'}{CB'}=x>0$, $\frac{AB'}{B'C}=y>0$, $\frac{BA'}{A'C}=z>0$, then we obtain

$$E(P)=\left(x+\frac{1}{x}\right)+\left(y+\frac{1}{y}\right)+\left(z+\frac{1}{z}\right)\geq 2+2+2=6$$

The minimum value will be obtained when $x=y=z=1$, therefore when $P$ is the weight center of the triangle $ABC$.

Multiplying the three relations we obtain:

$$E(P)=\left(x+\frac{1}{x}\right)+\left(y+\frac{1}{y}\right)+\left(z+\frac{1}{z}\right)+\frac{yz}{x}+\frac{x}{yz}\geq 8$$

The minimum value will be obtained when $x=y=z=1$, therefore when $P$ is the weight center of the triangle $ABC$.

**72.**

i) We apply the Menelaus' theorem in the triangle $A_1B_1C_1$ for the transversals $P_1-Q_3-R_2,\ P_2-Q_1-R_3,\ P_3-Q_2-R_1$, we obtain

$$\frac{P_1B_1}{P_1C_1}\cdot\frac{R_2C_1}{R_2A_1}\cdot\frac{Q_3A_1}{Q_3B_1}=1 \quad (1)$$
$$\frac{P_2C_1}{P_2A_1}\cdot\frac{Q_1B_1}{Q_2C_1}\cdot\frac{R_3A_1}{R_3B_1}=1 \quad (2)$$
$$\frac{P_3A_1}{P_3B_1}\cdot\frac{R_1B_1}{R_1C_1}\cdot\frac{Q_2C_1}{Q_2A_1}=1 \quad (3)$$

Multiplying relations (1), (2), and (3) side by side we obtain the proposed relation.

ii) If $A_1A_2; B_1B_2; C_1C_2$ are concurrent then $P_1, P_2, P_3$ are collinear and using the Menelaus theorem in the triangle $A_1B_1C_1$ gives

$$\frac{P_1B_1}{P_1C_1}\cdot\frac{P_2C_1}{P_2A_1}\cdot\frac{P_3A_1}{P_3B_1}=1 \quad (4)$$

Taking into account the relation from i) and (4) it will result relation that we are looking for.

Reciprocal

If the relation from ii) takes place, then substituting it in i) we obtain (4). Therefore, $A_1B_1C_1$, $A_2B_2C_2$ are homological, and their homological axis being $P_1P_2P_3$.

**73**.

The fixed point is the harmonic conjugate of the point $C$ in rapport with $A, B$ because $T_1T_2$ is the polar of $C$ in rapport with the circle which passes through $A, B$.

**74.**

We suppose that $AA_1$ is median, therefore $A_1B = A_1C$ and from the given relation in hypothesis, we obtain

$$A_1C^2 + C_1A^2 = AB_1^{\ 2} + BC_1^{\ 2}$$

From the concurrency of the Cevians $AA_1, BB_1, CC_1$ and the Ceva's theorem we retain that

$$\frac{AB_1}{B_1C} = \frac{AC_1}{C_1B} \quad (2)$$

We'll note $\frac{AB_1}{B_1C} = k$, $k > 0$, then

$$B_1C^2 + k^2C_1B^2 = k^2BC_1^{\ 2} + BC_1^{\ 2}$$

It follows that

$$\left(k^2 - 1\right)C_1B^2 - BC_1^{\ 2} = 0.$$

In order o have equality it results that $k = 1$ or $C_1B = B_1C$.

If $k = 1$ then $AB_1 = B_1C$ and $AC_1 = C_1B$, which means that $BB_1, CC_1$ are median.

If $C_1B = B_1C$ then from (2) we obtain $AB_1 = AC_1$ with the consequence that $AB = AC$, and the triangle $ABC$ is isosceles.

**75.**
We note

$$m(\sphericalangle C_1AB) = m(\sphericalangle B_1AC) = \alpha$$
$$m(\sphericalangle C_1BA) = m(\sphericalangle A_1BC) = \beta$$
$$m(\sphericalangle A_1CB) = m(\sphericalangle B_1CA) = \gamma$$

We have

$$\frac{BA'}{CA'} = \frac{Aria(ABA_1)}{Aria(ACA_1)} = \frac{AB \cdot BA_1 \cdot \sin(B+\beta)}{AC \cdot CA_1 \cdot \sin(c+\gamma)} = \frac{c \cdot \sin\gamma \cdot \sin(B+\beta)}{b \cdot \sin\beta \cdot \sin(c+\gamma)}$$

Similarly we compute $\frac{B'C}{B'A}, \frac{C'A}{C'B}$ and use the Ceva's theorem.

**76.**
Let $MNPQ$ a square inscribed in the triangle $ABC$.

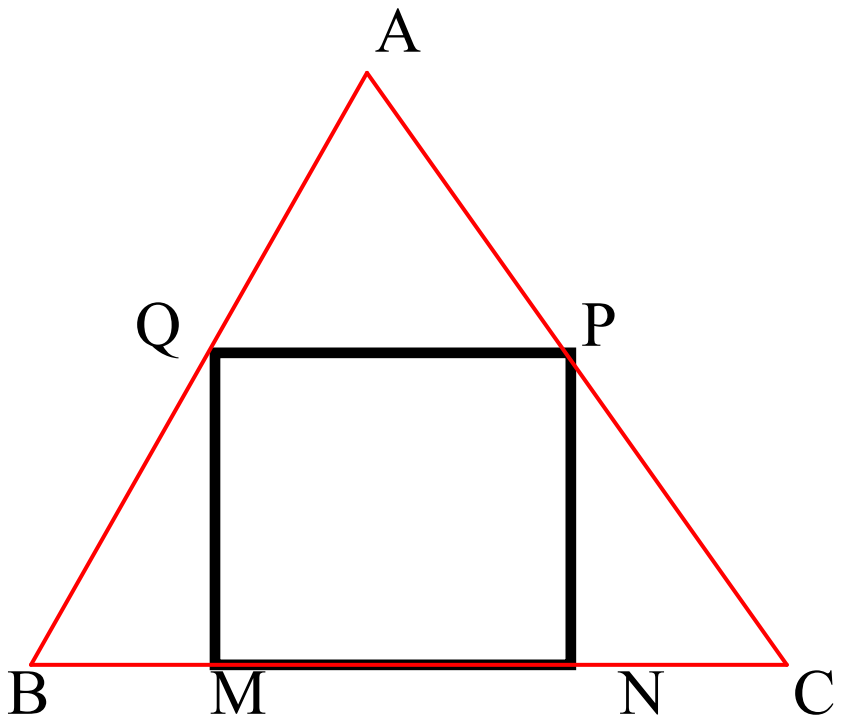

If we note the square's side with $x$ and the height from $A$ with $h_a$, we have

$$\Delta APQ \sim \Delta ACB$$

Therefore

$$\frac{x}{a} = \frac{h_a - x}{h_a}$$

But $ah_a = 2S$, where $S$ is the aria of the triangle $ABC$.
Therefore

$$x = \frac{2S}{a+h_a}.$$

Similarly we find $y = \frac{2S}{b+h_b}$, $z = \frac{2S}{c+h_c}$ the size of the inscribed squares.

From $x = y = z$ it results

$$a + h_a = b + h_b = c + h_c \qquad (1)$$

By squaring the relation we obtain:

$$a^2 + h_a^{\ 2} + 4S = b^2 + h_b^{\ 2} + 4S = c^2 + h_c^{\ 2} + 4S$$

Taking away $6S$ we have:

$$\left(a + h_a\right)^2 = \left(b + h_b\right)^2 = \left(c + h_c\right)^2$$

From here it results

$$\left|a - h_a\right| = \left|b - h_b\right| = \left|c - h_c\right| \qquad (2)$$

From (1) and (2) we find $a = b = c$, therefore the triangle is equilateral

**77.**

The line $NP$ is the polar of the point of intersection of the lines $AM$, $BC$, therefore $NP$ passes through the pole of the side $BC$ (the intersection of the tangents constructed in $B, C$ to the circumscribed circle of the triangle $ABC$)

**78.**

Let $\{P\} = AB \cap CD$. The line $MN$ is the polar of the point $P$ and passes through the fixed point $Q$ which is the harmonic conjugate of the point $P$ in rapport with $C, D$.

**79.**

It is known that $ON^2 = 9R^2 - \left(a^2 + b^2 + c^2\right)$

We observe that $C = 60°$. From the sinus' theorem we have that $R^2 = \frac{c^2}{3}$ and from the cosine's theorem $c^2 = a^2 + b^2 - ab$. By substituting we obtain that $ON = c - b$

**80.**

i)

$$\overrightarrow{OA} + \overrightarrow{OB} = 2\overrightarrow{OP}$$
$$\overrightarrow{OC} + \overrightarrow{OD} = 2\overrightarrow{OR}$$
$$\overrightarrow{OP} + \overrightarrow{OR} = \overrightarrow{OM}$$

The quadrilateral $OPMR$ is a parallelogram of center $T$, the middle of $PR$, therefore $2\overrightarrow{OT} = \overrightarrow{OM}$.

ii)

$OPMR$ is a parallelogram, it results that $PM \parallel OR$ and because $OR \perp CD$ it results that $PM \perp CD$, therefore $M \in PP'$, Similarly we show that $M \in QQ'$, $M \in RR'$, $M \in SS'$

**81.**

If $O_1, O_2, O_3$ are the centers of the three congruent circles, it result that $O_1O_2$, $O_2O_3$, $O_3O_1$ are parallel with the sides of the triangle
$ABC$.

It results that the triangles $ABC$, $O_1O_2O_3$ can be obtained one from the other through a homothety conveniently chosen.

Because $O_1, O_2, O_3$ belong to the interior bisectrices of the triangle $ABC$, it means that the homothety center is the point $T$, which is the center of the inscribed circle in the triangle $ABC$.

The common point of the three given circles, noted with $O'$. This point will be the center of the circumscribed circle to the triangle $O_1O_2O_3$ ($O'O_1 = O'O_2 = O'O_3$).

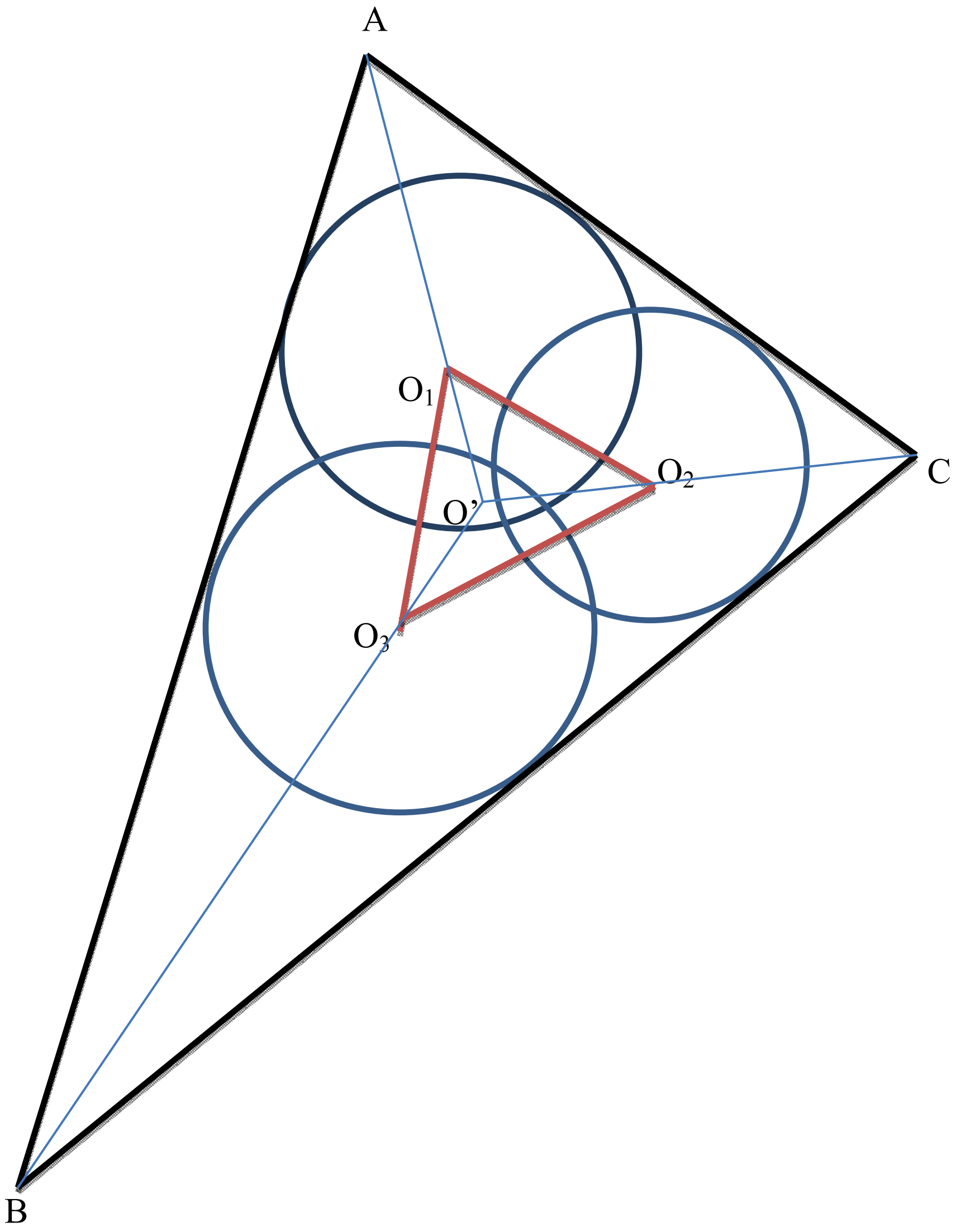

The homothety of center $I$, which transforms the triangle $ABC$ in $O_1O_2O_3$ will transform the center of the circumscribed circle $O$ of the triangle $ABC$ in the center $O'$ of the circumscribed circle to the triangle $O_1O_2O_3$, therefore $O, I, O'$ are collinear points.

**82.**

$OA$ is perpendicular on $B'C'$, the triangles $A'B'C'$, $ABC$ are homological, their homology axes is $A''B''C''$ the orthic axis of the triangle $ABC$.

The quadrilateral $AA_1A'A"$ is inscribable, the center of its circumscribed circle being $O_1$, the middle of the segment $(AA")$.

Similarly, the centers of the circumscribed circles to the triangles $BB_1B', CC_1C'$ will be $O_2, O_3$, the middle points of the segments $(BB"), (CC")$. The points $O_1, O_2, O_3$ are collinear because these are the middle of the diagonals of the complete quadrilateral $ABA"B"CC"$ (the Newton-Gauss line of the quadrilateral).

**83.**

The symmetric of the Cevian $AA_1$ in rapport with $BC$ is $A'A_1$ where $A'$ is the vertex of the anti-complementary triangle $A'B'C'$ of the triangle $ABC$ (the triangle formed by the parallels to the sides of the triangle $ABC$ constructed through $A, B, C$).

We'll use the reciprocal of Ceva's theorem.

**84.**

We will transform the configuration from the given data through an inversion of pole $O$ and of rapport $r^2$. The circles circumscribed to the triangles $B'OC'$, $C'OA'$, $A'OB'$ which transforms the sides $BC, CA, AB$ of a triangle $ABC$. The three given congruent circles are transformed respectively in the tangents constructed in $A, B, C$ to the circumscribed circle to the triangle $ABC$. The points $A_1^{'}, B_1^{'}, C_1^{'}$ have as inverse the intersection points $A_1, B_1, C_1$ of the tangents to the circumscribed circle to the triangle $ABC$ constructed in $A, B, C$ with $BC, CA, AB$. These point are, in conformity with a theorem of Carnot, collinear (the Lemoine line of the triangle $ABC$), therefore the given points from the problem $A_1^{'}, B_1^{'}, C_1^{'}$ are on a circle that contains the inversion pole $O$.

**85.**

Let $\{P\} = A'B \cap AB'$ and the middle point of $(AA')$; the points *P, M, W* are collinear.

If $\{Q\} = AC' \cap A'C$, we have that $V, M, W$ are collinear, and if we note $\{R\} = BC' \cap B'C$ and $N$ the middle point of $(CC')$ we have that $R$, $U$, $N$ are collinear, on the other side $V, M, N$ are collinear belonging to the median from the vertex $V$ of the triangle VCC'. From these we find that $U, V, W, M, N$ are collinear. If we note $S, S'$ the middle points of the segments $(BC), (B'C')$, we have $\{G\} = A'S' \cap AS$. The quadrilateral $AA'SS'$ is a trapeze and $MG$ passes through the middle of the segment $(SS')$. Through a similar rational and by noting $T, T'$ the middle points of the segments $(AB), (A'B')$ we have that the points $N, G$ and the middle of the segment $(TT')$ are collinear. Also, if we note $X, X'$ the middle points of the segments $(AC), (A'C')$, we find that $G$, the middle point of $(BB')$ and the middle point of $(NN')$ are collinear.

Because $M, N$ and the middle of $(BB')$ are collinear, it results that $G$ belongs to the line $MN$, on which we noticed that are placed also the points $U, V, W$.

**86.**

If we note $S_a, S_b, S_c$ the first triangle of Sharygin of the triangle $ABC$, $\{A_o\} = S_bS_c \cap BC$ and $A_1$ the intersection point of the external bisectrix of the angle $A$ with $BC$, through some computations we'll find

$$A_1B = \frac{ac}{b-c}$$
$$A_1A' = \frac{2abc}{b^2-c^2}$$
$$A_oB = \frac{ac^2}{b^2-c^2}$$
$$A_oC = \frac{ab^2}{b^2-c^2}$$

From $\frac{A_oB}{A_oC} = \frac{c^2}{b^2}$ it results that $A_o$ is the intersection point of the tangent from $A$ to the circumscribed circle to the triangle $ABC$ with $BC$, therefore the foot of the exterior symmedian of the vertex $A$.

**87.**

See Pascal's theorem.

**88.**

$SDA \sim \Delta ABC$ implies

$$\frac{SD}{AB} = \frac{SE}{AC} = \frac{DE}{BC} \qquad (1)$$

From $\sphericalangle ABC \equiv \sphericalangle SDE$ it results that $DE$ is the tangent to the circumscribed circle to the triangle $BSD$, therefore, $DE, BC$ are anti-polar, and $\sphericalangle AED \equiv \sphericalangle ABC$.

It results that $SD \parallel AC$. Similarly we'll obtain $SE \parallel AB$.

$$\frac{BS}{BC} = \frac{BD}{BA} = \frac{DS}{AC} \qquad (2)$$
$$\frac{CS}{CB} = \frac{CE}{CA} = \frac{SE}{BA} \qquad (3)$$

From (2) and (3) we deduct

$$\frac{BS}{CS} = \frac{AB}{AC} \cdot \frac{DS}{SE} \qquad (4)$$

From (1) we retain

$$\frac{DS}{SE} = \frac{AB}{AC} \qquad (5)$$

The relations (4) and(5) give us the relation

$$\frac{BS}{CS}=\left(\frac{AB}{AC}\right)^2$$

**89.**

If we note

$$m\sphericalangle\left(PC_aC_c\right)=\alpha$$
$$m\sphericalangle\left(PC_aC_b\right)=\alpha'$$

we have

$$\sphericalangle AC_cA_1=\alpha,\ \sphericalangle AC_bA_1=\alpha'.$$

Then

$$\frac{\sin BAA_1}{\sin CAA_1}=\left(\frac{\sin\alpha}{\sin\alpha'}\right)^2.$$

From the concurrency of $C_aP, C_bP, C_cP$ we have:

$$\frac{\sin\alpha}{\sin\alpha'}\cdot\frac{\sin\beta}{\sin\beta'}\cdot\frac{\sin\gamma}{\sin\gamma'}=1$$

We obtain that

$$\frac{\sin BAA_1}{\sin CAA_1}\cdot\frac{\sin CBB_1}{\sin ABB_1}\cdot\frac{\sin BCC_1}{\sin ACC_1}=1$$

**90.**

Considering the inversion $i_o^{R^2}$, we observe that the image of the line $BC$ through this inversion is the circumscribed circle to the triangle $BOC$.

The circle $\mathcal{C}\left(O_1,O_1A\right)$ is the image through the same inversion of the height $AA'$. Indeed, the height being perpendicular on BC, the images of $BC$ and $AA'$ will be orthogonal circles, and the circle $\mathcal{C}\left(O_1,O_1A\right)$ has the radius $O_1O$ perpendicular on the radius that passes through $O$ of the circumscribed circle of the triangle $BOC$. It results the image of $A'$ through the considered inversion will be the intersection point $D$ of the line $OA'$ with the circle $\mathcal{C}\left(O_1,O_1A\right)$ and this coincide with the second point of intersection of the circle $\left(BOC\right)$ and $\mathcal{C}\left(O_1,O_1A\right)$.

**91.**

We'll consider $k>1$ and we'll note $A_1,B_1,C_1$ the projections of the vertexes of the triangle $ABC$ on its opposite sides.

We have

$$BA_1=c\cdot\cos B$$
$$CA_1=b\cdot\cos C$$
$$AA_1=h_a$$
$$BA'=p-b;\ \ A'A''=\left(k-1\right)r$$

Also, we'll note $D,E,F$ the intersection points of the lines $AA'',BB'',CC''$ respectively with the sides $BC,CA,AB$.

From $\Delta AA_1D \sim \Delta A''A'D$ we find $\frac{DA_1}{DA'}=\frac{h_a}{(k-1)r}$ and further more:

$$DA'=\frac{(k-1)r\cdot(p-b-\cos B)}{h_a+(k-1)r}$$

$$\frac{BD}{DC}=\frac{BA'-DA'}{CA'+DA'}=\frac{(p-b)\cdot h_a+c\cdot r(k-1)\cos B}{(p-c)\cdot h_a+a\cdot r(k-1)-c\cdot r(k-1)\cos B}$$

We substitute $h_a=\frac{2S}{a}, r=\frac{S}{p}$, and we obtain

$$\frac{BD}{DC}=\frac{2p(p-b)+(k-1)ac\cos B}{2p(p-c)+(k-1)(a^2-ac\cos B)}$$

Taking into account of

$$a^2-ac\cos B=a(a-c\cos B)=a(a-BA_1)=a\cdot CA_1=ab\cos C$$

and that

$$2p(p-b)=a^2+c^2-b^2+2ac;\ 2p(p-c)=a^2+b^2-c^2+2ab,$$

along of the cosine's theorem, we obtain

$$\frac{BD}{DC}=\frac{c(1+k\cos B)}{b(1+k\cos C)}$$

Similarly we find

$$\frac{EC}{EA}=\frac{a(1+k\cos C)}{c(1+k\cos A)},$$

$$\frac{FA}{FB}=\frac{b(1+k\cos A)}{a(1+k\cos B)}$$

We observe that

$$\frac{BD}{DC}\cdot\frac{EC}{EA}\cdot\frac{FA}{FB}=1$$

which shows that the lines $AA'',BB'',CC''$ are concurrent.

**92.**

i) $m(\sphericalangle AIB)=90^\circ+\frac{C}{2};\ m(\sphericalangle AIA_1)=90^\circ$. It results $m(\sphericalangle A_1IB)=\frac{C}{2}$

$$\frac{A_1B}{A_1C}=\frac{Aria(A_1IB)}{Aria(A_1IC)}=\frac{IA_1\cdot IB-\sin\frac{C}{2}}{IA_1\cdot IC-\sin(\sphericalangle A_1IC)}$$

But,

$$A_1IC=90^\circ+\frac{A}{2}+\frac{C}{2};\sin(\sphericalangle A_1IC)=\sin\frac{B}{2}$$

From the sinus' theorem applied in the triangle $BIC$ we have:

$$\frac{IB}{\sin\frac{C}{2}} = \frac{IC}{\sin\frac{B}{2}}$$

therefore

$$\frac{A_1B}{A_1C} = \left(\frac{\sin\frac{C}{2}}{\sin\frac{B}{2}}\right)^2.$$

Similarly

$$\frac{B_1C}{B_1A} = \left(\frac{\sin\frac{A}{2}}{\sin\frac{C}{2}}\right)^2$$

$$\frac{C_1A}{C_1B} = \left(\frac{\sin\frac{B}{2}}{\sin\frac{A}{2}}\right)^2$$

We'll then obtain

$$\frac{A_1B}{A_1C} \cdot \frac{B_1C}{B_1A} \cdot \frac{C_1A}{C_1B} = 1$$

Therefore, $A_1, B_1, C_1$ are collinear.

ii)

We note $m\left(\sphericalangle B'A'A_1'\right) = \alpha$, then we have $\sphericalangle B'A'A_1' = \sphericalangle A_2IB$ as angles with perpendicular sides.

$$\frac{A_2B}{A_2C} = \frac{\sin\alpha \cdot IB}{\sin\left(\sphericalangle BIC + \alpha\right) \cdot IC}$$

But

$$m\left(\sphericalangle BIC\right) = 90^\circ + \frac{A}{2};\ \ m\left(\sphericalangle A_2IC\right) = 90^\circ + \frac{A}{2} + \alpha$$

$$\alpha = m\left(\sphericalangle B'A'A_1'\right);$$

$$\sin\left(90^\circ + \frac{A}{2} + \alpha\right) = \sin\left(180^\circ - A_2IC\right) = \sin\left(90^\circ - \frac{A}{2} - \alpha\right)$$

We find that

$$\sin\left(\sphericalangle A_2IC\right) = \sin\left(A' - \alpha\right);\ A' = \sphericalangle B'A'C'$$

Therefore

$$\frac{A_2B}{A_2C} = \frac{\sin\alpha}{\sin\left(A' - \alpha\right)} \cdot \frac{\sin\frac{C}{2}}{\sin\frac{B}{2}}$$

If we note $m\left(\sphericalangle C'B'B_1'\right)=\beta$, $m\left(\sphericalangle A'C'C_1'\right)=\gamma$, similarly we find:

$$\frac{B_2C}{B_2A}=\frac{\sin\beta}{\sin\left(B'-\beta\right)}\cdot\frac{\sin\frac{A}{2}}{\sin\frac{C}{2}}$$

$$\frac{C_2A}{C_2B}=\frac{\sin\gamma}{\sin\left(C'-\gamma\right)}\cdot\frac{\sin\frac{B}{2}}{\sin\frac{A}{2}}$$

Observation:
This problem can be resolved by duality transformation of the "Euler's line".

**93.**

It is shown that if $AD_1$ is the external bisectrix from $A$ and $A_1H \parallel AD_1$, $A_1\in BC$, we have that $A_1H$ is the exterior bisectrix in the triangle $BHC$, therefore

$$\frac{A_1B}{A_1C}=\frac{BH}{CH}$$

**94.**

Consider the inversion $i_o^{R^2}$, where $R$ is the radius of the given circle. Through this inversion the given circle remains invariant, and the circumscribed circles to triangles $AOB$, $BOC$, $COD$, $DOE$, $EOF$, $FOA$ are transformed on the lines $AB, BC, CD, DE, EF$ of the inscribed hexagon $ABCDEF$.

The points $A_1, B_1, C_1$ have as inverse the intersection points $A_1', B_1' C_1'$ of the pairs of opposite sides of the hexagon.

Considering the Pascal's theorem, these points are collinear, and it results that the initial points $A_1, B_1, C_1$ are situated on a circle that passes through $O$.

**95.**

The quadrilateral $DCFP$ is inscribable. It result

$$\sphericalangle DCP\equiv\sphericalangle DFP \text{ and } \sphericalangle PCF\equiv\sphericalangle PDF \qquad (1)$$

On the other side $AD\perp BC$ and $\sphericalangle DCP\equiv\sphericalangle DAC$

Let $\{Q\}=BD\cap FP$. The quadrilateral $BQPE$ is inscribable, it results

$$\sphericalangle EQP\equiv\sphericalangle EBP \qquad (2)$$

From (1) and (2) we deduct that $\sphericalangle EQP\equiv\sphericalangle PDF$, therefore the quadrilateral $EQDF$ is inscribable which implies that

$$\sphericalangle FED\equiv\sphericalangle FQD \qquad (3)$$

But $FQPE$ inscribable implies that

$$\sphericalangle BEF\equiv\sphericalangle FED \qquad (4)$$

From (3) and (4) we'll retain $\sphericalangle BED \equiv \sphericalangle FED$ which shows that $(ED$ is the bisectrix of the angle $BEF$. Because $(AD$ is the interior bisectrix in the triangle $EAF$ we'll find that $D$ is the center of the circle A-ex-inscribed in the triangle $AEF$. If $DT \perp EF$ we have $DT \perp EF / DT = DC = DB, TF = FC, TE = BE$. From $EF = TF + TE$ and the above relations we'll find that $EF = BE + CF$.

**96.**

It can be proved without difficulties that $I$ is the orthocenter for the triangle $A'B'C'$, also it is known that the symmetries of the orthocenter in rapport to the triangle's sides are on the circumscribed circle to the triangle, therefore $A$ is the symmetric of the point $I$ in rapport to $B'C'$. That means that the quadrilateral $IOO_1A$ is an isosceles trapeze.

We have $IO_1 = AO$ and $\sphericalangle IO_1O \equiv \sphericalangle IAO$. On the other side $\sphericalangle IAO \equiv \sphericalangle OA'I$, we deduct that $\sphericalangle IO_1O \equiv \sphericalangle OA'I$. Having also $IA' \| OO_1$ (are perpendicular on $B'C'$) we obtain that $OO_1IA'$ is a parallelogram. It results that $OA' \| IO_1$, but $OA' \perp BC$ leads us to $O_1I$ perpendicular on *BC*.

Similarly, it results that $IO_1 = IO_2 = R$ and $IO_2 \perp AC$, $IO_3 \perp AB$, consequently the lines $AO_1, BO_2, CO_3$ are concurrent in a Kariya point of the triangle $ABC$.

**97.**

Let $M, N, P$ the middle point of the diagonals $(AC), (BD), (EF)$ of the complete quadrilateral $ABCDEF$ and $R, S, T$ respectively the middle points of the diagonals $(FB), (ED), (AG)$ of the complete quadrilateral $EFDBAG$. The Newton-Gauss line $M-N-P$ respectively $R-S-T$ are perpendicular because are diagonals in the rhomb $PRNS$.

**98.**

It is known that $AH = 2R|\cos A|$ and $OM_a = \frac{1}{2}AH$, therefore $OA' = kR|\cos A|$. We'll note $\{P\} = OH \cap AA'$ ($H$ is the orthocenter of the triangle $ABC$). From the similarity of the triangles $APH, A'PO$ it results that

$$\frac{AH}{A'O} = \frac{HP}{OP} = \frac{2}{k} = const.$$

Consequently, $P$ is a fixed point on $OH$. Similarly it will be shown that $BB', CC'$ pass through $P$.

**99.**

Because the points $B_1, B_2, C_1, C_2$ are concyclic and $AC_1 \cdot AC_2 = AB_1 \cdot AB_2$, it results that the point $A$ belongs to the radical axis of the circles with the centers in $H$ and $P$. The radical axis being perpendicular on the line of the centers $H$ and $P$, which is parallel to $BC$ will lead us to the fact that the radical axis is the height from the point $A$.

A similar rational shows that the height from the vertex $B$ of the triangle $ABC$ is a radical axis for the circles of centers $M, P$. The intersection of the radical axes, which is the orthocenter of the triangle $ABC$ is the radical center of the constructed circles.

**100.**

Let $A_1$ the intersection of the tangent in $M$ to the circumscribed circle to the triangle $BMC$ with $BC$. Because $MA_1$ is exterior symmedian in the triangle $BMC$, we have

$$\frac{A_1C}{A_1B} = \left(\frac{MC}{MB}\right)^2.$$

Similarly, we note $B_1, C_1$ the intersection points of the tangents in $M$ to the circles $CMA, AMB$ with $AC$ respectively with $AB$.

We find

$$\frac{B_1A}{B_1C} = \left(\frac{MA}{MC}\right)^2 \text{ and } \frac{C_1B}{C_1A} = \left(\frac{MB}{MA}\right)^2$$

The reciprocal Menelaus' theorem leads us to the collinearity of the points $A_1, B_1, C_1$.

**Online papers**:

Comments: 13 pages, 10 geometrical figures
Subjects: General Mathematics (math.GM)